\documentclass[12pt]{amsbook}
\usepackage{amssymb}
\usepackage[all]{xy}

\usepackage[colorlinks=true,linkcolor=magenta,citecolor=magenta]{hyperref}
\makeindex

\newtheorem{theorem}{Theorem}[chapter]
\newtheorem{answer}[theorem]{Answer}
\newtheorem{cat}[theorem]{Cat}
\newtheorem{claim}[theorem]{Claim}
\newtheorem{comment}[theorem]{Comment}
\newtheorem{conclusion}[theorem]{Conclusion}
\newtheorem{convention}[theorem]{Convention}
\newtheorem{definition}[theorem]{Definition}
\newtheorem{example}[theorem]{Example}
\newtheorem{exercise}[theorem]{Exercise}
\newtheorem{fact}[theorem]{Fact}
\newtheorem{lemma}[theorem]{Lemma}
\newtheorem{plan}[theorem]{Plan}
\newtheorem{principle}[theorem]{Principle}
\newtheorem{problem}[theorem]{Problem}
\newtheorem{proposition}[theorem]{Proposition}
\newtheorem{question}[theorem]{Question}
\newtheorem{questions}[theorem]{Questions}
\newtheorem{rat}[theorem]{Rat}
\newtheorem{recipe}[theorem]{Recipe}
\newtheorem{remark}[theorem]{Remark}
\newtheorem{thought}[theorem]{Thought}

\begin{document}

\title{Continuous random variables}

\author{Teo Banica}
\address{Department of Mathematics, University of Cergy-Pontoise, F-95000 Cergy-Pontoise, France. {\tt teo.banica@gmail.com}}


\begin{abstract}
This is an introduction to probability, written with a quantum idea in mind, namely that the semicircle law comes first. We first discuss discrete probability, notably with the binomial and hypergeometric laws, positive and negative, and the Poisson and compound Poisson laws. Then we get into the continuous case, with the basics of the theory explained, and with as starting examples the exponential, semicircle and beta distributions. Afterwards, we investigate the central limits and normal variables, both real and complex, and with a look at Rayleigh laws, and hyperspherical laws too. Finally, we discuss a number of more specialized distributions, and more specialized techniques too, and we end with an introduction to random matrices and freeness.
\end{abstract}

\maketitle

\chapter*{Preface}

One of the fascinating aspects of probability theory, that you can enjoy both as a beginner, or as someone initiated, is the variety of probability laws available. There are so many of them, all interesting, and with relations between them, and with their mathematics sharing many common features. Put a bit poetically, while disciplines like basic algebra or analysis or physics or chemistry look very structured, say as mountains, when contemplating basic probability, you rather feel being in front of an ocean.

\bigskip

Which is certainly relaxing, and no wonder we are so many, to love probability. This being said, for basic learning purposes, and for some peace of mind later too, a bit of hierarchy between that various ocean beasts can help you a lot, in your quest.

\bigskip

So, what comes first? Here, based on both teaching and research beliefs, the general agreement is that the discrete laws come first, and in a quite precise order, namely Bernoulli, then binomial, then negative binomial including geometric, and then, upon taste, either hypergeometric, positive and negative, or Poisson and compound Poisson. Followed of course by the many other interesting discrete laws which are available.

\bigskip

However, passed this, which remains close to the shore, here we are now in the deep blue ocean, meaning continuous probability, the land of whales and sharks. And there, passed the uniform and exponential laws, which again, based on both teaching and research beliefs, come first, we have two big beasts fighting for supremacy, namely:
$$g_1=\frac{1}{\sqrt{2\pi}}\,e^{-x^2/2}dx\quad,\quad\gamma_1=\frac{1}{2\pi}\sqrt{4-x^2}\,dx$$

To be more precise, on the left we have the ubiquitous normal law of Gauss, coming via central limits, which are something fundamental. In a word, big beast that we have here. However, this creature has an Achilles heel, coming from the formula $\int_\mathbb Re^{-x^2}dx=\sqrt{\pi}$, which is something quite complicated, needed for its definition. So, this is the situation, and personally, I always think at $g_1$ as being some kind of a big whale.

\bigskip

As for the beast on the right, that is the semicircle law of Wigner, coming from nowhere, I mean the mass 1 property comes from the fact that the unit circle has area $\pi$. This law is central in quantum physics and related mathematical areas, such as operator algebras, random matrices and quantum groups. As for the weaknesses, none I would say, quantum being supposedly at the origin of everything. This $\gamma_1$ is a true shark.

\bigskip

So, this was for the philosophy, or at least for the philosophy that I personally believe in. Getting now to the present book, this will be an introduction to probability theory, mainly interested in the continuous case, and based on the above belief, that the semicircle law, being quantum, comes first. The book is organized in 4 parts, as follows:

\bigskip

I. We first discuss discrete probability, notably with the binomial and hypergeometric laws, positive and negative, and the Poisson and compound Poisson laws. 

\bigskip

II. Then we get into the continuous case, with the basics of the theory explained, and with as starting examples the exponential, semicircle and beta distributions.

\bigskip

III. Afterwards, we investigate the central limits and normal variables, both real and complex, and with a look at Rayleigh laws, and hyperspherical laws too.

\bigskip

IV. Finally, we discuss more specialized distributions, and more specialized techniques too, and we end with an introduction to random matrices and freeness.

\bigskip

In the hope that you will find this book useful, whether for learning the basics, or fine-tuning some previous knowledge. Of course, the philosophy behind it remains a bit speculative, and I would have loved to have a theorem stating that this whole world, including $g_1$, comes from some minuscule elementary particles, obeying to $\gamma_1$. But this is something that we don't have yet, still far from that, you will have to trust me here.

\bigskip

This book, although being not exactly beginner level, is quite self-contained, and normally no prior reading of anything is needed. However, as some references, for more on discrete probability you have my book \cite{ba1}, for more regarding the material at the end you have my book \cite{ba2}, and for some help with the physics you have my book \cite{ba3}.

\bigskip

Finally, writing a probability book is always a pleasure, save sometimes for certain axiomatic aspects, and I would like to thank here my cats, for some help with this.

\bigskip

\

{\em Cergy, August 2026}

\smallskip

{\em Teo Banica}

\baselineskip=15.95pt
\tableofcontents
\baselineskip=14pt

\part{Basic probability}

\ \vskip50mm

\begin{center}
{\em Bow down mister

Hare Rama, Hare Krishna

Bow down mister

We say Radhe Shyam}
\end{center}

\chapter{Counting, estimates}

\section*{1a. Binomial formula}

Welcome to probability, and as good news, we will start this book, the first 100 pages, with a study in the discrete case, meant to be introductory and relaxing. Expect lots of interesting computations featuring coins, dice, cards and the like, which are certainly the best way, in order to get introduced to probability. As for the continuous case, which is more complicated, we will have the remainder of the book, 300 pages, for that.

\bigskip

As bad news, however, the discrete probability that we will be doing will always lead us into factorials $n!$, which are, after all, quite abstract numbers. In order to convert our theorems into useful findings, we will need the Stirling formula, stating that:
$$n!\simeq\left(\frac{n}{e}\right)^n\sqrt{2\pi n}$$

And the point is that this formula, featuring both numbers $e,\pi$, certainly means good knowledge in real analysis. In fact, worse is true, because in order to prove this formula, we will need in fact to learn full calculus, meaning derivatives and integrals, and with a passage into 2 dimensions too, for proving Gauss, $\int_\mathbb Re^{-x^2}=\sqrt{\pi}$, used by Stirling.

\bigskip

Well, so this is the situation, we are right away into a logical loop, doing discrete probability, in order to avoid the continuous case, naturally drags us into calculus, via the Stirling formula, which calculus itself means more or less continuous probability.

\bigskip

Nevermind. Stopping complaining now, here is a plan for the present book:

\begin{plan}
We will talk in this first chapter about counting and estimates,
\begin{enumerate}
\item Factorials and binomials, and their basic applications,

\item Sequences and series, and related estimates for factorials and binomials,

\item Derivatives and integrals, and better estimates for factorials and binomials,

\item Stirling formula above, and best estimates for factorials and binomials,
\end{enumerate}
and afterwards, using this knowledge, we will do probability, discrete, then continuous.
\end{plan}

Getting started now, at the beginning of everything in probability, and in mathematics in general, we have the binomial formula, that you surely know well, as follows:

\index{binomial coefficient}
\index{factorials}

\begin{theorem}
The number of possibilities of choosing $k$ objects among $n$ objects is
$$\binom{n}{k}=\frac{n!}{k!(n-k)!}$$
called binomial number, where $n!=1\cdot2\cdot3\ldots(n-2)(n-1)n$, called ``factorial $n$''.
\end{theorem}

\begin{proof}
This is something very standard, the idea being as follows:

\medskip

(1) Imagine a set consisting of $n$ objects. We have then $n$ possibilities for choosing our 1st object, then $n-1$ possibilities for choosing our 2nd object, out of the $n-1$ objects left, and so on up to $n-k+1$ possibilities for choosing our $k$-th object, out of the $n-k+1$ objects left. Since the possibilities multiply, the total number of choices is:
$$N=n(n-1)\ldots(n-k+1)=\frac{n!}{(n-k)!}$$

(2) However, thiking well, the number $N$ that we computed is in fact the number of possibilities of choosing $k$ ordered objects among $n$ objects. Thus, we must divide everything by the number $M$ of orderings of the $k$ objects that we chose:
$$\binom{n}{k}=\frac{N}{M}$$

(3) In order to compute now the missing number $M$, imagine a set consisting of $k$ objects. There are $k$ choices for the object to be designated $\#1$, then $k-1$ choices for the object to be designated $\#2$, and so on up to 1 choice for the object to be designated $\#k$. We conclude that the missing number $M$ is given by the following formula:
$$M=k(k-1)\ldots 2\cdot 1=k!$$

(4) Now by getting back to what we wanted to compute, we are led to:
$$\binom{n}{k}
=\frac{N}{M}
=\frac{n!/(n-k)!}{k!}
=\frac{n!}{k!(n-k)!}$$

(5) Finally, for our theory to be complete, we must set $0!=1$, as for the following computation of $\binom{n}{n}=1$, based on the formula that we found, to work fine:
$$\binom{n}{n}=\frac{n!}{n!0!}=\frac{n!}{n!\times1}=1$$

And with this, and of what was to be said, on this subject.
\end{proof}

Coming next, we have the following result, that you surely know too:

\index{binomial formula}

\begin{theorem}
We have the binomial formula
$$(a+b)^n=\sum_{k=0}^n\binom{n}{k}a^kb^{n-k}$$
valid for any two numbers $a,b\in\mathbb R$. 
\end{theorem}

\begin{proof}
The question to be solved is that of computing the following quantity:
$$(a+b)^n=\underbrace{(a+b)(a+b)\ldots(a+b)}_{n\ terms}$$

When expanding, we obtain a certain sum of products of $a,b$ variables, with each such product being a quantity of type $a^kb^{n-k}$. Thus, we have a formula as follows:
$$(a+b)^n=\sum_{k=0}^nC_ka^kb^{n-k}$$

But, according to our product formula above, $C_k$ is the number of choices for the $k$ needed $a$ variables among the $n$ available $a$ variables. Thus $C_k=\binom{n}{k}$, as desired.
\end{proof}

Theorem 1.3 is something quite interesting, so let us doublecheck it with some numerics. At small values of $n$ we obtain the following formulae, which are all correct:
$$(a+b)^0=1$$
$$(a+b)^1=a+b$$
$$(a+b)^2=a^2+2ab+b^2$$
$$(a+b)^3=a^3+3a^2b+3ab^2+b^3$$
$$(a+b)^4=a^4+4a^3b+6a^2b^2+4ab^3+b^4$$
$$(a+b)^5=a^5+5a^4b+10a^3b^2+10a^2b^3+5ab^4+b^5$$
$$\vdots$$

Now observe that in these formulae, say for memorization purposes, the powers of the $a,b$ variables are something very simple, that can be recovered right away. What matters are the coefficients, which are the binomial coefficients $\binom{n}{k}$, which form a triangle. So, it is enough to memorize this triangle, and this can be done by using:

\index{Pascal triangle}
\index{binomial coefficients}

\begin{theorem}
The Pascal triangle, formed by the binomial coefficients $\binom{n}{k}$,
$$1$$
$$1\ \ ,\ \ 1$$
$$1\ \ ,\ \ 2\ \ ,\ \ 1$$
$$1\ \ ,\ \ 3\ \ ,\ \ 3\ \ ,\ \ 1$$
$$1\ \ ,\ \ 4\ \ ,\ \ 6\ \ ,\ \ 4\ \ ,\ \ 1$$
$$1\ \ ,\ \ 5\ \ ,\ \ 10\ \ ,\ \ 10\ \ ,\ \ 5\ \ ,\ \ 1$$
$$\vdots$$
has the property that each entry is the sum of the two entries above it.
\end{theorem}

\begin{proof}
In practice, the theorem states that the following formula holds:
$$\binom{n}{k}=\binom{n-1}{k-1}+\binom{n-1}{k}$$

There are many ways of proving this formula, all instructive, as follows:

\medskip

(1) Brute-force computation. We have indeed, as desired:
\begin{eqnarray*}
\binom{n-1}{k-1}+\binom{n-1}{k}
&=&\frac{(n-1)!}{(k-1)!(n-k)!}+\frac{(n-1)!}{k!(n-k-1)!}\\
&=&\frac{(n-1)!}{(k-1)!(n-k-1)!}\left(\frac{1}{n-k}+\frac{1}{k}\right)\\
&=&\frac{(n-1)!}{(k-1)!(n-k-1)!}\cdot\frac{n}{k(n-k)}\\
&=&\binom{n}{k}
\end{eqnarray*}

(2) Algebraic proof. We have the following formula, to start with:
$$(a+b)^n=(a+b)^{n-1}(a+b)$$

By using the binomial formula, this formula becomes:
$$\sum_{k=0}^n\binom{n}{k}a^kb^{n-k}=\left[\sum_{r=0}^{n-1}\binom{n-1}{r}a^rb^{n-1-r}\right](a+b)$$

Now let us perform the multiplication on the right. We obtain a certain sum of terms of type $a^kb^{n-k}$, and to be more precise, each such $a^kb^{n-k}$ term can either come from the $\binom{n-1}{k-1}$ terms $a^{k-1}b^{n-k}$ multiplied by $a$, or from the $\binom{n-1}{k}$ terms $a^kb^{n-1-k}$ multiplied by $b$. Thus, the coefficient of $a^kb^{n-k}$ on the right is $\binom{n-1}{k-1}+\binom{n-1}{k}$, as desired.

\medskip

(3) Combinatorics. Let us count $k$ objects among $n$ objects, with one of the $n$ objects having a hat on top. Obviously, the hat has nothing to do with the count, and we obtain $\binom{n}{k}$. On the other hand, we can say that there are two possibilities. Either the object with hat is counted, and we have $\binom{n-1}{k-1}$ possibilities here, or the object with hat is not counted, and we have $\binom{n-1}{k}$ possibilities here. Thus $\binom{n}{k}=\binom{n-1}{k-1}+\binom{n-1}{k}$, as desired.
\end{proof}

In practice, in order to get some applications out of our counting results, we will need some estimates for the binomials. And here, we have the following basic result:

\begin{theorem}
We have the following estimates for the binomial coefficients,
$$\left(\frac{n}{k}\right)^k\leq\binom{n}{k}\leq\frac{n^k}{k!}$$
both coming from definitions.
\end{theorem}

\begin{proof}
Regarding the first estimate, this can be established as follows:
$$\binom{n}{k}
=\frac{n}{k}\cdot\frac{n-1}{k-1}\ldots\frac{n-k+1}{1}
\geq\frac{n}{k}\cdot\frac{n}{k}\ldots\frac{n}{k}
=\left(\frac{n}{k}\right)^k$$

The proof of the second estimate is similar, as follows:
$$\binom{n}{k}
=\frac{n}{k}\cdot\frac{n-1}{k-1}\ldots\frac{n-k+1}{1}
\leq\frac{n}{k}\cdot\frac{n}{k-1}\ldots\frac{n}{1}
=\frac{n^k}{k!}$$

Thus, we are led to the conclusions in the statement.
\end{proof}

The problem is now, how to do better? Leaving aside the fact that both our estimates above were quite rough, so job for us later to improve our method, what is quite bothering is the upper bound $n^k/k!$, and more specifically, the number $k!$ appearing there. Thus, estimating the factorials $k!$, as sharp as possible, will be our next task.

\section*{1b. Series, estimates}

In order to solve our approximation questions for factorials, the idea will be that of using the famous number $e=2.71828\ldots$ We will need the following basic fact:

\begin{proposition}
We have the following inequality, for any $a_1,\ldots,a_n\geq0$,
$$\frac{a_1+\ldots+a_n}{n}\geq\sqrt[n]{a_1\ldots a_n}$$
telling us that the arithmetic mean is bigger than the geometric mean.
\end{proposition}

\begin{proof}
To start with, the result holds indeed at $n=2$, with this coming from:
$$\frac{a+b}{2}\geq\sqrt{ab}\iff(\sqrt{a}-\sqrt{b})^2\geq0$$

But with this, we can prove our inequality at $n=4$ too, as follows:
$$\frac{a+b+c+d}{4}
\geq\frac{\sqrt{ab}+\sqrt{cd}}{2}
\geq\sqrt[4]{abcd}$$

And so on, the conclusion being that we have the result in the case $n=2^s$. In general now, given numbers $a_1,\ldots,a_n\geq0$, consider their arithmetic mean:
$$m=\frac{a_1+\ldots+a_n}{n}$$

Now pick $s\in\mathbb N$ such that $n\leq 2^s$, and let us complete our series $a_1,\ldots,a_n$ with $2^s-n$ copies of $m$. The arithmetic mean stays the same, and using the result at $2^s$, we get:
\begin{eqnarray*}
m\geq\sqrt[2^s]{a_1\ldots a_nm^{2^s-n}}
&\implies&m^{2^s}\geq a_1\ldots a_nm^{2^s-n}\\
&\implies&m^n\geq a_1\ldots a_n\\
&\implies&m\geq\sqrt[n]{a_1\ldots a_n}
\end{eqnarray*}

Thus, we are led to the conclusion in the statement. 
\end{proof}

Good news, we can talk now about a very interesting convergence, as follows:

\index{e}

\begin{theorem}
We have the following convergence
$$\left(1+\frac{1}{n}\right)^n\to e$$
where $e=2.71828\ldots$ is a certain number.
\end{theorem}

\begin{proof}
This is something quite tricky, as follows:

\medskip

(1) Our first claim is that the following sequence is increasing:
$$x_n=\left(1+\frac{1}{n}\right)^n$$ 

In order to prove this, we use the following arithmetic-geometric inequality:
$$\frac{1+\sum_{i=1}^n\left(1+\frac{1}{n}\right)}{n+1}\geq\sqrt[n+1]{1\cdot\prod_{i=1}^n\left(1+\frac{1}{n}\right)}$$

In practice, by raising to the power $n+1$ we obtain, as desired:
$$\left(1+\frac{1}{n+1}\right)^{n+1}\geq\left(1+\frac{1}{n}\right)^n$$

(2) Normally we are left with proving that $x_n$ is bounded from above, but this is non-trivial, and we have to use a trick. Consider the following sequence:
$$y_n=\left(1+\frac{1}{n}\right)^{n+1}$$ 

We will prove that this sequence $y_n$ is decreasing, and together with the fact that we have $x_n/y_n\to1$, this will give the result. So, this will be our plan. 

\medskip

(3) In order to prove now that $y_n$ is decreasing, we use, a bit as before:
$$\frac{1+\sum_{i=1}^n\left(1-\frac{1}{n}\right)}{n+1}\geq\sqrt[n+1]{1\cdot\prod_{i=1}^n\left(1-\frac{1}{n}\right)}$$

In practice, by raising to the power $n+1$, we obtain from this:
$$\left(1-\frac{1}{n+1}\right)^{n+1}\geq\left(1-\frac{1}{n}\right)^n$$

And, by inverting this inequality that we found, we obtain, as desired:
$$\left(1+\frac{1}{n}\right)^{n+1}\leq\left(1+\frac{1}{n-1}\right)^n$$

(4) But with this, we can now finish. Indeed, the sequence $x_n$ is increasing, the sequence $y_n$ is decreasing, and we have $x_n<y_n$, as well as:
$$\frac{y_n}{x_n}=1+\frac{1}{n}\to1$$

Thus, both sequences $x_n,y_n$ converge to a certain number $e$. Finally, regarding the numerics for $e$, these follow from $x_n<e<y_n$, with a bit of patience, or a computer.
\end{proof}

More generally now, we have the following result, regarding the exponential:

\index{e}
\index{exponential}

\begin{theorem}
We have the following formula,
$$\left(1+\frac{x}{n}\right)^n\to e^x$$
valid for any $x\in\mathbb R$.
\end{theorem}

\begin{proof}
We know that this holds at $x=1$, by definition of $e$, and by inverting, we have it at $x=-1$ too. But then, when $x\in\mathbb R$ is arbitrary, we can proceed as follows:
$$\left(1+\frac{x}{n}\right)^n=\left[\left(1+\frac{x}{n}\right)^{n/x}\right]^x
\to e^x$$

Thus, we are led to the conclusion in the statement.
\end{proof}

All this is very nice, but we have as well an alternative approach to $e$, as follows:

\begin{theorem}
We have the following formula,
$$e=\sum_{k=0}^\infty\frac{1}{k!}$$
which can stand as an alternative definition for $e=2.71828\ldots$
\end{theorem}

\begin{proof}
This is something very standard, the idea being as follows:

\medskip

(1) In practice, we want to prove that we have the following equality:
$$\sum_{k=0}^\infty\frac{1}{k!}=\lim_{n\to\infty}\left(1+\frac{1}{n}\right)^n$$

For this purpose, the first observation is that we have the following estimate:
$$2<\sum_{k=0}^\infty\frac{1}{k!}<\sum_{k=0}^\infty\frac{1}{2^{k-1}}=3$$

In order to prove now that this limit is indeed $e$, observe that we have:
$$\left(1+\frac{1}{n}\right)^n
=\sum_{k=0}^n\binom{n}{k}\cdot\frac{1}{n^k}
\leq\sum_{k=0}^n\frac{1}{k!}$$

Thus, with $n\to\infty$, we get that the limit of the series $\sum_{k=0}^\infty\frac{1}{k!}$ belongs to $[e,3)$.

\medskip

(2) For the reverse inequality, we have the following computation:
\begin{eqnarray*}
\sum_{k=0}^n\frac{1}{k!}-\left(1+\frac{1}{n}\right)^n
&=&\sum_{k=2}^n\frac{n^k-n(n-1)\ldots(n-k+1)}{n^kk!}\\
&\leq&\sum_{k=2}^n\frac{n^k-(n-k)^k}{n^kk!}\\
&=&\sum_{k=2}^n\frac{1-\left(1-\frac{k}{n}\right)^k}{k!}
\end{eqnarray*}

Next, we can use the following trivial inequality, valid for any number $x\in(0,1)$:
$$1-x^k=(1-x)(1+x+x^2+\ldots+x^{k-1})\leq(1-x)k$$

Indeed, we can use this with $x=1-k/n$, and we obtain in this way:
$$\sum_{k=0}^n\frac{1}{k!}-\left(1+\frac{1}{n}\right)^n
\leq\sum_{k=2}^n\frac{\frac{k}{n}\cdot k}{k!}
\leq\frac{1}{n}\sum_{k=2}^n\frac{2}{2^{k-2}}
<\frac{4}{n}$$

Thus, we have our needed estimate, and this finishing the proof.
\end{proof}

As a last result now about $e$, which is something quite far-reaching, we have:

\begin{theorem}
We have the following formula,
$$e^x=\sum_{k=0}^\infty\frac{x^k}{k!}$$
valid for any $x\in\mathbb R$.
\end{theorem}

\begin{proof}
Let $f(x)=\sum_kx^k/k!$ the series in the statement. We have:
\begin{eqnarray*}
f(x+y)
&=&\sum_{k=0}^\infty\frac{(x+y)^k}{k!}\\
&=&\sum_{k=0}^\infty\sum_{s=0}^k\binom{k}{s}\cdot\frac{x^sy^{k-s}}{k!}\\
&=&\sum_{k=0}^\infty\sum_{s=0}^k\frac{x^sy^{k-s}}{s!(k-s)!}\\
&=&f(x)f(y)
\end{eqnarray*}

Next, observe that this shows that $f$ is continuous. Indeed, at $x=0$ we have:
$$\lim_{t\to0}f(t)=\lim_{t\to0}\left(1+t\sum_{k=1}^\infty\frac{t^{k-1}}{k!}\right)=1$$

But from this, we get $f(x+t)=f(x)f(t)\to f(x)$ with $t\to0$, at any $x$. Thus, as a conclusion, our function $f$ is continuous, and satisfies the following conditions:
$$f(x+y)=f(x)f(y)\quad,\quad f(1)=e$$

But with this, we can finish. Indeed, by iterating, we have $f(nx)=f(x)^n$ for any $n\in\mathbb N$. Then, by extracting roots, we have $f(rx)=f(x)^r$ for any $r\in\mathbb Q$. Thus $f(r)=e^r$ for any $r\in\mathbb Q$, and by continuity we obtain $f(x)=e^x$ for any $x\in\mathbb R$, as desired.
\end{proof}

Good news, we can go back now to our previous estimates from Theorem 1.5, for the binomials, and improve them into something truly useful, as follows:

\begin{theorem}
We have the following estimates for the binomials,
$$\left(\frac{n}{k}\right)^k\leq\binom{n}{k}<\left(\frac{ne}{k}\right)^k$$
with $e=2.71828\ldots$ being the usual constant from analysis.
\end{theorem}

\begin{proof}
We have indeed the following estimate, based on Theorem 1.10:
$$\frac{k^k}{k!}<\sum_{n=0}^\infty\frac{k^n}{n!}=e^k$$

But with this, we can go back to the upper bound in Theorem 1.5, and we have:
$$\frac{n^k}{k!}<\frac{n^k}{k^k/e^k}=\left(\frac{ne}{k}\right)^k$$

Thus, we are led to the conclusions in the statement.
\end{proof}

\section*{1c. Derivatives, integrals}

The problem is now, how to improve Theorem 1.11? In order to solve this question, which is something non-trivial, we need to learn some calculus. Let us start with:

\index{derivative}

\begin{definition}
A function $f:\mathbb R\to\mathbb R$ is called differentiable at $x$ when
$$f'(x)=\lim_{t\to0}\frac{f(x+t)-f(x)}{t}$$
called derivative of $f$ at that point $x$, exists.
\end{definition}

As a first remark, in order for $f$ to be differentiable at $x$, that is to say, in order for the above limit to converge, the numerator must go to $0$, as the denominator $t$ does:
$$\lim_{t\to0}\left[f(x+t)-f(x)\right]=0$$

Thus, $f$ must be continuous at $x$. However, the converse is not true, a basic counterexample being $f(x)=|x|$ at $x=0$. At the level of the general theory, we have:

\index{derivative}
\index{locally affine}

\begin{proposition}
Assuming that $f$ is differentiable at $x$, we have:
$$f(x+t)\simeq f(x)+f'(x)t$$
In other words, $f$ is, approximately, locally affine at $x$.
\end{proposition}

\begin{proof}
This is indeed something self-explanatory, appearing as a reformulation of Definition 1.12, by multiplying the formula there by $t$, and concluding.
\end{proof}

As a first illustration now, the derivatives of the power functions are as follows:

\begin{proposition}
We have the differentiation formula
$$(x^p)'=px^{p-1}$$
valid for any exponent $p\in\mathbb R$.
\end{proposition}

\begin{proof}
In the case $p\in\mathbb N$ we can use the binomial formula, which gives:
$$(x+t)^p
=x^p+px^{p-1}t+\ldots+t^p
\simeq x^p+px^{p-1}t$$

Next, for $p\in\mathbb Q$, we can write $p=m/n$, with $m\in\mathbb N$ and $n\in\mathbb Z$, and we have: 
\begin{eqnarray*}
(x+t)^{m/n}-x^{m/n}
&=&\frac{(x+t)^m-x^m}{(x+t)^{m(n-1)/n}+\ldots+x^{m(n-1)/n}}\\
&\simeq&\frac{mx^{m-1}t}{nx^{m(n-1)/n}}\\
&=&\frac{m}{n}\cdot x^{m/n-1}\cdot t
\end{eqnarray*}

But then, the general case, $p\in\mathbb R$, follows too, via a continuity argument.
\end{proof}

Here are some further computations, for other basic functions that we know:

\index{sin}
\index{cos}
\index{exp}
\index{log}

\begin{theorem}
We have the following results:
\begin{enumerate}
\item $(\sin x)'=\cos x$.

\item $(\cos x)'=-\sin x$.

\item $(e^x)'=e^x$.

\item $(\log x)'=x^{-1}$.
\end{enumerate}
\end{theorem}

\begin{proof}
This is quite tricky, as always when computing derivatives, as follows:

\medskip

(1) Regarding $\sin$, the computation here goes as follows:
\begin{eqnarray*}
(\sin x)'
&=&\lim_{t\to0}\frac{\sin(x+t)-\sin x}{t}\\
&=&\lim_{t\to0}\frac{\sin x\cos t+\cos x\sin t-\sin x}{t}\\
&=&\lim_{t\to0}\sin x\cdot\frac{\cos t-1}{t}+\cos x\cdot\frac{\sin t}{t}\\
&=&\cos x
\end{eqnarray*}

(2) The computation for $\cos$ is similar, as follows:
\begin{eqnarray*}
(\cos x)'
&=&\lim_{t\to0}\frac{\cos(x+t)-\cos x}{t}\\
&=&\lim_{t\to0}\frac{\cos x\cos t-\sin x\sin t-\cos x}{t}\\
&=&\lim_{t\to0}\cos x\cdot\frac{\cos t-1}{t}-\sin x\cdot\frac{\sin t}{t}\\
&=&-\sin x
\end{eqnarray*}

(3) For the exponential, the derivative can be computed as follows:
$$(e^x)'
=\left(\sum_{k=0}^\infty\frac{x^k}{k!}\right)'
=\sum_{k=0}^\infty\frac{kx^{k-1}}{k!}
=e^x$$

(4) As for the logarithm, the computation here is as follows, using $\log(1+y)\simeq y$ for $y\simeq 0$, which follows from $e^y\simeq 1+y$ that we found in (3), by taking the logarithm:
\begin{eqnarray*}
(\log x)'
&=&\lim_{t\to0}\frac{\log(x+t)-\log x}{t}\\
&=&\lim_{t\to0}\frac{\log(1+t/x)}{t}\\
&=&\frac{1}{x}
\end{eqnarray*}

Thus, we are led to the formulae in the statement.
\end{proof}

Back now to general theory, we have the following very useful result:

\index{Leinbitz rule}
\index{chain rule}

\begin{theorem}
The derivatives are subject to the following rules:
\begin{enumerate}
\item Leibnitz rule: $(fg)'=f'g+fg'$.

\item Chain rule: $(f\circ g)'=f'(g)g'$.
\end{enumerate}
\end{theorem}

\begin{proof}
Both formulae follow from the definition of the derivative, as follows:

\medskip

(1) Regarding the products, we have the following computation:
\begin{eqnarray*}
(fg)(x+t)
&=&f(x+t)g(x+t)\\
&\simeq&(f(x)+f'(x)t)(g(x)+g'(x)t)\\
&\simeq&f(x)g(x)+(f'(x)g(x)+f(x)g'(x))t
\end{eqnarray*}

(2) Regarding compositions, we have the following computation:
\begin{eqnarray*}
(f\circ g)(x+t)
&=&f(g(x+t))\\
&\simeq&f(g(x)+g'(x)t)\\
&\simeq&f(g(x))+f'(g(x))g'(x)t
\end{eqnarray*}

Thus, we are led to the conclusions in the statement.
\end{proof}

There are many applications of the derivative, summarized as follows:

\index{local minimum}
\index{local maximum}
\index{Rolle theorem}
\index{mean value theorem}
\index{vanishing derivative}

\begin{theorem}
Given a differentiable function $f:[a,b]\to\mathbb R$, we have:
\begin{enumerate}
\item The local minima and maxima of $f$ appear at the points where $f'(x)=0$.

\item Rolle theorem: if $f(a)=f(b)$, we must have $f'(c)=0$, for some $c\in(a,b)$. 

\item Mean value theorem: $\frac{f(b)-f(a)}{b-a}=f'(c)$, for some $c\in(a,b)$.

\item Main theorem: if $f'=0$ then $f$ must be constant.
\end{enumerate}
\end{theorem}

\begin{proof}
This is something very standard, the idea being as follows:

\medskip

(1) This is clear from the formula $f(x+t)\simeq f(x)+f'(x)t$, from Proposition 1.13.

\medskip

(2) This is something clear too, coming from (1), in the obvious way.

\medskip

(3) This comes indeed from (2), applied to the following function:
$$g(x)=f(x)-\frac{f(b)-f(a)}{b-a}\cdot x$$

(4) This assertion, which is extremely useful in practice, comes from (3).
\end{proof}

At a more advanced level now, we can talk about second derivatives, and we have:

\index{Taylor formula}

\begin{theorem}
Any twice differentiable $f:\mathbb R\to\mathbb R$ is approximately locally quadratic,
$$f(x+t)\simeq f(x)+f'(x)t+\frac{f''(x)}{2}\,t^2$$
with $f''(x)$ being the derivative of the function $f':\mathbb R\to\mathbb R$ at the point $x$.
\end{theorem}

\begin{proof}
This is something quite intuitive, when thinking geometrically. In practice, we can use L'H\^opital's rule, stating that the $0/0$ type limits can be computed as:
$$\frac{f(x)}{g(x)}\simeq\frac{f'(x)}{g'(x)}$$

Observe that this formula holds indeed, as an application of Proposition 1.13. Now by using this, if we denote by $\varphi(t)\simeq P(t)$ the formula to be proved, we have:
\begin{eqnarray*}
\frac{\varphi(t)-P(t)}{t^2}
&\simeq&\frac{\varphi'(t)-P'(t)}{2t}\\
&\simeq&\frac{\varphi''(t)-P''(t)}{2}\\
&=&\frac{f''(x)-f''(x)}{2}\\
&=&0
\end{eqnarray*}

Thus, we are led to the conclusion in the statement.
\end{proof}

We can further develop the above method, at order 3, at order 4, and so on, the ultimate result on the subject, called Taylor formula, being as follows:

\index{Taylor formula}

\begin{theorem}
Assuming that $f:\mathbb R\to\mathbb R$ is $n$ times differentiable, we have
$$f(x+t)\simeq\sum_{k=0}^n\frac{f^{(k)}(x)}{k!}\,t^k$$
where $f^{(k)}(x)$ are the higher derivatives of $f$ at the point $x$.
\end{theorem}

\begin{proof}
We use the same method as in the proof of Theorem 1.18. Indeed, if we denote by $\varphi(t)\simeq P(t)$ the approximation to be proved, we have:
\begin{eqnarray*}
\frac{\varphi(t)-P(t)}{t^n}
&\simeq&\frac{\varphi'(t)-P'(t)}{nt^{n-1}}\\
&\simeq&\frac{\varphi''(t)-P''(t)}{n(n-1)t^{n-2}}\\
&\vdots&\\
&\simeq&\frac{\varphi^{(n)}(t)-P^{(n)}(t)}{n!}\\
&=&0
\end{eqnarray*}

Thus, we are led to the conclusion in the statement.
\end{proof}

As a basic application of the derivatives and Taylor formula, we have:

\index{binomial formula}
\index{generalized binomial formula}

\begin{theorem}
We have the generalized binomial formula
$$(1+t)^p=\sum_{k=0}^\infty\binom{p}{k}t^k$$
with the generalized binomial coefficients being given by
$$\binom{p}{k}=\frac{p(p-1)\ldots(p-k+1)}{k!}$$
for any $p\in\mathbb R$, and any $|t|<1$. With $p\in\mathbb N$, we recover the usual binomial formula.
\end{theorem}

\begin{proof}
The Taylor series assertion is clear, by differentiating. Regarding now the fact that the formula is indeed exact, if $f$ is the series in the statement, we have:
$$(1+t)f'(t)=pf(t)$$

Now by using this formula, we have the following computation:
$$\left((1+t)^{-p}f(t)\right)'=-p(1+t)^{-p-1}f(t)+(1+t)^{-p}f'(t)=0$$

Thus we have $f(t)=c(1+t)^p$, with $c=f(0)=1$, as desired.
\end{proof}

As another basic application of the derivatives and Taylor formula, we have:

\index{Taylor formula}
\index{trigonometric functions}
\index{exp and log}

\begin{theorem}
We have the following formulae,
$$\sin t=\sum_{l=0}^\infty(-1)^l\frac{t^{2l+1}}{(2l+1)!}\quad,\quad 
\cos t=\sum_{l=0}^\infty(-1)^l\frac{t^{2l}}{(2l)!}$$
as well as the following formulae,
$$e^t=\sum_{k=0}^\infty\frac{t^k}{k!}\quad,\quad 
\log(1+t)=\sum_{k=0}^\infty(-1)^{k+1}\frac{t^k}{k}$$
as Taylor series, and in general as well, with $|t|<1$ needed for $\log$.
\end{theorem}

\begin{proof}
There are several statements here, the proofs being as follows:

\medskip

(1) Regarding $\sin$ and $\cos$, we can use here the following formulae:
$$(\sin x)'=\cos x\quad,\quad 
(\cos x)'=-\sin x$$

Thus, we can differentiate $\sin$ and $\cos$ as many times as we want to, and so compute the corresponding Taylor series, and we obtain the formulae in the statement.

\medskip

(2) Regarding $\exp$ and $\log$, here the needed formulae, which lead to the formulae in the statement for the corresponding Taylor series, are as follows:
$$(e^x)'=e^x\quad,\quad 
(\log x)'=x^{-1}\quad,\quad
(x^p)'=px^{p-1}$$

(3) Finally, the fact that the Taylor formulae in the statement are exact is something standard too. Indeed, for $\exp$ this is clear, and for $\sin$, $\cos$, $\log$ this is something that can be deduced by using tricks like the one in the proof of Theorem 1.20.
\end{proof}

With this discussed, done with calculus? You must be kidding. Indeed, as a next piece for input for our considerations, we have the following key definition:

\index{Riemann integration}

\begin{definition}
The integral of a function $f$ on $[a,b]$ is its normalized average
$$\int_a^bf(x)dx=(b-a)\times A(f)$$
that is, average $A(f)$ times the length of the interval $b-a$. Equivalently, 
$$\int_a^bf(x)dx=(b-a)\times\lim_{n\to\infty}\frac{1}{n}\sum_{k=1}^nf\left(a+\frac{b-a}{n}\cdot k\right)$$
which can stand as a formal definition for the integral.
\end{definition}

As a first observation, in more geometric terms, the integral of a function $f:[a,b]\to\mathbb R$ is the signed area below its graph. Indeed, we can compute this area by approximating with rectangles, in the obvious way, and we are led to the limiting formula above. 

\bigskip

As an illustration now for this, with some arithmetic know-how, for the computation of sums of type $1^p+2^p+\ldots+N^p$, we have the following formula, for $p\in\mathbb N$:
$$\int_0^1x^pdx=
\lim_{N\to\infty}\frac{1^p+2^p+\ldots+N^p}{N^{p+1}}=\frac{1}{p+1}$$

However, such things remain a bit amateurish. At the more advanced level, the point is that the derivatives and integrals are related in several subtle ways, as follows:

\index{change of variable}
\index{partial integration}

\begin{theorem}
We have the following formulae, called fundamental theorem of calculus, integration by parts formula, and change of variable formula,
$$\int_a^bF'(x)dx=\Big[F\Big]_a^b$$
$$\int_a^b(f'g+fg')(x)dx=\Big[fg\Big]_a^b$$
$$\int_a^bf(x)dx=\int_{\varphi^{-1}(a)}^{\varphi^{-1}(b)}f(\varphi(t))\varphi'(t)dt$$
with the convention $[F]_a^b=F(b)-F(a)$, for the first two formulae.
\end{theorem}

\begin{proof}
To start with, given a continuous function $f:[a,b]\to\mathbb R$, by integrating $\min f\leq f\leq \max f$ we obtain the following formula, called mean value property:
$$\exists c\in[a,b]\quad,\quad \int_a^bf(x)dx=(b-a)f(c)$$

Next, this mean value property shows that we have the following implication:
$$I(x)=\int_a^xf(s)ds\implies I'=f$$

Now given $F:\mathbb R\to\mathbb R$ as in the statement, by using this with $f=F'$, we obtain $I'=F'$. Since $I(a)=0$, this reads $F(x)=I(x)+F(a)$, and with $x=b$ we get:
$$F(b)=\int_a^bF'(x)dx+F(a)$$

Thus, first formula proved, and the second and third formulae follow as well.
\end{proof}

And with this, good news, end of our basic learning of calculus. We can get back now to our estimates in Theorem 1.11, and improve them via calculus, as follows:

\begin{theorem}
We have the following estimate,
$$\left(\frac{N}{e}\right)^Ne<N!<\left(\frac{N}{e}\right)^Ne(N+1)$$
valid for any $N\in\mathbb N$.
\end{theorem}

\begin{proof}
This is something quite tricky, the idea being as follows:

\medskip

(1) In order to get our estimate, the idea is to integrate the logarithm. Since $\exp$ is increasing, $\log$ is increasing too, and by drawing rectangles, we obtain:
$$\log 1+\ldots+\log(N-1)<\int_1^N\log x\,dx<\log 2+\ldots+\log N$$

Which is good news, because we get the following estimate, involving factorials:
$$\log((N-1)!)<\int_1^N\log x\,dx<\log(N!)$$

Equivalently, by exponentiating, we have the following estimate:
$$(N-1)!<\exp\left(\int_1^N\log x\,dx\right)<N!$$

(2) So, in the end it all comes to finding a function $F$ satisfying $F'=\log$. Now for this purpose, let us start with a formula that we know from Theorem 1.15, namely:
$$(\log x)'=\frac{1}{x}$$

In order to make now appear $\log$ on the right, the idea is quite clear, namely multiplying on the left by $x$. We obtain in this way the following formula:
$$(x\log x)'=1\cdot\log x+x\cdot\frac{1}{x}=\log x+1$$

We are almost there, all we have to do is to substract $x$ from the left, as to get:
$$(x\log x-x)'=\log x$$

(3) Summarizing, integration problem solved, and by using this, we have:
\begin{eqnarray*}
\int_1^N\log x\,dx
&=&(N\log N-N)-(1\log 1-1)\\
&=&N\log N-N+1
\end{eqnarray*}

But with this, by going back to the estimate found in (1), we obtain:
$$(N-1)!<\left(\frac{N}{e}\right)^Ne<N!$$

Equivalently, with $N\to N+1$ on the left, we have the following estimate:
$$\left(\frac{N}{e}\right)^Ne<N!<\left(\frac{N+1}{e}\right)^{N+1}e$$

(4) Now observe that the upper bound on the right satisfies:
$$\left(\frac{N+1}{e}\right)^{N+1}e
=\left(\frac{N+1}{N}\right)^N\left(\frac{N}{e}\right)^N(N+1)
<\left(\frac{N}{e}\right)^Ne(N+1)$$

Thus, we are led to the formula in the statement.
\end{proof}

Coming next, as a consequence of Theorem 1.24, let us formulate:

\begin{theorem}
We have the following estimate, in the $N\to\infty$ limit,
$$N!\approx \left(\frac{N}{e}\right)^Ne$$
with the $\approx$ sign meaning that the quotient is much smaller than both factors.
\end{theorem}

\begin{proof}
This is clear from what we have in Theorem 1.24, which tells us that the quotient of the functions in the statement satisfies the following estimate:
$$N!\Big/\left(\frac{N}{e}\right)^Ne\in[1,N+1]$$

Thus, with the above convention for $\approx$, we are led to the formula in the statement.
\end{proof}

Good news, we can now approximate the binomial coefficients, as follows:

\index{binomial coefficients}
\index{central binomial coefficients}

\begin{theorem}
We have the following estimate for the binomial coefficients,
$$\binom{N}{K}\approx\left(\frac{1}{t^t(1-t)^{1-t}}\right)^N$$
in the $K\simeq tN\to\infty$ limit, with $t\in(0,1]$. In particular we have
$$\binom{2N}{N}\approx4^N$$
in the $N\to\infty$ limit, for the central binomial coefficients.
\end{theorem}

\begin{proof}
We know from Theorem 1.25 that we have $N!\approx \left(\frac{N}{e}\right)^Ne$, and in practice, with our convention for $\approx$, this can be replaced by $N!\approx \left(\frac{N}{e}\right)^N$. Now by using this:

\medskip

(1) The first estimate follows from the definition of the binomial coefficients:
\begin{eqnarray*}
\binom{N}{K}
&=&\frac{N!}{K!(N-K)!}\\
&\approx&\left(\frac{N}{e}\right)^N\left(\frac{e}{K}\right)^K\left(\frac{e}{N-K}\right)^{N-K}\\
&=&\frac{N^N}{K^K(N-K)^{N-K}}\\
&=&\frac{N^N}{(tN)^{tN}((1-t)N)^{(1-t)N}}\\
&=&\left(\frac{1}{t^t(1-t)^{1-t}}\right)^N
\end{eqnarray*}

(2) The second estimate follows from a similar computation, as follows:
$$\binom{2N}{N}
=\frac{(2N)!}{N!N!}
\approx\left(\frac{2N}{e}\right)^{2N}\left(\frac{e}{N}\right)^{2N}
=4^N$$

Alternatively, we can take $t=1/2$ in (1), then rescale via $N\to 2N$.
\end{proof}

\section*{1d. Stirling formula}

Shall we stop here, and go get a drink? Not yet. Indeed, the $\approx$ sign that we have in Theorems 1.25 and 1.26 is quite annoying, so let us improve that. We will need:

\begin{theorem}
We have the following formula,
$$\int_\mathbb Re^{-x^2}dx=\sqrt{\pi}$$
called Gauss integral formula.
\end{theorem}

\begin{proof}
This is something truly magic, the idea being as follows:

\medskip

(1) To start with, we can certainly integrate $e^{-x^2}$ by using the formula of the exponential series, and the primitive which is worth 0 at $x=0$ is given by:
$$\int e^{-x^2}=\sum_{k=0}^\infty(-1)^k\frac{x^{2k+1}}{(2k+1)k!}$$

However, this series is not computable, in terms of the known, familiar series.

\medskip

(2) Thus, no primitive, but we can still ask for the computation of $\int_\mathbb Re^{-x^2}dx$, who knows. And here, another surprise awaits us, this is simply undoable, with bare hands, I mean all the formulae and tricks that we learned so far fail, for this integral.

\medskip

(3) Which seems to send our problem to the trash can. However, and here comes the magic, the Gauss integral can be computed by using two dimensions, as follows:
\begin{eqnarray*}
\int_\mathbb R\int_\mathbb Re^{-x^2-y^2}dxdy
&=&4\int_0^\infty\int_0^\infty e^{-x^2-y^2}dxdy\\
&=&4\int_0^\infty\int_0^\infty e^{-t^2y^2-y^2}ydtdy\\
&=&4\int_0^\infty\int_0^\infty ye^{-y^2(1+t^2)}dydt\\
&=&2\int_0^\infty\int_0^\infty\left(-\frac{e^{-y^2(1+t^2)}}{1+t^2}\right)'dydt\\
&=&2\int_0^\infty\frac{dt}{1+t^2}\\
&=&2\int_0^\infty(\arctan t)'dt\\
&=&\pi
\end{eqnarray*}

(4) To be more precise, we have used here at the end the fact that the derivative of $\arctan=\tan^{-1}$, with $\tan=\sin/\cos$ as usual, is the function $1/(1+t^2)$, and with this coming from the various formulae and rules that we have, good exercise for you.
\end{proof}

We can now establish a key result, the Stirling formula, as follows:

\index{Stirling formula}
\index{Riemann sum}
\index{Gauss formula}
\index{Laplace method}
\index{factorial}

\begin{theorem}
We have the Stirling formula
$$N!\simeq\left(\frac{N}{e}\right)^N\sqrt{2\pi N}$$
valid in the $N\to\infty$ limit.
\end{theorem}

\begin{proof}
This is something quite tricky, the idea being as follows:

\medskip

(1) We recall, from the proof of Theorem 1.24, that with rectangles we get:
$$\log(N!)
\approx\int_1^N\log x\,dx
=N\log N-N+1$$

By exponentiating, this gives the following estimate, which is not bad:
$$N!\approx\left(\frac{N}{e}\right)^N\cdot e$$

(2) We can improve this by replacing the rectangles with trapezoids. We get:
\begin{eqnarray*}
\log(N!)
&\approx&\int_1^N\log x\,dx+\frac{\log 1+\log N}{2}\\
&=&N\log N-N+1+\frac{\log N}{2}
\end{eqnarray*}

By exponentiating, this gives the following estimate, which gets us closer:
$$N!\approx\left(\frac{N}{e}\right)^N\cdot e\cdot\sqrt{N}$$

(3) In order to conclude, we must take some kind of mathematical magnifier, and carefully estimate the error made in (2). Fortunately, this mathematical magnifier exists, called Euler-Maclaurin formula, and after some computations, this leads to Stirling.

\medskip

(4) Alternatively, here is another approach to (3), which better does the job, explaining where the $\sqrt{2\pi}$ factor comes from. First, by iterated partial integration we have:
$$N!=\int_0^\infty x^Ne^{-x}dx$$

(5) Since the integrand is sharply peaked at $x=N$, as you can see by computing the derivative of $\log(x^Ne^{-x})$, this suggests writing $x=N+y$, and we obtain:
\begin{eqnarray*}
\log(x^Ne^{-x})
&=&N\log x-x\\
&=&N\log(N+y)-(N+y)\\
&=&N\log N+N\log\left(1+\frac{y}{N}\right)-(N+y)\\
&\simeq&N\log N+N\left(\frac{y}{N}-\frac{y^2}{2N^2}\right)-(N+y)\\
&=&N\log N-N-\frac{y^2}{2N}
\end{eqnarray*}

(6) By exponentiating, we obtain from this the following estimate:
$$x^Ne^{-x}\simeq\left(\frac{N}{e}\right)^Ne^{-y^2/2N}$$

(7) Now by integrating, we obtain from this the following estimate, as desired:
\begin{eqnarray*}
N!
&=&\int_0^\infty x^Ne^{-x}dx\\
&\simeq&\int_{-N}^N\left(\frac{N}{e}\right)^Ne^{-y^2/2N}\,dy\\
&\simeq&\left(\frac{N}{e}\right)^N\int_\mathbb Re^{-y^2/2N}\,dy\\
&=&\left(\frac{N}{e}\right)^N\sqrt{2N}\int_\mathbb Re^{-z^2}\,dz\\
&=&\left(\frac{N}{e}\right)^N\sqrt{2\pi N}
\end{eqnarray*}

(8) And exercise of course for you to learn more about this, as much as you can.
\end{proof}

Good news, we can now professionally estimate the binomials, as follows:

\index{binomial coefficients}
\index{central binomial coefficients}

\begin{theorem}
We have the following estimate for binomial coefficients,
$$\binom{N}{K}\simeq\left(\frac{1}{t^t(1-t)^{1-t}}\right)^N\frac{1}{\sqrt{2\pi t(1-t)N}}$$
in the $K\simeq tN\to\infty$ limit, with $t\in(0,1]$. In particular we have
$$\binom{2N}{N}\simeq\frac{4^N}{\sqrt{\pi N}}$$
in the $N\to\infty$ limit, for the central binomial coefficients.
\end{theorem}

\begin{proof}
All this is very standard, by using the Stirling formula, as follows:

\medskip

(1) The first estimate comes from the definition of the binomials, as follows:
\begin{eqnarray*}
\binom{N}{K}
&=&\frac{N!}{K!(N-K)!}\\
&\simeq&\left(\frac{N}{e}\right)^N\sqrt{2\pi N}\left(\frac{e}{K}\right)^K\frac{1}{\sqrt{2\pi K}}\left(\frac{e}{N-K}\right)^{N-K}\frac{1}{\sqrt{2\pi(N-K)}}\\
&=&\frac{N^N}{K^K(N-K)^{N-K}}\sqrt{\frac{N}{2\pi K(N-K)}}\\
&\simeq&\frac{N^N}{(tN)^{tN}((1-t)N)^{(1-t)N}}\sqrt{\frac{N}{2\pi tN(1-t)N}}\\
&=&\left(\frac{1}{t^t(1-t)^{1-t}}\right)^N\frac{1}{\sqrt{2\pi t(1-t)N}}
\end{eqnarray*}

(2) The second estimate follows from a similar computation, as follows:
$$\binom{2N}{N}
=\frac{(2N)!}{N!N!}
\simeq\left(\frac{2N}{e}\right)^{2N}\sqrt{4\pi N}\left(\frac{e}{N}\right)^{2N}\frac{1}{2\pi N}
=\frac{4^N}{\sqrt{\pi N}}$$

Alternatively, we can take $t=1/2$ in (1), then rescale via $N\to 2N$.
\end{proof}

\section*{1e. Exercises}

We learned many things about counting and analysis, and as exercises, we have:

\begin{exercise}
Further meditate on $0!=1$, but not too much.
\end{exercise}

\begin{exercise}
Memorize as many binomials and factorials you can.
\end{exercise}

\begin{exercise}
Look up and learn other interesting formulae relating binomials.
\end{exercise}

\begin{exercise}
Learn some more, about the arithmetic of the binomials.
\end{exercise}

\begin{exercise}
Learn also about the Catalan numbers, and their properties.
\end{exercise}

\begin{exercise}
Review if needed the basic theory of sequences and convergence.
\end{exercise}

\begin{exercise}
Review what we said above, about derivatives and integrals.
\end{exercise}

\begin{exercise}
Learn other proofs, and improvements, of the Stirling formula.
\end{exercise}

As bonus exercise, learn some programming too. You will need this, later.

\chapter{Basic probability}

\section*{2a. Random variables}

Welcome to probability, which will be of discrete type, to start with. Inspired by what happens with simple games, from the real life, let us have as starting point:

\index{probability space}
\index{discrete space}
\index{events}
\index{probability function}
\index{probability 1}

\begin{definition}
A discrete probability space is a set $X$, usually finite or countable, whose elements $x\in X$ are called events, together with a function
$$P:X\to[0,\infty)$$
called probability function, which is subject to the condition
$$\sum_{x\in X}P(x)=1$$
telling us that the overall probability for something to happen is $1$.
\end{definition}

Which sounds quite neat, this is definitely a rigorous mathematical definition, that we can build our theory upon. Here are a few comments, in relation with this:

\bigskip

(1) As a first comment, our condition $\sum_{x\in X}P(x)=1$ perfectly makes sense, and this even if $X$ is uncountable, because the sum of positive numbers is always defined, as a number in $[0,\infty]$, and this no matter how many positive numbers we have.

\bigskip

(2) As a second comment, we have chosen in the above not to assume that $X$ is finite or countable, and this for instance because we want to be able to regard any probability function on $\mathbb N$ as a probability function on $\mathbb R$, by setting $P(x)=0$ for $x\notin\mathbb N$. 

\bigskip

(3) As a third comment, standing as a complement to Definition 2.1, once we have a probability function $P:X\to[0,\infty)$, given a  set of events $Y\subset X$, we can say that the probability for the events of type $y\in Y$ to happen is $P(Y)=\sum_{y\in Y}P(y)$.

\bigskip

As our first result now, making the link with our intuition from real life, we have:

\begin{proposition}
In the simplest case, that where $X$ is finite and $P$ is constant,
$$P(Y)=\frac{{\rm number\ of\ times}\ Y\ {\rm happens}}{{\rm total\ number\ of\ possibilities}}$$
for any set of events $Y\subset X$.
\end{proposition}

\begin{proof}
Assume indeed that our probability function $P:X\to[0,\infty)$ is constant. But then, $\sum_{x\in X}P(x)=1$ implies that $X$ must be finite, and that we must have:
$$P(x)=\frac{1}{|X|}\quad,\quad\forall x\in X$$

Now by summing this over the elements $y\in Y$ of a subset $Y\subset X$, we obtain:
$$P(Y)=\frac{|Y|}{|X|}\quad,\quad\forall Y\subset X$$

We are therefore led to the conclusion in the statement.
\end{proof}

Time now for some examples? The basic examples will come of course from basic games, and at the beginning of everything, we certainly have flipping coins:

\begin{example}
Flipping coins.
\end{example}

Here things are simple and clear, because when you flip a coin the corresponding discrete probability space, together with its probability measure, is as follows:
$$X=\big\{{\rm heads},\,{\rm tails}\big\}\quad,\quad P({\rm heads})=P({\rm tails})=\frac{1}{2}$$

In the case where the coin is biased, as to land on heads with probability $2/3$, and on tails with probability $1/3$, the corresponding probability space is as follows:
$$X=\big\{{\rm heads},\,{\rm tails}\big\}\quad,\quad P({\rm heads})=\frac{2}{3}\quad,\quad P({\rm tails})=\frac{1}{3}$$

More generally, given any number $p\in[0,1]$, we have an abstract probability space as follows, where we have replaced heads and tails by win and lose:
$$X=\big\{{\rm win},\,{\rm lose}\big\}\quad,\quad P({\rm win})=p\quad,\quad P({\rm lose})=1-p$$

Things become even more interesting when flipping a coin, biased or not, several times in a row. In the case of a biased coin, landing on heads with probability $p_1$ and on tails with probability $p_2=1-p_2$, thrown $n$ times, the probability space is:
$$X=\big\{1,2\}^n\quad,\quad P(i_1\ldots i_n)=p_{i_1}\ldots p_{i_n}\quad,\quad p_i\geq0\ ,\ p_1+p_2=1$$

Finally, in relation with this latter example, let us check that the sum 1 condition in Definition 2.1 is indeed satisfied. But this is indeed the case, as shown by:
$$\sum_{i\in X}P(i)
=\sum_{i_1,\ldots,i_n}p_{i_1}\ldots p_{i_n}
=\sum_{i_1}p_{i_1}\ldots\sum_{i_n}p_{i_n}
=1$$

Summarizing, coins reasonably understood, or at least the mathematical formalism of the coin game, reasonably understood. However, in what regards the mathematics itself of the coin game, many interesting things can be said. More on this later.

\bigskip

Moving on, at a more advanced level we can talk about rolling dice, as follows:

\begin{example}
Rolling dice.
\end{example}

Again, things here are simple and clear, because when you throw a die the corresponding probability space, together with its probability measure, is as follows:
$$X=\big\{1,\ldots,6\big\}\quad,\quad P(i)=\frac{1}{6}\ ,\ \forall i$$

As before with the coins, we can further complicate this by assuming that the die is biased, say landing on face $i$ with probability $p_i\in[0,1]$. In this case the corresponding probability space, together with its probability measure, is as follows:
$$X=\big\{1,\ldots,6\big\}\quad,\quad P(i)=p_i\quad,\quad p_i\geq0\ ,\ \sum_ip_i=1$$

Also as before with coins, things become more interesting when throwing a die several times in a row, or equivalently, when throwing several identical dice at the same time. In this latter case, with $n$ identically biased dice, the probability space is as follows:
$$X=\big\{1,\ldots,6\big\}^n\quad,\quad P(i_1\ldots i_n)=p_{i_1}\ldots p_{i_n}\quad,\quad p_i\geq0\ ,\ \sum_ip_i=1$$

Observe that the sum 1 condition in Definition 2.1 is indeed satisfied, and with this proving that our dice modeling is bug-free, due to the following computation:
$$\sum_{i\in X}P(i)
=\sum_{i_1,\ldots,i_n}p_{i_1}\ldots p_{i_n}
=\sum_{i_1}p_{i_1}\ldots\sum_{i_n}p_{i_n}
=1$$

In fact, this computation is identical to the one that we did before, for coins.

\bigskip

Getting back now to theory, in the general context of Definition 2.1, we can see that what we have there is very close to the biased die, from Example 2.4. Indeed, in the general context of Definition 2.1, we can say that what happens is that we have a die with $|X|$ faces, which is biased such that it lands on face $i$ with probability $P(i)$. 

\bigskip

Which is something quite interesting, theoretically speaking, in relation with the discrete probability theory that we want to develop in this Part I, because we have now some intuition on all this. So, as a conclusion, let us record this finding as follows:

\begin{principle}
Discrete probability can be understood as being about throwing a general die, having an arbitrary number of faces, and which is arbitrarily biased too:
$$\xymatrix@R=7pt@C=7pt{
\ar@{-}[rr]\ar@{-}[dd]&&\ar@{-}[dd]&&\ar@{-}[rr]\ar@{-}[dd]&&\ar@{-}[dd]&&\ar@{-}[rr]\ar@{-}[dd]&&\ar@{-}[dd]&&\ar@{-}[rr]\ar@{-}[dd]&&\ar@{-}[dd]\\
&p_x&&&&p_y&&&&p_z&&&&p_t&&&&\ldots\\
\ar@{-}[rr]&&&&\ar@{-}[rr]&&&&\ar@{-}[rr]&&&&\ar@{-}[rr]&&}$$
To be more precise, when writing our probability space as $X=\{x,y,z,t,\ldots\}$, with the corresponding probabilities being $p_x,p_y,p_z,p_t,\ldots\,$, the picture is the one above.
\end{principle}

Getting now to what we want to do in this chapter, namely further building on Definition 2.1, let us formulate an intuitive and down-to-earth question, as follows:

\begin{question}
As usual when doing probability computations, we are mainly interested in winning. But, winning what?
\end{question}

And good question this is, because Definition 2.1 as stated does not provide an answer to it. So, let us look first at the examples. In case we are rolling a die, what we win is what the die says. And, with this convention, we are led to the following result:

\begin{proposition}
When rolling a die, what we win on the average is
$$E=\sum_ii\times p_i$$
which for a usual die reads $E=3.5$. More generally, when $X\subset\mathbb R$, what we win is
$$E=\sum_ii\times P(i)$$
with the sum being over all the elements $i\in X$ of our probability space.
\end{proposition}

\begin{proof}
Indeed, when rolling a usual die, what we win on the average is:
$$E=\frac{1+2+3+4+5+6}{6}=3.5$$

As for the other formulae above, these are theoretical, coming from this.
\end{proof}

With this understood, what about coins? Here, before doing any computation, we must assign some numbers to our events, and a standard choice here is as follows:
$$f:\big\{{\rm heads},\,{\rm tails}\big\}\to\mathbb R\quad,\quad f({\rm heads})=1\quad,\quad f({\rm tails})=0$$ 

With this choice made, what we can expect to win is the following quantity:
\begin{eqnarray*}
E(f)
&=&f({\rm heads})\times P({\rm heads})+f({\rm tails})\times P({\rm tails})\\
&=&1\times\frac{1}{2}+0\times\frac{1}{2}\\
&=&\frac{1}{2}
\end{eqnarray*}

When dealing with a biased coin, as in Example 2.3, the computation becomes:
\begin{eqnarray*}
E(f)
&=&f({\rm heads})\times P({\rm heads})+f({\rm tails})\times P({\rm tails})\\
&=&1\times p+0\times(1-p)\\
&=&p
\end{eqnarray*}

However, the story with coins is not over here, because with a different convention for $f$, we will reach to a different outcome. So, let us try to understand this. We have:

\begin{proposition}
When flipping a coin, and with our winning function being
$$f:\big\{{\rm heads},\,{\rm tails}\big\}\to\mathbb R\quad,\quad f({\rm heads})=a\quad,\quad f({\rm tails})=b$$
with $a,b\in\mathbb R$, what we can expect to win, on average, is the quantity
$$E(f)=ap+b(1-p)$$
in the case of a biased coin. For a normal coin, this reads $E(f)=(a+b)/2$.
\end{proposition}

\begin{proof}
With the above convention for the winning function, we have indeed:
\begin{eqnarray*}
E(f)
&=&f({\rm heads})\times P({\rm heads})+f({\rm tails})\times P({\rm tails})\\
&=&a\times p+b\times(1-p)
\end{eqnarray*}

Thus, we are led to the conclusions in the statement.
\end{proof}

Summarizing, in order to do some mathematics, in the context of Definition 2.1, we need a random variable $f:X\to\mathbb R$, and the math will consist in computing the expectation of this variable, $E(f)\in\mathbb R$. Let us axiomatize this situation as follows:

\index{random variable}
\index{average gain}
\index{expectation}

\begin{definition}
A random variable on a probability space $X$ is a function
$$f:X\to\mathbb R$$
and the expectation of such a random variable is the quantity
$$E(f)=\sum_{x\in X}f(x)P(x)$$
which is best thought as being the average gain, when the game is played.
\end{definition}

Here the word ``game'' refers to the probability space interpretation from Principle 2.5. Indeed, in that context, with our discrete set of events $X$ being thought of as corresponding to a generalized die, and by thinking of $f:X\to\mathbb R$ as representing some sort of money, the above quantity $E(f)$ is what we win, on average, when playing the game.

\bigskip

Getting back now to the examples, as our next result, we have:

\begin{theorem}
When flipping a usual coin $n$ times in a row:
\begin{enumerate}
\item $E({\rm heads})=n/2$.

\item $E({\rm tails})=n/2$.

\item $E(a\times{\rm heads}+b\times{\rm tails})=n(a+b)/2$.
\end{enumerate} 
\end{theorem}

\begin{proof}
Let us first fix some mathematical notation, in the spirit of what we have in Definition 2.9. Consider the variables in the statement, namely:
$$f:\big\{{\rm heads},\,{\rm tails}\big\}^n\to\mathbb R\quad,\quad f={\rm number\ of\ heads}$$ 
$$g:\big\{{\rm heads},\,{\rm tails}\big\}^n\to\mathbb R\quad,\quad g={\rm number\ of\ tails}$$ 
$$h:\big\{{\rm heads},\,{\rm tails}\big\}^n\to\mathbb R\quad,\quad h=a\times{\rm heads}+b\times{\rm tails}$$ 

With a bit of probabilistic know-how, we can compute the expectations of these variables, by using several methods, which are all instructive, as follows:

\medskip

(1) \underline{Symmetry proof}. Observe that, by using the obvious symmetry between heads and tails, at the level of the corresponding expectations, we have:
$$E(f)=E(g)$$

On the other hand, we also have $f+g=n$, which gives the following formula:
$$E(f)+E(g)=E(f+g)=E(n)=n$$

Thus, we reach to the following conclusion, which is the one in the statement:
$$E(f)=E(g)=\frac{n}{2}$$

As for the general variable $h$, no need here for new computations, because we have:
\begin{eqnarray*}
E(h)
&=&E(af+bg)\\
&=&aE(f)+bE(g)\\
&=&an/2+bn/2\\
&=&n(a+b)/2
\end{eqnarray*}

Thus, we are led to the conclusions in the statement.

\medskip

(2) \underline{Additive proof}. We can argue here that the formulae in the statement are all obvious, by using the fact that, when flipping our coin several times in a row, the expectations will sum up. That is, if we denote by $f_n,g_n,h_n$ our variables, we have:
$$E(f_n)=nE(f_1)\quad,\quad E(g_n)=nE(g_1)\quad,\quad E(h_n)=nE(h_1)$$

Indeed, by using these formulae, along with our computations before at $n=1$, telling us that we have $E(f_1)=E(g_1)=1/2$ and $E(h_1)=(a+b)/2$, we obtain, as desired:
$$E(f_n)=E(g_n)=\frac{n}{2}\quad,\quad E(h_n)=\frac{n(a+b)}{2}$$

(3) \underline{Algebraic proof}. Here is as well a third proof, which is my favorite one, having the advantage of using 0 thinking and neurons, or almost. We have indeed:
\begin{eqnarray*}
E(f)
&=&\sum_{x\in\{{\rm heads,\, tails}\}^n}\#({\rm heads}\in x)\times\frac{1}{2^n}\\
&=&\sum_{s=1}^n\sum_{\#({\rm heads}\in x)=s}\#({\rm heads}\in x)\times\frac{1}{2^n}\\
&=&\sum_{s=1}^n\binom{n}{s}\times s\times\frac{1}{2^n}
\end{eqnarray*}

But with this in hand, we can kill the problem with some algebra, as follows:
\begin{eqnarray*}
E(f)
&=&\frac{1}{2^n}\sum_{s=1}^n\frac{n!}{(s-1)!(n-s)!}\\
&=&\frac{n}{2^n}\sum_{t=0}^{n-1}\binom{n-1}{t}\\
&=&\frac{n}{2^n}\times 2^{n-1}\\
&=&\frac{n}{2}
\end{eqnarray*}

Similarly, $E(g)=n/2$, and then $E(h)=n(a+b)/2$ comes by linearity, as in (1).
\end{proof}

Getting now to the more general case of a biased coin, the result here is as follows:

\index{biased coin}

\begin{theorem}
When flipping a biased coin $n$ times in a row:
\begin{enumerate}
\item $E({\rm heads})=np$.

\item $E({\rm tails})=n(1-p)$.

\item $E(a\times{\rm heads}+b\times{\rm tails})=n(ap+b(1-p))$.
\end{enumerate} 
\end{theorem}

\begin{proof}
As before with the usual coins, if we denote by $f,g,h$ our random variables, the formulae in the statement for $E(f)$ and $E(g)$ give the formula for $E(h)$, by linearity. Thus, we are left with computing $E(f)$ and $E(g)$, and regarding the 3 possible approaches here, from the proof of Theorem 2.10, the situation with them is as follows:

\medskip

(1) \underline{Symmetry proof}. Normally, the idea here would be to exploit the symmetry between heads and tails, with the goal of reaching to the following formula:
$$(1-p)E(f)=pE(g)$$

Indeed, on the other hand, we have $f+g=n$, which gives the following formula:
$$E(f)+E(g)=E(f+g)=E(n)=n$$

Thus, we would reach to the following conclusion, which is the one in the statement:
$$E(f)=np\quad,\quad E(g)=n(1-p)$$

The problem, however, is that the symmetry study, leading to $(1-p)E(f)=pE(g)$, is not obvious when $p\neq1/2$. So, relax, and we should better run away from this. 

\medskip

(2) \underline{Additive proof}. As before, we can argue that the formulae in the statement are all obvious, by using the fact that, when flipping our coin several times in a row, the expectations will sum up. That is, if we denote by $f_n,g_n,h_n$ our variables, we have:
$$E(f_n)=nE(f_1)\quad,\quad E(g_n)=nE(g_1)\quad,\quad E(h_n)=nE(h_1)$$

Indeed, using these formulae, along with our computations before at $n=1$, telling us that we have $E(f_1)=p$, $E(g_1)=1-p$, $E(h_1)=ap+b(1-p)$, we obtain, as desired:
$$E(f_n)=np\quad,\quad E(g_n)=n(1-p)\quad,\quad E(h_n)=n(ap+b(1-p))$$

(3) \underline{Algebraic proof}. And here is as well the algebraic proof, certainly a bit less elegant, but having the advantage of using 0 thinking and neurons, or almost:
\begin{eqnarray*}
E(f)
&=&\sum_{s=1}^n\binom{n}{s}sp^s(1-p)^{n-s}\\
&=&np\sum_{s=1}^n\frac{(n-1)!}{(s-1)!(n-s)!}\,p^{s-1}(1-p)^{n-s}\\
&=&np\sum_{t=0}^{n-1}\binom{n-1}{t}p^t(1-p)^{n-t-1}\\
&=&np(p+1-p)^{n-1}\\
&=&np
\end{eqnarray*}

Indeed, a similar computation gives $E(g)=n(1-p)$, and then the general formula, $E(h)=n(ap+b(1-p))$, comes from this, by linearity, as explained before.
\end{proof}

Next, come the dice. Here there are several possible natural choices for our random variable, and regarding the essentials, in the unbiased die case, we have:

\begin{theorem}
When rolling a usual die $n$ times in a row:
\begin{enumerate}
\item For $f_i=\#i$, we have $E(f_i)=n/6$.

\item For $g_i=i\#i$, we have $E(g_i)=ni/6$.

\item For $h={\rm sum}$, we have $E(h)=3.5n$.
\end{enumerate} 
\end{theorem}

\begin{proof}
The story here is quite similar to that of a usual coin, and skipping some straightforward details, and going directly to the three possible methods for proving such things, from the proof of Theorem 2.10, the situation is as follows:

\medskip

(1) \underline{Symmetry proof}. The variables $f_i=\#i$, counting the number of occurrences of $i\in\{1,\ldots,6\}$ during our $n$ throws, obviously satisfy the following condition:
$$E(f_1)=\ldots=E(f_6)$$

On the other hand, we have as well the following trivial computation:
$$f_1+\ldots+f_6=n\implies E(f_1)+\ldots+E(f_6)=n$$

Thus we have $E(f_i)=n/6$ for any $i$, as stated. Then, by multiplying by $i$, we have $E(g_i)=ni/6$. And finally, by summing over $i$ we obtain $E(h)=3.5n$.

\medskip

(2) \underline{Additive proof}. We can argue here that the formulae in the statement are all obvious, using the intuitive fact that, when rolling our die several times in a row, the expectations will sum up. That is, if we denote by $f_i^n,g_i^n,h^n$ our variables, we have:
$$E(f_i^n)=nE(f_i^1)\quad,\quad E(g_i^n)=nE(g_i^1)\quad,\quad E(h^n)=nE(h^1)$$

Indeed, by using these formulae, along with our computations before at $n=1$, telling us that we have $E(f_i^1)=1/6$, $E(g_i^1)=i/6$, $E(h^1)=3.5$, we obtain, as desired:
$$E(f_i^n)=\frac{n}{6}\quad,\quad E(g_i^n)=\frac{ni}{6}\quad,\quad E(h^n)=3.5n$$

(3) \underline{Algebraic proof}. Things are a bit more tricky here, depending on the exact variables that we have in mind, the situation being as follows:

\medskip

-- Normally, we can compute indeed $E(f_i)$ algebrically, a bit as before for the coins, by using binomial coefficients and some summing, and I will leave this as an exercise for you, and then pass to $E(g_i)$, and then to $E(h)$, as explained in (1), by linearity.

\medskip

-- However, if we want to compute for instance $E(h)$ directly, algebrically, we are a bit in trouble, because, as you can convince yourself by doing some computations, finding the formula of the distribution of the sum $h$ is a quite delicate question.
\end{proof}

Finally, in the case of a biased die, the result is as follows:

\index{biased die}

\begin{theorem}
When rolling a biased die $n$ times in a row:
\begin{enumerate}
\item For $f_i=\#i$, we have $E(f_i)=np_i$.

\item For $g_i=i\#i$, we have $E(g_i)=nip_i$.

\item For $h={\rm sum}$, we have $E(h)=n\sum_iip_i$.
\end{enumerate} 
\end{theorem}

\begin{proof}
This is something hybrid between Theorem 2.11 and Theorem 2.12, and skipping the various straightforward details here, the situation is as follows:

\medskip

(1) \underline{Symmetry proof}. This won't quite work for a biased die, for the technical reasons explained in the proof of Theorem 2.11, dealing with a biased coin.

\medskip

(2) \underline{Additive proof}. Good news here, this works, as usual, because if we denote by $f_i^n,g_i^n,h^n$ our variables, it is clear that we have the following formulae:
$$E(f_i^n)=nE(f_i^1)\quad,\quad E(g_i^n)=nE(g_i^1)\quad,\quad E(h^n)=nE(h^1)$$

Indeed, by using these formulae, along with our computations before at $n=1$, telling us that we have $E(f_i^1)=p_i$, $E(g_i^1)=ip_i$, $E(h^1)=\sum_iip_i$, we obtain, as desired:
$$E(f_i^n)=np_i\quad,\quad E(g_i^n)=nip_i\quad,\quad E(h^n)=n\sum_iip_i$$

(3) \underline{Algebraic proof}. This won't work well for a biased die, for the technical reasons explained in the proof of Theorem 2.12, dealing with a usual die.
\end{proof}

\section*{2b. Variance, moments}

We have already seen some good illustrations for Definition 2.9, so time now to get into more delicate aspects. Imagine that you want to set up some sort of business, with your variable $f:X\to\mathbb R$. You are of course mostly interested in the expectation $E(f)\in\mathbb R$, but passed that, the way this expectation comes in matters too. For instance:

\bigskip

-- When your variable is constant, $f=c$, you certainly have $E(f)=c$, and your business will run smoothly, with not so many surprises on the way.

\bigskip

-- On the opposite, for a complicated variable satisfying $E(f)=c$, your business will be more bumpy, with wins and loses on the way, depending on your skills.

\bigskip

In short, and extrapolating now from business to mathematics, physics, chemistry and everything else, we must complement Definition 2.9 with something finer, regarding the ``quality'' of the expectation $E(f)\in\mathbb R$ appearing there. And the first thought here, which is the correct one, goes to the following definition, of the variance of our variable:

\index{variance}

\begin{definition}
The variance of a variable $f:X\to\mathbb R$ is the quantity
$$V(f)=E\left((f-E(f))^2\right)$$
intuitively measuring how far is $f$ from a constant variable.
\end{definition}

As a first observation, which is in tune with what we were saying in the above, the variance is 0 precisely when the variable is constant, equal to its expectation:
\begin{eqnarray*}
V(f)=0
&\iff&f={\rm constant}\\
&\iff&f=E(f)
\end{eqnarray*}

Many other things, more specialized, can be said, as a continuation of this. In general now, the variance is best computed by using the following result:

\index{variance formula}
\index{Cauchy-Schwarz}

\begin{theorem}
The variance of a variable $f:X\to\mathbb R$ is given by
$$V(f)=E(f^2)-E(f)^2$$
with $E$ denoting as usual the expectation.
\end{theorem}

\begin{proof}
This is something straightforward, the idea being as follows:

\medskip

(1) To start with, our formula holds indeed, as shown by the following computation:
\begin{eqnarray*}
V(f)
&=&E\left((f-E(f))^2\right)\\
&=&E\left(f^2-2E(f)f+E(f)^2\right)\\
&=&E(f^2)-2E(f)^2+E(f)^2\\
&=&E(f^2)-E(f)^2
\end{eqnarray*}

(2) However, as an interesting consequence, we have the following inequality, and the question is whether we can prove directly this inequality, with bare hands:
$$E(f^2)\geq E(f)^2$$

In practice, using Definition 2.9, this amounts in proving an inequality as follows, for any numbers $p_i\in[0,1]$ satisfying $\sum_ip_i=1$, and any real numbers $f_i\in\mathbb R$:
$$\sum_if_i^2p_i\geq\left(\sum_if_ip_i\right)^2$$

Now in order to prove this, let us write this inequality as follows:
$$\sum_if_i^2p_i\sum_ip_i\geq\left(\sum_if_ip_i\right)^2$$

(3) Our claim now is that this latter inequality, which might look a bit mysterious, is in fact the very familiar Cauchy-Schwarz inequality, which is as follows:
$$\sum_ia_i^2\sum_ib_i^2\geq\left(\sum_ia_ib_i\right)^2$$

Indeed, the passage from this to the inequality in (3) can be done by setting:
$$a_i=f_i\sqrt{p_i}\quad,\quad b_i=\sqrt{p_i}$$

As for the passage from (3) to Cauchy-Schwarz, this can be done too, as follows:
$$f_i=\frac{a_i}{b_i}\quad,\quad p_i=b_i^2$$

(4) Summarizing, saved by Cauchy-Schwarz, we can say now that we have a good understanding of $V(f)\geq0$. However, for finishing this discussion, actually, where does Cauchy-Schwarz come from? In answer, consider the following function:
\begin{eqnarray*}
\varphi(t)
&=&||at+b||^2\\
&=&<at+b,at+b>\\
&=&<a,a>t^2+2<a,b>t+<b,b>\\
&=&||a||^2t^2+2<a,b>t+||b||^2
\end{eqnarray*}

Since we have $\varphi\geq0$, the discriminant must be negative, $\Delta\leq0$, which reads:
$$4<a,b>^2-4||a||^2||b||^2\leq 0$$

Thus, we have the following inequality, which is the Cauchy-Schwarz one above:
$$||a||\cdot ||b||\geq|<a,b>|$$

And with this, we have now a truly good understanding of $V(f)\geq0$.
\end{proof}

As a complement to Definition 2.14, let us formulate as well:

\index{standard deviation}

\begin{definition}
The standard deviation of a variable $f:X\to\mathbb R$ is the quantity
$$\sigma(f)=\sqrt{E\left((f-E(f))^2\right)}$$
so that $V(f)=\sigma(f)^2$, intuitively measuring too how far is $f$ from a constant variable.
\end{definition}

Moving on, but still regarding expectations and variances of variables, as constructed above, at the abstract level, following Markov, we have the following result:

\index{Markov inequality}

\begin{theorem}
We have the following Markov inequality,
$$P\big(|f|\geq a\big)\leq\frac{E(f)}{a}$$
valid for any random variable $f:X\to\mathbb R$, and any $a>0$.
\end{theorem}

\begin{proof}
This is something trivial, coming from definitions. Indeed, we have:
\begin{eqnarray*}
E(f)
&\geq&\sum_{|f(x)|\geq a}f(x)P(x)\\
&\geq&\sum_{|f(x)|\geq a}aP(x)\\
&=&a\sum_{|f(x)|\geq a}P(x)\\
&=&aP\big(|f|\geq a\big)
\end{eqnarray*}

Thus, we are led to the Markov inequality in the statement.
\end{proof}

Next, following Chebycheff, we have the following key estimate:

\index{Chebycheff inequality}

\begin{theorem}
We have the following Chebycheff inequality,
$$P\big(|g-E|\geq b\big)=\frac{V}{b^2}$$
valid for any random variable $g:X\to\mathbb R$, having mean $E$ and variance $V$.
\end{theorem}

\begin{proof}
Given a random variable $g:X\to\mathbb R$, having mean $E$ and variance $V$, we can use the Markov inequality above with $f=(g-E)^2$, and we obtain in this way:
\begin{eqnarray*}
P\big(|g-E|\geq b\big)
&=&P\big((g-E)^2\geq b^2\big)\\
&\leq&\frac{E((g-E)^2)}{b^2}\\
&=&\frac{V}{b^2}
\end{eqnarray*}

Thus, we are led to the Chebycheff inequality in the statement.
\end{proof}

Getting back now to the real-life considerations from the beginning of this section, good that we have the notion of variance, for understanding how to run our business, but let us not stop here. For a total control of your business, be that of financial, mathematical, physical or chemical type, you will certainly want to know more about your variable $f:X\to\mathbb R$. Which leads us into general moments, constructed as follows:

\index{moments}

\begin{definition}
The moments of a variable $f:X\to\mathbb R$ are the numbers
$$M_k=E(f^k)$$
which satisfy $M_0=1$, then $M_1=E(f)$, and then $V(f)=M_2-M_1^2$.
\end{definition}

And with this, good news, we have now all the needed tools in our bag for doing some good business, in what follows. To put things in a very compacted way:

\medskip

-- $M_0$ is about foundations.

\medskip

-- $M_1$ is about running some business.

\medskip

-- $M_2$ is about running that business well.

\medskip

-- $M_3$ and higher are advanced level, about ruining all the competing businesses. 

\medskip

As a last piece of discussion, still regarding the moments, we can formulate the following version of Definition 2.19, making a more clear link with the variance:

\index{central moments}

\begin{definition}
The central moments of a variable $f:X\to\mathbb R$ are the numbers
$$M_k'=E((f-E)^k)$$
with $E=E(f)$, which satisfy $M_0'=1$, $M_1'=0$, $M_2'=V(f)$.
\end{definition}

In other words, the central moments of a variable $f:X\to\mathbb R$ are the moments of the shifted variable $f-E:X\to\mathbb R$, with the shift being there for having 0 mean. Which is something quite nice, and we will sometimes use this, instead of Definition 2.20.

\bigskip

Along the same lines, we can do in fact a bit better, as follows:

\index{normalized central moments}

\begin{definition}
The normalized central moments of $f:X\to\mathbb R$ are the numbers
$$M''_k=E\left(\left(\frac{f-E}{\sigma}\right)^k\right)$$
with $\sigma=\sigma(f)$, which satisfy $M''_0=1$, $M''_1=0$, $M''_2=1$.
\end{definition}

To be more precise here, we can define indeed the normalized central moments by the above formula, assuming of course $\sigma>0$, and we have then $M''_0=1$, $M''_1=0$. As for the second normalized central moment, this follows to be 1 indeed, according to:
$$M''_2=E\left(\left(\frac{f-E}{\sigma}\right)^2\right)=\frac{V(f)}{V(f)}=1$$

The normalization made in Definition 2.21 is something quite clever, and of particular interest are the cases $k=3,4$, where we can formulate:

\index{skewness}
\index{kurtosis}

\begin{definition}
The third and fourth normalized central moments
$$\gamma=E\left(\left(\frac{f-E}{\sigma}\right)^3\right)\quad,\quad 
\kappa=E\left(\left(\frac{f-E}{\sigma}\right)^4\right)$$
are called skewness and kurtosis of the variable $f:X\to\mathbb R$.
\end{definition}

And we will end our discussion regarding the moments with this. As a conclusion, a random variable $f:X\to\mathbb R$ is quite reasonably described by the data $(E,\sigma,\gamma,\kappa)$ consisting of the expectation, standard deviation, skewness and kurtosis, which are all quite intuitive quantities, whose knowledge is the same as that of $M_1,M_2,M_3,M_4$. And for more delicate questions, we still have the higher moments, $M_k$ with $k\geq5$.

\section*{2c. Dirac masses, laws} 

Done with the theory? Not really, because we still have to formulate one more key result, which is actually the most important one, in probability, as follows:

\index{law of variable}

\begin{theorem}
Each random variable $f:X\to\mathbb R$ has a law,
$$\mu=\sum_{x\in X}P(x)\delta_{f(x)}$$
linear combination of Dirac masses, regarded as probability measure on $\mathbb R$.
\end{theorem} 

\begin{proof}
Given a random variable $f:X\to\mathbb R$, we can certainly talk about its law $\mu$, as being the formal linear combination of Dirac masses in the statement. Our claim is that this is a probability measure on $\mathbb R$, in the sense of Definition 2.1:

\medskip

(1) Indeed, according to our conventions above, the weight of each point $y\in\mathbb R$ is the following quantity, which is positive, as it should:
$$d\mu(y)=\sum_{f(x)=y}P(x)$$

(2) Moreover, the total mass of this measure is $1$, as it should, due to:
\begin{eqnarray*}
\sum_{y\in\mathbb R}d\mu(y)
&=&\sum_{y\in\mathbb R}\sum_{f(x)=y}P(x)\\
&=&\sum_{x\in X}P(x)\\
&=&1
\end{eqnarray*}

Thus, we are led to the conclusions in the statement.
\end{proof}

Still talking basics, let us record as well the following alternative formula for the law, which is clear from definitions, and that we will often use, in what follows:
$$\mu=\sum_{y\in\mathbb R}P(f=y)\delta_y$$

But you might probably ask, what is the law good for? In answer, many things, and as a first illustration here, the law encodes the sequence of moments, as shown by:

\index{discrete measure}
\index{discrete law}
\index{Dirac mass}
\index{discrete integration}
\index{discrete probability}

\begin{theorem}
We have the following formula, for any $k\in\mathbb N$,
$$E(f^k)=\int_\mathbb Ry^kd\mu(y)$$
with the usual convention that each Dirac mass integrates up to $1$.
\end{theorem}

\begin{proof}
We have indeed the following computation, for the moments:
\begin{eqnarray*}
E(f^k)
&=&\sum_{x\in X}P(x)f(x)^k\\
&=&\sum_{y\in\mathbb R}y^k\sum_{f(x)=y}P(x)\\
&=&\int_\mathbb Ry^kd\mu(y)
\end{eqnarray*}

Thus, we are led to the conclusions in the statement.
\end{proof}

Along the same lines, we have as well the following straightforward generalization of Theorem 2.24, dealing with other functions that can be applied to our variables:

\begin{theorem}
Given a random variable $f:X\to\mathbb R$ having law $\mu$, we have
$$E(\phi(f))=\int_\mathbb R\phi(y)d\mu(y)$$
with the usual convention that each Dirac mass integrates up to $1$.
\end{theorem}

\begin{proof}
We have indeed the following straighforward computation:
\begin{eqnarray*}
E(\phi(f))
&=&\sum_{x\in X}P(x)\phi(f(x))\\
&=&\sum_{y\in\mathbb R}\phi(y)\sum_{f(x)=y}P(x)\\
&=&\int_\mathbb R\phi(y)d\mu(y)
\end{eqnarray*}

Thus, we are led to the conclusion in the statement.
\end{proof}

Getting now to more abstract aspects, the above results are quite interesting, theoretically speaking, and in view of them, it is tempting to formulate:

\index{discrete measure}
\index{Dirac mass}
\index{probability measure}

\begin{definition}
Given a set $X$, which can be finite, countable, or even uncountable, a discrete probability measure on it is a linear combination as follows,
$$\mu=\sum_{x\in X}\lambda_x\delta_x$$
with the coefficients $\lambda_i\in\mathbb R$ satisfying $\lambda_i\geq0$ and $\sum_i\lambda_i=1$. For $f:X\to\mathbb R$ we set
$$\int_Xf(x)d\mu(x)=\sum_{x\in X}\lambda_xf(x)$$
with the convention that each Dirac mass integrates up to $1$.
\end{definition}

Observe that, with this, we are now into pure mathematics. However, and we insist on this, it is basic probability, as developed before, which is behind all this.

\bigskip

Now by staying abstract a bit more, with Definition 2.26 in hand, we can recover our previous basic notions, from Definition 2.1 and from Theorem 2.23, as follows:

\index{push-forward}

\begin{theorem}
With the above notion of discrete probability measure in hand:
\begin{enumerate}
\item A discrete probability space is simply a space $X$, with a discrete probability measure on it $\nu$. In this picture, the probability function is $P(x)=d\nu(x)$.

\item Each random variable $f:X\to\mathbb R$ has a law, which is a discrete probability measure on $\mathbb R$. This law is given by $\mu=f_*\nu$, push-forward of $\nu$ by $f$.
\end{enumerate}
\end{theorem}

\begin{proof}
This might look a bit scary, but is in fact a collection of trivialities, coming straight from definitions, the details being as follows:

\medskip

(1) Nothing much to say in relation with the first assertion, with this being something which is plainly clear, just by comparing Definition 2.1 and Definition 2.26. 

\medskip

(2) As a interesting comment, however, in the general context of Definition 2.26, observe that a probability measure $\mu=\sum_{x\in X}\lambda_x\delta_x$ as there depends only on the following function, called density of our probability measure:
$$\varphi:X\to\mathbb R\quad,\quad \varphi(x)=\lambda_x$$

And, with this notion in hand, our equation $P(x)=d\nu(x)$ simply says that the probability function $P$ is the density of $\nu$. Which is something which is good to know.

\medskip

(3) Pretty much the same story with the first assertion in (2), which this being something clear, just by comparing Theorem 2.23 and Definition 2.26. 

\medskip

(4) As for the very last assertion in (2), consider more generally a probability space $(X,\nu)$, and a function $f:X\to Y$. We can construct then a probability measure $\mu=f_*\nu$ on $Y$, called push-forward of $\nu$ by $f$, according to the following formula:
$$\nu=\sum_{x\in X}\lambda_x\delta_x\implies
\mu=\sum_{y\in Y}\left(\sum_{f(x)=y}\lambda_x\right)\delta_y$$

Alternatively, at the level of the corresponding measures of the parts $Z\subset Y$, we have the following abstract formula, which looks more conceptual:
$$\mu(Z)=\nu(f^{-1}(Z))$$

In any case, one way or another we can talk about push-forward measures $\mu=f_*\nu$, and in the case of a random variable $f:X\to\mathbb R$, we obtain in this way the law of $f$.
\end{proof}

Getting back to more concrete things, and to Theorems 2.23 and 2.24 as stated, quite remarkably, the sequence of moments uniquely determines the law, as shown by:

\begin{theorem}
The sequence of moments of a variable $f:X\to\mathbb R$,
$$M_k=E(f^k)$$
uniquely determines the law of the variable.
\end{theorem}

\begin{proof}
Assume indeed that the law of our variable $f:X\to\mathbb R$ is as follows:
$$\mu=\sum_i\lambda_i\delta_{x_i}$$

The sequence of moments is then given by the following formula:
$$M_k=\sum_i\lambda_ix_i^k$$

But it is then a standard question to recover the numbers $\lambda_i,x_i\in\mathbb R$, and so the measure $\mu$, out of the sequence of numbers $M_k$. Indeed, assuming that the numbers $x_i$ are $0<x_1<\ldots<x_n$ for simplifying, with $k\to\infty$ we have the following estimate:
$$M_k\sim\lambda_nx_n^k$$

Thus, we got the parameters $\lambda_n,x_n\in\mathbb R$ of our measure $\mu$, and by substracting them and doing an obvious recurrence, we get the other parameters $\lambda_i,x_i\in\mathbb R$ as well.
\end{proof}

Moving on, the law of a variable $f:X\to\mathbb R$ is not the only main quantity related to $f$ that we can talk about. As a rival to this notion, we have:

\begin{definition}
Given a random variable $f:X\to\mathbb R$, the function
$$F(x)=P(f\leq x)$$
is called its cumulative distribution function (CDF).
\end{definition}

Observe that the cumulative distribution function $F:\mathbb R\to[0,1]$ must be by definition increasing, and also, that it must have the following limiting properties:
$$F(-\infty)=0\quad,\quad F(\infty)=1$$

In relation with the law of the variable, as previously axiomatized, in Theorem 2.23, the situation is very simple, as shown by the following result:

\begin{theorem}
Assuming that a random variable $f:X\to\mathbb R$ has as law 
$$\mu=\sum_i\lambda_i\delta_{x_i}\quad,\quad x_1<\ldots<x_N\quad,\quad\lambda_i>0,\ \sum_i\lambda_i=1$$
its cumulative distribution function is the following step function:
$$F(x)=\begin{cases}
0&:\ x\in(-\infty,x_1]\\
\lambda_1&:\ x\in(x_1,x_2]\\
\lambda_1+\lambda_2&:\ x\in(x_2,x_3]\\
\vdots&\\
\lambda_1+\ldots+\lambda_{N-1}&:\ x\in(x_{N-1},x_N]\\
1&:\ x\in(x_N,\infty)
\end{cases}$$
Conversely, any such step function $F$ produces a law $\mu$, in the obvious way.
\end{theorem}

\begin{proof}
The first assertion is clear, coming from definitions. As for the converse, given a function $F$ as in the statement, if we set $\lambda_i=\mu_i-\mu_{i-1}$, with the convention $\mu_0=0$, then $\mu=\sum_i\lambda_i\delta_{x_i}$ is a probability measure, having $F$ as CDF, as desired.
\end{proof}

\section*{2d. Independence}

As a continuation of the above, and at a more advanced level, inspired by what happens for coins, dice and cards, let us formulate the following key definition: 

\index{independence}
\index{independent variables}

\begin{definition}
We say that two variables $f,g:X\to\mathbb R$ are independent when
$$P(f=x,g=y)=P(f=x)P(g=y)$$
happens, for any $x,y\in\mathbb R$.
\end{definition}

As mentioned, this is something quite intuitive, inspired by what happens for coins, dice and cards. More on this, and on other examples available, later. For the moment, let us develop some theory. As our first result regarding independence, we have:

\index{mixed moments}

\begin{theorem}
Assuming that $f,g:X\to\mathbb R$ are independent, we have
$$E(f^kg^l)=E(f^k)E(g^l)$$
and the converse holds, in the sense that this implies the independence of $f,g$.
\end{theorem}

\begin{proof}
This is something very standard, the idea being as follows:

\medskip

(1) In one sense, we have the following computation, for the mixed moments:
\begin{eqnarray*}
E(f^kg^l)
&=&\sum_{xy}x^ky^lP(f=x,g=y)\\
&=&\sum_{xy}x^ky^lP(f=x)P(g=y)\\
&=&\sum_xx^kP(f=x)\sum_yy^lP(g=y)\\
&=&E(f^k)E(g^l)
\end{eqnarray*}

(2) Conversely, the moment condition $E(f^kg^l)=E(f^k)E(g^l)$ reformulates as:
$$\sum_{xy}c(x,y)x^ky^l=0\quad,\quad 
c(x,y)=P(f=x,g=y)-P(f=x)P(g=y)$$

So, let us examine this. Our variables $f,g:X\to\mathbb R$ being discrete, let $x_1<\ldots<x_n$ and $y_1<\ldots<y_m$ be their images. With $c_{ij}=c(x_i,x_j)$, our formula reads:
$$\sum_{ij}c_{ij}x_i^ky_j^l=0\quad,\quad\forall k,l\in\mathbb N$$

Assuming now $x_i>0$ for simplifying, we have, with $k\to\infty$, and with $l$ fixed:
$$\sum_{ij}c_{ij}x_i^ky_j^l\sim x_n^k\sum_jc_{nj}y_j^l$$

Thus $\sum_jc_{nj}y_j^l=0$, so $c_{nj}=0$. But with this, again with $k\to\infty$, and $l$ fixed:
$$\sum_{ij}c_{ij}x_i^ky_j^l\sim x_{n-1}^k\sum_jc_{n-1,j}y_j^l$$

And so on, the idea being that we obtain by recurrence $c_{ij}=0$, as desired.
\end{proof}

In order to talk now about sums of independent variables, we will need:

\index{convolution}
\index{semigroup}
\index{group}

\begin{definition}
Given a space $X$ with a sum operation $+$, we can define the convolution of any two discrete probability measures on it,
$$\mu=\sum_ia_i\delta_{x_i}\quad,\quad \nu=\sum_jb_j\delta_{y_j}$$
as being the discrete probability measure given by the following formula:
$$\mu*\nu=\sum_{ij}a_ib_j\delta_{x_i+y_j}$$
That is, the convolution operation $*$ is defined by $\delta_x*\delta_y=\delta_{x+y}$, and linearity.
\end{definition}

Observe that our $*$ operation is indeed well-defined, the weights $a_ib_j$ being positive, and summing up to 1. In relation now with independence, we have:

\begin{theorem}
Assuming that $f,g:X\to\mathbb R$ are independent, we have
$$\mu_{f+g}=\mu_f*\mu_g$$
where $*$ is the convolution of real probability measures.
\end{theorem}

\begin{proof}
We have indeed the following straightforward computation:
\begin{eqnarray*}
\mu_{f+g}
&=&\sum_{x\in\mathbb R}P(f+g=x)\delta_x\\
&=&\sum_{y,z\in\mathbb R}P(f=y,g=z)\delta_{y+z}\\
&=&\sum_{y,z\in\mathbb R}P(f=y)P(g=z)\delta_y*\delta_z\\
&=&\left(\sum_{y\in\mathbb R}P(f=y)\delta_y\right)*\left(\sum_{z\in\mathbb R}P(g=z)\delta_z\right)\\
&=&\mu_f*\mu_g
\end{eqnarray*}

Thus, we are led to the conclusion in the statement.
\end{proof}

Time for some illustrations, and applications, of all this? Let us start with:

\begin{proposition}
When flipping a biased coin, with the probabilities being
$$P({\rm heads})=p\quad,\quad P({\rm tails})=1-p$$
with $p\in[0,1]$, your winning law, for heads, is the following probability measure,
$$b_p=(1-p)\delta_0+p\delta_1$$
called Bernoulli law of parameter $p\in[0,1]$.
\end{proposition}

\begin{proof}
This is indeed something trivial, more standing as a definition for $b_p$.
\end{proof}

Coming next, we have the following result, which is something clear too:

\index{binomial law}
\index{iterated coin throws}

\begin{proposition}
When flipping a biased coin $n$ times in a row, we have
$$P(s\ {\rm heads}, n-s\ {\rm tails})=\binom{n}{s}p^s(1-p)^{n-s}$$
and your winning law, for heads, is the following probability measure,
$$b_{np}=\sum_{s=0}^n\binom{n}{s}p^s(1-p)^{n-s}\delta_s$$
called binomial law of parameter $p\in[0,1]$.
\end{proposition}

\begin{proof}
This is again clear, because in order to have $s$ heads and $n-s$ tails, we must select the $s$ occurrences of heads, among the $n$ total throws, and there are $\binom{n}{s}$ choices here, and then sum over these $\binom{n}{s}$ choices the common quantity $p^s(1-p)^{n-s}$.
\end{proof}

Here is now our theorem about this, using our independence technology:

\begin{theorem}
The following happen, in the context of the biased coin game:
\begin{enumerate}
\item The Bernoulli laws $b_p$ produce the binomial laws $b_{np}$, by iterating the coin game $n\in\mathbb N$ times, via the independence of the throws.

\item We have in fact $b_{np}=b_p^{*n}$, with $*$ being the convolution operation for the real probability measures, given by $\delta_x*\delta_y=\delta_{x+y}$, and linearity.
\end{enumerate}
\end{theorem}

\begin{proof}
By using the binomial formula, with numbers formally replaced by Dirac masses, which is something that we can do, we have indeed the following computation:
\begin{eqnarray*}
b_p^{*n}
&=&\big((1-p)\delta_0+p\delta_1\big)^{*n}\\
&=&\sum_{s=0}^np^s(1-p)^{n-s}\binom{n}{s}\delta_0^{*(n-s)}*\delta_1^{*s}\\
&=&\sum_{s=0}^np^s(1-p)^{n-s}\binom{n}{s}\delta_s\\
&=&b_{np}
\end{eqnarray*}

Thus, up to the general the independence theory developed above, and with everything here doing well, as desired, we are led to the various conclusions in the statement.
\end{proof}

We have as well an interesting application of independence to dice, as follows:

\begin{theorem}
When rolling a die $n$ times, the distribution of the sum is
$$P(s)=\frac{1}{6^n}\sum_{k=0}^{[\frac{s-n}{6}]}(-1)^k\binom{n}{k}\binom{s-6k-1}{n-1}$$
with this coming from independence, and the binomial formula for negative exponents.
\end{theorem}

\begin{proof}
This is something quite tricky, the idea being as follows:

\medskip

(1) According to our independence theory above, the law of the sum is:
$$\mu_n=\frac{1}{6^n}(\delta_1+\delta_2+\delta_3+\delta_4+\delta_5+\delta_6)^{*n}$$

(2) Equivalently, when thinking a bit, we have the following formula for the law of the sum, and I will leave it to you, to clarify the details here:
$$\sum_{s=n}^{6n}P(s)x^s=\frac{1}{6^n}(x+x^2+x^3+x^4+x^5+x^6)^n$$

(3) So, let us compute the function on the right. This is given by:
\begin{eqnarray*}
f(x)
&=&\frac{1}{6^n}(x+x^2+x^3+x^4+x^5+x^6)^n\\
&=&\frac{x^n}{6^n}(1+x+x^2+x^3+x^4+x^5)^n\\
&=&\frac{x^n}{6^n}\cdot\frac{(1-x^6)^n}{(1-x)^n}
\end{eqnarray*}

(4) By using the usual binomial formula for the numerator $(1-x^6)^n$, and the generalized  binomial formula from chapter 1 for the fraction $1/(1-x)^n$, we obtain:
\begin{eqnarray*}
f(x)
&=&\frac{x^n}{6^n}\sum_{k=0}^n(-1)^k\binom{n}{k}x^{6k}\sum_{l=0}^\infty\binom{n+l-1}{n-1}x^l\\
&=&\frac{1}{6^n}\sum_{k=0}^n\sum_{l=0}^\infty(-1)^k\binom{n}{k}\binom{n+l-1}{n-1}x^{n+6k+l}\\
&=&\frac{1}{6^n}\sum_{s=0}^\infty\sum_{k=0}^n(-1)^k\binom{n}{k}\binom{s-6k-1}{n-1}x^s
\end{eqnarray*}

(5) Now remember that $f$ was a polynomial of degree $6n$, divisible by $x^n$. Thus, the true range of the first summing index is $s=n,\ldots,6n$. Also, from $s=n+6k+l$ in the above computation, we get $s\geq n+6k$, and so $k\leq[(s-n)/6]$, integer part. Thus, our formula in (4), written by ignoring zero coefficients appearing from cancellation, is:
$$f(x)=\frac{1}{6^n}\sum_{s=n}^{6n}\sum_{k=0}^{[\frac{s-n}{6}]}(-1)^k\binom{n}{k}\binom{s-6k-1}{n-1}x^s$$

But with this, we are led via (2) to the formula in the statement.
\end{proof} 

As a last thing to be discussed, in relation with independence, we have:

\begin{theorem}
Given variables $f,g:X\to\mathbb R$, we can talk about their covariance
\begin{eqnarray*}
cov(f,g)
&=&E[(f-E(f))(g-E(g))]\\
&=&E(fg)-E(f)E(g)
\end{eqnarray*}
with this quantity having the following properties:
\begin{enumerate}
\item $cov(f,f)=V(f)$.

\item If $f,g$ are independent, $cov(f,g)=0$.

\item The converse of this latter fact is, unfortunately, not true.
\end{enumerate}
\end{theorem}

\begin{proof}
This is something very standard. To start with, we have the following computation, which shows that our two definitions of the covariance are equivalent:
\begin{eqnarray*}
&&E[(f-E(f))(g-E(g))]\\
&=&E[fg-E(f)g-fE(g)+E(f)E(g)]\\
&=&E(fg)-E(f)E(g)-E(f)E(g)+E(f)E(g)\\
&=&E(fg)-E(f)E(g)
\end{eqnarray*}

As for the various assertions in the statement, their proof goes as follows:

\medskip

(1) This is indeed clear, with either of the definitions for the variance $V(f)$.

\medskip

(2) This is clear too, coming from the multiplicativity formula from Theorem 2.32.

\medskip

(3) To start with, the converse of (2) has obviously zero chances to hold, because we perfectly know from Theorem 2.32 that the independence of $f,g:X\to\mathbb R$ amounts in having the following formula, for the mixed moments, meaning for any $k,l\in\mathbb N$:
$$E(f^kg^l)=E(f^k)E(g^l)$$

(4) This being said, an explicit counterexample will certainly not hurt. Consider two Bernoulli variables $f,g:X\to\mathbb R$, taking their values according to the following table:
\vskip-7mm
$$\xymatrix@R=1pt@C=4pt{
&&\ar@{-}[ddddddd]&&&\ar@{-}[ddddddd]\\
&g\backslash f&\ar@{-}[dd]&a&b\\
\ar@{-}[rrrrrrr]&&&&&&&\\
&c&&pq+x&(1-p)q-x&&q\\
&d&&p(1-q)-x&(1-p)(1-q)+x&&1-q\\
\ar@{-}[rrrrrrr]&&&&&&&\\
&&&p&1-p&\\
&&&&&}$$
\vskip-2mm

The covariance of $f,g$ can be then computed as follows:
\begin{eqnarray*}
&&cov(f,g)\\
&=&E(fg)-E(f)E(g)\\
&=&[ac(pq+x)+ad(p(1-q)-x)+bc((1-p)q-x)+bd((1-p)(1-q)+x)]\\
&&-[ap+b(1-p)][cq+d(1-q)]\\
&=&acx-adx-bcx+bdx\\
&=&(a-b)(c-d)x
\end{eqnarray*}

Which is bad news, because $cov(f,g)=0$ means $x=0$, which in turn means that $f,g$ must be independent. Damn, so instead of our counterexample, we have a theorem here, stating that the converse of (2) holds for the Bernoulli variables. Good to know.

\medskip

(5) Counterexample, take two. In view of the above, consider the simplest possible situation which is not (Bernoulli, Bernoulli), which is as follows, with $f$ following a uniform trivalent distribution on $[-1,1]$, and with $g$ being taken to be its square, $g=f^2$: 
\vskip-7mm
$$\xymatrix@R=1pt@C=6pt{
&&\ar@{-}[ddddddd]&&&&\ar@{-}[ddddddd]\\
&g\backslash f&\ar@{-}[dd]&-1&0&1\\
\ar@{-}[rrrrrrrr]&&&&&&&&\\
&0&&0&1/3&0&&1/3\\
&1&&1/3&0&1/3&&2/3\\
\ar@{-}[rrrrrrrr]&&&&&&&&\\
&&&1/3&1/3&1/3\\
&&&&&&}$$
\vskip-2mm

Our variables are then not independent, because the upper left 0 does not come as $0=1/3\cdot1/3$. On the other hand we have $E(f)=0$, $E(g)=1/3$, $E(fg)=0$, so the covariance is $cov(f,g)=0-0=0$, and we have our counterexample, as desired.
\end{proof}

And we will end our preliminary introduction to probability theory with this. All the above was of course quite quick, and for a more detailed treatement of all this, including more phenomenology, numerics, tables and so on, you can check my previous book \cite{ba1}, where all this material is explained over a full 6 chapters, worth 150 pages.

\section*{2e. Exercises}

This was a quite crowded theoretical chapter, and as exercises, we have:

\begin{exercise}
Further develop the symmetry approach, for coins and dice.
\end{exercise}

\begin{exercise}
Further develop the additive approach, for coins and dice.
\end{exercise}

\begin{exercise}
Further develop the algebraic approach, for coins and dice.
\end{exercise}

\begin{exercise}
Work out some further examples of expectations.
\end{exercise}

\begin{exercise}
Work out some further examples of variances.
\end{exercise}

\begin{exercise}
Learn the various versions of the Cauchy-Schwarz inequality.
\end{exercise}

\begin{exercise}
Clarify what we said, in relation with push-forward measures.
\end{exercise}

\begin{exercise}
Compute more covariances, for pairs of variables of your choice.
\end{exercise}

As bonus exercise, have a look at measure theory, which is related to all this.

\chapter{Binomial laws}

\section*{3a. Binomial laws}

We already met the binomial laws in chapter 2, and time now to get into a more detailed study of these laws, and their various versions. Let us start our study with:

\index{binomial law}
\index{Bernoulli law}
\index{biased coin}

\begin{theorem}
When playing with a biased coin, your winning law is
$$b_p=(1-p)\delta_0+p\delta_1$$
called Bernoulli law of parameter $p\in[0,1]$. When playing $n$ times, your winning law is
$$b_{np}=\sum_{s=0}^n\binom{n}{s}p^s(1-p)^{n-s}\delta_s$$
called binomial law of parameter $p\in[0,1]$. We have the formula $b_{np}=b_p^{*n}$.
\end{theorem}

\begin{proof}
This is something that we know from chapter 2, the idea being as follows:

\medskip

(1) The first assertion, regarding the Bernoulli law, is something which is plainly trivial, and with this standing more for a definition, for this law.

\medskip

(2) When playing $n$ times, in order to have $s$ wins, we must select the $s$ occurrences of heads, among the $n$ throws, and there are $\binom{n}{s}$ choices here. Then, we must sum over these $\binom{n}{s}$ choices the common probability $p^s(1-p)^{n-s}$, for $s$ heads and $n-s$ tails.
\medskip

(3) The formula $b_{np}=b_p^{*n}$ comes from the fact that the Bernoulli laws $b_p$ produce the binomial laws $b_{np}$, by iterating the game $n$ times, via the independence of the throws. 

\medskip

(4) Equivalently, this latter formula comes from the following formal computation:
\begin{eqnarray*}
b_p^{*n}
&=&\big((1-p)\delta_0+p\delta_1\big)^{*n}\\
&=&\sum_{s=0}^np^s(1-p)^{n-s}\binom{n}{s}\delta_0^{*(n-s)}*\delta_1^{*s}\\
&=&\sum_{s=0}^np^s(1-p)^{n-s}\binom{n}{s}\delta_s\\
&=&b_{np}
\end{eqnarray*}

Thus, we are led to the conclusions in the statement.
\end{proof}

Let us do now some computations. As a first concrete result about the binomial laws, regarding their mean and variance, as introduced in chapter 2, we have:

\index{mean of binomial law}
\index{variance of binomial law}
\index{moments of binomial law}

\begin{theorem}
The binomial law $b_{np}$ has the following properties:
\begin{enumerate}
\item The mean is $E=np$.

\item The variance is $V=np(1-p)$.
\end{enumerate}
\end{theorem}

\begin{proof}
This is very standard, by using the binomial formula, as follows:

\medskip

(1) Regarding the mean of a binomial variable $f:X\to\mathbb R$, the computation, based on a suitable application of the binomial formula, is as follows:
\begin{eqnarray*}
E(f)
&=&\sum_{s=1}^ns\binom{n}{s}p^s(1-p)^{n-s}\\
&=&\sum_{s=1}^n\frac{n!}{(s-1)!(n-s)!}p^s(1-p)^{n-s}\\
&=&np\sum_{s=1}^n\frac{(n-1)!}{(s-1)!(n-s)!}p^{s-1}(1-p)^{n-s}\\
&=&np\sum_{t=0}^{n-1}\binom{n-1}{t}p^t(1-p)^{n-t-1}\\
&=&np(p+1-p)^{n-1}\\
&=&np
\end{eqnarray*}

(2) Coming next, with the same trick, we can compute the following quantity:
\begin{eqnarray*}
E(f^2)-E(f)
&=&\sum_{s=2}^n(s^2-s)\binom{n}{s}p^s(1-p)^{n-s}\\
&=&\sum_{s=2}^n\frac{k!}{(s-2)!(n-s)!}p^s(1-p)^{n-s}\\
&=&n(n-1)p^2\sum_{s=2}^n\frac{(n-2)!}{(s-2)!(n-s)!}p^{s-2}(1-p)^{n-s}\\
&=&n(n-1)p^2\sum_{t=0}^{n-2}\binom{n-2}{t}p^t(1-p)^{n-t-2}\\
&=&n(n-1)p^2(p+1-p)^{n-2}\\
&=&n(n-1)p^2
\end{eqnarray*}

(3) Now by using the above formulae, we conclude that the variance of a binomial variable $V(f)=E(f^2)-E(f)^2$ is given by the following formula:
\begin{eqnarray*}
V(f)
&=&[E(f^2)-E(f)]+E(f)-E(f)^2\\
&=&n(n-1)p^2+np-(np)^2\\
&=&np(1-p)
\end{eqnarray*}

We are therefore led to the conclusions in the statement.
\end{proof}

Regarding now the higher moments, things are quite similar, and we have:

\begin{theorem}
The moments of the binomial law $b_{np}$ are given by
$$M_1=np$$
$$M_2=n(n-1)p^2+np$$
$$M_3=n(n-1)(n-2)p^3+3n(n-1)p^2+np$$
$$M_4=n(n-1)(n-2)(n-3)p^4+6n(n-1)(n-2)p^3+7n(n-1)p^2+np$$
$$\vdots$$
with the computation in general being possible by recurrence.
\end{theorem}

\begin{proof}
The moments of the binomial law $b_{np}$ of parameters $n\in\mathbb N$ and $p\in[0,1]$ are by definition given by the following formula:
$$M_k=\sum_{s=1}^ns^k\binom{n}{s}p^s(1-p)^{n-s}$$

We already know the first two formulae from Theorem 3.2 and its proof, and the continuation is by using the same trick as there, as follows:

\medskip

(1) Regarding the third moment, we can use here the following formula, that you can come upon by some reverse engineering, and I will leave this to you, as an exercise:
\begin{eqnarray*}
s(s-1)(s-2)+3s(s-1)+s
&=&s[(s-1)(s-2)+3(s-1)+1]\\
&=&s(s^2-3s+2+3s-3+1)\\
&=&s^3
\end{eqnarray*}

Indeed, this gives the following formula for the third moment:
\begin{eqnarray*}
M_3
&=&\sum_{s=1}^ns^3\binom{n}{s}p^s(1-p)^{n-s}\\
&=&\sum_{s=1}^n\left[s(s-1)(s-2)+3s(s-1)+s\right]\binom{n}{s}p^s(1-p)^{n-s}\\
&=&\sum_{s=3}^ns(s-1)(s-2)\binom{n}{s}p^s(1-p)^{n-s}\\
&&+3\sum_{s=2}^ns(s-1)\binom{n}{s}p^s(1-p)^{n-s}\\
&&+\sum_{s=1}^ns\binom{n}{s}p^s(1-p)^{n-s}\\
&=&n(n-1)(n-2)\sum_{s=3}^n\binom{n-3}{s-3}p^s(1-p)^{n-s}\\
&&+3n(n-1)\sum_{s=2}^n\binom{n-2}{s-2}p^s(1-p)^{n-s}\\
&&+n\sum_{s=1}^n\binom{n-1}{s-1}p^s(1-p)^{n-s}\\
&=&n(n-1)(n-2)p^3+3n(n-1)p^2+np
\end{eqnarray*}

(2) Regarding the fourth moment, we can use here the following formula, which again can be discovered via some routine work, say good exercise for you:
\begin{eqnarray*}
&&s(s-1)(s-2)(s-3)+6s(s-1)(s-2)+7s(s-1)+s\\
&=&s[(s-1)(s^2-5s+6+6s-12+7)+1]\\
&=&s[(s-1)(s^2+s+1)+1]\\
&=&s[s^3-1+1]\\
&=&s^4
\end{eqnarray*}

But this gives the formula in the statement for the fourth moment, via a computation which is very similar to the one in (1) above, again good exercise for you.

\medskip

(3) As a continuation of this, it is quite clear that we can compute $M_k$ by recurrence. We will be back to this later, with more results on the subject.
\end{proof}

Let us record as well a result about the corresponding central moments:

\index{central moments}

\begin{theorem}
The central moments of the binomial law are
$$M_1'=0$$
$$M_2'=np(1-p)$$
$$M_3'=np(1-p)(1-2p)$$
$$M_4'=np(1-p)(1+(3n-6)p(1-p))$$
$$\vdots$$
with the computation in general being possible by recurrence.
\end{theorem}

\begin{proof}
We recall from chapter 2 that the central moments of a variable $f:X\to\mathbb R$, designed for having $M_1'=0$, $M_2'=V$, are the following numbers, with $E=E(f)$:
$$M_k'=E((f-E)^k)$$

Thus $M_1'=0$, and the formula in the statement for $M_2'=V$ is the one from Theorem 3.2. Regarding now the higher moments, by the binomial formula, we have:
\begin{eqnarray*}
M_k'
&=&E\left(\sum_{r=0}^k(-1)^r\binom{k}{r}E^rf^{k-r}\right)\\
&=&\sum_{r=0}^k(-1)^r\binom{k}{r}E^rM_{k-r}\\
&=&\sum_{r=0}^{k-2}(-1)^r\binom{k}{r}E^rM_{k-r}+(-1)^{k-1}(k-1)E^k
\end{eqnarray*}

Now with $k=3$, by using the formulae from Theorem 3.3, we obtain:
\begin{eqnarray*}
M_3'
&=&M_3-3EM_2+2E^3\\
&=&[n(n-1)(n-2)p^3+3n(n-1)p^2+np]-3np(n(n-1)p^2+np)+2(np)^3\\
&=&np(1-p)(1-2p)
\end{eqnarray*}

Similarly, at $k=4$, again by using the formulae from Theorem 3.3, we have:
\begin{eqnarray*}
&&M_4'\\
&=&M_4-4EM_3+6E^2M_2-3E^4\\
&=&n(n-1)(n-2)(n-3)p^4+6n(n-1)(n-2)p^3+7n(n-1)p^2+np\\
&&\!\!-4np[n(n-1)(n-2)p^3+3n(n-1)p^2+np]+6(np)^2(n(n-1)p^2+np)-3(np)^4\\
&=&np(1-p)(1+(3n-6)p(1-p))
\end{eqnarray*}

As for the last assertion, we will be back to it later, with details.
\end{proof}

Let us end this discussion with the following result, capturing the essentials:

\index{mean}
\index{variance}
\index{skewness}
\index{kurtosis}

\begin{theorem}
The mean, variance, skewness and kurtosis of $b_{np}$ are
$$E=np\quad,\quad V=npq\quad,\quad\gamma=\frac{q-p}{\sqrt{npq}}\quad,\quad\kappa=3+\frac{1-6pq}{npq}$$
with the convention $q=1-p$.
\end{theorem}

\begin{proof}
We know the first two formulae, from Theorem 3.2. Regarding the last two formulae, we can use here the central moment formulae from Theorem 3.4. Observe first that with the convention $q=1-p$, our formulae for $M_3',M_4'$ there read:
$$M_3'=npq(q-p)\quad,\quad M_4'=npq(1+(3n-6)pq)$$

Now by dividing respectively by $\sigma^3=V\sqrt{V}$ and $\sigma^4=V^2$, we obtain the formulae in the statement for $\gamma$ and $\kappa$. Needless to say, the higher normalized central moments can be computed too, and we will back to this later, with some explicit results.
\end{proof}

In order to discuss now a more conceptual approach to the computation of the moments of the binomial laws, let us start with something seemingly unrelated, namely:

\index{partitions}
\index{partitions of set}

\begin{definition}
We denote by $P(k)$ the set of partitions of $\{1,\ldots,k\}$, with these partitions $\pi\in P(k)$ being most conveniently being drawn as diagrams,
$$
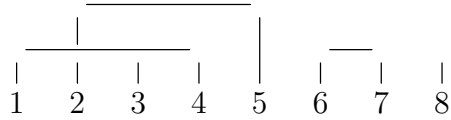
\xymatrix@R=10pt@C=10pt{
&\ar@{-}[rrr]&&&&&&\\
\ar@{-}[rrr]&\ar@{-}[u]&&&&\ar@{-}[r]&&\\
1\ar@{-}[u]&2\ar@{-}[u]&3\ar@{-}[u]&4\ar@{-}[u]&5\ar@{-}[uu]&6\ar@{-}[u]&7\ar@{-}[u]&8\ar@{-}[u]
}$$
with the strings joining the numbers belonging to the same block of $\pi$. That is, the above diagram represents the partition $\{1,\ldots,8\}=\{1,3,4\}\cup\{2,5\}\cup\{6,7\}\cup\{8\}$.
\end{definition}

Observe that there is a bit of care to be taken with this convention, with respect to the crossings. We can either proceed as above, with the $\{2,5\}$ block being respresented ``under'' the block $\{1,3,4\}$, or use different types of strings, as for instance:
$$
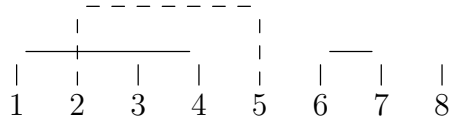
\xymatrix@R=10pt@C=10pt{
&\ar@{--}[rrr]&&&&&&\\
\ar@{-}[rrr]&&&&&\ar@{-}[r]&&\\
1\ar@{-}[u]&2\ar@{--}[uu]&3\ar@{-}[u]&4\ar@{-}[u]&5\ar@{--}[uu]&6\ar@{-}[u]&7\ar@{-}[u]&8\ar@{-}[u]
}$$

Both conventions are good, and we will be mostly using here the first one, that from Definition 3.6, which in practice leads to quicker drawings.

\bigskip

Now back to the binomial laws, we have the following result, about them:

\index{number of blocks}

\begin{theorem}
The moments of the binomial law $b_{np}$ are given by
$$M_k=\sum_{\pi\in P(k)}\frac{n!}{(n-|\pi|)!}\,p^{|\pi|}$$
with $|.|$ being the number of blocks.
\end{theorem}

\begin{proof}
This is something very standard, the idea being as follows:

\medskip

(1) Some numerics first. At $k=1$ we only have one partition, namely $|$\,, and the formula in the statement holds indeed, as shown by the following computation:
\begin{eqnarray*}
\sum_{\pi\in P(1)}\frac{n!}{(n-|\pi|)!}\,p^{|\pi|}
&=&\frac{n!}{(n-|{\ }_|\;|)!}\,p^{|{\ }_|\;|}\\
&=&\frac{n!}{(n-1)!}\,p^1\\
&=&np
\end{eqnarray*}

(2) At $k=2$ now, we have two partitions, namely $|\,|$ and $\sqcap$, and the formula in the statement holds again, as shown by the following computation:
\begin{eqnarray*}
\sum_{\pi\in P(2)}\frac{n!}{(n-|\pi|)!}\,p^{|\pi|}
&=&\frac{n!}{(n-|{\ }_{|\,|}\;|)!}\,p^{|{\ }_{|\,|}\;|}+\frac{n!}{(n-|\sqcap|)!}\,p^{|\sqcap|}\\
&=&\frac{n!}{(n-2)!}\,p^2+\frac{n!}{(n-1)!}\,p^1\\
&=&n(n-1)p^2+np
\end{eqnarray*}

(3) At $k=3$ we have 5 partitions, namely $|\ |\ |$\,, $\sqcap\ |$\,, $\sqcap\hskip-3.2mm{\ }_|$\,\,, $|\ \sqcap$\,, $\sqcap\hskip-0.7mm\sqcap$, and the formula in the statement holds again, as shown by the following computation:
\begin{eqnarray*}
\sum_{\pi\in P(3)}\frac{n!}{(n-|\pi|)!}\,p^{|\pi|}
&=&\frac{n!}{(n-|{\ }_{|\,|\,|}\;|)!}\,p^{|{\ }_{|\,|\,|}\;|}+\ldots+\frac{n!}{(n-|\sqcap\hskip-1.6mm\sqcap|)!}\,p^{|\sqcap\hskip-0.57mm\sqcap|}\\
&=&\frac{n!}{(n-3)!}\,p^3+3\cdot\frac{n!}{(n-2)!}\,p^2+\frac{n!}{(n-1)!}\,p^1\\
&=&n(n-1)(n-2)p^3+3n(n-1)p^2+np
\end{eqnarray*}

(4) At $k=4$ we have 15 partitions, with the list of relevant partitions, and their multiplicities up to permutations, which is self-explanatory, being as follows:
$$|\ |\ |\ |\times1\quad,\quad \sqcap\ |\ |\times 6\quad,\quad \sqcap\sqcap\times 3
\quad,\quad \sqcap\hskip-1.6mm\sqcap|\times4\quad,\quad \sqcap\hskip-1.6mm\sqcap\hskip-1.6mm\sqcap\times1$$

Now by doing the computation, as before, the formula in the statement holds again:
\begin{eqnarray*}
&&\sum_{\pi\in P(4)}\frac{n!}{(n-|\pi|)!}\,p^{|\pi|}\\
&=&\frac{n!}{(n-|{\ }_{|\,|\,|\,|}\;|)!}\,p^{|{\ }_{|\,|\,|\,|}\;|}+\ldots\ldots+\frac{n!}{(n-|\sqcap\hskip-1.6mm\sqcap\hskip-1.6mm\sqcap|)!}\,p^{|\sqcap\hskip-0.57mm\sqcap\hskip-0.57mm\sqcap|}\\
&=&\frac{n!}{(n-4)!}\,p^4+6\cdot\frac{n!}{(n-3)!}\,p^3+(3+4)\cdot\frac{n!}{(n-2)!}\,p^2+\frac{n!}{(n-1)!}\,p^1\\
&=&n(n-1)(n-2)(n-3)p^4+6n(n-1)(n-2)p^3+7n(n-1)p^2+np
\end{eqnarray*}

(5) In general now, in order to prove the result, let us go back to our computations from the proof of Theorem 3.3. The idea there was that of using a formula as follows:
\begin{eqnarray*}
s^k
&=&s(s-1)\ldots(s-k+1)\\
&&+c_1s(s-1)\ldots (s-k+2)\\
&&\vdots\\
&&+c_{k-2}s(s-1)\\
&&+s
\end{eqnarray*}

(6) Indeed, assuming that we have such a magic formula, we obtain:
\begin{eqnarray*}
M_k
&=&\sum_{s=k}^ns(s-1)\ldots(s-k+1)\binom{n}{s}p^s(1-p)^{n-s}\\
&&+c_1\sum_{s=k-1}^ns(s-1)\ldots (s-k+2)\binom{n}{s}p^s(1-p)^{n-s}\\
&&\vdots\\
&&+c_{k-2}\sum_{s=2}^ns(s-1)\binom{n}{s}p^s(1-p)^{n-s}\\
&&+\sum_{s=1}^ns\binom{n}{s}p^s(1-p)^{n-s}\\
&=&n(n-1)\ldots(n-k+1)\sum_{s=k}^n\binom{n-k}{s-k}p^s(1-p)^{n-s}\\
&&+c_1n(n-1)\ldots (n-k+2)\sum_{s=k-1}^n\binom{n-k+1}{s-k+1}p^s(1-p)^{n-s}\\
&&\vdots\\
&&+c_{k-2}n(n-1)\sum_{s=2}^n\binom{n-2}{s-2}p^s(1-p)^{n-s}\\
&&+n\sum_{s=1}^n\binom{n-1}{s-1}p^s(1-p)^{n-s}\\
&
\end{eqnarray*}

(7) Thus, we are led to the following formula, for the $k$-th moment:
\begin{eqnarray*}
M_k
&=&n(n-1)\ldots(n-k+1)p^k\\
&&+c_1n(n-1)\ldots (n-k+2)p^{k-1}\\
&&\vdots\\
&&+c_{k-2}n(n-1)p^2\\
&&+np
\end{eqnarray*}

(8) The problem is now, how to find that magic formula from (5)? And the answer here comes from partitions. Indeed, with the standard convention that for a multi-index $i=(i_1,\ldots,i_k)$, its kernel is the partition $\ker i\in P(k)$ joining equal indices, we have:
\begin{eqnarray*}
s^k
&=&\sum_{i_1=1}^s\ldots\sum_{i_k=1}^s1\\
&=&\sum_{\pi\in P(k)}\#\left\{(i_1,\ldots,i_k)\in\{1,\ldots,s\}^k\Big|\ker i=\pi\right\}\\
&=&\sum_{\pi\in P(k)}\frac{s!}{(s-|\pi|)!}
\end{eqnarray*}

(9) Thus, good news, we have the magic formula that we need as data in (5), and according now to (7), we are led to the following formula, for the moment:
$$M_k=\sum_{\pi\in P(k)}\frac{n!}{(n-|\pi|)!}\,p^{|\pi|}$$

Thus, we have indeed the formula in the statement.
\end{proof}

\section*{3b. Negative binomials}

By using the generalized binomial formula established in chapter 1, at negative integer values of the exponent, we can formulate the following definition:

\index{Pascal law}
\index{negative binomial law}

\begin{definition}
The negative binomial law $c_{mp}$, with $m\in\mathbb N$ and $p\in[0,1]$, is:
$$P(s)=\binom{s+m-1}{s}(1-p)^sp^m$$
We also call this Pascal law of parameters $m\in\mathbb N$ and $p\in[0,1]$.
\end{definition}

The relation with the negative binomial coefficients comes from:
$$\binom{s+m-1}{s}=(-1)^s\binom{-m}{s}$$

Indeed, by using this formula, we can write our law as follows:
$$P(s)=(-1)^s\binom{-m}{s}(1-p)^sp^m=\binom{-m}{s}(p-1)^sp^m$$

As a consequence of this, we have indeed a probability measure, as shown by:
\begin{eqnarray*}
\sum_{s\geq0}P(s)
&=&\sum_{s\geq0}\binom{-m}{s}(p-1)^sp^m\\
&=&p^m\sum_{s\geq0}\binom{-m}{s}(p-1)^s\\
&=&p^m[1+(p-1)]^{-m}\\
&=&p^mp^{-m}\\
&=&1
\end{eqnarray*}

Let us record this finding, which is something quite interesting, as follows:

\index{generalized binomial coefficient}

\begin{proposition}
The negative binomial law $c_{mp}$ is given by
$$P(s)=\binom{-m}{s}(p-1)^sp^m$$
with $\binom{-m}{s}$ being a generalized binomial coefficient.
\end{proposition}

\begin{proof}
This follows indeed from the above discussion.
\end{proof}

In practice now, the negative binomial laws appear quite naturally, as follows:

\begin{theorem}
When flipping a biased coin until reaching to $m$ heads,
$$P(s\ {\rm tails})=\binom{s+m-1}{s}(1-p)^sp^m$$
and so, $c_{mp}$ is the law of the number of failures, when playing this game.
\end{theorem}

\begin{proof}
When computing $P(s\ {\rm tails})$, we certainly have the $(1-p)^sp^m$ factor, coming from the $m$ heads and $s$ tails. As for the multiplicity, this is the binomial coefficient in the statement, coming from choosing the positions of the $s$ tails, among the total of $s+m-1$ attempts, up to the last one, which does not count, because it must be heads. 
\end{proof}

Regarding now the basic properties of the negative binomial laws, we have:

\begin{theorem}
The negative binomial law $c_{mp}$ has the following properties:
\begin{enumerate}
\item The mean is $E=\frac{m(1-p)}{p}$.

\item The variance is $V=\frac{m(1-p)}{p^2}$.
\end{enumerate}
\end{theorem}

\begin{proof}
Regarding the mean, this can be computed as follows:
\begin{eqnarray*}
E
&=&\sum_{s\geq1}sP(s)\\
&=&\sum_{s\geq1}s\binom{-m}{s}(p-1)^sp^m\\
&=&\sum_{s\geq1}-m\binom{-m-1}{s-1}(p-1)^sp^m\\
&=&-mp^m(p-1)\sum_{s\geq1}\binom{-m-1}{s-1}(p-1)^{s-1}\\
&=&mp^m(1-p)[1+(p-1)]^{-m-1}\\
&=&mp^m(1-p)p^{-m-1}\\
&=&\frac{m(1-p)}{p}
\end{eqnarray*}

With a similar trick, we can compute the difference $M_2-M_1$, as follows:
\begin{eqnarray*}
M_2-M_1
&=&\sum_{s\geq2}(s^2-s)P(s)\\
&=&\sum_{s\geq2}s(s-1)\binom{-m}{s}(p-1)^sp^m\\
&=&\sum_{s\geq2}m(m+1)\binom{-m-2}{s-2}(p-1)^sp^m\\
&=&m(m+1)p^m(p-1)^2\sum_{s\geq2}\binom{-m-2}{s-2}(p-1)^{s-2}\\
&=&m(m+1)p^m(1-p)^2[1+(p-1)]^{-m-2}\\
&=&m(m+1)p^m(1-p)^2p^{-m-2}\\
&=&m(m+1)\,\frac{(1-p)^2}{p^2}
\end{eqnarray*}

We conclude that the second moment is given by the following formula:
$$M_2
=m\,\frac{1-p}{p^2}(1+m-mp)$$

Thus, we are led to the variance formula in the statement.
\end{proof}

Moving on, of particular interest is the exponent $m=1$, where we have:

\index{geometric law}

\begin{proposition}
The negative binomial laws $c_p=c_{1p}$ of exponent $m=1$ are
$$P(s)=(1-p)^sp$$
with these being called geometric laws.
\end{proposition}

\begin{proof}
This is indeed clear from definitions, and with the name of these laws coming from the geometric series which produces them, namely $\sum_{s\geq0}(1-p)^sp=1$.
\end{proof}

Regarding now the basic properties of the geometric laws, these are as follows:

\begin{proposition}
The geometric law $c_p$ has the following properties:
\begin{enumerate}
\item The mean is $E=\frac{(1-p)}{p}$.

\item The variance is $V=\frac{(1-p)}{p^2}$.
\end{enumerate}
\end{proposition}

\begin{proof}
This comes indeed as a consequence of Theorem 3.11, at $m=1$.
\end{proof}

Finally, we have the following result, similar to our previous formula $b_{np}=b_p^{*n}$:

\begin{theorem}
We have the following formula, for any $m\in\mathbb N$,
$$c_{mp}=c_p^{*m}$$
relating the geometric laws and the negative binomial laws.
\end{theorem}

\begin{proof}
We have indeed the following straightforward computation:
\begin{eqnarray*}
c_p^{*m}
&=&\left(\sum_{a_1\geq0}(1-p)^{a_1}p\delta_{a_1}\right)*\ldots*\left(\sum_{a_m\geq0}(1-p)^{a_m}p\delta_{a_m}\right)\\
&=&\sum_{a_1,\ldots,a_m\geq0}(1-p)^{a_1+\ldots+a_m}p^m\delta_{a_1+\ldots+a_m}\\
&=&\sum_{s\geq0}\#\left\{a_1,\ldots,a_m\geq0\Big|a_1+\ldots+a_m=s\right\}(1-p)^sp^m\delta_s\\
&=&\sum_{s\geq0}\#\left\{b_1,\ldots,b_m\geq1\Big|b_1+\ldots+b_m=s+m\right\}(1-p)^sp^m\delta_s\\
&=&\sum_{s\geq0}\binom{s+m-1}{m-1}(1-p)^sp^m\delta_s\\
&=&\sum_{s\geq0}\binom{s+m-1}{s}(1-p)^sp^m\delta_s\\
&=&c_{mp}
\end{eqnarray*}

To be more precise here, the count at the end comes from the fact that, in order to partition $s+m$, best thought as a sequence of $s+m$ boxes, as $s+m=b_1+\ldots+b_m$, we have to insert $m-1$ markers between these boxes, at different locations, at the $s+m-1$ positions available. Thus, we are led to the conclusion in the statement.
\end{proof}

Let us discuss now the computation of the higher moments. We have here:

\begin{theorem}
The moments of the negative binomial law $c_{mp}$ are
$$M_1=\frac{m(1-p)}{p}$$
$$M_2=m\,\frac{1-p}{p^2}(1+m-mp)$$
$$M_3=m\,\frac{1-p}{p^3}[(m+1)(m+2)(1-p)^2+3(m+1)(1-p)p+p^2]$$
$$M_4=m\,\frac{1-p}{p^4}[(m+1)(m+2)(m+3)(1-p)^3
+6(m+1)(m+2)(1-p)^2p+7(m+1)(1-p)p^2+p^3]$$
$$\vdots$$
with the computation in general being possible by recurrence.
\end{theorem}

\begin{proof}
We already know the first two formulae from Theorem 3.11 and its proof, and the continuation is by using the same trick, as follows:

\medskip

(1) Regarding the third moment, we can use here the following formula:
$$s^3=s(s-1)(s-2)+3s(s-1)+s$$

Indeed, this gives the formula in the statement for the third moment.

\medskip

(2) Regarding the fourth moment, we can use here the following formula:
$$s^4=s(s-1)(s-2)(s-3)+6s(s-1)(s-2)+7s(s-1)+s$$

Indeed, this gives the formula in the statement for the fourth moment.

\medskip

(3) As a continuation of this, it is quite clear that we can compute $M_k$ by recurrence. We will be back to this later, with more results on the subject.
\end{proof}

Let us record as well a result about the corresponding central moments:

\index{central moments}

\begin{theorem}
The central moments of $c_{mp}$ are given by the formulae
$$M_1'=0$$
$$M_2'=\frac{m(1-p)}{p^2}$$
$$M_3'=\frac{m(1-p)(2-p)}{p^3}$$
$$M_4'=\frac{m(1-p)((3m+6)(1-p)+p^2)}{p^4}$$
$$\vdots$$
with the computation in general being possible by recurrence.
\end{theorem}

\begin{proof}
This follows indeed from the various formulae that we have in Theorem 3.15, and we will leave the computations here as an easy exercise.
\end{proof}

Let us end this discussion with the following result, capturing the essentials:

\begin{theorem}
For the negative binomial law $c_{mp}$,
$$E=\frac{m(1-p)}{p}\ ,\  V=\frac{m(1-p)}{p^2}\ ,\ \gamma=\frac{2-p}{\sqrt{m(1-p)}}\ ,\ \kappa=3+\frac{6}{m}+\frac{p^2}{m(1-p)}$$
are the mean, variance, skewness and kurtosis.
\end{theorem}

\begin{proof}
Regarding the formulae of $E,V$, we know these from Theorem 3.11. As for the formulae of $\gamma,\kappa$, these follow from the formulae of $M_3',M_4'$ from Theorem 3.16.
\end{proof}

Again inspired from what we did before for the usual binomial laws, we have:

\begin{theorem}
The moments of $c_{mp}$ are given by the following formula,
$$M_k=\sum_{\pi\in P(k)}\frac{(m+|\pi|-1)!}{(m-1)!}\left(\frac{1-p}{p}\right)^{|\pi|}$$
with $|.|$ standing as usual for the number of blocks.
\end{theorem}

\begin{proof}
In order to prove the result, let us go back to our computations from the proof of Theorem 3.15. The idea there was that of using a formula as follows:
\begin{eqnarray*}
s^k
&=&s(s-1)\ldots(s-k+1)\\
&&+c_1s(s-1)\ldots (s-k+2)\\
&&\vdots\\
&&+c_{k-2}s(s-1)\\
&&+s
\end{eqnarray*}

But, as explained in the proof of Theorem 3.7, we have the following formula:
$$s^k=\sum_{\pi\in P(k)}\frac{s!}{(s-|\pi|)!}$$

Thus, we have the formula that we need, and this leads to the following formula:
$$M_k=\sum_{\pi\in P(k)}\frac{(m+|\pi|-1)!}{(m-1)!}\left(\frac{1-p}{p}\right)^{|\pi|}$$

And with this, good news, we have reached to the formula in the statement.
\end{proof}

\section*{3c. Hypergeometric laws}

Let us go back to the various types of binomial laws studied before. As a variation of that constructions, we can talk about hypergeometric laws, as follows:

\index{hypergeometric law}
\index{Vandermonde formula}

\begin{theorem}
Given a population of size $N\in\mathbb N$, with $m\in\{0,\ldots,N\}$ of these objects having a certain feature, the probability of having $s$ successes, when performing $p\in\{0,\ldots,N\}$ draws without replacement, and looking for that feature, is
$$P(s)=\frac{\binom{m}{s}\binom{N-m}{p-s}}{\binom{N}{p}}$$
with this being called hypergeometric law of parameters $(N,m,p)$. When doing the same thing, but with replacement, we obtain a binomial law of parameter $m/N$.
\end{theorem}

\begin{proof}
Many things can be said here, the idea being as follows:

\medskip

(1) To start with, the main assertion is clear, because we have a total of $\binom{N}{p}$ possibilities for our $p$ draws without replacement, and among these draws, those amounting to $s$ successes come by multiplying $\binom{m}{s}$, standing for the $s$ successes, and $\binom{N-m}{p-s}$, standing for the $p-s$ fails. Thus, we are led to the formula of $P(s)$ in the statement.

\medskip

(2) Next, observe that the number of successes $s\in\mathbb N$ cannot be anything, because it must be subject to the condition that all binomial coefficients appearing in the formula of $P(s)$ must be well-defined, as usual binomial coefficients. Thus, we must have:
$$0\leq s\leq m\quad,\quad 0\leq p-s\leq N-m$$

In other words, the number of successes $s\in\mathbb N$, where $P(s)>0$, is subject to:
$$s\in\{\max(0,p+m-N),\ldots,\min(m,p)\}$$

(3) Regarding now the fact that we have indeed a probability measure, this is something which is plainly clear, because when performing our $p$ draws, something must happen, I mean we have a well-defined number of successes $s\in\mathbb N$, and $\sum_sP(s)=1$.

\medskip

(4) However, let us comment more on this. The formula $\sum_sP(s)=1$ looks in practice as follows, with the bounds for the summing parameter $s$ being those found in (2):
$$\sum_s\binom{m}{s}\binom{N-m}{p-s}=\binom{N}{p}$$

But this is something which is clear indeed, because when drawing $p$ objects out of $N$ objects, with $m$ of the objects having a certain feature, we have $\binom{N}{p}$ choices on one hand, and $\sum_s\binom{m}{s}\binom{N-m}{p-s}$ choices on the other hand, by looking at this feature.

\medskip

(5) As a further comment on this, $\sum_sP(s)=1$, this is related as well to the Vandermonde formula, which is as follows, coming itself from $(1+x)^m(1+x)^n=(1+x)^{m+n}$, by looking at the coefficient of $x^p$, on the left and on the right:
$$\sum_{s=0}^p\binom{m}{s}\binom{n}{p-s}=\binom{m+n}{p}$$

To be more precise, the formula from (4) is of this type, obtained by looking at the coefficient of $x^p$ inside $(1+x)^m(1+x)^{N-m}=(1+x)^N$, which gives:
$$\sum_s\binom{m}{s}\binom{N-m}{p-s}=\binom{N}{p}$$

(6) Finally, in what regards the second assertion of the theorem, assume that we are doing our $p$ draws, but this time with replacement. In this case, for having $s$ successes, we have $\binom{p}{s}$ choices for the occurences of these successes, and then $(m/N)^s$ probability for the successes, and $(1-m/N)^{n-s}$ probability for the fails, so we obtain:
$$P(s)=\binom{p}{s}\left(\frac{m}{N}\right)^s\left(1-\frac{m}{N}\right)^{n-s}$$

But this is a binomial law of parameter $m/N$, as stated.
\end{proof}

Still at the theoretical level, as an interesting observation about our hypergeometric laws, which is not exactly obvious when looking at their definition, we have:

\begin{proposition}
The hypergeometric law of parameters $(N,m,p)$, given by
$$P(s)=\frac{m!p!(N-m)!(N-p)!}{s!(m-s)!(p-s)!N!(N-m-p+s)!}$$
is symmetric in $m,p$, so is equal to the hypergeometric law of parameters $(N,p,m)$.
\end{proposition}

\begin{proof}
This is indeed something self-explanatory, coming from the definition of the binomial laws. Thus, we are led to the conclusion in the statement.
\end{proof}

Getting now to our usual business, computation of moments of low order, regarding the mean and variance of the hypergeometric laws, these are as follows:

\begin{theorem}
The mean and variance of the hypergeometric laws are
$$E=\frac{mp}{N}\quad,\quad V=\frac{mp(N-m)(N-p)}{N^2(N-1)}$$
with $(N,m,p)$ standing as usual for the parameters.
\end{theorem}

\begin{proof}
Regarding the mean, we can compute it is as follows, by making use of the Vandermonde formula discussed in the proof of Theorem 3.19:
\begin{eqnarray*}
E
&=&\sum_ss\,\frac{\binom{m}{s}\binom{N-m}{p-s}}{\binom{N}{p}}\\
&=&\frac{1}{\binom{N}{p}}\sum_ss\binom{m}{s}\binom{N-m}{p-s}\\
&=&\frac{1}{\binom{N}{p}}\sum_sm\binom{m-1}{s-1}\binom{N-m}{p-s}\\
&=&\frac{m}{\binom{N}{p}}\sum_s\binom{m-1}{s-1}\binom{N-m}{p-s}\\
&=&\frac{m}{\binom{N}{p}}\binom{N-1}{p-1}\\
&=&m\cdot\frac{p!(N-p)!}{N!}\cdot\frac{(N-1)!}{(p-1)!(N-p)!}\\
&=&\frac{mp}{N}
\end{eqnarray*}

Next, with the same trick we can compute $M_2-M_1$, and we obtain:
$$M_2-M_1=\frac{mp(m-1)(p-1)}{N(N-1)}$$

We conclude that the second moment is given by the following formula:
$$M_2=\frac{mp(mp-m-p+N)}{N(N-1)}$$

Thus, we are led to the variance formula in the statement.
\end{proof}

In order to get beyond this, we can use some algebra. Consider the standard matrix coordinates over the symmetric group $S_N\subset O_N$, which form a matrix, as follows:
$$u=\begin{pmatrix}
u_{11}&\ldots&u_{1p}&\ldots&\ldots&u_{1N}\\
\vdots&\bigstar&\vdots&&&\vdots\\
u_{m1}&\ldots&u_{mp}&\ldots&\ldots&u_{mN}\\
\vdots&&\vdots&&&\vdots\\
\vdots&&\vdots&&&\vdots\\
u_{N1}&\ldots&u_{Np}&\ldots&\ldots&u_{NN}
\end{pmatrix}$$

The point now is that the sum of the entries of the rectangular matrix $\bigstar$ appearing in the upper left corner is hypergeometric, as shown by the following result:

\index{hypergeometric variable}

\begin{theorem}
The following variable over the symmetric group $S_N\subset O_N$,
$$S_{mp}=\sum_{i=1}^m\sum_{j=1}^pu_{ij}$$
is hypergeometric, having parameters $(N,m,p)$.
\end{theorem}

\begin{proof}
We know that the coordinates of the symmetric group, viewed as group of permutation matrices, $S_N\subset O_N$, are given by the following formula:
$$u_{ij}=\chi\left(\sigma\in S_N\Big|\sigma(j)=i\right)$$

Thus, the variable in the statement is given by the following formula:
$$S_{mp}=\sum_{i=1}^m\sum_{j=1}^p\chi\left(\sigma\in S_N\Big|\sigma(j)=i\right)$$

Let us compute now the law of this variable. We have the following formula:
\begin{eqnarray*}
P(S_{mp}=s)
&=&\frac{1}{N!}\,\#\left(\sigma\in S_N\Big|S_{mp}(\sigma)=s\right)\\
&=&\frac{1}{N!}\,\#\left(\sigma\in S_N\Big|\#\left\{i\leq m,j\leq p\big|\sigma(j)=i\right\}=s\right)\\
&=&\frac{1}{N!}\,\#\left(\sigma\in S_N\Big|\#\left\{j\leq p\big|\sigma(j)\leq m\right\}=s\right)
\end{eqnarray*}

Now let us try to count the permutations on the right. Such a permutation $\sigma\in S_N$ comes from 5 operations involved, as follows:

\medskip

(1) We first have to pick a subset $X\subset\{1,\ldots,p\}$, with $|X|=s$.

\medskip

(2) We also have to pick a subset $Y\subset\{1,\ldots,m\}$, with $|Y|=s$.

\medskip

(3) Then, we have to bijectively map $X\to Y$.

\medskip

(4) Next, we have to map injectively $\{1,\ldots,p\}-X\to\{m+1,\ldots,N\}$.

\medskip

(5) And finally, we have to bijectively map $\{p+1,\ldots,N\}$ to what is left.

\medskip

Regarding now the precise numbers of choices for these 5 operations involved, these are easy to compute, as follows:

\medskip

(1) Here we have $\binom{p}{s}$ choices.

\medskip

(2) Here we have $\binom{m}{s}$ choices.

\medskip

(3) Here we have $s!$ choices.

\medskip

(4) Here we have $\frac{(N-m)!}{(N-m-p+s)!}$ choices.

\medskip

(5) Here we have $(N-p)!$ choices.

\medskip

We can now finish our probability computation started above, as follows:
\begin{eqnarray*}
P(S_{mp}=s)
&=&\frac{1}{N!}\binom{p}{s}\binom{m}{s}s!\,\frac{(N-m)!}{(N-m-p+s)!}(N-p)!\\
&=&\frac{1}{N!}\cdot\frac{p!}{s!(p-s)!}\cdot\frac{m!}{s!(m-s)!}\cdot s!\cdot\frac{(N-m)!}{(N-m-p+s)!}(N-p)!\\
&=&\frac{p!m!(N-m)!(N-p)!}{N!s!(p-s)!(m-s)!(N-m-p+s)!}\\
&=&\frac{m!}{s!(m-s)!}\cdot\frac{(N-m)!}{(p-s)!(N-m-p+s)!}\cdot\frac{p!(N-p)!}{N!}\\
&=&\frac{\binom{m}{s}\binom{N-m}{p-s}}{\binom{N}{p}}
\end{eqnarray*}

We are therefore led to the conclusion in the statement.
\end{proof}

As a first observation, with the above in hand, we have now a better explanation for the fact, that we knew from Proposition 3.20, coming there via a quite opaque proof, that the hypergeometric law of parameters $(N,m,p)$ is symmetric in $m,p$. Nice.
 
\bigskip

In order to discuss now the moment problematics, using our group theory approach coming from Theorem 3.22, we will need the following standard fact:

\begin{proposition}
The integrals over the symmetric group $S_N$ are given by
$$\int_{S_N}u_{i_1j_1}\ldots u_{i_kj_k}=\begin{cases}
\frac{(N-|\ker i|)!}{N!}&{\rm if}\ \ker i=\ker j\\
0&{\rm otherwise}
\end{cases}$$
where $\ker i$ denotes as usual the partition of $\{1,\ldots,k\}$ whose blocks collect the equal indices of $i$, and where $|.|$ denotes the number of blocks.
\end{proposition}

\begin{proof}
The first assertion follows from Stone-Weierstrass, because the standard coordinates $g_{ij}$ separate the points of $S_N$. Regarding now the second assertion, according to the definition of the coordinates $u_{ij}$, the integrals in the statement are given by:
$$\int_{S_N}u_{i_1j_1}\ldots u_{i_kj_k}=\frac{1}{N!}\#\left\{\sigma\in S_N\Big|\sigma(j_1)=i_1,\ldots,\sigma(j_k)=i_k\right\}$$

Now observe that the existence of $\sigma\in S_N$ as above requires $i_m=i_n\iff j_m=j_n$, so the above integral vanishes when the following condition is satisfied:
$$\ker i\neq\ker j$$

Regarding now the case $\ker i=\ker j$, if we denote by $b\in\{1,\ldots,k\}$ the number of blocks of this partition $\ker i=\ker j$, we have $N-b$ points to be sent bijectively to $N-b$ points, and so $(N-b)!$ solutions, and the integral is $\frac{(N-b)!}{N!}$, as claimed.
\end{proof}

Getting back now to the hypergeometric variables, as a continuation of our previous work on the moments, and standing as a final result on the subject, we have:

\index{moments of hypergeometric law}

\begin{theorem}
For the hypergeometric law of parameters $(N,m,p)$ we have
$$M_k=\sum_{\pi\in P(k)}\frac{m!}{(m-|\pi|)!}\cdot\frac{p!}{(p-|\pi|)!}\cdot\frac{(N-|\pi|)!}{N!}$$
where $|.|$ denotes as usual the number of blocks.
\end{theorem}

\begin{proof}
We have indeed, by some magic, the following computation, using the interpretation from Theorem 3.22, and the integration formula in Proposition 3.23:
\begin{eqnarray*}
M_k
&=&\int_{S_N}S_{mp}^k\\
&=&\int_{S_N}\left(\sum_{i=1}^m\sum_{j=1}^pu_{ij}\right)^k\\
&=&\int_{S_N}\sum_{i_1=1}^m\ldots\sum_{i_k=1}^m\sum_{j_1=1}^p\ldots\sum_{j_k=1}^pu_{i_1j_1}\ldots u_{i_kj_k}\\
&=&\sum_{i_1=1}^m\ldots\sum_{i_k=1}^m\sum_{j_1=1}^p\ldots\sum_{j_k=1}^p\int_{S_N}u_{i_1j_1}\ldots u_{i_kj_k}\\
&=&\sum_{i_1=1}^m\ldots\sum_{i_k=1}^m\sum_{j_1=1}^p\ldots\sum_{j_k=1}^p\delta_{\ker i,\ker j}\frac{(N-|\ker i|)!}{N!}\\
&=&\sum_{\pi\in P(k)}\sum_{i_1=1}^m\ldots\sum_{i_k=1}^m\sum_{j_1=1}^p\ldots\sum_{j_k=1}^p\delta_{\ker i,\ker j,\pi}\frac{(N-|\pi|)!}{N!}\\
&=&\sum_{\pi\in P(k)}\left(\sum_{i_1=1}^m\ldots\sum_{i_k=1}^m\delta_{\ker i,\pi}\right)
\left(\sum_{j_1=1}^p\ldots\sum_{j_k=1}^p\delta_{\ker j,\pi}\right)\frac{(N-|\pi|)!}{N!}\\
&=&\sum_{\pi\in P(k)}\frac{m!}{(m-|\pi|)!}\cdot\frac{p!}{(p-|\pi|)!}\cdot\frac{(N-|\pi|)!}{N!}
\end{eqnarray*}

Thus, we are led to the conclusion in the statement, and I would leave to you the remaining discussion, regarding the central moments, skewness and kurtosis.
\end{proof}

\section*{3d. Negative versions}

As a variation of what we did so far in this chapter, regarding the hypergeometric laws, we can talk as well about negative hypergeometric laws, as follows:

\index{negative hypergeometric law}
\index{negative binomial law}

\begin{theorem}
Given a population of size $N\in\mathbb N$, with $m\in\{0,\ldots,N\}$ of these objects having a certain feature, the probability of having $s$ successes, when performing draws without replacement, and stopping after $r\in\{0,\ldots,N-m\}$ failures, is
$$P(s)=\frac{\binom{s+r-1}{s}\binom{N-r-s}{m-s}}{\binom{N}{m}}$$
with this being called negative hypergeometric law of parameters $(N,m,r)$. When doing this with replacement, we obtain a negative binomial law of parameter $m/N$.
\end{theorem}

\begin{proof}
Many things can be said here, the idea being as follows:

\medskip

(1) Observe first the similarity with the definition of the usual hypergeometric laws, from Theorem 3.19. To be more precise, the common part is that, in both these cases, we have a population of size $N\in\mathbb N$, with $m\in\{0,\ldots,N\}$ of these objects having a certain feature, and we are looking for the probability of having $s$ successes, when performing draws with replacement. The difference comes from the stopping convention, for the hypergeometric laws this being after $p\in\{0,\ldots,N\}$ draws, and for the negative hypergeometric laws this being after $r\in\{0,\ldots,N-m\}$ failures.

\medskip

(2) In practice, based on this observation, we can develop the theory of negative hypergeometric laws by using the previously developed theory of the usual hypergeometric laws. Indeed, to start with, the probability $P(s)$ that we are interested in comes by multiplying a probability coming from a usual hypergeometric law, corresponding to the $m$ succeses in the first $m+r-1$ draws, and a basic fraction, corresponding to the failure required at the $(m+r)$-th draw. In practice, this gives the following formula:
\begin{eqnarray*}
P(s)
&=&\frac{\binom{m}{s}\binom{N-m}{s+r-1-s}}{\binom{N}{s+r-1}}\cdot\frac{N-m-(r-1)}{N-(s+r-1)}\\
&=&\frac{\binom{m}{s}\binom{N-m}{r-1}}{\binom{N}{s+r-1}}\cdot\frac{N-m-r+1}{N-s-r+1}\\
&=&\frac{m!(N-m)!(s+r-1)!(N-s-r+1)!(N-m-r+1)}{s!(m-s)!(r-1)!N!(N-m-r+1)!(N-s-r+1)}\\
&=&\frac{m!(N-m)!(s+r-1)!(N-s-r)!}{s!(m-s)!(r-1)!N!(N-m-r)!}\\
&=&\frac{(s+r-1)!(N-r-s)!(N-m)!m!}{s!(r-1)!(m-s)!(N-r-m)!N!}\\
&=&\frac{\binom{s+r-1}{s}\binom{N-r-s}{m-s}}{\binom{N}{m}}
\end{eqnarray*}

Thus, we have indeed the formula of $P(s)$ in the statement.

\medskip

(3) Regarding now the fact that we have indeed a probability measure, $\sum_sP(s)=1$, this is something which is plainly clear from definitions, but which can be deduced as well from the Vandermonde formula, which is as follows, coming from $(1+x)^m(1+x)^n=(1+x)^{m+n}$, by looking at the coefficient of $x^p$, on the left and on the right:
$$\sum_s\binom{m}{s}\binom{n}{p-s}=\binom{m+n}{p}$$

Indeed, by using this formula, along with the identity $\binom{n}{k}=(-1)^k\binom{k-n-1}{k}$ for the generalized binomial coefficients, which is something trivial, we have:
\begin{eqnarray*}
\sum_s\binom{s+r-1}{s}\binom{N-r-s}{m-s}
&=&\sum_s(-1)^s\binom{-r}{s}(-1)^{m-s}\binom{m+r-N-1}{m-s}\\
&=&(-1)^m\sum_s\binom{-r}{s}\binom{m+r-N-1}{m-s}\\
&=&(-1)^m\binom{m-N-1}{m}\\
&=&\binom{N}{m}
\end{eqnarray*}

(4) Finally, in what regards the second assertion of the theorem, assume that we are doing our draws, but this time with replacement. In this case, for having $s$ successes, we have $\binom{s+r-1}{s}$ choices for the occurences of these successes, and then $(m/N)^s$ probability for the successes, and $(1-m/N)^r$ probability for the fails, so we obtain:
$$P(s)=\binom{s+r-1}{s}\left(\frac{m}{N}\right)^s\left(1-\frac{m}{N}\right)^r$$

But this is a negative binomial law of parameter $m/N$, as stated.
\end{proof}

Regarding the mean and variance of the negative hypergeometric laws, these can be computed a bit as for the usual hypergeometric laws, the result being as follows:

\begin{theorem}
The mean and variance of the negative hypergeometric laws are
$$E=\frac{mr}{N-m+1}\quad,\quad V=\frac{mr(N+1)(N-m-r+1)}{(N-m+1)^2(N-m+2)}$$
with $(N,m,r)$ standing as usual for the parameters.
\end{theorem}

\begin{proof}
The computations are quite similar to those from before, as follows:

\medskip

(1) Regarding the mean, we can compute it with a trick, namely:
$$s=(s+r)-r$$

Indeed, by using this trick, and the summation formula based on Vandermonde discussed in (3) in the proof of Theorem 3.25, we obtain the following formula:
\begin{eqnarray*}
E
&=&\sum_ss\,\frac{\binom{s+r-1}{s}\binom{N-r-s}{m-s}}{\binom{N}{m}}\\
&=&\sum_s[(s+r)-r]\,\frac{\binom{s+r-1}{s}\binom{N-r-s}{m-s}}{\binom{N}{m}}\\
&=&\left[\sum_s(s+r)\,\frac{\binom{s+r-1}{s}\binom{N-r-s}{m-s}}{\binom{N}{m}}\right]-r\\
&=&\frac{r}{\binom{N}{m}}\left[\sum_s\frac{s+r}{r}\binom{s+r-1}{s}\binom{N-r-s}{m-s}\right]-r\\
&=&\frac{r}{\binom{N}{m}}\left[\sum_s\binom{s+r}{s}\binom{N-r-s}{m-s}\right]-r\\
&=&\frac{r}{\binom{N}{m}}\binom{N+1}{m}-r\\
&=&r\cdot\frac{m!(N-m)!}{N!}\cdot\frac{(N+1)!}{m!(N-m+1)!}-r\\
&=&\frac{(N+1)r}{N-m+1}-r\\
&=&\frac{mr}{N-m+1}
\end{eqnarray*}

(2) In order to compute now $M_2-M_1$, we can use a similar trick, as follows:
$$s^2-s=(s+r)(s+r+1)-2(r+1)(s+r)+r(r+1)$$

Indeed, by summing as in (1), using Vandermonde, we get the following formula:
$$M_2-M_1=\frac{m(m-1)r(r+1)}{(N-m+1)(N-m+2)}$$

(3) We conclude that the second moment is given by the following formula:
$$M_2=\frac{m(m-1)r(r+1)}{(N-m+1)(N-m+2)}+\frac{mr}{N-m+1}$$

(4) But this gives the formula of the variance in the statement.
\end{proof}

As a main result now regarding the negative hypergeometric laws, we have:

\index{moments of negative hypergeometric law}

\begin{theorem}
For the negative hypergeometric law of parameters $(N,m,r)$,
$$M_k=\sum_{\pi\in P(k)}\frac{m!}{(m-|\pi|)!}\cdot\frac{(r+|\pi|-1)!}{(r-1)!}\cdot\frac{(N-m)!}{(N-m+|\pi|)!}$$
where $|.|$ denotes as usual the number of blocks.
\end{theorem}

\begin{proof}
This is something quite standard, the idea being that the formula in the statement follows by recurrence, by using the tricks from the proof of Theorem 3.26. We will leave the various computations here, as well as the remaining discussion, regarding the central moments, skewness and kurtosis, as an instructive exercise. Enjoy.
\end{proof}

And with this, end of our study of the binomial laws, and their generalizations. All the above was of course quite quick, and for a more detailed treatement of all this, including more phenomenology, numerics, tables and so on, you can check my previous book \cite{ba1}, where all this material is explained over a full 3 chapters, worth 75 pages.

\section*{3e. Exercises}

This was a quite technical chapter, as we like them, and as exercises, we have:

\begin{exercise}
Further meditate on the symmetry of the hypergeometric law in $m,p$.
\end{exercise}

\begin{exercise}
Do the missing hypergeometric law computations, at $k=3$.
\end{exercise}

\begin{exercise}
Do the missing hypergeometric law computations, at $k=4$.
\end{exercise}

\begin{exercise}
Compute the higher moments by a direct recurrence.
\end{exercise}

\begin{exercise}
Work out the asymptotics of the hypergeometric laws.
\end{exercise}

\begin{exercise}
Learn more about the negative hypergeometric laws.
\end{exercise}

\begin{exercise}
Try finding a group-theoretical model for these latter laws.
\end{exercise}

\begin{exercise}
Do the missing moment computations, for these latter laws.
\end{exercise}

As bonus exercise, quite scary, learn a bit about hypergeometric functions.

\chapter{Poisson laws}

\section*{4a. Poisson limits}

We would like to discuss in this chapter the central objects in discrete probability, which are the Poisson laws $p_t$, appearing via the Poisson Limit Theorem (PLT). Let us start with the following definition, which is something quite natural:

\index{Poisson law}
\index{Poisson variable}

\begin{definition}
The Poisson law of parameter $1$ is the following measure,
$$p_1=\frac{1}{e}\sum_{k\geq0}\frac{\delta_k}{k!}$$
and the Poisson law of parameter $t>0$ is the following measure,
$$p_t=e^{-t}\sum_{k\geq0}\frac{t^k}{k!}\,\delta_k$$
with the letter ``p'' standing for Poisson.
\end{definition}

As a first observation, our laws have indeed mass 1, as they should, due to the following formula, which is actually the key formula in analysis, as we learned in chapter 1:
$$e^t=\sum_{k\geq0}\frac{t^k}{k!}$$

We will see in the moment why the above laws appear a bit everywhere, in discrete contexts, the reasons for this coming from the Poisson Limit Theorem. In the meantime, as a piece of advertisement for this, Poisson laws and limits, let us formulate:

\begin{fact}
The number of chewing gum pieces on a single tile of a sidewalk
$$\xymatrix@R=5pt@C=7pt{
\ar@{-}[rrrrrrrrrr]\ar@{-}[dddddd]&&\ar@{-}[dddddd]&&\ar@{-}[dddddd]&&\ar@{-}[dddddd]&&\ar@{-}[dddddd]&&\ar@{-}[dddddd]\\
&\bullet&&&&\bullet&&\ \,\bullet\bullet\,\ &&&&\ldots\\
\ar@{-}[rrrrrrrrrr]&&&&&&&&&&\\
&&&\,\bullet\bullet\bullet\,&&&&\bullet&&\ \,\bullet\bullet\,\ &&\ldots\\
\ar@{-}[rrrrrrrrrr]&&&&&&&&&&\\
&\ \ \, \bullet\ \ \,&&&&\ \,\bullet\bullet\,\ &&\bullet&&&&\ldots\\
\ar@{-}[rrrrrrrrrr]&&&&&&&&&&
}$$
follows a Poisson law.
\end{fact}

But more on such things later, once we will have tools for proving them. Let us first develop some straightforward general theory. We first have the following result:

\index{mean of Poisson law}
\index{variance of Poisson law}

\begin{proposition}
The mean and variance of $p_t$ are given by:
$$E=t\quad,\quad V=t$$
In particular for the Poisson law $p_1$ we have $E=1,V=1$.
\end{proposition}

\begin{proof}
Regarding the mean, this can be computed as follows:
\begin{eqnarray*}
E
&=&e^{-t}\sum_{k\geq0}\frac{t^k}{k!}\cdot k\\
&=&e^{-t}\sum_{k\geq1}\frac{t^k}{(k-1)!}\\
&=&e^{-t}\sum_{l\geq0}\frac{t^{l+1}}{l!}\\
&=&te^{-t}\sum_{l\geq0}\frac{t^l}{l!}\\
&=&t
\end{eqnarray*}

For the variance, we first compute the second moment, as follows:
\begin{eqnarray*}
M_2
&=&e^{-t}\sum_{k\geq0}\frac{t^k}{k!}\cdot k^2\\
&=&e^{-t}\sum_{k\geq1}\frac{t^kk}{(k-1)!}\\
&=&e^{-t}\sum_{l\geq0}\frac{t^{l+1}(l+1)}{l!}\\
&=&te^{-t}\sum_{l\geq0}\frac{t^ll}{l!}+te^{-t}\sum_{l\geq0}\frac{t^l}{l!}\\
&=&te^{-t}\sum_{l\geq1}\frac{t^l}{(l-1)!}+t\\
&=&t^2e^{-t}\sum_{m\geq0}\frac{t^m}{m!}+t\\
&=&t^2+t
\end{eqnarray*}

Thus the variance is given by $V=(t^2+t)-t^2=t$, as claimed.
\end{proof}

Next, we have the following quite conceptual result, in relation with convolution:

\index{convolution semigroup}

\begin{theorem}
We have the following formula, for any $s,t>0$,
$$p_s*p_t=p_{s+t}$$
so the Poisson laws form a convolution semigroup.
\end{theorem}

\begin{proof}
By using $\delta_k*\delta_l=\delta_{k+l}$ and the binomial formula, we obtain:
\begin{eqnarray*}
p_s*p_t
&=&e^{-s}\sum_{k\geq0}\frac{s^k}{k!}\,\delta_k*e^{-t}\sum_{l\geq0}\frac{t^l}{l!}\,\delta_l\\
&=&e^{-s-t}\sum_{n\geq0}\delta_n\sum_{k+l=n}\frac{s^kt^l}{k!l!}\\
&=&e^{-s-t}\sum_{n\geq0}\frac{\delta_n}{n!}\sum_{k+l=n}\frac{n!}{k!l!}s^kt^l\\\
&=&e^{-s-t}\sum_{n\geq0}\frac{(s+t)^n}{n!}\,\delta_n\\
&=&p_{s+t}
\end{eqnarray*}

Thus, we are led to the conclusion in the statement.
\end{proof}

Next in line, we have the following result, which is fundamental as well:

\index{convolution exponential}
\index{formal exponential}

\begin{theorem}
The Poisson laws appear as formal exponentials
$$p_t=\sum_{k\geq0}\frac{t^k(\delta_1-\delta_0)^{*k}}{k!}$$
with respect to the convolution of measures $*$.
\end{theorem}

\begin{proof}
By using the binomial formula, the measure on the right is:
\begin{eqnarray*}
\mu
&=&\sum_{k\geq0}t^k\sum_{r+s=k}(-1)^s\frac{\delta_r}{r!s!}\\
&=&\sum_{r\geq0}\frac{t^r\delta_r}{r!}\sum_{s\geq0}\frac{(-1)^s}{s!}\\
&=&\frac{1}{e}\sum_{r\geq0}\frac{t^r\delta_r}{r!}\\
&=&p_t
\end{eqnarray*}

Thus, we are led to the conclusion in the statement.
\end{proof}

Getting now to what we wanted to do, Poisson limits, it is convenient to use, as our main tool for dealing with independence, the Fourier transform, coming via:

\index{independence}
\index{Fourier transform}

\begin{theorem}
Assuming that $f,g:X\to\mathbb R$ are independent, we have
$$F_{f+g}=F_fF_g$$
where $F_f(x)=E(e^{ixf})$ is the Fourier transform.
\end{theorem}

\begin{proof}
We have indeed the following computation, using $\mu_{f+g}=\mu_f*\mu_g$:
\begin{eqnarray*}
F_{f+g}(x)
&=&\int_Xe^{ixz}d\mu_{f+g}(z)\\
&=&\int_Xe^{ixz}d(\mu_f*\mu_g)(z)\\
&=&\int_{X\times X}e^{ix(z+t)}d\mu_f(z)d\mu_g(t)\\
&=&\int_Xe^{ixz}d\mu_f(z)\int_Xe^{ixt}d\mu_g(t)\\
&=&F_f(x)F_g(x)
\end{eqnarray*}

Thus, we are led to the conclusion in the statement.
\end{proof}

As an illustration, for the basic measures that we know, from chapter 3, we have:

\begin{proposition}
The Fourier transform of the binomial law $b_{np}$ is given by:
$$F_f(x)=(1-p+pe^{ix})^n$$
Also, the Fourier transform of the negative binomial law $c_{mp}$ is given by:
$$F_f(x)=\left(\frac{p}{1+(p-1)e^{ix}}\right)^m$$
In particular, in both cases we can see that $\log F$ is linear in $n/m$, as it should.
\end{proposition}

\begin{proof}
For the binomial law, we have indeed the following computation:
\begin{eqnarray*}
F_f(x)
&=&\sum_{s=0}^n\binom{n}{s}p^s(1-p)^{n-s}e^{ixs}\\
&=&(1-p)^n\sum_{s=0}^n\binom{n}{s}\left(\frac{p}{1-p}\,e^{ix}\right)^s\\
&=&(1-p)^n\left(1+\frac{p}{1-p}\,e^{ix}\right)^n\\
&=&(1-p+pe^{ix})^n
\end{eqnarray*}

For the negative binomial law the computation is quite similar, as follows:
\begin{eqnarray*}
F_f(x)
&=&\sum_{s\geq0}\binom{-m}{s}(p-1)^sp^me^{ixs}\\
&=&p^m\sum_{s\geq0}\binom{-m}{s}[(p-1)e^{ix}]^s\\
&=&p^m(1+(p-1)e^{ix})^{-m}\\
&=&\left(\frac{p}{1+(p-1)e^{ix}}\right)^m
\end{eqnarray*}

As for the last assertion, this is just an observation, to be related, via Theorem 4.6, to the formulae $b_{np}=b_p^{*n}$ and $c_{mp}=c_p^{*m}$, that we know from chapter 3.
\end{proof}

Regarding now the Fourier transform computation for $p_t$, this is as follows:

\index{Fourier transform}

\begin{theorem}
The Fourier transform of $p_t$ is given by
$$F_{p_t}(x)=\exp\left((e^{ix}-1)t\right)$$
for any $t>0$, with its logarithm being linear in $t$, as it should.
\end{theorem}

\begin{proof}
We have indeed the following computation, which gives the result:
\begin{eqnarray*}
F_{p_t}(x)
&=&e^{-t}\sum_{k\geq0}\frac{t^k}{k!}F_{\delta_k}(x)\\
&=&e^{-t}\sum_{k\geq0}\frac{t^k}{k!}\,e^{ikx}\\
&=&e^{-t}\sum_{k\geq0}\frac{(e^{ix}t)^k}{k!}\\
&=&\exp(-t)\exp(e^{ix}t)\\
&=&\exp\left((e^{ix}-1)t\right)
\end{eqnarray*}

As for the last assertion, this is just a remark, to be related to Theorem 4.4.
\end{proof}

Good news, we can now establish the Poisson Limit Theorem, as follows:

\index{PLT}
\index{Poisson Limit Theorem}
\index{Bernoulli laws}
\index{Poisson limit}

\begin{theorem}[PLT]
We have the following convergence, in moments,
$$\left(\left(1-\frac{t}{n}\right)\delta_0+\frac{t}{n}\,\delta_1\right)^{*n}\to p_t$$
for any $t>0$. Equivalently, $b_{np}\to p_t$, when $p=t/n$ with $t>0$ fixed.
\end{theorem}

\begin{proof}
Let us denote by $\mu_n$ the Bernoulli law under the convolution sign:
$$\mu_n=\left(1-\frac{t}{n}\right)\delta_0+\frac{t}{n}\,\delta_1$$

We have the following computation, for the Fourier transform of the limit: 
\begin{eqnarray*}
F_{\delta_r}(x)=e^{irx}
&\implies&F_{\mu_n}(x)=\left(1-\frac{t}{n}\right)+\frac{t}{n}\,e^{ix}\\
&\implies&F_{\mu_n^{*n}}(x)=\left(\left(1-\frac{t}{n}\right)+\frac{t}{n}\,e^{ix}\right)^n\\
&\implies&F_{\mu_n^{*n}}(x)=\left(1+\frac{(e^{ix}-1)t}{n}\right)^n\\
&\implies&F(x)=\exp\left((e^{ix}-1)t\right)
\end{eqnarray*}

Thus, we obtain indeed the Fourier transform of $p_t$, as desired.
\end{proof}

We have as well a negative Poisson Limit Theorem, as follows:

\index{PLT}
\index{Poisson Limit Theorem}
\index{Bernoulli laws}
\index{Poisson limit}
\index{negative PLT}
\index{NPLT}
\index{geometric law}

\begin{theorem}[negative PLT]
We have the following convergence, in moments,
$$\left(\sum_{k\geq0}\left(\frac{t}{m+t}\right)^k\frac{m}{m+t}\,\delta_k\right)^{*m}\to p_t$$
for any $t>0$. Equivalently, $c_{mp}\to p_t$, when $p=m/(m+t)$ with $t>0$ fixed.
\end{theorem}

\begin{proof}
Let us denote by $\nu_m$ the geometric law under the convolution sign:
$$\nu_m=\sum_{k\geq0}\left(\frac{t}{m+t}\right)^k\frac{m}{m+t}\,\delta_k$$

We have the following computation, for the Fourier transform of the limit: 
\begin{eqnarray*}
F_{\delta_k}(x)=e^{ikx}
&\implies&F_{\nu_m}(x)=\sum_{k\geq0}\left(\frac{t}{m+t}\right)^k\frac{m}{m+t}\,e^{ikx}\\
&\implies&F_{\nu_m}(x)=\frac{m}{m+t-te^{ix}}\\
&\implies&F_{\nu_m^{*m}}(x)=\left(\frac{m}{m+t-te^{ix}}\right)^m\\
&\implies&F_{\nu_m^{*m}}(x)=\left(1+\frac{t(e^{ix}-1)}{m+t-te^{ix}}\right)^m\\
&\implies&F(x)=\exp\left(t(e^{ix}-1)\right)
\end{eqnarray*}

Thus, we obtain indeed the Fourier transform of $p_t$, as desired.
\end{proof}

In practice, Theorems 4.9 and 4.10 lead to many occurrences of the Poisson laws, in relation with continuous phenomena from the real life. For instance, we have:

\begin{theorem}
The number of chewing gum pieces on a single tile of a sidewalk
$$\xymatrix@R=5pt@C=7pt{
\ar@{-}[rrrrrrrrrr]\ar@{-}[dddddd]&&\ar@{-}[dddddd]&&\ar@{-}[dddddd]&&\ar@{-}[dddddd]&&\ar@{-}[dddddd]&&\ar@{-}[dddddd]\\
&\bullet&&&&\bullet&&\ \,\bullet\bullet\,\ &&&&\ldots\\
\ar@{-}[rrrrrrrrrr]&&&&&&&&&&\\
&&&\,\bullet\bullet\bullet\,&&&&\bullet&&\ \,\bullet\bullet\,\ &&\ldots\\
\ar@{-}[rrrrrrrrrr]&&&&&&&&&&\\
&\ \ \, \bullet\ \ \,&&&&\ \,\bullet\bullet\,\ &&\bullet&&&&\ldots\\
\ar@{-}[rrrrrrrrrr]&&&&&&&&&&
}$$
follows a Poisson law.
\end{theorem}

\begin{proof}
It is pretty much clear, by thinking a bit, that the law that we are looking for, that of the number of chewing gum pieces on a single tile of a sidewalk, must be a discrete law, supported by $\mathbb N$, and appearing as a suitable limit of binomial laws:
$$b_{np}\to\sum_{k\geq 0}c_k\delta_k$$

But, how does this limit exactly appear? In answer, since the total number of pieces of chewing gum must not blow up with $n\to\infty$, and more specifically, must remain approximately constant, we must impose the condition $p=t/n$, for some fixed $t>0$. But with this, we are led into the Poisson limiting formula from Theorem 4.9, namely:
$$\left(\left(1-\frac{t}{n}\right)\delta_0+\frac{t}{n}\,\delta_1\right)^{*n}\to p_t$$

Thus, we are led to the conclusion in the statement. And to the above picture too, the data there coming, up to an epsilon, from the Poisson law of parameter $t=1$.
\end{proof}

\section*{4b. Bell and Stirling}

Getting back now to the basics of the Poisson laws, let us discuss now the higher moments of these laws. As we will soon discover, things are quite tricky here. Let us start with something that we are familiar with, from chapter 3, namely:

\index{partitions}
\index{partitions of set}

\begin{definition}
We denote by $P(k)$ the set of partitions of $\{1,\ldots,k\}$, with these partitions $\pi\in P(k)$ being most conveniently being drawn as diagrams,
$$
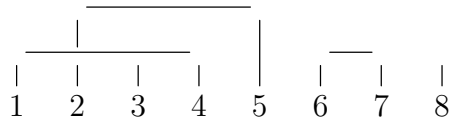
\xymatrix@R=10pt@C=10pt{
&\ar@{-}[rrr]&&&&&&\\
\ar@{-}[rrr]&\ar@{-}[u]&&&&\ar@{-}[r]&&\\
1\ar@{-}[u]&2\ar@{-}[u]&3\ar@{-}[u]&4\ar@{-}[u]&5\ar@{-}[uu]&6\ar@{-}[u]&7\ar@{-}[u]&8\ar@{-}[u]
}$$
with the strings joining the numbers belonging to the same block of $\pi$. That is, the above diagram represents the partition $\{1,\ldots,8\}=\{1,3,4\}\cup\{2,5\}\cup\{6,7\}\cup\{8\}$.
\end{definition}

Now, let us study these partitions. And here, surprise, instead of pulling a theorem, as you would expect, we must formulate something quite modest, as follows:

\index{number of partitions}
\index{Bell numbers}

\begin{proposition}
The Bell numbers $B_k=|P(k)|$ satisfy the recurrence relation
$$B_{k+1}=\sum_s\binom{k}{s}B_{k-s}$$
with initial data $B_0=1$, $B_1=1$, and are numerically as follows:
$$1,\ 1,\ 2,\ 5,\ 15,\ 52,\ 203,\ 877,\ 4140,\ 21147,\ 115975,\ 678570,\ \ldots$$
However, there is no mathematical formula for $B_k$.
\end{proposition}

\begin{proof}
There are several things going on here, the idea being as follows:

\medskip

(1) Experiments first, before anything, let us compute a few Bell numbers. Obviously $B_1=1$, and then we have $B_2=2$, the partitions being as follows:
$$|\ |\quad,\quad\sqcap$$

Next, we have $B_3=5$, the partitions at $k=3$ being as follows:
$$|\ |\ |\quad,\quad\sqcap\ |\quad,\quad\sqcap\hskip-3.2mm{\ }_|\quad,\quad|\ \sqcap\quad,\quad\sqcap\hskip-0.7mm\sqcap$$

At $k=4$ now, things become more complex, and it is better to trick. We can count the partitions up to permutations of the corresponding diagrams, and with this convention made, here are the relevant partitions and their multiplicities, leading to $B_4=15$:
$$|\ |\ |\ |\times1\quad,\quad \sqcap\ |\ |\times 6\quad,\quad \sqcap\sqcap\times 3
\quad,\quad \sqcap\hskip-1.6mm\sqcap|\times4\quad,\quad \sqcap\hskip-1.6mm\sqcap\hskip-1.6mm\sqcap\times1$$

The same method works at $k=5$, with the block distributions and multiplicities, which are simpler to draw than partitions, being as follows, leading to $B_5=52$:
$$11111\to1\ ,\ 2111\to10\ ,\ 221\to15\ ,\ 311\to10\ ,\ 32\to10\ ,\ 41\to5\ ,\ 5\to1$$

As for the case $k=6$, where $B_6=203$, we will leave this as an instructive exercise.

\medskip

(2) Let us try now to find a recurrence for these Bell numbers. Since a partition of $\{1,\ldots,k+1\}$ appears by choosing $s$ partners for $1$, among the $k$ numbers available, and then partitioning the $k-s$ elements left, we have the following formula:
$$B_{k+1}=\sum_s\binom{k}{s}B_{k-s}$$

Observe that this formula forces us to talk about $B_0=1$, as done in the statement.

\medskip

(3) As for the last assertion, regarding the non-computability of the Bell numbers, take this as a physics fact. Mankind has tried to find a formula for these numbers, had not found anything, and we are reporting here this finding, which is of course rock-solid.
\end{proof}

In relation now with the Poisson laws, and their higher moments, we first have the following result, dealing with the simplest case, where the parameter is $t=1$:

\index{Bell numbers}
\index{partitions}
\index{moments of Poisson laws}

\begin{theorem}
The moments of $p_1$ are the Bell numbers,
$$M_k(p_1)=|P(k)|$$
where $P(k)$ is the set of partitions of $\{1,\ldots,k\}$.
\end{theorem}

\begin{proof}
The moments of $p_1$ are given by the following formula:
$$M_k=\frac{1}{e}\sum_{n\geq1}\frac{n^k}{n!}$$

We therefore have the following recurrence formula for these moments:
\begin{eqnarray*}
M_{k+1}
&=&\frac{1}{e}\sum_{n\geq1}\frac{n^{k+1}}{n!}\\
&=&\frac{1}{e}\sum_{m\geq0}\frac{(m+1)^k}{m!}\\
&=&\frac{1}{e}\sum_{m\geq0}\frac{m^k}{m!}\left(1+\frac{1}{m}\right)^k\\
&=&\frac{1}{e}\sum_{m\geq0}\frac{m^k}{m!}\sum_{s=0}^k\binom{k}{s}m^{-s}\\
&=&\sum_{s=0}^k\binom{k}{s}\cdot\frac{1}{e}\sum_{m\geq0}\frac{m^{k-s}}{m!}\\
&=&\sum_{s=0}^k\binom{k}{s}M_{k-s}
\end{eqnarray*}

Thus, our moments $M_k$ satisfy the same recurrence as the Bell numbers $B_k=|P(k)|$. Regarding now the initial values, in what concerns the first moment of $p_1$, we have:
$$M_1=\frac{1}{e}\sum_r\frac{r}{r!}=1$$

Also, by using the above recurrence for the numbers $M_k$, we obtain from this:
$$M_2
=\sum_s\binom{1}{s}M_{k-s}
=1+1
=2$$

On the other hand, $B_1=1$ and $B_2=2$. Thus we obtain $M_k=B_k$, as claimed.
\end{proof}

Summarizing, the moments of the Poisson law $p_1$ are not computable. However, it is possible to say many theoretical things about them. We will be back to this.

\bigskip

As a continuation now of the above study, more generally, we have the following result, dealing with the case of the arbitrary Poisson laws, $p_t$ with $t>0$:

\index{number of blocks}
\index{Stirling numbers}

\begin{theorem}
The moments of $p_t$ with $t>0$ are given by
$$M_k(p_t)=\sum_{\pi\in P(k)}t^{|\pi|}$$
where $|.|$ is the number of blocks.
\end{theorem}

\begin{proof}
The moments of the Poisson law $p_t$ with $t>0$ are given by:
$$M_k=e^{-t}\sum_{n\geq1}\frac{t^nn^k}{n!}$$

We have the following recurrence formula for these moments:
\begin{eqnarray*}
M_{k+1}
&=&e^{-t}\sum_{n\geq1}\frac{t^nn^{k+1}}{n!}\\
&=&e^{-t}\sum_{m\geq0}\frac{t^{m+1}(m+1)^k}{m!}\\
&=&e^{-t}\sum_{m\geq0}\frac{t^{m+1}m^k}{m!}\left(1+\frac{1}{m}\right)^k\\
&=&e^{-t}\sum_{m\geq0}\frac{t^{m+1}m^k}{m!}\sum_{s=0}^k\binom{k}{s}m^{-s}\\
&=&\sum_{s=0}^k\binom{k}{s}\cdot e^{-t}\sum_{m\geq0}\frac{t^{m+1}m^{k-s}}{m!}\\
&=&t\sum_{s=0}^k\binom{k}{s}M_{k-s}
\end{eqnarray*}

Regarding now the initial values, the first moment of $p_t$ is given by:
$$M_1
=e^{-t}\sum_r\frac{t^rr}{r!}
=e^{-t}\sum_r\frac{t^r}{(r-1)!}
=t$$

Now by using the above recurrence, we deduce from this that we have:
$$M_2
=t\sum_s\binom{1}{s}M_{k-s}
=t(1+t)
=t+t^2$$

On the other hand, consider the numbers in the statement, namely:
$$S_k=\sum_{\pi\in P(k)}t^{|\pi|}$$

Since a partition of $\{1,\ldots,k+1\}$ appears by choosing $s$ neighbors for $1$, among the $k$ numbers available, and then partitioning the $k-s$ elements left, we have:
$$S_{k+1}=t\sum_s\binom{k}{s}S_{k-s}$$

As for the initial values of these numbers, these are as follows:
$$S_1=t\quad,\quad S_2=t+t^2$$

Thus the initial values coincide, and so by recurrence we obtain that these latter numbers are the moments of the Poisson law $p_t$, as stated.
\end{proof}

As an illustration now for the above result, numerically, we have:

\index{low order moments}

\begin{theorem}
The low order moments of $p_t$ are the numbers
$$M_1=t$$
$$M_2=t+t^2$$
$$M_3=t+3t^2+t^3$$
$$M_4=t+7t^2+6t^3+t^4$$
$$M_5=t+15t^2+25t^3+10t^4+t^5$$
$$\vdots$$
which at $t=1$ are the Bell numbers $1,2,5,15,52,\ldots$.
\end{theorem}

\begin{proof}
We can use here input from the proof of Proposition 4.13, as follows:

\medskip

(1) Obviously we have $B_1=1$, and $M_1=t$.

\medskip

(2) Then $B_2=2$, the partitions being as follows, leading to $M_2=t+t^2$:
$$|\ |\quad,\quad\sqcap$$

(3) Next, $B_3=5$, the partitions being as follows, leading to $M_3=t+3t^2+t^3$:
$$|\ |\ |\quad,\quad\sqcap\ |\quad,\quad\sqcap\hskip-3.2mm{\ }_|\quad,\quad|\ \sqcap\quad,\quad\sqcap\hskip-0.7mm\sqcap$$

(4) At $k=4$ now, things become more complex, and it is better to trick. We can count the partitions up to permutations of the corresponding diagrams, and with this convention made, here are the relevant partitions and their multiplicities, leading to $B_4=15$:
$$|\ |\ |\ |\times1\quad,\quad \sqcap\ |\ |\times 6\quad,\quad \sqcap\sqcap\times 3
\quad,\quad \sqcap\hskip-1.6mm\sqcap|\times4\quad,\quad \sqcap\hskip-1.6mm\sqcap\hskip-1.6mm\sqcap\times1$$

But this leads to the formula in the statement for the fourth moment, namely:
$$M_4=t+7t^2+6t^3+t^4$$

(5) The same method works at $k=5$, with the block distributions and multiplicities, which are simpler to draw than partitions, being as follows, leading to $B_5=52$:
$$11111\to1\ ,\ 2111\to10\ ,\ 221\to15\ ,\ 311\to10\ ,\ 32\to10\ ,\ 41\to5\ ,\ 5\to1$$

But this leads to the formula in the statement for the fifth moment, namely:
$$M_5=t+15t^2+25t^3+10t^4+t^5$$

(6) Next, in what regards the case $k=6$, where $B_6=203$, or even higher, $k=7,8$ and so on, we will leave some computations here as an instructive exercise.

\medskip

(7) Finally, let us mention that the coefficients appearing in the statement are called Stirling numbers. And exercise of course for you, to learn more about them.
\end{proof}

Regarding now the central moments, the result here is as follows: 

\begin{theorem}
The low order central moments of $p_t$ are the numbers
$$M_1'=0$$
$$M_2'=t$$
$$M_3'=t$$
$$M_4'=t+3t^2$$
$$M_5'=t+10t^2$$
$$\vdots$$
which at $t=1$ are the numbers $0,1,1,4,11,\ldots$
\end{theorem}

\begin{proof}
According to our various formulae above, and using $E=t$, we have:
$$M_k'=\sum_{s=0}^k(-1)^s\binom{k}{s}t^sM_{k-s}$$

In particular, at $t=1$ we have the following formula, $B_s$ being the Bell numbers:
$$M_k'=\sum_{s=0}^k(-1)^s\binom{k}{s}B_{k-s}$$

But with this, in practice, we are led to the formulae in the statement.
\end{proof}

Let us end this discussion with the following result, capturing the essentials:

\begin{theorem}
For the Poisson law $p_t$ of parameter $t>0$,
$$E=t\quad,\quad  V=t\quad ,\quad \gamma=\frac{1}{\sqrt{t}}\quad ,\quad \kappa=3+\frac{1}{t}$$
are the mean, variance, skewness and kurtosis.
\end{theorem}

\begin{proof}
Regarding the formulae of $E,V$, we know these from Proposition 4.3. As for the formulae of $\gamma,\kappa$, these follow from Theorem 4.17. Indeed, we first have:
$$\gamma=\frac{M_2'}{V\sqrt{V}}=\frac{t}{t\sqrt{t}}=\frac{1}{\sqrt{t}}$$

As for the kurtosis, again by using Theorem 4.17, this is given by:
$$\kappa=\frac{M_3'}{V^2}=\frac{t+3t^2}{t^2}=3+\frac{1}{t}$$

Thus, we are led to the conclusions in the statement.
\end{proof}

\section*{4c. Derangements}

What is next? Some more advanced combinatorics, I guess, in relation with the Poisson laws. Let us start with the following result, which is something very useful:

\begin{proposition}
We have the following formula,
$$\left|\left(\bigcup_iA_i\right)^c\right|=|A|-\sum_i|A_i|+\sum_{i<j}|A_i\cap A_j|-\sum_{i<j<k}|A_i\cap A_j\cap A_k|+\ldots$$
for any subsets $A_i$ of a certain set $A$, called inclusion-exclusion principle.
\end{proposition}

\begin{proof}
This is indeed quite clear, by thinking a bit, as follows:

\medskip

(1) In order to count $(\cup_iA_i)^c$, we certainly have to start with $|A|$.

\medskip

(2) Then, we obviously have to remove each $|A_i|$, and so remove $\sum_i|A_i|$.

\medskip

(3) But then, we have to put back each $|A_i\cap A_j|$, and so put back $\sum_{i<j}|A_i\cap A_j|$.

\medskip

$\vdots$

\medskip

(4) And so on, which leads to the formula in the statement.
\end{proof}

With this discussed, here is now the application of the inclusion-exclusion principle that we were having in mind, making appear $e$, in a nice combinatorial way:

\index{permutation}
\index{derangement}
\index{fixed points}
\index{random permutation}

\begin{theorem}
The probability for a random permutation $\sigma\in S_N$ to be a derangement, that is, to have no fixed points, is given by the following formula:
$$P=1-\frac{1}{1!}+\frac{1}{2!}-\frac{1}{3!}+\ldots+(-1)^N\frac{1}{N!}$$
Thus we have the estimate $P\simeq1/e$, in the $N\to\infty$ limit.
\end{theorem}

\begin{proof}
Consider the sets $S_N^i=\{\sigma\in S_N|\sigma(i)=i\}$. According to the inclusion-exclusion principle, the probability that we are interested in is given by:
\begin{eqnarray*}
P
&=&\frac{1}{N!}\left(|S_N|-\sum_i|S_N^i|+\sum_{i<j}|S_N^i\cap S_N^j|-\sum_{i<j<k}|S_N^i\cap S_N^j\cap S_N^k|+\ldots\right)\\
&=&\frac{1}{N!}\sum_{k=0}^N(-1)^k\sum_{i_1<\ldots<i_k}|S_N^{i_1}\cap\ldots\cap S_N^{i_k}|\\
&=&\frac{1}{N!}\sum_{k=0}^N(-1)^k\binom{N}{k}(N-k)!\\
&=&\sum_{k=0}^N\frac{(-1)^k}{k!}
\end{eqnarray*}

Since at the end we have the expansion of $1/e$, we obtain the result.
\end{proof}

Let us discuss now, more generally, what happens when counting permutations having exactly $k$ fixed points. The result here, extending Theorem 4.20, is as follows:

\begin{theorem}
The probability for a random permutation $\sigma\in S_N$ to have exactly $k$ fixed points is given by the following formula:
$$P=\frac{1}{k!}\left(1-\frac{1}{1!}+\frac{1}{2!}-\ldots+(-1)^{N-1}\frac{1}{(N-1)!}+(-1)^N\frac{1}{N!}\right)$$
Thus we have the estimate $P\simeq1/(ek!)$, in the $N\to\infty$ limit.
\end{theorem}

\begin{proof}
We already know, from Theorem 4.20, that this formula holds at $k=0$. In the general case now, we have to count the permutations $\sigma\in S_N$ having exactly $k$ points. Since having such a permutation amounts in choosing $k$ points among $1,\ldots,N$, and then permuting the $N-k$ points left, without fixed points allowed, we have:
\begin{eqnarray*}
\#\left\{\sigma\in S_N\Big|\chi(\sigma)=k\right\}
&=&\binom{N}{k}\#\left\{\sigma\in S_{N-k}\Big|\chi(\sigma)=0\right\}\\
&=&\frac{N!}{k!(N-k)!}\#\left\{\sigma\in S_{N-k}\Big|\chi(\sigma)=0\right\}\\
&=&N!\times\frac{1}{k!}\times\frac{\#\left\{\sigma\in S_{N-k}\Big|\chi(\sigma)=0\right\}}{(N-k)!}
\end{eqnarray*}

Now by dividing everything by $N!$, we obtain from this the following formula:
$$\frac{\#\left\{\sigma\in S_N\Big|\chi(\sigma)=k\right\}}{N!}=\frac{1}{k!}\times\frac{\#\left\{\sigma\in S_{N-k}\Big|\chi(\sigma)=0\right\}}{(N-k)!}$$

By using now the computation at $k=0$, from Theorem 4.20, it follows that with $N\to\infty$ we have the following estimate, $\chi$ being the number of fixed points:
$$P(\chi=k)
\simeq\frac{1}{k!}\cdot P(\chi=0)
\simeq\frac{1}{k!}\cdot\frac{1}{e}$$

Thus, we are led to the conclusion in the statement.
\end{proof}

Now observe that Theorem 4.21 reformulates into something very simple, namely:

\begin{theorem}
The number of fixed points of permutations,
$$\chi(\sigma)=\#\left\{i\in\{1,\ldots,N\}\Big|\sigma(i)=i\right\}$$ 
when regarded as random variable $\chi:S_N\to\mathbb N$, satisfies the formula
$$\chi\sim p_1$$
that is, follows the Poisson law of parameter $1$, in the $N\to\infty$ limit.
\end{theorem}

\begin{proof}
We know from Theorem 4.21 that with $N\to\infty$ we have, for any $k\in\mathbb N$:
$$P(\chi=k)\simeq\frac{1}{ek!}$$

But this tells us precisely that $\chi:S_N\to\mathbb N$ is asymptotically Poisson (1), as stated.
\end{proof}

More generally now, and quite remarkably, we have in fact the following result:

\index{symmetric group}
\index{Poisson law}

\begin{theorem}
The number of partial fixed points of permutations,
$$\chi_t(\sigma)=\#\left\{i\in\{1,\ldots,[tN]\}\Big|\sigma(i)=i\right\}$$ 
when regarded as random variable $\chi_t:S_N\to\mathbb N$, satisfies the formula
$$\chi_t\sim p_t$$
that is, follows the Poisson law of parameter $t\in(0,1]$, in the $N\to\infty$ limit.
\end{theorem}

\begin{proof}
We already know from Theorem 4.22 that the result holds at $t=1$. In general, the proof is similar, the idea being as follows:

\medskip

(1) Consider first, as in the proof of Theorem 4.20, the sets $S_N^i=\{\sigma\in S_N|\sigma(i)=i\}$. The set of permutations having no fixed points among $1,\ldots,[tN]$ is then:
$$X_N=\left(\bigcup_{i\leq[tN]}S_N^i\right)^c$$

In order to compute now the cardinality $|X_N|$, consider as well the following sets, depending on indices $i_1<\ldots<i_k$, obtained by taking intersections:
$$S_N^{i_1\ldots i_k}=S_N^{i_1}\cap\ldots\cap S_N^{i_k}$$

As before in the proof of Theorem 4.20, we obtain by inclusion-exclusion that:
\begin{eqnarray*}
P(\chi_t=0)
&=&\frac{1}{N!}\sum_{k=0}^{[tN]}(-1)^k\sum_{i_1<\ldots<i_k<[tN]}|S_N^{i_1\ldots i_k}|\\
&=&\frac{1}{N!}\sum_{k=0}^{[tN]}(-1)^k\sum_{i_1<\ldots<i_k<[tN]}(N-k)!\\
&=&\frac{1}{N!}\sum_{k=0}^{[tN]}(-1)^k\binom{[tN]}{k}(N-k)!\\
&=&\sum_{k=0}^{[tN]}\frac{(-1)^k}{k!}\cdot\frac{[tN]!(N-k)!}{N!([tN]-k)!}
\end{eqnarray*}

But with $N\to\infty$, we obtain from this the following estimate, as desired:
$$P(\chi_t=0)
\simeq\sum_{k=0}^{[tN]}\frac{(-1)^k}{k!}\cdot t^k
=\sum_{k=0}^{[tN]}\frac{(-t)^k}{k!}
\simeq e^{-t}$$

(2) More generally now, by counting the permutations $\sigma\in S_N$ having exactly $k$ fixed points among $1,\ldots,[tN]$, as in the proof of Theorem 4.21, our claim is that we get:
$$P(\chi_t=k)\simeq\frac{t^k}{k!e^t}$$

Indeed, we already know from (1) that this formula holds at $k=0$. In general now, we have to count the permutations $\sigma\in S_N$ having exactly $k$ fixed points among $1,\ldots,[tN]$. Since having such a permutation amounts in choosing $k$ points among $1,\ldots,[tN]$, and then permuting the $N-k$ points left, without fixed points among $1,\ldots,[tN]$ allowed, we obtain the following formula, where $s\in(0,1]$ is such that $[s(N-k)]=[tN]-k$:
\begin{eqnarray*}
\#\left\{\sigma\in S_N\Big|\chi_t(\sigma)=k\right\}
&=&\binom{[tN]}{k}\#\left\{\sigma\in S_{N-k}\Big|\chi_s(\sigma)=0\right\}\\
&=&\frac{[tN]!}{k!([tN]-k)!}\#\left\{\sigma\in S_{N-k}\Big|\chi_s(\sigma)=0\right\}\\
&=&\frac{1}{k!}\times\frac{[tN]!(N-k)!}{([tN]-k)!}\times\frac{\#\left\{\sigma\in S_{N-k}\Big|\chi_s(\sigma)=0\right\}}{(N-k)!}
\end{eqnarray*}

Now by dividing everything by $N!$, we obtain from this the following formula:
$$\frac{\#\left\{\sigma\in S_N\Big|\chi_t(\sigma)=k\right\}}{N!}=\frac{1}{k!}\times\frac{[tN]!(N-k)!}{N!([tN]-k)!}\times\frac{\#\left\{\sigma\in S_{N-k}\Big|\chi_s(\sigma)=0\right\}}{(N-k)!}$$

By using now the computation at $k=0$, that we already have, from (1) above, it follows that with $N\to\infty$ we have the following estimate:
\begin{eqnarray*}
P(\chi_t=k)
&\simeq&\frac{1}{k!}\times\frac{[tN]!(N-k)!}{N!([tN]-k)!}\cdot P(\chi_s=0)\\
&\simeq&\frac{t^k}{k!}\cdot P(\chi_s=0)\\
&\simeq&\frac{t^k}{k!}\cdot\frac{1}{e^s}
\end{eqnarray*}

Now recall that the parameter $s\in(0,1]$ was chosen in the above such that:
$$[s(N-k)]=[tN]-k$$

Thus in the $N\to\infty$ limit we have $s=t$, and so we obtain, as claimed:
$$P(\chi_t=k)\simeq\frac{t^k}{k!}\cdot\frac{1}{e^t}$$

It follows that we obtain in the limit a Poisson law of parameter $t$, as stated.
\end{proof}

And with this, end of the story? Not yet, because with a bit of group theory know-how, we can formulate the following definition, which is something quite conceptual:

\begin{definition}
Given a group of matrices $G\subset O_N$, the function
$$\chi:G\to\mathbb C\quad,\quad 
\chi(g)=\sum_{i=1}^Ng_{ii}$$
is called main character of $G$. More generally, the function
$$\chi_t:G\to\mathbb C\quad,\quad 
\chi_t(g)=\sum_{i=1}^{[tN]}g_{ii}$$
is called main truncated character of $G$, of parameter $t\in(0,1]$.
\end{definition}

These notions are quite interesting in relation with group theory, and the point now is that, in terms of them, our main result, Theorem 4.23, reformulates as follows:

\begin{theorem}
For the symmetric group $S_N\subset O_N$ we have
$$\chi_t\sim p_t$$
in the $N\to\infty$ limit, for any $t\in(0,1]$.
\end{theorem} 

\begin{proof}
As we already know from chapter 3, the standard coordinates on $S_N\subset O_N$, which are the coefficients of the permutation matrices, are given by:
$$g_{ij}=\delta_{\sigma(j)i}$$

Thus the sums of diagonal coordinates $g_{ii}$ count the fixed points, and with this observation in hand, Theorem 4.23 reformulates as in the statement. Great.
\end{proof}

\section*{4d. Compound Poisson}

As a continuation of the above material, Definition 4.24 and Theorem 4.25 suggest doing some computations for other interesting groups $G\subset O_N$ that we know. And here, I don't know about you, but my favorite example is the hyperoctahedral group:

\begin{theorem}
Consider the hyperoctahedral group $H_N$, which appears as the symmetry group of the $N$-cube, or the symmetry group of the $N$ coordinate axes of $\mathbb R^N$:
$$S_N\subset H_N\subset O_N$$
In matrix terms, $H_N$ consists of the permutation-type matrices having $\pm1$ as nonzero entries, and we have a wreath product decomposition as follows:
$$H_N=\mathbb Z_2\wr S_N$$
In this picture, the main character counts the signed number of fixed points, among the coordinate axes, and its truncations count the truncations of such numbers.
\end{theorem}

\begin{proof}
Many things going on here, the idea being as follows:

\medskip

(1) To start with, we can certainly talk about the hyperoctahedral group $H_N$, as being the symmetry group of the $N$-cube, and in small dimensions we have $H_1=\mathbb Z_2$, obvious, then $H_2=D_4=\mathbb Z_4\rtimes\mathbb Z_2$, the symmetry group of the square, and then $H_3=S_4\times\mathbb Z_2$, with the copy of $S_4$ best viewed as permuting the main diagonals of the cube.

\medskip

(2) Next, by centering the $N$-cube at the origin, its symmetry group $H_N$ appears as the symmetry group of the $N$ coordinate axes of $\mathbb R^N$. Which is a quite fruitful viewpoint, proving right away $H_N=M_N(-1,0,1)\cap O_N$, and giving as well the cardinality formula $|H_N|=2^NN!$, which after a bit more work gives $H_N=\mathbb Z_2\wr S_N$, as stated.

\medskip

(3) Finally, in what regards the assertions at the end, concerning the main character and its truncations, these are clear from definitions, exactly as previously for $S_N$.
\end{proof}

Regarding now the character laws, we can compute them by using the same method as for the symmetric group $S_N$, namely inclusion-exclusion, and we have:

\index{main character}

\begin{theorem}
For the hyperoctahedral group $H_N\subset O_N$, the law of the variable
$$\chi_t=\sum_{i=1}^{[tN]}u_{ii}$$
with $u_{ij}$ being the standard coordinates, becomes in the $N\to\infty$ limit the measure
$$b^2_t=e^{-t}\sum_{k=-\infty}^\infty\delta_k\sum_{p=0}^\infty \frac{(t/2)^{|k|+2p}}{(|k|+p)!p!}$$ 
where $\delta_k$ is the Dirac mass at $k\in\mathbb Z$.
\end{theorem}

\begin{proof}
We can regard the group $H_N$ as being the symmetry group of the graph $I_N=\{I^1,\ldots ,I^N\}$ formed by $N$ segments. The diagonal coefficients are given by:
$$u_{ii}(g)=\begin{cases}
\ 0\ \mbox{ if $g$ moves $I^i$}\\
\ 1\ \mbox{ if $g$ fixes $I^i$}\\
-1\mbox{ if $g$ returns $I^i$}
\end{cases}$$

We denote by $\uparrow g,\downarrow g$ the number of segments among $\{I^1,\ldots ,I^s\}$ which are fixed, respectively returned by an element $g\in H_N$. With this notation, we have:
$$u_{11}+\ldots+u_{ss}=\uparrow g-\downarrow g$$

Let us denote by $P_N$ probabilities computed over the group $H_N$. The density of the law of $u_{11}+\ldots+u_{ss}$ at a point $k\geq 0$ is then given by the following formula:
\begin{eqnarray*}
D(k)
&=&P_N(\uparrow g-\downarrow g=k)\\
&=&\sum_{p=0}^\infty P_N(\uparrow g=k+p, \downarrow g=p)
\end{eqnarray*}

Assume first that we have $t=1$. We use the fact, that we know well from before, that the probability of $\sigma\in S_N$ to have no fixed points is asymptotically given by:
$$P_0=\frac{1}{e}$$

Thus the probability of $\sigma\in S_N$ to have $m$ fixed points is asymptotically given by:
$$P_m=\frac{1}{em!}$$

In terms of probabilities over $H_N$, we obtain from this, as desired:
\begin{eqnarray*}
\lim_{N\to\infty}D(k)
&=&\lim_{N\to\infty}\sum_{p=0}^\infty(1/2)^{k+2p}\begin{pmatrix}k+2p\\ k+p\end{pmatrix} P_N(\uparrow g+\downarrow g=k+2p)\\ 
&=&\sum_{p=0}^\infty(1/2)^{k+2p}\begin{pmatrix}k+2p\\
k+p\end{pmatrix}\frac{1}{e(k+2p)!}\\
&=&\frac{1}{e}\sum_{p=0}^\infty \frac{(1/2)^{k+2p}}{(k+p)!p!}
\end{eqnarray*}

As for the general case $0<t\leq 1$, here the result follows by performing some modifications in the above computation. The asymptotic density is computed as follows:
\begin{eqnarray*}
\lim_{N\to\infty}D(k)
&=&\lim_{N\to\infty}\sum_{p=0}^\infty(1/2)^{k+2p}\begin{pmatrix}k+2p\\ k+p\end{pmatrix} P_N(\uparrow g+\downarrow g=k+2p)\\
&=&\sum_{p=0}^\infty(1/2)^{k+2p}\begin{pmatrix}k+2p\\
k+p\end{pmatrix}\frac{t^{k+2p}}{e^t(k+2p)!}\\
&=&e^{-t}\sum_{p=0}^\infty \frac{(t/2)^{k+2p}}{(k+p)!p!}
\end{eqnarray*}

Together with $D(-k)=D(k)$, this gives the formula in the statement.
\end{proof}

The above result is quite interesting, and based on it, we can formulate:

\index{Bessel function}
\index{Bessel law}

\begin{definition}
The Bessel law of parameter $t>0$ is the measure
$$b^2_t=e^{-t}\sum_{k=-\infty}^\infty\delta_k\,f_k(t/2)$$
with the density being the following function, depending on $k\in\mathbb Z$,
$$f_k(t)=\sum_{p=0}^\infty \frac{t^{|k|+2p}}{(|k|+p)!p!}$$
which is the Bessel function of the first kind, from classical analysis.
\end{definition}

Let us study now these Bessel laws that we found. In analogy with what we know about the Poisson laws $p_t$, from before, we first have the following result:

\index{convolution}
\index{convolution semigroup}

\begin{theorem}
The Bessel laws $b^2_t$ have the property
$$b^2_s*b^2_t=b^2_{s+t}$$
so they form a truncated one-parameter semigroup with respect to convolution.
\end{theorem}

\begin{proof}
We use the formula in Definition 4.28. The Fourier transform is:
$$F(x)=e^{-t}\sum_{k=-\infty}^\infty e^{ikx}\,f_k(t/2)$$

We can compute the derivative of $F$ with respect to $t$, as follows:
\begin{eqnarray*}
F(x)'
&=&-e^{-t}\sum_{k=-\infty}^\infty e^{ikx}\,f_k(t/2)+\frac{e^{-t}}{2}\sum_{k=-\infty}^\infty e^{ikx}\,f_k'(t/2)\\
&=&-F(x)+\frac{e^{-t}}{2}\sum_{k=-\infty}^\infty e^{ikx}\,f_k'(t/2)
\end{eqnarray*}

On the other hand, the derivative of $f_k$ with $k\geq 1$ is given by:
\begin{eqnarray*}
f_k'(t)
&=&\sum_{p=0}^\infty \frac{(k+2p)t^{k+2p-1}}{(k+p)!p!}\\
&=&\sum_{p=0}^\infty \frac{(k+p)t^{k+2p-1}}{(k+p)!p!}+\sum_{p=0}^\infty\frac{p\,t^{k+2p-1}}{(k+p)!p!}\\
&=&\sum_{p=0}^\infty \frac{t^{k+2p-1}}{(k+p-1)!p!}+\sum_{p=1}^\infty\frac{t^{k+2p-1}}{(k+p)!(p-1)!}\\
&=&\sum_{p=0}^\infty \frac{t^{(k-1)+2p}}{((k-1)+p)!p!}+\sum_{p=1}^\infty\frac{t^{(k+1)+2(p-1)}}{((k+1)+(p-1))!(p-1)!}\\
&=&f_{k-1}(t)+f_{k+1}(t)
\end{eqnarray*}

This computation works in fact for any $k\in\mathbb Z$, so we get:
\begin{eqnarray*}
F(x)'
&=&-F(x)+\frac{e^{-t}}{2}
\sum_{k=-\infty}^\infty e^{ikx} (f_{k-1}(t/2)+f_{k+1}(t/2))\\
&=&-F(x)+\frac{e^{-t}}{2} \sum_{k=-\infty}^\infty
e^{i(k+1)x}f_{k}(t/2)+e^{i(k-1)x}f_{k}(t/2)\\
&=&-F(x)+\frac{e^{ix}+e^{-ix}}{2}\,F(x)\\
&=&\left(\frac{e^{ix}+e^{-ix}}{2}-1\right)F(x)
\end{eqnarray*}

Thus, the Fourier transform is subject to the following equation:
$$\frac{F(x)'}{F(x)}=\frac{e^{ix}+e^{-ix}}{2}-1$$

Now by integrating, it follows that the Fourier transform is given by:
$$F(x)=\exp\left(\left(\frac{e^{ix}+e^{-ix}}{2}-1\right)t\right)$$

Thus the log of the Fourier transform is linear in $t$, and we get the result.
\end{proof}

In order to unify now the theory of the Poisson and Bessel laws, we have the following notion, extending the Poisson limit theory from the beginning of this chapter:

\index{compound Poisson law}

\begin{definition}
Associated to any compactly supported positive measure $\nu$ on $\mathbb C$ is the probability measure
$$p_\nu=\lim_{n\to\infty}\left(\left(1-\frac{c}{n}\right)\delta_0+\frac{1}{n}\nu\right)^{*n}$$
where $c=mass(\nu)$, called compound Poisson law.
\end{definition}

In other words, what we are doing here is to generalize the construction in the Poisson Limit Theorem, by allowing the only parameter there, which was the positive real number $t>0$, to be replaced by a certain probability measure $\nu$, of arbitrary mass $c>0$.

\bigskip

In what follows we will be mostly interested in the case where $\nu$ is discrete, as is for instance the case for the measure $\nu=t\delta_1$ with $t>0$, which produces the Poisson laws. To be more precise, the Poisson limit theorem from before reformulates as:

\begin{theorem}[PLT]
The Poisson laws $p_t$ are compound Poisson laws,
$$p_t=p_{t\varepsilon}$$
with $\varepsilon=\delta_1$, Dirac mass at $1$.
\end{theorem}

\begin{proof}
We recall from before that the PLT states that we have:
$$\lim_{n\to\infty}\left(\left(1-\frac{t}{n}\right)\delta_0+\frac{t}{n}\,\delta_1\right)^{*n}=p_t$$

Thus, we are led to the conclusion in the statement.
\end{proof}

Getting back now to our comments on Definition 4.30, as stated, as mentioned in the above, we will be mostly interested here in the case where the measure $\nu$ is discrete. In fact, we will be mainly interested in the case where $\nu$ is a multiple of the uniform measure on the $s$-th roots of unity. But more on this later, once we will be more advanced.

\bigskip

Finally, as a last comment on Definition 4.30, we allow there the measure $\nu$ to be complex, instead of being real, and this precisely in order to have as examples, later, the multiples of the uniform measures on the $s$-th roots of unity. But, more on this later.

\bigskip

Getting now to some theory, for the compound Poisson laws as introduced in Definition 4.30, we first have the following result, allowing us to detect such laws:

\index{Fourier transform}

\begin{proposition}
For a discrete measure, $\nu=\sum_{i=1}^sc_i\delta_{z_i}$ with $c_i>0$ and $z_i\in\mathbb C$, we have the formula
$$F_{p_\nu}(x)=\exp\left(\sum_{i=1}^sc_i(e^{ixz_i}-1)\right)$$
where $F$ denotes as usual the Fourier transform.
\end{proposition}

\begin{proof}
In order to prove the formula in the statement, consider the measure $\mu_n$ appearing in Definition 4.30, under the convolution sign, namely:
$$\mu_n=\left(1-\frac{c}{n}\right)\delta_0+\frac{1}{n}\nu$$

We have the following computation, in the context of Definition 4.30:
\begin{eqnarray*}
&&F_{\mu_n}(x)=\left(1-\frac{c}{n}\right)+\frac{1}{n}\sum_{i=1}^sc_ie^{ixz_i}\\
&\implies&F_{\mu_n^{*n}}(x)=\left(\left(1-\frac{c}{n}\right)+\frac{1}{n}\sum_{i=1}^sc_ie^{ixz_i}\right)^n\\
&\implies&F_{p_\nu}(x)=\exp\left(\sum_{i=1}^sc_i(e^{ixz_i}-1)\right)
\end{eqnarray*}

Thus, we have obtained the formula in the statement.
\end{proof}

Next, we have the following key result, providing an alternative to Definition 4.30, and which will be our formulation here of the Compound Poisson Limit Theorem:

\index{compound Poisson Limit theorem}
\index{CPLT}
\index{compound PLT}
\index{compound Poisson limit}

\begin{theorem}[CPLT]
For a discrete measure, $\nu=\sum_{i=1}^sc_i\delta_{z_i}$ with $c_i>0$ and $z_i\in\mathbb C$, we have the following formula,
$$p_\nu={\rm law}\left(\sum_{i=1}^sz_i\alpha_i\right)$$
with the variables $\alpha_i$ being Poisson $(c_i)$, and independent.
\end{theorem}

\begin{proof}
Let $\alpha$ be the sum of Poisson variables in the statement:
$$\alpha=\sum_{i=1}^sz_i\alpha_i$$

By using the Fourier transform formulae from before, we have:
\begin{eqnarray*}
F_{\alpha_i}(x)=\exp(c_i(e^{ix}-1))
&\implies&F_{z_i\alpha_i}(x)=\exp(c_i(e^{ixz_i}-1))\\
&\implies&F_\alpha(x)=\exp\left(\sum_{i=1}^sc_i(e^{ixz_i}-1)\right)
\end{eqnarray*}

Thus we have the same formula as in Proposition 4.32, as desired.
\end{proof}

Getting back now to the Bessel laws, we have the following result:

\index{Bessel law}

\begin{theorem}
The Bessel laws $b^2_t$ are compound Poisson laws, given by
$$b^2_t=p_{t\varepsilon}$$
where $\varepsilon=\frac{1}{2}(\delta_{-1}+\delta_1)$ is the uniform measure on $\mathbb Z_2$.
\end{theorem}

\begin{proof}
We recall from the proof of Theorem 4.29 that the Fourier transform of the Bessel law $b^2_t$, with respect to a variable $x$, is given by the following formula:
$$F(x)=\exp\left(\left(\frac{e^{ix}+e^{-ix}}{2}-1\right)t\right)$$

But, according to Proposition 4.32, this is exactly the Fourier transform of the compound Poisson law $p_{t\varepsilon}$ in the statement. Thus, we have $b^2_t=p_{t\varepsilon}$, as claimed.
\end{proof}

Quite nice all this, and as a continuation, we can talk, more generally, about complex reflection groups $G\subset U_N$, which according to a theorem of Shephard and Todd consist of a two-parameter series $H_N^{sd}$, and a number of exceptional examples. And, in what regards the series $H_N^{sd}$, the character computations can be done, a bit as before for $S_N,H_N$, with the conclusion that the parameter $d$ is irrelevant, and that we obtain as limiting laws the compound Poisson laws $p_{t\varepsilon}$, with $\varepsilon$ being the uniform measure on the $s$-th roots of unity. For more on all this, story and details, you can check my previous book \cite{ba1}.

\section*{4e. Exercises}

This was a key theoretical chapter, and as exercises on this, we have:

\begin{exercise}
Learn more about the precise convergence in the PLT.
\end{exercise}

\begin{exercise}
Work out some asymptotics for the Bell numbers.
\end{exercise}

\begin{exercise}
Learn more about the Stirling numbers, and their properties.
\end{exercise}

\begin{exercise}
Learn more about permutation matrices, and their properties.
\end{exercise}

\begin{exercise}
Learn about the signature, and the alternating group $A_N\subset S_N$.
\end{exercise}

\begin{exercise}
Do the character computations for the alternating group $A_N\subset S_N$.
\end{exercise}

\begin{exercise}
Learn about the classification of the complex reflection groups.
\end{exercise}

\begin{exercise}
Do the probabilistic computations for the groups $H_N^{sd}$.
\end{exercise}

As bonus exercise, and no surprise here, read some group theory.

\part{Continuous laws}

\ \vskip50mm

\begin{center}
{\em Cuando sei Maria Dolores

Cuando sei quei mal d'amore

Cuando sei quei mal a su vera

Cuando sei me va al dottore}
\end{center}

\chapter{Variables, laws}

\section*{5a. Measured spaces}

Welcome to probability theory, take two. When looking around, in the real life, phenomena can be of both discrete and continuous nature, and our goal here, in this chapter, will be that of axiomatizing the general case, covering both the discrete and continuous situations. And then afterwards, in the remainder of the book, we will have a detailed look at the continuous case, phenomena and distributions appearing there.

\bigskip

So, how to get started? We already know some good mathematics in the discrete case, and we will start with a review of the basics there, that we know from chapter 2. In a very compacted form, our axiomatic knowledge from chapter 2 is as follows:

\index{probability space}
\index{discrete space}
\index{events}
\index{probability function}
\index{probability 1}

\begin{fact}
Discrete probability theory can be axiomatized as follows:
\begin{enumerate}
\item A discrete probability space is a set $X$ with a probability measure on it, $\nu(Y)=\sum_{y\in Y}P(y)$ for $Y\subset X$, with $P:X\to[0,\infty)$ satisfying $\sum_{x\in X}P(x)=1$.

\item A discrete random variable is a function $f:X\to\mathbb R$, and the law of such a variable is the real probability measure $\mu=f_*\nu$, push-forward of $\nu$ by $f$.
\end{enumerate}
\end{fact}

All this is of course, and we insist, something very compact, and in practice, many comments can be made on both (1) and (2). We refer to the material in chapter 2 for more on all this, and to the material from chapters 3-4, for illustrations and more.

\bigskip

Getting now to what we wanted to do, general axiomatization of probability, covering both the discrete and continuous situations, what we have to do is, obviously, to find the correct generalizations of (1) and (2) in Fact 5.1. And since both that (1) and (2) are about measures, we will attempt to shoot both rabbits at the same time:

\begin{question}
What is a measured space? Because, with an answer to this, we can easily axiomatize both the probability spaces, and the random variable laws.
\end{question}

As a first comment here, you might wonder why I chose to talk about general measured spaces, instead of the probability spaces, that is, the measured spaces having mass 1, that we are truly interested in. Well, this is because basic spaces like $\mathbb R$ itself, and more generally $\mathbb R^N$, that we need all the time when doing probability computations, are measured spaces, but not of mass 1. So, we want in fact to shoot 3 rabbits at a time, the third one being the clarification of integration theory over spaces like $\mathbb R$, or $\mathbb R^N$.

\bigskip

Which sounds quite good, nice and ambitious plan that we have here. This being said I must admit that, although being a rural being for years, here in France, and an Eastern European before, I still get my meat from the supermarket, not being that familiar with hunting. So, time I guess to ask the expert. And, here is what she says:

\begin{cat}
Sounds good your plan, to shoot 3 rabbits at the same time. I would start with some training, catching a small mouse would be not bad.
\end{cat}

Humm, okay cat, so shall I understand that all this will be quite difficult. Nevermind. So, thanks for the advice, but I really have this feeling that I am quite good at mathematics, and that I can shoot my 3 rabbits, as planned. I'll just go for that.

\bigskip

Getting started now with our study, as a first difficulty on the way, we have:

\begin{proposition}
In the continuous setting, some sets are not measurable.
\end{proposition}

\begin{proof}
Obviously, this cannot be something rigorous, because we have not talked yet about the precise definition of the measured spaces. But listen to what I have to say, all this will be very useful, precisely for this goal, axiomatizing the measured spaces. Consider the unit circle $\mathbb T$, measured as usual, the total mass being $\mu(\mathbb T)=2\pi$, and then consider the group $G$ consisting of the rational rotations of this unit circle $\mathbb T$:
$$G=\left\{z\to e^{2r\pi i}z\Big|r\in\mathbb Q\cap[0,1)\right\}$$

The action $G\curvearrowright\mathbb T$ partitions then $\mathbb T$ into a certain number of orbits, which is uncountable, and this because $G$ is obviously countable, while $\mathbb T$ itself is uncountable. Now if we choose $A\subset\mathbb T$ containing exactly one point from each orbit, we have:
$$\mathbb T=\bigsqcup_{r\in\mathbb Q\cap[0,1)}e^{2r\pi i}A$$  

But this shows that $A$ is not measurable. Indeed, with the convention that we list our group as $\mathbb Q\cap[0,1)=\{r_1,r_2,r_3,\ldots\}$, we have two possible cases, as follows:

\medskip

(1) Assuming $\mu(A)=c>0$, we would get from this $\mu(\mathbb T)=\infty$, which is contradiction, as shown by the following computation, which is something quite self-explanatory:
$$\mu(\mathbb T)=\sum_{n=1}^\infty\mu(e^{2r_n\pi i}A)=\sum_{n=1}^\infty c=c\cdot\infty=\infty$$

(2) Assuming $\mu(A)=0$, we would get from this $\mu(\mathbb T)=0$, which is a contradiction again, as shown by the following estimate, self-explanatory too, valid for any $\varepsilon>0$:
$$\mu(\mathbb T)=\sum_{n=1}^\infty\mu(e^{2r_n\pi i}A)<\sum_{n=1}^\infty\frac{\varepsilon}{2^n}=\varepsilon$$

Summarizing, save for some discussion of the axiomatics, and more on this in a moment, we are led, just via some common sense, to the conclusion in the statement.
\end{proof}

In view of Proposition 5.4, as a first task in our axiomatization work, we must start with a set $X$, and come up with some axioms for the measurable subsets $Y\subset X$. And here, a bit of thinking and common sense lead to the following definition:

\index{measurable set}

\begin{definition}
A measurable space is a set $X$, given with a set of subsets $M\subset P(X)$, called measurable sets, which form an algebra, in the sense that:
\begin{enumerate}
\item $\emptyset,X\in M$.

\item $E\in M\implies E^c\in M$.

\item $M$ is stable under countable unions and intersections.
\end{enumerate}
\end{definition}

Observe that some of the axioms above are redundant, because assuming (2) we have $\emptyset\in M\iff X\in M$, which can help in verifying (1). The same goes for (3), with only one of the conditions there being in need to be verified, due to $\left(\bigcup_iE_i\right)^c=\bigcap_iE_i^c$.

\bigskip

At the level of examples, these usually come via the following simple fact:

\begin{proposition}
Given a set of subsets $S\subset P(X)$, there is a smallest algebra
$$M=\bar{S}$$
containing it, called algebra generated by $X$.
\end{proposition}

\begin{proof}
This can be viewed in two possible ways, as follows:

\medskip

(1) According to the axioms in Definition 5.5, what we have to do in order to construct $M=\bar{S}$ is to start with $S$, then add $\emptyset,X$ to it, along with the complements $E^c$ of all the sets $E\in S$, and then take countable unions and intersections of such sets. And, some elementary verifications show that what we get in this way is indeed an algebra.

\medskip

(2) Alternatively, we can define $M=\bar{S}$ as being the intersection of all algebras containing $S$, and with the remark that we have at least one such algebra, namely $P(X)$ itself. It is then clear from definitions that $M$ is an algebra, as desired.
\end{proof}

Getting now to the concrete examples of measurable spaces, we have:

\index{Borel set}

\begin{theorem}
Any metric space $X$ is automatically a measurable space, with the algebra of measurable sets being 
$$B=\bar{O}$$
that is, the smallest algebra containing the open sets, called Borel algebra of $X$.
\end{theorem}

\begin{proof}
Many things can be said here, the idea being as follows:

\medskip

(1) To start with, a metric space is by definition a set $X$ with a distance function $d:X\times X\to[0,\infty)$, subject to the following conditions, which are all quite natural:

\medskip

--  $d(x,y)>0$ if $x\neq y$, and $d(x,x)=0$.

\smallskip

-- $d(x,y)=d(y,x)$.

\smallskip

-- $d(x,y)\leq d(x,z)+d(y,z)$.

\medskip

(2) As a basic example, we have $\mathbb R^N$, as well as any of its subsets $X\subset\mathbb R^N$. Indeed, the first two axioms are clear, and for the third axiom, we must prove that:
$$\sqrt{\sum_i(a_i+b_i)^2}\leq\sqrt{\sum_i a_i^2}+\sqrt{\sum_i b_i^2}$$

But this latter inequality is one of the many equivalent formulations of the Cauchy-Schwarz inequality, that we know well from chapter 2, as corresponding to $V(f)\geq 0$.

\medskip

(3) Next, once we have a metric space as in (1), we can talk about open and closed sets, in the obvious way, by saying that $U\subset X$ is open if for any $x\in U$ we have a ball $B_x(r)\subset U$, and that $F\subset X$ is called closed if its complement $F^c\subset X$ is open.

\medskip

(4) Getting now to what the statement says, this is something more or less trivial, based on Proposition 5.6. Thus, we have our measurable space, as desired.

\medskip

(5) Finally, observe that the Borel sets include all open sets, all closed sets, as well as all countable unions of closed sets, and all countable intersections of open sets. As an example here, in the case $X=\mathbb R$, all kinds of intervals are Borel sets:
$$(a,b),\ [a,b],\ (a,b],\ [a,b)\in B$$

Indeed, the first interval is open, and the second one is closed, so these are certainly Borel sets. As for the third and fourth intervals, these appear as countable unions of closed intervals, or as countable intersections of open intervals, so they are Borel too.
\end{proof}

\bigskip

Getting back now to the general case, as a complement to Definition 5.5, we have:

\begin{definition}
Given a measurable space $(X,M)$, a measure on it is a function $\mu:M\to[0,\infty]$ which is countably additive, in the sense that
$$\mu\left(\bigcup_{i=1}^\infty E_i\right)=\sum_{i=1}^\infty\mu(E_i)$$
for any countable family of disjoint measurable sets $E_i\in M$. In this case, we say that $(X,M,\mu)$ is a measured space.
\end{definition}

Observe that we did not assume that our measures are finite, and this for including spaces like $X=\mathbb R$. In the case where $\mu$ happens to be bounded, up to a rescaling we can assume $\mu(X)=1$, and we say in this case that we have a probability measure on $X$.

\bigskip

Looking at what we have so far, Definition 5.5, Theorem 5.7 and Definition 5.8, many natural questions appear, and leaving aside anything too specialized, we are led to:

\begin{question}
Given a metric space $X$, such as $X=\mathbb R$, or $X=\mathbb R^N$, how to construct a measure on it? Also, once we have such a measure, how to integrate the functions $f:X\to\mathbb R$, or $f:X\to\mathbb C$, with respect to this measure?
\end{question}

As you can see, we have two questions here, and none is trivial:

\medskip

(1) Indeed, regarding the first question, this is something that we do not know yet how to solve, even for simple spaces like $X=\mathbb R$. Indeed, while we certainly know how to measure the real intervals, by setting $\mu(a,b)=b-a$, and then unions of such intervals too, nothing guarantees that we can measure any Borel set $E\in B$, in this way. 

\medskip

(2) As for the second question, experience with the usual Riemann integral, that we know well from calculus, shows that such things can be quite tricky too. Observe also that this second question is more general than the first one, because we have $\mu(E)=\int\chi_E$ for any $E\in M$, so if we know how to integrate, we know how to measure. 

\medskip

In short, we face non-trivial questions here, and we need a plan. Here that plan is:

\begin{plan}
We will jointly develop measure and integration theory, as follows:
\begin{enumerate}
\item We will first keep staying abstract, and understand how the functions $f:X\to\mathbb R$, or $f:X\to\mathbb C$, can be integrated, with respect to an abstract measure.

\item With this understood, we will conclude that the integration over $X=\mathbb R$ can only be a straightforward extension of the usual Riemann integration.

\item Which, once done, will enable us to both measure all the Borel sets $E\subset\mathbb R$, and integrate the measurable functions $f:\mathbb R\to\mathbb C$.

\item Finally, we will discuss how to deal with $X=\mathbb R^N$ too, and with other product spaces $X=Y\times Z$, both measure theory and integration.
\end{enumerate}
\end{plan}

Getting started now, let us first talk about measurable functions, in the general context of Definition 5.5, with no measure involved. We have here the following notion:

\begin{definition}
Given a measurable space $X$, and a metric space $Y$, a function $f:X\to Y$ is called measurable when it satisfies the following condition:
$$U\in O\implies f^{-1}(U)\in M$$
When $X$ comes with a measure, we also call such functions integrable.
\end{definition}

Obviously, this is something simplified, the idea being that, in view of Theorem 5.7, the case where the target space $Y$ is a metric space is the one that really matters. In practice, many things can be said about the measurable functions, and we first have:

\begin{proposition}
The measurable functions have the following properties:
\begin{enumerate}
\item If $f:X\to Y$ is measurable and $g:Y\to Z$ is continuous, $g\circ f$ is measurable.

\item If $f,g:X\to\mathbb R$ are both measurable and $h:\mathbb R^2\to Y$ is continuous, then the function $x\to h(f(x),g(x))$ is measurable. 

\item $f:X\to\mathbb C$ is measurable precisely when $Re(f),Im(f):X\to\mathbb R$ are measurable. In this case, the function $|f|:X\to\mathbb R$ is measurable too.

\item If $f,g:X\to\mathbb C$ are measurable, then so are $f+g,fg:X\to\mathbb C$. 
\end{enumerate}
\end{proposition}

\begin{proof}
This is something very standard, the idea being as follows:

\medskip

(1) This is something quite clear, say easy exercise for you, because we have:
$$U\in O\implies g^{-1}(U)\in O\implies f^{-1}(g^{-1}(U))\in M$$

(2) By using (1), it is enough to check that the function $k(x)=(f(x),g(x))$ is measurable. But this follows by writing any open set $U\subset\mathbb R^2$ as a union of rectangles. Indeed, the preminage of each rectangle $R=I\times J$ is measurable, as shown by:
$$k^{-1}(R)
=(f,g)^{-1}(I\times J)
=f^{-1}(I)\cap g^{-1}(J)
\in M$$

But then, with this in hand, since any open set $U\subset\mathbb R^2$ can be written as a  union of rectangles, it follows that $k^{-1}(U)$ is measurable, as desired.

\medskip

(3) This follows indeed by using (2), with the following functions:
$$h(z)=z,Re(z),Im(z),|z|$$

(4) In the real case, $f,g:X\to\mathbb R$, the result follows by using (2), with:
$$h(x,y)=x+y,xy$$

Then, the result can be extended to the complex case, $f,g:X\to\mathbb C$, by using (3).
\end{proof}

Coming next, we have the following useful characterization of the real measurable functions, that we will use many times, in what follows:

\begin{proposition}
A function $f:X\to[-\infty,\infty]$ is measurable precisely when 
$$f^{-1}(I)\in M$$
for any interval of type $I=(a,\infty]$, with $a\in\mathbb R$. 
\end{proposition}

\begin{proof}
Consider the following set, which is easily seen to be an algebra:
$$\Omega=\left\{E\subset[-\infty,\infty]\Big|f^{-1}(E)\in M\right\}$$

We want to prove that $\Omega$ contains all open sets, and this can be done as follows:

\medskip

(1) Pick $a\in\mathbb R$, and then pick an increasing sequence $a_n\to a$. We have then the following equality, which shows that we have $[-\infty,a)\in\Omega$:
$$[-\infty,a)=\bigcup_n[-\infty,a_n)=\bigcup_n(a_n,\infty]^c\in\Omega$$

(2) But with this in hand, we obtain, for any $a<b$, that we have:
$$(a,b)=[-\infty,b)\cap(a,\infty]\in\Omega$$

(3) Now since any open set $U\subset[-\infty,\infty]$ can be written as a union of open intervals, we conclude that we have $U\in\Omega$, as desired.
\end{proof}

Getting now to limits of measurable functions, we have the following result:

\begin{proposition}
The measurable functions $f:X\to[-\infty,\infty]$ are as follows:
\begin{enumerate}
\item If $f_n$ are measurable, so are $g=\sup_nf_n$, and $h=\limsup_nf_n$.

\item If $f_n$ are measurable, so are $k=\inf_nf_n$, and $l=\liminf_nf_n$.

\item If $f_n\to f$ and $f_n$ are measurable, then $f$ is measurable.

\item If $f,g$ are measurable, so are $h=\min(f,g)$ and $k=\max(f,g)$.

\item If $f$ is measurable, so are $f^+=\max(f,0)$ and $f^-=-\min(f,0)$.
\end{enumerate}
\end{proposition}

\begin{proof}
This is again something very standard, the idea being as follows:

\medskip

(1) For the function $g=\sup_nf_n$ we can use the measurability criterion from Proposition 5.13, along with the following elementary fact, valid for any $a\in\mathbb R$:
$$g=\sup_nf_n\implies g^{-1}(a,\infty]=\bigcup_nf_n^{-1}(a,\infty]$$

By symmetry we obtain that the function $k=\inf_nf_n$ is measurable as well. But with these results in hand, the last assertion follows too, by using the following formula:
$$\limsup_nf_n=\inf_k\left(\sup_{n\geq k}f_n\right)$$

(2) Here the fact that $k=\inf_nf_n$ is measurable was already proved in the above, and for $l=\liminf_nf_n$ we can use the same argument, symmetry, or the following formula:
$$\liminf_nf_n=\sup_k\left(\inf_{n\geq k}f_n\right)$$

(3) This follows indeed from (1), or from (2), and from the following fact:
$$f_n\to f\implies f=\inf_nf_n=\sup_nf_n$$

(4) This is a trivial application of (1) and (2).

\medskip

(5) This follows from (4), and from the fact that if $f$ is measurable, so is $-f$.
\end{proof}

As a main result now regarding the measurable functions, with the convention that a step function means a function which takes finitely many values, we have:

\begin{theorem}
Given a measurable function $f:X\to[0,\infty]$, we can write
$$f(x)=\lim_{n\to\infty}\varphi_n(x)$$
with $0\leq\varphi_1\leq\varphi_2\leq\ldots\leq f$, and with each $\varphi_i$ being a measurable step function.
\end{theorem}

\begin{proof}
As a first observation, the converse holds too, thanks to the results in Proposition 5.14. Regarding now the proof, this goes as follows:

\medskip

(1) First, it is clear by drawing a picture that we can approximate the identity of $[0,\infty]$ with step functions as in the statement. That is, we can obviously find an increasing sequence of measurable step functions $0\leq\psi_1\leq\psi_2\leq\ldots\leq id$, satisfying:
$$\lim_{n\to\infty}\psi_n(x)=x$$ 

(2) Now let us set $\varphi_n=\psi_n\circ f$. According to the limiting formula above, we have:
$$\lim_{n\to\infty}\varphi_n(x)=f(x)$$

On the other hand, by using the measurability criterion from Proposition 5.13, it follows that our truncation functions $\varphi_n=\psi_n\circ f$ are measurable, as desired.
\end{proof}

\section*{5b. Lebesgue, Fatou}

Good news, with the above general theory understood, we can now integrate functions, by following the good old method of Riemann, that we know well from one-variable calculus. To be more precise, we have the following result, which is of course stated a bit informally, with some of the details being left to you, as a routine exercise:

\begin{theorem}
We can integrate the measurable functions $f:X\to\mathbb R_+$ by setting
$$\int_Xf(x)d\mu(x)=\sup_{0\leq\varphi\leq f}\int_X\varphi(x)d\mu(x)$$
with sup over measurable step functions, then extend this by linearity.
\end{theorem}

\begin{proof}
We can certainly integrate the step functions $\varphi:X\to\mathbb R_+$, by writing each such function as a linear combination of characteristic functions, as follows:
$$\varphi=\sum_i\lambda_i\chi_{E_i}$$

Indeed, with this formula in hand, we can integrate our function $\varphi$, as follows:
$$\int_X\varphi(x)\,d\mu(x)=\sum_i\lambda_i\mu(E_i)$$

Next, consider an arbitrary measurable function $f:X\to[0,\infty]$. We know from Theorem 5.15 that we can write this function as an increasing limit, as follows, with $0\leq\varphi_1\leq\varphi_2\leq\ldots\leq f$, and with each $\varphi_i$ being a measurable step function:
$$f(x)=\lim_{n\to\infty}\varphi_n(x)$$

But this suggests to define the integral of $f$ by the formula in the statement, namely:
$$\int_Xf(x)d\mu(x)=\sup_{0\leq\varphi\leq f}\int_X\varphi(x)d\mu(x)$$

Indeed, we can see that the integral constructed in this way has all the linearity and positivity properties that you might expect, and behaves well with respect to limits.
\end{proof}

More in detail now, here are some basic properties of the integrals:

\begin{proposition}
The integrals of measurable functions $f:X\to\mathbb R_+$ have the following properties:
\begin{enumerate}
\item $f\leq g$ implies $\int f\leq\int g$.

\item $E\subset F$ implies $\int_Ef\leq\int_Ff$.

\item $\int f+g=\int f+\int g$.

\item $\int\lambda f=\lambda\int f$.
\end{enumerate}
\end{proposition}

\begin{proof}
All this is indeed very standard, all routine verifications.
\end{proof}

Summarizing, we know how to integrate the real positive functions, in our abstract measure theory setting. In general, we can use the following formula:
$$\int_X(f-g)(x)\,d\mu(x)=\int_Xf(x)d\mu(x)-\int_Xg(x)d\mu(x)$$

We can integrate as well the complex functions, by setting:
$$\int_X(f+ig)(x)\,d\mu(x)=\int_Xf(x)d\mu(x)+i\int_Xg(x)d\mu(x)$$

All this is, indeed, very standard, exactly as for the Riemann integral. Let us record these findings as an upgrade of Theorem 5.16, as follows, once again with the statement being a bit informal, and with some of the details being left as a routine exercise:

\begin{theorem}
We can integrate the measurable functions $f:X\to\mathbb C$ by setting
$$\int_Xf(x)d\mu(x)=\sup_{0\leq\varphi\leq f}\int_X\varphi(x)d\mu(x)$$
with sup over measurable step functions, for $f:X\to\mathbb R_+$ then extend this by linearity.
\end{theorem}

\begin{proof}
This follows indeed from what we already know from Theorem 5.16, and from the above discussion, exactly as for the Riemann integral.
\end{proof}

Many other things can be said here, following Lebesgue, Fatou and others. We first have the following result, regarding the monotone convergence, due to Lebesgue:

\begin{theorem}[Lebesgue]
Given an increasing sequence of measurable functions 
$$0\leq f_1\leq f_2\leq\ldots\leq\infty$$
which converges pointwise, $f_n\to f$, their limit is measurable, and we have
$$\int_Xf_n(x)d\mu(x)\to\int_Xf(x)d\mu(x)$$
for any positive measure on $X$.
\end{theorem}

\begin{proof}
This is indeed something very standard, the idea being as follows:

\medskip

(1) We first have the following obvious implication, showing that the sequence of integrals on the right converges, to a certain number in $[0,\infty]$: 
$$f_1\leq f_2\leq\ldots\implies\int_Xf_1(x)d\mu(x)\leq\int_Xf_2(x)d\mu(x)\leq\ldots$$

Moreover, since we have $f_n\to f$, it follows that we have the following inequality:
$$\lim_{n\to\infty}\int f_n(x)d\mu(x)\leq\int_Xf(x)d\mu(x)$$

(2) In order to prove now the reverse inequality, pick a measurable simple function $0\leq\varphi\leq f$, pick also a number $c\in(0,1)$, and consider the following sets:
$$E_n=\left\{x\in X\Big|f_n(x)\geq c\varphi(x)\right\}$$

These sets are then measurable, we have $E_1\subset E_2\subset\ldots$, and our claim is that:
$$X=\bigcup_nE_n$$

Indeed, given $x\in X$, if $f(x)=0$ then $x\in E_1$ and things fine. Otherwise, we have $f(x)>0$, and so $f(x)>c\varphi(x)$, since $c<1$, and so $x\in E_n$ for some $n$, as desired.

\medskip

(3) Now observe that, with $0\leq\varphi\leq f$ and $c\in(0,1)$ as above, we have:
$$\int_Xf_n(x)d\mu(x)\geq\int_{E_n}f_n(x)d\mu(x)\geq c\int_{E_n}\varphi(x)d\mu(x)$$

By taking the limit of this estimate, with $n\to\infty$, we obtain:
$$\lim_{n\to\infty}\int f_n(x)d\mu(x)\geq c\int_X\varphi(x)d\mu(x)$$

Now with $c\to1$, we obtain from this the following estimate:
$$\lim_{n\to\infty}\int f_n(x)d\mu(x)\geq\int_X\varphi(x)d\mu(x)$$

But this being true for any simple function $0\leq\varphi\leq f$, we conclude that we have:
$$\lim_{n\to\infty}\int f_n(x)d\mu(x)\geq\int_X f(x)d\mu(x)$$

Thus we have the reverse of the estimate found in (1), which finishes the proof.
\end{proof}

Regarding now the series of functions, again following Lebesgue, we have:

\begin{theorem}
Given measurable functions $f_n:X\to[0,\infty]$, their sum $f(x)=\sum_{n=1}^\infty f_n(x)$
is measurable, and we have the formula
$$\int_Xf(x)d\mu(x)=\sum_{n=1}^\infty\int_Xf_n(x)d\mu(x)$$
for any positive measure on $X$.
\end{theorem}

\begin{proof}
This follows indeed from Theorem 5.19, applied to the partial sums.
\end{proof}

Following now Fatou, we have as well the following key result:

\begin{theorem}[Fatou]
Given measurable functions $f_n:X\to[0,\infty]$, their limit $f(x)=\liminf_nf_n(x)$ is measurable, and we have the formula
$$\int_Xf(x)d\mu(x)\leq\liminf_n\int_Xf_n(x)d\mu(x)$$
for any positive measure on $X$.
\end{theorem}

\begin{proof}
The first assertion is something that we know, from Proposition 5.14. As for the second assertion, this can be proved by using the same trick as there, namely: 
$$\liminf_nf_n=\sup_k\left(\inf_{n\geq k}f_n\right)$$

Indeed, we can apply Theorem 5.19 to the following sequence of functions:
$$g_k(x)=\inf_{n\geq k}f_n(x)$$

Thus, we are led to the inequality in the statement. Let us also mention that this inequality can be strict, a basic example here being $f_{2k}=\chi_E$, and $f_{2k+1}=1-\chi_E$.
\end{proof}

Finally, as a last basic result on the subject, due to Lebesgue again, we have:

\begin{theorem}[Lebesgue]
Assuming $f_n\to f$ pointwise, and assuming that we have $|f_n|\leq g$, for a certain function $g$ satisfying $\int g<\infty$, the following happen:
\begin{enumerate}
\item $\int|f|<\infty$.

\item $\int|f_n-f|\to0$.

\item $\int f_n\to\int f$.
\end{enumerate}
\end{theorem}

\begin{proof}
This is something very standard, using Fatou, the idea being as follows:

\medskip

(1) This follows from $|f_n|\leq g$, which in the limit gives $|f|<g$, so $\int|f|<\infty$.

\medskip

(2) Since we have $|f_n-f|\leq 2g$, we can apply Theorem 5.21 to the sequence of functions $h_n=2g-|f_n-f|$. We obtain in this way the following estimate:
\begin{eqnarray*}
\int_X2g
&\leq&\liminf_n\int_X2g-|f_n-f|\\
&=&\int_X2g+\liminf_n\left(-\int_X|f_n-f|\right)\\
&=&\int_X2g-\limsup_n\int_X|f_n-f|
\end{eqnarray*}

Now by substracting $\int_X2g$, this estimate gives the following formula:
$$\limsup_n\int_X|f_n-f|\leq0$$

We conclude that this limit must be 0, as claimed in (2).

\medskip

(3) This assertion follows indeed from (2).
\end{proof}

\section*{5c. Main theorems}

Getting back now to our Plan 5.10, done with (1) there, and we are left with (2,3,4), namely measure and integration on $X=\mathbb R$ and variations, which in turn corresponds to the third rabbit in our original hunting plan, made after Question 5.2. As for the first two rabbits, namely probability spaces, and variables and laws, we can see after.

\bigskip

And tough hunt this will be. The idea, however, is very simple, as follows:

\begin{principle}
In order to construct our measures, that is, measure the arbitrary measurable sets $E\subset X$, the idea will be to squeeze these sets $E$ as 
$$K\subset E\subset U$$
between a compact set $K\subset X$, and an open set $U\subset X$.
\end{principle}

So, let us go for this. Regarding the compact sets, needed in the above, you have certainly heard about them, in the case $X=\mathbb R$. In general, their theory is as follows:

\begin{theorem}
Given a metric space $X$, let us call $K\subset X$ compact when any open cover $K\subset\cup_iU_i$ has a finite subcover. Then, the following happen:
\begin{enumerate}
\item Compact implies closed and bounded.

\item The converse of this fact is not true.

\item Closed inside compact is compact.

\item If $K_1\supset K_2\supset K_3\supset\ldots$ are non-empty compacts, then $\cap_iK_i\neq\emptyset$.

\item If $K$ is compact, any sequence $\{x_n\}\subset K$ has a limit point in $K$.
\end{enumerate}
Also, coming as good news, in the case $X=\mathbb R^N$, compact means closed and bounded.
\end{theorem}

\begin{proof}
This is something very standard, the idea being as follows:

\medskip

(1-2) The first assertion is easy, good exercise for you, and for the second assertion, the standard counterexample is the standard basis $K=\{e_i\}$ of the space $X=l^2(\mathbb N)$, consisting of the square-summable sequences, with metric $d(x,y)=\sqrt{\sum_i(x_i-y_i)^2}$.

\medskip

(3-5) The first assertion, which is technically very useful, is easy, good exercise for you, and the next two assertions follow from it, again exercise for you.

\medskip

(6) We must prove here that $K\subset\mathbb R^N$ closed and bounded is compact, and in view of (3), it is enough to do this for the unit cube. But here, the best is to proceed by contradiction. Indeed, assuming that $K$ has a cover $\cup_iU_i$ having no finite subcover, the same would be true for one of the subcubes $K_1$ obtained by cutting $K$ into $2^N$ smaller cubes, in the obvious way. But then, the same would be true for one of the subcubes $K_2$ obtained by cutting $K_1$ into $2^N$ smaller cubes. And so on, and by setting now $\cap_iK_i=\{x\}$, and picking an open set $U_i$ containing $x$, we have $K_r\subset U_i$ for $r>>0$ , contradiction.
\end{proof}

Getting back now to Principle 5.23, with the convention that the support of a function $f$ is the closure of the set $\{x|f(x)\neq0\}$, the statement that we will need is as follows:

\begin{lemma}[Urysohn]
Given $K\subset U$, compact inside open, we can find
$$\chi_K\leq f\leq\chi_U$$
with $f$ being continuous, and compactly supported.
\end{lemma}

\begin{proof}
Beware of lemmas, as we say in mathematics, the idea being as follows:

\medskip

(1) Our first claim is that given $K\subset U$ as in the statement, compact inside open, we can find inclusions as follows, with the set $V$ being open, with compact closure $\bar{V}$:
$$K\subset V\subset\bar{V}\subset U$$

Indeed, for any $x\in K$ pick a neighborhood $W_x$ having compact closure. Since $K$ is compact, we can find finitely many such neighborhoods, covering it:
$$K\subset W_{x_1}\cup\ldots\cup W_{x_n}$$

In order to prove our claim, observe that if we are in the case $U=X$ we are done, because here we can simply take $V=W$, the set on the right. So, assume $U\neq X$. Given a point $x\in U^c$, since this point is compact, and our set $K$ not containing it is compact as well, we can find an open set $V_x$ having the following properties:
$$K\subset V_x\quad,\quad x\notin\bar{V}_x$$

Now consider the following family of sets, with $W$ being the open set above:
$$\left\{K\cap\bar{W}\cap\bar{V}_x\Big|x\in U^c\right\}$$

This is then a family of compact sets having empty intersection, so we can find a finite subfamily having empty intersection too. That is, we can find points $x_i\in U^c$ such that:
$$K\cap\bar{W}\cap\bar{V}_{x_1}\cap\ldots\cap\bar{V}_{x_m}=\emptyset$$

But with this done, $V=W\cap V_{x_1}\cap\ldots\cap V_{x_m}$ is the open set we are looking for.

\medskip

(2) Getting now to what our statement says, given $K\subset U$ as there, compact inside open, by (1) we can find an open set $V_0$, having compact closure $\bar{V}_0$, such that:
$$K\subset V_0\subset\bar{V}_0\subset U$$

But then, by using again (1), applied this time to the inclusion $K\subset V_0$, we can find a second open set $V_1$, having compact closure $\bar{V}_1$, such that:
$$K\subset V_1\subset\bar{V}_1\subset V_0\subset\bar{V}_0\subset U$$

And so on, would be the idea. In practice now, by using the fact that the rational numbers are countable, we can construct in this way a whole family of open sets $V_r$, having compact closures $\bar{V}_r$, one for each rational number $r\in[0,1]$, such that:
$$r<s\implies \bar{V}_s\subset V_r$$

Time now to construct our function $f$. Let us set, for any $r\in[0,1]$ rational:
$$f_r(x)=\begin{cases}
r&{\rm if}\ x\in V_r\\
0&{\rm otherwise}
\end{cases}$$

We can then set $f=\sup_rf_r$, which does indeed the job, say good exercice for you.
\end{proof}

Good news, eventually, the above is all that we need, as abstract preliminaries. We can now formulate the key result in abstract measure theory, as follows:

\begin{theorem}[Riesz]
Any positive functional $I:C_c(X)\to\mathbb R$ comes by integrating with respect to a certain measure on $X$,
$$I(f)=\int_Xf(x)d\mu(x)$$
which is unique, modulo null sets. In addition, the following happen:
\begin{enumerate}
\item $\mu(K)<\infty$, for any $K\subset X$ compact.

\item $\mu(E)=\inf\{\mu(U)|E\subset U\ {\rm open}\}$, for any $E$ measurable.

\item $\mu(E)=\sup\{\mu(K)|K\subset E\ {\rm compact}\}$, for any $E$ open, or of finite measure.
\end{enumerate}
\end{theorem}

\begin{proof}
This is something quite tricky, the idea being as follows:

\medskip

(1) Let us first prove the uniqueness. Assuming that our measure $\mu$ produces $I$, and has the properties (1,2,3) in the statement, it is clear from (2,3) that $\mu$ is uniquely determined by its values on the compact sets $K\subset X$. Thus, we must show that given two measures $\mu_1,\mu_2$ as in the statement, and a compact set $K\subset X$, we have:
$$\mu_1(K)=\mu_2(K)$$

(2) For this purpose, let us pick $\varepsilon>0$. By using the properties (1,3) in the statement, for the measure $\mu_2$, we can find $K\subset U$ open such that:
$$\mu_2(U)<\mu_2(K)+\varepsilon$$

Now by using the Urysohn lemma for the inclusion $K\subset U$, we obtain a certain continuous, compactly supported function $f$, such that:
$$\chi_K\leq f\leq\chi_U$$

(3) But, with this choice of $f$, we have the following computation:
$$\mu_1(K)
=\int_X\chi_K(x)d\mu_1(x)
\leq\int_Xf(x)d\mu_1(x)
=I(f)$$

On the other hand, we have as well the following computation:
$$I(f)
=\int_Xf(x)d\mu_2(x)
\leq\int_X\chi_U(x)d\mu_2(x)
=\mu_2(U)
<\mu_2(K)+\varepsilon$$

Now all this being true for any $\varepsilon>0$, we conclude that we have:
$$\mu_1(K)\leq\mu_2(K)$$

But by interchanging $\mu_1,\mu_2$, in this reasoning, we have the reverse inequality as well. Thus we obtain, as desired, proving the uniqueness statement:
$$\mu_1(K)=\mu_2(K)$$

(4) With this done, let us get now to the main part, and discuss the construction of the algebra $M$, that we will choose to be saturated with respect to null sets, and of the measure $\mu$. What we have as data is a positive functional, as follows:
$$I:C_c(X)\to\mathbb R$$

We must first measure the Borel sets $E\subset X$. And here, to start with, in order to measure the open sets $U\subset X$, the formula is straightforward, namely:
$$\mu(U)=\sup\left\{I(f)\Big| f\leq\chi_U\right\}$$

Now observe that with this definition in hand, for the open sets, we have:
$$U_1\subset U_2\implies\mu(U_1)\leq\mu(U_2)$$

We deduce that the following happens, for the open sets:
$$\mu(E)=\inf\left\{\mu(U)\Big| E\subset U\ {\rm open}\right\}$$

But this can serve as the definition of $\mu$, for all the subsets $E\subset X$.

\medskip

(5) Regarding now the measurable sets, let us first consider the class $N$ of subsets $E\subset X$ which are of finite measure, $\mu(E)<\infty$, and satisfy the following condition:
$$\mu(E)=\sup\left\{\mu(K)\Big| K\subset E\ {\rm compact}\right\}$$

Then, we can define the class of measurable sets $M\subset P(X)$ as follows:
$$M=\left\{E\subset X\Big|E\cap K\in N,\forall K\subset X\ {\rm compact}\right\}$$

(6) Summarizing, done with definitions, and it remains now to prove that $M$ is an algebra, and that $\mu$ is a measure on it, along with the other things claimed in the statement. But this can be done indeed, by making a heavy use of the Urysohn lemma, and for the details, taking a few pages, you can find them in the book of Rudin \cite{rud}.
\end{proof}

As a consequence of the Riesz theorem, we can now formulate, following Lebesgue:

\begin{theorem}[Lebesgue]
We can measure all the Borel sets 
$$E\subset\mathbb R$$
and the corresponding integration theory extends the usual Riemann integral.
\end{theorem}

\begin{proof}
This follows indeed from Theorem 5.26, by using the Riemann integral, that we know well from basic calculus, as the positive functional needed there:
$$f\to\int_\mathbb Rf(x)dx$$

Indeed, it follows from basic calculus that the Riemann integral satisfies all the needed properties in Theorem 5.26, and this gives the result. For details, see Rudin \cite{rud}.
\end{proof}

Goind ahead with more theory, let us discuss now the case of $\mathbb R^N$. And here, we are a bit in trouble, because do we really know how to integrate on $\mathbb R^N$, as shown by:

\begin{proposition}
The Fubini formula, namely 
$$\int_\mathbb R\int_\mathbb Rf(x,y)dxdy=\int_\mathbb R\int_\mathbb Rf(x,y)dydx$$
can fail, for certain suitably chosen functions.
\end{proposition}

\begin{proof}
Listen carefully. We have the following computation:
$$\int_0^1\int_0^1\frac{y^2-x^2}{(x^2+y^2)^2}\,dxdy
=\int_0^1\left[\frac{x}{x^2+y^2}\right]_0^1dy
=\int_0^1\frac{1}{1+y^2}\,dy
=\frac{\pi}{4}$$

On the other hand, by using this, and symmetry, we have as well:
$$\int_0^1\int_0^1\frac{y^2-x^2}{(x^2+y^2)^2}\,dydx
=-\int_0^1\int_0^1\frac{y^2-x^2}{(x^2+y^2)^2}\,dxdy
=-\frac{\pi}{4}$$

Thus Fubini fails. Please pass this knowledge to your children, and grandchildren.
\end{proof}

Summarizing, our knowledge of multivariable integration from calculus is flawed, I mean we don't have a reliable definition for the Riemann integral there. Of course, this is not a problem for most questions appearing in the real life, but in the present setting, it is quite obvious that we cannot cheat, as we usually do. So, here is the fix:

\begin{proposition}
Given a compactly supported function $f:\mathbb R^N\to\mathbb C$, if we set
$$I_k(f)=\frac{1}{2^{Nk}}\sum_{x\in\mathbb Z^N/2^k}f(x)$$
with $\mathbb Z^N/2^k\subset\mathbb R^N$ being the points having as coordinates integer multiples of $1/2^k$, then
$$I(f)=\lim_{k\to\infty}I_k(f)$$
converges, and $f\to I(f)$ satisfies the assumptions of the Riesz theorem.
\end{proposition}

\begin{proof}
In order to prove the first assertion, we can assume that our function is real, $f:\mathbb R^N\to\mathbb R$. But then, we can use the fact that $f$ is uniformly continuous, and this gives the convergence in the statement. As for the second assertion, the fact that $f\to I(f)$ satisfies indeed the assumptions of the Riesz theorem is clear from definitions.
\end{proof}

In view of the above, we can indeed apply the Riesz theorem, and we obtain:

\begin{theorem}[Lebesgue]
There is a unique regular, translation-invariant measure on $\mathbb R^N$, normalized as for the volume of the unit cube to be $1$. This measure is called Lebesgue measure, and appears alternatively as the measure on the Borel sets 
$$E\subset\mathbb R$$
whose corresponding integration theory extends the usual Riemann integral, as constructed in Proposition 5.29, or equivalently, as in multivariable calculus.
\end{theorem}

\begin{proof}
This is something quite self-explanatory, and we leave some reading here, details on the notions in the statement, and proof, from Rudin \cite{rud}, as an exercise.
\end{proof}

Summarizing, we know how to rigorously integrate on $\mathbb R^N$. As a continuation of this, as a most pressing task, let us go back to the Fubini formula problematics. We have:

\begin{theorem}[Fubini, Tonelli]
Given a function $f:\mathbb R^N\to\mathbb R$ which is measurable and integrable, in the sense that the following integral is finite,
$$\int_{\mathbb R^N}|f(z)|dz<\infty$$
we have the following equalities, for any decomposition $N=N_1+N_2$:
$$\int_{\mathbb R^{N_1}}\int_{\mathbb R^{N_2}}f(x,y)dydx=\int_{\mathbb R^{N_2}}\int_{\mathbb R^{N_1}}f(x,y)dxdy=\int_{\mathbb R^N}f(z)dz$$
Moreover, the same holds when $f:\mathbb R^N\to\mathbb R$ is assumed positive, and measurable. 
\end{theorem}

\begin{proof}
In order to prove the first assertion, due to Fubini, consider the vector space $L$ consisting of the functions $f:\mathbb R^N\to\mathbb R$ which are measurable, and integrable:
$$L=\left\{f:\mathbb R^N\to\mathbb R\Big|\int_{\mathbb R^N}|f(z)|dz<\infty\right\}$$

Now fix a decomposition $N=N_1+N_2$, and consider the family $F\subset L$ of functions $f:\mathbb R^N\to\mathbb R$ which satisfy the formula in the statement, namely:
$$\int_{\mathbb R^{N_1}}\int_{\mathbb R^{N_2}}f(x,y)dydx=\int_{\mathbb R^{N_2}}\int_{\mathbb R^{N_1}}f(x,y)dxdy=\int_{\mathbb R^N}f(z)dz$$

With these conventions, we want to prove that the inclusion $F\subset L$ is an equality. But this comes from the following sequence of observations, which all come from definitions, and from the Lebesgue dominated convergence theorem, that is, Theorem 5.22:

\medskip

-- $F$ is stable under taking finite linear combinations.

\medskip

-- $F$ is stable under taking pointwise limits of increasing sequences.

\medskip

-- $1_E\in F$, for any set of finite measure of type $E=E_1\times E_2\subset\mathbb R^N$.

\medskip

-- $1_E\in F$, for any open set of finite measure $E\subset\mathbb R^N$.

\medskip

-- $1_E\in F$, for any countable intersection of open sets, of finite measure $E\subset\mathbb R^N$.

\medskip

-- $1_E\in F$, for any set of null measure $E\subset\mathbb R^N$.

\medskip

-- $1_E\in F$, for any set of finite measure $E\subset\mathbb R^N$.

\medskip

Indeed, once we have the last assertion, the result follows, because we can approximate any given function $f\in L$ by step functions, and we obtain $f\in F$. Thus, Fubini proved, and the second assertion in the statement, Tonelli, follows from this. See Rudin \cite{rud}.
\end{proof}

\section*{5d. Random variables}

Good news, time to get back to our regular business in this book, probability theory. Based on what we have, we can quickly axiomatize everything, as folows:

\begin{definition}
The axioms for probability theory are as follows:
\begin{enumerate}
\item A probability space is a measured space $(X,M,\nu)$ of mass one, $\nu(X)=1$.

\item A random variable on $X$ is a real measurable function $f:X\to\mathbb R$.

\item The expectation of $f:X\to\mathbb R$ is its integral, $E(f)=\int_Xf(x)d\nu(x)$.

\item The law of $f$ is the real probability measure $\mu=f_*\nu$, push-forward of $\nu$ by $f$.
\end{enumerate}
\end{definition}

In short, what we are doing here is to extend our previous discrete probability theory, axiomatized as in Fact 5.1, by using our measure theory knowledge. Which looks quite straightforward, but in practice, there are a few subtleties with this, as follows:

\bigskip

(1) To start with, when the measure $\nu$ is discrete, meaning that $M=P(X)$, and that $\nu(Y)=\sum_{y\in Y}\varphi(y)$, for some function $\varphi:X\to[0,1)$, we recover discrete probability theory, as axiomatized in Fact 5.1. Indeed, this comes from our results in chapter 2.

\bigskip

(2) Next, as a new phenomenon appearing, in our generalized setting, the expectation $E(f)=\int_Xf(x)d\nu(x)$ can now take infinite values, $\pm\infty$, or even be undefined, due to the lack of convergence of the integral. So, in a word, always careful with this.

\bigskip

(3) Finally, in what regards the formula $\mu=f_*\nu$, we are using here the standard fact that given a measured space $(X,M,\nu)$ and a measurable function $f:X\to\mathbb R$, the formula $\mu(Z)=\nu(f^{-1}(Z))$ defines a measure $\mu=f_*\nu$ on $\mathbb R$, called push-forward of $\nu$ by $f$.

\bigskip

In practice now, Definition 5.32 remains something a bit abstract, and in what regards the random variables and their laws, which is what matters the most, we have:

\begin{theorem}
Given a random variable $f:X\to\mathbb R$ having law $\mu$, we have the following formula, for any measurable function $\phi:\mathbb R\to\mathbb R$,
$$E(\phi(f))=\int_\mathbb R\phi(x)d\mu(x)$$
and this can stand, via the Riesz theorem, as an alternative definition for $\mu$.
\end{theorem}

\begin{proof}
This is something quite self-explanatory, the idea being as follows:

\medskip

(1) To start with, the formula in the statement for the expectation $E(\phi(f))$ comes indeed from our definition of the law as a push-forward measure, $\mu=f_*\nu$.

\medskip

(2) As for the last assertion, this comes from the Riesz theorem, the point being that a probability measure $\mu$ is uniquely determined by the integration over it. 

\medskip

(3) However, there are more things that can be said here, for instance by using the monomials $\phi(x)=x^k$, producing the moments $M_k=E(f^k)$. More on this later.
\end{proof}

As in the discrete probability setting, the law of a variable $f:X\to\mathbb R$ is not the only main quantity related to $f$ that we can talk about. As a rival to it, we have:

\begin{definition}
Given a random variable $f:X\to\mathbb R$, the function
$$F(x)=P(f\leq x)$$
is called its cumulative distribution function (CDF).
\end{definition}

Observe that the cumulative distribution function $F:\mathbb R\to[0,1]$ must be by definition increasing, and also, that it must have the following limiting properties:
$$F(-\infty)=0\quad,\quad F(\infty)=1$$

With a bit more care, by looking into continuity properties as well, we are led to the following result, regarding the CDF of the general random variables $f:X\to\mathbb R$:

\begin{theorem}
The cumulative distribution function $F:\mathbb R\to[0,1]$ must be
\begin{enumerate}
\item Increasing,

\item Right continuous,

\item Such that $F(-\infty)=0$,

\item Such that $F(\infty)=1$,
\end{enumerate}
and conversely, any function satisfying these properties appears as a CDF.
\end{theorem}

\begin{proof}
Here the first assertion is quite clear, and I will leave it to you to check the right continuity property, after remembering what that property exactly is. As for the second assertion, this is something which is clear too, good exercise for you.
\end{proof}

Moving on, as something which is new, of genuine continuous nature, we have:

\begin{theorem}
Assume that a random variable $f:X\to\mathbb R$ has density $\varphi:\mathbb R\to\mathbb R$, in the sense that its law is given by $d\mu(x)=\varphi(x)dx$.
\begin{enumerate}
\item The density must be positive, $\varphi\geq0$, and of mass one, $\int_\mathbb R\varphi(x)dx=1$.

\item $\mu(Z)=\int_Z\varphi(x)dx$, which can stand as an alternative definition for $\varphi$.

\item $E(\phi(f))=\int_\mathbb R\phi(x)\varphi(x)dx$, which can stand as well as a definition for $\varphi$.

\item The CDF is $F(x)=\int_{-\infty}^x\varphi(t)dt$, and conversely, we have $\varphi(x)=F'(x)$.
\end{enumerate}
\end{theorem}

\begin{proof}
As before with other such things, this is quite self-explanatory, as follows:

\medskip

(1) This is something clear, coming from $d\mu(x)=\varphi(x)dx$.

\medskip

(2) This is again clear, also coming from $d\mu(x)=\varphi(x)dx$.

\medskip

(3) This comes indeed via the general theory from Theorem 5.33.

\medskip

(4) This comes from the fundamental theorem of calculus.
\end{proof}

Quite nice all this, we are learning new things here, with respect to what we previously knew from Part I. However, no idea if you noticed, but there is an elephant in the room, coming from our lack of knowledge of the probability measures on $\mathbb R$. I mean, personally I would love to have a theorem stating that such measures appear as follows:
$$\mu=\varphi(x)dx+\sum_ic_i\delta_{x_i}$$

And, is this possible. Not clear, the few computations that I have leading nowhere, and I will have I guess to disturb again the cat. Which cat, after a big yawn, declares:

\begin{cat}
Your axioms seem to be wrong, indeed, but I forgot the good ones, we felines don't bother with such things. Ask rat, the new one is quite smart.
\end{cat}

Damn, so we have a rat at the house, you should have told me earlier, sweetie. Anyway, thanks for the answer, and I'll think some more, about both axioms and pests.

\bigskip

So, getting back now to Rudin \cite{rud}, here is something interesting, from there:

\begin{theorem}
Any probability measure on $\mathbb R^N$ is of the following form,
$$d\mu(x)=\varphi(x)dx+d\eta(x)$$
with $\varphi:\mathbb R^N\to\mathbb R_+$ being a certain density, and $\eta$ being supported on a measure $0$ set.
\end{theorem}

\begin{proof}
Long story with this, the idea being that this is a particular case of a general theorem of Radon-Nikodym, dealing with the comparison of pairs of measures. For more on all this, which is definitely something worth learning, have a look at Rudin \cite{rud}.
\end{proof}

Very nice all this, there is still hope, but coming now as bad news, we have:

\begin{theorem}
Exotic probability measures on $\mathbb R$ exist, with a standard example being the Cantor distribution, supported by the Cantor set $K\subset\mathbb R$.
\end{theorem}

\begin{proof}
This is Devil's contribution to mathematics. Consider the Cantor set $K\subset\mathbb R$, appearing by definition as the intersection of the following unions of intervals:
$$\xymatrix@R=10pt@C=10pt{
K_0:&0\ar@{-}[rrrrrrrrr]&&&&&&&&&1\\
K_1:&0\ar@{-}[rrr]&&&\frac{1}{3}&&&\frac{2}{3}\ar@{-}[rrr]&&&1\\
K_2:&0\ar@{-}[r]&\frac{1}{9}&\frac{2}{9}\ar@{-}[r]&\frac{1}{3}&&&\frac{2}{3}\ar@{-}[r]&\frac{7}{9}&\frac{8}{9}\ar@{-}[r]&1\\
&&&&&\vdots&\vdots}$$

As a first observation, the Lebesgue measure of this Cantor set $K$ is given by:
\begin{eqnarray*}
\lambda(K)
&=&1-\frac{1}{3}-\frac{2}{9}-\frac{4}{27}-\ldots\\
&=&1-\frac{1}{3}\left(1+\frac{2}{3}+\frac{4}{9}+\ldots\right)\\
&=&1-\frac{1}{3}\cdot\frac{1}{1-2/3}\\
&=&0
\end{eqnarray*}

The point now is that we can construct a natural probability measure on $K$, which is however exotic, meaning that it is not a linear combination of Dirac masses, by stating that $\mu$ must be uniform with respect to the 2 intervals of $K_1$, then with respect to the 4 intervals of $K_2$, and so on. And, exercise for you, to learn more about all this.
\end{proof}

Well, not fun all this, and personally I feel a bit lost, at this point. Are our probability axioms wrong, as suggested by Cat 5.37, and confirmed by Theorem 5.39? What are then the correct axioms? And, in general, how to get out of the mud that we are in?

\bigskip

Guess I will have to ask the rat. While working on Theorems 5.38 and 5.39, in my kitchen as usual, I got to know him a bit, guy is really smart, I doubt that cat will ever catch him, and with me of course not being able to catch the cat. So, let's see what the smartest animal in the house has to say, about the difficulties of measure theory:

\begin{rat}
Mathematics is a natural science, go with the axioms most adapted to the phenomena that you are interested in. In harmony with your inner self.
\end{rat}

Well, looks like this rat fellow is far more advanced than us, both me and cat, at mathematics, life and everything. Must come from food I guess, adding a few vegetables, and some junk too, to your diet, most likely makes you smarter. I'll tell cat about this.

\bigskip

In any case, going now for his advice, $100\%$ of the probability measures that I know, from my daily quantum physics work, are measures having densities, as in Theorem 5.36, with some Dirac masses added. So, let us make the following convention:

\begin{convention}
We will study random variables $f:X\to\mathbb R$ having laws of type
$$\mu=\varphi(x)dx+\sum_ic_i\delta_{x_i}$$
with $\varphi:\mathbb R\to[0,\infty)$ piecewise continuous, $c_i>0$, and $\int_\mathbb R\varphi+\sum_ic_i=1$.
\end{convention}

And with this, end of our discussion regarding the axiomatics. Tough learning this was, and with a Mobutu type axiom added at the end, but we made it. Great.

\section*{5e. Exercises}

This was a quite crowded theoretical chapter, and as exercises, we have:

\begin{exercise}
Get familiar with open and closed sets, via theory and practice.
\end{exercise}

\begin{exercise}
Do the same with the compact sets, and the connected sets too.
\end{exercise}

\begin{exercise}
Read the full proof of the Riesz theorem.
\end{exercise}

\begin{exercise}
Read the proof of the Lebesgue theorem.
\end{exercise}

\begin{exercise}
Read the proof of the Fubini theorem.
\end{exercise}

\begin{exercise}
Read the proof of the Radon-Nikodym theorem.
\end{exercise}

\begin{exercise}
Fill in the details, in what we said about variables and laws.
\end{exercise}

\begin{exercise}
Further meditate on the axiomatization of probability theory.
\end{exercise}

As bonus exercise, and no surprise here, purchase a copy of Rudin \cite{rud}.

\chapter{Exponential laws}

\section*{6a. Uniform laws}

Good news, with the axiomatics of continuous probability discussed, we can now start doing some concrete work, on one hand by attempting to model the various continuous phenomena, appearing in the real life, in probabilistic terms, and on the other hand with lots of computations, as we like them, for the various laws that we find.

\bigskip

We will use of course some inspiration from discrete probability, which is a quite elementary business, that we master well, after all that learning from Part I. 

\bigskip

So, what is at the beginning? In discrete probability, we perfectly know from Part I that at the beginning we have activities like flipping coins, or rolling dice. Which in probabilistic terms means that we are dealing with uniform laws, as follows:
$$\mu=\frac{1}{n}\left(\delta_{x_1}+\ldots+\delta_{x_n}\right)$$

The point now is that we can talk about uniform laws in the continuous setting as well, with these being something quite elementary, introduced as follows:

\begin{definition}
We say that a variable $f:X\to\mathbb R$ is uniform if its law has a density, taking only 2 values. That is, the law must be of the following type,
$$\mu=\frac{1}{\lambda(E)}\,\chi_E(x)dx$$
for a certain subset $E\subset\mathbb R$, with $\lambda$ being the Lebesgue measure. In the case of an interval $E=[a,b]$, the corresponding law must be given by the following formula,
$$\mu=\frac{1}{b-a}\,\chi_{[a,b]}(x)dx$$
with the $\chi$ symbols standing in all this, as usual, for characteristic functions.
\end{definition}

Looking at what we have here, it is pretty much clear that the mathematics of the uniform laws is something very simple. However, before getting into more complicated laws, let us spend some time on these, with the aim of reviewing the material from chapter 5, which was something quite abstract, now that we have some examples available.

\bigskip

We will focus on the uniform laws on intervals $[a,b]\subset\mathbb R$, which are the most important ones, and leave the extensions involving more general subsets $E\subset\mathbb R$ as exercises. As our first result, which is something trivial, coming from definitions, we have:

\begin{proposition}
A variable $f:X\to\mathbb R$ follows the uniform law of $[a,b]$ precisely when the following happens, for any measurable function $\phi:\mathbb R\to\mathbb R$:
$$E(\phi(f))=\frac{1}{b-a}\int_a^b\phi(x)dx$$
Equivalently, we must have $f:X\to[a,b]$, up to a set of null measure, and given a measurable set $Y\subset[a,b]$, the following must happen, $\lambda$ being the Lebesgue measure:
$$P(f\in Y)=\frac{\lambda(Y)}{b-a}$$
Also, in this latter statement, we can restrict the attention to the intervals $Y\subset[a,b]$.
\end{proposition}

\begin{proof}
This is indeed something self-explanatory, the idea being as follows:

\medskip

(1) Regarding the formula of $E(\phi(f))$, and the converse assertion too, that comes from our definition of the law of a random variable, as axiomatized in chapter 5.

\medskip

(2) Regarding the formula of $P(f\in Y)$, that comes from (1) with $\phi=\chi_Y$, and the converse is clear too, by using the measure theory developed in chapter 5.

\medskip

(3) Finally, the last assertion, stating that we can restrict the attention to the intervals $Y\subset[a,b]$, is clear too, again by using the measure theory developed in chapter 5.
\end{proof}

Let us record as well a result regarding the cumulative distribution function:

\begin{proposition}
Given a variable $f:X\to\mathbb R$ which follows the uniform law of $[a,b]$, its cumulative distribution function (CDF) is as follows,
$$
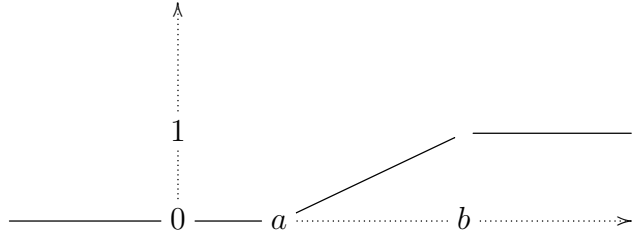
\xymatrix@R=18pt@C=25pt{
&&&\\
&&&&&\\
&&1\ar@.[uu]&&&\ar@{-}[rr]&&\\
\ar@{-}[rr]&&0\ar@{.}[u]\ar@{-}[r]&a\ar@{.}[rr]\ar@{-}[urr]&&b\ar@.[rr]&&}$$
and the converse holds too, in the sense that having this CDF means that our variable $f:X\to\mathbb R$ must follow the uniform law of $[a,b]$.
\end{proposition}

\begin{proof}
Again, this is something self-explanatory, and trivial, coming from the definition and general theory of the CDF, as developed in chapter 5.
\end{proof}

What is next? Some moment computations I guess, because, as we know well from the discrete case, from Part I, moments are where the mathematical fun happens. For our uniform laws the moment computations are quite simple, and we are led to:

\begin{theorem}
Given a variable $f:X\to\mathbb R$ which follows the uniform law of $[a,b]$, its moments are given by the following formula,
$$M_k=\frac{1}{k+1}\cdot\frac{b^{k+1}-a^{k+1}}{b-a}$$
with the small order moments being as follows:
$$M_1=\frac{a+b}{2}$$
$$M_2=\frac{a^2+ab+b^2}{3}$$
$$M_3=\frac{a^3+a^2b+ab^2+b^3}{4}$$
$$M_4=\frac{a^4+a^3b+a^2b^2+ab^3+b^4}{5}$$
In particular, the mean is $E=(a+b)/2$, and the variance is $V=(b-a)^2/12$.
\end{theorem}

\begin{proof}
As mentioned, this is something elementary, the idea being as follows:

\medskip

(1) We use the general integration formula from Proposition 6.2, namely:
$$E(\phi(f))=\frac{1}{b-a}\int_a^b\phi(x)dx$$

With $\phi(x)=x^k$ we obtain in this way the formula of the $k$-th moment, as follows:
\begin{eqnarray*}
M_k
&=&\frac{1}{b-a}\int_a^bx^kdx\\
&=&\frac{1}{b-a}\left[\frac{x^{k+1}}{k+1}\right]_a^b\\
&=&\frac{1}{b-a}\cdot\frac{b^{k+1}-a^{k+1}}{k+1}\\
&=&\frac{1}{k+1}\cdot\frac{b^{k+1}-a^{k+1}}{b-a}
\end{eqnarray*}

(2) Next, the explicit formulae for the small order moments are all clear, and clear as well is the assertion regarding the mean, which is the first moment, $E=M_1$.

\medskip

(3) Regarding now the variance, this can be computed as follows:
\begin{eqnarray*}
V
&=&M_2-M_1^2\\
&=&\frac{a^2+ab+b^2}{3}-\left(\frac{a+b}{2}\right)^2\\
&=&\frac{4(a^2+ab+b^2)-3(a+b)^2}{12}\\
&=&\frac{a^2+b^2-2ab}{12}\\
&=&\frac{(b-a)^2}{12}
\end{eqnarray*}

Thus, we are led to the conclusions in the statement.
\end{proof}

Still following the material from Part I, we can try next to compute the central moments, and then the normalized central moments, with the final aim of computing the skewness $\gamma$ and kurtosis $\kappa$, as to have all the parameters $(E,V,\gamma,\kappa)$. We first have:

\begin{theorem}
Given a variable $f:X\to\mathbb R$ which follows the uniform law of $[a,b]$, its even central moments are given by the following formula,
$$M_k'=\frac{1}{k+1}\left(\frac{b-a}{2}\right)^k$$
with the small order even central moments being as follows:
$$M_2'=\frac{(b-a)^2}{12}\quad,\quad 
M_4'=\frac{(b-a)^4}{80}\quad,\quad
M_6'=\frac{(b-a)^6}{448}\quad,\quad\ldots$$
As for the odd central moments, these all vanish.
\end{theorem}

\begin{proof}
We recall from chapter 2 that the central moments of a variable $f:X\to\mathbb R$, designed for having $M_1'=0$, $M_2'=V$, are the following numbers, with $E=E(f)$:
$$M_k'=E((f-E)^k)$$

Thus, we are in need of a remake of the computation from the proof of Theorem 6.4. As there, we will use the general integration formula from Proposition 6.2, namely:
$$E(\phi(f))=\frac{1}{b-a}\int_a^b\phi(x)dx$$

With the function $\phi(x)=(x-E)^k$, where $E=(a+b)/2$, as computed in Theorem 6.4, we obtain in this way a formula for the $k$-th central moment, as follows:
\begin{eqnarray*}
M_k'
&=&\frac{1}{b-a}\int_a^b\left(x-\frac{a+b}{2}\right)^kdx\\
&=&\frac{1}{b-a}\left[\frac{1}{k+1}\left(x-\frac{a+b}{2}\right)^{k+1}\right]_a^b\\
&=&\frac{1}{b-a}\cdot\frac{1}{k+1}\left[\left(b-\frac{a+b}{2}\right)^{k+1}-\left(a-\frac{a+b}{2}\right)^{k+1}\right]\\
&=&\frac{1}{b-a}\cdot\frac{1}{k+1}\left[\left(\frac{b-a}{2}\right)^{k+1}-\left(\frac{a-b}{2}\right)^{k+1}\right]
\end{eqnarray*}

We can see that this quantity vanishes when $k$ is odd. As for the case where $k$ is even, here we obtain the formula in the statement, as follows:
$$M_k'=\frac{1}{b-a}\cdot\frac{1}{k+1}\cdot 2\left(\frac{b-a}{2}\right)^{k+1}
=\frac{1}{k+1}\left(\frac{b-a}{2}\right)^k$$

Thus, we are led to the conclusions in the statement.
\end{proof}

Getting now to the normalized central moments, we have here:

\begin{theorem}
Given a variable $f:X\to\mathbb R$ which follows the uniform law of $[a,b]$, its even normalized central moments are given by the following formula,
$$M_k''=\frac{3^{k/2}}{k+1}$$
with the small order even normalized central moments being as follows:
$$M_2''=1\quad,\quad 
M_4''=\frac{9}{5}\quad,\quad
M_6''=\frac{27}{7}\quad,\quad\ldots$$
As for the odd normalized central moments, these all vanish.
\end{theorem}

\begin{proof}
We recall from chapter 2 that the normalized central moments of a variable $f:X\to\mathbb R$, which are designed for having $M_1''=0$, $M_2''=1$, are the following numbers, with $E=E(f)$ being the expectation, and $\sigma=\sqrt{V}$ being the standard deviation:
$$M''_k=E\left(\left(\frac{f-E}{\sigma}\right)^k\right)$$

In other words, we have $M_k''=M_k'/\sigma^k$, with $\sigma=\sqrt{V}$, which in our case is:
$$\sigma=\sqrt{\frac{(b-a)^2}{12}}=\frac{b-a}{2\sqrt{3}}$$

Now by using the formula of $M_k'$ from Theorem 6.5, for $k$ even we obtain:
$$M_k''=\frac{1}{k+1}\left(\frac{b-a}{2}\right)^k\Big/\left(\frac{b-a}{2\sqrt{3}}\right)^k
=\frac{3^{k/2}}{k+1}$$

As for the case where $k$ is odd, here, again by Theorem 6.5, we obtain $M_k''=0$. 
\end{proof}

Let us end this discussion with the following result, capturing the essentials:

\begin{theorem}
For a variable $f:X\to\mathbb R$ which follows the uniform law of $[a,b]$,
$$E=\frac{a+b}{2}\quad,\quad V=\frac{(b-a)^2}{12}\quad,\quad\gamma=0\quad,\quad\kappa=\frac{9}{5}$$
are the mean, variance, skewness and kurtosis. 
\end{theorem}

\begin{proof}
Here the formulae of the mean $E$ and variance $V$ are those from Theorem 6.4, and the formulae of the skewness $\gamma=M_3''$ and kurtosis $\kappa=M_4''$ are those from Theorem 6.6. Thus, we are led to the conclusions in the statement.
\end{proof}

Regarding now the Fourier transform computation, this is as follows:

\begin{theorem}
For a variable $f:X\to\mathbb R$ which follows the uniform law of $[a,b]$,
$$F_f(x)=\frac{e^{ibx}-e^{iax}}{i(b-a)x}$$
is its Fourier transform, $F_f(x)=E(e^{ifx})$. 
\end{theorem}

\begin{proof}
This is again a standard computation, based on the general integration formula in Proposition 6.2, by using the function $\phi(t)=e^{itx}$ there, which gives:
\begin{eqnarray*}
F_f(x)
&=&E(e^{ifx})\\
&=&\frac{1}{b-a}\int_a^be^{itx}dt\\
&=&\frac{1}{b-a}\left[\frac{e^{itx}}{ix}\right]_a^b\\
&=&\frac{1}{b-a}\left(\frac{e^{ibx}-e^{iax}}{ix}\right)\\
&=&\frac{e^{ibx}-e^{iax}}{i(b-a)x}
\end{eqnarray*}

Thus, we are led to the conclusion in the statement.
\end{proof}

And with this, end of our discussion of the uniform laws. There are of course more things that can be said, for instance in regards with summing independent uniform variables, but more on such things later, when discussing more in detail independence.

\section*{6b. Exponential laws} 

What is next? More complicated distributions I guess, and with me writing this book, passioned by particle physics, I will choose the exponential laws, which are related to particle decay, as we will soon discover. By the way, these laws are as well something very simple, mathematically speaking, and according to a recent poll among probability teachers at the university, they ranked $\#2$, right after the uniform laws.

\bigskip

But let us begin with some philosophy. Forgetting about applications, which is what probability theory is good for, if you want to construct probability measures, just for the sake of constructing probability measures, and studying them, here is the recipe:

\begin{recipe}
For constructing a probability measure with density, $d\mu(x)=\varphi(x)dx$, what you need is a function $\varphi:\mathbb R\to[0,\infty)$ satisfying $\int_\mathbb R\varphi(x)dx=1$. And here:
\begin{enumerate}
\item In the case of bounded support, you need a function $\varphi$ as follows,
$$
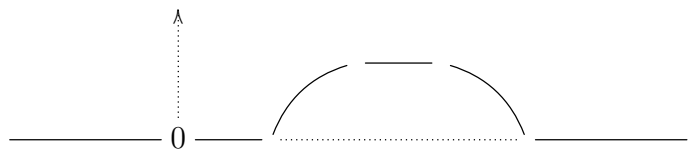
\xymatrix@R=18pt@C=25pt{
&&&&&\\
&&&&\ar@{-}[r]&\ar@/^/@{-}[dr]&&\\
\ar@{-}[rr]&&0\ar@.[uu]\ar@{-}[r]&\ar@{.}[rrr]\ar@/^/@{-}[ur]&&&\ar@{-}[rr]&&}$$
the simplest example being the rectangle, leading to the uniform laws.

\item In the case of positive support, you need a function $\varphi$ as follows,
$$
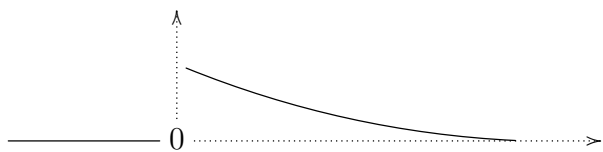
\xymatrix@R=18pt@C=25pt{
&&&&&\\
&&\ar@/_/@{-}[drrrr]&&&&&\\
\ar@{-}[rr]&&0\ar@.[uu]\ar@.[rrrrr]&&&&&}$$
the simplest example being the exponential, due to $\int_0^\infty e^{-x}dx=1$.

\item In the case of full support on $\mathbb R$, you need a bell-shaped curve, as follows,
$$
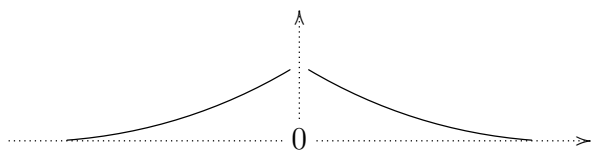
\xymatrix@R=18pt@C=15pt{
&&&&&&&&\\
&&&&&\ar@/_/@{-}[drrrr]\ar@/^/@{-}[dllll]&&&&&\\
\ar@{.}[rrrrr]&&&&&0\ar@.[uu]\ar@.[rrrrr]&&&&&}$$
and you can use here for instance the functions in (2), symmetrized.
\end{enumerate}
\end{recipe}

Now in view of this, what should be our next distributions? By leaving (3) aside, we are left between choosing between (1) with more complicated shapes, such as semicircles, and (2) with the exponential there. Which is a tough choice, believe me, and we will choose here (2) with the exponential, leaving the semicircles for the next chapter.

\bigskip

Getting started now, with the exponentials, as always when constructing probability distributions, it is convenient to make use of a real parameter, say $\lambda$, if that parameter looks like something quite natural to be added. We are led in this way to:

\begin{definition}
We say that a random variable $f:X\to\mathbb R$ is exponential, with parameter $\lambda>0$, if its law is as follows:
$$\mu(x)=\lambda e^{-\lambda x}\chi_{[0,\infty]}(x)dx$$
That is, the density of $f$ must be $\varphi(x)=\lambda e^{-\lambda x}$, over the positive numbers.
\end{definition}

As a first observation, our measures are indeed of mass 1, as shown by the following computation, which is arguably the simplest integral computation yielding 1:
\begin{eqnarray*}
\mu(\mathbb R)
&=&\int_\mathbb R\lambda e^{-\lambda x}\chi_{[0,\infty]}(x)dx\\
&=&\int_0^\infty \lambda e^{-\lambda x}dx\\
&=&\left[-e^{-\lambda x}\right]_0^\infty\\
&=&0-(-1)\\
&=&1
\end{eqnarray*}

So long for the mathematical phenomenology of these laws, hope we agree on this, the exponential laws are something very simple, naturally coming after the uniform ones. In what regards now the real-life phenomenology, many things can be said here, and as already mentioned, my own favorite example comes from particle physics:

\begin{fact}
The exponential laws model the decay of radioactive particles, which decay is, notoriously, exponential.
\end{fact}

And more on this later, the point being that it is better to get into physics once you know some mathematics, matter of being able to prove some theorems there. Getting now to the mathematics, we will follow our previous treatment of the uniform laws. As our first result, which is something trivial, coming from definitions, we have:

\begin{proposition}
A variable $f:X\to\mathbb R$ follows the exponential law of parameter $\lambda>0$ precisely when the following happens, for any measurable function $\phi:\mathbb R\to\mathbb R$:
$$E(\phi(f))=\lambda\int_0^\infty\phi(x)e^{-\lambda x}dx$$
Equivalently, we must have $f:X\to[0,\infty)$, up to a set of null measure, and given a measurable set $Y\subset[0,\infty)$, the following must happen:
$$P(f\in Y)=\lambda\int_Ye^{-\lambda x}dx$$
Also, in this latter statement, we can restrict the attention to the intervals $Y\subset[0,\infty)$.
\end{proposition}

\begin{proof}
This is indeed something self-explanatory, the idea being as follows:

\medskip

(1) Regarding the formula of $E(\phi(f))$, and the converse assertion too, that comes from our definition of the law of a random variable, as axiomatized in chapter 5.

\medskip

(2) Regarding the formula of $P(f\in Y)$, that comes from (1) with $\phi=\chi_Y$, and the converse is clear too, by using the measure theory developed in chapter 5.

\medskip

(3) Finally, the last assertion, stating that we can restrict the attention to the intervals $Y\subset[a,b]$, is clear too, again by using the measure theory developed in chapter 5.
\end{proof}

Let us record as well a result regarding the cumulative distribution function:

\begin{proposition}
Given a variable $f:X\to\mathbb R$ which follows the exponential law of parameter $\lambda>0$, its cumulative distribution function (CDF) is given by
$$F(x)=1-e^{-\lambda x}$$
with the picture of this CDF being as follows,
\vskip-5mm
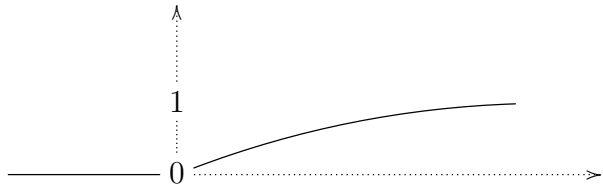
$$\xymatrix@R=12pt@C=25pt{
&&&&&\\
&&&&&\\
&&1\ar@.[uu]&&&&&\\
\ar@{-}[rr]&&0\ar@/^/@{-}[urrrr]\ar@{.}[u]\ar@.[rrrrr]&&&&&}$$
and the converse holds too, in the sense that having this CDF means that our variable $f:X\to\mathbb R$ must follow the exponential law of parameter $\lambda>0$.
\end{proposition}

\begin{proof}
We have indeed the following computation, for any $x>0$:
\begin{eqnarray*}
F(x)
&=&\int_{-\infty}^x\lambda e^{-\lambda t}\chi_{[0,\infty]}(t)dt\\
&=&\int_0^x\lambda e^{-\lambda t}dt\\
&=&\left[-e^{-\lambda t}\right]_0^x\\
&=&1-e^{-\lambda x}
\end{eqnarray*}

Thus, we are led to the various conclusions in the statement.
\end{proof}

Getting now to the moments, the result here is as follows:

\begin{theorem}
Given a variable $f:X\to\mathbb R$ which follows the exponential law of parameter $\lambda>0$, its moments are given by the following formula:
$$M_k=\frac{k!}{\lambda^k}$$
In particular, the mean is $E=1/\lambda$, and the variance is $V=1/\lambda^2$.
\end{theorem}

\begin{proof}
We use the general integration formula from Proposition 6.12, namely:
$$E(\phi(f))=\lambda\int_0^\infty\phi(x)e^{-\lambda x}dx$$

With the function $\phi(x)=x^k$ we obtain in this way, by using a straightforward partial integration, a recurrence formula for the $k$-th moment, as follows:
\begin{eqnarray*}
M_k
&=&\lambda\int_0^\infty x^ke^{-\lambda x}dx\\
&=&\lambda\int_0^\infty kx^{k-1}\cdot\frac{e^{-\lambda x}}{\lambda}\,dx\\
&=&\frac{k}{\lambda}\cdot\lambda\int_0^\infty x^{k-1}e^{-\lambda x}dx\\
&=&\frac{k}{\lambda}\cdot M_{k-1}
\end{eqnarray*}

Now since the $0$-th moment is $M_0=1$, we obtain by recurrence, as claimed:
$$M_k=\frac{k!}{\lambda^k}$$

Finally, the expectation is $E=M_1=1/\lambda$, and the variance is given by:
$$V=M_2-M_1^2=\frac{2}{\lambda^2}-\frac{1}{\lambda^2}=\frac{1}{\lambda^2}$$

Thus, we are led to the conclusions in the statement.
\end{proof}

Good news, with the mathematics that we have we can do some physics:

\begin{theorem}
In the context of decay, the quantity to look at is the decay rate $\lambda$, which is the probability per unit time that the particle will disintegrate. With this:
\begin{enumerate}
\item The number of particles remaining at time $t>0$ is $N_t=e^{-\lambda t}N_0$.

\item The mean lifetime of a particle is $\tau=1/\lambda$.

\item The half-life of the substance is $t_{1/2}=(\log 2)/\lambda$. 
\end{enumerate}
\end{theorem}

\begin{proof}
This can be deduced indeed from our knowledge of the exponential laws:

\medskip

(1) In mathematical terms, our definition of the decay rate reads:
$$\frac{dN}{dt}=-\lambda N$$

By integrating, we are led to the formula in the statement, namely:
$$N_t=e^{-\lambda t}N_0$$

(2) Let us convert now what we have into a probability law. We have:
$$\int_0^\infty N_tdt=\int_0^\infty N_0e^{-\lambda t}dt=\frac{N_0}{\lambda}$$

We conclude that the density of the probability decay function is given by:
$$\varphi(t)=\frac{\lambda}{N_0}\cdot N_0e^{-\lambda t}=\lambda e^{-\lambda t}$$

But this is the density of the exponential law of parameter $\lambda$, and the mean lifetime therefore appears as the expectation of this law, which according to Theorem 6.14 is:
$$\tau=\frac{1}{\lambda}$$

(3) Finally, regarding the half-life, this is by definition the time $t_{1/2}$ required for the decaying quantity to fall to one-half of its initial value. Mathematically, this means:
$$N_t=2^{-\frac{t}{t_{1/2}}}N_0$$

Now by comparing with $N_t=e^{-\lambda t}N_0$, this gives $t_{1/2}=(\log 2)/\lambda$, as stated.
\end{proof}

Quite nice all this, and needless to say, other occurrences of the exponential laws abound, in the real life, and exercise for you, to learn more about all this. 

\bigskip

Getting back now to mathematics, as a continuation of Theorem 6.14, let us compute now the central moments as well. And here, surprise, while the central moments are meant in general to simplify the moment combinatorics, for our exponential laws these are in fact given by a more complicated formula, the result being as follows:

\begin{theorem}
Given a variable $f:X\to\mathbb R$ which follows the exponential law of parameter $\lambda>0$, its central moments are given by the following formula,
$$M_k'=\frac{k!}{\lambda^k}\sum_{s=0}^k\frac{(-1)^s}{s!}$$
with the small order central moments being as follows,
$$M_1'=0\quad,\quad 
M_2'=\frac{1}{\lambda^2}\quad,\quad
M_3'=\frac{2}{\lambda^3}\quad,\quad
M_4'=\frac{9}{\lambda^4}\quad,\quad\ldots$$
and with all this being more complicated than for the usual moments.
\end{theorem}

\begin{proof}
We recall from chapter 2 that the central moments of a variable $f:X\to\mathbb R$, designed for having $M_1'=0$, $M_2'=V$, are the following numbers, with $E=E(f)$:
$$M_k'=E((f-E)^k)$$

Thus, we are in need of a remake of the computation from the proof of Theorem 6.14. As there, we will use the general integration formula from Proposition 6.12, namely:
$$E(\phi(f))=\lambda\int_0^\infty\phi(x)e^{-\lambda x}dx$$

With the function $\phi(x)=(x-E)^k$, where $E=1/\lambda$, as computed in Theorem 6.14, we obtain in this way a recurrence formula for the $k$-th central moment, as follows:
\begin{eqnarray*}
M_k'
&=&\lambda\int_0^\infty\left(x-\frac{1}{\lambda}\right)^ke^{-\lambda x}dx\\
&=&\lambda\left[\int_0^\infty k\left(x-\frac{1}{\lambda}\right)^{k-1}\cdot\frac{e^{-\lambda x}}{\lambda}\,dx+\left(0-\frac{1}{\lambda}\right)^k\frac{e^{-\lambda\cdot0}}{\lambda}\right]\\
&=&\frac{k}{\lambda}\cdot\lambda\int_0^\infty\left(x-\frac{1}{\lambda}\right)^{k-1}e^{-\lambda x}dx+\left(-\frac{1}{\lambda}\right)^k\\
&=&\frac{k}{\lambda}\cdot M_{k-1}'+\frac{(-1)^k}{\lambda^k}
\end{eqnarray*}

Now let us see how this recurrence works. With $M_1'=0$ as input, we obtain:
$$M_2'=\frac{2}{\lambda}\cdot0+\frac{1}{\lambda^2}=\frac{1}{\lambda^2}$$
$$M_3'=\frac{3}{\lambda}\cdot\frac{1}{\lambda^2}-\frac{1}{\lambda^3}=\frac{2}{\lambda^3}$$
$$M_4'=\frac{4}{\lambda}\cdot\frac{2}{\lambda^3}+\frac{1}{\lambda^4}=\frac{9}{\lambda^4}$$

Thus, good news, we have the low order moments in the statement. In order to solve now the recurrence in general, with a bit of recurrence know-how, let us write:
$$M_k'=\frac{k!}{\lambda^k}\,N_k$$

With this convention, the recurrence formula found above takes the following form:
$$\frac{k!}{\lambda^k}\,N_k=\frac{k}{\lambda}\cdot\frac{(k-1)!}{\lambda^{k-1}}\,N_{k-1}+\frac{(-1)^k}{\lambda^k}$$

Now by multiplying by $\lambda^k/k!$, we are led to something very simple, namely:
$$N_k=N_{k-1}+\frac{(-1)^k}{k!}$$

The initial data being $N_1=0$, we conclude from this that we have:
$$N_k=\sum_{s=2}^k\frac{(-1)^s}{s!}=\sum_{s=0}^k\frac{(-1)^s}{s!}$$

To be more precise, the first formula is what we get by recurrence, at $k\geq2$, and the passage to the second formula, done by adding $1-1=0$, is there for our final formula to work at $k=1$ too. Thus, we are led to the formula of $M_k'$ in the statement.
\end{proof}

Getting now to the normalized central moments, we have here:

\begin{theorem}
Given a variable $f:X\to\mathbb R$ which follows the exponential law of parameter $\lambda>0$, its normalized central moments are given by the following formula,
$$M_k''=\sum_{s=0}^k(-1)^s\frac{k!}{s!}$$
with the small order normalized central moments being as follows,
$$M_1''=0\quad,\quad 
M_2''=1\quad,\quad
M_3''=2\quad,\quad
M_4''=9\quad,\quad\ldots$$
and with all this being, again, more complicated than for the usual moments.
\end{theorem}

\begin{proof}
We recall that the normalized central moments are the following numbers, with $E=E(f)$ being the expectation, and $\sigma=\sqrt{V}$ being the standard deviation:
$$M''_k=E\left(\left(\frac{f-E}{\sigma}\right)^k\right)=\frac{M_k'}{\sigma^k}$$

In our case the standard deviation is $\sigma=\sqrt{V}=1/\lambda$, and we obtain:
$$M_k''=\lambda^kM_k'=k!\sum_{s=0}^k\frac{(-1)^s}{s!}=\sum_{s=0}^k(-1)^s\frac{k!}{s!}$$

Thus, we are led to the conclusions in the statement. 
\end{proof}

Let us end this discussion with the following result, capturing the essentials:

\begin{theorem}
For a variable $f:X\to\mathbb R$ which is exponential of parameter $\lambda>0$,
$$E=\frac{1}{\lambda}\quad,\quad V=\frac{1}{\lambda^2}\quad,\quad\gamma=2\quad,\quad\kappa=9$$
are the mean, variance, skewness and kurtosis. 
\end{theorem}

\begin{proof}
Here the formulae of the mean $E$ and variance $V$ are those from Theorem 6.14, and the formulae of the skewness $\gamma=M_3''$ and kurtosis $\kappa=M_4''$ are those from Theorem 6.17. Thus, we are led to the conclusions in the statement.
\end{proof}

Regarding now the Fourier transform computation, this is as follows:

\begin{theorem}
For a variable $f:X\to\mathbb R$ which is exponential of parameter $\lambda>0$,
$$F_f(x)=\frac{\lambda}{\lambda-ix}$$
is its Fourier transform, $F_f(x)=E(e^{ixf})$. 
\end{theorem}

\begin{proof}
This is again a standard computation, based on the general integration formula in Proposition 6.12, by using the function $\phi(t)=e^{ixt}$ there, which gives:
\begin{eqnarray*}
F_f(x)
&=&E(e^{ixf})\\
&=&\lambda\int_0^\infty e^{ixt}e^{-\lambda t}dt\\
&=&\lambda\int_0^\infty e^{(ix-\lambda)t}dt\\
&=&\lambda\left[\frac{e^{(ix-\lambda)t}}{ix-\lambda}\right]_0^\infty\\
&=&\lambda\left(0-\frac{1}{ix-\lambda}\right)\\
&=&\frac{\lambda}{\lambda-ix}
\end{eqnarray*}

Thus, we are led to the conclusion in the statement.
\end{proof}

And we will stop our preliminary learning of the exponential laws here. As already mentioned, these laws appear in many situations, in the real life, and the various formulae established above are basically what you need, in order to study these situations.

\section*{6c. Memory loss}

We would like to discuss now an interesting property of the exponential laws, namely their loss of memory. Which might sound not very exciting to young guys like you, but listen to the old man here, memory loss is an interesting phenomenon, appearing in the real life, that we should know more about. Here is what we want to prove:

\begin{fact}
The following happen, in the context of general probability theory:
\begin{enumerate}
\item The exponential laws experience memory loss, in the sense that they accurately predict the future, without remembering the past. 

\item In fact, they are the only continuous laws having this property. As for the discrete case, the unique examples here are the geometric laws.
\end{enumerate}
\end{fact}

Which sounds quite good, even with a nice axiomatic flavor, for the laws involved above, namely exponential and geometric, hope you agree with me here. In practice now, looking at this, the problem is, what does (1) above exactly say? Not very clear, and with the key to such things being conditional probability, here will be our plan:

\begin{plan}
In order to understand all this, our plan will be:
\begin{enumerate}
\item Talk about conditional probability, something nice and very useful.

\item Prove that the exponential laws experience indeed memory loss.

\item Then discuss the converse, and the discrete case as well.
\end{enumerate}
\end{plan}

Getting started now, with conditional probability, that is something which is best learned in the discrete context, and should normally have been discussed in chapter 2, at the beginning, short after the definition of the discrete probability spaces.

\bigskip

This being said, we were busy with other things in chapter 2, and then in chapters 3-4 as well, and had not time for this. Nevermind. So, never too late for learning new things, and as a starting point, here is an interesting question, that I have for you:

\begin{question}
What is the probability of getting $12$ as total, when throwing $3$ dice, knowing that none of the dice falls on $1$?
\end{question}

And isn't this a bit perplexing, because despite our hard learning of probability theory so far, we are now at almost page 150 of the present book, having endured things like heavy calculus, heavy combinatorics, axiomatics and Rudin \cite{rud}, advice of cats and rats, and even the Cantor set $K$, we don't really have tools for solving such things.

\bigskip

Fortunately, there is a systematic approach to such questions, conditional probability. Indeed, by thinking a bit, we are led to the following simple, but very useful fact:

\index{conditional probability}

\begin{theorem}
The probability for $A$ to happen, knowing that $B$ happens, is
$$P(A|B)=\frac{P(A\cap B)}{P(B)}$$
with this being called conditional probability.
\end{theorem}

\begin{proof}
This is something quite obvious, by drawing a Venn diagram, as follows:
$$
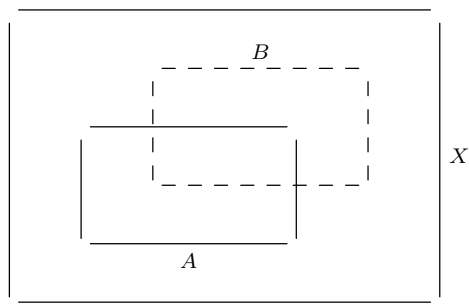
\xymatrix@R=15pt@C=20pt{
\ar@{-}[rrrrrr]\ar@{-}[ddddd]&&&&&&\ar@{-}[ddddd]^X\\
&&\ar@{--}[rrr]^B\ar@{--}[dd]&&&\ar@{--}[dd]\\
&\ar@{-}[rrr]\ar@{-}[dd]&&&\ar@{-}[dd]\\
&&\ar@{--}[rrr]&&&\\
&\ar@{-}[rrr]_A&&&\\
\ar@{-}[rrrrrr]&&&&&&}$$

Indeed, knowing that $B$ happens means that our computation takes place inside $B$. Next, we are interested in the probability for $A$ to happen, and since $B$ happens, this means that we are interested in the probability for $A\cap B$ to happen. Thus, we have:
$$P(A|B)=\frac{P(A\cap B)}{P(B)}$$

Summarizing, we are led to the conclusion in the statement.
\end{proof}

As an illustration for the above formula, let us go back to Question 6.22. By using Theorem 6.23, can now answer that question, majestically, as follows:

\begin{answer}
The probability of getting $12$, with none of our $3$ dice falling on $1$, is
$$P=\frac{19/216}{125/216}=\frac{19}{125}$$
according to our general formula above, for conditional probability.
\end{answer}

To be more precise here, we have indeed the following self-explanatory computation, based on the general formula from Theorem 6.23:
\begin{eqnarray*}
P
&=&\frac{P(246)+P(255)+P(336)+P(345)+P(444)}{P(\geq2,\geq2,\geq2)}\\
&=&\frac{6/216+3/216+3/216+6/216+1/216}{125/216}\\
&=&\frac{19/216}{125/216}\\
&=&\frac{19}{125}
\end{eqnarray*}

Very nice, all this. Still in relation with conditional probability questions, at a more theoretical level, we have the following well-known and useful formula, due to Bayes:

\index{Bayes formula}

\begin{theorem}
We have the following formula,
$$P(A|B)=\frac{P(B|A)P(A)}{P(B)}$$
involving conditional probabilities.
\end{theorem}

\begin{proof}
This is something which might look quite difficult to establish, just by thinking at probabilities, but which follows easily by using our mathematical formalism for probability, and for conditional probability in particular. Recall indeed that we have:
$$P(A|B)=\frac{P(A\cap B)}{P(B)}$$

Now observe that this equality can be written in the following way:
$$P(A|B)P(B)=P(A\cap B)$$

By symmetry reasons, the following formula must hold as well:
$$P(B|A)P(A)=P(A\cap B)$$

Now by forgetting about $P(A\cap B)$, we conclude that the following holds:
$$P(A|B)P(B)=P(B|A)P(A)$$

But this leads to the above Bayes formula, as stated.
\end{proof}

In practice, there are many interesting applications of the Bayes formula. Here is a first illustration, which is quite hard to come upon, with bare hands:

\begin{example}
Assume that we have $100$ dice, $75$ of which are normal, and $25$ of which are biased, with $P(6)=1/2$. If we pick a die, throw it, and get $6$, then this means that the probability of the die that we picked to be biased is $1/2$. 
\end{example}

Indeed, according to the Bayes formula, the probability that we are interested in, that of having a biased die, knowing that we got a 6 after throwing it, is:
$$P({\rm biased}|6)=\frac{P(6|{\rm biased})P({\rm biased})}{P(6)}$$

Now in this formula, the numbers on top are already known to us, given by:
$$P(6|{\rm biased})=\frac{1}{2}\quad,\quad P({\rm biased})=\frac{25}{100}=\frac{1}{4}$$

As for the number on the bottom, this can be computed as follows:
\begin{eqnarray*}
P(6)
&=&P(6,{\rm normal})+P(6,{\rm biased})\\
&=&P(6|{\rm normal})P({\rm normal})+P(6|{\rm biased})P({\rm biased})\\
&=&\frac{1}{6}\times\frac{3}{4}+\frac{1}{2}\times\frac{1}{4}\\
&=&\frac{1}{8}+\frac{1}{8}\\
&=&\frac{1}{4}
\end{eqnarray*}

Thus, we can finish our computation, and we obtain indeed probability $1/2$:
$$P({\rm biased}|6)
=\frac{P(6|{\rm biased})P({\rm biased})}{P(6)}
=\frac{1/2\times 1/4}{1/4}
=\frac{1}{2}$$

And with this, end of our general discussion regarding conditional probability. Getting back now to our regular business in this chapter, exponential laws, we have:

\begin{theorem}
An exponential variable $f:X\to\mathbb R$ experiences loss of memory:
$$P(f>s+t|f>s)=P(f>t)\quad,\quad\forall s,t>0$$
Conversely, the only continuous variables having this property are the exponential ones. 
\end{theorem}

\begin{proof}
We have two things to be proved, the idea being as follows:

\medskip

(1) Regarding the first assertion, this comes from a standard computation, based on the CDF formula from Proposition 6.13, as follows, $\lambda>0$ being the parameter:
\begin{eqnarray*}
P(f>s+t|f>s)
&=&\frac{P((f>s+t)\cap(f>s))}{P(f>s)}\\
&=&\frac{P(f>s+t)}{P(f>s)}\\
&=&\frac{e^{-\lambda(s+t)}}{e^{-\lambda s}}\\
&=&e^{-\lambda t}\\
&=&P(f>t)
\end{eqnarray*}

(2) Conversely, assume that $f:X\to\mathbb R$ experiences loss of memory. In practice, this means that the following formula must be satisfied, for any $s,t>0$:
$$P(f>s+t)=P(f>s)P(f>t)$$

But this tells us that $\psi(t)=P(f>t)$ has the property $\psi(s+t)=\psi(s)\psi(t)$, and by some routine analysis, involving integers, then rationals, and then reals, say good exercise for you, this function $\psi$ must be a power function, $\psi(t)=a^t$. But this is the same as saying that we have $\psi(t)=e^{-\lambda t}$ for some $\lambda\in\mathbb R$, which back to $f$, means:
$$P(f>t)=e^{-\lambda t}$$

And with this, done, because we must have $\lambda>0$, as to have the correct limits at $t=\pm\infty$, and then, the rest follows by using the CDF theory from Proposition 6.13.
\end{proof}

Let us record as well the result in the discrete case, which is as follows:

\begin{theorem}
A geometric variable $f:X\to\mathbb N$ experiences loss of memory:
$$P(f>s+t|f>s)=P(f>t)\quad,\quad\forall s,t>0$$
Conversely, the only discrete variables having this property are the geometric ones. 
\end{theorem}

\begin{proof}
It is most convenient here to shoot both rabbits at the same time. In order for $f:X\to\mathbb N$ to experience loss of memory, we must have, for any $s,t>0$:
$$P(f>s+t)=P(f>s)P(f>t)$$

As before in the continuous case, this tells us that $\psi(t)=P(f>t)$ has the property $\psi(s+t)=\psi(s)\psi(t)$, which in turn means that $\psi$ must be a power function, $\psi(t)=a^t$. Now by getting back to our variable $f$, the conclusion is that we must have:
$$P(f>t)=a^t$$

Now since the law was assumed to be supported by $\mathbb N$, this reads $P(t)=a^t(1-a)$, and so $f$ must follow a geometric law of parameter $p=1-a$, as desired.
\end{proof}

\section*{6d. Independence}

Looking at what we did so far in the continuous case, almost 50 pages, we certainly have some theory started, but when comparing to what we did before in Part I, in the discrete case, still long way to go, with far more distributions and phenomena in need to be investigated, as to get on par with that. But no worries here, still 250 pages to go.

\bigskip

This being said, as a most pressing issue, arising from this comparison, we must definitely talk about independence. So, we will do this now, first with various generalities, in the spirit of what we did before in chapter 2, in the discrete case, and then, time permitting, with a look at this, in the context of the uniform and exponential laws.

\bigskip

Let us start our study with the following straightforward definition:

\index{independence}
\index{independent variables}

\begin{definition}
We say that two variables $f,g:X\to\mathbb R$ are independent when
$$P(f\in I,g\in J)=P(f\in I)P(g\in J)$$
happens, for any two intervals $I,J\subset\mathbb R$.
\end{definition}

As a first observation, in the discrete case this is the same as our notion of independence from chapter 2, which was asking for the following equality, for any $x,y\in\mathbb R$:
$$P(f=x,g=y)=P(f=x)P(g=y)$$

Indeed, the equivalence comes from the following computation:
\begin{eqnarray*}
P(f\in I,g\in J)
&=&\sum_{x\in I}\sum_{y\in J}P(f=x,g=y)\\
&=&\sum_{x\in I}\sum_{y\in J}P(f=x)P(g=y)\\
&=&\sum_{x\in I}P(f=x)\sum_{y\in J}P(g=y)\\
&=&P(f\in I)P(g\in J)
\end{eqnarray*}

As our first result now, we can in fact recycle the above computation in the continuous setting, by replacing the points $x,y$ there by intervals, and we are led to:

\begin{proposition}
Two variables $f,g:X\to\mathbb R$ are independent precisely when
$$P(f\in K,g\in L)=P(f\in K)P(g\in L)$$
happens, for any two measurable sets $K,L\subset\mathbb R$.
\end{proposition}

\begin{proof}
In the case where our measurable sets $K,L\subset\mathbb R$ are disjoint unions of intervals, $K=\sqcup_rI_r$ and $L=\sqcup_sJ_s$, we have indeed the following computation:
\begin{eqnarray*}
P(f\in K,g\in L)
&=&\sum_r\sum_sP(f\in I_r,g\in J_s)\\
&=&\sum_r\sum_sP(f\in I_r)P(g\in J_s)\\
&=&\sum_rP(f\in I_r)\sum_sP(g\in J_s)\\
&=&P(f\in K)P(g\in L)
\end{eqnarray*}

As for the general case, this follows from this, by using the theory from chapter 5.
\end{proof}

Importantly now, Proposition 6.30 allows us to pass to expectations, as follows:

\begin{proposition}
Two variables $f,g:X\to\mathbb R$ are independent precisely when
$$E(\phi(f)\psi(g))=E(\phi(f))E(\psi(g))$$
happens, for any two measurable functions $\phi,\psi:\mathbb R\to\mathbb R$.
\end{proposition}

\begin{proof}
We know from Proposition 6.30 that the independence of $f,g:X\to\mathbb R$ is equivalent to the following condition, for any two measurable sets $K,L\subset\mathbb R$:
$$P(f\in K,g\in L)=P(f\in K)P(g\in L)$$

Now observe that in terms of characteristic functions, this condition reads:
$$E(\chi_K(f)\chi_L(g))=E(\chi_K(f))E(\chi_L(g))$$

But with this, we are done, by approximating $\phi,\psi:\mathbb R\to\mathbb R$ by step functions.
\end{proof}

As a third and main result now on the subject, generalizing what we knew from chapter 2, in the discrete case, we have the moment characterization of independence:

\index{mixed moments}

\begin{theorem}
Assuming that $f,g:X\to\mathbb R$ are independent, we have
$$E(f^kg^l)=E(f^k)E(g^l)$$
and the converse holds, in the sense that this implies the independence of $f,g$.
\end{theorem}

\begin{proof}
This folllows indeed from Proposition 6.31, by approximating the continuous functions, and so the measurable functions too, by polynomials. We will leave the details here, including some learning about polynomial approximation, as an exercise. 
\end{proof}

Still following the material from chapter 2, as a next task, we must talk about convolution of measures. In the present continuous case, this can be done as follows:

\index{convolution}
\index{semigroup}
\index{group}

\begin{definition}
Given a space $X$ with a sum operation $+$, we can define the convolution of any two probability measures on it by the following formula:
$$(\mu*\nu)(E)=(\mu\times\nu)\left((x,y)\Big|x+y\in E\right)$$
Equivalently, in terms of the associated integrals, we must have 
$$\int_Xf(x)d(\mu*\nu)(x)=\int_{X\times X}f(x+y)d\mu(x)d\nu(y)$$
for any measurable function $f:X\to\mathbb R$.
\end{definition}

As a first observation on this, which is something quite subtle, this fits with our previous notion of convolution, from chapter 2, for discrete laws, which was as follows:
$$\mu=\sum_ia_i\delta_{x_i}\quad,\quad \nu=\sum_jb_j\delta_{y_j}\quad\implies\quad 
\mu*\nu=\sum_{ij}a_ib_j\delta_{x_i+y_j}$$

Indeed, this amounts in defining $*$ by the formula $\delta_x*\delta_y=\delta_{x+y}$, and linearity. But, according to Definition 6.33, the convolution of two Dirac masses is, at is should:
$$(\delta_x*\delta_y)(E)=(\delta_x\times\delta_y)\left((x,y)\Big|x+y\in E\right)=\delta_{x+y}(E)$$

In general now, observe that our $*$ operation in Definition 6.33 is indeed well-defined, the total mass being 1. Also, the equivalent formulation, in terms of integrals, is clear by using characteristic functions, as shown by the following computation:
\begin{eqnarray*}
(\mu*\nu)(E)
&=&\int_X\chi_E(x)d(\mu*\nu)(x)\\
&=&\int_{X\times X}\chi_E(x+y)d\mu(x)d\nu(y)\\
&=&(\mu\times\nu)\left((x,y)\Big|x+y\in E\right)
\end{eqnarray*}

In relation now with the notion of independence, we have the following key result, extending what we previously knew from chapter 2, in the discrete case:

\begin{theorem}
Assuming that $f,g:X\to\mathbb R$ are independent, we have
$$\mu_{f+g}=\mu_f*\mu_g$$
where $*$ is the convolution of real probability measures.
\end{theorem}

\begin{proof}
There are several possible proofs here, and with our present knowledge, which is rather functional analysis oriented, the simplest is to use moments. Indeed, we have the following computation, using the formula for mixed moments in Theorem 6.32:
\begin{eqnarray*}
M_k(f+g)
&=&E((f+g)^k)\\
&=&\sum_r\binom{k}{r}E(f^rg^{k-r})\\
&=&\sum_r\binom{k}{r}M_r(f)M_{k-r}(g)
\end{eqnarray*}

On the other hand, if $h:X\to\mathbb R$ is a random variable following the law $\mu_f*\mu_g$, the corresponding moments are given by the following formula:
\begin{eqnarray*}
M_k(h)
&=&\int_Xx^kd(\mu_f*\mu_g)(x)\\
&=&\int_{X\times X}(x+y)^kd\mu_f(x)d\mu_g(y)\\
&=&\sum_r\binom{k}{r}\int_Xx^rd\mu_f(x)\int_Xy^{k-r}d\mu_g(y)\\
&=&\sum_r\binom{k}{r}M_r(f)M_{k-r}(g)
\end{eqnarray*}

Thus our laws have the same moments, so they coincide, as stated.
\end{proof}

Still in relation with the basic theory of independence, here is now a second result, coming as a continuation of Theorem 6.34, which is something more advanced:

\index{independence}
\index{Fourier transform}

\begin{theorem}
Assuming that $f,g:X\to\mathbb R$ are independent, we have
$$F_{f+g}=F_fF_g$$
where $F_f(x)=E(e^{ixf})$ is the Fourier transform.
\end{theorem}

\begin{proof}
We have the following computation, using Theorem 6.34:
\begin{eqnarray*}
F_{f+g}(x)
&=&\int_Xe^{ixz}d\mu_{f+g}(z)\\
&=&\int_Xe^{ixz}d(\mu_f*\mu_g)(z)\\
&=&\int_{X\times X}e^{ix(z+t)}d\mu_f(z)d\mu_g(t)\\
&=&\int_Xe^{ixz}d\mu_f(z)\int_Xe^{ixt}d\mu_g(t)\\
&=&F_f(x)F_g(x)
\end{eqnarray*}

Thus, we are led to the conclusion in the statement.
\end{proof}

As a comment here, which might actually come a bit late, because we know about Fourier since chapter 4, you might wonder what $i\in\mathbb C$ has to do with all this. Good point, and this is something quite subtle, more on this later, trust me in the meantime.

\bigskip

Let us also mention that the converses of Theorems 6.34 and 6.35 do not necessarily hold, with what comes from there being known as ``subindependence''. And, exercise of course for you to learn more about this, notably with some counterexamples.

\bigskip

And with this, end of our preliminary study of abstract independence. Getting back now to the laws that we investigated in this chapter, namely the uniform and exponential ones, many interesting things can be said, for instance in relation with the sums of independent uniform variables, or the sums of independent exponential variables.

\bigskip

However, such questions, which can be a bit technical, can be investigated in a straightforward way, by using the Fourier transform theory from Theorem 6.35, and the explicit Fourier transform formulae for our laws, from Theorem 6.8 and Theorem 6.19.

\bigskip

Instead, we would like to talk about an exciting topic, which is of main theoretical interest, namely the infinite divisibility of probability measures. Let us start with: 

\begin{definition}
We call a real probability measure $\mu$ infinitely divisible when
$$\underbrace{\nu*\ldots*\nu}_{n\ terms}=\mu$$
has solutions, for any $n\in\mathbb N$. That is, $\mu$ must admit $n$-th roots, for any $n\in\mathbb N$.
\end{definition}

Obviously, this is something quite subtle, and in relation with this, many things can be said, of various levels of difficulty. Here are a few basic observations, about this:

\bigskip

(1) As a first remark, once we have a $n$-th root as above, by convolving this $n$-th root $m$ times with itself, we obtain a measure that we would like to denote as follows:
$$\mu^{*m/n}=\underbrace{\nu*\ldots*\nu}_{m\ terms}$$

(2) However, this remains a bit vague, because we still have uniqueness issues to be discussed. Then, assuming these uniqueness issues solved, next comes the problem of talking about arbitrary multiples of our measure, $\mu^{*r}$ with $r>0$ real. And so on.

\bigskip

(3) As another comment, the Fourier transform technology from Theorem 6.35 provides a key for the study of infinite divisibility, due to the following implication:
$$\nu^{*n}=\mu\implies F_\nu(x)^n=F_\mu(x)$$

(4) Summarizing, all this leads us into the question of recapturing probability measures from their Fourier transforms. And more on such things, later in this book. 

\bigskip

In relation now with the various laws that we know, from this book so far, we have:

\begin{theorem}
The following happen:
\begin{enumerate}
\item The Dirac masses are infinitely divisible.

\item The measures supported by $2\leq n<\infty$ points are not infinitely divisible.

\item In particular, the binomial laws are not infinitely divisible.

\item The negative binomial laws, including the geometric ones, are infinitely divisible.

\item The Poisson laws are infinitely divisible, and so are their compound versions.

\item The uniform laws are not infinitely divisible, while the exponential laws are.
\end{enumerate}
\end{theorem}

\begin{proof}
This is a mixture of trivial and non-trivial facts, as follows:

\medskip

(1-3) The first assertion comes from the formula $\delta_{x/n}^{*n}=\delta_x$. The second assertion comes from some straightforward cardinality considerations, that we will leave as an exercise. As for the third assertion, this comes right away from the second one.

\medskip

(4-6) Here things are elementary for the negative binomial laws, trivial via the PLT for the Poisson and compound Poisson laws, and less trivial for the continuous laws. We will leave some study here as an exercise, and come back to this, later in this book.
\end{proof}

\section*{6e. Exercises}

Welcome, eventually, to concrete continuous laws, and as exercises, we have:

\begin{exercise}
Work out formulae, for uniform measures on unions of intervals.
\end{exercise}

\begin{exercise}
Work out formulae, for uniform measures on arbitrary sets.
\end{exercise}

\begin{exercise}
Learn about the physical mechanisms behind particle decay.
\end{exercise}

\begin{exercise}
Learn about the other occurences of the exponential laws.
\end{exercise}

\begin{exercise}
Compute more conditional probabilities, of your choice.
\end{exercise}

\begin{exercise}
Further meditate at the Bayes formula, is that trivial or not?
\end{exercise}

\begin{exercise}
Clarify the things that we said, in relation with independence.
\end{exercise}

\begin{exercise}
Have some work done on the notion of infinite divisibility.
\end{exercise}

As bonus exercise, read a physics book. This is where laws come from.

\chapter{Semicircle laws}

\section*{7a. Random walks}

We have seen in the previous chapter that some interesting theory can be developed for the uniform and exponential laws. And with these being, respectively, the simplest continuous laws which are compactly supported, and non-compactly supported.

\bigskip

What is next? More complicated laws I guess, again compactly or non-compactly supported, and kept as simple as possible. And with this being something subjective, we will need someone to make a choice, and that will be most likely job for me, author.

\bigskip

So, based on years of experience in quantum physics, we will get soon to that sub-quark particles making the world go round, that is promised, and on years of teaching too, here is my claim, which will direct the material from the present chapter: 

\begin{claim}
The next simplest continuous law is the Wigner semicircle law
$$\gamma_1=\frac{1}{2\pi}\sqrt{4-x^2}\,dx$$
which is something very simple, both mathematically and physically speaking. 
\end{claim}

Which might sound of course a bit bold, and here is the rationale, behind this:

\bigskip

(1) As a first comment, mathematically speaking, when looking for laws having a density, $d\mu(x)=\varphi(x)dx$, what you need is a function $\varphi\geq0$ integrating up to 1. And with the rectangles being already studied, leading to the uniform laws, it is a no-brainer to look next at the semicircles, which produce the Wigner law. Agreed, I hope.

\bigskip

(2) As a next piece of evidence, still mathematically speaking, as we will soon discover, the combinatorics of the Wigner law is something very interesting, involving the Catalan numbers $C_k=\frac{1}{k+1}\binom{2k}{k}$, which are at the core of advanced combinatorics. In short, the study here is both foundational and fun, and always a pleasure to teach this in class. 

\bigskip

(3) Finally, regarding the physics, the Wigner law is known to appear in a myriad of interesting situations, involving random walks, statistical mechanics, random matrices, phase transitions, operator algebras, quantum groups, and all sorts of microscopic beasts, ranging from heavy atoms to hypothetical very light particles. Very nice too.

\bigskip

Getting to work now, we have two things to do, as usual when meeting a new law, namely talking about its phenomenology, and then developing the math, with moments and everything. So, let us begin with phenomenology. And here, looking at what we have in (3), random walks are the simplest topic, and here will be our starting question:

\begin{question}
Given a graph $X$, with a distinguished vertex $*$:
\begin{enumerate}
\item What is the number $L_k$ of length $k$ loops on $X$, based at $*$? 

\item Equivalently, what is the measure $\mu$ having $L_k$ as moments?
\end{enumerate}
\end{question}

To be more precise, we are mainly interested in the first question, counting loops on graphs, with this being notoriously related to many applied mathematics questions, of discrete type. Think for instance percolation, with the water droplets making their way through coffee particles, that obviously corresponds to a random walk on a graph. And by making the technical assumption that the paths are actually loops, which helps with the mathematics, we are led into counting length $k$ loops on a graph, as in (1). 

\bigskip

As for question (2), this is a technical, useful probabilistic reformulation of the first question, that we will usually prefer, in what follows. There are of course many things that can be said here, and once you are reasonably advanced at this, you can even have this as starting question, and (1) as secondary question. However, let us not get here into this, which can be something quite technical, and keep Question 7.2 as stated, with (1), which is quite natural, as main question, and with (2) being a useful version of it.

\bigskip

Let us start our study with some graph theory basics. We first have:

\index{adjacency matrix}
\index{symmetric matrix}

\begin{proposition}
A graph $X$, with vertices labeled $1,\ldots,N$, is uniquely determined by its adjacency matrix, which is the matrix $d\in M_N(0,1)$ given by:
$$d_{ij}=\begin{cases}
1&{\rm if}\ i-j\\
0&{\rm if}\ {i\not\!\!-\,j}
\end{cases}$$
Moreover, the matrices $d\in M_N(0,1)$ which can appear in this way, from graphs, are precisely those which are symmetric, and have $0$ on the diagonal.
\end{proposition}

\begin{proof}
We have two things to be proved, the idea being as follows:

\medskip

(1) Given a graph $X$, we can construct a matrix $d\in M_N(0,1)$ as in the statement, and this matrix is obviously symmetric, and has $0$ on the diagonal.

\medskip

(2) Conversely, given a matrix $d\in M_N(0,1)$ which is symmetric, and has $0$ on the diagonal, we can associate to it the graph $X$ having as vertices the numbers $1,\ldots,N$, and with the edges between these vertices being defined as follows:
$$i-j\iff d_{ij}=1$$

It is then clear that the adjacency matrix of this graph $X$ is the matrix $d$ itself. Thus, we have established a correspondence $X\leftrightarrow d$, as in the statement.
\end{proof}

Getting now to the walks on graphs, this brings us into the mathematics of the corresponding adjacency matrix $d$, thanks to the following key observation:

\index{paths on graph}
\index{power of matrix}

\begin{proposition}
For a graph $X$, with adjacency matrix $d\in M_N(0,1)$, we have:
$$(d^k)_{ij}=\#\Big\{i=i_0-i_1-\ldots-i_{k-1}-i_k=j\Big\}$$
That is, the $k$-th power of $d$ describes the length $k$ paths on $X$.
\end{proposition}

\begin{proof}
According to the usual rule of matrix multiplication, the formula for the powers of the adjacency matrix $d\in M_N(0,1)$ is as follows:
\begin{eqnarray*}
(d^k)_{i_0i_k}
&=&\sum_{i_1,\ldots,i_{k-1}}d_{i_0i_1}d_{i_1i_2}\ldots d_{i_{k-1}i_k}\\
&=&\sum_{i_1,\ldots,i_{k-1}}\delta_{i_0-i_1}\delta_{i_1-i_2}\ldots\delta_{i_{k-1}-i_k}\\
&=&\sum_{i_1,\ldots,i_{k-1}}\delta_{i_0-i_1-\ldots-i_{k-1}-i_k}\\
&=&\#\Big\{i_0-i_1-\ldots-i_{k-1}-i_k\Big\}
\end{eqnarray*}

Thus, we are led to the conclusion in the statement.
\end{proof}

Of particular interest are the paths which begin and end at the same point. These are called loops, and in the case of loops, the formula in Proposition 7.4 leads to:

\index{loops on graph}
\index{diagonalization}

\begin{theorem}
For a graph $X$, with adjacency matrix $d\in M_N(0,1)$, we have:
$$(d^k)_{ii}=\#\Big\{k-{\rm loops\ based\ at\ }i\in I\Big\}$$
Also, the total number of $k$-loops on $X$, at various vertices, is the number
$$Tr(d^k)=\sum_i(d^k)_{ii}$$
which can be computed by diagonalizing $d$.
\end{theorem}

\begin{proof}
There are several things going on here, the idea being as follows:

\medskip

(1) The first assertion follows from Proposition 7.4, which at $i=j$ gives the following formula, which translates into the first formula in the statement:
$$(d^k)_{ii}=\#\Big\{i=i_0-i_1-\ldots-i_{k-1}-i_k=i\Big\}$$

(2) Regarding now the second assertion, this follows from the first one, simply by summing over all the vertices $i\in X$, which gives, as desired:
$$Tr(d^k)=\sum_i\#\Big\{k-{\rm loops\ based\ at\ }i\in I\Big\}$$

(3) Finally, the third assertion comes from linear algebra. Indeed, the adjacency matrix $d\in M_N(0,1)$ being symmetric, it diagonalizes over the reals, and we have:
\begin{eqnarray*}
d=U\lambda U^{-1}
&\implies&d^k=U\lambda ^kU^{-1}\\
&\implies&Tr(d^k)=Tr(\lambda ^k)
\end{eqnarray*}

Thus, if we denote by $\lambda_1,\ldots,\lambda_N\in\mathbb R$ the eigenvalues of $d$, we have:
$$Tr(d^k)=\lambda_1^k+\ldots+\lambda_N^k$$

Summarizing, we are led to the conclusions in the statement.
\end{proof}

Getting now to Question 7.2 as stated, we still have to talk about the associated probability measure $\mu$. So, here is the result that we will need, in what follows:

\index{adjacency matrix}
\index{random walk}

\begin{theorem}
Given a graph $X$, with adjacency matrix $d\in M_N(0,1)$, and with a distinguished vertex $*$, the number $L_k$ of loops of length $k$ based at $*$ is given by:
$$L_k=(d^k)_{**}$$
When writing $d=U\lambda U^{-1}$ with $U\in O_N$ and $\lambda=diag(\lambda_1,\ldots,\lambda_N)$ with $\lambda_i\in\mathbb R$, we have
$$L_k=\sum_iU_{*i}^2\lambda_i^k$$
and the real probability measure $\mu$ having these numbers as moments is given by
$$\mu=\sum_iU_{*i}^2\delta_{\lambda_i}$$
with the delta symbols standing as usual for Dirac masses.
\end{theorem}

\begin{proof}
The first formula comes from Theorem 7.5, and the second one can be deduced from it, by using the orthogonality of the passage matrix, $U^t=U^{-1}$, as follows:
\begin{eqnarray*}
L_k
&=&(d^k)_{**}\\
&=&(U\lambda^kU^t)_{**}\\
&=&\sum_iU_{*i}\lambda_i^k(U^t)_{i*}\\
&=&\sum_iU_{*i}^2\lambda_i^k
\end{eqnarray*}

Finally, the last assertion is clear from this, because the moments of the measure in the statement, $\mu=\sum_iU_{*i}^2\delta_{\lambda_i}$, are the following numbers, as they should:
$$M_k
=\int_\mathbb Rx^kd\mu(x)
=\sum_iU_{*i}^2\lambda_i^k
=L_k$$

Observe also that $\mu$ is indeed of mass 1, because all rows of $U\in O_N$ must be of norm 1, and so $\sum_iU_{*i}^2=1$. Thus, we are led to the conclusions in the statement.
\end{proof}

Getting now to the examples, what are the simplest non-trivial graphs $X$, that we can try to do some computations for? And here, we have 3 possible answers, as follows:

\index{circle graph}
\index{segment graph}

\begin{fact}
The following are graphs $X$, with a distinguished vertex $*\in X$:
\begin{enumerate}
\item The circle graph, having $N$ vertices, with $*$ being one of the vertices.

\item The segment graph, having $N$ vertices, with $*$ being the vertex at left.

\item The segment graph, having $2N+1$ vertices, with $*$ being in the middle.
\end{enumerate}
\end{fact}

So, let us start with these. For the circle, the computations are quite non-trivial, and you can try doing some, in order to understand what I am talking about. The problem comes from the fact that the loops of length $k=0,2,4,6,\ldots$ are quite easy to count, but then, once we pass $k=N$, the loops can turn around the circle or not, and they can even turn several times, and so on, and all this makes the count complicated. In addition, again due to loop turning, when $N$ is odd, we have as well loops of odd length.

\bigskip

As for the two segment graphs, here the computations look again complicated, and even more complicated than for the circle, because, again, once we pass $k=N$ many things can happen, and this makes the count quite complicated. As before, you can try doing some computations here, in order to understand what I am talking about. 

\bigskip

So, what to do? In the lack of a bright idea, let us pull an analysis trick, and formulate the following asymptotic result, which is of course something quite modest:

\begin{proposition}
For the circle graph, having $N$ vertices, the number of length $k$ loops based at one of the vertices is approximately
$$L_k\simeq\frac{2^k}{N}$$
in the $k\to\infty$ limit, when $N$ is odd, and is approximately
$$L_k\simeq\begin{cases}
\frac{2^{k+1}}{N}&(k\ {\rm even})\\
0&(k\ {\rm odd})
\end{cases}$$
also with $k\to\infty$, when $N$ is even.
\end{proposition}

\begin{proof}
This is something quite standard, the idea being as follows:

\medskip

(1) Consider the circle graph $X$, with vertices denoted $0,1,\ldots,N-1$. Since each vertex has valence 2, any length $k$ path based at 0 will consist of a binary choice at the beginning, then another binary choice afterwards, and so on up to a $k$-th binary choice at the end. Thus, there is a total of $2^k$ such paths, based at 0, and having length $k$. 

\medskip

(2) But now, based on the obvious ``uniformity'' of the circle, we can argue that, in the $k\to\infty$ limit, the endpoint of such a path will become random among the vertices $0,1,\ldots,N-1$. Thus, if we want this endpoint to be 0, as to have a loop, we have $1/N$ chances for this to happen, so the total number of loops is $L_k\simeq 2^k/N$, as stated.

\medskip

(3) With the remark, however, that the above argument works fine only when $N$ is odd. Indeed, when $N$ is even, the endpoint of a length $k$ path will be random among $0,2,\ldots,2N-2$ when $k$ is even, and random among $1,3,\ldots,2N-1$ when $k$ is odd. Thus for getting a loop we must assume that $k$ is even, and in this case the number of such loops is the total number of length $k$ paths, namely $2^k$, approximately divided by $N/2$, the number of points in $\{0,2,\ldots,2N-2\}$, which gives $L_k=2^k/(N/2)$, as stated.
\end{proof}

Getting now to the two segment graphs from Fact 7.7, in the context of the above proof, it is pretty much clear that for both, we lack the ``uniformity'' needed in (2), and this due to the 2 endpoints of the segment. In fact, thinking well, these graphs are no longer 2-valent, again due to the 2 endpoints, each having valence 1, and so even step (1) there must be fixed. And so, too bad, it is probably wiser to stop here.

\bigskip

This being said, as a matter of not leaving this subject empty-handed, let us record the following result, which is something more advanced, based on diagonalization:

\index{circle graph}
\index{Chebycheff polynomials}
\index{orthogonal polynomials}

\begin{theorem}
For the segment graphs having $N$ vertices, the characteristic polynomials $P_N$ of the adjacency matrices are subject to the following recurrence, 
$$P_0=1\quad,\quad P_1=x\quad,\quad P_{N+1}=xP_N-P_{N-1}$$
and are the well-known Chebycheff polynomials, enjoying many interesting properties. The corresponding numbers of loops and measures can be deduced from this.
\end{theorem} 

\begin{proof}
Obviously, many things going on here, the idea being as follows:

\medskip

(1) To start with, by computing determinants, we are led to the recurrence formula in the statement. Here is the proof at $N=4$, the general case being similar:
\begin{eqnarray*}
P_5
&=&\begin{vmatrix}
x&-1&0&0&0\\
-1&x&-1&0&0\\
0&-1&x&-1&0\\
0&0&-1&x&-1\\
0&0&0&-1&x
\end{vmatrix}\\
&=&x\begin{vmatrix}
x&-1&0&0\\
-1&x&-1&0\\
0&-1&x&-1\\
0&0&-1&x
\end{vmatrix}
+\begin{vmatrix}
-1&-1&0&0\\
0&x&-1&0\\
0&-1&x&-1\\
0&0&-1&x
\end{vmatrix}\\
&=&x\begin{vmatrix}
x&-1&0&0\\
-1&x&-1&0\\
0&-1&x&-1\\
0&0&-1&x
\end{vmatrix}
-\begin{vmatrix}
x&-1&0\\
-1&x&-1\\
0&-1&x
\end{vmatrix}\\
&=&xP_4-P_3
\end{eqnarray*}

(2) Regarding the intial values, these are normally $P_2=x^2-1$ and $P_3=x^3-2x$, but we can formally add the values $P_0=1$ and $P_1=x$, with this being justified by:
$$\ \ \ P_0=1\quad,\quad 
P_1=x$$
$$\implies P_2=xP_1-P_0=x^2-1$$
$$\implies P_3=xP_2-P_1=x^3-2x$$

(3) Thus, we have the recurrence formula in the statement, and with the initial values there, and an internet search, or some advanced calculus know-how, tells us that these are the well-known Chebycheff polynomials, enjoying lots of interesting properties.

\medskip

(4) Finally, in what regards the numbers of loops $L_k$, and the corresponding measures $\mu$, these can be in principle computed by using Theorem 7.6. However, the computations are quite tough here, and we will not choose this path. More on all this, later.

\medskip

(5) To be more precise, we will see, towards the end of the present book, that there is a more clever way of dealing with both the circle and segment graphs, substantially improving what we have in Proposition 7.8 and here. However, all this will be quite technical, not yet the moment for such things. More later, that is promised.
\end{proof}

\section*{7b. Catalan numbers}

Time perhaps for a conclusion? The loop counting business on graphs looks like a nice activity, but in practice, when trying this for the simplest possible graphs, we got into quite complicated mathematics. So, if there is a conclusion to this, we are in trouble.

\bigskip

Well, instead of giving up, let us look face-to-face at the difficulties that we met. We are led in this way, after analyzing the situation, to the following thought:

\begin{thought}
The difficulties that we met, with the circle and the two segments, come from the fact that our loops are not ``free to move'',
\begin{enumerate}
\item for the circle, because these can circle around the circle,

\item for the segments, obviously because of the endpoints,
\end{enumerate}
and so our difficulties will dissapear, and we will be able to do our exact loop count, once we find a graph $X$ where the loops are truly ``free to move''.
\end{thought}

Thinking some more, all this definitely buries the first interval graph, where the vertex $0$ is one of the endpoints. However, we can still try to recycle the circle, by unwrapping it, or extend our second interval graph up to $\infty$. But in both cases what we get is the graph $\mathbb Z$ formed by the integers. So, let us formulate the following definition:

\index{infinite graph}

\begin{definition}
An infinite graph is the same thing as a finite graph, but now with an infinity of vertices, $|X|=\infty$. As a basic example, we have $\mathbb Z$. We also have $\mathbb N$.
\end{definition}

Leaving aside now $\mathbb N$, which looks more complicated, let us try to count the length $k$ paths on $\mathbb Z$, based at $0$. At $k=1$ we have $2$ such paths, ending at $-1$ and $1$, and the count results can be pictured as follows, with everything being self-explanatory:

$$
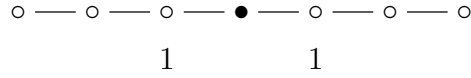
\xymatrix@R=5pt@C=15pt{
\circ\ar@{-}[r]&\circ\ar@{-}[r]&\circ\ar@{-}[r]&\bullet\ar@{-}[r]&\circ\ar@{-}[r]&\circ\ar@{-}[r]&\circ\\
&&1&&1
}$$

At $k=2$ now, we have 4 paths, one of which ends at $-2$, two of which end at 0, and one of which ends at 2. The results can be pictured as follows:
$$
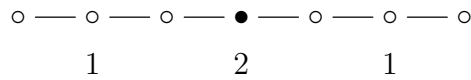
\xymatrix@R=5pt@C=15pt{
\circ\ar@{-}[r]&\circ\ar@{-}[r]&\circ\ar@{-}[r]&\bullet\ar@{-}[r]&\circ\ar@{-}[r]&\circ\ar@{-}[r]&\circ\\
&1&&2&&1
}$$

At $k=3$ now, we have 8 paths, the distribution of the endpoints being as follows:
$$
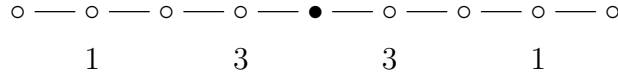
\xymatrix@R=5pt@C=15pt{
\circ\ar@{-}[r]&\circ\ar@{-}[r]&\circ\ar@{-}[r]&\circ\ar@{-}[r]&\bullet\ar@{-}[r]&\circ\ar@{-}[r]&\circ\ar@{-}[r]&\circ\ar@{-}[r]&\circ\\
&1&&3&&3&&1
}$$

As for $k=4$, here we have 16 paths, the distribution of the endpoints being as follows:
$$
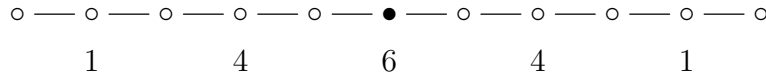
\xymatrix@R=5pt@C=15pt{
\circ\ar@{-}[r]&\circ\ar@{-}[r]&\circ\ar@{-}[r]&\circ\ar@{-}[r]&\circ\ar@{-}[r]&\bullet\ar@{-}[r]&\circ\ar@{-}[r]&\circ\ar@{-}[r]&\circ\ar@{-}[r]&\circ\ar@{-}[r]&\circ\\
&1&&4&&6&&4&&1
}$$

And good news, we can see in the above the Pascal triangle. Thus, getting back now to Question 7.2, we can answer it for the graph $\mathbb Z$, the result being as follows:

\index{central binomial coefficients}
\index{Pascal triangle}

\begin{theorem}
The paths on $\mathbb Z$ are counted by the binomial coefficients. In particular, the $2k$-paths based at $0$ are counted by the central binomial coefficients,
$$L_{2k}=\binom{2k}{k}$$
and $\mu$ is the centered measure having these numbers as even moments.
\end{theorem}

\begin{proof}
This basically follows from the above discussion, as follows:

\medskip

(1) In what regards the count, we certainly have the Pascal triangle, as discovered above, and the rest is just a matter of finishing. There are many possible ways here, a straightforward one being that of arguing that the number $C_k^l$ of length $k$ loops $0\to l$  is subject, due to the binary choice at the end, to the following recurrence relation:
$$C_k^l=C_{k-1}^{l-1}+C_{k-1}^{l+1}$$

But this is exactly the recurrence for the Pascal triangle, so done with the count. 

\medskip

(2) As for the second assertion, the first part, regarding $L_{2k}$, is clear from this, and the second part is more of an empty statement, with $\mu$ still remaining to be computed.
\end{proof}

As a second illustration, let us try now to count the loops of $\mathbb N$, based at 0. This is something less obvious, and at the experimental level, the result is as follows:

\index{Catalan numbers}

\begin{proposition}
The numbers $C_k$ counting the loops on $\mathbb N$ based at $0$,
$$C_k=\#\Big\{0-i_1-\ldots-i_{2k-1}-0\Big\}$$
are numerically $1,2,5,14,42,132,429,1430,4862,16796,58786,\ldots$
\end{proposition}

\begin{proof}
To start with, we have indeed $C_1=1$, the only loop here being $0-1-0$. Then we have $C_2=2$, due to two possible loops, namely:
$$0-1-0-1-0$$
$$0-1-2-1-0$$

Then we have $C_3=5$, the possible loops here being as follows:
$$0-1-0-1-0-1-0$$
$$0-1-0-1-2-1-0$$
$$0-1-2-1-0-1-0$$
$$0-1-2-1-2-1-0$$
$$0-1-2-3-2-1-0$$

And so on, and I will leave some more computations to you, as an exercise.
\end{proof}

Obviously, computing the numbers $C_k$ is no easy task, and finding the formula of $C_k$, out of the data that we have, does not look as an easy task either. So, let us look for other objects counted by the same numbers $C_k$. With a bit of luck, among these objects some will be easier to count than the others, and this will eventually compute $C_k$.

\bigskip

This was for the strategy. In practice now, we first have the following result:

\begin{theorem}
The following objects are counted by the same numbers $C_k$, called Catalan numbers, given by $C_0=C_1=1$ and $C_{k+1}=\sum_{a+b=k}C_aC_b$:
\begin{enumerate}
\item The length $2k$ loops on $\mathbb N$, based at $0$.

\item The noncrossing pairings of $1,\ldots,2k$.

\item The noncrossing partitions of $1,\ldots,k$.

\item The length $2k$ Dyck paths in the plane.
\end{enumerate}
\end{theorem}

\begin{proof}
All this is standard combinatorics, the idea being as follows:

\medskip

(1) To start with, in what regards the various objects involved, the length $2k$ loops on $\mathbb N$ are the usual length $2k$ loops on $\mathbb N$, and the same goes for the noncrossing pairings of $1,\ldots,2k$, and for the noncrossing partitions of $1,\ldots,k$, the idea here being that you must be able to draw the pairing or partition in a noncrossing way. 

\medskip

(2) Regarding now the length $2k$ Dyck paths in the plane, these are by definition the paths from $(0,0)$ to $(k,k)$, marching North-East over the integer lattice $\mathbb Z^2\subset\mathbb R^2$, by staying inside the square $[0,k]\times[0,k]$, and staying as well under the diagonal of this square. As an example, here are the 5 possible Dyck paths at $n=3$:
$$
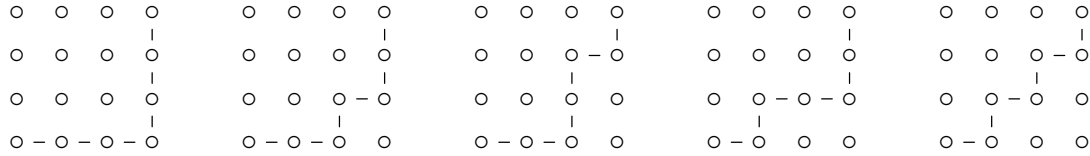
\xymatrix@R=4pt@C=4pt
{\circ&\circ&\circ&\circ\\
\circ&\circ&\circ&\circ\ar@{-}[u]\\
\circ&\circ&\circ&\circ\ar@{-}[u]\\
\circ\ar@{-}[r]&\circ\ar@{-}[r]&\circ\ar@{-}[r]&\circ\ar@{-}[u]}
\qquad
\xymatrix@R=4pt@C=4pt
{\circ&\circ&\circ&\circ\\
\circ&\circ&\circ&\circ\ar@{-}[u]\\
\circ&\circ&\circ\ar@{-}[r]&\circ\ar@{-}[u]\\
\circ\ar@{-}[r]&\circ\ar@{-}[r]&\circ\ar@{-}[u]&\circ}
\qquad
\xymatrix@R=4pt@C=4pt
{\circ&\circ&\circ&\circ\\
\circ&\circ&\circ\ar@{-}[r]&\circ\ar@{-}[u]\\
\circ&\circ&\circ\ar@{-}[u]&\circ\\
\circ\ar@{-}[r]&\circ\ar@{-}[r]&\circ\ar@{-}[u]&\circ}
\qquad
\xymatrix@R=4pt@C=4pt
{\circ&\circ&\circ&\circ\\
\circ&\circ&\circ&\circ\ar@{-}[u]\\
\circ&\circ\ar@{-}[r]&\circ\ar@{-}[r]&\circ\ar@{-}[u]\\
\circ\ar@{-}[r]&\circ\ar@{-}[u]&\circ&\circ}
\qquad
\xymatrix@R=4pt@C=4pt
{\circ&\circ&\circ&\circ\\
\circ&\circ&\circ\ar@{-}[r]&\circ\ar@{-}[u]\\
\circ&\circ\ar@{-}[r]&\circ\ar@{-}[u]&\circ\\
\circ\ar@{-}[r]&\circ\ar@{-}[u]&\circ&\circ}
$$

(3) Thus, we have definitions for all objects involved, and in each case, if you start counting them, you end up with the same sequence, namely that in Proposition 7.13:
$$1,1,2,5,14,42,132,429,1430,4862,\ldots$$

(4) In order to prove now that (1-4) produce indeed the same numbers, many things can be said. The idea is that, leaving aside mathematical brevity, and more specifically abstract reasonings of type $a=b,b=c\implies a=c$, what we have to do, in order to fully understand what is going on, is to etablish $\binom{4}{2}=6$ equalities, via bijective proofs.

\medskip

(5) But this can be done, indeed. As an example here, the noncrossing pairings of $1,\ldots,2k$ from (2) are in bijection with the noncrossing partitions of $1,\ldots,k$ from (3), via  fattening the pairings and shrinking the partitions. Consider a noncrossing partition:
$$
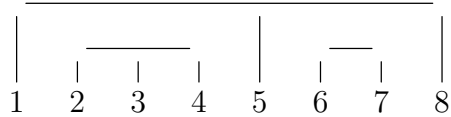
\xymatrix@R=10pt@C=10pt{
\ar@{-}[rrrrrrr]&&&&&&&\\
&\ar@{-}[rr]&&&&\ar@{-}[r]&&\\
1\ar@{-}[uu]&2\ar@{-}[u]&3\ar@{-}[u]&4\ar@{-}[u]&5\ar@{-}[uu]&6\ar@{-}[u]&7\ar@{-}[u]&8\ar@{-}[uu]
}$$

Now let us ``fatten'' this partition, by doubling everything, as follows:
$$
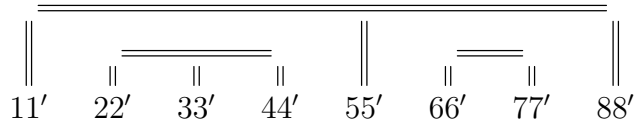
\xymatrix@R=10pt@C=10pt{
\ar@{=}[rrrrrrr]&&&&&&&\\
&\ar@{=}[rr]&&&&\ar@{=}[r]&&\\
11'\ar@{=}[uu]&22'\ar@{=}[u]&33'\ar@{=}[u]&44'\ar@{=}[u]&55'\ar@{=}[uu]&66'\ar@{=}[u]&77'\ar@{=}[u]&88'\ar@{=}[uu]
}$$

We can see emerging here a noncrossing pairing, and by relabeling the points $1,\ldots,16$, and properly redrawing the picture, what we have is indeed a noncrossing pairing:
$$
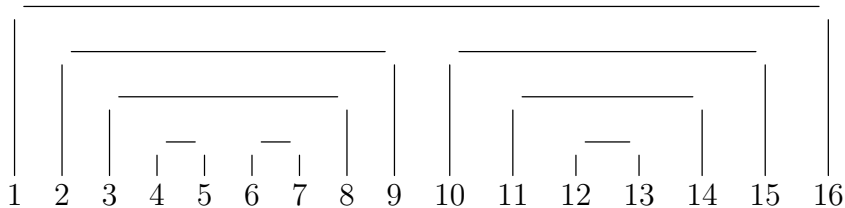
\xymatrix@R=10pt@C=5pt{
\ar@{-}[rrrrrrrrrrrrrrr]&&&&&&&&&&&&&&&\\
&\ar@{-}[rrrrrrr]&&&&&&&&\ar@{-}[rrrrr]&&&&&&\\
&&\ar@{-}[rrrrr]&&&&&&&&\ar@{-}[rrr]&&&&&\\
&&&\ar@{-}[r]&&\ar@{-}[r]&&&&&&\ar@{-}[r]&&&&\\
1\ar@{-}[uuuu]&2\ar@{-}[uuu]&3\ar@{-}[uu]&4\ar@{-}[u]&5\ar@{-}[u]&6\ar@{-}[u]&7\ar@{-}[u]&8\ar@{-}[uu]&9\ar@{-}[uuu]&10\ar@{-}[uuu]&11\ar@{-}[uu]&12\ar@{-}[u]&13\ar@{-}[u]&14\ar@{-}[uu]&15\ar@{-}[uuu]&16\ar@{-}[uuuu]
}$$

As for the reverse operation, that is obviously obtained by ``shrinking'' our pairing, by collapsing pairs of consecutive neighbors, that is, by identifying $1=2$, then $3=4$, then $5=6$, and so on. We will leave the details here as an instructive exercise, and exercise as well, to add the loops in (1) and the Dyck paths in (4) to the bijective picture.

\medskip

(6) This being said, as a matter of having our theorem formally proved, I mean by me professor and not by you student, here is a less elegant argument, which is however quick, and does the job. The point is that, in each of the cases (1-4) under consideration, the numbers $C_k$ that we get are easily seen to be subject to the following recurrence:
$$C_{k+1}=\sum_{a+b=k}C_aC_b$$ 

Now the initial data being the same, namely $C_0=C_1=1$, in each of the cases (1-4) under consideration, we get indeed the same numbers, as desired.
\end{proof}

We can pass now to the second step in our plan, namely selecting in the above list the objects that we find the most convenient to count, and count them. This leads to:

\index{Catalan numbers}

\begin{theorem}
The Catalan numbers are given by the formula
$$C_k=\frac{1}{k+1}\binom{2k}{k}$$
with this being best seen by counting the length $2k$ Dyck paths in the plane.
\end{theorem}

\begin{proof}
This is something quite tricky, the idea being as follows:

\medskip

(1) Let us count indeed the Dyck paths in the plane. For this purpose, we use a trick. Indeed, if we ignore the assumption that our path must stay under the diagonal of the square, we have $\binom{2k}{k}$ such paths. And among these, we have the ``good'' ones, those that we want to count, and then the ``bad'' ones, those that we want to ignore.

\medskip

(2) So, let us count the bad paths, those crossing the diagonal of the square, and reaching the higher diagonal next to it, the one joining $(0,1)$ and $(k,k+1)$. In order to count these, the trick is to ``flip'' their bad part over that higher diagonal, as follows:
$$
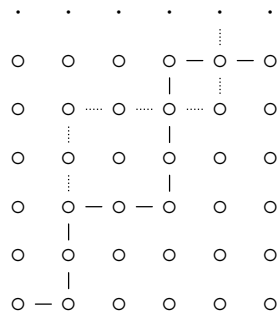
\xymatrix@R=6pt@C=6pt
{\cdot&\cdot&\cdot&\cdot&\cdot&\cdot\\
\circ&\circ&\circ&\circ\ar@{-}[r]&\circ\ar@{-}[r]\ar@{.}[u]&\circ\\
\circ&\circ\ar@{.}[r]&\circ\ar@{.}[r]&\circ\ar@{.}[r]\ar@{-}[u]&\circ\ar@{.}[u]&\circ\\
\circ&\circ\ar@{.}[u]&\circ&\circ\ar@{-}[u]&\circ&\circ\\
\circ&\circ\ar@{-}[r]\ar@{.}[u]&\circ\ar@{-}[r]&\circ\ar@{-}[u]&\circ&\circ\\
\circ&\circ\ar@{-}[u]&\circ&\circ&\circ&\circ\\
\circ\ar@{-}[r]&\circ\ar@{-}[u]&\circ&\circ&\circ&\circ}$$

(3) Now observe that, as it is obvious on the above picture, due to the flipping, the flipped bad path will no longer end in $(k,k)$, but rather in $(k-1,k+1)$. Moreover, more is true, in the sense that, by thinking a bit, we see that the flipped bad paths are precisely those ending in $(k-1,k+1)$. Thus, we can count these flipped bad paths, and so the bad paths, and so the good paths too, and so good news, we are done.

\medskip

(4) To be more precise, by putting everything together, we obtain:
\begin{eqnarray*}
C_k
&=&\binom{2k}{k}-\binom{2k}{k-1}\\
&=&\binom{2k}{k}-\frac{k}{k+1}\binom{2k}{k}\\
&=&\frac{1}{k+1}\binom{2k}{k}
\end{eqnarray*}

Thus, we are led to the formula in the statement.
\end{proof}

Finally, still talking Catalan numbers, all the above was sweet combinatorics, that must be learned, but we can in fact kill as well the problem with calculus, as follows: 

\index{Catalan numbers}

\begin{theorem}
The Catalan numbers have the following properties:
\begin{enumerate}
\item They satisfy $C_{k+1}=\sum_{a+b=k}C_aC_b$.

\item The series $f(z)=\sum_{k\geq0}C_kz^k$ satisfies $zf^2-f+1=0$.

\item This series is given by $f(z)=\frac{1-\sqrt{1-4z}}{2z}$.

\item We have the formula $C_k=\frac{1}{k+1}\binom{2k}{k}$.
\end{enumerate}
\end{theorem}

\begin{proof}
This is best viewed by using noncrossing pairings, as follows: 

\medskip

(1) Let us count the noncrossing pairings of $\{1,\ldots,2k+2\}$. Such a pairing appears by pairing 1 to an odd number, $2a+1$, and then inserting a noncrossing pairing of $\{2,\ldots,2a\}$, and a noncrossing pairing of $\{2a+2,\ldots,2k+2\}$. Thus we have, as claimed:
$$C_{k+1}=\sum_{a+b=k}C_aC_b$$ 

(2) Consider now the generating series of the Catalan numbers, $f(z)=\sum_{k\geq0}C_kz^k$. In terms of this generating series, the above recurrence gives, as desired:
\begin{eqnarray*}
zf^2
&=&\sum_{a,b\geq0}C_aC_bz^{a+b+1}\\
&=&\sum_{k\geq1}\sum_{a+b=k-1}C_aC_bz^k\\
&=&\sum_{k\geq1}C_kz^k\\
&=&f-1
\end{eqnarray*}

(3) By solving now the equation $zf^2-f+1=0$ found above, and choosing the solution which is bounded at $z=0$, we obtain the following formula, as claimed:
$$f(z)=\frac{1-\sqrt{1-4z}}{2z}$$ 

(4) In order to compute this function, we use the generalized binomial formula, which is as follows, with $p\in\mathbb R$ being an arbitrary exponent, and with $|t|<1$:
$$(1+t)^p=\sum_{k=0}^\infty\binom{p}{k}t^k$$

Indeed, for the exponent $p=1/2$, the generalized binomial coefficients are:
\begin{eqnarray*}
\binom{1/2}{k}
&=&\frac{1/2(-1/2)(-3/2)\ldots(3/2-k)}{k!}\\
&=&(-1)^{k-1}\frac{1\cdot 3\cdot 5\ldots(2k-3)}{2^kk!}\\
&=&(-1)^{k-1}\frac{(2k-2)!}{2^{k-1}(k-1)!2^kk!}\\
&=&\frac{(-1)^{k-1}}{2^{2k-1}}\cdot\frac{1}{k}\binom{2k-2}{k-1}\\
&=&-2\left(\frac{-1}{4}\right)^k\cdot\frac{1}{k}\binom{2k-2}{k-1}
\end{eqnarray*}

Thus the generalized binomial formula at exponent $p=1/2$ reads:
$$\sqrt{1+t}=1-2\sum_{k=1}^\infty\frac{1}{k}\binom{2k-2}{k-1}\left(\frac{-t}{4}\right)^k$$

With $t=-4z$ we obtain from this the following formula:
$$\sqrt{1-4z}=1-2\sum_{k=1}^\infty\frac{1}{k}\binom{2k-2}{k-1}z^k$$

Now back to our series $f$, we obtain the following formula for it:
\begin{eqnarray*}
f(z)
&=&\frac{1-\sqrt{1-4z}}{2z}\\
&=&\sum_{k=1}^\infty\frac{1}{k}\binom{2k-2}{k-1}z^{k-1}\\
&=&\sum_{k=0}^\infty\frac{1}{k+1}\binom{2k}{k}z^k
\end{eqnarray*}

Thus, the Catalan numbers are indeed given by $C_k=\frac{1}{k+1}\binom{2k}{k}$, as claimed.
\end{proof}

Quite interesting all this, and as a byproduct of our study, let us record the following result, allowing us to extract square roots, which is something very useful:

\begin{theorem}
We have the following formula,
$$\sqrt{1+t}=1-2\sum_{k=1}^\infty C_{k-1}\left(\frac{-t}{4}\right)^k$$
with $C_k=\frac{1}{k+1}\binom{2k}{k}$ being the Catalan numbers. Also, we have
$$\frac{1}{\sqrt{1+t}}=\sum_{k=0}^\infty D_k\left(\frac{-t}{4}\right)^k$$
with $D_k=\binom{2k}{k}$ being the central binomial coefficients.
\end{theorem}

\begin{proof}
Here the first formula comes from our considerations above, namely generalized binomial formula with exponent $p=1/2$. As for the second formula, this comes in a similar way, from the generalized binomial formula with exponent $p=-1/2$.
\end{proof}

Getting back now to graphs, we can complement Theorem 7.12 with:

\begin{theorem}
The loops on $\mathbb N$ are counted by the Catalan numbers,
$$L_{2k}=\frac{1}{k+1}\binom{2k}{k}$$
and $\mu$ is the centered measure having these numbers as even moments.
\end{theorem}

\begin{proof}
This follows indeed from what we have, namely Theorem 7.14 for the recurrence, and then Theorem 7.15 or Theorem 7.16 for the solution of the recurrence.
\end{proof}

Summarizing, loop count question solved for both the graphs $\mathbb Z$ and $\mathbb N$. Which is very nice, and in what follows we will keep working on these two graphs, $\mathbb Z$ and $\mathbb N$.

\section*{7c. Semicircle laws}

According to what we have in Theorems 7.12 and 7.18, the problem is now, how to recover a real probability measure out of its moments. And things here are quite subtle, with the main result, which is something quite technical, being as follows:

\index{Stieltjes inversion}
\index{Cauchy transform}
\index{moment problem}

\begin{theorem}
The density of a continuous probability measure $\mu$ can be recaptured from the sequence of moments $\{M_k\}_{k\geq0}$ via the Stieltjes inversion formula
$$d\mu (x)=\lim_{t\searrow 0}-\frac{1}{\pi}\,Im\left(G(x+it)\right)\cdot dx$$
where the function on the right, given in terms of moments by
$$G(\xi)=\xi^{-1}+M_1\xi^{-2}+M_2\xi^{-3}+\ldots$$
is the Cauchy transform of the measure $\mu$.
\end{theorem}

\begin{proof}
To start with, the Cauchy transform of our measure $\mu$ is given by:
\begin{eqnarray*}
G(\xi)
&=&\xi^{-1}\sum_{k=0}^\infty M_k\xi^{-k}\\\
&=&\int_\mathbb R\frac{\xi^{-1}}{1-\xi^{-1}y}\,d\mu(y)\\
&=&\int_\mathbb R\frac{1}{\xi-y}\,d\mu(y)
\end{eqnarray*}

Now with $\xi=x+it$, we obtain from this the following formula:
\begin{eqnarray*}
Im(G(x+it))
&=&\int_\mathbb RIm\left(\frac{1}{x-y+it}\right)d\mu(y)\\
&=&\int_\mathbb R\frac{1}{2i}\left(\frac{1}{x-y+it}-\frac{1}{x-y-it}\right)d\mu(y)\\
&=&-\int_\mathbb R\frac{t}{(x-y)^2+t^2}\,d\mu(y)
\end{eqnarray*}

By integrating over $[a,b]$ we obtain, with the change of variables $x=y+tz$:
\begin{eqnarray*}
\int_a^bIm(G(x+it))dx
&=&-\int_\mathbb R\int_a^b\frac{t}{(x-y)^2+t^2}\,dx\,d\mu(y)\\
&=&-\int_\mathbb R\int_{(a-y)/t}^{(b-y)/t}\frac{t}{(tz)^2+t^2}\,t\,dz\,d\mu(y)\\
&=&-\int_\mathbb R\int_{(a-y)/t}^{(b-y)/t}\frac{1}{1+z^2}\,dz\,d\mu(y)\\
&=&-\int_\mathbb R\left(\arctan\frac{b-y}{t}-\arctan\frac{a-y}{t}\right)d\mu(y)
\end{eqnarray*}

Now observe that with $t\searrow0$ we have the following formulae:
$$\lim_{t\searrow0}\left(\arctan\frac{b-y}{t}-\arctan\frac{a-y}{t}\right)
=\begin{cases}
\frac{\pi}{2}-\frac{\pi}{2}=0& (y<a)\\
\frac{\pi}{2}-0=\frac{\pi}{2}& (y=a)\\
\frac{\pi}{2}-(-\frac{\pi}{2})=\pi& (a<y<b)\\
0-(-\frac{\pi}{2})=\frac{\pi}{2}& (y=b)\\
-\frac{\pi}{2}-(-\frac{\pi}{2})=0& (y>b)
\end{cases}$$

We therefore obtain the following formula:
$$\lim_{t\searrow0}\int_a^bIm(G(x+it))dx=-\pi\left(\mu(a,b)+\frac{\mu(a)+\mu(b)}{2}\right)$$

Thus, we are led to the conclusion in the statement.
\end{proof}

As a first comment on this, with some suitable modifications, the Stieltjes inversion formula applies as well to measures featuring Dirac masses. However, for our purposes here, where most measures will be continuous, what we have in Theorem 7.19 will do.

\bigskip

As a second comment, Theorem 7.19 and its generalizations do not fully solve the moment problem, because we still have the question of understanding when a sequence of numbers $M_1,M_2,M_3,\ldots$ can be the moments of a measure $\mu$. And here, we have:

\index{Hankel determinant}

\begin{theorem}
A sequence of numbers $M_0,M_1,M_2,M_3,\ldots\in\mathbb R$, with $M_0=1$, is the sequence of moments of a real probability measure $\mu$ precisely when:
$$\begin{vmatrix}M_0\end{vmatrix}\geq0\quad,\quad 
\begin{vmatrix}
M_0&M_1\\
M_1&M_2
\end{vmatrix}\geq0\quad,\quad 
\begin{vmatrix}
M_0&M_1&M_2\\
M_1&M_2&M_3\\
M_2&M_3&M_4\\
\end{vmatrix}\geq0\quad,\quad 
\ldots$$
That is, the associated Hankel determinants must be all positive.
\end{theorem}

\begin{proof}
As a first observation, by linear algebra, the positivity conditions in the statement tell us that the following associated linear forms must be positive:
$$\sum_{i,j=1}^nc_i\bar{c}_jM_{i+j}\geq0$$

But with this observation in hand, in one sense the result is elementary, coming from the following computation, which shows that we have positivity indeed:
\begin{eqnarray*}
\int_\mathbb R\left|\sum_{i=1}^nc_ix^i\right|^2d\mu(x)
&=&\int_\mathbb R\sum_{i,j=1}^nc_i\bar{c}_jx^{i+j}d\mu(x)\\
&=&\sum_{i,j=1}^nc_i\bar{c}_jM_{i+j}
\end{eqnarray*}

As for the other sense, here the result comes once again from the above formula, this time via some standard functional analysis. So, this is the situation, welcome to the moment problem, and of course exercise for you, to learn a bit more about all this.
\end{proof}

Getting back now to the Stieltjes formula from Theorem 7.19, let us try to solve the moment problem for the Catalan numbers $C_k$, and for the central binomial coefficients $D_k$. We first have the following result, coming as a complement to Theorem 7.18:

\index{semicircle law}
\index{Wigner law}

\begin{theorem}
The real measure having as even moments the Catalan numbers, $C_k=\frac{1}{k+1}\binom{2k}{k}$, and having all odd moments $0$ is the measure
$$\gamma_1=\frac{1}{2\pi}\sqrt{4-x^2}dx$$
called Wigner semicircle law on $[-2,2]$.
\end{theorem}

\begin{proof}
In order to apply the inversion formula, our starting point will be the formula from Theorem 7.16 for the generating series of the Catalan numbers, namely:
$$\sum_{k=0}^\infty C_kz^k=\frac{1-\sqrt{1-4z}}{2z}$$

By using this formula with $z=\xi^{-2}$, we obtain the following formula:
\begin{eqnarray*}
G(\xi)
&=&\xi^{-1}\sum_{k=0}^\infty C_k\xi^{-2k}\\
&=&\xi^{-1}\cdot\frac{1-\sqrt{1-4\xi^{-2}}}{2\xi^{-2}}\\
&=&\frac{\xi}{2}\left(1-\sqrt{1-4\xi^{-2}}\right)\\
&=&\frac{\xi}{2}-\frac{1}{2}\sqrt{\xi^2-4}
\end{eqnarray*}

Now let us apply Theorem 7.19. The study here goes as follows:

\medskip

(1) According to the general philosophy of the Stieltjes formula, the first term, namely $\xi/2$, which is ``trivial'', will not contribute to the density. 

\medskip

(2) As for the second term, which is something non-trivial, this will contribute to the density, the rule here being that the square root $\sqrt{\xi^2-4}$ will be replaced by the ``dual'' square root $\sqrt{4-x^2}\,dx$, and that we have to multiply everything by $-1/\pi$. 

\medskip

(3) As a conclusion, by Stieltjes inversion we obtain the following density:
$$d\mu(x)
=-\frac{1}{\pi}\cdot-\frac{1}{2}\sqrt{4-x^2}\,dx
=\frac{1}{2\pi}\sqrt{4-x^2}dx$$

Thus, we have obtained the mesure in the statement, and we are done.
\end{proof}

Next, we have the following interesting version of the above result, whose precise interest will become clear later in this chapter, and later in this book as well:

\index{Marchenko-Pastur law}

\begin{theorem}
The real measure having as sequence of moments the Catalan numbers, $C_k=\frac{1}{k+1}\binom{2k}{k}$, is the measure
$$\pi_1=\frac{1}{2\pi}\sqrt{4x^{-1}-1}\,dx$$
called Marchenko-Pastur law on $[0,4]$.
\end{theorem}

\begin{proof}
As before, we use the standard formula for the generating series of the Catalan numbers. With $z=\xi^{-1}$ in that formula, we obtain the following formula:
\begin{eqnarray*}
G(\xi)
&=&\xi^{-1}\sum_{k=0}^\infty C_k\xi^{-k}\\
&=&\xi^{-1}\cdot\frac{1-\sqrt{1-4\xi^{-1}}}{2\xi^{-1}}\\
&=&\frac{1}{2}\left(1-\sqrt{1-4\xi^{-1}}\right)\\
&=&\frac{1}{2}-\frac{1}{2}\sqrt{1-4\xi^{-1}}
\end{eqnarray*}

With this in hand, let us apply now the Stieltjes inversion formula, from Theorem 7.19. We obtain, a bit as before in Theorem 7.21, the following density:
$$d\mu(x)
=-\frac{1}{\pi}\cdot-\frac{1}{2}\sqrt{4x^{-1}-1}\,dx
=\frac{1}{2\pi}\sqrt{4x^{-1}-1}\,dx$$

Thus, we are led to the conclusion in the statement.
\end{proof}

Getting now to the central binomial coefficients, we have two computations to do for them, in the spirit of those done in Theorems 7.21 and 7.22. The result is:

\index{arcsine law}

\begin{theorem}
The real probability measures having as even moments/moments the central binomial coefficients, $D_k=\binom{2k}{k}$, are the measures
$$\alpha_1=\frac{1}{\pi\sqrt{4-x^2}}\,dx\quad,\quad\alpha_1'=\frac{1}{\pi\sqrt{x(4-x)}}\,dx$$
called arcsine laws on $[-2,2]$, and on $[0,4]$.
\end{theorem}

\begin{proof}
For the first measure, by using Theorem 7.17, the Cauchy transform is:
$$G(\xi)
=\xi^{-1}\sum_{k=0}^\infty D_k\xi^{-2k}
=\frac{1}{\xi}\cdot\frac{1}{\sqrt{1-4/\xi^2}}
=\frac{1}{\sqrt{\xi^2-4}}$$

As for the second measure, again by using Theorem 7.17, the Cauchy transform is:
$$G(\xi)
=\xi^{-1}\sum_{k=0}^\infty D_k\xi^{-k}
=\frac{1}{\xi}\cdot\frac{1}{\sqrt{1-4/\xi}}
=\frac{1}{\sqrt{\xi(\xi-4)}}$$

But this gives the densities in the statement, via Theorem 7.19. 
\end{proof}

Finally, we have the following useful technical version of the above result:

\index{modified arcsine law}
\index{middle binomial coefficients}

\begin{theorem}
The real probability measure having as moments the middle binomial coefficients, $E_k=\binom{k}{[k/2]}$, is the following law on $[-2,2]$,
$$\sigma_1=\frac{1}{2\pi}\sqrt{\frac{2+x}{2-x}}\,dx$$
called modified arcsine law on $[-2,2]$.
\end{theorem}

\begin{proof}
In terms of the central binomial coefficients $D_k$, we have:
$$E_{2k}=\binom{2k}{k}=\frac{(2k)!}{k!k!}=D_k$$
$$E_{2k-1}=\binom{2k-1}{k}=\frac{(2k-1)!}{k!(k-1)!}=\frac{D_k}{2}$$

Standard calculus based on the Taylor formula for $(1+t)^{-1/2}$ gives:
$$\frac{1}{2x}\left(\sqrt{\frac{1+2x}{1-2x}}-1\right)=\sum_{k=0}^\infty E_kx^k$$

With $x=\xi^{-1}$ we obtain the following formula for the Cauchy transform:
\begin{eqnarray*}
G(\xi)
&=&\xi^{-1}\sum_{k=0}^\infty E_k\xi^{-k}\\
&=&\frac{1}{\xi}\left(\sqrt{\frac{1+2/\xi}{1-2/\xi}}-1\right)\\
&=&\frac{1}{\xi}\left(\sqrt{\frac{\xi+2}{\xi-2}}-1\right)
\end{eqnarray*}

By Stieltjes inversion we obtain the density in the statement.
\end{proof}

And with this, end of our discussion regarding the random walks on basic infinite graphs, and related topics. Good learning this was, and more on such things, later.

\section*{7d. Parametric versions}

Getting back now to the plan made in the beginning of this chapter, we have certainly understood quite well the phenomenology of the Wigner law, at least from the random walk viewpoint, and as a bonus, we came across a number of other interesting measures as well, namely Marchenko-Pastur, and the arcsine and modified arcsine laws.

\bigskip

What is next? Second part of that plan, namely doing the mathematics for the Wigner law. And here, it makes sense to include as well some of the new measures that we found. More precisely, forgetting about the arcsine and modified arcsine laws, which are something quite secondary, related to the central binomial coefficients $D_k=\binom{2k}{k}$, we would like to include in our study the Marchenko-Pastur law, which is crucially related to the Catalan numbers $C_k=\frac{1}{k+1}\binom{2k}{k}$, as the Wigner semicircle law itself is.

\bigskip

In practice now, here is the question that we would like to solve:

\begin{question}
Do the Wigner and Marchenko-Pastur laws, namely
$$\gamma_1=\frac{1}{2\pi}\sqrt{4-x^2}dx\quad,\quad 
\pi_1=\frac{1}{2\pi}\sqrt{4x^{-1}-1}\,dx$$
have interesting parametric versions $\gamma_t,\pi_t$, and what is the mathematics of these?
\end{question}

In answer now, let us start with the Wigner law. That is by definition the semicircle law on $[-2,2]$, and it is a no-brainer to look, more generally, at the semicircle laws on arbitrary intervals $[-a,a]$. We are led in this way to the following result:

\begin{theorem}
Given $t>0$, the real measure having as even moments the numbers $M_{2k}=t^kC_k$ and having all odd moments $0$ is the measure
$$\gamma_t=\frac{1}{2\pi t}\sqrt{4t-x^2}dx$$
called Wigner semicircle law on $[-2\sqrt{t},2\sqrt{t}]$.
\end{theorem}

\begin{proof}
This looks quite straightforward, out of what we already have at $t=1$, save perhaps for our choice of the parametrization, favoring the moments over the support, which certainly needs some discussion. The idea with all this is as follows:

\medskip

(1) A first approach is by redoing our previous Stieltjes inversion computation, with a parameter $t>0$ added. To be more precise, as before for $\gamma_1$, the starting point is the formula in Theorem 7.16 for the generating series of the Catalan numbers, namely:
$$\sum_{k=0}^\infty C_kz^k=\frac{1-\sqrt{1-4z}}{2z}$$

Indeed, by using this formula with $z=t\xi^{-2}$, we obtain the following formula, for the Cauchy transform of the real measure that we want to compute:
\begin{eqnarray*}
G(\xi)
&=&\xi^{-1}\sum_{k=0}^\infty t^kC_k\xi^{-2k}\\
&=&\xi^{-1}\cdot\frac{1-\sqrt{1-4t\xi^{-2}}}{2t\xi^{-2}}\\
&=&\frac{\xi}{2t}\left(1-\sqrt{1-4t\xi^{-2}}\right)\\
&=&\frac{\xi}{2t}-\frac{1}{2t}\sqrt{\xi^2-4t}
\end{eqnarray*}

Thus, by Stieltjes inversion we obtain the following density, as claimed:
$$d\mu(x)=\frac{1}{2\pi t}\sqrt{4t-x^2}\,dx$$

(2) We have as well a second approach, based on what we have at $t=1$, and a change of variables. Indeed, the even moments of the measure in the statement are given by:
\begin{eqnarray*}
M_{2k}
&=&\frac{1}{2\pi t}\int_{-2\sqrt{t}}^{2\sqrt{t}}\sqrt{4t-x^2}\,x^{2k}\,dx\\
&=&\frac{1}{2\pi t}\int_{-1}^1\sqrt{4t-ty^2}\,(\sqrt{t}y)^{2k}\,\sqrt{t}\,dy\\
&=&\frac{t^k}{2\pi}\int_{-1}^1\sqrt{4-y^2}\,y^{2k}\,dy\\
&=&t^kC_k
\end{eqnarray*}

As for the odd moments, these all vanish, the density of $\gamma_t$ being an even function.

\medskip

(3) We have in fact a third approach as well. Indeed, if we want to get everything via standard calculus only, we first have to go back to $\gamma_1$, and compute its moments. But this can be done with the change of variable $x=2\cos t$, and we obtain, with some trigonometric know-how, namely Wallis formula, that we talk about later in this book:
\begin{eqnarray*}
M_{2k}
&=&\frac{1}{\pi}\int_0^2\sqrt{4-x^2}x^{2k}\,dx\\
&=&\frac{1}{\pi}\int_0^{\pi/2}\sqrt{4-4\cos^2t}\,(2\cos t)^{2k}2\sin t\,dt\\
&=&\frac{4^{k+1}}{\pi}\int_0^{\pi/2}\cos^{2k}t\sin^2t\,dt\\
&=&\frac{4^{k+1}}{\pi}\cdot\frac{\pi}{2}\cdot\frac{(2k)!!2!!}{(2k+3)!!}\\
&=&2\cdot 4^k\cdot\frac{(2k)!/2^kk!}{2^{k+1}(k+1)!}\\
&=&C_k
\end{eqnarray*}

As for the odd moments, these all vanish, the density of $\gamma_1$ being an even function. But with this done, we can then upgrade to $t>0$, again via calculus only, as in (2).

\medskip

(4) Summarizing, one way or another, we are led to the conclusion in the statement. This being said, things are not over here, because we still have to talk about our parametrization, in terms of $t>0$ as above, favoring the moments over the support.

\medskip

(5) And here, I am afraid that I have no simple explanation, and that we will need some more phenomenology. As a first good explanation for our parametrization, coming from the original work of Wigner himself, the point is that, in the random matrix context, no one cares about the support, and what comes out from computations is the moment formula $M_{2k}=t^kC_k$, leading to $\gamma_t$ as above. And more on this, in chapter 16.

\medskip

(6) As a second explanation, which is something advanced too, but is a bit more conceptual, recall from chapter 4 that our main parametric laws so far, namely the Poisson ones, were subject to the following formula, which obviously dictates the scaling:
$$p_s*p_t=p_{s+t}$$

 In what regards now the Wigner laws, these certainly cannot be subject to such a convolution formula, and this no matter the parametrization, as shown for instance by some quick Fourier transform computations, that we will leave as an exercise. However, and here comes the point, the convolution operation $*$ has a certain free analogue $\boxplus$, due to Voiculescu, and we have the following formula, which dictates the scaling:
$$\gamma_s\boxplus\gamma_t=\gamma_{s+t}$$

And more on such things, which are actually related to the random matrix considerations in (5) too, via a certain key phenomenon of asymptotic freeness for the Wigner matrices, and other random matrices, towards the end of the present book.
\end{proof}

Getting now to the mathematics of the Wigner laws, we already have the moments, our laws being centered the central moments are the same, and finally the Fourier transform computation is irrelevant, these laws being of ``free'' nature, as just explained above. So, good news, done with all work already, and let us just record the following fact:

\begin{theorem}
For the Wigner law $\gamma_t$ of parameter $t>0$,
$$E=0\quad,\quad V=t\quad,\quad \gamma=0\quad,\quad \kappa=2$$
are the expectation, variance, skewness and kurtosis.
\end{theorem}

\begin{proof}
We know from Theorem 7.26 that the odd moments vanish, and that the even moments are given by $M_{2k}=t^kC_k$. Thus the expectation is $E=0$, and the variance is $V=M_2=t$. Now by using this, we conclude that the odd normalized central moments vanish, and that the even normalized central moments are given by:
$$M_{2k}''=\frac{M_{2k}}{V^k}=\frac{t^kC_k}{t^k}=C_k$$

Thus the skewness is $\gamma=M_3''=0$, and the kurtosis is $\kappa=M_4''=C_2=2$, as stated.
\end{proof}

Regarding now the Marchenko-Pastur law, the parametric problematics here is substantially more tricky, leading among others to an atom at $0$, as follows:

\begin{theorem}
The Marchenko-Pastur law of parameter $t>0$ is given by
$$\pi_t=\max(1-t,0)\delta_0+\frac{\sqrt{4t-(x-1-t)^2}}{2\pi x}\,dx$$
the support being $[(1-\sqrt{t})^2,(1+\sqrt{t})^2]$, 
and the moments of this law are
$$M_k=\sum_{\pi\in NC(k)}t^{|\pi|}$$
$NC(k)$ being the noncrossing partitions of $\{1,\ldots,k\}$, and $|.|$, the number of blocks.
\end{theorem}

\begin{proof}
Obviously, some magic trickery going on here, the idea being as follows:

\medskip

(1) To start with, the ``classical'' phenomenology of $\pi_1$, including the random walks on $\mathbb N$, discussed above, and other things that we have not discussed, such as a relation of type $1+\chi\sim\pi_1$ with the group $SO_3$, in the spirit of the general theory from chapter 4, do not provide a clue on what the parametric version $\pi_t$ with $t>0$ should be.

\medskip

(2) Thus, we must go into the ``free'' phenomenology of $\pi_1$, in the spirit of what was said about the Wigner laws, in the proof of Theorem 7.26. And here, two things happen. First, and following the original approach of Marchenko and Pastur themselves, what comes out of random matrix computations is the moment formula $M_k=\sum_{\pi\in NC(k)}t^{|\pi|}$, which itself leads, via some standard calculus, to the law $\pi_t$ in the statement.

\medskip

(3) Alternatively, following Voiculescu, the PLT leading to the Poisson laws $p_t$, which satisfy $p_s*p_t=p_{s+t}$, has a free analogue, the FPLT, leading to the above laws $\pi_t$, which are also called free Poisson laws, and which satisfy $\pi_s\boxplus\pi_t=\pi_{s+t}$. So, this was for the story, and more on all this, random matrices and freeness, later in this book.

\medskip

(4) Finally, let us verify that what we have works well at $t=1$. In what regards the density, at $t=1$ we have no atom, and we obtain the formula from Theorem 7.22:
$$\pi_1=\frac{\sqrt{4-(x-2)^2}}{2\pi x}\,dx
=\frac{\sqrt{4x-x^2}}{2\pi x}\,dx
=\frac{1}{2\pi}\sqrt{4x^{-1}-1}\,dx$$

As for the moments, at $t=1$ we obtain $M_k=|NC(k)|=C_k$, as needed.
\end{proof}

Getting now to the mathematics of the Marchenko-Pastur laws, we already have the moments, and the Fourier transform computation is irrelevant, these laws being of ``free'' nature, as explained above. So, we are left with studying the central moments, normalized or not, and then computing the skewness and kurtosis, and we have here:

\begin{theorem}
For the Marchenko-Pastur law $\pi_t$ of parameter $t>0$,
$$E=t\quad,\quad V=t\quad,\quad \gamma=\frac{1}{\sqrt{t}}\quad,\quad \kappa=2+\frac{1}{t}$$
are the expectation, variance, skewness and kurtosis.
\end{theorem}

\begin{proof}
This follows from the moment formula in Theorem 7.28, but the best is to do this with a trick, by recycling our previous computations for the Poisson laws, from chapter 4. Indeed, let us put the moment formulae for $p_t$ and $\pi_t$ side by side:
$$M_k(p_t)=\sum_{\pi\in P(k)}t^{|\pi|}\quad,\quad 
M_k(\pi_t)=\sum_{\pi\in NC(k)}t^{|\pi|}$$

Now since $P(k)=NC(k)$ at $k=1,2,3$, the first three moments of $p_t,\pi_t$ coincide. As for the 4th moment, since $P(4)=NC(4)\cup\{\cap\!\!\cap\}$, the passage $p_t\to\pi_t$ amounts in substracting a $t^2$ term. We are therefore led to the following formulae, for $\pi_t$:
$$M_1=t$$
$$M_2=t+t^2$$
$$M_3=t+3t^2+t^3$$
$$M_4=t+6t^2+6t^3+t^4$$

Thus the mean is $E=t$, the variance is $V=t$, and the skewness is $\gamma=1/\sqrt{t}$, exactly as for $p_t$. As for the kurtosis, what happens here is that we must substract a $1$ factor from the kurtosis $3+1/t$ previously computed for $p_t$, and we obtain $\kappa=2+1/t$.
\end{proof}

\section*{7e. Exercises}

This was a quite tricky chapter, and as exercises on this, we have:

\begin{exercise}
Work out loop count results for various products of graphs.
\end{exercise}

\begin{exercise}
Learn more about Catalan numbers, the more, the better.
\end{exercise}

\begin{exercise}
Clarify what happens to the atoms, in the Stieltjes formula.
\end{exercise}

\begin{exercise}
Learn more about Hankel determinants, and the moment problem.
\end{exercise}

\begin{exercise}
Learn more about the Wigner laws, and their various properties.
\end{exercise}

\begin{exercise}
Do the missing computations for the Marchenko-Pastur laws.
\end{exercise}

\begin{exercise}
Learn more about the arcsine laws, both plain and modified.
\end{exercise}

\begin{exercise}
Conjecture something about $g_t$, the classical version of $\gamma_t$.
\end{exercise}

As bonus exercise, you can start reading some random matrix theory.

\chapter{Beta distributions}

\section*{8a. Gamma function}

What is next? More specialized distributions I guess, and thinking first at the compactly supported case, what we already studied as densities are the rectangles, leading to the uniform laws, and the semicircles, leading to the Wigner laws. Of course, we came upon some other distributions too, such as Marchenko-Pastur and the arcsine laws, but philosophically, we are basically still at this stage, rectangles and semicircles.

\bigskip

In short, our question is, what comes after rectangles and semicircles. And here, a quite natural idea is that of looking at densities of type $f(x)=x^p(1-x)^q$, over $[0,1]$. Which proves to be a very fruitful idea, and here is our claim, about this:

\begin{claim}
The laws of the following type, suitably normalized as to have mass $1$,
$$\mu=x^p(1-x)^qdx$$ 
are very useful in statistics, and will lead us into a lot of fundamental mathematics too.
\end{claim} 

Getting started now, we must first find the normalization constant, making our measures of mass 1. When $p,q$ are positive integers, this can be done as follows:

\begin{proposition}
We have the following formula, valid for any $p,q\in\mathbb N$,
$$\int_0^1x^p(1-x)^qdx=\frac{p!q!}{(p+q+1)!}$$
providing us with the missing normalization constant, in order to have a density.
\end{proposition}

\begin{proof}
Our integral $I_{pq}$ can be computed by partial integration, as follows:
\begin{eqnarray*}
I_{pq}
&=&\int_0^1x^p(1-x)^qdx\\
&=&\int_0^1px^{p-1}\frac{(1-x)^{q+1}}{q+1}\,dx\\
&=&\frac{p}{q+1}\int_0^1x^{p-1}(1-x)^{q+1}dx\\
&=&\frac{p}{q+1}\,I_{p-1,q+1}
\end{eqnarray*}

Indeed, by using this recurrence, we obtain the following formula:
\begin{eqnarray*}
I_{pq}
&=&\frac{p}{q+1}\,I_{p-1,q+1}\\
&=&\frac{p}{q+1}\cdot\frac{p-1}{q+2}\,I_{p-2,q+2}\\
&=&\frac{p}{q+1}\cdot\frac{p-1}{q+2}\cdot\frac{p-2}{q+3}\,I_{p-3,q+3}\\
&\vdots&\\
&=&\frac{p}{q+1}\cdot\frac{p-1}{q+2}\cdot\frac{p-2}{q+3}\ldots\frac{1}{q+p}\,I_{0,q+p}
\end{eqnarray*}

But $I_{0n}=1/(n+1)$, as we know well from calculus, and this gives:
\begin{eqnarray*}
I_{pq}
&=&\frac{p}{q+1}\cdot\frac{p-1}{q+2}\cdot\frac{p-2}{q+3}\ldots\frac{1}{q+p}\cdot\frac{1}{q+p+1}\\
&=&\frac{p!q!}{(q+p+1)!}
\end{eqnarray*}

Thus, we are led to the conclusion in the statement.
\end{proof}

With the above discussed, before getting into the study of the corresponding probability measures, let us have as well a look at the general case, where $p,q$ are no longer integers. In order for our integral to converge we must have $p,q>-1$, and it is standard to write $p=a-1,q=b-1$ with $a,b>0$. Thus, we are led to the following question:

\begin{question}
How does our integration formula above, which reads
$$\int_0^1x^{a-1}(1-x)^{b-1}dx=\frac{(a-1)!(b-1)!}{(a+b-1)!}$$
generalize to the case of non-integer exponents, $a,b>0$?
\end{question}

And good question this is, because the integral on the left looks uncomputable, in general, at least with elementary means. However, there is a solution to this, which is actually something extremely clever, producing a major earthquake over the whole mathematics that we know, ruins everywhere, no less than that, according to:

\begin{principle}
We know that we can talk about generalized binomial coefficients,
$$\binom{n}{k}=\frac{n(n-1)\ldots(n-k+1)}{k!}$$
with $n\in\mathbb R$. However, more is true, because we can talk about factorials themselves,
$$n!=\int_0^\infty x^ne^{-x}\,dx$$
for $n>-1$. And with this, all mathematics that we know potentially generalizes.
\end{principle}

Quite interesting all this, hope you agree with me. So, let us get now into understanding what this fact exactly says. As for Question 8.3, that can wait, I mean with mathematics destroyed, down to the foundations, answering that will be an easy task.

\bigskip

At the beginning of everything we have the following key result:

\index{gamma function}
\index{generalized factorial}

\begin{theorem}
The following integral converges, for any real number $s>0$,
$$\Gamma(s)=\int_0^\infty x^{s-1}e^{-x}\,dx$$
satisfies $\Gamma(s+1)=s\Gamma(s)$, and extends the usual factorial, $\Gamma(n)=(n-1)!$ for $n\in\mathbb N$.
\end{theorem}

\begin{proof}
There are several assertions here, the idea being as follows:

\medskip

(1) Regarding the convergence claim, the integral converges at $\infty$, due to the presence of the exponential, and converges at $0$ too, because $x^{s-1}$ with $s>0$ is integrable at $0$.

\medskip

(2) Next, by partial integration we have the following formula, as claimed:
\begin{eqnarray*}
\Gamma(s+1)
&=&\int_0^\infty x^se^{-x}\,dx\\
&=&\int_0^\infty sx^{s-1}e^{-x}\,dx\\
&=&s\Gamma(s)
\end{eqnarray*}

(3) Regarding now the case $n\in\mathbb N$, for the initial value $n=1$ we have:
$$\Gamma(1)=\int_0^\infty e^{-x}dx=1$$

Thus, by recurrence we obtain $\Gamma(n)=(n-1)!$ for any $n\in\mathbb N$, as claimed.
\end{proof}

The above result is something quite magic, because it allows to talk about the factorial of any real number $n>-1$, according to the following formula:
$$n!=\Gamma(n+1)$$

Many interesting things can be said about the gamma function, notably with an explicit computation at the half-integers, with the result here being as follows:

\index{half-integers}
\index{double factorial}

\begin{theorem}
The gamma function is given at half-integers by
$$\Gamma(n)=(n-1)!\quad,\quad\Gamma\left(n+\frac{1}{2}\right)=\frac{(2n)!!}{2^n}\,\sqrt{\pi}$$
and we have the following almost uniform formula for it, valid for any $N\in\mathbb N$,
$$\Gamma\left(\frac{N}{2}\right)=\frac{(N-1)!!}{2^{(N-1)/2}}\,c$$
with $M!!=(M-1)(M-3)\ldots$ as usual, and $c=\sqrt{2},\sqrt{\pi}$ for $N$ even, odd.
\end{theorem}

\begin{proof}
There are several things going on here, the idea being as follows:

\medskip

(1) Regarding the first formula, values of the gamma function at the positive integers, we know indeed, from Theorem 8.5, that for $n\in\mathbb N$ we have:
$$\Gamma(n)=(n-1)!$$

(2) Regarding now the half-integers, we first have the following computation:
\begin{eqnarray*}
\Gamma\left(\frac{1}{2}\right)
&=&\int_0^\infty x^{-1/2}e^{-x}\,dx\\
&=&\int_0^\infty y^{-1}e^{-y^2}2y\,dy\\
&=&2\int_0^\infty e^{-y^2}dy\\
&=&2\times\frac{\sqrt{\pi}}{2}\\
&=&\sqrt{\pi}
\end{eqnarray*}

Next, by using the formula $\Gamma(s+1)=s\Gamma(s)$, we have the following computations:
$$\Gamma\left(\frac{3}{2}\right)=\frac{1}{2}\,\Gamma\left(\frac{1}{2}\right)=\frac{1}{2}\,\sqrt{\pi}$$
$$\Gamma\left(\frac{5}{2}\right)=\frac{3}{2}\,\Gamma\left(\frac{3}{2}\right)=\frac{3}{4}\,\sqrt{\pi}$$
$$\Gamma\left(\frac{7}{2}\right)=\frac{5}{2}\,\Gamma\left(\frac{5}{2}\right)=\frac{15}{8}\,\sqrt{\pi}$$
$$\vdots$$

Thus, we can solve the problem by recurrence, and we obtain in this way, with our usual convention $N!!=(N-1)(N-3)(N-5)\ldots$ for the double factorials:
\begin{eqnarray*}
\Gamma\left(n+\frac{1}{2}\right)
&=&\frac{2n-1}{2}\cdot\frac{2n-3}{2}\ldots\frac{5}{2}\cdot\frac{3}{2}\cdot\frac{1}{2}\,\Gamma\left(\frac{1}{2}\right)\\
&=&\frac{1.3.5\ldots(2n-1)}{2^n}\,\sqrt{\pi}\\
&=&\frac{(2n)!!}{2^n}\,\sqrt{\pi}
\end{eqnarray*}

(3) Regarding now the unification of the formulae that we have in (1) and (2), let us first rewrite the formula that we just found in (2), in the following way:
$$\Gamma\left(\frac{2n+1}{2}\right)=\frac{(2n)!!}{2^n}\,\sqrt{\pi}$$

Which looks good and in final form, no question about this, so for a unification we are left with refurbishing the nice formula found in (1), in the following way:
\begin{eqnarray*}
\Gamma\left(\frac{2n}{2}\right)
&=&\Gamma(n)\\
&=&(n-1)!\\
&=&\frac{2.4.6\ldots(2n-2)}{2^{n-1}}\\
&=&\frac{(2n-1)!!}{2^{n-1}}\\
&=&\frac{(2n-1)!!}{2^{n-1/2}}\,\sqrt{2}
\end{eqnarray*}

And with this, job done, we are led to the uniform formula in the statement.
\end{proof}

Quite interestingly, the gamma function at half-integers is related as well to geometry, and more specifically, to the core formulae of geometry in $N$ dimensions, as follows:

\index{volume of sphere}
\index{area of sphere}

\begin{theorem}
The volume and area of the unit sphere in $\mathbb R^N$, given by
$$V=\left(\frac{\pi}{2}\right)^{[N/2]}\frac{2^N}{(N+1)!!}\quad,\quad 
A=\left(\frac{\pi}{2}\right)^{[N/2]}\frac{2^N}{(N-1)!!}$$
can be expressed in terms of the gamma function, the formulae being
$$V=\left(\frac{\pi}{2}\right)^{[N/2]}\frac{2^{(N-1)/2}c}{\Gamma(N/2+1)}\quad,\quad 
A=\left(\frac{\pi}{2}\right)^{[N/2]}\frac{2^{(N+1)/2}c}{\Gamma(N/2)}$$
with $c=\sqrt{2},\sqrt{\pi}$ for $N$ even, odd, as before.
\end{theorem}

\begin{proof}
There are several things going on here, the idea being as follows:

\medskip

(1) In order to compute the volume of the sphere $V$, we will need the Wallis formula, which is as follows, with $\varepsilon(p)=1$ if $p$ is even, and $\varepsilon(p)=0$ if $p$ is odd:
$$\int_0^{\pi/2}\cos^pt\,dt=\left(\frac{\pi}{2}\right)^{\varepsilon(p)}\frac{p!!}{(p+1)!!}$$

For proving this, we use partial integration. We have the following formula:
\begin{eqnarray*}
(\cos^pt\sin t)'
&=&p\cos^{p-1}t(-\sin t)\sin t+\cos^pt\cos t\\
&=&p\cos^{p+1}t-p\cos^{p-1}t+\cos^{p+1}t\\
&=&(p+1)\cos^{p+1}t-p\cos^{p-1}t
\end{eqnarray*}

By integrating this between $0$ and $\pi/2$, we obtain, for our Wallis integral $I_p$:
$$(p+1)I_{p+1}=pI_{p-1}$$

Thus we can compute $I_p$ by recurrence, and we obtain in this way:
\begin{eqnarray*}
I_p
&=&\frac{p-1}{p}\,I_{p-2}\\
&=&\frac{p-1}{p}\cdot\frac{p-3}{p-2}\,I_{p-4}\\
&=&\frac{p-1}{p}\cdot\frac{p-3}{p-2}\cdot\frac{p-5}{p-4}\,I_{p-6}\\
&\vdots&\\
&=&\frac{p!!}{(p+1)!!}\,I_{1-\varepsilon(p)}
\end{eqnarray*}

Now the initial data being $I_0=\pi/2$ and $I_1=1$, we obtain the result. 

\medskip

(2) Next, if we denote by $V_N$ the volume of the unit sphere in $\mathbb R^N$, we have:
\begin{eqnarray*}
V_N
&=&\int_{-1}^1(1-x^2)^{(N-1)/2}dx\cdot V_{N-1}\\
&=&2V_{N-1}\int_0^1(1-x^2)^{(N-1)/2}dx\\
&=&2V_{N-1}\int_0^{\pi/2}(1-\sin^2 t)^{(N-1)/2}\cos tdt\\
&=&2V_{N-1}\int_0^{\pi/2}\cos^{N-1}t\cos tdt\\
&=&2V_{N-1}\int_0^{\pi/2}\cos^Ntdt
\end{eqnarray*}

Now by recurrence, and by using the Wallis formula in (1), we obtain:
\begin{eqnarray*}
V_N
&=&2^N\int_0^{\pi/2}\cos^Ntdt\int_0^{\pi/2}\cos^{N-1}tdt\ldots\ldots\int_0^{\pi/2}\cos tdt\\
&=&2^N\left(\frac{\pi}{2}\right)^{\varepsilon(N)+\varepsilon(N-1)+\ldots+\varepsilon(1)}\frac{N!!}{(N+1)!!}\cdot\frac{(N-1)!!}{N!!}\ldots\frac{1!!}{2!!}\\
&=&\left(\frac{\pi}{2}\right)^{\varepsilon(N)+\varepsilon(N-1)+\ldots+\varepsilon(1)}\frac{2^N}{(N+1)!!}\\
&=&\left(\frac{\pi}{2}\right)^{[N/2]}\frac{2^N}{(N+1)!!}
\end{eqnarray*}

Thus, we have the formula for the volume of the sphere $V$ in the statement.

\medskip

(3) Regarding now the area of the sphere $A$, this is related, a bit as in 2D where $2\pi=2\cdot\pi$, to the volume of the sphere $V$ by the formula $A=N\cdot V$, and this gives: 
$$A=N\left(\frac{\pi}{2}\right)^{[N/2]}\frac{2^N}{(N+1)!!}=\left(\frac{\pi}{2}\right)^{[N/2]}\frac{2^N}{(N-1)!!}$$

(4) Getting now to the second series of formulae, involving gamma, recall from Theorem 8.6 that we have the following formula, with $c=\sqrt{2},\sqrt{\pi}$ for $N$ even, odd:
$$\Gamma\left(\frac{N}{2}\right)=\frac{(N-1)!!}{2^{(N-1)/2}}\,c$$

Thus, the double factorials can be expresed as follows, in terms of gamma:
$$(N-1)!!=\frac{2^{(N-1)/2}}{c}\,\Gamma\left(\frac{N}{2}\right)$$

(5) But with this in hand, the previous volume formula for the sphere becomes:
\begin{eqnarray*}
V
&=&\left(\frac{\pi}{2}\right)^{[N/2]}\frac{2^N}{(N+1)!!}\\
&=&\left(\frac{\pi}{2}\right)^{[N/2]}2^N\cdot\frac{c}{2^{(N+1)/2}}\cdot\frac{1}{\Gamma(N/2+1)}\\
&=&\left(\frac{\pi}{2}\right)^{[N/2]}\frac{2^{(N-1)/2}c}{\Gamma(N/2+1)}
\end{eqnarray*}

As for the formula for the area of the sphere, this becomes, again as claimed:
\begin{eqnarray*}
A
&=&\left(\frac{\pi}{2}\right)^{[N/2]}\frac{2^N}{(N-1)!!}\\
&=&\left(\frac{\pi}{2}\right)^{[N/2]}2^N\cdot\frac{c}{2^{(N-1)/2}}\cdot\frac{1}{\Gamma(N/2)}\\
&=&\left(\frac{\pi}{2}\right)^{[N/2]}\frac{2^{(N+1)/2}c}{\Gamma(N/2)}
\end{eqnarray*}

(6) Time perhaps for some numerics? At $N=2$ we have $c=\sqrt{2}$, and we get:
$$V=\frac{\pi}{2}\cdot\frac{\sqrt{2}\cdot\sqrt{2}}{\Gamma(2)}=\pi
\quad,\quad A=\frac{\pi}{2}\cdot\frac{2\sqrt{2}\cdot\sqrt{2}}{\Gamma(1)}=2\pi$$

Also, at $N=3$ we have $c=\sqrt{\pi}$, and we obtain, again as we should:
$$V=\frac{\pi}{2}\cdot\frac{2\cdot\sqrt{\pi}}{3\sqrt{\pi}/4}=\frac{4\pi}{3}
\quad,\quad A=\frac{\pi}{2}\cdot\frac{4\cdot\sqrt{\pi}}{\sqrt{\pi}/2}=4\pi$$

(7) Thus, formulae proved and doublechecked. And with the remark, however, that the original formulae, in terms of double factorials, are quite often better in practice.
\end{proof}

As a last piece of mathematics regarding the gamma function, we have:

\index{Stirling formula}
\index{generalized Stirling formula}

\begin{theorem}
The gamma function is given at half-integers by
$$\Gamma(n)\simeq\left(\frac{n}{e}\right)^n\sqrt{\frac{2\pi}{n}}
\quad,\quad\Gamma\left(n+\frac{1}{2}\right)\simeq\left(\frac{n}{e}\right)^n\sqrt{2\pi}$$
and we have the following uniform approximation formula for it,
$$\Gamma\left(\frac{N}{2}\right)\simeq\left(\frac{N}{2e}\right)^{N/2}\sqrt{\frac{4\pi}{N}}$$
valid for $N\in\mathbb N$, in the $N\to\infty$ limit.
\end{theorem}

\begin{proof}
This is very standard, based on the formulae in Theorem 8.6:

\medskip

(1) At the usual integers $n\in\mathbb N$ we have the following Stirling estimate:
\begin{eqnarray*}
\Gamma(n)
&=&(n-1)!\\
&\simeq&\left(\frac{n-1}{e}\right)^{n-1}\sqrt{2\pi(n-1)}\\
&=&\left(\frac{n}{e}\right)^{n-1}\left(\frac{n-1}{n}\right)^{n-1}\sqrt{2\pi(n-1)}\\
&\simeq&\left(\frac{n}{e}\right)^{n-1}\frac{\sqrt{2\pi n}}{e}\\
&=&\left(\frac{n}{e}\right)^n\sqrt{\frac{2\pi}{n}}
\end{eqnarray*}

(2) At the half-integers, $n+1/2$ with $n\in\mathbb N$, we have the following estimate:
\begin{eqnarray*}
\Gamma\left(n+\frac{1}{2}\right)
&=&\frac{(2n)!!}{2^n}\,\sqrt{\pi}\\
&=&\frac{(2n)!}{4^nn!}\,\sqrt{\pi}\\
&\simeq&\left(\frac{2n}{e}\right)^{2n}\sqrt{4\pi n}\left(\frac{e}{n}\right)^n\frac{1}{\sqrt{2\pi n}}\cdot\sqrt{\pi}\\
&=&\left(\frac{n}{e}\right)^n\sqrt{2\pi}
\end{eqnarray*}

(3) Regarding the uniform formula, with $N=2n$ the estimate in (1) reads:
$$\Gamma\left(\frac{N}{2}\right)\simeq\left(\frac{N}{2e}\right)^{N/2}\sqrt{\frac{4\pi}{N}}$$

(4) Also, with $N=2n+1$, the estimate that we found in (2) reads:
\begin{eqnarray*}
\Gamma\left(\frac{N}{2}\right)
&\simeq&\left(\frac{N-1}{2e}\right)^{(N-1)/2}\sqrt{2\pi}\\
&=&\left(\frac{N}{2e}\right)^{(N-1)/2}\left(\frac{N-1}{N}\right)^{(N-1)/2}\sqrt{2\pi}\\
&\simeq&\left(\frac{N}{2e}\right)^{(N-1)/2}\sqrt{\frac{2\pi}{e}}\\
&\simeq&\left(\frac{N}{2e}\right)^{N/2}\sqrt{\frac{4\pi}{N}}
\end{eqnarray*}

We are therefore led to the uniform formula in the statement.

\medskip

(5) Finally, let us mention that the uniform approximation formula that we found works in fact in general, with the formula being as follows, for $s>>0$:
$$\Gamma(s)\simeq \left(\frac{s}{e}\right)^s\sqrt{\frac{2\pi}{s}}$$

And we will leave some learning here, generalized Stirling formula, as an exercise.
\end{proof}

Many other things can be said about the gamma function, some in relation with the usual factorial of the integers, and some other being analysis specific. We will be back to this, in what follows. For the moment, what we have in the above will do.

\bigskip

Finally, let us mention that $\Gamma(s)$ can be defined in fact at any $s\in\mathbb C$ with $Re(s)>0$, and this because in the convergence discussion from the proof of Theorem 8.5, the phase does not matter. In what follows, we will only need gamma over the reals.

\section*{8b. Beta distributions}

With the above discussed, time to get back to probability? In relation with Question 8.3, that we started with, we have the following result, answering that question:

\begin{theorem}
We have the following formula, valid for any $a,b>0$,
$$\int_0^1x^{a-1}(1-x)^{b-1}dx=\frac{\Gamma(a)\Gamma(b)}{\Gamma(a+b)}$$
generalizing our previous computation for $a,b\in\mathbb N$, which was as follows:
$$\int_0^1x^{a-1}(1-x)^{b-1}dx=\frac{(a-1)!(b-1)!}{(a+b-1)!}$$
The above quantity is denoted $B(a,b)$, and called beta function. 
\end{theorem}

\begin{proof}
As a first comment, you might wonder why calling B beta, but the thing is, the letter B is originally Greek, upper-case $\beta$. As for the proof, this goes as follows:

\medskip

(1) Let us attempt to compute $\Gamma(a)\Gamma(b)$. By using the definition of gamma, this amounts in computing a certain integral over $\mathbb R^2$, as follows:
\begin{eqnarray*}
\Gamma(a)\Gamma(b)
&=&\int_0^\infty x^{a-1}e^{-x}dx\int_0^\infty y^{b-1}e^{-y}dy\\
&=&\int_0^\infty\int_0^\infty x^{a-1}y^{b-1}e^{-x-y}dxdy
\end{eqnarray*}

(2) Now comes the trick. Let us perform the following change of variables:
$$\begin{cases}x=st\\ y=s(1-t)\end{cases}\iff
\begin{cases}s=x+y\\ t=x/(x+y)\end{cases}$$ 

We have then $dxdy=sdsdt$, coming from the following Jacobian computation, and more on such things later, when discussing multivariable calculus:
$$\frac{dxdy}{dsdt}
=\begin{vmatrix}dx/ds&dx/dt\\ \\ dy/ds&dy/dt\end{vmatrix}
=\begin{vmatrix}t&s\\ \\1-t&-s\end{vmatrix}=-s$$

(3) We conclude that our integral over $\mathbb R^2$ above takes the following form:
\begin{eqnarray*}
\Gamma(a)\Gamma(b)
&=&\int_0^1\int_0^\infty(st)^{a-1}(s(1-t))^{b-1}e^{-s}sdsdt\\
&=&\int_0^1s^{a+b-1}e^{-s}ds\int_0^\infty t^{a-1}(1-t)^{b-1}dt\\
&=&\Gamma(a+b)\int_0^\infty t^{a-1}(1-t)^{b-1}dt
\end{eqnarray*}

Thus, we are led to the conclusions in the statement.
\end{proof}

Before getting into probability, along the lines of Claim 8.1, let us study a bit more the beta function. According to Theorem 8.9, this is given by the following formula:
$$B(a,b)=\frac{\Gamma(a)\Gamma(b)}{\Gamma(a+b)}$$

We also know from Theorem 8.9 that for integer variables $a,b\in\mathbb N$ we have:
$$B(a,b)=\frac{(a-1)!(b-1)!}{(a+b-1)!}$$

Thus, the beta function is some sort of generalized binomial coefficient, written upside down, and coming with an extra factor, or a factor missing, upon taste. As a basic result now about this function, inspired from what happens for the usual binomials, we have:

\begin{theorem}
We have the following formula,
$$\sum_{s=0}^n\binom{n}{s}B(s+a,n-s+b)=B(a,b)$$
generalizing the usual Vandermonde formula for binomials.
\end{theorem}

\begin{proof}
This is very standard, as for the usual Vandermonde formula, as follows:

\medskip

(1) At $n=0$ the formula to be proved is $B(a,b)=B(a,b)$, true.

\medskip

(2) At $n=1$ now, the formula to be proved is as follows, called Pascal formula:
$$B(a,b+1)+B(a+1,b)=B(a,b)$$

But this can be proved by using the identity $\Gamma(s+1)=s\Gamma(s)$, as follows:
\begin{eqnarray*}
B(a,b+1)+B(a+1,b)
&=&\frac{\Gamma(a)\Gamma(b+1)}{\Gamma(a+b+1)}+\frac{\Gamma(a+1)\Gamma(b)}{\Gamma(a+b+1)}\\
&=&\frac{b\Gamma(a)\Gamma(b)}{(a+b)\Gamma(a+b)}+\frac{a\Gamma(a)\Gamma(b)}{(a+b)\Gamma(a+b)}\\
&=&\frac{\Gamma(a)\Gamma(b)}{\Gamma(a+b)}\\
&=&B(a,b)
\end{eqnarray*}

(3) At $n=2$ now, the formula to be proved is as follows:
$$B(a,b+2)+2B(a+1,b+1)+B(a+2,b)=B(a,b)$$

But this can be proved by using the above Pascal formula three times, as follows:
\begin{eqnarray*}
&&B(a,b+2)+2B(a+1,b+1)+B(a+2,b)\\
&=&[B(a,b+2)+B(a+1,b+1)]+[B(a+1,b+1)+B(a+2,b)]\\
&=&B(a,b+1)+B(a+1,b)\\
&=&B(a,b)
\end{eqnarray*}

(4) And so on, with the general case following by iterating Pascal, as above.
\end{proof}

Many other things can be said about the beta function, using the theory of the gamma function developed before, notably with formulae at half-integers, other special values, Stirling estimates, and complex extensions too. In what concerns us, what we have in the above will do. In fact, save for Theorem 8.9 making the link with Claim 8.1, and for Theorem 8.10 which looks really good in beta formulation, our policy will be always that of writing everything in terms of gamma, where the theory available is just fine.

\bigskip

Time, eventually, to do some probability? Based on what we have, we can now formulate the following key definition, which is something quite far-reaching:

\begin{definition}
The beta distribution of parameters $a,b>0$ is
$$\mu=\frac{x^{a-1}(1-x)^{b-1}}{B(a,b)}\,dx$$
supported on $[0,1]$, with $B(a,b)=\Gamma(a)\Gamma(b)/\Gamma(a+b)$ being the beta function.
\end{definition}

To be more precise, this is exactly what we wanted to do, in Claim 8.1, with the mathematics for the mass 1 condition being provided by Theorem 8.9.

\bigskip

In practice now, many things can be said about the beta distribution, the idea being that this is something very useful in statistics, when a sampling over $[0,1]$, more sophisticated than the uniform one, is needed. For more on this, we recommend a solid statistics book. In what concerns us, we will be back to this, but still in a theoretical form, a bit later in this chapter, with an application of such things to the discrete setting.

\bigskip

Getting now to moment computations, the result here is as follows:

\begin{theorem}
The moments of the beta distribution are given by
$$M_k=\frac{a(a+1)\ldots(a+k-1)}{(a+b)(a+b+1)\ldots(a+b+k-1)}$$
and in particular, the mean and variance are given by the following formulae,
$$E=\frac{a}{a+b}\quad,\quad V=\frac{ab}{(a+b)^2(a+b+1)}$$
and the skewness and kurtosis are given by the following formulae,
$$\gamma=\frac{2(b-a)\sqrt{a+b+1}}{(a+b+2)\sqrt{ab}}\quad,\quad
\kappa=3+\frac{6[(a-b)^2(a+b+1)-ab(a+b+2)]}{ab(a+b+2)(a+b+3)}$$
with $a,b>0$ being as usual the parameters.
\end{theorem}

\begin{proof}
The moments can be indeed computed as follows, using Theorem 8.9:
\begin{eqnarray*}
M_k
&=&\frac{1}{B(a,b)}\int_0^1x^{a+k-1}(1-x)^{b-1}dx\\
&=&\frac{B(a+k,b)}{B(a,b)}\\
&=&\frac{\Gamma(a+k)\Gamma(b)}{\Gamma(a+b+k)}\cdot\frac{\Gamma(a+b)}{\Gamma(a)\Gamma(b)}\\
&=&\frac{\Gamma(a+k)}{\Gamma(a)}\cdot\frac{\Gamma(a+b)}{\Gamma(a+b+k)}\\
&=&\frac{a(a+1)\ldots(a+k-1)}{(a+b)(a+b+1)\ldots(a+b+k-1)}
\end{eqnarray*}

In particular, we have the formula of $E$. Getting now to the variance, this is:
\begin{eqnarray*}
V
&=&M_2-M_1^2\\
&=&\frac{a(a+1)}{(a+b)(a+b+1)}-\left(\frac{a}{a+b}\right)^2\\
&=&\frac{a}{(a+b)^2(a+b+1)}[(a+1)(a+b)-a(a+b+1)]\\
&=&\frac{ab}{(a+b)^2(a+b+1)}
\end{eqnarray*}

As for the computations of the skewness and kurtosis, these are similar.
\end{proof}

Finally, observe that at $a=b=3/2$ the beta distribution has density $\frac{4}{\pi}\sqrt{x(1-x)}$ on $[0,1]$, so is a shifted semicircle law. Thus, what we are doing here is, at least from an abstract point of view, a natural continuation of the material from chapter 7.

\section*{8c. Beta binomials}

As a next topic of discussion, that we will soon discover to be very interesting as well, related to the material regarding the negative hypergeometric laws from chapter 3, we can talk about beta binomial distributions, constructed in the following way:

\index{beta binomial law}
\index{beta function}

\begin{definition}
The beta binomial distributions are the discrete laws given by
$$P(s)=\binom{n}{s}\frac{B(s+a,n-s+b)}{B(a,b)}$$
with $B(a,b)=\Gamma(a)\Gamma(b)/\Gamma(a+b)$ being as usual the beta function.
\end{definition}

As a first job, let us verify that we have indeed probability measures, of mass 1. But this comes indeed from the Vandermonde formula in Theorem 8.10, as follows:
$$\sum_{s=0}^nP(s)
=\frac{1}{B(a,b)}\sum_{s=0}^n\binom{n}{s}B(s+a,n-s+b)
=1$$

Next, we have the following fact, which is the key to understanding our laws:

\begin{proposition}
At integer values of the parameters, $a,b\in\mathbb N$, we have
$$P(s)=\frac{\binom{s+a-1}{s}\binom{n-s+b-1}{n-s}}{\binom{n+a+b-1}{n}}$$
and our beta binomial distribution is a negative hypergeometric law.
\end{proposition}

\begin{proof}
At integer values of the parameters $a,b$, the beta function is given by:
$$B(a,b)=\frac{\Gamma(a)\Gamma(b)}{\Gamma(a+b)}=\frac{(a-1)!(b-1)!}{(a+b-1)!}$$

By using this, the formula of the beta binomial distribution is as follows:
\begin{eqnarray*}
P(s)
&=&\binom{n}{s}\frac{B(s+a,n-s+b)}{B(a,b)}\\
&=&\frac{n!}{s!(n-s)!}\cdot\frac{(s+a-1)!(n-s+b-1)!}{(n+a+b-1)!}\cdot\frac{(a+b-1)!}{(a-1)!(b-1)!}\\
&=&\frac{n!(a+b-1)!}{(n+a+b-1)!}\cdot\frac{(s+a-1)!}{s!(a-1)!}\cdot\frac{(n-s+b-1)!}{(n-s)!(b-1)!}\\
&=&\frac{\binom{s+a-1}{s}\binom{n-s+b-1}{n-s}}{\binom{n+a+b-1}{n}}
\end{eqnarray*}

But this is a negative hypergeometric law, of parameters as follows:
$$(N,m,r)=(n+a+b-1,n,a)$$

Thus, we are led to the conclusion in the statement.
\end{proof}

We can see now the interest in Definition 8.13, the point being that, in certain situations where the negative hypergeometric laws do not produce accurate results, we can use the beta binomial distributions, instead. It is possible to be a bit more precise here, the idea being that the uniform sampling used in chapter 3 in order to define the negative hypergeometric laws can be replaced with a sampling based on a continuous beta distribution, and this is what produces the beta binomial laws from Definition 8.13.

\bigskip

In what follows we will be mostly interested in the abstract mathematics of the beta binomial laws, coming as a continuation of our previous work from chapter 3, on the negative hypergeometric laws. So, as our next task, coming as a continuation of Proposition 8.14, let us see what happens when $a,b$ are half-integers. We have here:

\begin{theorem}
At half-integer parameters, $a,b\in\mathbb N/2$, we still have the formula
$$P(s)=\frac{\binom{s+a-1}{s}\binom{n-s+b-1}{n-s}}{\binom{n+a+b-1}{n}}$$
with the following convention for the quantities in brackets,
$$\binom{n}{k}=\frac{(2n+1)!!}{(2k+1)!!(2n-2k+1)!!}$$
where $M!!=(M-1)(M-3)\ldots$ for any integer $M\in\mathbb N$, as usual.
\end{theorem}

\begin{proof}
This is something standard, based on Theorem 8.6, as follows:

\medskip

(1) To start with, we have the following formula, valid for any $n\in\mathbb N$:
$$n!=\frac{(2n+1)!!}{2^n}$$

Thus, at $n,k\in\mathbb N$, the quantities in brackets in the statement are indeed binomials:
$$\binom{n}{k}=\frac{n!}{k!(n-k)!}=\frac{(2n+1)!!}{(2k+1)!!(2n-2k+1)!!}$$

(2) Next, we recall from Theorem 8.6 that we have the following formula for the gamma function at the half-integers, with $c=\sqrt{2},\sqrt{\pi}$ for $N$ even, odd:
$$\Gamma\left(\frac{N}{2}\right)=\frac{(N-1)!!}{2^{(N-1)/2}}\,c$$

Thus, at half-integer values of the parameters, $a,b\in\mathbb N/2$, say $a=A/2$ and $b=B/2$ with $A,B\in\mathbb N$, the beta function is given by the following formula:
\begin{eqnarray*}
B(a,b)
&=&\frac{\Gamma(a)\Gamma(b)}{\Gamma(a+b)}\\
&=&\frac{\Gamma(A/2)\Gamma(B/2)}{\Gamma((A+B)/2)}\\
&=&\frac{(A-1)!!}{2^{(A-1)/2}}\,c_A\cdot\frac{(B-1)!!}{2^{(B-1)/2}}\,c_B
\cdot\frac{2^{(A+B-1)/2}}{(A+B-1)!!}\cdot\frac{1}{c_{A+B}}\\
&=&\frac{(A-1)!!(B-1)!!}{(A+B-1)!!}\cdot\frac{\sqrt{2}c_Ac_B}{c_{A+B}}
\end{eqnarray*}

(3) By using now (1) and (2), the formula of the beta binomial distribution is:
\begin{eqnarray*}
P(s)
&=&\binom{n}{s}\frac{B(s+a,n-s+b)}{B(a,b)}\\
&=&\frac{(2n+1)!!}{(2s+1)!!(2n-2s+1)!!}\cdot\frac{(2s+A-1)!!(2n-2s+B-1)!!}{(2n+A+B-1)!!}\cdot\frac{\sqrt{2}c_Ac_B}{c_{A+B}}\\
&&\cdot\frac{(A+B-1)!!}{(A-1)!!(B-1)!!}\cdot\frac{c_{A+B}}{\sqrt{2}c_Ac_B}\\
&=&\frac{(2n+1)!!(2s+A-1)!!(2n-2s+B-1)!!(A+B-1)!!}
{(2s+1)!!(2n-2s+1)!!(2n+A+B-1)!!(A-1)!!(B-1)!!}\\
&=&\frac{(2n+1)!!(A+B-1)!!}{(2n+A+B-1)!!}\cdot\frac{(2s+A-1)!!}{(2s+1)!!(A-1)!!}\cdot\frac{(2n-2s+B-1)!!}{(2n-2s+1)!!(B-1)!!}\\
&=&\frac{\binom{s+a-1}{s}\binom{n-s+b-1}{n-s}}{\binom{n+a+b-1}{n}}
\end{eqnarray*}

Thus, we are led to the conclusion in the statement.
\end{proof}

Regarding now the mean and variance of the beta binomial laws, we have:

\begin{theorem}
The mean and variance of the beta binomial laws are
$$E=\frac{na}{a+b}\quad,\quad V=\frac{nab(n+a+b)}{(a+b)^2(a+b+1)}$$
in terms of the parameters $(n,a,b)$.
\end{theorem}

\begin{proof}
This is something quite straightforward, the idea being as follows:

\medskip

(1) We recall from chapter 3 that the computation of the mean for the negative hypergeometric laws was based on the following trick, and the Vandermonde formula:
$$s=(s+r)-r$$

To be more precise, in terms of our new parameters $(n,a,b)$ here, which are given by $(N,m,r)=(n+a+b-1,n,a)$ the computation from chapter 3 was as follows:
\begin{eqnarray*}
E
&=&\sum_ss\,\frac{\binom{s+a-1}{s}\binom{n+b-s-1}{n-s}}{\binom{n+a+b-1}{n}}\\
&=&\sum_s[(s+a)-a]\,\frac{\binom{s+a-1}{s}\binom{n+b-s-1}{n-s}}{\binom{n+a+b-1}{n}}\\
&=&\left[\sum_s(s+a)\,\frac{\binom{s+a-1}{s}\binom{n+b-s-1}{n-s}}{\binom{n+a+b-1}{n}}\right]-a\\
&=&\frac{a}{\binom{n+a+b-1}{n}}\left[\sum_s\frac{s+a}{a}\binom{s+a-1}{s}\binom{n+b-s-1}{n-s}\right]-a\\
&=&\frac{a}{\binom{n+a+b-1}{n}}\left[\sum_s\binom{s+a}{s}\binom{n+b-s-1}{n-s}\right]-a\\
&=&\frac{a}{\binom{n+a+b-1}{n}}\binom{n+a+b}{n}-a\\
&=&a\cdot\frac{n!(a+b-1)!}{(n+a+b-1)!}\cdot\frac{(n+a+b)!}{n!(a+b)!}-a\\
&=&\frac{(n+a+b)a}{a+b}-a\\
&=&\frac{na}{a+b}
\end{eqnarray*}

But the same trick works in general, at any $a,b>0$, by using the generalized Vandermonde formula from Theorem 8.10 for the summation, as desired.

\medskip

(2) Next, again by following the material in chapter 3, corresponding to the case $a,b\in\mathbb N$, in order to compute $M_2-M_1$, we can use a similar trick, as follows:
$$s^2-s=(s+r)(s+r+1)-2(r+1)(s+r)+r(r+1)$$

Indeed, as before in (1), by following the computation from chapter 3, we obtain:
$$M_2-M_1=\frac{n(n-1)a(a+1)}{(a+b)(a+b+1)}$$

(3) We conclude that the second moment is given by the following formula:
\begin{eqnarray*}
M_2
&=&(M_2-M_1)+M_1\\
&=&\frac{n(n-1)a(a+1)}{(a+b)(a+b+1)}+\frac{na}{a+b}\\
&=&\frac{na(na+n+b)}{(a+b)(a+b+1)}
\end{eqnarray*}

(4) Getting now to the variance, this is given by the following formula:
\begin{eqnarray*}
V
&=&M_2-M_1^2\\
&=&\frac{na(na+n+b)}{(a+b)(a+b+1)}-\frac{n^2a^2}{(a+b)^2}\\
&=&\frac{nab(n+a+b)}{(a+b)^2(a+b+1)}
\end{eqnarray*}

Thus, we are led to the formulae in the statement.
\end{proof}

Regarding the higher moments, again in view of what we found in chapter 3 for the negative hypergeometric laws, things are more complicated. In order to deal with this question, let us introduce the following notion, which is something quite subtle:

\begin{definition}
The factorial moments of $f:X\to\mathbb R$ are the numbers
$$\widetilde{M}_k=E[f(f-1)\ldots(f-k+1)]$$
with the convention $\widetilde{M}_0=1$.
\end{definition}

As a first comment here, you would probably say, we already have the usual moments $M_k$, the central moments $M_k'$, and the normalized central moments $M_k''$, so why complicating things with one more definition. In answer, there is some magic in Definition 8.17, the idea being that the normalization there ``kills the underlying partitions''.

\bigskip

As an illustration for this principle, for the Poissons law $p_t$, which are the central objects in discrete probability theory, we have the following remarkable result:

\begin{theorem}
The factorial moments of the Poisson law $p_t$ are
$$\widetilde{M_k}=t^k$$
and with this being something infinitely simpler than the formula of $M_k$. 
\end{theorem}

\begin{proof}
A first approach is via limits of binomial laws. Indeed, for the binomial law $b_{np}$ of parameters $n\in\mathbb N$ and $p\in[0,1]$, we have the following formula:
\begin{eqnarray*}
\widetilde{M}_k(b_{np})
&=&\sum_{s=k}^ns(s-1)\ldots(s-k+1)\binom{n}{s}p^s(1-p)^{n-s}\\
&=&\sum_{s=k}^n\frac{s!}{(s-k)!}\cdot\frac{n!}{s!(n-s)!}\,p^s(1-p)^{n-s}\\
&=&\sum_{s=k}^n\frac{n!}{(s-k)!(n-s)!}\,p^s(1-p)^{n-s}\\
&=&\sum_{r=0}^{n-k}\frac{n!}{r!(n-k-r)!}\,p^{k+r}(1-p)^{n-k-r}\\
&=&\frac{n!p^k}{(n-k)!}\sum_{r=0}^{n-k}\frac{(n-k)!}{r!(n-k-r)!}\,p^r(1-p)^{n-k-r}\\
&=&\frac{n!p^k}{(n-k)!}\sum_{r=0}^{n-k}\binom{n-k}{r}p^r(1-p)^{n-k-r}\\
&=&\frac{n!p^k}{(n-k)!}[p+(1-p)]^{k-r}\\
&=&\frac{n!p^k}{(n-k)!}
\end{eqnarray*}

Now recall from the Poisson Limit Theorem from chapter 4 that we have $b_{np}\to p_t$, when $p=t/n$ with $t>0$ fixed. Thus, by doing a Poisson limit we obtain, as claimed:
\begin{eqnarray*}
\widetilde{M_k}(p_t)
&=&\widetilde{M_k}\left(\lim_{n\to\infty}b_{n,t/n}\right)\\
&=&\lim_{n\to\infty}\widetilde{M_k}(b_{n,t/n})\\
&=&\lim_{n\to\infty}\frac{n!(t/n)^k}{(n-k)!}\\
&=&t^k\lim_{n\to\infty}\frac{n(n-1)\ldots(n-k+1)}{n^k}\\
&=&t^k
\end{eqnarray*}

Alternatively, we can use the moment formula for the Poisson laws that we established in chapter 4, which was as follows, with $|.|$ being as usual the number of blocks:
$$M_k(p_t)=\sum_{\pi\in P(k)}t^{|\pi|}$$

Now in order to convert this into a formula regarding the factorial moments, we need a general conversion formula, connecting the usual and the factorial moments. So, recall from chapter 3 the following magic formula, valid for any number $s\in\mathbb N$:
\begin{eqnarray*}
s^k
&=&\sum_{i_1=1}^s\ldots\sum_{i_k=1}^s1\\
&=&\sum_{\pi\in P(k)}\#\left\{(i_1,\ldots,i_k)\in\{1,\ldots,s\}^k\Big|\ker i=\pi\right\}\\
&=&\sum_{\pi\in P(k)}\frac{s!}{(s-|\pi|)!}
\end{eqnarray*}

Now since this is an equality of polynomials in $s$, we can replace our $s\in\mathbb N$ by anything, and in particular, by an arbitrary random variable $f:X\to\mathbb R$. We obtain:
$$f^k=\sum_{\pi\in P(k)}f(f-1)\ldots(f-|\pi|+1)$$

By taking now the expectation, on both sides, we obtain the following conversion formula, connecting the usual and the factorial moments of a variable $f:X\to\mathbb R$:
$$M_k=\sum_{\pi\in P(k)}\widetilde{M}_{|\pi|}$$

But with this formula in hand, back to the Poisson law, we can convert our formula above for the moments of $p_t$ into a formula for the factorial moments, and we obtain:
$$\widetilde{M_k}(p_t)=t^k$$

Thus, one way or another, we are led to the formula in the statement. We will be back to such things later in this book, when discussing partitions and cumulants.
\end{proof}

Getting now to the negative hypergeometric laws, the result for them is:

\begin{theorem}
The factorial moments of the negative hypergeometric laws are
$$\widetilde{M}_k=\frac{m!}{(m-k)!}\cdot\frac{r!}{(r-k)!}\cdot\frac{(N-k)!}{N!}$$
with $(N,m,r)$ being as usual the parameters.
\end{theorem}

\begin{proof}
We recall from chapter 3 that the moments of the negative hypergeometric law of parameters $(N,m,r)$ are given by the following formula:
$$M_k=\sum_{\pi\in P(k)}\frac{m!}{(m-|\pi|)!}\cdot\frac{(r+|\pi|-1)!}{(r-1)!}\cdot\frac{(N-m)!}{(N-m+|\pi|)!}$$

Now recall from the proof of Theorem 8.18 that we have the following formula:
$$M_k=\sum_{\pi\in P(k)}\widetilde{M}_{|\pi|}$$

Thus, we are led via some combinatorics to the formula in the statement, namely:
$$\widetilde{M}_k=\frac{m!}{(m-k)!}\cdot\frac{r!}{(r-k)!}\cdot\frac{(N-k)!}{N!}$$

We will leave the details here as an exercise, and we will be back to such things later in this book, when discussing partitions and cumulants.
\end{proof}

Getting now to what we wanted to do, for the beta binomial distribution, we have:

\index{factorial moments}

\begin{theorem}
The factorial moments of the beta binomial distribution are
$$\widetilde{M}_k=\frac{n!}{(n-k)!}\cdot\frac{B(a+k,b)}{B(a,b)}$$
with $B$ being as usual the beta function.
\end{theorem}

\begin{proof}
Many things can be said here, the idea being as follows:

\medskip

(1) To start with, let us record the following more digest formulation of the moment formula in the statement, in terms of the gamma function only:
\begin{eqnarray*}
\widetilde{M}_k
&=&\frac{n!}{(n-k)!}\cdot\frac{B(a+k,b)}{B(a,b)}\\
&=&\frac{n!}{(n-k)!}\cdot\frac{\Gamma(a+k)\Gamma(b)}{\Gamma(a+b+k)}\cdot\frac{\Gamma(a+b)}{\Gamma(a)\Gamma(b)}\\
&=&\frac{n!}{(n-k)!}\cdot\frac{\Gamma(a+k)}{\Gamma(a)}\cdot\frac{\Gamma(a+b)}{\Gamma(a+b+k)}
\end{eqnarray*}

(2) Equivalently, by using the formula $\Gamma(s+1)=s\Gamma(s)$ from Theorem 8.5, we have the following formula, directly in terms of the parameters $a,b>0$:
$$\widetilde{M_k}=\frac{n(n-1)\ldots(n-k+1)a(a+1)\ldots(a+k-1)}{(a+b)(a+b+1)\ldots(a+b+k-1)}$$

(3) But with this we can see, as a first observation, that our formula fits with the moment formulae at $k=1,2$ from Theorem 8.16 and its proof, which were as follows:
$$\widetilde{M}_1=\frac{na}{a+b}\quad,\quad 
\widetilde{M}_2=\frac{n(n-1)a(a+1)}{(a+b)(a+b+1)}$$

(4) Next, in the case of integer parameters, what we have fits with the formula in Theorem 8.19, for the negative hypergeometric laws. Indeed, let us set:
$$(N,m,r)=(n+a+b-1,n,a)$$

But with this, we obtain indeed the formulae from Theorem 8.19.

\medskip

(5) In general now, this is something quite standard, by using the various tricks from the proof of Theorem 8.16, as for the negative hypergeometric laws. We will leave this as an exercise, and come back to it later, when discussing partitions and cumulants.
\end{proof}

As a main application, we can now compute the skewness and kurtosis:

\begin{theorem}
The first factorial moments of the beta binomial distribution are
$$\widetilde{M}_1=\frac{na}{a+b}$$
$$\widetilde{M}_2=\frac{n(n-1)a(a+1)}{(a+b)(a+b+1)}$$
$$\widetilde{M}_3=\frac{n(n-1)(n-2)a(a+1)(a+2)}{(a+b)(a+b+1)(a+b+2)}$$
$$\widetilde{M}_4=\frac{n(n-1)(n-2)(n-3)a(a+1)(a+2)(a+3)}{(a+b)(a+b+1)(a+b+2)(a+b+3)}$$
and this allows computing $(E,V,\gamma,\kappa)$, with $E,V$ being as before, the skewness being 
$$\gamma=\frac{(2n+a+b)(b-a)}{a+b+2}\sqrt{\frac{a+b+1}{nab(n+a+b)}}$$
and with the formula of the kurtosis looking fairly bad.
\end{theorem}

\begin{proof}
This is indeed something self-explanatory, based on the general moment formula from Theorem 8.20, and I would leave the details to you, as an exercise.
\end{proof}

\section*{8d. Negative versions}

We would like to end this chapter, and present Part II, and first half of this book, with some philosophy, meaning general discussion, and some difficult theorems coming without proof. Of course, I know that you are a bit tired, after all our recent computations, and believe me, so am I, but this is precisely the point, it is when you're tired that it is best to sit down, relax, and have some thinking at the difficulties of life.

\bigskip

The material in this chapter was dictated by Principle 8.4, telling us that we can talk about factorials $n!$ with $n>-1$ real, which is something of key importance in modern mathematics. Indeed, with such beasts in hand, we can have combinatorics of continuous  nature, and then probability and other mathematics of continuous nature too. 

\bigskip

Moreover, as briefly explained when discussing the gamma function, we can talk in fact about factorials $n!$ with $n$ complex satisfying $Re(n)>-1$, with this being something elementary. And there is even more, because with a bit of complex analytic know-how, we can talk in fact about factorials $n!$, with $n$ complex satisfying $Re(n)\notin\!-\mathbb N$.

\bigskip

Quite amazing all this, but in practice, when getting into this continuous mathematics, be that of real or complex nature, things quickly become complicated, with the complexity of your functions drastically increasing with the number of parameters. And with this phenomenon, called ``special function blues'', being something that we already started to experience, in the context of our computations from this chapter.

\bigskip

So, looking for more blues? Nothing better than that, mathematical modesty is a key asset, that must be acquired from somewhere, and special functions are here for that, for every mathematician on this planet. So, getting back to probability, let us formulate:

\index{negative beta binomial law}

\begin{definition}
The negative beta binomial distribution of parameters $r,a,b>0$ is
$$P(s)=\frac{\Gamma(s+b)}{s!\Gamma(b)}\cdot\frac{B(s+r,a+b)}{B(r,a)}$$
with $B=\Gamma(a)\Gamma(b)/\Gamma(a+b)$ being as usual the beta function.
\end{definition}

Obviously, this is something quite tricky. As a first observation, we can express the density in terms of the beta function only, the formula being as follows:
$$P(s)=\frac{1}{s}\cdot\frac{B(s+r,a+b)}{B(r,a)B(s,b)}$$

Alternatively, in terms of the gamma function only, the formula is as follows:
$$P(s)=\frac{\Gamma(r+s)\Gamma(a+b)\Gamma(r+a)\Gamma(s+b)}{\Gamma(r+s+a+b)\Gamma(r)\Gamma(a)\Gamma(s+1)\Gamma(b)}$$

Observe that, by rearranging terms, we have as well the following formula:
$$P(s)=\frac{1}{s}\cdot\frac{B(r+a,s+b)}{B(a,b)B(r,s)}$$

Finally, when all parameters are integers, $r,a,b\in\mathbb N$, our gamma formula reads:
$$P(s)=\frac{(r+s-1)!(a+b-1)!(r+a-1)!(s+b-1)!}{(r+s+a+b-1)!(r-1)!(a-1)!s!(b-1)!}$$

In practice, the simplest case is $r=1$, and we have here the following result:

\begin{proposition}
The negative beta binomial distribution at $r=1$ is given by
$$P(s)=\frac{B(a+1,b+s)}{B(a,b)}$$
and with this being called beta geometric law of parameters $a,b>0$.
\end{proposition}

\begin{proof}
This comes from our second beta-only formula above, which reads:
$$P(s)=\frac{1}{s}\cdot\frac{B(a+1,s+b)}{B(a,b)B(1,s)}=\frac{\Gamma(s+1)}{s\Gamma(s)}\cdot\frac{B(a+1,b+s)}{B(a,b)}$$

Thus, we are led to the formula in the statement.
\end{proof}

Generally speaking, the negative beta binomial distributions are best thought as coming in 3 stages, namely the beta geometric laws above, corresponding to the case $r=1$, then their generalizations with $r\in\mathbb N$, and then their further generalizations with $r>0$. Many other things can be here, and as before for the usual beta binomial distributions, for more on the phenomenology of these laws, we refer to a solid statistics book.

\bigskip

Getting now to the mathematics of these laws, as a central result, we have:

\begin{theorem}
The factorial moments of the negative beta binomial law are
$$\widetilde{M}_k=\frac{\Gamma(r+k)\Gamma(a-k)\Gamma(b+k)}{\Gamma(r)\Gamma(a)\Gamma(b)}$$
with $r,a,b>0$ being as usual the parameters.
\end{theorem}

\begin{proof}
This is something quite complicated, that we will not attempt to prove here, but feel free of course to have a look at it, for getting that special function blues I was talking about. This being said, a few comments. As a first observation, by using $\Gamma(s+1)=s\Gamma(s)$, the formula in the statement can be written as follows:
$$\widetilde{M}_k=\frac{(r+k-1)(r+k-2)\ldots r(b+k-1)(b+k-2)\ldots b}{(a-1)(a-2)\ldots(a-k)}$$

In practice, at $k=0$ this tells us that we have $\widetilde{M}_0=1$, with this showing that our measures are indeed of mass one, as they should. Next, at $k=1,2$ we get:
$$\widetilde{M}_1=\frac{rb}{a-1}\quad,\quad \widetilde{M}_2=\frac{(r+1)r(b+1)b}{(a-1)(a-2)}$$

And so on, the idea being that the formula of $\widetilde{M}_k$ is something simple and useful.
\end{proof}

Regarding the standard parameters $(E,V,\gamma,\kappa)$ of our distributions, we have:

\begin{theorem}
The mean and variance of the negative beta binomial law are
$$E=\frac{rb}{a-1}\quad,\quad V=\frac{rb(r+a-1)(a+b-1)}{(a-1)^2(a-2)}$$
and the corresponding skewness is given by the following formula,
$$\gamma=\frac{(2r+a-1)(2b+a-1)}{(a-3)\sqrt{\frac{rb(r+a-1)(a+b-1)}{a-2}}}$$
assuming $a>1,2,3$ respectively. As for $\kappa$, this is given by a more complicated formula.
\end{theorem}

\begin{proof}
This follows indeed from Theorem 8.24, good exercise for you.
\end{proof}

\section*{8e. Exercises}

This was a quite technical chapter, and as exercises, all technical, we have:

\begin{exercise}
Learn about the gamma function at $z\in\mathbb C$, with $Re(z)>0$ or not.
\end{exercise}

\begin{exercise}
Learn the proof of the generalized Stirling formula, for gamma.
\end{exercise}

\begin{exercise}
Try proving the formula of beta, without using advanced calculus.
\end{exercise}

\begin{exercise}
Further meditate on the relation between beta and semicircle laws.
\end{exercise}

\begin{exercise}
Learn more about the phenomenology of the beta binomial laws.
\end{exercise}

\begin{exercise}
Compute the skewness and kurtosis of the beta binomial laws.
\end{exercise}

\begin{exercise}
Learn about the phenomenology of negative beta binomial laws.
\end{exercise}

\begin{exercise}
Compute the skewness and kurtosis of negative beta binomial laws.
\end{exercise}

As bonus exercise, read about Ramanujan, the master of special functions.

\part{Limiting results}

\ \vskip50mm

\begin{center}
{\em You move a little fast

Now wait a little while

If you want to last

You show a little style}
\end{center}

\chapter{Advanced calculus}

\section*{9a. Linear algebra}

Welcome back to probability, for the second half of this book, hope you enjoyed the holidays, and on the menu today, calculus. The point indeed is that, save for a few straightforward applications of Fubini, and some unwanted complications in the previous chapter, we managed to find our way by using one-variable calculus only. However, with the present Part III being of higher dimensional nature, we will have to learn multivariable calculus, in order to have tools, for solving the various problems on our way.

\bigskip

Before starting, in answer to a question that you might have, what about the normal laws, that we have not talked about yet? In answer, sure these are 1D, but their introduction needs the Gauss formula, which is something 2D. So, it is in fact the 2D normal laws that are the ``central'' objects, in the normal law galaxy, and with the further remark that these 2D laws, viewed as 1D complex laws, are what you need, when doing things like quantum physics. And finally, as a further twist to the story, it is in fact in our usual 3D that the normal laws really shine, solving the heat equation. More on this later.

\bigskip

Getting started now, with multivariable calculus, the idea is very simple:

\begin{principle}
The simplest functions $f:\mathbb R^N\to\mathbb R^M$ are those which are linear,
$$f(x+y)=f(x)+f(y)\quad,\quad f(\lambda x)=\lambda f(x)$$
with the arbitrary ones $f:\mathbb R^N\to\mathbb R^M$ locally appearing as perturbations of these.
\end{principle}

To be more precise here, in what regards the first assertion, that is just common sense, assuming a bit of familiarity with vectors. As for the second assertion, that is the higher dimensional analogue of the fundamental formula $f(x+t)\simeq f(x)+f'(x)t$, with the linear map in question being here $t\to f'(x)t$, that we will learn later in this chapter.

\bigskip

In view of Principle 9.1, let us first have a look at the linear maps $f:\mathbb R^N\to\mathbb R^M$. And for everything else, meaning arbitrary functions $f:\mathbb R^N\to\mathbb R^M$, and their approximation by various linear means, we can see later, once we will be truly familiar with the linear case. As a starting result for our theory, that you surely know, we have:

\index{linear map}
\index{matrix}
\index{rectangular matrix}
\index{scalar product}

\begin{theorem}
The linear maps $f:\mathbb R^N\to\mathbb R^M$ are those of the form
$$f(x)=Ax$$
with the matrix $A\in M_{M\times N}(\mathbb R)$ coming via the formula $A_{ij}=<f(e_j),e_i>$.
\end{theorem}

\begin{proof}
A linear map $f:\mathbb R^N\to\mathbb R^M$ must send a vector $x\in\mathbb R^N$ to a certain vector $f(x)\in\mathbb R^M$, all whose components are linear combinations of the components of $x$, so:
$$f\begin{pmatrix}
x_1\\
\vdots\\
x_N
\end{pmatrix}
=\begin{pmatrix}
A_{11}x_1+\ldots+A_{1N}x_N\\
\vdots\\
A_{M1}x_1+\ldots+A_{MN}x_N
\end{pmatrix}$$

Now the parameters $A_{ij}\in\mathbb R$ can be regarded as being the entries of a rectangular matrix $A\in M_{M\times N}(\mathbb R)$, and we have $f(x)=Ax$, as claimed. As for the second assertion, if we denote by $e_1,e_2,e_3,\ldots$ the standard bases of our vector spaces, we have:
$$<f(e_j),e_i>
=<Ae_j,e_i>
=(Ae_j)_i
=A_{ij}$$

Thus, we are led to the conclusions in the statement.
\end{proof}

At the level of examples, we have the following statement, which is a must-know:

\index{rotation}
\index{symmetry}
\index{projection}
\index{rank 1 projection}

\begin{proposition}
The following happen, in the plane:
\begin{enumerate}
\item The rotation of angle $t\in\mathbb R$ is given by the following matrix:
$$R_t=\begin{pmatrix}\cos t&-\sin t\\ \sin t&\cos t\end{pmatrix}$$

\item The symmetry with respect to the $Ox$ axis rotated by $t/2\in\mathbb R$ is given by:
$$S_t=\begin{pmatrix}\cos t&\sin t\\ \sin t&-\cos t\end{pmatrix}$$

\item The projection on the $Ox$ axis rotated by $t/2\in\mathbb R$ is given by:
$$P_t=\frac{1}{2}\begin{pmatrix}1+\cos t&\sin t\\ \sin t&1-\cos t\end{pmatrix}$$

\item However, the non-trivial translations are not linear maps, in our sense.
\end{enumerate}
\end{proposition}

\begin{proof}
All this is well-known and elementary, the idea being as follows:

\medskip

(1) We can guess the rotation matrix $R_t=\binom{a\ b}{c\ d}$ via its action on the standard coordinate vectors $\binom{1}{0}$ and $\binom{0}{1}$. Indeed, a quick picture shows that we must have:
$$\begin{pmatrix}a&b\\ c&d\end{pmatrix}\begin{pmatrix}1\\ 0\end{pmatrix}=
\begin{pmatrix}\cos t\\ \sin t\end{pmatrix}\quad,\quad 
\begin{pmatrix}a&b\\ c&d\end{pmatrix}\begin{pmatrix}0\\ 1\end{pmatrix}=
\begin{pmatrix}-\sin t\\ \cos t\end{pmatrix}$$

Recovering now the matrix is not complicated, because the first equation gives us the first column, and the second equation gives us the second column:
$$\binom{a}{c}=\begin{pmatrix}\cos t\\ \sin t\end{pmatrix}\quad,\quad 
\binom{b}{d}=\begin{pmatrix}-\sin t\\ \cos t\end{pmatrix}$$

Thus, we can just put together these two vectors, and we obtain our matrix.

\medskip

(2) This is again clear on a picture, drawn with $t/2$ instead of $t$, as indicated.

\medskip

(3) Again, picture drawn with $t/2$, as indicated, plus some easy trigonometry.

\medskip

(4) Indeed, the non-trivial translations satisfy $f(0)\neq0$, so they are not linear.
\end{proof}

As a next piece of general theory, we have the following key definition:

\index{eigenvalues}
\index{eigenvectors}
\index{diagonalization}

\begin{definition}
A linear map $f:\mathbb R^N\to\mathbb R^N$ is called diagonalizable if there is a basis  $v_1,\ldots,v_N\in\mathbb R^N$ such that $f$ multiplies by $\lambda_i$ in the direction $v_i$:
$$f(v_i)=\lambda_iv_i$$
In terms of the writing $f(x)=Ax$, we say that the corresponding matrix $A\in M_N(\mathbb R)$ is diagonalizable, with eingenvectors $v_i$ and eigenvalues $\lambda_i$.
\end{definition}

Here the system of vectors $v_1,\ldots,v_N\in\mathbb R^N$ being a basis means by definition that each vector $x\in\mathbb R^N$ must decompose uniquely as $x=\sum_ic_iv_i$, and we can see that the diagonalizability assumption uniquely determines $f$, which follows to be given by:
$$f\left(\sum_ic_iv_i\right)=\sum_i\lambda_ic_iv_i$$

Obviously, being diagonalizable means to be ``good'', and being not diagonalizable means to be ``bad''. In order to understand how this works, here are some examples:

\begin{proposition}
The following happen:
\begin{enumerate}
\item The rotation $R_t$ is not diagonalizable, unless at $t=0$ where it is the identity, $R_0=1$, and at $t=\pi$ where it is minus the identity, $R_\pi=-1$.

\item The symmetry $S_t$ is diagonalizable, with eigenvectors on the symmetry axis, and on its orthogonal, with respective eigenvalues $1,-1$.

\item The projection $P_t$ is diagonalizable, with the eigenvectors exactly as for the symmetry $S_t$, this time with respective eigenvalues $1,0$.

\item In fact, any projection is diagonalizable, with eigenvectors on it image, and on the orthogonal of its image, with respective eigenvalues $1,0$.
\end{enumerate}
\end{proposition}

\begin{proof}
All this is self-explanatory, with no need for any computation, and exercise for you, to figure out how all this works, just by drawing pictures, and thinking.
\end{proof}

Still in relation with diagonalization, at the general level now, we have:

\index{diagonal form}
\index{passage matrix}

\begin{theorem}
Assuming that a matrix $A\in M_N(\mathbb R)$ is diagonalizable, with eigenvectors $v_1,\ldots,v_N$ and corresponding eigenvalues $\lambda_1,\ldots,\lambda_N$, we have
$$A=PDP^{-1}$$
with the matrices $P,D\in M_N(\mathbb R)$ being given by the formulae
$$P=[v_1,\ldots,v_N]\quad,\quad 
D=diag(\lambda_1,\ldots,\lambda_N)$$
and respectively called passage matrix, and diagonal form of $A$.
\end{theorem}

\begin{proof}
We have $Pe_i=v_i$, where $\{e_i\}$ is the standard basis of $\mathbb R^N$, and so:
$$APe_i
=Av_i
=\lambda_iv_i$$

On the other hand, once again by using $Pe_i=v_i$, we have as well:
$$PDe_i
=P\lambda_ie_i
=\lambda_iPe_i
=\lambda_iv_i$$

Thus we have $AP=PD$, and so $A=PDP^{-1}$, as claimed.
\end{proof}

As an illustration for this, you can have some fun with the matrices from Proposition 9.3, with in each case the corresponding diagonalization formula $A=PDP^{-1}$ coming without much pain, because Proposition 9.5 tells us what both $P,D$ are, in each case, and the only piece of work remaining is that of figuring out what $P^{-1}$ is. Enjoy.

\bigskip

Summarizing, the diagonalizable matrices are the ``good'' ones, and their diagonalization is quite often a matter of doing some geometry. Regarding the non-diagonalizable matrices, these actually fall into two classes, ``bad'' and ``evil''. The bad ones are those which diagonalize over $\mathbb C$, with a main example here being the rotation $R_t$:

\begin{proposition}
The rotation of angle $t\in\mathbb R$ in the plane diagonalizes as:
$$R_t=\frac{1}{2}\begin{pmatrix}1&1\\i&-i\end{pmatrix}
\begin{pmatrix}e^{-it}&0\\0&e^{it}\end{pmatrix}
\begin{pmatrix}1&-i\\1&i\end{pmatrix}$$
Over the real numbers this is impossible, unless $t=0,\pi$.
\end{proposition}

\begin{proof}
The last assertion is something clear, that we already know, coming from the fact that at $t\neq0,\pi$ our rotation is a ``true'' rotation, having no eigenvectors in the plane. Regarding the first assertion, the point is that we have the following computation:
$$R_t\binom{1}{i}
=\begin{pmatrix}\cos t&-\sin t\\ \sin t&\cos t\end{pmatrix}\binom{1}{i}
=\binom{\cos t-i\sin t}{i\cos t+\sin t}
=e^{-it}\binom{1}{i}$$

We have as well a second eigenvector, as follows:
$$R_t\binom{1}{-i}
=\begin{pmatrix}\cos t&-\sin t\\ \sin t&\cos t\end{pmatrix}\binom{1}{-i}
=\binom{\cos t+i\sin t}{-i\cos t+\sin t}
=e^{it}\binom{1}{-i}$$

Thus, by using Theorem 9.6, we are led to the conclusion in the statement.
\end{proof}

As for the evil matrices, these are truly evil, a basic example being as follows:

\index{non-diagonalizable}
\index{Jordan block}

\begin{proposition}
The following matrix is not diagonalizable,
$$J=\begin{pmatrix}0&1\\0&0\end{pmatrix}$$
because it has only $1$ eigenvector.
\end{proposition}

\begin{proof}
The eigenvector/eigenvalue equation $Jv=\lambda v$ reads $\binom{y}{0}=\lambda\binom{x}{y}$, and we have two cases, depending on $\lambda$, which are as follows, and which give the result:

\medskip

(1) For $\lambda\neq0$ we must have $y=0$, coming from the second row, and so $x=0$ as well, coming from the first row, so we have no nontrivial eigenvectors. 

\medskip

(2) As for the case $\lambda=0$, here we must have $y=0$, coming from the first row, and so the eigenvectors here are the vectors of the form $\binom{x}{0}$, that is, the multiples of $\binom{1}{0}$. 
\end{proof}

At a more advanced level now, in order to study the diagonalization problem, the idea is that the eigenvectors can be grouped into linear spaces, called eigenspaces. By doing this directly over $\mathbb C$, as suggested by Proposition 9.7, the result is as follows:

\begin{theorem}
Let $A\in M_N(\mathbb C)$, and for any eigenvalue $\lambda\in\mathbb C$ define the corresponding eigenspace as being the vector space formed by the corresponding eigenvectors:
$$E_\lambda=\left\{v\in\mathbb C^N\Big|Av=\lambda v\right\}$$
These eigenspaces $E_\lambda$ are then in a direct sum position, in the sense that given vectors $v_1\in E_{\lambda_1},\ldots,v_k\in E_{\lambda_k}$ corresponding to different eigenvalues $\lambda_1,\ldots,\lambda_k$, we have:
$$\sum_ic_iv_i=0\implies c_i=0$$
In particular, we have $\sum_\lambda\dim(E_\lambda)\leq N$, with the sum being over all the eigenvalues, and our matrix is diagonalizable precisely when we have equality.
\end{theorem}

\begin{proof}
We prove the first assertion by recurrence on $k\in\mathbb N$. Assume by contradiction that we have a formula as follows, with the scalars $c_1,\ldots,c_k$ being not all zero:
$$c_1v_1+\ldots+c_kv_k=0$$

By dividing by one of these scalars, we can assume that our formula is:
$$v_k=c_1v_1+\ldots+c_{k-1}v_{k-1}$$

Now let us apply $A$ to this vector. On the left we obtain:
$$Av_k
=\lambda_kv_k
=\lambda_kc_1v_1+\ldots+\lambda_kc_{k-1}v_{k-1}$$

On the right we obtain something different, as follows:
\begin{eqnarray*}
A(c_1v_1+\ldots+c_{k-1}v_{k-1})
&=&c_1Av_1+\ldots+c_{k-1}Av_{k-1}\\
&=&c_1\lambda_1v_1+\ldots+c_{k-1}\lambda_{k-1}v_{k-1}
\end{eqnarray*}

We conclude from this that the following equality must hold:
$$\lambda_kc_1v_1+\ldots+\lambda_kc_{k-1}v_{k-1}=c_1\lambda_1v_1+\ldots+c_{k-1}\lambda_{k-1}v_{k-1}$$

On the other hand, we know by recurrence that the vectors $v_1,\ldots,v_{k-1}$ must be linearly independent. Thus, the coefficients must be equal, at right and at left:
$$\lambda_kc_1=c_1\lambda_1$$
$$\vdots$$
$$\lambda_kc_{k-1}=c_{k-1}\lambda_{k-1}$$

Now since at least one of the numbers $c_i$ must be nonzero, from $\lambda_kc_i=c_i\lambda_i$ we obtain $\lambda_k=\lambda_i$, which is a contradiction. Thus our proof by recurrence of the first assertion is complete. As for the second assertion, this follows from the first one.
\end{proof}

Finally, back to the examples, here is a particularly sweet one, in $N$ dimensions:

\index{all-one matrix}
\index{flat matrix}
\index{Fourier matrix}

\begin{theorem}
The all-one, or flat matrix, diagonalizes according to
$$\frac{1}{N}\begin{pmatrix}
1&\ldots&1\\
\vdots&&\vdots\\
1&\ldots&1\end{pmatrix}=\frac{1}{N}\,F_N
\begin{pmatrix}
1\\
&0\\
&&\ddots\\
&&&0\end{pmatrix}
F_N^*$$
where $F_N=(w^{ij})_{ij}$ with $w=e^{2\pi i/N}$ is the Fourier matrix, and $F_N^*=(w^{-ij})_{ij}$.
\end{theorem}

\begin{proof}
This is something quite tricky, worth a detailed discussion, as follows:

\medskip

(1) Let us denote by $\mathbb I_N$ the all-1 matrix, set $P_N=\mathbb I_N/N$, and consider as well the all-1 vector $\xi=(1)_i\in\mathbb R^N$. The matrices $\mathbb I_N,P_N$ act on the vectors $v\in\mathbb R^N$ as follows:
$$\mathbb I_Nv=\sum_iv_i\cdot\xi\quad,\quad P_Nv=\frac{\sum_iv_i}{N}\cdot\xi$$

Now the point is that the formula on the right shows that $P_N$ is the orthogonal projection on $\xi$, and so, as for any rank 1 projection, its diagonal form must be:
$$P_N\sim
\begin{pmatrix}
1\\
&0\\
&&\ddots\\
&&&0
\end{pmatrix}$$

(2) Next, how to diagonalize $P_N$? And the problem here is that, when looking for 0-eigenvectors, we must explicitly solve the following equation:
$$x_1+\ldots+x_N=0$$

Which is not an easy task, if we want a nice basis for the space of solutions. I mean, try this over the reals, at $N=3,4$, in order to understand what I am talking about.

\medskip

(3) Fortunately, the complex numbers, and more specifically the $N$-th roots of unity, come to the rescue. Indeed, with $w=e^{2\pi i/N}$ as in the statement, we have the following magic formula, coming from an obvious barycenter computation, in the plane:
$$\frac{1}{N}\sum_{k=0}^{N-1}w^{ks}=\delta_{N|s}$$

But this shows that the columns of the Fourier matrix $F_N=(w^{ij})_{ij}$ are the desired eigenvectors for $P_N$, and with the standard discrete Fourier analysis convention that the matrix indices are $i,j=0,1,\ldots,N-1$, we are led to the formula in the statement.

\medskip

(4) To be more precise here, we have the diagonal form from (1), and the passage matrix from (3), and the point is that the root of unity formula in (3) shows that the rescaling $U=F_N/\sqrt{N}$ is unitary, in the sense that $U^{-1}=U^*$, with $U^*=(\bar{U}_{ji})_{ij}$. Thus, by putting everything together, we are led to the formula in the statement.
\end{proof}

Finally, as a question that you might have, does the Fourier matrix $F_N$ above has something to do with our Fourier transform $F_f(x)=e^{ifx}$ from probability. In answer, yes, the idea being that, when doing probability over $\{1,\ldots,N\}$, the matrix of the Fourier transform $F_f(x)=e^{ifx}$ is the above Fourier matrix $F_N$. And, more on this, later.

\section*{9b. Spectral theorems}

In order to reach now to a more advanced diagonalization theory, we need to have a more geometric look at the matrices $A\in M_N(\mathbb C)$, and at linear algebra in general. We first have the following result, regarding the space $\mathbb C^N$ itself, which is straightforward: 

\begin{theorem}
We can talk about scalar products and lengths in $\mathbb C^N$,
$$<x,y>=\sum_ix_i\bar{y}_i\qquad,\qquad ||x||=\sqrt{\sum_i|x_i|^2}$$
which are related by the following conversion formulae,
$$||x||=\sqrt{<x,x>}\quad,\quad <x,y>
=\frac{||x+y||^2-||x-y||^2+i||x+iy||^2-i||x-iy||^2}{4}$$
and the following happen:
\begin{enumerate}
\item $<\lambda x,y>=<x,\bar{\lambda}y>=\lambda<x,y>$.

\item $<x+y,z>=<x,z>+<y,z>$.

\item $<x,y+z>=<x,z>+<y,z>$. 

\item $<x,y>=\overline{<y,x>}$.

\item $||\lambda x||=|\lambda|\cdot||x||$.

\item $|<x,y>|\leq||x||\cdot||y||$.

\item $||x+y||\leq||x||+||y||$.

\item $x\perp y\iff<x,y>=0$, by definition.
\end{enumerate}
\end{theorem}

\begin{proof}
We can certainly talk about scalar products and lengths, as above, and with the second conversion formula, called polarization identity, coming from:
\begin{eqnarray*}
&&||x+y||^2-||x-y||^2+i||x+iy||^2-i||x-iy||^2\\
&=&||x||^2+||y||^2-||x||^2-||y||^2+i||x||^2+i||y||^2-i||x||^2-i||y||^2\\
&&+2Re(<x,y>)+2Re(<x,y>)+2iIm(<x,y>)+2iIm(<x,y>)\\
&=&4<x,y>
\end{eqnarray*}

By the way, talking useful general identities, we have as well a parallelogram rule, that I forgot to mention in the above statement, which is as follows:
\begin{eqnarray*}
||x+y||^2+||x-y||^2
&=&<x+y,x+y>+<x-y,x-y>\\
&=&2||x||^2+2||y||^2+2Re(<x,y>)-2Re(<x,y>)\\
&=&2(||x||^2+||y||^2)
\end{eqnarray*}

As for the various claims in the statement, their proof is as follows:

\medskip

(1-5) All the verifications here are trivial, easy exercise for you.

\medskip

(6) Given $x,y\in\mathbb C^N$, consider the following real function, depending on $w\in\mathbb T$:
\begin{eqnarray*}
f(t)
&=&||wx+ty||^2\\
&=&<wx+ty,wx+ty>\\
&=&||x||^2+2tRe(w<x,y>)+t^2||y||^2
\end{eqnarray*}

We can see that $f$ is a degree 2 polynomial in its variable $t\in\mathbb R$, and since this polynomial is positive, its discriminant must be negative, $\Delta\leq0$. But, if we pick the parameter $w\in\mathbb T$ such that $w<x,y>=|<x,y>|$, the discriminant is given by:
$$\Delta=4|<x,y>|^2-4||x||^2||y||^2$$

Thus we have the Cauchy-Schwarz inequality, $|<x,y>|\leq||x||\cdot||y||$, as claimed. 

\medskip

(7) This can be proved by raising to the square and simplifying, as follows:
\begin{eqnarray*}
&&||x+y||\leq||x||+||y||\\
&\iff&||x+y||^2\leq||x||^2+||y||^2+2||x||\cdot||y||\\
&\iff&||x||^2+||y||^2+2Re(<x,y>)\leq||x||^2+||y||^2+2||x||\cdot||y||\\
&\iff&Re(<x,y>)\leq||x||\cdot||y||
\end{eqnarray*}

Indeed, the last inequality holds by (6), so we have our triangle inequality.

\medskip

(8) This is a definition, inspired from the fact that in $\mathbb R^N$ two vectors are orthogonal precisely when their scalar product is zero, that we included here for convenience.
\end{proof}

Back now to the complex matrices, the point is that these do not come alone, but rather in pairs $(A,A^*)$, with the basics of the $A\to A^*$ operation being as follows:

\index{adjoint operator}
\index{adjoint matrix}
\index{scalar product}
\index{unitary}
\index{projection}

\begin{theorem}
Given $A\in M_N(\mathbb C)$, define its adjoint matrix by $(A^*)_{ij}=\bar{A}_{ji}$.
\begin{enumerate}
\item $<Ax,y>=<x,A^*y>$, for any two vectors $x,y\in\mathbb C^N$.

\item $x\to Ux$ with $U\in M_N(\mathbb C)$ is an isometry precisely when $U^*=U^{-1}$.

\item $x\to Px$ with $P\in M_N(\mathbb C)$ is a projection precisely when $P^2=P^*=P$.
\end{enumerate}
\end{theorem}

\begin{proof}
This is something very standard, the idea being as follows:

\medskip

(1) The formula $<Ax,y>=<x,A^*y>$ is indeed clear from definitions on the basis vectors, $x,y\in\{e_1,\ldots,e_N\}$, and then by linearity, it must hold in general, as stated.

\medskip

(2) Given a matrix $U\in M_N(\mathbb C)$, we have indeed the following equivalences, with the first one coming from the polarization identity, and the other ones being clear:
\begin{eqnarray*}
||Ux||=||x||
&\iff&<Ux,Uy>=<x,y>\\
&\iff&<x,U^*Uy>=<x,y>\\
&\iff&U^*Uy=y\\
&\iff&U^*U=1\\
&\iff&U^*=U^{-1}
\end{eqnarray*}

(3) Given a matrix $P\in M_N(\mathbb C)$, in order for $x\to Px$ to be an oblique projection, we must have $P^2=P$. Now observe that this projection is orthogonal when:
\begin{eqnarray*}
<Px-x,Py>=0
&\iff&<P^*Px-P^*x,y>=0\\
&\iff&P^*Px-P^*x=0\\
&\iff&P^*P-P^*=0\\
&\iff&P^*P=P^*
\end{eqnarray*}

The point now is that by conjugating the last formula, we get $P^*P=P$. Now by comparing with $P^*P=P^*$, we obtain $P=P^*$, and this gives the result. 
\end{proof}

Back now to our usual diagonalization business, we have the following key result:

\begin{theorem}
Any matrix $A\in M_N(\mathbb C)$ which is self-adjoint, $A=A^*$, is diagonalizable, with the diagonalization being of the following type,
$$A=UDU^*$$
with $U\in U_N$, and with $D\in M_N(\mathbb R)$ diagonal. The converse holds too.
\end{theorem}

\begin{proof}
As a first remark, the converse trivially holds, because if we take a matrix of the form $A=UDU^*$, with $U$ unitary and $D$ diagonal and real, then we have:
$$A^*
=(UDU^*)^*
=UD^*U^*
=UDU^*
=A$$

In the other sense now, assume that $A$ is self-adjoint, $A=A^*$.  Our first claim is that the eigenvalues are real. Indeed, assuming $Av=\lambda v$, we have:
\begin{eqnarray*}
\lambda<v,v>
&=&<Av,v>\\
&=&<v,Av>\\
&=&\bar{\lambda}<v,v>
\end{eqnarray*}

Thus we obtain $\lambda\in\mathbb R$, as claimed. Our next claim now is that the eigenspaces corresponding to different eigenvalues are pairwise orthogonal. Assume indeed that:
$$Av=\lambda v\quad,\quad 
Aw=\mu w$$

We have then the following computation, by using the fact that we have $\mu\in\mathbb R$:
\begin{eqnarray*}
\lambda<v,w>
&=&<Av,w>\\
&=&<v,Aw>\\
&=&\mu<v,w>
\end{eqnarray*}

Thus $\lambda\neq\mu$ implies $v\perp w$, as claimed. In order now to finish the proof, it remains to prove that the eigenspaces of $A$ span the whole space $\mathbb C^N$. For this purpose, we will use a recurrence method. Let us pick an eigenvector of our matrix:
$$Av=\lambda v$$

Assuming now that we have a vector $w$ orthogonal to it, $v\perp w$, we have:
\begin{eqnarray*}
<Aw,v>
&=&<w,Av>\\
&=&<w,\lambda v>\\
&=&\lambda<w,v>\\
&=&0
\end{eqnarray*}

Thus, if $v$ is an eigenvector, then the vector space $v^\perp$ is invariant under $A$. Moreover, since a matrix $A$ is self-adjoint precisely when $<Av,v>\in\mathbb R$ for any vector $v\in\mathbb C^N$, as one can see by expanding the scalar product, the restriction of $A$ to the subspace $v^\perp$ is self-adjoint. Thus, we can proceed by recurrence, and we obtain the result.
\end{proof}

As an interesting consequence of the above result, in the real case, we have:

\begin{theorem}
Any matrix $A\in M_N(\mathbb R)$ which is symmetric, $A=A^t$, is diagonalizable, with the diagonalization being of the following type,
$$A=UDU^t$$
with $U\in O_N$, and with $D\in M_N(\mathbb R)$ diagonal. The converse holds too.
\end{theorem}

\begin{proof}
As before, the converse trivially holds, because if we take a matrix of the form $A=UDU^t$, with $U$ orthogonal and $D$ diagonal and real, then we have $A^t=A$. In the other sense now, this follows from Theorem 9.13, and its proof.
\end{proof}

An important class of self-adjoint matrices, which includes for instance all the orthogonal projections, are the positive matrices. The theory here is as follows:

\begin{theorem}
For a matrix $A\in M_N(\mathbb C)$ the following conditions are equivalent, and if they are satisfied, we say that $A$ is positive:
\begin{enumerate}
\item $A=B^2$, with $B=B^*$.

\item $A=CC^*$, for some $C\in M_N(\mathbb C)$.

\item $<Ax,x>\geq0$, for any vector $x\in\mathbb C^N$.

\item $A=A^*$, and the eigenvalues are positive, $\lambda_i\geq0$.

\item $A=UDU^*$, with $U\in U_N$ and with $D\in M_N(\mathbb R_+)$ diagonal.
\end{enumerate}
\end{theorem}

\begin{proof}
This is something very standard, the idea being as follows:

\medskip

$(1)\implies(2)$ This is clear, because we can take $C=B$.

\medskip

$(2)\implies(3)$ This follows from the following computation:
$$<Ax,x>
=<CC^*x,x>
=<C^*x,C^*x>
\geq0$$

$(3)\implies(4)$ By using the fact that $<Ax,x>$ is real, we have:
$$<Ax,x>
=<x,A^*x>
=<A^*x,x>$$

Thus we have $A=A^*$, and the remaining assertion, regarding the eigenvalues, follows from the following computation, assuming $Ax=\lambda x$:
$$<Ax,x>
=<\lambda x,x>
=\lambda<x,x>
\geq0$$

$(4)\implies(5)$ This follows indeed by using Theorem 9.13.

\medskip

$(5)\implies(1)$ Assuming $A=UDU^*$ with $U\in U_N$, and with $D\in M_N(\mathbb R_+)$ diagonal, we can set $B=U\sqrt{D}U^*$. Then $B$ is self-adjoint, and we have $B^2=A$, as desired.
\end{proof}

Finally, let us record the following technical version of the above result:

\begin{theorem}
For a matrix $A\in M_N(\mathbb C)$ the following conditions are equivalent, and if they are satisfied, we say that $A$ is strictly positive, and write $A>0$:
\begin{enumerate}
\item $A=B^2$, with $B=B^*$, invertible.

\item $A=CC^*$, for some $C\in M_N(\mathbb C)$ invertible.

\item $<Ax,x>>0$, for any nonzero vector $x\in\mathbb C^N$.

\item $A=A^*$, and the eigenvalues are strictly positive, $\lambda_i>0$.

\item $A=UDU^*$, with $U\in U_N$ and with $D\in M_N(\mathbb R_+^*)$ diagonal.
\end{enumerate}
\end{theorem}

\begin{proof}
This follows either from Theorem 9.15, by adding the extra assumptions in the statement, or from the proof of Theorem 9.15, by modifying where needed.
\end{proof}

With the above discussed, done with spectral theorems, and time to get a beer? You must be kidding. Indeed, we know from Proposition 9.7 that the basic rotations are diagonalizable too, so let us have a look now into this. Here is our general result:

\begin{theorem}
Any matrix $U\in M_N(\mathbb C)$ which is unitary, $U^*=U^{-1}$, is diagonalizable, with the eigenvalues on the unit circle $\mathbb T$. More precisely we have
$$U=VDV^*$$
with $V\in U_N$, and with $D\in M_N(\mathbb T)$ diagonal. The converse holds too.
\end{theorem}

\begin{proof}
As a first remark, the converse trivially holds, because given a matrix of type $U=VDV^*$, with $V\in U_N$, and with $D\in M_N(\mathbb T)$ being diagonal, we have:
$$U^*
=VD^*V^*
=(V^*)^{-1}D^{-1}V^{-1}
=U^{-1}$$

Let us prove now the first assertion. To start with, assuming $Uv=\lambda v$, we have:
\begin{eqnarray*}
<v,v>
&=&<U^*Uv,v>\\
&=&<Uv,Uv>\\
&=&<\lambda v,\lambda v>\\
&=&|\lambda|^2<v,v>
\end{eqnarray*}

Thus we obtain $\lambda\in\mathbb T$, as desired. Our claim now is that the eigenspaces corresponding to different eigenvalues are pairwise orthogonal. Assume indeed that:
$$Uv=\lambda v\quad,\quad 
Uw=\mu w$$

We have then the following computation, using $U^*=U^{-1}$ and $\lambda,\mu\in\mathbb T$:
\begin{eqnarray*}
\lambda<v,w>
&=&<Uv,w>\\
&=&<v,U^*w>\\
&=&<v,U^{-1}w>\\
&=&<v,\mu^{-1}w>\\
&=&\mu<v,w>
\end{eqnarray*}

Thus $\lambda\neq\mu$ implies $v\perp w$, as claimed. In order now to finish the proof, it remains to prove that the eigenspaces of $U$ span the whole space $\mathbb C^N$. For this purpose, let us pick an eigenvector of our matrix, $Uv=\lambda v$. Assuming that we have $v\perp w$, we obtain:
\begin{eqnarray*}
<Uw,v>
&=&<w,U^*v>\\
&=&<w,U^{-1}v>\\
&=&<w,\lambda^{-1}v>\\
&=&\lambda<w,v>\\
&=&0
\end{eqnarray*}

Thus, if $v$ is an eigenvector, then the vector space $v^\perp$ is invariant under $U$. Now since $U$ is an isometry, so is its restriction to this space $v^\perp$. Thus this restriction is a unitary, and so we can proceed by recurrence, and we obtain the result.
\end{proof}

Well done all this, and I have this feeling that you are thirsty again. Well, not yet, because we still have to unify Theorem 9.13, dealing with the self-adjoints, with Theorem 9.17, dealing with the unitaries. And fortunately, this can be done, as follows:

\begin{theorem}
Any matrix $A\in M_N(\mathbb C)$ which is normal, $AA^*=A^*A$, is diagonalizable, with the diagonalization being of the following type,
$$A=UDU^*$$
with $U\in U_N$, and with $D\in M_N(\mathbb C)$ diagonal. The converse holds too.
\end{theorem}

\begin{proof}
As a first remark, the converse trivially holds, because if we take a matrix of the form $A=UDU^*$, with $U$ unitary and $D$ diagonal, then we have:
$$AA^*
=UDD^*U^*
=UD^*DU^*
=A^*A$$

In the other sense now, this is something more technical. Our first claim is that a matrix $A$ is normal precisely when the following happens, for any vector $v$:
$$||Av||=||A^*v||$$

Indeed, the above equality can be written as follows:
$$<AA^*v,v>=<A^*Av,v>$$

But this is equivalent to $AA^*=A^*A$, by expanding the scalar products. Our claim now is that $A,A^*$ have the same eigenvectors, with conjugate eigenvalues:
$$Av=\lambda v\implies A^*v=\bar{\lambda}v$$

Indeed, this follows from the following computation, and from the trivial fact that if $A$ is normal, then so is any matrix of type $A-\lambda 1_N$:
$$||(A^*-\bar{\lambda}1_N)v||
=||(A-\lambda 1_N)^*v||
=||(A-\lambda 1_N)v||
=0$$

Let us prove now, by using this, that the eigenspaces of $A$ are pairwise orthogonal. Assume that we have two eigenvectors, corresponding to different eigenvalues, $\lambda\neq\mu$:
$$Av=\lambda v\quad,\quad 
Aw=\mu w$$

We have the following computation, which shows that $\lambda\neq\mu$ implies $v\perp w$:
\begin{eqnarray*}
\lambda<v,w>
&=&<\lambda v,w>\\
&=&<Av,w>\\
&=&<v,A^*w>\\
&=&<v,\bar{\mu}w>\\
&=&\mu<v,w>
\end{eqnarray*}

In order to finish, it remains to prove that the eigenspaces of $A$ span the whole $\mathbb C^N$. This is something that we have already seen for the self-adjoint matrices, and for unitaries, and we will use here these results, in order to deal with the general normal case. As a first observation, given an arbitrary matrix $A$, the matrix $AA^*$ is self-adjoint:
$$(AA^*)^*=AA^*$$

Thus, we can diagonalize this matrix $AA^*$, as follows, with the passage matrix being a unitary, $V\in U_N$, and with the diagonal form being real, $E\in M_N(\mathbb R)$:
$$AA^*=VEV^*$$

Now observe that, for matrices of type $A=UDU^*$, which are those that we supposed to deal with, we have the following formulae:
$$V=U\quad,\quad 
E=D\bar{D}$$

In particular, the matrices $A$ and $AA^*$ have the same eigenspaces. So, this will be our idea, proving that the eigenspaces of $AA^*$ are eigenspaces of $A$. In order to do so, let us pick two eigenvectors $v,w$ of the matrix $AA^*$, corresponding to different eigenvalues, $\lambda\neq\mu$. The eigenvalue equations are then as follows:
$$AA^*v=\lambda v\quad,\quad 
AA^*w=\mu w$$

We have the following computation, using the normality condition $AA^*=A^*A$, and the fact that the eigenvalues of $AA^*$, and in particular $\mu$, are real:
\begin{eqnarray*}
\lambda<Av,w>
&=&<\lambda Av,w>\\
&=&<A\lambda v,w>\\
&=&<AAA^*v,w>\\
&=&<AA^*Av,w>\\
&=&<Av,AA^*w>\\
&=&<Av,\mu w>\\
&=&\mu<Av,w>
\end{eqnarray*}

We conclude that we have $<Av,w>=0$. But this reformulates as follows:
$$\lambda\neq\mu\implies A(E_\lambda)\perp E_\mu$$

Now since the eigenspaces of $AA^*$ are pairwise orthogonal, and span the whole $\mathbb C^N$, we deduce from this that these eigenspaces are invariant under $A$:
$$A(E_\lambda)\subset E_\lambda$$

But with this result in hand, we can finish the proof of the theorem. Indeed, we can decompose the problem, and the matrix $A$ itself, following these eigenspaces of $AA^*$, which in practice amounts in saying that we can assume that we only have 1 eigenspace. By rescaling, this is the same as assuming that we have $AA^*=1$, and so we are now into the unitary case, that we know how to solve, as explained in Theorem 9.17.
\end{proof}

And with this, curtain, good work that we did, and time now to celebrate.

\section*{9c. Derivatives, Taylor}

Getting now to functions and analysis, in several variables, as a first result here, which is something fundamental, getting us into matrices and linear algebra, we have:

\index{derivative}
\index{partial derivative}

\begin{theorem}
A function $f:\mathbb R^N\to\mathbb R^M$ is continuously differentiable,
$$f(x+t)\simeq f(x)+f'(x)t$$
with $f'(x)$ linear, and $x\to f'(x)$ continuous, precisely when it has partial derivatives,
$$\frac{df_i}{dx_j}(x)=\lim_{t\to 0}\frac{f_i(x+te_j)-f_i(x)}{t}$$
which depend continuously on $x$. In this case the derivative is
$$f'(x)=\left(\frac{df_i}{dx_j}(x)\right)_{ij}\in M_{M\times N}(\mathbb R)$$ 
acting on the vectors $t\in\mathbb R^N$ by usual multiplication.
\end{theorem}

\begin{proof}
This is something very standard, the idea being as follows:

\medskip

(1) As a first observation, the formula in the statement makes sense indeed, as an equality, or rather approximation, of vectors in $\mathbb R^M$, as follows:
$$f\begin{pmatrix}x_1+t_1\\ \vdots\\ x_N+t_N\end{pmatrix}
\simeq f\begin{pmatrix}x_1\\ \vdots\\ x_N\end{pmatrix}
+\begin{pmatrix}
\frac{df_1}{dx_1}(x)&\ldots&\frac{df_1}{dx_N}(x)\\
\vdots&&\vdots\\
\frac{df_M}{dx_1}(x)&\ldots&\frac{df_M}{dx_N}(x)
\end{pmatrix}\begin{pmatrix}t_1\\ \vdots\\ t_N\end{pmatrix}$$

(2) Getting now to the proof, in the simplest possible case, $N=M=1$, what we have is a usual 1-variable function $f:\mathbb R\to\mathbb R$, and the formula in the statement is something that we know well from 1-variable calculus, namely the definition of the derivative:
$$f(x+t)\simeq f(x)+f'(x)t$$

(3) Next, at $N=2,M=1$ the result holds as well, by using (2) twice, as follows:
\begin{eqnarray*}
f\binom{x_1+t_1}{x_2+t_2}
&\simeq&f\binom{x_1+t_1}{x_2}+\frac{df}{dx_2}(x)t_2\\
&\simeq&f\binom{x_1}{x_2}+\frac{df}{dx_1}(x)t_1+\frac{df}{dx_2}(x)t_2\\
&=&f\binom{x_1}{x_2}+\begin{pmatrix}\frac{df}{dx_1}(x)&\frac{df}{dx_2}(x)\end{pmatrix}\binom{t_1}{t_2}
\end{eqnarray*}

(4) More generally, we can deal in this way with the general case $M=1$, with the formula here, obtained via a straightforward recurrence, being as follows:
\begin{eqnarray*}
f\begin{pmatrix}x_1+t_1\\ \vdots\\ x_N+t_N\end{pmatrix}
&\simeq&f\begin{pmatrix}x_1\\ \vdots\\ x_N\end{pmatrix}+\frac{df}{dx_1}(x)t_1+\ldots+\frac{df}{dx_N}(x)t_N\\
&=&f\begin{pmatrix}x_1\\ \vdots\\ x_N\end{pmatrix}+
\begin{pmatrix}\frac{df}{dx_1}(x)&\ldots&\frac{df}{dx_N}(x)\end{pmatrix}
\begin{pmatrix}t_1\\ \vdots\\ t_N\end{pmatrix}
\end{eqnarray*}

(5) But this gives the result in the case where both $N,M\in\mathbb N$ are arbitrary too. Indeed, consider a function $f:\mathbb R^N\to\mathbb R^M$, and let us write it as follows:
$$f=\begin{pmatrix}f_1\\ \vdots\\ f_M\end{pmatrix}$$

We can apply (4) to each of the components $f_i:\mathbb R^N\to\mathbb R$, and we get:
$$f_i\begin{pmatrix}x_1+t_1\\ \vdots\\ x_N+t_N\end{pmatrix}
\simeq f_i\begin{pmatrix}x_1\\ \vdots\\ x_N\end{pmatrix}+
\begin{pmatrix}\frac{df_i}{dx_1}(x)&\ldots&\frac{df_i}{dx_N}(x)\end{pmatrix}
\begin{pmatrix}t_1\\ \vdots\\ t_N\end{pmatrix}$$

(6) But this collection of $M$ formulae, taken altogether, tells us precisely that the formula in (1) happens indeed. Thus, we are led to the conclusion in the statement.
\end{proof}

Generally speaking, Theorem 9.19 is all you need to know, for extending to several variables the basic results from one-variable calculus. As a basic illustration, we have:

\index{chain rule}

\begin{theorem}
We have the chain derivative formula
$$(f\circ g)'(x)=f'(g(x))\cdot g'(x)$$
as an equality of matrices.
\end{theorem}

\begin{proof}
This is something standard in one variable, that we know from chapter 1, and in several variables the proof is similar, by using the notion of derivative coming from Theorem 9.19. To be more precise, consider a composition of functions, as follows:
$$f:\mathbb R^N\to\mathbb R^M\quad,\quad 
g:\mathbb R^K\to\mathbb R^N\quad,\quad 
f\circ g:\mathbb R^K\to\mathbb R^M$$

According to Theorem 9.19, the derivatives of these functions are certain linear maps, corresponding to certain rectangular matrices, as follows:
$$f'(g(x))\in M_{M\times N}(\mathbb R)\quad,\quad 
g'(x)\in M_{N\times K}(\mathbb R)\quad\quad
(f\circ g)'(x)\in M_{M\times K}(\mathbb R)$$

Thus, our formula makes sense indeed. As for proof, this comes from:
\begin{eqnarray*}
(f\circ g)(x+t)
&=&f(g(x+t))\\
&\simeq&f(g(x)+g'(x)t)\\
&\simeq&f(g(x))+f'(g(x))g'(x)t
\end{eqnarray*}

Thus, we are led to the conclusion in the statement.
\end{proof}

At a more advanced level now, at order 2, we have the following result:

\index{second derivative}
\index{Taylor formula}
\index{Hessian matrix}

\begin{theorem}
Given a twice differentiable function $f:\mathbb R^N\to\mathbb R$, we have
$$f(x+t)\simeq f(x)+f'(x)t+\frac{<f''(x)t,t>}{2}$$
with $f'(x)\in M_{1\times N}(\mathbb R)$ being a row vector, and with $f''(x)\in M_N(\mathbb R)$, given by 
$$f''(x)=\left(\frac{d^2f}{dx_idx_j}\right)_{ij}(x)$$
being the Hessian matrix of $f$, at the point $x\in\mathbb R^N$.
\end{theorem}

\begin{proof}
This is something quite tricky, the idea being as follows:

\medskip

(1) As a first remark, at $N=1$ the Hessian matrix is the $1\times1$ matrix having as entry the usual second derivative $f''(x)\in\mathbb R$, and the formula in the statement is something that we know well from 1-variable calculus, namely the Taylor formula at order 2:
$$f(x+t)\simeq f(x)+f'(x)t+\frac{f''(x)t^2}{2}$$

(2) In general now, this is in fact something which does not need a whole new proof, because it follows from the one-variable formula above, applied to the restriction of $f$ to the following segment in $\mathbb R^N$, which can be regarded as being a one-variable interval:
$$I=[x,x+t]$$

To be more precise, let $y\in\mathbb R^N$, and consider the following function, with $r\in\mathbb R$:
$$g(r)=f(x+ry)$$

We know from (1) that the Taylor formula for $g$, at the point $r=0$, reads:
$$g(r)\simeq g(0)+g'(0)r+\frac{g''(0)r^2}{2}$$

And our claim is that, with $t=ry$, this is precisely the formula in the statement.

\medskip

(3) So, let us see if our claim is correct. By using the chain rule, we have the following formula, with on the right, as usual, a row vector multiplied by a column vector:
$$g'(r)=f'(x+ry)\cdot y$$

By using again the chain rule, we can compute the second derivative as well:
\begin{eqnarray*}
g''(r)
&=&(f'(x+ry)\cdot y)'\\
&=&\left(\sum_i\frac{df}{dx_i}(x+ry)\cdot y_i\right)'\\
&=&\sum_i\sum_j\frac{d^2f}{dx_idx_j}(x+ry)\cdot\frac{d(x+ry)_j}{dr}\cdot y_i\\
&=&\sum_i\sum_j\frac{d^2f}{dx_idx_j}(x+ry)\cdot y_iy_j\\
&=&<f''(x+ry)y,y>
\end{eqnarray*}

(4) Time now to conclude. We know that we have $g(r)=f(x+ry)$, and according to our various computations above, we have the following formulae:
$$g(0)=f(x)\quad,\quad 
g'(0)=f'(x)\quad,\quad 
g''(0)=<f''(x)y,y>$$

Buit with this data in hand, the usual Taylor formula for our one variable function $g$, at order 2, at the point $r=0$, takes the following form, with $t=ry$:
\begin{eqnarray*}
f(x+ry)
&\simeq&f(x)+f'(x)ry+\frac{<f''(x)y,y>r^2}{2}\\
&=&f(x)+f'(x)t+\frac{<f''(x)t,t>}{2}
\end{eqnarray*}

Thus, we have obtained the formula in the statement.
\end{proof}

As an application of our methods, we can now study the local extrema, as follows:

\begin{theorem}
Given a twice differentiable function $f:\mathbb R^N\to\mathbb R$, in order for $x\in\mathbb R^N$ to be a local extremum, we must have $f'(x)=0$, and:
\begin{enumerate}
\item $f''(x)\geq0$ is needed for $x$ to be a local minimum.

\item $f''(x)>0$ guarantees that $x$ is a local minimum.

\item $f''(x)\leq0$ is needed for $x$ to be a local maximum.

\item $f''(x)<0$ guarantees that $x$ is a local maximum.
\end{enumerate}
\end{theorem}

\begin{proof}
The first assertion comes, as in one variable, from the following formula, the point being that the error term $f'(x)t$ must vanish, which means $f'(x)=0$:
$$f(x+t)\simeq f(x)+f'(x)t$$

Next, consider the Taylor formula at order 2, which in our case $f'(x)=0$ reads:
$$f(x+t)\simeq f(x)+\frac{<f''(x)t,t>}{2}$$

But with this, by using Theorems 9.15 and 9.16, we obtain the results.
\end{proof}

\section*{9d. Change of variables}

As a complement to the above material, let us discuss now integration. The key result here is the change of variable formula, and in order to discuss this, let us start with:

\begin{proposition}
We have the change of variable formula
$$\int_a^bf(x)dx=\int_c^df(\varphi(t))\varphi'(t)dt$$
where $c=\varphi^{-1}(a)$ and $d=\varphi^{-1}(b)$.
\end{proposition}

\begin{proof}
This is something that we know well, coming with $f=F'$ from:
$$(F\varphi)'(t)=F'(\varphi(t))\varphi'(t)$$

Indeed, by integrating this between $c$ and $d$, we obtain the result.
\end{proof}

In several variables now, we can only expect the above $\varphi'(t)$ factor to be replaced by something similar, a sort of ``derivative of $\varphi$, arising as a real number''. But this can only be the quantity $\det(\varphi'(t))$, called Jacobian, and with this in mind, we are led to:

\begin{theorem}
Given a transformation $\varphi=(\varphi_1,\ldots,\varphi_N)$, we have
$$\int_Ef(x)dx=\int_{\varphi^{-1}(E)}f(\varphi(t))|J_\varphi(t)|dt$$
with the $J_\varphi$ quantity, called Jacobian, being given by
$$J_\varphi(t)=\det\left[\left(\frac{d\varphi_i}{dx_j}(x)\right)_{ij}\right]$$ 
and with this generalizing the formula from Proposition 9.23.
\end{theorem}

\begin{proof}
This is something quite tricky, the idea being as follows:

\medskip

(1) To start with, the determinant of a matrix $A\in M_N(\mathbb R)$, that I forgot to talk about before, is the signed volume of the parallelepiped made by its column vectors:
$$A=\begin{bmatrix}v_1&\ldots&v_N\end{bmatrix}\quad\implies\quad \det A=\pm vol<v_1,\ldots,v_N>$$

To be more precise, the sign here is $+$ when we can continuously pass from the standard basis $e_1,\ldots,e_N$ to the system of vectors $v_1,\ldots,v_N$, and is $-$ otherwise.

\medskip

(2) Next, many interesting things can be said about the determinant, notably with the following key formula, with $\varepsilon:S_N\to\{\pm1\}$ being the signature of permutations:
$$\det A=\sum_{\sigma\in S_N}\varepsilon(\sigma)A_{1\sigma(1)}\ldots A_{N\sigma(N)}$$

So, exercise I guess for you, to learn all this, from a good linear algebra book.

\medskip

(3) Now back to our theorem, observe first that our formula generalizes indeed the change of variable formula in 1 dimension, the point being that the absolute value on the derivative appears as to compensate for the lack of explicit bounds for the integral.

\medskip 

(4) As a second observation, we can assume if we want, by linearity, that we are dealing with the constant function $f=1$. For this function, our formula reads:
$$vol(E)=\int_{\varphi^{-1}(E)}|J_\varphi(t)|dt$$

In terms of $D={\varphi^{-1}(E)}$, this amounts in proving that we have:
$$vol(\varphi(D))=\int_D|J_\varphi(t)|dt$$

And here, as a first remark, our formula is clear for the linear maps $\varphi$, by using the definition of the determinant of real matrices, given in (1), as a signed volume.

\medskip

(5) However, the extension of this to the case of non-linear maps $\varphi$ is something non-trivial, so we will not follow this path. In order to prove now the result, as stated, our first claim is that the validity of the theorem is stable under taking compositions of transformations $\varphi$. In order to prove this claim, consider a composition, as follows:
$$\varphi:E\to F\quad,\quad 
\psi:D\to E\quad,\quad 
\varphi\circ\psi:D\to F$$

Assuming that the theorem holds for $\varphi,\psi$, we deduce that we have, as desired:
\begin{eqnarray*}
\int_Ff(x)dx
&=&\int_Ef(\varphi(s))|J_\varphi(s)|ds\\
&=&\int_Df(\varphi\circ\psi(t))|J_\varphi(\psi(t))|\cdot|J_\psi(t)|dt\\
&=&\int_Df(\varphi\circ\psi(t))|J_{\varphi\circ\psi}(t)|dt
\end{eqnarray*}

(6) Next, as a key ingredient, let us examine the case where we are in $N=2$ dimensions, and our transformation $\varphi$ has one of the following special forms:
$$\varphi(x,y)=(\psi(x,y),y)\quad,\quad\varphi(x,y)=(x,\psi(x,y))$$

By symmetry, it is enough to deal with the first case. Here the Jacobian is $d\psi/dx$, and by replacing if needed $\psi\to-\psi$, we can assume that this Jacobian is positive, $d\psi/dx>0$. Now by assuming as before that $D=\varphi^{-1}(E)$ is a rectangle, $D=[a,b]\times[c,d]$, we can prove our formula by using the change of variables in 1 dimension, as follows:
\begin{eqnarray*}
\int_Ef(s)ds
&=&\int_{\varphi(D)}f(x,y)dxdy\\
&=&\int_c^d\int_{\psi(a,y)}^{\psi(b,y)}f(x,y)dxdy\\
&=&\int_c^d\int_a^bf(\psi(x,y),y)\frac{d\psi}{dx}\,dxdy\\
&=&\int_Df(\varphi(t))J_\varphi(t)dt
\end{eqnarray*}

(7) But with this, we can now prove the theorem, in $N=2$ dimensions. Indeed, given a transformation $\varphi=(\varphi_1,\varphi_2)$, consider the following two transformations:
$$\phi(x,y)=(\varphi_1(x,y),y)\quad,\quad \psi(x,y)=(x,\varphi_2\circ\phi^{-1}(x,y))$$

We have then $\varphi=\psi\circ\phi$, and by using (4) for $\psi,\phi$, which are of the special form there, and then (5) for composing, we conclude that the theorem holds indeed for $\varphi$, as desired. Thus, theorem proved in $N=2$ dimensions, and the extension of the above proof to arbitrary $N$ dimensions is straightforward, that we will leave here as an exercise.
\end{proof}

In order to discuss now some basic applications, in 2D, we will need:

\begin{proposition}
We have polar coordinates in $2$ dimensions,
$$\begin{cases}
x\!\!\!&=\ r\cos t\\
y\!\!\!&=\ r\sin t
\end{cases}$$
the corresponding Jacobian being $J=r$.
\end{proposition}

\begin{proof}
This is something elementary, the Jacobian being computed as follows:
$$J
=\begin{vmatrix}
\frac{d(r\cos t)}{dr}&&\frac{d(r\cos t)}{dt}\\
\\
\frac{d(r\sin t)}{dr}&&\frac{d(r\sin t)}{dt}
\end{vmatrix}
=\begin{vmatrix}
\cos t&-r\sin t\\
\sin t&r\cos t
\end{vmatrix}
=r$$

Thus, we are led to the conclusions in the statement.
\end{proof}

We can now compute the Gauss integral, which is the best calculus formula ever:

\begin{theorem}
We have the following formula,
$$\int_\mathbb Re^{-x^2}dx=\sqrt{\pi}$$
called Gauss integral formula.
\end{theorem}

\begin{proof}
This is something quite magic, the point being that all 1D methods that we know fail, for this integral. However, we can do this using two dimensions, as follows:
\begin{eqnarray*}
\left(\int_\mathbb Re^{-x^2}dx\right)^2
&=&\int_\mathbb R\int_\mathbb Re^{-x^2-y^2}dxdy\\
&=&\int_0^{2\pi}\int_0^\infty e^{-r^2}rdrdt\\
&=&2\pi\int_0^\infty\left(-\frac{e^{-r^2}}{2}\right)'dr\\
&=&2\pi\left[0-\left(-\frac{1}{2}\right)\right]\\
&=&\pi
\end{eqnarray*}

Amazing all this, isn't it. There are of course some other known proofs for the Gauss formula, as for instance the one that we met in chapter 1, when first discussing the Gauss integral, but all using two dimensions, appearing as variations of the above.
\end{proof}

\section*{9e. Exercises}

This was a very enjoyable chapter, all good old things, and as exercises, we have:

\begin{exercise}
Learn more about the Fourier matrix $F_N$, and its properties.
\end{exercise}

\begin{exercise}
Learn about the polar decomposition theorem.
\end{exercise}

\begin{exercise}
Learn as well about the singular value theorem.
\end{exercise}

\begin{exercise}
Solve some physics equations using the chain rule.
\end{exercise}

\begin{exercise}
Clarify what we said, regarding the local extrema.
\end{exercise}

\begin{exercise}
Try to find the multivariable Taylor formula, at order 3.
\end{exercise}

\begin{exercise}
Learn about the determinant, from a good linear algebra book.
\end{exercise}

\begin{exercise}
Try computing the Gauss integral without using 2D.
\end{exercise}

As bonus exercise, with basic calculus learned, you are good for Rudin \cite{rud}.

\chapter{Normal laws}

\section*{10a. Large numbers}

Welcome to convergence of random variables, and this will be a key chapter. We already met convergence in chapter 4, when talking about Poisson laws, with the Poisson Limit Theorem (PLT) and its negative version stating that we have:
$$b_{np}\to p_t\quad,\quad p=t/n\quad ,\quad \forall t>0$$
$$c_{mp}\to p_t\quad,\quad p=m/(m+t)\quad,\quad \forall t>0$$

It is of course possible to get beyond this, with other such convergence results, for instance for the hypergeometric laws from chapter 3, and their beta versions from chapter 8. However, it is not about such things, which remain of discrete nature, that we would like to talk here. The question that we would like to solve is as follows:

\begin{question}
Given a sequence of independent identically distributed (i.i.d.) variables $f_1,f_2,f_3,\ldots\,$, what can we say about their partial sums
$$S_n=f_1+\ldots+f_n$$
with $n\to\infty$? Do we have convergence, under which normalization, and in which exact sense? What are the mean, variance of the limit? Is the limiting law continuous?
\end{question}

So, this will be our question, which is something quite subtle. More in detail now, here are a few comments about this, meant to clarify what our question exactly says:

\bigskip

(1) As a first observation, the above-mentioned PLT results do not fall in this class, due to the conditionings $p=t/n$ and $p=m/(m+t)$ used for the i.i.d. variables appearing there. And so we are here, right from the beginning, into uncharted territory.

\bigskip

(2) As a second observation, due to our i.i.d. assumption on our variables, by obvious mean and variance reasons, the above sums $S_n$ cannot converge as such, without a suitable normalization, or an extra assumption on our variables. More on this later.

\bigskip

(3) Finally, we must talk about the meaning of the convergence too. Indeed, even if our variables $f_1,f_2,f_3,\ldots\,$ are discrete, their sums $S_n$, suitably normalized, can converge to bell-shaped curves. So, we need some convergence theory for such situations.

\bigskip

Quite interesting all this, hope you agree with me, we are now into delicate analysis questions, that will take us some time, to untangle. But we have this whole chapter, for this. So, let us start with the beginning, which as usual in probability, means study of the expectation $E$. We have here the following soft version of Question 10.1:

\begin{question}[soft version]
Given a sequence of i.i.d. variables $f_1,f_2,f_3,\ldots\,$, intuition tells us that we should have a convergence as follows,
$$\frac{f_1+\ldots+f_n}{n}\to E$$
 with $E=E(f_i)$ being the common expectation of our variables. So, do we really have this, and in which exact sense, in what regards the convergence?
\end{question}

Which sounds very good and concrete, as a starting question, so let us get now to work. Obviously, we must first talk about convergence of random variables. And here, we have several possible notions. As a first notion, which is quite natural, we have:

\index{convergence in law}
\index{convergence in moments}

\begin{definition}
Given variables $f_1,f_2,f_3,\ldots\,$, we say that $f_n\to f$ in law when 
$$E(\varphi(f_n))\to E(\varphi(f))$$
for any continuous function $\varphi:\mathbb R\to\mathbb R$.
\end{definition}

As a first observation, by linearity we can assume that we are dealing with the power functions $\varphi(x)=x^k$, so the following condition must be satisfied, for any $k\in\mathbb N$:
$$E(f_n^k)\to E(f^k)$$

In other words, with $M_k$ denoting as usual the moments, we must have:
$$M_k(f_n)\to M_k(f)$$

Alternatively, again by linearity, and by using some standard analysis, we can assume that we are dealing with characteristic functions of intervals, which leads to:

\begin{proposition}
We have $f_n\to f$ in law precisely when
$$P(f_n\leq x)\to P(f\leq x)$$
at all continuity points of the function $x\to P(f\leq x)$.
\end{proposition}

\begin{proof}
This is indeed very standard, coming from definitions, as explained above, by using linearity, and approximation by step functions, say good exercise for you.
\end{proof}

Moving on, as a second notion of convergence for random variables, we have:

\index{convergence in probability}

\begin{definition}
We say that we have $f_n\to f$ in probability when 
$$P(|f_n-f|\geq\varepsilon)\to 0$$
for any $\varepsilon>0$.
\end{definition}

As before with the convergence in law, there are many illustrations for this notion, and further things that can be said. Before everything, however, let us formulate:

\begin{theorem}
We have the following implication,
$$\begin{pmatrix}f_n\to f\\ {\rm in\ probability}\end{pmatrix}\implies
\begin{pmatrix}f_n\to f\\ {\rm in\ law}\end{pmatrix}$$
and the reverse implication does not necessarily hold.
\end{theorem}

\begin{proof}
As a technical ingredient, we will need the following estimate, valid for any two random variables $f,g:X\to\mathbb R$, and any two numbers $a\in\mathbb R$ and $\varepsilon>0$:
\begin{eqnarray*}
P(g\leq a)
&=&P(g\leq a,f\leq a+\varepsilon)+P(g\leq a,f>a+\varepsilon)\\
&\leq&P(f\leq a+\varepsilon)+P(g\leq a,f>a+\varepsilon)\\
&=&P(f\leq a+\varepsilon)+P(g-f\leq a-f<-\varepsilon)\\
&\leq&P(f\leq a+\varepsilon)+P(g-f<-\varepsilon)\\
&\leq&P(f\leq a+\varepsilon)+P(|g-f|>\varepsilon)
\end{eqnarray*}

In order to prove now the result, the most convenient is to use the criterion for convergence in law from Proposition 10.4. That is, given a continuity point $a\in\mathbb R$ of the function $x\to P(f\leq x)$, we would like to prove that the following happens:
$$P(f_n\leq a)\to P(f\leq a)$$

So, let $\varepsilon>0$. By using the estimate found above, we have:
$$P(f_n\leq a)\leq P(f\leq a+\varepsilon)+P(|f_n-f|>\varepsilon)$$
$$P(f\leq a-\varepsilon)\leq P(f_n\leq a)+P(|f_n-f|>\varepsilon)$$

We conclude from this that we have the following estimate:
\begin{eqnarray*}
P(f\leq a+\varepsilon)-P(|f_n-f|>\varepsilon)
&\leq&P(f_n\leq a)\\
&\leq&P(f\leq a+\varepsilon)+P(|f_n-f|>\varepsilon)
\end{eqnarray*}

Now with $n\to\infty$ we obtain from this, in terms of $F(a)=P(x\leq a)$:
$$F(a-\varepsilon)\leq\lim_{n\to\infty}P(f_n\leq a)\leq F(a+\varepsilon)$$

Since $F$ was assumed to be continuous at $a$, with $\varepsilon\to 0$ this gives:
$$\lim_{n\to\infty}P(f_n\leq a)=P(f\leq a)$$

Thus we have indeed the convergence $f_n\to f$ in law, as stated. As for the counterexamples at the end, which are elementary, we will leave them as an exercise.
\end{proof}

\index{convergence almost surely}
\index{a.s.}
\index{almost surely}
\index{convergence a.s.}
\index{measure zero}

Moving on, as a third and final notion of convergence, we have:

\begin{definition}
We say that we have $f_n\to f$ almost surely when 
$$f_n(x)\to f(x)$$
outside a measure zero set.
\end{definition}

Again, there are many illustrations for this notion, and further things that can be said, about this. At the theoretical level, we have the following key result:

\begin{theorem}
We have the following implication
$$\begin{pmatrix}f_n\to f\\ {\rm almost\ surely}\end{pmatrix}\implies
\begin{pmatrix}f_n\to f\\ {\rm in\ probability}\end{pmatrix}$$
and the reverse implication does not necessarily hold.
\end{theorem}

\begin{proof}
Assume that $f_n\to f$ in probability in our sense, namely:
$$P(|f_n-f|\geq\varepsilon)\to 0\quad,\quad\forall\varepsilon>0$$

Given a continuous function $\varphi:\mathbb R\to\mathbb R$, we can write:
\begin{eqnarray*}
|E(\varphi(f_n))-E(\varphi(f))|
&=&|E(\varphi(f_n)-\varphi(f))|
\end{eqnarray*}

The point now is that we can decompose the expectation on the right over the sets where $|f_n-f|\geq\varepsilon$ and $|f_n-f|<\varepsilon$, and we obtain zero in the limit, as desired. As for the counterexamples at the end, which are elementary, we will leave them as an exercise.
\end{proof}

As a conclusion to all this, we have the following grand result:

\begin{theorem}
We have the following implications,
$$\begin{pmatrix}f_n\to f\\ {\rm almost\ surely}\end{pmatrix}\implies
\begin{pmatrix}f_n\to f\\ {\rm in\ probability}\end{pmatrix}\implies
\begin{pmatrix}f_n\to f\\ {\rm in\ law}\end{pmatrix}$$
and the reverse implications do not necessarily hold.
\end{theorem}

\begin{proof}
This comes indeed by putting together what we have, namely Theorem 10.6 and Theorem 10.8. Thus, we are led to the conclusion in the statement.
\end{proof}

Well, looks like we have some good theory going on here, reinforced concrete foundations for convergence, but before going further, let us listen as well to what cat has to say. Cat indeed, that we last met in chapter 5, in relation with some other axiomatic aspects, is constantly meowing since I started this chapter. So, what is it, katie?

\begin{cat}
In quantum physics there are no well-defined probability spaces $X$, and the convergence is always in law. Why bothering with more.
\end{cat}

Well, thanks cat, this sounds quite deep, as a comment. This being said, I should have told you earlier, I'm not writing here a quantum physics book, and I feel genuinely interested in concepts like almost sure convergence. So, thanks again, and good luck with your mice chasing in law, I'm pretty much sure that you'll have $P=1$ catches.

\bigskip

Getting now to what we wanted to do, namely convergence of an average of i.i.d. variables towards their common mean, we have the following result, about this:

\index{law of large numbers}
\index{large numbers}

\begin{theorem}[Weak law of large numbers]
Given i.i.d. variables $f_1,f_2,f_3,\ldots\,$, with common mean $E(f_i)=E$, we have the following convergence in probability,
$$\frac{f_1+\ldots+f_n}{n}\to E$$
and as a consequence, this convergence holds as well in law.
\end{theorem}

\begin{proof}
This is something very standard, the idea being as follows:

\medskip

(1) Let us first establish the weaker result, namely the convergence in law, for our averages. Consider the sequence of averages of our variables:
$$g_n=\frac{f_1+\ldots+f_n}{n}$$

Now let us look at the Fourier transforms of our variables. We have, for any $i$, the following estimate, with $E=E(f_i)$ being the common expectation of our variables:
$$F_{f_i}(t)\simeq 1+iEt$$

Now by using independence, we deduce from this that we have, with $n\to\infty$:
$$F_{g_n}(t)
=\left[F_{f_i}\left(\frac{t}{n}\right)\right]^n
\simeq\left[1+iE\,\frac{t}{n}\right]^n
\simeq e^{iEt}$$

Now since the limiting function that we found $e^{iEt}$ is the Fourier transform of the constant variable $E$, we conclude from this that we have $g_n\to E$ in law.

\medskip

(2) In order to establish now the convergence in probability, as stated, we recall from chapter 2 that given a random variable $h:X\to\mathbb R$, having mean $E$ and variance $V$, the Chebycheff inequality states that we have, for any $a>0$:
$$P\big(|h-E|\geq a\big)\leq\frac{V}{a^2}$$

To be more precise, this was indeed something established in chapter 2, when doing discrete probability at that time, and the proof is identical in the general, continuous case. In order to prove now our result, consider the following sequence of variables:
$$h_n=\frac{f_1+\ldots+f_n}{n}-E$$

The mean and variance of these variables are then as follows, with $v=V(f_i)$:
$$E(h_n)=0\quad,\quad V(h_n)=\frac{v}{n}$$

Indeed, the mean formula $E(h_n)=0$ is clear, and by using this, we have as well:
\begin{eqnarray*}
V(h_n)
&=&E(h_n^2)\\
&=&E\left(\frac{\sum_{ij}f_if_j}{n^2}-\frac{2E\sum_if_i}{n}+E^2\right)\\
&=&E\left(\frac{\sum_if_i^2}{n^2}+\frac{\sum_{i\neq j}f_if_j}{n^2}-\frac{2E\sum_if_i}{n}+E^2\right)\\
&=&\frac{\sum_iE(f_i^2)}{n^2}+\frac{\sum_{i\neq j}E(f_i)E(f_j)}{n^2}-\frac{2E\sum_iE(f_i)}{n}+E^2\\
&=&\frac{v+E^2}{n}+n(n-1)\cdot\frac{E^2}{n^2}-2E^2+E^2\\
&=&\frac{v}{n}+\left(n+n(n-1)-n^2\right)\frac{E^2}{n^2}\\
&=&\frac{v}{n}
\end{eqnarray*}

Now let us apply the Chebycheff inequality above, to the variable $h_n$. This gives:
$$P(|h_n|\geq\varepsilon)\leq\frac{v/n}{\varepsilon^2}$$

Thus, we have the following inequality, in terms of the original variables $f_i$:
$$P\left(\left|\frac{f_1+\ldots+f_n}{n}-E\right|\geq\varepsilon\right)\leq\frac{v}{\varepsilon^2 n}$$

But this gives, as desired, the convergence in probability, namely:
$$\lim_{n\to\infty}P\left(\left|\frac{f_1+\ldots+f_n}{n}-E\right|\geq\varepsilon\right)=0$$

Thus, we are led to the conclusion in the statement.
\end{proof}

The above result is not the end of the story with the law of large numbers, because, quite remarkably, the convergence there holds almost surely. In order to explain this, which is something quite technical, we will need some preliminaries. Let us start with:

\index{i.o.}
\index{infinitely often}

\begin{definition}
Given a sequence of events $A_1,A_2,A_3,\ldots\,$, we set
$$\limsup_nA_n=\bigcap_{m\in\mathbb N}\bigcup_{n=m}^\infty A_n$$
with the limit notation coming from the decreasing intersection. Equivalently, we set
$$\limsup_nA_n=\left\{w\in X\Big|w\in A_n\ {\rm i.o.}\right\}$$
with i.o. standing for infinitely often.
\end{definition}

To be more precise here, given a sequence of events $A_1,A_2,A_3,\ldots$ as above, we certainly have a decreasing sequence of events, as follows:
$$\bigcup_{n=1}^\infty A_n\supset\bigcup_{n=2}^\infty A_n\supset\bigcup_{n=3}^\infty A_n\supset\ldots$$

Thus, we can consider the intersection of this sequence, and in analogy with what we know about numbers and limits, it makes sense to denote this intersection as above: 
$$\limsup_nA_n=\lim_{m\to\infty}\bigcup_{n=m}^\infty A_n$$

Finally, it is clear from definitions that we have the following equivalent formula for our intersection, with i.o. standing for infinitely often:
$$\limsup_nA_n=\left\{w\in X\Big|w\in A_n\ {\rm i.o.}\right\}$$

We can now state a key technical result, the Borel-Cantelli lemma:

\index{Borel-Cantelli}

\begin{lemma}[Borel-Cantelli]
Given events $A_1,A_2,A_3,\ldots\,$, we have
$$\sum_{n=1}^\infty P(A_n)<\infty\implies P(w\in A_n\ {\rm i.o.})=0$$
with i.o. standing as usual for infinitely often.
\end{lemma}

\begin{proof}
We use the interpretation of $P(w\in A_n\ {\rm i.o.})$ coming from Definition 10.12. In that context, since the intersection there is decreasing, we have, for any $m\in\mathbb N$:
$$\limsup_nA_n\subset\bigcup_{n=m}^\infty A_n$$

We conclude from this that we have the following estimate, for any $m\in\mathbb N$:
$$P(\limsup_nA_n)\leq\sum_{n=m}^\infty P(A_n)$$

On the other hand, the assumption in the statement is that the series $\sum_{n=1}^\infty P(A_n)$ converges, and this can be stated in the following way:
$$\lim_{m\to\infty}\sum_{n=m}^\infty P(A_n)=0$$

Thus $P(\limsup_nA_n)$ is arbitrarily small, and we conclude that we have:
$$P(\limsup_nA_n)=0$$

But this is the same as saying that $P(w\in A_n\ {\rm i.o.})=0$, as stated.
\end{proof}

As a second technical ingredient, we will need the following straightforward generalization of the Chebycheff inequality from chapter 2, featuring an exponent $p$:

\index{Chebycheff inequality}

\begin{proposition}
We have the following Chebycheff inequality,
$$P\big(|g|\geq a\big)\leq\frac{E(g^p)}{a^p}$$
valid for any random variable $g:X\to\mathbb R$, and any $p\in\mathbb N$.
\end{proposition}

\begin{proof}
We have indeed the following computation, $\varphi$ being the density:
\begin{eqnarray*}
E(|g|^p)
&=&\int_\mathbb R|x|^p\varphi(x)dx\\
&=&\int_{|x|\geq a}|x|^p\varphi(x)dx\\
&\geq&a^p\int_{|x|\geq a}\varphi(x)dx\\
&=&a^pP(|g|\geq a)
\end{eqnarray*}

Thus, we are led to the Chebycheff inequality in the statement.
\end{proof}

We can now improve the weak law of large numbers, into something final, namely:

\index{strong law of large numbers}

\begin{theorem}[Strong law of large numbers]
Given i.i.d. variables $f_1,f_2,f_3,\ldots\,$, with common mean $E(f_i)=E$, we have the convergence
$$\frac{f_1+\ldots+f_n}{n}\to E$$
happening almost surely.
\end{theorem}

\begin{proof}
This is something quite technical, the idea being as follows:

\medskip

(1) By replacing the variables $f_i$ by their centered versions $f_i-E$, we can assume:
$$E=0$$

Consider now the sums of our variables, which have also zero mean:
$$g_n=f_1+\ldots+f_n$$

We would like to prove that the following convergence is almost sure:
$$\frac{g_n}{n}\to0$$

(2) In order to do this, consider an arbitrary event $w\in X$ such that:
$$\lim_{n\to\infty}\frac{g_n(w)}{n}\neq0$$

In this situation, there exists $\varepsilon>0$ such that, for infinitely many $n$:
$$\left|\frac{g_n(w)}{n}\right|>\varepsilon$$

Thus, in order to prove our theorem, we must show that we have, for any $\varepsilon>0$:
$$P(|g_n|>\varepsilon n\ {\rm i.o.})=0$$

(3) But for this latter purpose, we can use the Borel-Cantelli lemma, applied to:
$$A_n=\left\{w\in X\Big||g_n(w)|\geq n\varepsilon\right\}$$

Indeed, let us first verify that the convergence assumption in the Borel-Cantelli lemma is indeed satisfied. For this purpose, we will use the generalized Chebycheff inequality from Proposition 10.14, with exponent $p=4$. By using our assumption $E=0$, we have the following computation, with $M_2,M_4$ being the second and fourth moments of $f_i$:
\begin{eqnarray*}
E(g_n^4)
&=&E\left(\sum_{ijkl}f_if_jf_kf_k\right)\\
&=&E\left(\sum_if_i^4\right)+6E\left(\sum_{i<j}f_i^2f_j^2\right)\\
&=&nM_4+6\binom{n}{2}M_2^2\\
&=&nM_4+3n(n-1)M_2^2
\end{eqnarray*}

But this shows that for $n>>0$, we have a constant $C$ such that:
$$E(g_n^4)<Cn^2$$

(4) Thus, the generalized Chebycheff inequality from Proposition 10.14 applies to our situation, with exponent $p=4$, and gives the following estimate:
$$P(|g_n|\geq n\varepsilon)\leq\frac{E(g_n^4)}{(n\varepsilon)^4}\leq\frac{C}{\varepsilon^4n^2}$$

And with this, we are almost there. Indeed, if we denote by $N\in\mathbb N$ the smallest integer where the inequality found in (3) holds, we have:
$$\sum_{n\geq N}P(|g_n|\geq n\varepsilon)\leq\sum_{n\geq N}\frac{C}{\varepsilon^4n^2}<\infty$$

Thus, the Borel-Cantelli lemma applies, and shows that we have:
$$P(|g_n|>\varepsilon n\ {\rm i.o.})=0$$

But this ends the proof of our theorem, as explained in (2) above.
\end{proof}

\section*{10b. Central limits}

At a more advanced level now, on the same topics, you have certainly heard about bell-shaped curves, and perhaps even observed them in physics or chemistry class, because any routine measurement leads to such curves. Mathematically, here is the question that we would like to solve, coming as a refined version of our previous Question 10.1:

\index{CLT}
\index{Central Limit theorem}
\index{normal variable}
\index{Gaussian variable}

\begin{question}
Given i.i.d. random variables $f_1,f_2,f_3,\ldots\,$, assumed to be centered, and having common variance $t>0$, do we have
$$\frac{f_1+\ldots+f_n}{\sqrt{n}}\sim g_t$$
in the $n\to\infty$ limit, for some bell-shaped density $g_t$? And, what is the formula of $g_t$? 
\end{question}

Obviously, this is something quite tricky, notably with the above $1/\sqrt{n}$ factor needing some explanations. So, here is the idea. In order for our sums $S_n=f_1+\ldots+f_n$ to have a chance to converge, we must arrange our normalizations as for the expectation $E(S_n)$ not to blow up, and as for the variance $V(S_n)$ to not blow up either. But this is impossible in general, because the expectation requirement leads to a $1/n$ normalization factor, and with this, we are into the law of large numbers, and away from bell-shaped curves.

\bigskip

So, what to do? Waive the expectation requirement, of course, by assuming that our variables are centered, $E(f_i)=0$, as said in Question 10.16. Now with this done, we still have to find the correct normalization factor, making the variance $V(S_n)$ not to blow up. And here, the good factor is indeed $1/\sqrt{n}$, as stated above, as shown by:

\begin{proposition}
Given i.i.d. random variables $f_1,f_2,f_3,\ldots\,$, assumed to be centered, and having common variance $t>0$, we have
$$V(f_1+\ldots+f_n)=nt$$
so in order for the partial sums $S_n=f_1+\ldots+f_n$ to have a chance to converge,
we must normalize them by a $1/\sqrt{n}$ factor, as stated in Question 10.16.
\end{proposition}

\begin{proof}
The variables $f_i$ being centered, so are their sums $S_n$, and we have the following computation, based on this, and on the independence of our variables:
\begin{eqnarray*}
V(S_n)
&=&E(S_n^2)\\
&=&\sum_{ij}E(f_if_j)\\
&=&\sum_iE(f_i^2)+\sum_{i\neq j}E(f_i)E(f_j)\\
&=&\sum_iE(f_i^2)\\
&=&nt
\end{eqnarray*}

Thus, we are in need indeed of a normalization factor $\alpha$. But, but repeating the computation with $S_n$ replaced by $\alpha S_n$, or simply by using $V(\alpha f)=\alpha^2V(f)$, we have:
$$V(\alpha S_n)=\alpha^2nt$$

Thus the good normalization factor is indeed $\alpha=1/\sqrt{n}$, as claimed.
\end{proof}

With Question 10.16 understood, what is next? Solving it I guess, and by using our linearization results for convolution from chapter 6, we first have:

\begin{theorem}
Given i.i.d. random variables $f_1,f_2,f_3,\ldots\,$, assumed to be centered, and having common variance $t>0$, we have
$$\frac{f_1+\ldots+f_n}{\sqrt{n}}\sim g_t$$
with $n\to\infty$, with $g_t$ being the law having $F(x)=e^{-tx^2/2}$ as Fourier transform.
\end{theorem}

\begin{proof}
There are several things going on here, the idea being as follows:

\medskip

(1) Observe first that in terms of moments, the Fourier transform of an arbitrary random variable $f:X\to\mathbb R$ is given by the following formula:
\begin{eqnarray*}
F_f(x)
&=&E(e^{ixf})\\
&=&E\left(\sum_{k=0}^\infty\frac{(ixf)^k}{k!}\right)\\
&=&\sum_{k=0}^\infty\frac{(ix)^kE(f^k)}{k!}\\
&=&\sum_{k=0}^\infty\frac{i^kM_k(f)}{k!}\,x^k
\end{eqnarray*}

(2) In particular, in the case of a centered variable, $E(f)=0$, as those that we are interested in, the Fourier transform formula that we get is as follows:
$$F_f(x)=1-\frac{M_2(f)}{2}\cdot x^2-i\,\frac{M_3(f)}{6}\cdot x^3+\ldots$$

Moreover, by further assuming that the Fourier variable is small, $x\simeq0$, the Fourier transform formula that we get, that we will use in what follows, becomes:
$$F_f(x)=1-\frac{M_2(f)}{2}\cdot x^2+O(x^2)$$

(3) In addition to this, we will also need to know what happens to the Fourier transform when rescaling. But the formula here is very easy to find, as follows:
$$F_{\alpha f}(x)
=E(e^{ix\alpha f})
=E(e^{i\alpha xf})
=F_f(\alpha x)$$

(4) Good news, we can now do our computation. By using the above formulae in (2) and (3), the Fourier transform of the variable in the statement is given by:
\begin{eqnarray*}
F(x)
&=&\left[F_f\left(\frac{x}{\sqrt{n}}\right)\right]^n\\
&=&\left[1-\frac{M_2(f)}{2}\cdot\frac{x^2}{n}+O(n^{-2})\right]^n\\
&=&\left[1-\frac{tx^2}{2n}+O(n^{-2})\right]^n\\
&\simeq&\left[1-\frac{tx^2}{2n}\right]^n\\
&\simeq&e^{-tx^2/2}
\end{eqnarray*}

(5) We are therefore led to the conclusion in the statement, modulo the fact that we do not know yet that a density $g_t$ having as Fourier transform $F(x)=e^{-tx^2/2}$ really exists, plus perhaps some other abstract issues to be discussed too. More on this later.
\end{proof}

Moving on, as a continuation, we have the following result:

\begin{theorem}
The ``normal'' law $g_t$, having as Fourier transform 
$$F(x)=e^{-tx^2/2}$$
must have all odd moments zero, and its even moments must be the numbers
$$M_k(g_t)=t^{k/2}\times k!!$$
where $k!!=(k-1)(k-3)(k-5)\ldots\,$, as usual. 
\end{theorem}

\begin{proof}
We must reformulate the equation $F(x)=e^{-tx^2/2}$, in terms of moments. We know from the proof of Theorem 10.18 that we have:
$$F(x)=\sum_{k=0}^\infty\frac{i^kM_k(g_t)}{k!}\,x^k$$

On the other hand, we have the following formula, for the exponential:
$$e^{-tx^2/2}
=\sum_{r=0}^\infty
(-1)^r\frac{t^rx^{2r}}{2^rr!}$$

We conclude that our equation $F(x)=e^{-tx^2/2}$ takes the following form:
$$\sum_{k=0}^\infty\frac{i^kM_k(g_t)}{k!}\,x^k=\sum_{r=0}^\infty
(-1)^r\frac{t^rx^{2r}}{2^rr!}$$

Thus the odd moments must vanish, and the even moments must be given by:
\begin{eqnarray*}
M_k(g_t)
&=&k!\times\frac{t^{k/2}}{2^{k/2}(k/2)!}\\
&=&t^{k/2}\times\frac{2\cdot3\cdot4\ldots (k-1)\cdot k}{2\cdot 4\cdot 6\ldots (k-2)\cdot k}\\
&=&t^{k/2}\times 3\cdot 5\ldots(k-3)(k-1)\\
&=&t^{k/2}\times k!!
\end{eqnarray*}

And with this, we are led to the conclusion in the statement.
\end{proof}

We can now finish our work, with some integration know-how, as follows:

\index{double factorial}

\begin{theorem}
The normal laws, appearing from central limits, are
$$g_t=\frac{1}{\sqrt{2\pi t}}\,e^{-x^2/2t}dx$$
with their mass $1$ property coming from the Gauss formula.
\end{theorem}

\begin{proof}
According to Theorem 10.19, the moments of the normal law $g_t$ must satisfy the following recurrence formula, with initial data $M_0=1$, $M_1=0$:
$$M_k=t(k-1)M_{k-2}$$

Now observe that by partial integration, we have the following formula:
\begin{eqnarray*}
\int_\mathbb Rx^ke^{-x^2/2t}dx
&=&\int_\mathbb R\frac{x^{k+1}}{k+1}\cdot\frac{x}{t}\cdot e^{-x^2/2t}dx\\
&=&\frac{1}{t(k+1)}\int_\mathbb Rx^{k+2}e^{-x^2/2t}dx
\end{eqnarray*}

Thus our integrals $I_k$ satisfy $I_{k+2}=t(k+1)I_k$, which is the same as $I_k=t(k-1)I_{k-2}$, same recurrence as for the numbers $M_k$. We also have $I_1=0$, because $e^{-x^2/2t}$ is an even function. So, it remains to suitably rescale our integrals, as to have $I_0=1$. But:
\begin{eqnarray*}
I_0
&=&\int_\mathbb Re^{-x^2/2t}\,dx\\
&=&\int_\mathbb Re^{-2ty^2/2t}\,\sqrt{2t}dy\\
&=&\sqrt{2t}\int_\mathbb Re^{-y^2}dy\\
&=&\sqrt{2\pi t}
\end{eqnarray*}

Thus, we are led to the formula of $g_t$ in the statement.
\end{proof}

\section*{10c. Normal variables}

Good work that we did in the above, I mean main question about limits solved, no question about this, but for better understanding, and for further reference too, it is customary at this point to make tabula rasa, and review our findings, in reverse order. 

\bigskip

So, let us get now into this, remake of the previous section, on the grounds that each good story is worth being told many times, with slight changes. Let us start with:

\begin{definition}
The normal law of parameter $1$ is the following measure:
$$g_1=\frac{1}{\sqrt{2\pi}}\,e^{-x^2/2}dx$$
More generally, the normal law of parameter $t>0$ is the following measure:
$$g_t=\frac{1}{\sqrt{2\pi t}}\,e^{-x^2/2t}dx$$
These are also called Gaussian distributions, with ``g'' standing for Gauss.
\end{definition}

Observe that the above laws have indeed mass 1, as they should. This follows indeed from the Gauss formula, which gives, with $x=\sqrt{2t}\,y$:
$$\int_\mathbb R e^{-x^2/2t}dx
=\int_\mathbb R e^{-y^2}\sqrt{2t}\,dy
=\sqrt{2t}\int_\mathbb R e^{-y^2}dy
=\sqrt{2\pi t}$$

Generally speaking, the normal laws appear as bit everywhere, in real life. The reasons behind this phenomenon come from the Central Limit Theorem (CLT), that we will explain in a moment, after developing some general theory. As a first result, we have:

\begin{proposition}
We have the variance formula
$$V(g_t)=t$$
valid for any $t>0$.
\end{proposition}

\begin{proof}
The first moment is 0, because our normal law $g_t$ is centered. As for the second moment, this can be computed by partial integration, as follows:
\begin{eqnarray*}
M_2
&=&\frac{1}{\sqrt{2\pi t}}\int_\mathbb Rx^2e^{-x^2/2t}dx\\
&=&\frac{1}{\sqrt{2\pi t}}\int_\mathbb R(tx)\left(-e^{-x^2/2t}\right)'dx\\
&=&\frac{1}{\sqrt{2\pi t}}\int_\mathbb Rte^{-x^2/2t}dx\\
&=&t
\end{eqnarray*}

We conclude from this that the variance is $V=M_2=t$, as stated.
\end{proof}

We can in fact compute all moments, by using the same trick, and we obtain:

\begin{theorem}
The even moments of the normal law are the numbers
$$M_k(g_t)=t^{k/2}\times k!!$$
where $k!!=(k-1)(k-3)(k-5)\ldots\,$, and the odd moments vanish. 
\end{theorem}

\begin{proof}
We have the following computation, valid for any integer $k\in\mathbb N$:
\begin{eqnarray*}
M_k
&=&\frac{1}{\sqrt{2\pi t}}\int_\mathbb Ry^ke^{-y^2/2t}dy\\
&=&\frac{1}{\sqrt{2\pi t}}\int_\mathbb R(ty^{k-1})\left(-e^{-y^2/2t}\right)'dy\\
&=&\frac{1}{\sqrt{2\pi t}}\int_\mathbb Rt(k-1)y^{k-2}e^{-y^2/2t}dy\\
&=&t(k-1)\times\frac{1}{\sqrt{2\pi t}}\int_\mathbb Ry^{k-2}e^{-y^2/2t}dy\\
&=&t(k-1)M_{k-2}
\end{eqnarray*}

Thus by recurrence, we are led to the formula in the statement.
\end{proof}

Regarding now the Fourier transform computation, this is as follows:

\begin{theorem}
We have the following formula, valid for any $t>0$:
$$F_{g_t}(x)=e^{-tx^2/2}$$
In particular, the normal laws satisfy $g_s*g_t=g_{s+t}$, for any $s,t>0$.
\end{theorem}

\begin{proof}
The Fourier transform formula can be established as follows:
\begin{eqnarray*}
F_{g_t}(x)
&=&\frac{1}{\sqrt{2\pi t}}\int_\mathbb Re^{-y^2/2t+ixy}dy\\
&=&\frac{1}{\sqrt{2\pi t}}\int_\mathbb Re^{-(y/\sqrt{2t}-\sqrt{t/2}ix)^2-tx^2/2}dy\\
&=&\frac{1}{\sqrt{2\pi t}}\int_\mathbb Re^{-z^2-tx^2/2}\sqrt{2t}dz\\
&=&\frac{1}{\sqrt{\pi}}e^{-tx^2/2}\int_\mathbb Re^{-z^2}dz\\
&=&\frac{1}{\sqrt{\pi}}e^{-tx^2/2}\cdot\sqrt{\pi}\\
&=&e^{-tx^2/2}
\end{eqnarray*}

As for the last assertion, this follows from the fact that $\log F_{g_t}$ is linear in $t$.
\end{proof}

Good news, we are now ready to state and prove the CLT, as follows:

\begin{theorem}[CLT]
Given i.i.d. random variables $f_1,f_2,f_3,\ldots\,$, assumed to be centered, and having common variance $t>0$, we have, with $n\to\infty$, in moments,
$$\frac{f_1+\ldots+f_n}{\sqrt{n}}\sim g_t$$
where $g_t$ is the Gaussian law of parameter $t$.
\end{theorem}

\begin{proof}
In terms of moments, the Fourier transform is given by:
$$F_f(x)
=E\left(\sum_{k=0}^\infty\frac{(ixf)^k}{k!}\right)
=\sum_{k=0}^\infty\frac{(ix)^kE(f^k)}{k!}
=\sum_{k=0}^\infty\frac{i^kM_k(f)}{k!}\,x^k$$

We conclude that the Fourier transform of the variable in the statement is:
$$F(x)
=\left[F_f\left(\frac{x}{\sqrt{n}}\right)\right]^n
\simeq\left[1-\frac{tx^2}{2n}\right]^n
\simeq e^{-tx^2/2}$$

But this latter function being the Fourier transform of $g_t$, we obtain the result.
\end{proof}

And with this, material from the previous section fully reviewed, and end of our logical loop. By the way, talking logical loops, cat just came back from her hunt, with a mouse in her claws, so now that she's home again, let's see what she has to say about this:

\begin{cat}
Just caught this young fellow, with a quantum logical loop. Guy must have been trained in pure mathematics or something, was not familiar with such tricks.
\end{cat}

Okay, thanks cat, we are in full agreement I guess, logical loops are what makes this world go round, and mastering them is essential for survival, and other good things.

\bigskip

Back to our work now, what is next? Normally, some more mathematics for the normal laws, in relation with their moments, as we usually do, when meeting a new law. And we first have here, coming as an interesting reformulation of Theorem 10.23:

\begin{proposition}
The moments of the normal law are the numbers
$$M_k(g_t)=t^{k/2}|P_2(k)|$$
where $P_2(k)$ is the set of pairings of $\{1,\ldots,k\}$.
\end{proposition}

\begin{proof}
Let us count the pairings of $\{1,\ldots,k\}$. In order to have such a pairing, we must pair $1$ with one of the numbers $2,\ldots,k$, and then use a pairing of the remaining $k-2$ numbers. Thus, we have the following recurrence formula:
$$|P_2(k)|=(k-1)|P_2(k-2)|$$

As for the initial data, this is $P_1=0$, $P_2=1$. Thus, we are led to the result.
\end{proof}

We are not done yet with moments, here being one more version of what we have:

\begin{theorem}
The moments of the normal law are the numbers
$$M_k(g_t)=\sum_{\pi\in P_2(k)}t^{|\pi|}$$
where $P_2(k)$ is the set of pairings of $\{1,\ldots,k\}$, and $|.|$ is the number of blocks.
\end{theorem}

\begin{proof}
This follows indeed from Proposition 10.27, because the number of blocks of a pairing of $\{1,\ldots,k\}$ is trivially $k/2$, independently of the pairing.
\end{proof}

As a comment here, the above result makes a link with the material from chapter 7 regarding the semicircle laws, and more specifically with the following formula for the moments, which can be easily deduced, a bit as above, from our formulae there:
$$M_k(\gamma_t)=\sum_{\pi\in NC_2(k)}t^{|\pi|}$$

In more advanced terms, we can say that $g_t\to\gamma_t$ is a liberation operation, a bit like $p_t\to\pi_t$ was, as explained in chapter 7. And, more on such things, later in this book.

\bigskip

Now back to our regular business, various moment manipulations, as a final result on the subject, complementing what we already have, let us formulate:

\begin{theorem}
The even normalized central moments of the normal law are
$$M_k''(g_t)=k!!$$
with $k!!=(k-1)(k-3)\ldots\,$ as usual, and the odd such moments vanish. Also,
$$E=0\quad,\quad V=t\quad,\quad \gamma=0\quad,\quad \kappa=3$$
are the expectation, variance, skewness and kurtosis.
\end{theorem}

\begin{proof}
The normal laws being centered, their central moments equal the plain ones, $M_k'=M_k$. Since the variance is $V=t$, the normalized central moments are given by:
$$M_k''=\frac{M'_k}{t^{k/2}}=\frac{M_k}{t^{k/2}}=\frac{\delta_{2|k}t^{k/2}\times k!!}{t^{k/2}}=\delta_{2|k}k!!$$

Finally, at $k=3$ this gives $\gamma=0$, and at $k=4$ this gives $\kappa=4!!=3$.
\end{proof}

And with this, end our story regarding the normal laws? Not exactly, because we have some interesting comments regarding the kurtosis, as follows:

\begin{comment}
Regarding the kurtosis of the normal laws, $\kappa=3$:
\begin{enumerate}
\item This explains, via variations of the CLT, why kurtoses often appear as $\kappa=3+\kappa'$,

\item With a similar phenomenon happening in the free case, namely $\kappa=2+\kappa'$,

\item And with the $3-2=1$ factor coming from the fact that $\cap\!\!\cap$ is crossing.
\end{enumerate}
\end{comment}

Which sounds quite good and conceptual, and we will end with this, and exercise of course for you, to learn more about such things. By the way, I would like to thank here my colleague Paul Dario from Cergy, for some help with this, when I was in need.

\section*{10d. Complex versions} 

Let us discuss now the complex analogues of all the above, with a notion of complex normal, or Gaussian law. To start with, we have the following definition:

\begin{definition}
A complex random variable is a variable $f:X\to\mathbb C$. In the discrete case, the law of such a variable is the complex probability measure
$$\mu=\sum_i\alpha_i\delta_{z_i}\quad,\quad\alpha_i\geq0\quad,\quad\sum_i\alpha_i=1\quad,\quad z_i\in\mathbb C$$
given by the following formula, with $P$ being the probability over $X$,
$$\mu=\sum_{z\in\mathbb C}P(f=z)\delta_z$$
with the sum being finite or countable, as per our discretness assumption. In general, the law of such a variable $f:X\to\mathbb C$ is a probability measure on the complex plane.
\end{definition}

Many things can be said, regarding the complex random variables, but going now straight to the point, namely probabilistic limiting theorems, let us discuss the complex analogue of the CLT. We have the following statement, to start with:

\begin{theorem}
Given complex random variables $f_1,f_2,f_3,\ldots$ whose real and imaginary parts are i.i.d., centered, and with common variance $t>0$, we have
$$\frac{f_1+\ldots+f_n}{\sqrt{n}}\sim C_t$$
with $n\to\infty$, in moments, where $C_t$ is the law of a complex variable whose real and imaginary parts are independent, and each following the law $g_t$.
\end{theorem}

\begin{proof}
This follows indeed from the real CLT, established in Theorem 10.25, simply by taking the real and imaginary parts of all the variables involved.
\end{proof}

It is tempting at this point to call Theorem 10.32 the complex CLT, but before doing that, let us study a bit more all this. We would like to have a better understanding of the limiting law $C_t$ at the end, and for this purpose, let us look at a sum as follows, with $a,b$ being real independent variables, both following the normal law $g_t$:
$$c=a+ib$$

To start with, this variable is centered, in a complex sense, because we have:
\begin{eqnarray*}
E(c)
&=&E(a+ib)\\
&=&E(a)+iE(b)\\
&=&0+i\cdot 0\\
&=&0
\end{eqnarray*}

However, regarding the variance, things are more complicated, because the usual variance formula from the real case, which is $V(c)=E(c^2)$ in the centered case, will not provide us with a positive number, in the case where our variable is not real. So, in order to have a variance which is real, and positive too, we must rather use a formula of type $V(c)=E(|c|^2)$, in the complex centered case. And, with this convention for the variance, we have then the following computation, for the variance of the above variable $c$:
\begin{eqnarray*}
V(c)
&=&E(|c|^2)\\
&=&E(a^2+b^2)\\
&=&E(a^2)+E(b^2)\\
&=&V(a^2)+V(b^2)\\
&=&t+t\\
&=&2t
\end{eqnarray*}

But this suggests to divide everything by $\sqrt{2}$, as to have in the end a variable having complex variance $t$, in our sense, and we are led in this way into:

\index{complex normal law}
\index{complex Gaussian law}

\begin{definition}
The complex normal, or Gaussian law of parameter $t>0$ is
$$G_t=law\left(\frac{a+ib}{\sqrt{2}}\right)$$
where $a,b$ are real and independent, each following the law $g_t$.
\end{definition}

In short, the complex normal laws appear as natural complexifications of the real normal laws. As in the real case, these measures form convolution semigroups:

\index{convolution semigroup}

\begin{proposition}
The complex Gaussian laws have the property
$$G_s*G_t=G_{s+t}$$
for any $s,t>0$, and so they form a convolution semigroup.
\end{proposition}

\begin{proof}
This follows indeed from the real result, namely $g_s*g_t=g_{s+t}$, established in Theorem 10.24, simply by taking real and imaginary parts.
\end{proof}

We have as well the following complex analogue of the CLT:

\index{CCLT}
\index{Complex CLT}
\index{Complex Central Limit Theorem}
\index{complex variables}
\index{complex normal law}

\begin{theorem}[CCLT]
Given complex random variables $f_1,f_2,f_3,\ldots$ whose real and imaginary parts are i.i.d. and centered, and having variance $t>0$ as above,
$$\frac{f_1+\ldots+f_n}{\sqrt{n}}\sim G_t$$
with $n\to\infty$, in moments.
\end{theorem}

\begin{proof}
This follows indeed from our previous CCLT result, the one from Theorem 10.32, by dividing everything there by $\sqrt{2}$, as explained in the above.
\end{proof}

Regarding now the moments, the situation here is more complicated than in the real case, because in order to have good results, we have to deal with both the complex variables, and their conjugates. Let us formulate the following definition:

\index{colored integers}
\index{colored moments}

\begin{definition}
The moments a complex variable $f\in L^\infty(X)$ are the numbers
$$M_k=E(f^k)$$
depending on colored integers $k=\circ\bullet\bullet\circ\ldots\,$, with the conventions
$$f^\emptyset=1\quad,\quad f^\circ=f\quad,\quad f^\bullet=\bar{f}$$
and multiplicativity, in order to define the colored powers $f^k$.
\end{definition}

As an illustration for this notion, which is something very intuitive, here are the formulae of the four possible order 2 moments of a complex variable $f$:
$$M_{\circ\circ}=E(f^2)\quad,\quad M_{\circ\bullet}=E(f\bar{f})$$
$$M_{\bullet\circ}=E(\bar{f}f)\quad,\quad M_{\bullet\bullet}=E(\bar{f}^2)$$

Observe that, since $f,\bar{f}$ commute, we have the following identity, which shows that there is a bit of redundancy in our above definition, as formulated:
$$M_{\circ\bullet}=M_{\bullet\circ}$$

In fact, again since $f,\bar{f}$ commute, we can permute terms, in the general context of Definition 10.36, and restrict the attention to exponents of the following type:
$$k=\ldots\circ\circ\circ\bullet\bullet\bullet\bullet\ldots$$

However, the point is that our results to come regarding the complex Gaussian laws, and some other complex laws too, later on, not to talk about laws of matrices, random matrices and other noncommuting variables, that will appear later too, in this book, will look better without doing this. So, we will use Definition 10.36 as stated. 

\bigskip

Getting to work now, we have the following result, regarding the normal laws:

\begin{theorem}
The moments of the complex normal law are given by
$$M_k(G_t)=\begin{cases}
t^pp!&(k\ {\rm uniform, of\ length}\ 2p)\\
0&(k\ {\rm not\ uniform})
\end{cases}$$
where $k=\circ\bullet\bullet\circ\ldots$ is called uniform when it contains the same number of $\circ$ and $\bullet$.
\end{theorem}

\begin{proof}
We must compute the moments, with respect to colored integer exponents $k=\circ\bullet\bullet\circ\ldots$ as above, of the variable from Definition 10.33, namely:
$$f=\frac{a+ib}{\sqrt{2}}$$

We can assume that we are in the case $t=1$, and the proof here goes as follows:

\medskip

(1) As a first observation, in the case where our exponent $k=\circ\bullet\bullet\circ\ldots$ is not uniform, a standard rotation argument shows that the corresponding moment of $f$ vanishes. To be more precise, the variable $f'=wf$ is complex Gaussian too, for any complex number $w\in\mathbb T$, and from $M_k(f)=M_k(f')$ we obtain $M_k(f)=0$, in this case.

\medskip

(2) In the uniform case now, where the exponent $k=\circ\bullet\bullet\circ\ldots$ consists of $p$ copies of $\circ$ and $p$ copies of $\bullet$\,, the corresponding moment can be computed as follows:
\begin{eqnarray*}
M_k
&=&E\left((f\bar{f})^p\right)\\
&=&\frac{1}{2^p}E\left((a^2+b^2)^p\right)\\
&=&\frac{1}{2^p}\sum_r\binom{p}{r}E(a^{2r})E(b^{2p-2r})\\
&=&\frac{1}{2^p}\sum_r\binom{p}{r}(2r)!!(2p-2r)!!\\
&=&\frac{1}{2^p}\sum_r\frac{p!}{r!(p-r)!}\cdot\frac{(2r)!}{2^rr!}\cdot\frac{(2p-2r)!}{2^{p-r}(p-r)!}\\
&=&\frac{p!}{4^p}\sum_r\binom{2r}{r}\binom{2p-2r}{p-r}
\end{eqnarray*}

(3) In order to finish now the computation, let us recall that we have the following formula, coming from the generalized binomial formula, or from the Taylor formula:
$$\frac{1}{\sqrt{1+t}}=\sum_q\binom{2q}{q}\left(\frac{-t}{4}\right)^q$$

By taking the square of this series, we obtain the following formula:
$$\frac{1}{1+t}
=\sum_p\left(\frac{-t}{4}\right)^p\sum_r\binom{2r}{r}\binom{2p-2r}{p-r}$$

Now by looking at the coefficient of $t^p$ on both sides, we conclude that the sum on the right equals $4^p$. Thus, we can finish the moment computation in (2), as follows:
$$M_k=\frac{p!}{4^p}\times 4^p=p!$$

We are therefore led to the conclusion in the statement.
\end{proof}

In general, regarding the complex normal laws, what we have in Theorem 10.37 about them is usually what is needed in practice, when dealing with moments. But, as before with the real Gaussian laws, or even further before with the Poisson and Bessel laws, a better-looking statement regarding the moments is in terms of partitions. 

\bigskip

In order to explain this, let us make the following convention, which is standard:

\begin{definition}
Given a colored integer $k=\circ\bullet\bullet\circ\ldots\,$, we denote by
$$\mathcal P_2(k)\subset P_2(k)$$
the set of matching pairings, that is, pairings pairing $\circ-\bullet$ symbols. 
\end{definition}

With this convention, we can reformulate Theorem 10.37 into something better, in tune with what we previously had regarding some other laws, as follows:

\index{matching pairings}

\begin{theorem}
The moments of the complex normal law are the numbers
$$M_k(G_t)=\sum_{\pi\in\mathcal P_2(k)}t^{|\pi|}$$
where $\mathcal P_2(k)$ are the matching pairings of $\{1,\ldots,k\}$, and $|.|$ is the number of blocks.
\end{theorem}

\begin{proof}
This is a standard reformulation of Theorem 10.37, as follows:

\medskip

(1) In the case $t=1$, we know from Theorem 10.37 that the moments are:
$$M_k=\begin{cases}
(|k|/2)!&(k\ {\rm uniform})\\
0&(k\ {\rm not\ uniform})
\end{cases}$$

(2) On the other hand, the numbers $|\mathcal P_2(k)|$ are given by exactly the same formula. Indeed, in order to have a matching pairing of $k$, our exponent $k=\circ\bullet\bullet\circ\ldots$ must be uniform, consisting of $p$ copies of $\circ$ and $p$ copies of $\bullet$, with $p=|k|/2$. 

\medskip

(3) But then the matching pairings of $k$ correspond to the permutations of the $\bullet$ symbols, as to be matched with $\circ$ symbols, and so we have $p!$ such pairings. Thus, we have the same formula as for the moments of $f$, and we are led to the conclusion in the statement. As for the extension to the general case $t>0$, this is straightforward.
\end{proof}

In practice, we also need to know how to compute joint moments. We have here:

\index{joint moments}
\index{Wick formula}

\begin{theorem}[Wick formula]
Given independent variables $f_i$, each following the complex normal law $G_t$, with $t>0$ being a fixed parameter, we have the formula
$$E\left(f_{i_1}^{k_1}\ldots f_{i_s}^{k_s}\right)=t^{s/2}\#\left\{\pi\in\mathcal P_2(k)\Big|\pi\leq\ker i\right\}$$
where $k=k_1\ldots k_s$ and $i=i_1\ldots i_s$, for the joint moments of these variables, where $\pi\leq\ker i$ means that the indices of $i$ must fit into the blocks of $\pi$, in the obvious way.
\end{theorem}

\begin{proof}
This is something well-known, which can be proved as follows:

\medskip

(1) Let us first discuss the case where we have a single variable $f$, which amounts in taking $f_i=f$ for any $i$ in the formula in the statement. What we have to compute here are the moments of $f$, with respect to colored integer exponents $k=\circ\bullet\bullet\circ\ldots\,$, and the formula in the statement tells us that these moments must be:
$$E(f^k)=t^{|k|/2}|\mathcal P_2(k)|$$

But this is the formula in Theorem 10.39, so we are done with this case.

\medskip

(2) In general now, when expanding the product $f_{i_1}^{k_1}\ldots f_{i_s}^{k_s}$ and rearranging the terms, we are left with doing a number of computations as in (1), and then making the product of the expectations that we found. But this amounts in counting the partitions in the statement, with the condition $\pi\leq\ker i$ there standing for the fact that we are doing the various type (1) computations independently, and then making the product.
\end{proof}

As a variation of the Wick formula, which is quite useful in practice, we have:

\index{Wick formula}

\begin{theorem}[Wick formula 2]
Given independent variables $f_i$, each following the complex normal law $G_t$, with $t>0$ being a fixed parameter, we have the formula
$$E\left(f_{i_1}\ldots f_{i_k}f_{j_1}^*\ldots f_{j_k}^*\right)=t^k\#\left\{\pi\in S_k\Big|i_{\pi(r)}=j_r,\forall r\right\}$$
for the non-vanishing joint moments of these variables.
\end{theorem}

\begin{proof}
With some changes in the indices and notations, the Wick formula reads:
$$E\left(f_{I_1}^{K_1}\ldots f_{I_s}^{K_s}\right)=t^{s/2}\#\left\{\sigma\in\mathcal P_2(K)\Big|\sigma\leq\ker I\right\}$$

Now observe that we have $\mathcal P_2(K)=\emptyset$, unless the colored integer $K=K_1\ldots K_s$ is uniform, in the sense that it contains the same number of $\circ$ and $\bullet$ symbols. Up to permutations, the non-trivial case, where the moment is non-vanishing, is the case where the colored integer $K=K_1\ldots K_s$ is of the following special form:
$$K=\underbrace{\circ\circ\ldots\circ}_k\ \underbrace{\bullet\bullet\ldots\bullet}_k$$

So, let us focus on this case, which is the non-trivial one. Here we have $s=2k$, and we can write the multi-index $I=I_1\ldots I_s$ in the following way:
$$I=i_1\ldots i_k\ j_1\ldots j_k$$

With these changes made, the above usual Wick formula reads:
$$E\left(f_{i_1}\ldots f_{i_k}f_{j_1}^*\ldots f_{j_k}^*\right)=t^k\#\left\{\sigma\in\mathcal P_2(K)\Big|\sigma\leq\ker(ij)\right\}$$

The point now is that the matching pairings $\sigma\in\mathcal P_2(K)$, with $K=\circ\ldots\circ\bullet\ldots\bullet\,$, of length $2k$, as above, correspond to the permutations $\pi\in S_k$, in the obvious way. With this identification made, the above modified usual Wick formula becomes:
$$E\left(f_{i_1}\ldots f_{i_k}f_{j_1}^*\ldots f_{j_k}^*\right)=t^k\#\left\{\pi\in S_k\Big|i_{\pi(r)}=j_r,\forall r\right\}$$

Thus, we have reached to the formula in the statement, and we are done.
\end{proof}

Finally, here is one more formulation of the Wick formula, useful as well:

\index{Wick formula}

\begin{theorem}[Wick formula 3]
Given independent variables $f_i$, each following the complex normal law $G_t$, with $t>0$ being a fixed parameter, we have the formula
$$E\left(f_{i_1}f_{j_1}^*\ldots f_{i_k}f_{j_k}^*\right)=t^k\#\left\{\pi\in S_k\Big|i_{\pi(r)}=j_r,\forall r\right\}$$
for the non-vanishing joint moments of these variables.
\end{theorem}

\begin{proof}
This follows from our second Wick formula, from Theorem 10.41, simply by permuting the terms, as to have an alternating sequence of plain and conjugate variables. Alternatively, we can start with Theorem 10.40, and then perform the same manipulations as in the proof of Theorem 10.41, but with the exponent being this time as follows: 
$$K=\underbrace{\circ\bullet\circ\bullet\ldots\ldots\circ\bullet}_{2k}$$

Thus, we are led to the conclusion in the statement.
\end{proof}

And with this, end of our learning on the complex normal variables. The above formulae are very useful when dealing with random matrices, and more on this, later.

\section*{10e. Exercises}

This was a quite standard chapter, all classical results, and as exercises, we have:

\begin{exercise}
Find some counterexamples, for the 3 notions of convergence.
\end{exercise}

\begin{exercise}
Try recovering $g_t$ from its moments, via Stieltjes inversion.
\end{exercise}

\begin{exercise}
Find a genius method for dealing with the Gauss integral.
\end{exercise}

\begin{exercise}
Perform some real-life verifications, of your choice, for the CLT.
\end{exercise}

\begin{exercise}
Do as well some computer simulations for the CLT.
\end{exercise}

\begin{exercise}
How to best integrate $g_t$, over a finite interval $[a,b]$?
\end{exercise}

\begin{exercise}
Learn more about complex variables, theory and examples.
\end{exercise}

\begin{exercise}
Get familiar with matching pairings, and the Wick formula.
\end{exercise}

As bonus exercise, with the Wick formula learned, have a look at random matrices.

\chapter{Rayleigh laws}

\section*{11a. Shifted bells}

We have seen in the previous chapter the basic theory of the normal laws, real and complex, and this is certainly what is essentially needed, for dealing with various continuous probabilistic phenomena, originally appearing from central limits.

\bigskip

At a more advanced level, however, far more things can be said, the idea being that the real normal laws $g_t$ and their complex versions $G_t$ are surrounded by a myriad of other interesting laws, appearing as variations of them, which must be mastered as well. We will explore this ``normal law galaxy'' in this chapter, and in the next one too.

\bigskip

To start with, let us talk a bit more about Gaussian bells. The ones $g_t$ that we studied in the previous chapter appeared in a purely mathematical way, from central limits, and as the word ``central'' indicates, they are of course centered. However, this was theory, and in practice, say coming via various real-life measurements, the Gaussian bells quite often appear shifted by a parameter $a\in\mathbb R$, I mean with symmetry axis there at $a$. 

\bigskip

So, we must study these versions too, and regarding them, we have:

\index{shifted normal law}

\begin{theorem}
The normal law $g_t$, shifted to the right by $a\in\mathbb R$, has density
$$g_t^a=\frac{1}{\sqrt{2\pi t}}\,e^{-(x-a)^2/2t}\,dx$$
and the basic properties of this shifted law $g_t^a$ are as follows:
\begin{enumerate}
\item The mean is $E=a$, the variance is $V=t$.

\item The moments satisfy $M_k=aM_{k-1}+(k-1)tM_{k-2}$.

\item The central moments are $M_k'=\delta_{2|k}t^{k/2}k!!$.

\item The normalized central moments are $M_k''=\delta_{2|k}k!!$.

\item The skewness is $\gamma=0$, the kurtosis is $\kappa=3$.

\item The Fourier transform is $F(x)=\exp(iax-tx^2/2)$.

\item We have the convolution formula $g_s^a*g_t^b=g_{s+t}^{a+b}$.
\end{enumerate}
\end{theorem}

\begin{proof}
These are things that we know well at $a=0$, from the previous chapter, and the general case basically follows from this, with a few computations added. To start with, let us draw a picture, that will help us understanding what is going on:
$$
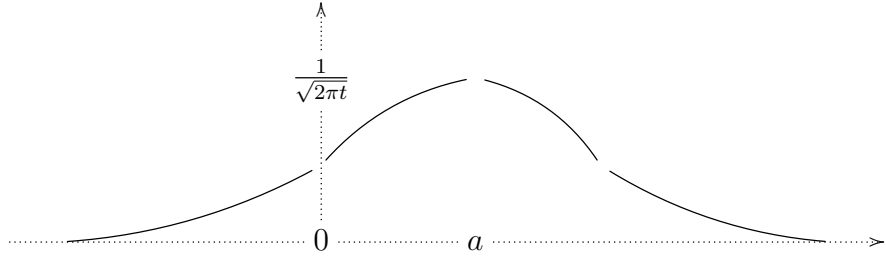
\xymatrix@R=18pt@C=15pt{
&&&&&\\
&&&&&\frac{1}{\sqrt{2\pi t}}\ar@.[u]&&\ar@/^/@{-}[drr]\ar@/_/@{-}[dll]&\\
&&&&&\ar@/^/@{-}[dllll]&&&&\ar@/_/@{-}[drrrr]&\\
\ar@{.}[rrrrr]&&&&&0\ar@{.}[uu]\ar@{.}[rr]&&a\ar@.[rrrrrrr]&&&&&&&}$$

(1) Regarding the mean, that is clearly $E=a$ on the picture, because the symmetry axis lies there, at $a\in\mathbb R$. Mathematically, we can deduce this from the following fact, which extends the basic observation that if the density is even, the mean is zero:
$$\varphi(a-x)=\varphi(a+x)\implies\int_\mathbb Rx\varphi(x)dx=a$$

To be more precise, this latter implication can be established as follows:
\begin{eqnarray*}
\int_\mathbb Rx\varphi(x)dx
&=&\int_{-\infty}^ax\varphi(x)dx+\int_a^\infty x\varphi(x)dx\\
&=&\int_0^\infty (a-y)\varphi(a-y)dy+\int_0^\infty(a+y)\varphi(a+y)dy\\
&=&a-\int_0^\infty y\varphi(a-y)dy+\int_0^\infty y\varphi(a+y)dy\\
&=&a
\end{eqnarray*}

Regarding now the variance, again we can trick here, this time with a shifting argument, but the simplest is perhaps just to compute. By partial integration, we have:
\begin{eqnarray*}
M_2-aM_1
&=&\frac{1}{\sqrt{2\pi t}}\int_\mathbb R(x^2-ax)e^{-(x-a)^2/2t}dx\\
&=&\frac{1}{\sqrt{2\pi t}}\int_\mathbb R(tx)\left(-e^{-(x-a)^2/2t}\right)'dx\\
&=&\frac{1}{\sqrt{2\pi t}}\int_\mathbb Rte^{-(x-a)^2/2t}dx\\
&=&t
\end{eqnarray*}

Thus $M_2=a^2+t$, and the variance is $V=(a^2+t)-a^2=t$, as stated.

\medskip

(2) Getting now to the higher moments of $g_t^a$, we can certainly compute them by partial integration, a bit as above, with the next computation being as follows:
\begin{eqnarray*}
M_3-aM_2
&=&\frac{1}{\sqrt{2\pi t}}\int_\mathbb R(x^3-ax^2)e^{-(x-a)^2/2t}dx\\
&=&\frac{1}{\sqrt{2\pi t}}\int_\mathbb R(tx^2)\left(-e^{-(x-a)^2/2t}\right)'dx\\
&=&\frac{1}{\sqrt{2\pi t}}\int_\mathbb R2tx\cdot e^{-(x-a)^2/2t}dx\\
&=&2tM_1
\end{eqnarray*}

Indeed, we conclude from this that the third moment is given by:
$$M_3=aM_2+2tM_1=a(a^2+t)+2at=a^3+3at$$

Next, it is pretty much clear that the same trick will work in general. And this is indeed the case, and we obtain the following recurrence, for our moments:
\begin{eqnarray*}
M_k-aM_{k-1}
&=&\frac{1}{\sqrt{2\pi t}}\int_\mathbb R(x^k-ax^{k-1})e^{-(x-a)^2/2t}dx\\
&=&\frac{1}{\sqrt{2\pi t}}\int_\mathbb R(tx^{k-1})\left(-e^{-(x-a)^2/2t}\right)'dx\\
&=&\frac{1}{\sqrt{2\pi t}}\int_\mathbb R(k-1)tx^{k-2}\cdot e^{-(x-a)^2/2t}dx\\
&=&(k-1)tM_{k-2}
\end{eqnarray*}

In practice, here is the list of the first few moments, obtained in this way:
$$M_1=a$$
$$M_2=a^2+t$$
$$M_3=a^3+3at$$
$$M_4=a^4+6a^2t+3t^2$$
$$M_5=a^5+10a^3t+15at^2$$
$$M_6=a^6+15a^4t+45at^2+15t^3$$
$$M_7=a^7+21a^5t+105a^3t^2+105at^3$$
$$M_8=a^8+28a^6t+210a^4t^2+420at^3+105t^4$$
$$\vdots$$

Which is nice, but since the coefficients do not look very good, we will not go further this way. More on this later in this book, when discussing orthogonal polynomials.

\medskip

(3) Instead, let us look at the central moments. And here, fortunately, we obtain the same moments as for $g_t$, as shown by the following computation:
\begin{eqnarray*}
M_k'(g_t^a)
&=&\frac{1}{\sqrt{2\pi t}}\int_\mathbb R(x-a)^ke^{-(x-a)^2/2t}dx\\
&=&\frac{1}{\sqrt{2\pi t}}\int_\mathbb Rx^ke^{-x^2/2t}dx\\
&=&M_k(g_t)
\end{eqnarray*}

Thus, by using now the formula from chapter 10 for the moments of $g_t$, we have:
$$M_k'(g_t^a)=\delta_{2|k}t^{k/2}k!!$$

(4) Next, since the central moments are the same for $g_t^a$ and $g_t$, and since the variances are also the same, the normalized central moments must all coincide:
$$M_k''(g^t_a)=M_k''(g_t)$$

(5) In particular, we deduce from this formula that the skewness is $\gamma=0$, and that the kurtosis is $\kappa=3$, exactly as for the laws $g_t$ in chapter 10, as claimed.

\medskip

(6) Regarding now the Fourier transform computation, we have:
\begin{eqnarray*}
F_{g_t^a}(x)
&=&\frac{1}{\sqrt{2\pi t}}\int_\mathbb Re^{-(y-a)^2/2t+ixy}dy\\
&=&\frac{1}{\sqrt{2\pi t}}\int_\mathbb Re^{-y^2/2t+ix(y+a)}dy\\
&=&\frac{e^{iax}}{\sqrt{2\pi t}}\int_\mathbb Re^{-y^2/2t+ixy}dy\\
&=&e^{iax}F_{g_t}(x)
\end{eqnarray*}

Now by using the formula from chapter 10 for the Fourier transform of $g_t$, we get:
$$F_{g_t^a}(x)=e^{iax}e^{-tx^2/2}=e^{iax-tx^2/2}$$
 
(7) Observe now that this Fourier transform that we just computed satisfies:
\begin{eqnarray*}
F_{g_s^a}(x)F_{g_t^b}(x)
&=&e^{iax-sx^2/2}e^{ibx-tx^2/2}\\
&=&e^{i(a+b)x-(s+t)x^2/2}\\
&=&F_{g_{s+t}^{a+b}}(x)
\end{eqnarray*}

Thus, at the level of the corresponding laws, we have the following formula:
$$g_s^a*g_t^b=g_{s+t}^{a+b}$$

Summarizing, theorem proved, but as a bonus, let us try as well to establish the formula $g_s^a*g_t^b=g_{s+t}^{a+b}$ directly, I mean without Fourier, with this being an instructive computation. We have, based on the definition of the convolution operation:
\begin{eqnarray*}
&&g_s^a*g_t^b\\
&=&\left(\frac{1}{\sqrt{2\pi s}}\,e^{-(x-a)^2/2s}\,dx\right)*
\left(\frac{1}{\sqrt{2\pi t}}\,e^{-(x-b)^2/2t}\,dx\right)\\
&=&\frac{1}{2\pi\sqrt{st}}\left(e^{-(x-a)^2/2s}\,dx\right)*\left(e^{-(x-b)^2/2t}\,dx\right)\\
&=&\frac{1}{2\pi\sqrt{st}}\left(\int_\mathbb Re^{-(x-y-a)^2/2s}e^{-(y-b)^2/2t}dy\right)dx\\
&=&\frac{1}{2\pi\sqrt{st}}\left(\int_\mathbb R\exp\left[-\frac{(x-y-a)^2}{2s}-\frac{(y-b)^2}{2t}\right]dy\right)dx\\
&=&\frac{1}{2\pi\sqrt{st}}\left(\int_\mathbb R\exp\left[-\frac{t(x-y-a)^2+s(y-b)^2}{2st}\right]dy\right)dx\\
&=&\frac{1}{2\pi\sqrt{st}}\left(\int_\mathbb R\exp\left[-\frac{(s+t)y^2+2[t(a-x)-sb]y+t(x-a)^2+sb^2}{2st}\right]dy\right)dx\\
&=&\frac{1}{2\pi\sqrt{st}}\left(\int_\mathbb R\exp\left[-\frac{1}{2st}\left(\left[\sqrt{s+t}\cdot y+\frac{t(a-x)-sb}{\sqrt{s+t}}\right]^2+\frac{st(x-a-b)^2}{s+t}\right)\right]dy\right)dx\\
&=&\frac{1}{2\pi\sqrt{st}}\left(\int_\mathbb R\exp\left[-\frac{1}{2st}\left(z^2+\frac{st(x-a-b)^2}{s+t}\right)\right]\frac{dz}{\sqrt{s+t}}\right)dx\\
&=&\frac{1}{2\pi\sqrt{st(s+t)}}\left(e^{-(x-a-b)^2/2(s+t)}\int_\mathbb Re^{-z^2/2st}dz\right)dx\\
&=&\frac{1}{2\pi\sqrt{st(s+t)}}e^{-(x-a-b)^2/2(s+t)}\sqrt{2\pi st}\,dx\\
&=&\frac{1}{\sqrt{2\pi(s+t)}}\,e^{-(x-a-b)^2/2(s+t)}\,dx\\
&=&g_{s+t}^{a+b}
\end{eqnarray*}

Cool stuff, hope you agree with me. This is actually my favorite computation in mathematics, when my low year students get nervous, for one reason or another, I give them this, and unless someone has heard about Fourier before, $1/2$ hour of silence.
\end{proof}

And with this, end of our study of the shifted normal laws $g_t^a$. Needless to say, the above computations are very useful, when dealing with real-life probability questions.

\bigskip

As a conclusion now to what we did, a bit philosophical, focusing on the difficulties that we met, rather than on our successes, let us formulate things as follows:

\begin{conclusion}
The theory of $g_t^a$ is quite straightforward, save for:
\begin{enumerate}
\item The computation of the CDF, which amounts in fighting with $\int e^{-x^2}$.

\item The moments at $a\neq0$, most likely a hypergeometric function business.
\end{enumerate}
\end{conclusion}

And with of course, exercise for you to learn more, about all this. In fact, thinking a bit, for various real-life problems involving the laws $g_t^a$, the CDF is what you need.

\section*{11b. Chi squared}

Getting now to our main objective in this chapter, namely talking about the Rayleigh laws, and their variations, let us formulate the following key definition:

\index{chi}
\index{chi squared}
\index{Rayleigh law}

\begin{definition}
The chi squared law, depending on $n\in\mathbb N$, is given by
$$\chi_n^2=law\Big(f_1^2+\ldots+f_n^2\Big)$$
with $f_1,\ldots,f_n$ being independent, each following $g_1$. We can talk as well about
$$\chi_n=law\left(\sqrt{f_1^2+\ldots+f_n^2}\right)$$
called chi law, whose simplest $n>1$ instance, $\chi_2$, is called Rayleigh law.
\end{definition}

So, these will be our objects of interest, in this chapter. As a first observation, all these laws basically come from $\chi_2$, which itself comes from Lord Rayleigh from physics, so expect in what follows, besides our usual calculus, an interesting mixture of geometry and physics. Let us also mention that all these laws are extremely useful as well in statistics. In short, good work that we will be doing here, for all kinds of purposes.

\bigskip

As a second comment, we already met in fact the Rayleigh law in chapter 10, when discussing the complex normal variables. However, our study here will be mostly real, and this because the chi and chi squared distributions, while certainly interpretable in terms of complex numbers at $n\in 2\mathbb N$, are defined in odd dimensions $n\in2\mathbb N+1$ too.

\bigskip

Getting now to work, what we have in Definition 11.3 is quite interesting, but as a question that you might have, and that I have too, where are all the parameters gone? I mean, we can certainly add a few useful parameters there, as summarized in:

\begin{remark}
We can add parameters to the $\chi_n$ and $\chi_n^2$ laws, as follows:
\begin{enumerate}
\item A parameter $t>0$, by assuming that $f_1,\ldots,f_n$ follow $g_t$.

\item An extra parameter $a\in\mathbb R$, assuming that $f_1,\ldots,f_n$ follow $g_t^a$.

\item Or full house, $t_1,\ldots,t_n>0$ and $a_1,\ldots,a_n\in\mathbb R$, assuming $f_i\sim g_{t_i}^{a_i}$.
\end{enumerate}
\end{remark}

So, shall we do this, ot not? Which is not an easy question, and with cat gone, although I doubt that she could have been useful here, I mean I know her well, she just goes for the kill, not bothering much with parameters and tricks, I will have to ask the rat. And rat, that we haven't met since chapter 5, philosopher as usual, answers:

\begin{rat}
The parameters can be both your strength, and your weakness. Go your own way, and peace upon you.
\end{rat}

Well, thanks rat, I was sort of expecting this, knowing you a bit, but frankly, what you say is a bit too philosophical. I have this feeling that with your modest lifestyle, hiding under the stove, and eating a few leftovers, you reached to a very good understanding of life and mathematics. But why not sharing a bit more your knowledge, with less modest animals, like me and cat? Although, we perhaps don't really deserve this. Guess at some point we should do the same as you, hide under the stove, and eat a few leftovers.

\bigskip

Nevermind. So going now for the advice Rat 11.5 as stated, a bit vague, and knowing myself a bit, I mean not that good at complex mathematics, and presuming that you reader are in a similar situation, we will go with Definition 11.3 as stated. As for generalizations, involving further parameters, we can see later, once the basics done.

\bigskip

Time now for some calculus? As a first and main result regarding $\chi_n^2$, we have:

\begin{theorem}
The density of the chi squared law is given by
$$\chi_n^2=\frac{x^{n/2-1}e^{-x/2}}{2^{n/2}\Gamma(n/2)}\,dx$$
with $\Gamma$ being as usual the gamma function.
\end{theorem}

\begin{proof}
This is something quite tricky, the idea being as follows:

\medskip

(1) To start with, and standing as a useful complement to the above, let us recall from chapter 8 that the gamma function is given by the following formula:
$$\Gamma(s)=\int_0^\infty x^{s-1}e^{-x}\,dx$$

By partial integration we have $\Gamma(s+1)=s\Gamma(s)$, and since $\Gamma(1)=1$, trivially, and $\Gamma(1/2)=\sqrt{\pi}$, obtained via $x=y^2$, we have the following formulae, for $N\in\mathbb N$:
$$\Gamma(N)=(N-1)!\quad,\quad\Gamma\left(N+\frac{1}{2}\right)=\frac{(2N)!!}{2^N}\,\sqrt{\pi}$$

Alternatively, we have the following formula, with $c=\sqrt{2},\sqrt{\pi}$ for $n$ even, odd:
$$\Gamma\left(\frac{n}{2}\right)=\frac{(n-1)!!}{2^{(n-1)/2}}\,c$$

And for more on all this, gamma function and its properties, we refer to chapter 8.

\medskip

(2) Getting now to what our statement says, let us first have a look at the case $n=1$. Here our law is $\chi_1^2=law(f^2)$, with $f\sim g_1$, and the predicted density is:
$$\chi_1^2=\frac{x^{-1/2}e^{-x/2}}{\sqrt{2}\cdot\sqrt{\pi}}\,dx$$

And, is this correct or not? To be more precise, we must prove that we have:
$$f\sim\frac{e^{-x^2/2}}{\sqrt{2\pi}}\,dx\ \implies\ f^2\sim\frac{x^{-1/2}e^{-x/2}}{\sqrt{2\pi}}\,dx$$

In order to do so, we can use moments. Indeed, we have the following formula:
$$M_k(f^2)=M_{2k}(f)=(2k)!!$$

On the other hand, the moments of the predicted density are as follows:
\begin{eqnarray*}
I_k
&=&\frac{1}{\sqrt{2\pi}}\int_0^\infty x^{k-1/2}e^{-x/2}dx\\
&=&\frac{1}{\sqrt{2\pi}}\int_0^\infty\left(k-\frac{1}{2}\right)x^{k-3/2}\cdot 2e^{-x/2}dx\\
&=&(2k-1)\cdot\frac{1}{\sqrt{2\pi}}\int_0^\infty x^{k-3/2}e^{-x/2}dx\\
&=&(2k-1)I_{k-1}
\end{eqnarray*}

As for the initial value, for this recurrence, this is given by the following formula, which by the way shows that we have indeed a probability measure, as we should:
\begin{eqnarray*}
I_0
&=&\frac{1}{\sqrt{2\pi}}\int_0^\infty x^{-1/2}e^{-x/2}dx\\
&=&\frac{1}{\sqrt{2\pi}}\int_0^\infty y^{-1}e^{-y^2/2}\,2y\,dy\\
&=&\sqrt{\frac{2}{\pi}}\int_0^\infty e^{-y^2/2}dy\\
&=&\sqrt{\frac{2}{\pi}}\times\frac{\sqrt{2\pi}}{2}\\
&=&1
\end{eqnarray*}

We conclude that $I_k=(2k)!!$, which matches with $M_k(f^2)=(2k)!!$, so job done.

\medskip

(3) At $n=2$ now, the problematics is closely related to the one for the complex normal laws $G_t$, that we studied in chapter 10, and we can recycle the computations there. To be more precise, with $f,g$ being independent variables, both following $g_1$, we have:
\begin{eqnarray*}
M_k(\chi_2^2)
&=&E((f^2+g^2)^k)\\
&=&\sum_r\binom{k}{r}E(f^{2r})E(g^{2k-2r})\\
&=&\sum_r\binom{k}{r}(2r)!!(2k-2r)!!\\
&=&\sum_r\frac{k!}{r!(k-r)!}\cdot\frac{(2r)!}{2^rr!}\cdot\frac{(2k-2r)!}{2^{k-r}(k-r)!}\\
&=&\frac{k!}{2^k}\sum_r\binom{2r}{r}\binom{2k-2r}{k-r}
\end{eqnarray*}

In order to finish now this moment computation, let us recall from chapter 7 that we have the following formula, coming from the generalized binomial formula:
$$\frac{1}{\sqrt{1+t}}=\sum_{q=0}^\infty\binom{2q}{q}\left(\frac{-t}{4}\right)^q$$

By taking the square of this series, we obtain the following formula:
$$\frac{1}{1+t}
=\sum_k\left(\frac{-t}{4}\right)^k\sum_r\binom{2r}{r}\binom{2k-2r}{k-r}$$

Now by looking at the coefficient of $t^k$ on both sides, we conclude that the sum on the right equals $4^k$. Thus, we can finish our moment computation, as follows:
$$M_k(\chi_2^2)=\frac{k!}{2^k}\times 4^k=2^kk!$$

(4) Still talking $n=2$, let us turn now to the second part of the work which is needed, namely computation of the moments for the predicted density, which is:
$$\chi_2^2=\frac{e^{-x/2}}{2}\,dx$$

As usual by partial integration, the moments of this predicted density are:
\begin{eqnarray*}
I_k
&=&\frac{1}{2}\int_0^\infty x^ke^{-x/2}dx\\
&=&\frac{1}{2}\int_0^\infty kx^{k-1}\cdot 2e^{-x/2}dx\\
&=&2k\cdot\frac{1}{2}\int_0^\infty x^{k-1}e^{-x/2}dx\\
&=&2kI_{k-1}
\end{eqnarray*}

As for the initial value, for this recurrence, this is given by the following formula, which by the way shows that we have indeed a probability measure, as we should:
\begin{eqnarray*}
I_0
&=&\frac{1}{2}\int_0^\infty e^{-x/2}dx\\
&=&\frac{1}{2}\left[-2e^{-x/2}\right]_0^\infty\\
&=&1
\end{eqnarray*}

We conclude that $I_k=2^kk!$, which matches with $M_k(\chi_2^2)=2^kk!$, so job done.

\medskip

(5) After all these preliminaries, time now to get into the general case, $n\in\mathbb N$. With $f_1,\ldots,f_n$ being independent variables, all following the normal law $g_1$, we have:
\begin{eqnarray*}
M_k(\chi_n^2)
&=&E((f_1^2+\ldots+f_n^2)^k)\\
&=&\sum_{i_1+\ldots+i_n=k}\binom{k}{i_1,\ldots,i_n}E(f_1^{2i_1})\ldots E(f_n^{2i_n})\\
&=&\sum_{i_1+\ldots+i_n=k}\frac{k!}{i_1!\ldots i_n!}(2i_1)!!\ldots(2i_n)!!\\
&=&\sum_{i_1+\ldots+i_n=k}\frac{k!}{i_1!\ldots i_n!}\cdot\frac{(2i_1)!}{2^{i_1}i_1!}\cdots\frac{(2i_n)!}{2^{i_n}i_n!}\\
&=&\frac{k!}{2^k}\sum_{i_1+\ldots+i_n=k}\frac{(2i_1)!}{i_1!i_1!}\cdots\frac{(2i_n)!}{i_n!i_n!}\\
&=&\frac{k!}{2^k}\sum_{i_1+\ldots+i_n=k}\binom{2i_1}{i_1}\ldots\binom{2i_n}{i_n}
\end{eqnarray*}

In order to finish this computation, we can use, as in (3), the following formula:
$$\frac{1}{\sqrt{1-4z}}=\sum_{i=0}^\infty\binom{2i}{i}z^i$$

Indeed, by taking the $n$-th power of this series, we obtain the following formula:
$$\frac{1}{\sqrt{(1-4z)^n}}
=\sum_kz^k\sum_{i_1+\ldots+i_n=k}\binom{2i_1}{i_1}\ldots\binom{2i_n}{i_n}$$

Now by putting everything together, the conclusion is that we have:
$$\frac{1}{\sqrt{(1-4z)^n}}
=\sum_{k\geq0}\frac{(2z)^k}{k!}\,M_k(\chi_n^2)$$

But the function on the left can be computed by using the generalized binomial formula, with exponent $p=-n/2$, and we obtain in this way:
\begin{eqnarray*}
\frac{1}{\sqrt{(1-4z)^n}}
&=&\sum_{k\geq0}\binom{-n/2}{k}(-4z)^k\\
&=&\sum_{k\geq0}\frac{(-n/2)(-n/2-1)\ldots(-n/2-k+1)}{k!}(-4z)^k\\
&=&\sum_{k\geq0}\frac{n(n+2)\ldots(n+2k-2)}{k!}(2z)^k
\end{eqnarray*}

And with this, good news, done with our moment computation, which yields:
$$M_k(\chi_n^2)=n(n+2)\ldots(n+2k-2)$$

(6) However, not time to celebrate yet, because we still have to do the second part of the work, namely computation of the moments for the predicted density, which is:
$$\chi_n^2=\frac{x^{n/2-1}e^{-x/2}}{2^{n/2}\Gamma(n/2)}\,dx$$

By using the definition of the gamma function, as formulated in (1), and then several times the formula $\Gamma(s+1)=s\Gamma(s)$, the moments of this predicted density are:
\begin{eqnarray*}
I_k
&=&\frac{1}{2^{n/2}\Gamma(n/2)}\int_0^\infty x^{n/2+k-1}e^{-x/2}\,dx\\
&=&\frac{1}{2^{n/2}\Gamma(n/2)}\int_0^\infty (2y)^{n/2+k-1}e^{-y}\,2dy\\
&=&\frac{2^k}{\Gamma(n/2)}\int_0^\infty y^{n/2+k-1}e^{-y}\,dy\\
&=&2^k\cdot\frac{\Gamma(n/2+k)}{\Gamma(n/2)}\\
&=&2^k\cdot\frac{\Gamma(n/2+k)}{\Gamma(n/2+k-1)}\cdot\frac{\Gamma(n/2+k-1)}{\Gamma(n/2+k-2)}\ldots\frac{\Gamma(n/2+1)}{\Gamma(n/2)}\\
&=&2^k\left(\frac{n}{2}+k-1\right)\left(\frac{n}{2}+k-2\right)\ldots\left(\frac{n}{2}\right)\\
&=&n(n+2)\ldots(n+2k-2)
\end{eqnarray*}

We conclude that we have $I_k=M_k(\chi_n^2)$, so job done, and theorem proved.
\end{proof}

Still with me I hope, after all these computations, and looking retrospectively, the final proof was  extremely quick, 1.5 pages, which can be further compacted to 1 page, or even 0.5, assuming some central binomial and gamma know-how, and, wait I'm not done with the laudations, with the moments being computed too. That was good work.

\bigskip

In short, we are basically done with the mathematics of the chi squared laws. And by adding a few missing things, all easy, here is our grand result, about these laws:

\begin{theorem}
The chi squared laws have the following properties: 
\begin{enumerate}
\item $\chi_n^2=law(f_1^2+\ldots+f_n^2)$, with $f_i\sim g_1$, independent.

\item The density is given by $\chi_n^2=\frac{x^{n/2-1}e^{-x/2}}{2^{n/2}\Gamma(n/2)}\,dx$.

\item The moments are $M_k=n(n+2)\ldots(n+2k-2)$.

\item We have $E=n$, $V=2n$, $\gamma=\sqrt{8/n}$, $\kappa=3+12/n$.

\item The Fourier transform is $F(x)=(1-2ix)^{-n/2}$.
\end{enumerate}
\end{theorem}

\begin{proof}
Here (1) is the definition of the chi squared laws, and (2) and (3) are things that we know, from Theorem 11.6 and its proof, which I insist, was ultimately something quite quick. Regarding now (4), according to our moment formula in (3), we have:
$$M_1=n$$
$$M_2=n(n+2)$$
$$M_3=n(n+2)(n+4)$$
$$M_4=n(n+2)(n+4)(n+6)$$

Thus we have the expectation $E=n$, and the variance $V=(n^2+2n)-n^2=2n$. In order to compute $\gamma$ and $\kappa$, observe that the needed central moments are given by:
$$M_3'=M_3-3EM_2+2E^3=8n$$
$$M_4'=M_4-4EM_3+6E^2M_2-3E^4=12n(n+4)$$

Thus, the skewness and kurtotis are given by the following formulae:
$$\gamma=\frac{8n}{2n\sqrt{2n}}=\sqrt{\frac{8}{n}}\quad,\quad\kappa=\frac{12n(n+4)}{4n^2}=3+\frac{12}{n}$$

Finally, regarding (5), the Fourier transform computation is as follows:
\begin{eqnarray*}
F(x)
&=&\frac{1}{2^{n/2}\Gamma(n/2)}\int_0^\infty y^{n/2-1}e^{-y/2+ixy}\,dy\\
&=&\frac{1}{2^{n/2}\Gamma(n/2)}\int_0^\infty y^{n/2-1}e^{(2ix-1)y/2}\,dy\\
&=&\frac{1}{2^{n/2}\Gamma(n/2)}\int_0^\infty\left(\frac{z}{1-2ix}\right)^{n/2-1}e^{-z/2}\,\frac{dz}{1-2ix}\\
&=&(1-2ix)^{-n/2}\cdot\frac{1}{2^{n/2}\Gamma(n/2)}\int_0^\infty z^{n/2-1}e^{-z/2}\,dz\\
&=&(1-2ix)^{-n/2}
\end{eqnarray*}

Thus, we are led to the conclusions in the statement.
\end{proof}

Finally, in relation with Remark 11.4 and Rat 11.5, a quick examination of what we did, at the end of the proof of Theorem 11.6 and afterwards, shows that our computations work for any $n>0$ real. With the notation $n=2m$, with $m>0$, and by adding a scale parameter $s>0$ too, we are led to the following generalization of Theorem 11.7:

\begin{theorem}
The gamma distribution of parameters $m,s>0$, given by
$$\mu_{ms}=\frac{x^{m-1}e^{-x/s}}{s^m\Gamma(m)}\,dx$$
generalizes the chi squared laws, $\chi^2_n=\mu_{n/2,2}$, and has the following properties: 
\begin{enumerate}
\item The moments are $M_k=s^km(m+1)\ldots(m+k-1)$.

\item We have $E=sm$, $V=s^2m$, $\gamma=2/\sqrt{m}$, $\kappa=3+6/m$.

\item The Fourier transform is $F(x)=(1-six)^{-m}$.
\end{enumerate}
\end{theorem}

\begin{proof}
This is a straightforward remake of Theorem 11.7, as follows:

\medskip

(1) The moment computation is as follows, using the various properties of the gamma function, and with $M_0=1$ showing that we have indeed a probability measure:
\begin{eqnarray*}
M_k
&=&\frac{1}{s^m\Gamma(m)}\int_0^\infty x^{m+k-1}e^{-x/s}\,dx\\
&=&\frac{1}{s^m\Gamma(m)}\int_0^\infty (sy)^{m+k-1}e^{-y}\,sdy\\
&=&\frac{s^k}{\Gamma(m)}\int_0^\infty y^{m+k-1}e^{-y}\,dy\\
&=&s^k\cdot\frac{\Gamma(m+k)}{\Gamma(m)}\\
&=&s^km(m+1)\ldots(m+k-1)
\end{eqnarray*}

(2) According to our moment formula found above, the first few moments are:
$$M_1=sm$$
$$M_2=s^2m(m+1)$$
$$M_3=s^3m(m+1)(m+2)$$
$$M_4=s^4m(m+1)(m+2)(m+3)$$

Thus the mean is $E=sm$, and the variance is $V=s^2m(m+1)-(sm)^2=s^2m$. In order to compute $\gamma$ and $\kappa$, observe that the needed central moments are given by:
$$M_3'=M_3-3EM_2+2E^3=2s^3m$$
$$M_4'=M_4-4EM_3+6E^2M_2-3E^4=3s^4m(m+2)$$

Thus, the skewness and kurtotis are given by the following formulae:
$$\gamma=\frac{2s^3m}{s^3m\sqrt{m}}=\frac{2}{\sqrt{m}}\quad,\quad\kappa=\frac{3s^4m(m+2)}{s^4m^2}=3+\frac{6}{m}$$

(3) The Fourier transform computation is, as before, as follows:
\begin{eqnarray*}
F(x)
&=&\frac{1}{s^m\Gamma(m)}\int_0^\infty y^{m-1}e^{-y/s+ixy}\,dy\\
&=&\frac{1}{s^m\Gamma(m)}\int_0^\infty y^{m-1}e^{(six-1)y/s}\,dy\\
&=&\frac{1}{2^m\Gamma(m)}\int_0^\infty\left(\frac{z}{1-six}\right)^{m-1}e^{-z/s}\,\frac{dz}{1-six}\\
&=&(1-six)^{-m}\cdot\frac{1}{s^m\Gamma(m)}\int_0^\infty z^{m-1}e^{-z/s}\,dz\\
&=&(1-six)^{-m}
\end{eqnarray*}

Thus, we are led to the conclusions in the statement.
\end{proof}

However, things are not over with this, and several other interesting parameters can be added to the picture. And exercise of course for you, to learn more about all this.

\section*{11c. Rayleigh laws}

Time to get now into the second family of laws introduced in Definition 11.3, namely the chi ones. These are a bit more complicated than the chi squared ones, with their mathematics appearing by extracting the ``square root'' of the mathematics from the previous section. Fortunately, this can be done, with the essentials being as follows:

\begin{theorem}
The chi law of parameter $n\in\mathbb N$, appearing as
$$\chi_n=law\left(\sqrt{f_1^2+\ldots+f_n^2}\right)$$
with $f_1,\ldots,f_n$ being independent, each following $g_1$, has density given by
$$\chi_n=\frac{x^{n-1}e^{-x^2/2}}{2^{n/2-1}\Gamma(n/2)}\,dx$$
and its moments are given by the following formula, valid for any $k\in\mathbb N$,
$$M_k=2^{k/2}\frac{\Gamma\left(\frac{n+k}{2}\right)}{\Gamma\left(\frac{n}{2}\right)}$$
with $\Gamma$ being as usual the gamma function.
\end{theorem}

\begin{proof}
This is something quite standard, the idea being as follows:

\medskip

(1) In order to compute the density, the simplest is to use Theorem 11.6, along with the following fact, coming from the definition of the chi and chi squared laws:
$$f\geq0\quad,\quad f\sim\chi_n\iff f^2\sim\chi_n^2$$

Indeed, assume $f\sim\chi_n$, and consider an arbitrary function $\varphi:[0,\infty)\to\mathbb R$. If we define $\psi:[0,\infty)\to\mathbb R$ by $\psi(x)=\varphi(\sqrt{x})$, we have the following computation:
\begin{eqnarray*}
E(\varphi(f))
&=&E(\psi(f^2))\\
&=&\frac{1}{2^{n/2}\Gamma(n/2)}\int_0^\infty\psi(x)x^{n/2-1}e^{-x/2}\,dx\\
&=&\frac{1}{2^{n/2}\Gamma(n/2)}\int_0^\infty\varphi(y)y^{n-2}e^{-y^2/2}\,2ydy\\
&=&\frac{1}{2^{n/2-1}\Gamma(n/2)}\int_0^\infty\varphi(y)y^{n-1}e^{-y^2/2}\,dy
\end{eqnarray*}

We conclude from this that the density of $f\sim\chi_n$ is the one in the statement.

\medskip

(2) As for the moment computation, this is straightforward, as follows:
\begin{eqnarray*}
M_k
&=&\frac{1}{2^{n/2-1}\Gamma(n/2)}\int_0^\infty x^{n+k-1}e^{-x^2/2}\,dx\\
&=&\frac{1}{2^{n/2-1}\Gamma(n/2)}\int_0^\infty\left(\sqrt{2y}\right)^{n+k-1}e^{-y}\,\frac{dy}{\sqrt{2y}}\\
&=&\frac{2^{(n+k-1)/2}}{2^{n/2-1/2}\Gamma(n/2)}\int_0^\infty y^{(n+k)/2-1}e^{-y}\,dy\\
&=&\frac{2^{k/2}}{\Gamma(n/2)}\,\Gamma\left(\frac{n+k}{2}\right)
\end{eqnarray*}

(3) Finally, in what regards the moments, whose formula is something quite tricky, let us comment a bit on this. In what regards the even moments, these are the moments of the chi squared law, that we computed before, the formulae being as follows:
$$M_2=n$$
$$M_4=n(n+2)$$
$$M_6=n(n+2)(n+4)$$
$$M_8=n(n+2)(n+4)(n+6)$$
$$\vdots$$

Which is very nice, at least we know one thing, so done with these.

\medskip

(4) In what regards the odd moments, however, things are more complicated. According to the general moment formula that we found, at $k=1$, the mean is:
$$E=\frac{\sqrt{2}}{\Gamma(n/2)}\,\Gamma\left(\frac{n+1}{2}\right)$$

Now recall from chapter 8 that the gamma function is given at half-integers by the following almost uniform formula, with $c_N=\sqrt{2},\sqrt{\pi}$ for $N$ even, odd:
$$\Gamma\left(\frac{N}{2}\right)=\frac{(N-1)!!}{2^{(N-1)/2}}\,c_N$$

By using this, the above formula of the mean becomes more explicit, as follows:
\begin{eqnarray*}
E
&=&\sqrt{2}\cdot\Gamma\left(\frac{n+1}{2}\right)\Gamma\left(\frac{n}{2}\right)^{-1}\\
&=&\sqrt{2}\cdot\frac{n!!}{2^{n/2}}\,c_{n+1}\cdot\frac{2^{(n-1)/2}}{(n-1)!!}\cdot\frac{1}{c_n}\\
&=&\frac{n!!}{(n-1)!!}\cdot\frac{c_{n+1}}{c_n}
\end{eqnarray*}

Finally, observe that we have the following formula, for the term on the right:
$$\frac{c_{n+1}}{c_n}=\begin{cases}
\sqrt{\frac{\pi}{2}}&(n\ {\rm even})\\
\sqrt{\frac{2}{\pi}}&(n\ {\rm odd})
\end{cases}$$

Summarizing, we know what $E$ is, but this remains something quite complicated.

\medskip

(5) Regarding the higher odd moments, no need for new computations, because the formula $\Gamma(s+1)=s\Gamma(s)$ does the job. If we denote as usual by $E$ the mean, which in practice is the beast computed above, the sequence of odd moments is:
$$M_1=E$$
$$M_3=(n+1)E$$
$$M_5=(n+1)(n+3)E$$
$$M_7=(n+1)(n+3)(n+5)E$$
$$\vdots$$

(6) Next, with what we have, we can certainly compute $(E,V,\gamma,\kappa)$, but the formulae of the variance $V$, skewness $\gamma$ and kurtosis $\kappa$ look quite bad, all featuring the above alien quantity $E$. And I will leave the computations here to you, as an exercise.

\medskip

(7) Finally, the formula of the Fourier transform is something quite complicated too, and I will leave again some exploration here to you, as an instructive exercise.
\end{proof}

Next, let us add a scaling parameter $t>0$ to what we have. We are led to:

\begin{theorem}
The chi law of parameters $n\in\mathbb N$ and $t>0$, appearing as
$$\chi_{nt}=law\left(\sqrt{f_1^2+\ldots+f_n^2}\right)$$
with $f_1,\ldots,f_n$ being independent, each following $g_t$, has density given by
$$\chi_{nt}=\frac{2x^{n-1}e^{-x^2/2t}}{(2t)^{n/2}\Gamma(n/2)}\,dx$$
and its moments are given by the following formula, valid for any $k\in\mathbb N$,
$$M_k=(2t)^{k/2}\frac{\Gamma\left(\frac{n+k}{2}\right)}{\Gamma\left(\frac{n}{2}\right)}$$
with $\Gamma$ being as usual the gamma function.
\end{theorem}

\begin{proof}
This is a straightforward parametric remake of what we have:

\medskip

(1) Let us first discuss the parametric chi squared law, which appears as follows, with $f_1,\ldots,f_n$ being as in the statement, namely independent, each following $g_t$:
$$\chi_{nt}^2=law\left(f_1^2+\ldots+f_n^2\right)$$

Since the moments of the normal laws obey to $M_k(g_t)=t^{k/2}M_k(g_1)$, we have:
$$M_k(\chi_{nt}^2)=t^kM_k(\chi_n^2)=t^kn(n+2)\ldots(n+2k-2)$$

But this is precisely the moment formula from Theorem 11.8, for the gamma distributions $\mu_{ms}$ there, at $m=n/2$ and $s=2t$. Thus, we have the following formula:
$$\chi_{nt}=\mu_{n/2,2t}$$

Finally, regarding the density, this comes from Theorem 11.8, and is given by:
$$\chi_{nt}^2
=\frac{x^{n/2-1}e^{-x/2t}}{(2t)^{n/2}\Gamma(n/2)}\,dx$$

(2) Now let us turn to the chi laws themselves. Again, since the moments of the normal laws obey to $M_k(g_t)=t^{k/2}M_k(g_1)$, we conclude from this that we have:
$$M_k(\chi_{nt})=t^{k/2}M_k(\chi_n)=(2t)^{k/2}\frac{\Gamma\left(\frac{n+k}{2}\right)}{\Gamma\left(\frac{n}{2}\right)}$$

In order to compute the density, the simplest is to use the density found in (1), along with the following fact, coming from our definition of the chi and chi squared laws:
$$f\geq0\quad,\quad f\sim\chi_{nt}\iff f^2\sim\chi_{nt}^2$$

Indeed, assume $f\sim\chi_{nt}$, and consider an arbitrary function $\varphi:[0,\infty)\to\mathbb R$. If we define $\psi:[0,\infty)\to\mathbb R$ by $\psi(x)=\varphi(\sqrt{x})$, we have the following computation:
\begin{eqnarray*}
E(\varphi(f))
&=&E(\psi(f^2))\\
&=&\frac{1}{(2t)^{n/2}\Gamma(n/2)}\int_0^\infty\psi(x)x^{n/2-1}e^{-x/2t}\,dx\\
&=&\frac{1}{(2t)^{n/2}\Gamma(n/2)}\int_0^\infty\varphi(y)y^{n-2}e^{-y^2/2t}\,2ydy\\
&=&\frac{2}{(2t)^{n/2}\Gamma(n/2)}\int_0^\infty\varphi(y)y^{n-1}e^{-y^2/2t}\,dy
\end{eqnarray*}

We conclude from this that the density of $f\sim\chi_{nt}$ is the one in the statement.

\medskip

(3) As for the discussion of $(E,V,\gamma,\kappa)$ and of the Fourier transform, along the lines of our previous discussion, at $t=1$, this would be a good exercise for you.
\end{proof}

What is next? Back to our familiar $n=1,2,3$ dimensions, and here, we have:

\index{Maxwell-Boltzmann}

\begin{theorem}
The chi laws in small dimensions are as follows:
\begin{enumerate}
\item At $n=1$ the chi law and its moments are as follows:
$$\chi_{1t}=\sqrt{\frac{2}{t\pi}}\,e^{-x^2/2t}dx\quad,\quad 
M_k=\sqrt{\frac{(2t)^k}{\pi}}\,\Gamma\left(\frac{k+1}{2}\right)$$

\item At $n=2$ the chi law and its moments are as follows:
$$\chi_{2t}=\frac{1}{t}\,xe^{-x^2/2t}dx\quad,\quad 
M_k=\sqrt{(2t)^k}\,\Gamma\left(\frac{k}{2}+1\right)$$

\item At $n=3$ the chi law and its moments are as follows:
$$\chi_{3t}=\sqrt{\frac{2}{t^3\pi}}\,x^2e^{-x^2/2t}dx\quad,\quad 
M_k=2\sqrt{\frac{(2t)^k}{\pi}}\,\Gamma\left(\frac{k+3}{2}\right)$$
\end{enumerate}
The law $\chi_{2t}$ is called Rayleigh law, and the law $\chi_{3t}$ is called Maxwell-Boltzmann law.
\end{theorem}

\begin{proof}
We use the general formulae found in Theorem 11.10, namely:
$$\chi_{nt}=\frac{2x^{n-1}e^{-x^2/2t}}{(2t)^{n/2}\Gamma(n/2)}\,dx\quad,\quad 
M_k=(2t)^{k/2}\frac{\Gamma\left(\frac{n+k}{2}\right)}{\Gamma\left(\frac{n}{2}\right)}$$

(1) At $n=1$ the density formula is as follows, using $\Gamma(1/2)=\sqrt{\pi}$:
$$\chi_{1t}=\frac{2e^{-x^2/2t}}{(2t)^{1/2}\sqrt{\pi}}\,dx=\sqrt{\frac{2}{t\pi}}\,e^{-x^2/2t}dx$$

As for the moment formula, this is as follows, again using $\Gamma(1/2)=\sqrt{\pi}$:
$$M_k=(2t)^{k/2}\frac{\Gamma\left(\frac{k+1}{2}\right)}{\sqrt{\pi}}
=\sqrt{\frac{(2t)^k}{\pi}}\,\Gamma\left(\frac{k+1}{2}\right)$$

(2) At $n=2$ the density formula is as follows, using $\Gamma(1)=1$:
$$\chi_{2t}=\frac{2xe^{-x^2/2t}}{2t}\,dx=\frac{1}{t}\,xe^{-x^2/2t}dx$$

As for the moment formula, this is as follows, again using $\Gamma(1)=1$:
$$M_k=(2t)^{k/2}\,\Gamma\left(\frac{k}{2}+1\right)
=\sqrt{(2t)^k}\,\Gamma\left(\frac{k}{2}+1\right)$$

(3) At $n=3$ the density formula is as follows, using $\Gamma(3/2)=\sqrt{\pi}/2$:
$$\chi_{3t}=\frac{2x^2e^{-x^2/2t}}{(2t)^{3/2}\sqrt{\pi}/2}\,dx
=\sqrt{\frac{2}{t^3\pi}}\,x^2e^{-x^2/2t}dx$$

As for the moment formula, this is as follows, again using $\Gamma(3/2)=\sqrt{\pi}/2$:
$$M_k=(2t)^{k/2}\frac{\Gamma\left(\frac{k+3}{2}\right)}{\sqrt{\pi}/2}
=2\sqrt{\frac{(2t)^k}{\pi}}\,\Gamma\left(\frac{k+3}{2}\right)$$

(4) Thus, we have the result, but let us comment a bit more on this, with a number of useful formulae. In what regards the expectations, the formulae are as follows:
$$E(\chi_{1t})=\sqrt{\frac{2t}{\pi}}\quad,\quad 
E(\chi_{2t})=\sqrt{\frac{t\pi}{2}}\quad,\quad 
E(\chi_{3t})=\sqrt{\frac{8t}{\pi}}$$

Regarding now the second moments, these are given by the following formulae:
$$M_2(\chi_{1t})=t\quad,\quad 
M_2(\chi_{2t})=2t\quad,\quad 
M_2(\chi_{3t})=3t$$

We conclude that the corresponding variances are given by:
$$V(\chi_{1t})=\left(1-\frac{2}{\pi}\right)t\quad,\quad
V(\chi_{2t})=\left(2-\frac{\pi}{2}\right)t\quad,\quad
V(\chi_{3t})=\left(3-\frac{8}{\pi}\right)t$$

(5) Finally, regarding the computation of the corresponding skewnesses $\gamma$ and kurtoses $\kappa$, and the study of the Fourier transforms, I will leave this to you, as an exercise.
\end{proof}

At the level of the phenomenology, the Rayleigh law appears for instance in the modeling of wind. As for the Maxwell-Boltzmann law, this appears in the context of gases. Many interesting things can be said here, of both mathematical and physical nature, and exercise of course for you, to learn more about all this, all good science.

\section*{11d. The heat kernel}

We would like to end this chapter with something simple and refreshing, namely the heat kernel. Let us start with the following result, which is of key importance:

\index{heat equation}
\index{lattice model}

\begin{theorem}
Heat diffusion in $\mathbb R^N$ is described by the heat equation
$$\dot{\varphi}=\alpha\Delta\varphi$$
where $\alpha>0$ is the thermal diffusivity of the medium, and with $\Delta$ given by
$$\Delta\varphi=\sum_{i=1}^N\frac{d^2\varphi}{dx_i^2}$$
being the Laplace operator, playing the role of a numeric second derivative.
\end{theorem}

\begin{proof}
This is something quite intuitive, the idea being as follows:

\medskip

(1) To start with, as a quick explanation for our equation, since the second derivative $\varphi''$ in one dimension, or the quantity $\Delta\varphi$ in general, computes the average value of a function $\varphi$ around a point, minus the value of $\varphi$ at that point, the heat equation as formulated above tells us that the rate of change $\dot{\varphi}$ of the temperature of the material at any given point must be proportional, with proportionality factor $\alpha>0$, to the average difference of temperature between that given point and the surrounding material.

\medskip

(2) The point now is that we can recover this equation by using a basic lattice model. Indeed, let us first assume, for simplifying, that we are in the one-dimensional case, $N=1$. Here our model looks as follows, with distance $l>0$ between neighbors:
$$\xymatrix@R=10pt@C=20pt{
\ar@{-}[r]&\circ_{x-l}\ar@{-}[r]^l&\circ_x\ar@{-}[r]^l&\circ_{x+l}\ar@{-}[r]&
}$$

In order to model heat diffusion, we have to implement the intuitive mechanism explained above, namely ``the rate of change of the temperature of the material at any given point must be proportional, with proportionality factor $\alpha>0$, to the average difference of temperature between that given point and the surrounding material''.

\medskip

(3) In practice, this leads to a condition as follows, expressing the change of the temperature $\varphi$, over a small period of time $\delta>0$:
$$\varphi(x,t+\delta)=\varphi(x,t)+\frac{\alpha\delta}{l^2}\sum_{x\sim y}\left[\varphi(y,t)-\varphi(x,t)\right]$$

To be more precise, we have made several assumptions here, as follows:

\medskip

-- General heat diffusion assumption: the change of temperature at any given point $x$ is proportional to the average over neighbors, $y\sim x$, of the differences $\varphi(y,t)-\varphi(x,t)$ between the temperatures at $x$, and at these neighbors $y$.

\medskip

-- Infinitesimal time and length conditions: in our model, the change of temperature at a given point $x$ is proportional to small period of time involved, $\delta>0$, and is inverse proportional to the square of the distance between neighbors, $l^2$.

\medskip

(4) Regarding these latter assumptions, the one regarding the proportionality with the time elapsed $\delta>0$ is something quite natural, physically speaking, and mathematically speaking too, because we can rewrite our equation as follows, making it clear that we have here an equation regarding the rate of change of temperature at $x$:
$$\frac{\varphi(x,t+\delta)-\varphi(x,t)}{\delta}=\frac{\alpha}{l^2}\sum_{x\sim y}\left[\varphi(y,t)-\varphi(x,t)\right]$$

As for the second assumption that we made above, namely inverse proportionality with $l^2$, this can be justified on physical grounds too, but again, perhaps the best is to do the math, which will show right away where this proportionality comes from. 

\medskip

(5) So, let us do the math. In the context of our 1D model the neighbors of $x$ are the points $x\pm l$, and so the equation that we wrote above takes the following form:
$$\frac{\varphi(x,t+\delta)-\varphi(x,t)}{\delta}=\frac{\alpha}{l^2}\Big[(\varphi(x+l,t)-\varphi(x,t))+(\varphi(x-l,t)-\varphi(x,t))\Big]$$

Now observe that we can write this equation as follows:
$$\frac{\varphi(x,t+\delta)-\varphi(x,t)}{\delta}
=\alpha\cdot\frac{\varphi(x+l,t)-2\varphi(x,t)+\varphi(x-l,t)}{l^2}$$

(6) We recognize on the right the usual approximation of the second derivative, coming from calculus. Thus, when taking the continuous limit of our model, $l\to 0$, we get:
$$\frac{\varphi(x,t+\delta)-\varphi(x,t)}{\delta}
=\alpha\cdot\varphi''(x,t)$$

Now with $t\to0$, we are led in this way to the heat equation, namely:
$$\dot{\varphi}(x,t)=\alpha\cdot\varphi''(x,t)$$

(7) With this done, let us discuss now 2 dimensions. Here we can use a lattice model as follows, with all lengths being $l>0$, for simplifying:
$$
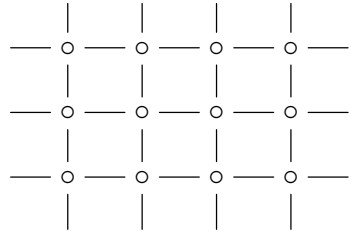
\xymatrix@R=12pt@C=15pt{
&\ar@{-}[d]&\ar@{-}[d]&\ar@{-}[d]&\ar@{-}[d]\\
\ar@{-}[r]&\circ\ar@{-}[r]\ar@{-}[d]&\circ\ar@{-}[r]\ar@{-}[d]&\circ\ar@{-}[r]\ar@{-}[d]&\circ\ar@{-}[r]\ar@{-}[d]&\\
\ar@{-}[r]&\circ\ar@{-}[r]\ar@{-}[d]&\circ\ar@{-}[r]\ar@{-}[d]&\circ\ar@{-}[r]\ar@{-}[d]&\circ\ar@{-}[r]\ar@{-}[d]&\\
\ar@{-}[r]&\circ\ar@{-}[r]\ar@{-}[d]&\circ\ar@{-}[r]\ar@{-}[d]&\circ\ar@{-}[r]\ar@{-}[d]&\circ\ar@{-}[r]\ar@{-}[d]&\\
&&&&&
}$$

(8) We have to implement now the physical heat diffusion mechanism, namely ``the rate of change of the temperature of the material at any given point must be proportional, with proportionality factor $\alpha>0$, to the average difference of temperature between that given point and the surrounding material''. In practice, this leads to a condition as follows, expressing the change of the temperature $\varphi$, over a small period of time $\delta>0$:
$$\varphi(x,y,t+\delta)=\varphi(x,y,t)+\frac{\alpha\delta}{l^2}\sum_{(x,y)\sim(u,v)}\left[\varphi(u,v,t)-\varphi(x,y,t)\right]$$

In fact, we can rewrite our equation as follows, making it clear that we have here an equation regarding the rate of change of temperature at $x$:
$$\frac{\varphi(x,y,t+\delta)-\varphi(x,y,t)}{\delta}=\frac{\alpha}{l^2}\sum_{(x,y)\sim(u,v)}\left[\varphi(u,v,t)-\varphi(x,y,t)\right]$$

(9) So, let us do the math. In the context of our 2D model the neighbors of $x$ are the points $(x\pm l,y\pm l)$, so the equation above takes the following form:
\begin{eqnarray*}
&&\frac{\varphi(x,y,t+\delta)-\varphi(x,y,t)}{\delta}\\
&=&\frac{\alpha}{l^2}\Big[(\varphi(x+l,y,t)-\varphi(x,y,t))+(\varphi(x-l,y,t)-\varphi(x,y,t))\Big]\\
&+&\frac{\alpha}{l^2}\Big[(\varphi(x,y+l,t)-\varphi(x,y,t))+(\varphi(x,y-l,t)-\varphi(x,y,t))\Big]
\end{eqnarray*}

Now observe that we can write this equation as follows:
\begin{eqnarray*}
\frac{\varphi(x,y,t+\delta)-\varphi(x,y,t)}{\delta}
&=&\alpha\cdot\frac{\varphi(x+l,y,t)-2\varphi(x,y,t)+\varphi(x-l,y,t)}{l^2}\\
&+&\alpha\cdot\frac{\varphi(x,y+l,t)-2\varphi(x,y,t)+\varphi(x,y-l,t)}{l^2}
\end{eqnarray*}

(10) We recognize on the right the usual approximations of the second derivative. Thus, when taking the continuous limit of our model, $l\to 0$, we get:
$$\frac{\varphi(x,y,t+\delta)-\varphi(x,y,t)}{\delta}
=\alpha\left(\frac{d^2\varphi}{dx^2}+\frac{d^2\varphi}{dy^2}\right)(x,y,t)$$

Now with $t\to0$, we are led in this way to the heat equation, namely:
$$\dot{\varphi}(x,y,t)=\alpha\cdot\Delta\varphi(x,y,t)$$

(11) Finally, in arbitrary $N$ dimensions the same argument carries over, namely a straightforward lattice model, and this gives the heat diffusion equation, as formulated in the statement. We will leave the details here as a routine exercise.
\end{proof}

Regarding now the resolution of the heat equation, we have here:

\index{heat kernel}
\index{convolution}

\begin{theorem}
The heat diffusion equation, $\dot{\varphi}=\alpha\Delta\varphi$ with $\alpha>0$, with initial condition $\varphi(x,0)=f(x)$, has as solution the function
$$\varphi(x,t)=\int_{\mathbb R^N}K_t(x-y)f(y)dy$$
where the function $K_t:\mathbb R^N\to\mathbb R$, called heat kernel, given by
$$K_t(x)=\frac{1}{\sqrt{(4\pi\alpha t)^N}}\,e^{-||x||^2/4\alpha t}$$
is the standard solution, coming from the initial data $f=\delta_0$, Dirac mass at $0$.
\end{theorem}

\begin{proof}
Many things can be said here, the idea being as follows:

\medskip

(1) Let us first discuss what happens in 1 dimension. We first have to check that the heat kernel $K_t$ is indeed a solution. The time derivative is:
$$\dot{K_t}=-\frac{1}{2t\sqrt{4\pi\alpha t}}\,e^{-x^2/4\alpha t}+\frac{x^2}{4\alpha t^2\sqrt{4\pi\alpha t}}\,e^{-x^2/4\alpha t}$$

Regarding the first space derivative, this is given by the following formula:
$$K_t'=-\frac{x}{2\alpha t\sqrt{4\pi\alpha t}}\,e^{-x^2/4\alpha t}$$

As for the second space derivative, this is given by the following formula:
$$K_t''=-\frac{1}{2\alpha t\sqrt{4\pi\alpha t}}\,e^{-x^2/4\alpha t}+\frac{x^2}{4\alpha^2t^2\sqrt{4\pi\alpha t}}\,e^{-x^2/4\alpha t}$$

Thus, we can see that the heat equation $\dot{\varphi}=\alpha\varphi''$ is indeed satisfied.

\medskip

(2) Next, let us convolve $K_t$ with an arbitrary function $f$, as follows:
$$\varphi(x,t)=\int_{\mathbb R^N}K_t(x-y)f(y)dy$$

The point now is that, when doing so, all the computations from (1) will perturb well, according to the following formulae, which are both clear:
$$\dot{\varphi}(x,t)=\int_\mathbb R\dot{K}_t(x-y)f(y)dy$$
$$\varphi''(x,t)=\int_\mathbb RK''_t(x-y)f(y)dy$$

Thus, we can see that the heat equation $\dot{\varphi}=\alpha\varphi''$ is again satisfied. 

\medskip

(3) Next, we must discuss initial conditions and normalization. And as a first remark here, thanks to the Gauss integral, our normalization for $K_t$ makes it of mass one:
\begin{eqnarray*}
\int_\mathbb RK_t(x)dx
&=&\frac{1}{\sqrt{4\pi\alpha t}}\int_\mathbb Re^{-x^2/4\alpha t}dx\\
&=&\frac{1}{\sqrt{4\pi\alpha t}}\int_\mathbb Re^{-y^2}\sqrt{4\alpha t}\,dy\\
&=&\frac{1}{\sqrt{4\pi\alpha t}}\cdot \sqrt{4\alpha t}\cdot\sqrt{\pi}\\
&=&1
\end{eqnarray*}

But with this in hand, it is routine to check that $K_t$ comes indeed from the simplest situation, that of a radiating point body placed at $0$, as said in the statement, and then that the solution found in (2) by convolving comes indeed from the initial data $f$.

\medskip

(4) Finally, everything extends well to arbitrary dimensions, and we will leave this as an exercise. As for the uniqueness issues, we will leave this as an exercise too.
\end{proof}

As a conclusion to this, the 3D normal laws, also called heat kernel, are smart enough for solving the heat diffusion equation, in our usual 3D. Which is quite remarkable.

\section*{11e. Exercises}

Tough chapter that we had here, and as exercises on this, we have:

\begin{exercise}
Experiment a bit with the CDF of the normal laws.
\end{exercise}

\begin{exercise}
Experiment too with the moments of the shifted normal laws $g_t^a$.
\end{exercise}

\begin{exercise}
Learn more about the gamma distribution, and related laws.
\end{exercise}

\begin{exercise}
Compute the variance, skewness and kurtosis of $\chi_n$.
\end{exercise}

\begin{exercise}
Then compute the variance, skewness and kurtosis of $\chi_{nt}$.
\end{exercise}

\begin{exercise}
Learn more about the occurences of the Rayleigh law $\chi_{2t}$.
\end{exercise}

\begin{exercise}
Learn more as well about the Maxwell-Boltzmann law $\chi_{3t}$.
\end{exercise}

\begin{exercise}
Learn about the physics corrections to the heat diffusion equation.
\end{exercise}

A bonus exercise, now that you master matematics, start reading a physics book.

\chapter{Hyperspherical laws}

\section*{12a. Circle, Wallis}

We would like to end the present Part III with some pure mathematics, all useful things, in relation with the spaces $\mathbb R^N$, and their unit spheres $S^{N-1}_\mathbb R\subset\mathbb R^N$. Which will eventually lead us, and no surprise here, into probability and normal laws.

\bigskip

Let us start with a discussion at $N=2$. I mean, back to school, and here is what we learned there, the hard way, and believe me, these things were tough for me too:

\begin{theorem}
The trigonometric functions of $t\in\mathbb R$ can be computed using
$$
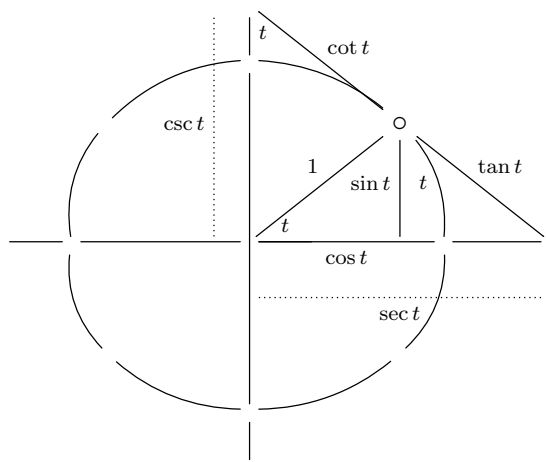
\xymatrix@R=14pt@C=6.4pt{
&&&&&&\ar@{.}[dddd]_{\csc t}&\ar@{-}[ddrrrr]^{\cot t}\ar@{-}[d]^t\\
&&&&&&&\ar@{-}[dddddd]\ar@{-}@/_/[dllll]&\\
&&&\ar@{-}@/_/[ddl]&&&&&&&&\circ\ar@{-}[ddrrrr]^{\tan t}\ar@{-}[dd]_{\sin t}\ar@{-}@/_/[ullll]\\
&&&&&&&&&&&\\
\ar@{-}[rr]&&\ar@{-}@/_/[ddr]\ar@{-}[rrrrr]&&&&&\ar@{-}[uurrrr]^1\ar@{-}[rrrrr]_{\cos t}\ar@{-}[rr]^t&&&&&\ar@{-}@/_/[uul]^t\ar@{-}[rrr]&&&\\
&&&&&&&\ar@{.}[rrrrrrrr]_{\sec t}&&&&&&&&\\
&&&\ar@{-}@/_/[drrrr]&&&&&&&&\ar@{-}@/_/[uur]\\
&&&&&&&\ar@{-}@/_/[urrrr]\ar@{-}[d]\\
&&&&&&&}$$
with due attention to the positives and negatives, for $t<0$ or $t>\pi/2$.
\end{theorem}

\begin{proof}
No comment here, this is something that we know, and are proud of.
\end{proof}

Coming as good news now, and away from these childhood nightmares, for our pure mathematics purposes here, we will only need a soft version of Theorem 12.1, namely:

\begin{principle}
The interplay between the plane and its unit circle
$$C\subset\mathbb R^2$$
is the alpha and omega of everything, in 2D mathematics.
\end{principle}

So, this will be the guiding principle for the present chapter, and of course with the above inclusion $C\subset\mathbb R^2$ being replaced by $S^{N-1}_\mathbb R\subset\mathbb R^N$, in higher dimensions.

\bigskip

This being said, do we really know well the circle $C\subset\mathbb R^2$? What we have in Theorem 12.1 is certainly first class mathematics, providing the key to everything 2D, but in practice, it is possible to further build with some interesting theorems, on that.

\bigskip

To be more precise, let us formulate the following questions:

\begin{questions}
In relation with the plane and its unit circle, $C\subset\mathbb R^2$:
\begin{enumerate}
\item How to parametrize these, and what can we get out of that?

\item How to integrate on the unit circle, any interesting laws there?

\item What about deformations of the unit circle, anything to do there?

\item What can we say about the king in all this, $\pi$, of analytic nature?
\end{enumerate}
\end{questions}

Getting started now, with question (1), we certainly know that we have polar coordinates in the plane, $x=r\cos t,y=r\sin t$, or $z=re^{it}$ in fancier notation, which restrict to the circle into $x=\cos t,y=\sin t$, or $z=e^{it}$. As basic applications of this, we have:

\index{Gauss integral}
\index{Fresnel integral}

\begin{theorem}
We have the following formula, due to Gauss,
$$\int_\mathbb Re^{-x^2}dx=\sqrt{\pi}$$
as well as the following formula, due to Fresnel, with $\sqrt{i}=e^{\pi i/4}$,
$$\int_\mathbb R e^{it^2}dt=\sqrt{\pi i}$$
which in practice gives $\int_\mathbb R\sin(t^2)dt=\int_\mathbb R\cos(t^2)dt=\sqrt{\pi/2}$.
\end{theorem}

\begin{proof}
Many things can be said here, the idea being as follows:

\medskip

(1) Regarding the Gauss formula, we already know it, but always a pleasure to talk about it again. The Jacobian of the polar coordinates is given by:
$$J=\begin{vmatrix}
\frac{d(r\cos t)}{dr}&&\frac{d(r\cos t)}{dt}\\
\\
\frac{d(r\sin t)}{dr}&&\frac{d(r\sin t)}{dt}
\end{vmatrix}\
=\begin{vmatrix}
\cos t&-r\sin t\\
\sin t&r\cos t
\end{vmatrix}
=r$$

Thus, the Gauss integral can be computed by using 2 dimensions, as follows:
$$\int_\mathbb R\int_\mathbb Re^{-x^2-y^2}dxdy
=\int_0^{2\pi}\int_0^\infty e^{-r^2}rdrdt
=2\pi\left[-\frac{e^{-r^2}}{2}\right]_0^\infty
=\pi$$

And with the remark, of course, that we already made several times in the above, that the Gauss integral cannot be computed in 1D alone. This is how life is.

\medskip

(2) Getting now to the Fresnel integral, due to $e^{\pi i/4}=(1+i)/\sqrt{2}$, the formula in the statement takes the following form, which makes the link with the last 2 formulae:
$$\int_0^\infty e^{it^2}dt=\sqrt{\frac{\pi}{2}}\cdot\frac{1+i}{2}$$

So, let us prove this. For this purpose, consider the following function:
$$f(t)=\int_0^\infty\frac{e^{(i-u^2)t^2}}{i-u^2}\,du$$

The derivative of this function is then given by the following formula:
$$f'(t)
=2t\int_0^\infty e^{(i-u^2)t^2}du
=2te^{it^2}\int_0^\infty e^{-u^2t^2}du
=\sqrt{\pi}e^{it^2}$$

Now let us integrate this derivative, from 0 to $\infty$. We obtain in this way:
$$\sqrt{\pi}\int_0^\infty e^{it^2}dt
=f(\infty)-f(0)
=\int_0^\infty\frac{1}{u^2-i}\,du$$

Summarizing, we have obtained the following formula, for the Fresnel integral:
$$\int_0^\infty e^{it^2}dt=\frac{1}{\sqrt{\pi}}\int_0^\infty\frac{1}{u^2-i}\,du$$

(3) In order to compute this latter integral, observe first that we have:
$$\int_0^\infty\frac{1}{u^2-i}\,du
=\int_0^\infty\frac{u^2}{u^4+1}\,du+i\int_0^\infty\frac{1}{u^4+1}\,du$$

Next, the two real integrals are equal, because with $u\to u^{-1}$ we obtain:
$$\int_0^\infty\frac{u^2}{u^4+1}\,du
=\int_0^\infty\frac{u^{-2}}{u^{-4}+1}\,u^{-2}du
=\int_0^\infty\frac{1}{u^4+1}\,du$$

Also, we can compute the sum of these integrals by using $t=u-u^{-1}$, as follows:
$$\int_0^\infty\frac{u^2+1}{u^4+1}\,du
=\int_0^\infty\frac{1+u^{-2}}{u^2+u^{-2}}\,du
=\int_0^\infty\frac{dt}{t^2+2}=\frac{\pi}{\sqrt{2}}$$

Thus, we are led to the following conclusion, regarding the integral in (2):
$$\int_0^\infty\frac{1}{u^2-i}\,du
=\frac{\pi}{2\sqrt{2}}+i\,\frac{\pi}{2\sqrt{2}}
=\frac{\pi}{2}\cdot\frac{1+i}{\sqrt{2}}$$

(4) Summarizing, we are done, and by putting everything together, we obtain:
$$\int_0^\infty e^{it^2}dt=\frac{1}{\sqrt{\pi}}\times\frac{\pi}{2}\cdot\frac{1+i}{\sqrt{2}}
=\sqrt{\frac{\pi}{2}}\cdot\frac{1+i}{2}$$

But this is exactly what we wanted, and this ends the proof of our result.
\end{proof}

And with this, done with Question 12.3 (1)? Not exactly, because as a rival to the polar coordinates, which are certainly very useful, and whose applications abound, we have the stereographic projection, which can be useful too, for various purposes:

\begin{theorem}
The stereographic projection from $i$, namely
$$
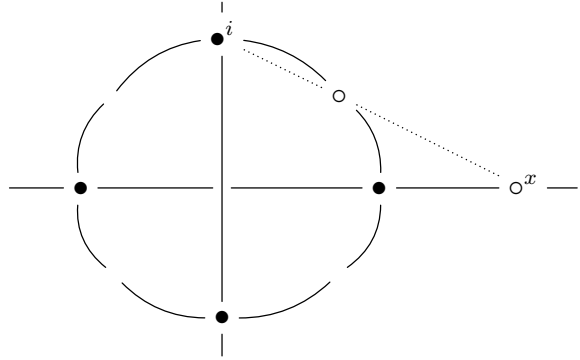
\xymatrix@R=7pt@C=2.1pt{
&&&&&&&&\ar@{-}[d]\\
&&&&&&&&\bullet^i\ar@{-}[dddddd]\ar@{-}@/_/[dllll]\ar@{.}[drrrr]&\\
&&&&\ar@{-}@/_/[ddl]&&&&&&&&\circ\ar@{-}@/_/[ullll]\ar@{.}[ddrrrrrr]\\
&&&&&&&&&&&&\\
\ar@{-}[rrr]&&&\bullet\ar@{-}@/_/[ddr]\ar@{-}[rrrrr]&&&&&\ar@{-}[rrrrr]&&&&&\bullet\ar@{-}@/_/[uul]\ar@{-}[rrrrr]&&&&&\circ^x\ar@{-}[rr]&&\\
&&&&&&&&&&&&\\
&&&&\ar@{-}@/_/[drrrr]&&&&&&&&\ar@{-}@/_/[uur]\\
&&&&&&&&\bullet\ar@{-}@/_/[urrrr]\ar@{-}[d]\\
&&&&&&&&}$$
makes correspond numbers $e^{it}\in\mathbb T$ to real numbers $x\in\mathbb R$ according to the formulae 
$$x=\frac{\cos t}{1-\sin t}\quad,\quad e^{it}=\frac{2x+i(x^2-1)}{x^2+1}$$
and can be used for various mathematical purposes.
\end{theorem}

\begin{proof}
By Thales, the formula of the stereographic projection is as follows:
$$\frac{x-\cos t}{x}=\frac{\sin t}{1}\implies x=\frac{\cos t}{1-\sin t}$$

As for the inverse, observe that with $x$ as above, we have the following formulae:
$$2x=\frac{2\cos t}{1-\sin t}\quad,\quad x^2-1=\frac{2\sin t}{1-\sin t}
\quad,\quad x^2+1=\frac{2}{1-\sin t}$$

We conclude that the angle $t$ can be recovered from $x$ via the following formulae:
$$\cos t=\frac{2x}{x^2+1}\quad,\quad\sin t=\frac{x^2-1}{x^2+1}$$

Thus, we are led to the formula of $e^{it}$ in the statement. As for the last assertion, the applications abound, and exercise for you, to learn a bit about all this.
\end{proof}

Getting now to Question 12.3 (2), the circle coordinates are $\cos t,\sin t$, and the first question is, how to integrate arbitrary products of these quantities $\cos t,\sin t$. And in answer here, we first have the following formula, for the arbitrary powers of $\cos t,\sin t$, that we already met in chapter 8, and that we reproduce here for convenience:

\index{double factorials}
\index{trigonometric integral}

\begin{theorem}[Wallis]
We have the following formulae,
$$\int_0^{\pi/2}\cos^pt\,dt=\int_0^{\pi/2}\sin^pt\,dt=\left(\frac{\pi}{2}\right)^{\varepsilon(p)}\frac{p!!}{(p+1)!!}$$
where $\varepsilon(p)=1$ if $p$ is even, and $\varepsilon(p)=0$ if $p$ is odd.
\end{theorem}

\begin{proof}
Let us first compute the integral on the left in the statement:
$$I_p=\int_0^{\pi/2}\cos^pt\,dt$$

We can do this by partial integration. Indeed, we have the following formula:
\begin{eqnarray*}
(\cos^pt\sin t)'
&=&p\cos^{p-1}t(-\sin t)\sin t+\cos^pt\cos t\\
&=&p\cos^{p+1}t-p\cos^{p-1}t+\cos^{p+1}t\\
&=&(p+1)\cos^{p+1}t-p\cos^{p-1}t
\end{eqnarray*}

By integrating now between $0$ and $\pi/2$, we obtain the following formula:
$$(p+1)I_{p+1}=pI_{p-1}$$

Thus we can compute $I_p$ by recurrence, and we obtain in this way:
\begin{eqnarray*}
I_p
&=&\frac{p-1}{p}\,I_{p-2}\\
&=&\frac{p-1}{p}\cdot\frac{p-3}{p-2}\,I_{p-4}\\
&=&\frac{p-1}{p}\cdot\frac{p-3}{p-2}\cdot\frac{p-5}{p-4}\,I_{p-6}\\
&\vdots&\\
&=&\frac{p!!}{(p+1)!!}\,I_{1-\varepsilon(p)}
\end{eqnarray*}

But $I_0=\frac{\pi}{2}$ and $I_1=1$, so we get the result. As for the second formula, this follows from the first one, with $t=\frac{\pi}{2}-s$. Thus, we have proved both formulae in the statement.
\end{proof}

More generally now, we have the following result, solving our integration problem:

\index{trigonometric integral}

\begin{theorem}[Wallis 2]
We have the following formula,
$$\int_0^{\pi/2}\cos^pt\sin^qt\,dt=\left(\frac{\pi}{2}\right)^{\varepsilon(p)\varepsilon(q)}\frac{p!!q!!}{(p+q+1)!!}$$
where $\varepsilon(p)=1$ if $p$ is even, and $\varepsilon(p)=0$ if $p$ is odd.
\end{theorem}

\begin{proof}
Let $I_{pq}$ be the integral in the statement. In order to do the partial integration, a bit as we previously did at $p=0$ or $q=0$, observe that we have:
\begin{eqnarray*}
(\cos^pt\sin^qt)'
&=&p\cos^{p-1}t(-\sin t)\sin^qt\\
&+&\cos^pt\cdot q\sin^{q-1}t\cos t\\
&=&-p\cos^{p-1}t\sin^{q+1}t+q\cos^{p+1}t\sin^{q-1}t
\end{eqnarray*}

By integrating between $0$ and $\pi/2$, we obtain, for $p,q>0$:
$$pI_{p-1,q+1}=qI_{p+1,q-1}$$

Thus, we can compute $I_{pq}$ by recurrence. When $q$ is even we obtain, in this way:
\begin{eqnarray*}
I_{pq}
&=&\frac{q-1}{p+1}\,I_{p+2,q-2}\\
&=&\frac{q-1}{p+1}\cdot\frac{q-3}{p+3}\,I_{p+4,q-4}\\
&=&\frac{q-1}{p+1}\cdot\frac{q-3}{p+3}\cdot\frac{q-5}{p+5}\,I_{p+6,q-6}\\
&\vdots&\\
&=&\frac{p!!q!!}{(p+q)!!}\,I_{p+q}
\end{eqnarray*}

But the last term comes from Theorem 12.6, and we obtain the result, for $q$ even:
$$I_{pq}
=\frac{p!!q!!}{(p+q)!!}\left(\frac{\pi}{2}\right)^{\varepsilon(p+q)}\frac{(p+q)!!}{(p+q+1)!!}
=\left(\frac{\pi}{2}\right)^{\varepsilon(p)\varepsilon(q)}\frac{p!!q!!}{(p+q+1)!!}$$

Observe that this gives the result for $p$ even as well, by symmetry. In the remaining case, where both $p,q$ are odd, we can use again $pI_{p-1,q+1}=qI_{p+1,q-1}$, which gives:
\begin{eqnarray*}
I_{pq}
&=&\frac{q-1}{p+1}\,I_{p+2,q-2}\\
&=&\frac{q-1}{p+1}\cdot\frac{q-3}{p+3}\,I_{p+4,q-4}\\
&=&\frac{q-1}{p+1}\cdot\frac{q-3}{p+3}\cdot\frac{q-5}{p+5}\,I_{p+6,q-6}\\
&\vdots&\\
&=&\frac{p!!q!!}{(p+q-1)!!}\,I_{p+q-1,1}
\end{eqnarray*}

But the last term is easy to compute, by using the primitive, as follows:
$$I_{p1}
=\int_0^{\pi/2}\cos^pt\sin t\,dt
=-\frac{1}{p+1}\int_0^{\pi/2}(\cos^{p+1}t)'\,dt
=\frac{1}{p+1}$$

We can therefore finish our computation in the case $p,q$ odd, and we obtain:
$$I_{pq}
=\frac{p!!q!!}{(p+q-1)!!}\cdot\frac{1}{p+q}
=\frac{p!!q!!}{(p+q+1)!!}$$

Thus, we are led to the formula in the statement, the exponent of $\pi/2$ appearing there being $\varepsilon(p)\varepsilon(q)=0\cdot 0=0$ in the present case, and this finishes the proof.
\end{proof}

Still in relation with Question 12.3 (2), we can now get to its second part, probability on the circle. And here, we actually already met the circle in chapter 5, when discussing the general probability axiomatics, with our observation there being as follows:

\begin{proposition}
The set $A\subset C$ obtained by picking one point from each orbit of
$$\mathbb Q/\mathbb Z\curvearrowright C\quad,\quad r(z)=e^{2r\pi i}z$$
cannot be measurable, with respect to the uniform measure on $C$.
\end{proposition}

\begin{proof}
This is indeed something that we know from chapter 5, the idea being that $C=\bigsqcup_rr(A)$ prevents $A$ from having measure 0, because this would lead to $\mu(C)=0$, and prevents $A$ from having measure $>0$ either, because this would lead to $\mu(C)=\infty$.
\end{proof}

However, it is not about such phenomena, duly analyzed in chapter 5 by myself, cat and rat, with some help from Rudin \cite{rud}, that we would like to talk about here. Indeed, based on Theorem 12.7, we can have some concrete mathematics going, as follows:

\begin{theorem}
With the notation $z=(x,y)$ for the points $z\in C$, and with respect to the uniform measure on $C$, both variables $x,y:C\to[-1,1]$ follow the law
$$\mu=\frac{1}{\pi\sqrt{1-x^2}}\,dx$$
that is, the arcsine law on $[-1,1]$, whose even the moments are as follows:
$$M_{2r}=\frac{1}{4^r}\binom{2r}{r}$$
Moreover, the joint moments of $(x,y)$ are given by the following formula,
$$E(x^{2r}y^{2s})=\frac{1}{4^{r+s}}\cdot\frac{(2r)!(2s)!}{r!s!(r+s)!}$$
with the other joint moments, containing an odd exponent, vanishing.
\end{theorem}

\begin{proof}
We certainly have variables $x,y:C\to[-1,1]$ as above, whose individual and joint behavior, meaning law and moments, we can investigate. But this can be done by using the Wallis formula from Theorem 12.7, the details being as follows:

\medskip

(1) With the convention in the statement, the Wallis formula from Theorem 12.7 takes the following form, in the non-vanishing case, where both exponents are even:
\begin{eqnarray*}
E(x^{2r}y^{2s})
&=&\frac{(2r)!!(2s)!!}{(2r+2s+1)!!}\\
&=&\frac{(2r)!/(2^rr!)\cdot (2s)!/(2^ss!)}{2^{r+s}(r+s)!}\\
&=&\frac{1}{4^{r+s}}\cdot\frac{(2r)!(2s)!}{r!s!(r+s)!}
\end{eqnarray*}

Thus, we are led to the various moment assertions in the statement. 

\medskip

(2) As for the formula of the density, this comes either by thinking, or from our previous knowledge of the arcsine law, from chapter 7, or simply as follows:
$$\frac{1}{\pi}\int_{-1}^1\frac{x^{2r}}{\sqrt{1-x^2}}\,dx=\frac{2}{\pi}\int_0^{\pi/2}\frac{\sin^{2r}t}{\cos t}\,\cos tdt=\frac{2}{\pi}\int_0^{\pi/2}\sin^{2r}t\,dt$$

Thus, the moments match indeed, and we have our density, as desired.
\end{proof}

\section*{12b. Hermite, Buffon}

Getting now to Question 12.3 (3), regarding deformations of the circle, and their theory, here we can talk about ellipses, whose basic theory is as follows:

\index{ellipse}
\index{conic}

\begin{theorem}
We can talk about ellipses, in several possible ways, as follows:
\begin{enumerate}
\item As the curves appearing via $d(z,p)+d(z,q)=l$, with $p,q\in\mathbb R^2$, and $l>0$.

\item Up to a translation and rotation, as the curves given by $(x/a)^2+(y/b)^2=1$.

\item Or equivalently, in polar coordinates, given by $x=a\cos t$, $y=b\sin t$.

\item As the conics, that is, the plane curves of degree $2$, which are compact. 

\item As the compact curves obtained by cutting a 3D cone with a plane.

\item As the pictures of a circle, taken from various 3D perspectives.

\item As the trajectories of the planets around the Sun.
\end{enumerate}
\end{theorem}

\begin{proof}
All this is first-class mathematics, the idea being as follows:

\medskip

(1) This is something practical, allowing you to draw ellipses, armed with a string.

\medskip

(2) The equivalence with (1) comes from some computations, good exercise for you.

\medskip

(3) This is indeed something which is clearly equivalent to (2).

\medskip

(4) This is indeed equivalent to (2), via some standard algebraic computations.

\medskip

(5) This is indeed equivalent to (4), via some standard 3D computations.

\medskip

(6) This is the same thing as (5), for the attention of Humanities students.

\medskip

(7) This is tough, due to Newton, check for instance chapter 1 in my book \cite{ba3}.
\end{proof}

Getting now to integration aspects, we have here the following result:

\begin{theorem}
Given an ellipse, parametrized as usual as
$$\left(\frac{x}{a}\right)^2+\left(\frac{y}{b}\right)^2=1$$
with $a,b>0$, its area is $A=\pi ab$. As for the length, this is given by the formula
$$L=4\int_0^{\pi/2}\sqrt{a^2\sin^2t+b^2\cos^2t}\,dt$$
and with this integral being generically not computable.
\end{theorem}

\begin{proof}
The formula of the area is standard, either by computing, or by saying that the passage from the unit circle to our ellipse is given by $(x,y)\to(ax,by)$, having Jacobian $J=ab$. As for the length assertion, which is something quite surprising:

\medskip

(1) To start with, what is the length of a curve $\gamma:[a,b]\to\mathbb R^2$? Good question, and in answer, a physicist would say that this is the quantity obtained by integrating the magnitude of the velocity vector over the curve, with respect to time, leading to:
$$L(\gamma)=\int_a^b||\gamma'(t)||dt$$

(2) Regarding now mathematicians, these would say that the length of a curve is the following quantity, with $(t_0=a,t_1,\ldots,t_{n-1},t_n=b)$ being a uniform division of $(a,b)$:
$$L(\gamma)=\lim_{n\to\infty}\sum_{i=1}^n||\gamma(t_i)-\gamma(t_{i-1})||$$

But, by using the fundamental theorem of calculus, this is the same as (2).

\medskip

(3) Getting back now to the ellipses, we can compute their length as follows:
\begin{eqnarray*}
L
&=&4\int_0^{\pi/2}\sqrt{\left(\frac{dx}{dt}\right)^2+\left(\frac{dy}{dt}\right)^2}\,dt\\
&=&4\int_0^{\pi/2}\sqrt{\left(\frac{da\cos t}{dt}\right)^2+\left(\frac{db\sin t}{dt}\right)^2}\,dt\\
&=&4\int_0^{\pi/2}\sqrt{a^2\sin^2t+b^2\cos^2t}\,dt
\end{eqnarray*}

(4) As for the last assertion, when $a=b=R$ we get of course $L=2\pi R$, as we should, but in general, when $a\neq b$, there is no trick for computing the above integral.
\end{proof}

Regarding now advanced integration and probability, in view of the above, we can certainly do some on the ellipses. But in practice, this will just bring some parameters $a,b>0$ in our Wallis formulae, and is this worth it. In addition, with the length of our ellipses being not computable, everything will be quite odd. So, time to stop here.

\bigskip

Moving on, and getting now to our final Question 12.3 (4), regarding the number $\pi$ itself, we have here the following famous 18th century result, due to Lambert: 

\begin{theorem}
The number $\pi=3.14159\ldots$ is irrational.
\end{theorem}

\begin{proof}
The original proof of Lambert was quite complicated, and we will explain here a more recent proof, from the 19th century, due to Hermite. Let us set:
$$I_n(t)=\int_{-1}^1(1-x^2)^n\cos(xt)dx$$

By double partial integration we obtain the following formula:
\begin{eqnarray*}
I_n(t)
&=&\int_{-1}^1(1-x^2)^n\cos(xt)dx\\
&=&\int_{-1}^12nx(1-x^2)^{n-1}\frac{\sin(xt)}{t}\,dx\\
&=&\frac{2n}{t}\int_{-1}^1x(1-x^2)^{n-1}\sin(xt)dx\\
&=&\frac{2n}{t}\int_{-1}^1\left[(1-x^2)^{n-1}-2(n-1)x^2(1-x^2)^{n-2}\right]\frac{\cos(xt)}{t}\,dx\\
&=&\frac{2n}{t^2}\int_{-1}^1(1-x^2)^{n-2}\left[1-x^2-2(n-1)x^2\right]\cos(xt)dx\\
&=&\frac{2n}{t^2}\int_{-1}^1(1-x^2)^{n-2}\left[1-(2n-1)x^2\right]\cos(xt)dx\\
&=&\frac{2n}{t^2}\int_{-1}^1(1-x^2)^{n-2}\left[(2n-1)(1-x^2)-(2n-2)\right]\cos(xt)dx\\
&=&\frac{2n}{t^2}\big[(2n-1)I_{n-1}(t)-(2n-2)I_{n-2}(t)\big]
\end{eqnarray*}

Thus, we have the following recurrence relation for our quantities:
$$t^2I_n(t)=2n(2n-1)I_{n-1}(t)-4n(n-1)I_{n-2}(t)$$

In terms of $J_n(t)=t^{2n+1}I_n(t)$, this recurrence formula becomes:
$$J_n(t)=2n(2n-1)J_{n-1}(t)-4n(n-1)t^2J_{n-2}(t)$$

Regarding now the initial data, for this latter recurrence, this is as follows:
$$J_0(t)=2\sin t\quad,\quad J_1(t)=-4t\cos t+4\sin t$$

We conclude from this that we must have a formula as follows, with $P_n,Q_n$ being certain polynomials of degree $\leq n$, with integer coefficients:
$$J_n(t)=n!(P_n(t)\sin t+Q_n(t)\cos t)$$

Now observe that with $t=\pi/2$, we obtain from this the following formula:
$$\left(\frac{\pi}{2}\right)^{2n+1}I_n\left(\frac{\pi}{2}\right)=n!P_n\left(\frac{\pi}{2}\right)$$

Assume now by contradiction that $\pi$ is rational, so that $\pi/2=a/b$ with $a,b\in\mathbb N$. We can rewrite the formula found above in the following more convenient way:
$$\frac{a^{2n+1}}{n!}\,I_n\left(\frac{a}{b}\right)=b^{2n+1}P_n\left(\frac{a}{b}\right)$$

But, by definition of the integrals $I_n$, we know that we have:
$$I_n\left(\frac{a}{b}\right)=I_n\left(\frac{\pi}{2}\right)\in(0,2)$$

Thus with $n>>0$ the number on the left belongs to $(0,1)$, which is contradictory, because the number on the right is an integer. And so $\pi$ is irrational, as claimed.
\end{proof}

Many other things can be said, as a continuation of the above, notably with the fact that $\pi$ is in fact transcendental, $P(\pi)\neq0$ for any $P\in\mathbb Z[X]$. For more on all this, which is quite fascinating mathematics, have a look at a basic number theory book.

\bigskip

As a last piece of mathematics, still in relation with the very nature of $\pi$ itself, we have the following remarkable probabilistic model for it, due to Buffon:

\begin{theorem}
The probability for a needle of length $1$, when trown on a grid of parallel $1$-spaced lines, to intersect one line, is:
$$P=\frac{2}{\pi}$$
Moreover, we have generalizations of this result, with needles of arbitrary length, thrown over a grid of parallel lines, with arbitrary spacing.
\end{theorem}

\begin{proof}
This is something quite tricky, and mandatory for properly learning probability theory, and science in general, because there are several possible modelings of the problem, leading, quite surprisingly, to different values of $P$. And, obviously, only one such modeling can be the correct one. So, here is how the correct computation goes:

\medskip

(1) Let us examine, right from the beginning, the general question involving throwing a needle of length $l$ on a grid of parallel $r$-spaced lines. We denote by $x\in[0,r/2]$ the distance from the middle of the needle to the closest line, and by $t\in[0,\pi/2]$ the acute angle formed by the needle and that line, according to the following picture:
$$
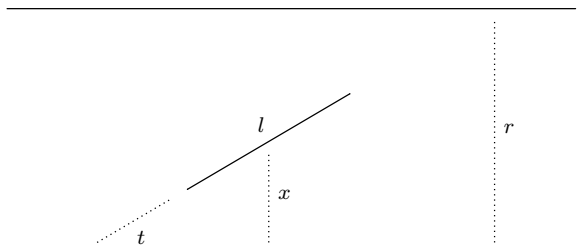
\xymatrix@R=3pt@C=10pt{
\ar@{-}[rrrrrrrrrrrrr]&&&&&&&&&&&\ar@{.}[ddddddddd]^r&&\\
\\
\\
&&&&&&&&\\
&&&&\\
&&&&&&\ar@{.}[dddd]^x\\
&&&\\
&&&&\ar@{.}[ddll]^t\ar@{-}[uuuurrrr]^l&\\
&&&&&\\
\ar@{-}[rrrrrrrrrrrrr]&&&&&&&&&&&&&}$$

(2) Now since the distance $x\in[0,r/2]$ varies uniformly, with density $2/r$, and the angle $t\in[0,\pi/2]$ varies uniformly too, with density $2/\pi$, and these variables are independent, in the case $l\leq r$ the probability for our needle to cross a line is given by:
\begin{eqnarray*}
P\left(x\leq\frac{l\sin t}{2}\right)
&=&\int_0^{\pi/2}\int_0^{l\sin t/2}\frac{4}{\pi r}\,dxdt\\
&=&\frac{4}{\pi r}\int_0^{\pi/2}\frac{l\sin t}{2}\,dt\\
&=&\frac{2l}{\pi r}
\end{eqnarray*}

(3) In particular, with $l=r=1$ the crossing probability is $P=2/\pi$, as stated.
\end{proof}

And with this, end of our study of the unit circle. Good learning this was, and in what follows we will try to understand what happens to this, in higher dimensions.

\section*{12c. Spherical integrals}

Getting now to higher dimensions, we potentially have many things to be done. To start with, as a straightforward extension of Principle 12.2, we have:

\begin{principle}
The interplay between the ambient space and its unit sphere
$$S^{N-1}_\mathbb R\subset\mathbb R^N$$
is the alpha and omega of everything, in $N$-dimensional mathematics.
\end{principle}

To be more precise, this is something that we already had a taste of, in chapter 8 when talking about the gamma function and related topics, and which remains to be further confirmed. Also, the notation $S^{N-1}_\mathbb R$, which might look a bit complicated, is justified by the fact that $N-1$ is the dimension of our sphere, viewed as manifold, and with the $\mathbb R$ subscript being something useful, because we will talk about $S^{N-1}_\mathbb C\subset\mathbb C^N$ too.

\bigskip

Getting now to the extension of Questions 12.3, many things can be said here, and going directly for the answers, the situation with this is as follows:

\begin{questions}
In relation with the space and its unit sphere, $S^{N-1}_\mathbb R\subset\mathbb R^N$:
\begin{enumerate}
\item We have spherical coordinates, and the stereographic projection too.

\item Using spherical coordinates and Wallis, we can integrate over $S^{N-1}_\mathbb R$.

\item We have ellipsoids, whose theory is similar to that of the ellipses.

\item Area and volumes teach us that there is no other $\pi$ than $\pi$.
\end{enumerate}
\end{questions}

More in detail now, in what regards the spherical coordinates, then integration via Wallis, and subsequent probability, which is related to our usual business in this book, continuous probability, this is something that we intend to discuss, right next.

\bigskip

As for the other topics, in what regards the stereographic projection and the ellipsoids, with these counting as basic differential and algebraic geometry in $\mathbb R^N$, we will leave some learning here as an exercise. As for the last assertion, this comes from our computations from chapter 8, which showed that, thanks God, the areas and volumes of the higher dimensional spheres $S^{N-1}_\mathbb R$ involve the good old $\pi$ only, and no other beasts. Nice.

\bigskip

Getting to work now, let us first talk about spherical coordinates. In 3D, to start with, you certainly know about latitude and longitude. Mathematically, the result is:

\index{spherical coordinates}

\begin{proposition}
We have spherical coordinates in $3$ dimensions,
$$\begin{cases}
x\!\!\!&=\ r\cos s\\
y\!\!\!&=\ r\sin s\cos t\\
z\!\!\!&=\ r\sin s\sin t
\end{cases}$$
the corresponding Jacobian being $J(r,s,t)=r^2\sin s$.
\end{proposition}

\begin{proof}
The Jacobian is indeed given by the following formula:
\begin{eqnarray*}
&&J(r,s,t)\\
&=&\begin{vmatrix}
\cos s&-r\sin s&0\\
\sin s\cos t&r\cos s\cos t&-r\sin s\sin t\\
\sin s\sin t&r\cos s\sin t&r\sin s\cos t
\end{vmatrix}\\
&=&r^2\sin s\sin t
\begin{vmatrix}\cos s&-r\sin s\\ \sin s\sin t&r\cos s\sin t\end{vmatrix}
+r\sin s\cos t\begin{vmatrix}\cos s&-r\sin s\\ \sin s\cos t&r\cos s\cos t\end{vmatrix}\\
&=&r\sin s\sin^2 t
\begin{vmatrix}\cos s&-r\sin s\\ \sin s&r\cos s\end{vmatrix}
+r\sin s\cos^2 t\begin{vmatrix}\cos s&-r\sin s\\ \sin s&r\cos s\end{vmatrix}\\
&=&r\sin s(\sin^2t+\cos^2t)\begin{vmatrix}\cos s&-r\sin s\\ \sin s&r\cos s\end{vmatrix}\\
&=&r\sin s\times 1\times r\\
&=&r^2\sin s
\end{eqnarray*}

Thus, we have indeed the formula in the statement.
\end{proof}

In general, the result, which generalizes those at $N=2,3$, is as follows:

\index{spherical coordinates}

\begin{theorem}
We have spherical coordinates in $N$ dimensions,
$$\begin{cases}
x_1\!\!\!&=\ r\cos t_1\\
x_2\!\!\!&=\ r\sin t_1\cos t_2\\
\vdots\\
x_{N-1}\!\!\!&=\ r\sin t_1\sin t_2\ldots\sin t_{N-2}\cos t_{N-1}\\
x_N\!\!\!&=\ r\sin t_1\sin t_2\ldots\sin t_{N-2}\sin t_{N-1}
\end{cases}$$
the corresponding Jacobian being given by the following formula,
$$J(r,t)=r^{N-1}\sin^{N-2}t_1\sin^{N-3}t_2\,\ldots\,\sin^2t_{N-3}\sin t_{N-2}$$
and with this generalizing the known formulae at $N=2,3$.
\end{theorem}

\begin{proof}
As before, the fact that we have spherical coordinates is clear. Regarding now the Jacobian, also as before, by developing over the last column, we have:
\begin{eqnarray*}
J_N
&=&r\sin t_1\ldots\sin t_{N-2}\sin t_{N-1}\times \sin t_{N-1}J_{N-1}\\
&+&r\sin t_1\ldots \sin t_{N-2}\cos t_{N-1}\times\cos t_{N-1}J_{N-1}\\
&=&r\sin t_1\ldots\sin t_{N-2}(\sin^2 t_{N-1}+\cos^2 t_{N-1})J_{N-1}\\
&=&r\sin t_1\ldots\sin t_{N-2}J_{N-1}
\end{eqnarray*}

Thus, we obtain the formula in the statement, by recurrence.
\end{proof}

As a comment here, the above convention for spherical coordinates is one among many, designed to best work in arbitrary $N$ dimensions. Also, in what regards the precise range of the angles $t_1,\ldots,t_{N-1}$, I will leave this to you, as an instructive exercise.

\bigskip

Good news, we can now integrate over the spheres, as follows:

\index{spherical integral}
\index{double factorials}

\begin{theorem}[Wallis 3]
The polynomial integrals over the unit sphere $S^{N-1}_\mathbb R\subset\mathbb R^N$, with respect to the uniform mass $1$ measure, are given by the formula
$$\int_{S^{N-1}_\mathbb R}x_1^{k_1}\ldots x_N^{k_N}\,dx=\frac{(N-1)!!k_1!!\ldots k_N!!}{(N+\Sigma k_i-1)!!}$$
valid when all exponents $k_i$ are even. If an exponent is odd, the integral vanishes.
\end{theorem}

\begin{proof}
This is something routine, generalizing Wallis 2, as follows:

\medskip

(1) Assume first that one of the exponents $k_i$ is odd. We can make then the following change of variables, which shows that the integral in the statement vanishes:
$$x_i\to-x_i$$

(2) Assume now that all the exponents $k_i$ are even. As a first observation, the result holds indeed at $N=2$, due to the Wallis formula from Theorem 12.7, which reads:
\begin{eqnarray*}
\int_0^{\pi/2}\cos^pt\sin^qt\,dt
&=&\left(\frac{\pi}{2}\right)^{\varepsilon(p)\varepsilon(q)}\frac{p!!q!!}{(p+q+1)!!}\\\
&=&\frac{p!!q!!}{(p+q+1)!!}
\end{eqnarray*}

Indeed, this formula computes the integral in the statement over the first quadrant. But since the exponents $p,q\in\mathbb N$ are assumed to be even, the integrals over the other quadrants are given by the same formula, so when averaging we obtain the result. 

\medskip

(3) In the general case now, where the dimension $N\in\mathbb N$ is arbitrary, the integral in the statement can be written in spherical coordinates, as follows:
$$I=\frac{2^N}{A}\int_0^{\pi/2}\ldots\int_0^{\pi/2}x_1^{k_1}\ldots x_N^{k_N}J\,dt_1\ldots dt_{N-1}$$

Here $A$ is the area of the sphere, $J$ is the Jacobian, and the $2^N$ factor comes from the restriction to the $1/2^N$ part of the sphere where all coordinates are positive. According to our formulae from chapter 8, the normalization constant in front of the integral is:
$$\frac{2^N}{A}=\left(\frac{2}{\pi}\right)^{[N/2]}(N-1)!!$$

As for the unnormalized integral, by using the formulae from Theorem 12.17, for the spherical coordinates and their Jacobian, this is given by the following formula:
\begin{eqnarray*}
I'=\int_0^{\pi/2}\ldots\int_0^{\pi/2}
&&(\cos t_1)^{k_1}
(\sin t_1\cos t_2)^{k_2}\\
&&\vdots\\
&&(\sin t_1\sin t_2\ldots\sin t_{N-2}\cos t_{N-1})^{k_{N-1}}\\
&&(\sin t_1\sin t_2\ldots\sin t_{N-2}\sin t_{N-1})^{k_N}\\
&&\sin^{N-2}t_1\sin^{N-3}t_2\ldots\sin^2t_{N-3}\sin t_{N-2}\\
&&dt_1\ldots dt_{N-1}
\end{eqnarray*}

(4) By rearranging the terms, we obtain the following formula:
\begin{eqnarray*}
I'
&=&\int_0^{\pi/2}\cos^{k_1}t_1\sin^{k_2+\ldots+k_N+N-2}t_1\,dt_1\\
&&\int_0^{\pi/2}\cos^{k_2}t_2\sin^{k_3+\ldots+k_N+N-3}t_2\,dt_2\\
&&\vdots\\
&&\int_0^{\pi/2}\cos^{k_{N-2}}t_{N-2}\sin^{k_{N-1}+k_N+1}t_{N-2}\,dt_{N-2}\\
&&\int_0^{\pi/2}\cos^{k_{N-1}}t_{N-1}\sin^{k_N}t_{N-1}\,dt_{N-1}
\end{eqnarray*}

Now by using the above-mentioned Wallis formula at $N=2$, this gives:
\begin{eqnarray*}
I'
&=&\frac{k_1!!(k_2+\ldots+k_N+N-2)!!}{(k_1+\ldots+k_N+N-1)!!}\left(\frac{\pi}{2}\right)^{\varepsilon(N-2)}\\
&&\frac{k_2!!(k_3+\ldots+k_N+N-3)!!}{(k_2+\ldots+k_N+N-2)!!}\left(\frac{\pi}{2}\right)^{\varepsilon(N-3)}\\
&&\vdots\\
&&\frac{k_{N-2}!!(k_{N-1}+k_N+1)!!}{(k_{N-2}+k_{N-1}+l_N+2)!!}\left(\frac{\pi}{2}\right)^{\varepsilon(1)}\\
&&\frac{k_{N-1}!!k_N!!}{(k_{N-1}+k_N+1)!!}\left(\frac{\pi}{2}\right)^{\varepsilon(0)}
\end{eqnarray*}

(5) Now let $F$ be the part involving the double factorials, and $P$ be the part involving the powers of $\pi/2$, so that $I'=F\cdot P$. Regarding $F$, by cancelling terms we have:
$$F=\frac{k_1!!\ldots k_N!!}{(\Sigma k_i+N-1)!!}$$

As in what regards $P$, by summing the exponents, we obtain $P=\left(\frac{\pi}{2}\right)^{[N/2]}$. We can now put everything together, and we obtain:
\begin{eqnarray*}
I
&=&\frac{2^N}{A}\times F\times P\\
&=&\left(\frac{2}{\pi}\right)^{[N/2]}(N-1)!!\times\frac{k_1!!\ldots k_N!!}{(\Sigma k_i+N-1)!!}\times\left(\frac{\pi}{2}\right)^{[N/2]}\\
&=&\frac{(N-1)!!k_1!!\ldots k_N!!}{(\Sigma k_i+N-1)!!}
\end{eqnarray*}

Thus, we are led to the conclusion in the statement.
\end{proof}

Next, we have the following useful generalization of the above formula:

\index{spherical integral}

\begin{theorem}[Wallis 3']
We have the following integration formula over the sphere $S^{N-1}_\mathbb R\subset\mathbb R^N$, with respect to the uniform measure, valid for any exponents $k_i\in\mathbb N$,
$$\int_{S^{N-1}_\mathbb R}|x_1^{k_1}\ldots x_N^{k_N}|\,dx=\left(\frac{2}{\pi}\right)^{\Sigma(k_1,\ldots,k_N)}\frac{(N-1)!!k_1!!\ldots k_N!!}{(N+\Sigma k_i-1)!!}$$
with $\Sigma=[odds/2]$ if $N$ is odd and $\Sigma=[(odds+1)/2]$ if $N$ is even, where ``odds'' denotes the number of odd numbers in the sequence $k_1,\ldots,k_N$.
\end{theorem}

\begin{proof}
As before, the formula holds at $N=2$, due to Theorem 12.7. In general, the integral in the statement can be written in spherical coordinates, as follows:
$$I=\frac{2^N}{A}\int_0^{\pi/2}\ldots\int_0^{\pi/2}x_1^{k_1}\ldots x_N^{k_N}J\,dt_1\ldots dt_{N-1}$$

Here $A$ is the area of the sphere, $J$ is the Jacobian, and the $2^N$ factor comes from the restriction to the $1/2^N$ part of the sphere where all the coordinates are positive. The normalization constant in front of the integral is, as before, given by:
$$\frac{2^N}{A}=\left(\frac{2}{\pi}\right)^{[N/2]}(N-1)!!$$

As for the unnormalized integral, this can be written, as before, as follows:
\begin{eqnarray*}
I'
&=&\int_0^{\pi/2}\cos^{k_1}t_1\sin^{k_2+\ldots+k_N+N-2}t_1\,dt_1\\
&&\int_0^{\pi/2}\cos^{k_2}t_2\sin^{k_3+\ldots+k_N+N-3}t_2\,dt_2\\
&&\vdots\\
&&\int_0^{\pi/2}\cos^{k_{N-2}}t_{N-2}\sin^{k_{N-1}+k_N+1}t_{N-2}\,dt_{N-2}\\
&&\int_0^{\pi/2}\cos^{k_{N-1}}t_{N-1}\sin^{k_N}t_{N-1}\,dt_{N-1}
\end{eqnarray*}

Now by using the Wallis 2 formula, from Theorem 12.7, we obtain:
\begin{eqnarray*}
I'
&=&\frac{\pi}{2}\cdot\frac{k_1!!(k_2+\ldots+k_N+N-2)!!}{(k_1+\ldots+k_N+N-1)!!}\left(\frac{2}{\pi}\right)^{\delta(k_1,k_2+\ldots+k_N+N-2)}\\
&&\frac{\pi}{2}\cdot\frac{k_2!!(k_3+\ldots+k_N+N-3)!!}{(k_2+\ldots+k_N+N-2)!!}\left(\frac{2}{\pi}\right)^{\delta(k_2,k_3+\ldots+k_N+N-3)}\\
&&\vdots\\
&&\frac{\pi}{2}\cdot\frac{k_{N-2}!!(k_{N-1}+k_N+1)!!}{(k_{N-2}+k_{N-1}+k_N+2)!!}\left(\frac{2}{\pi}\right)^{\delta(k_{N-2},k_{N-1}+k_N+1)}\\
&&\frac{\pi}{2}\cdot\frac{k_{N-1}!!k_N!!}{(k_{N-1}+k_N+1)!!}\left(\frac{2}{\pi}\right)^{\delta(k_{N-1},k_N)}
\end{eqnarray*}

In order to compute this quantity, let us denote by $F$ the part involving the double factorials, and by $P$ the part involving the powers of $\pi/2$, so that we have:
$$I'=F\cdot P$$

Regarding $F$, there are many cancellations there, and we end up with:
$$F=\frac{k_1!!\ldots k_N!!}{(\Sigma k_i+N-1)!!}$$

As in what regards $P$, the $\delta$ exponents on the right sum up to the following number:
$$\Delta(k_1,\ldots,k_N)=\sum_{i=1}^{N-1}\delta(k_i,k_{i+1}+\ldots+k_N+N-i-1)$$

In other words, with this notation, the above formula reads:
\begin{eqnarray*}
I'
&=&\left(\frac{\pi}{2}\right)^{N-1}\frac{k_1!!k_2!!\ldots k_N!!}{(k_1+\ldots+k_N+N-1)!!}\left(\frac{2}{\pi}\right)^{\Delta(k_1,\ldots,k_N)}\\
&=&\left(\frac{2}{\pi}\right)^{\Delta(k_1,\ldots,k_N)-N+1}\frac{k_1!!k_2!!\ldots k_N!!}{(k_1+\ldots+k_N+N-1)!!}\\
&=&\left(\frac{2}{\pi}\right)^{\Sigma(k_1,\ldots,k_N)-[N/2]}\frac{k_1!!k_2!!\ldots k_N!!}{(k_1+\ldots+k_N+N-1)!!}
\end{eqnarray*}

To be more precise, the formula relating $\Delta$ to $\Sigma$ follows from a number of simple observations, the first of which being the fact that, due to obvious parity reasons, the sequence of $\delta$ numbers appearing in the definition of $\Delta$ cannot contain two consecutive zeroes. Now together with $I=(2^N/V)I'$, this gives the formula in the statement.
\end{proof}

Finally, we have the following complex version of Theorems 12.18 and 12.19:

\index{spherical integral}

\begin{theorem}[Wallis 4]
We have the following integration formula over the complex sphere $S^{N-1}_\mathbb C\subset\mathbb C^N$, with respect to the uniform mass $1$ measure, 
$$\int_{S^{N-1}_\mathbb C}|z_1|^{2k_1}\ldots|z_N|^{2k_N}\,dz=\frac{(N-1)!k_1!\ldots k_n!}{(N+\sum k_i-1)!}$$
valid for any exponents $k_i\in\mathbb N$. As for the other polynomial integrals in $z_1,\ldots,z_N$ and their conjugates $\bar{z}_1,\ldots,\bar{z}_N$, these all vanish.
\end{theorem}

\begin{proof}
Consider a polynomial integral over $S^{N-1}_\mathbb C$, containing the same number of plain and conjugated variables, as to not vanish trivially, written as follows:
$$I=\int_{S^{N-1}_\mathbb C}z_{i_1}\bar{z}_{i_2}\ldots z_{i_{2k-1}}\bar{z}_{i_{2k}}\,dz$$

By using transformations of type $p\to\lambda p$ with $|\lambda|=1$, we see that this integral $I$ vanishes, unless each $z_a$ appears as many times as $\bar{z}_a$ does, and this gives the last assertion. So, assume now that we are in the non-vanishing case. Then the $k_a$ copies of $z_a$ and the $k_a$ copies of $\bar{z}_a$ produce by multiplication a factor $|z_a|^{2k_a}$, so we have:
$$I=\int_{S^{N-1}_\mathbb C}|z_1|^{2k_1}\ldots|z_N|^{2k_N}\,dz$$

Now by using the standard identification $S^{N-1}_\mathbb C\simeq S^{2N-1}_\mathbb R$, we obtain:
\begin{eqnarray*}
I
&=&\int_{S^{2N-1}_\mathbb R}(x_1^2+y_1^2)^{k_1}\ldots(x_N^2+y_N^2)^{k_N}\,d(x,y)\\
&=&\sum_{r_1\ldots r_N}\binom{k_1}{r_1}\ldots\binom{k_N}{r_N}\int_{S^{2N-1}_\mathbb R}x_1^{2k_1-2r_1}y_1^{2r_1}\ldots x_N^{2k_N-2r_N}y_N^{2r_N}\,d(x,y)
\end{eqnarray*}

By using the formula in Theorem 12.18 for integrating, we obtain:
\begin{eqnarray*}
&&I\\
&=&\sum_{r_1\ldots r_N}\binom{k_1}{r_1}\ldots\binom{k_N}{r_N}\frac{(2N-1)!!(2r_1)!!\ldots(2r_N)!!(2k_1-2r_1)!!\ldots (2k_N-2r_N)!!}{(2N+2\sum k_i-1)!!}\\
&=&\sum_{r_1\ldots r_N}\binom{k_1}{r_1}\ldots\binom{k_N}{r_N}\frac{2^{N-1}(N-1)!\prod(2r_i)!/(2^{r_i}r_i!)\prod(2k_i-2r_i)!/(2^{k_i-r_i}(k_i-r_i)!)}{2^{N+\sum k_i-1}(N+\sum k_i-1)!}\\
&=&\sum_{r_1\ldots r_N}\binom{k_1}{r_1}\ldots\binom{k_N}{r_N}
\frac{(N-1)!(2r_1)!\ldots (2r_N)!(2k_1-2r_1)!\ldots (2k_N-2r_N)!}{4^{\sum k_i}(N+\sum k_i-1)!r_1!\ldots r_N!(k_1-r_1)!\ldots (k_N-r_N)!}
\end{eqnarray*}

Now observe that can rewrite this quantity in the following way:
\begin{eqnarray*}
&&I\\
&=&\sum_{r_1\ldots r_N}\frac{k_1!\ldots k_N!(N-1)!(2r_1)!\ldots (2r_N)!(2k_1-2r_1)!\ldots (2k_N-2r_N)!}{4^{\sum k_i}(N+\sum k_i-1)!(r_1!\ldots r_N!(k_1-r_1)!\ldots (k_N-r_N)!)^2}\\
&=&\sum_{r_1}\binom{2r_1}{r_1}\binom{2k_1-2r_1}{k_1-r_1}\ldots\sum_{r_N}\binom{2r_N}{r_N}\binom{2k_N-2r_N}{k_N-r_N}\frac{(N-1)!k_1!\ldots k_N!}{4^{\sum k_i}(N+\sum k_i-1)!}\\
&=&4^{k_1}\times\ldots\times 4^{k_N}\times\frac{(N-1)!k_1!\ldots k_N!}{4^{\sum k_i}(N+\sum k_i-1)!}\\
&=&\frac{(N-1)!k_1!\ldots k_N!}{(N+\sum k_i-1)!}
\end{eqnarray*}

To be more precise, we have used here the following formula, that we know well from chapter 10, from our study there of the complex normal variables:
$$\sum_r\binom{2r}{r}\binom{2l-2r}{l-r}=4^l$$

Thus, we have obtained the formula in the statement, as desired.
\end{proof}

Summarizing, we have complete results for the integration over the spheres, with the answers involving various multinomial type coefficients, defined in terms of factorials, or of double factorials. All these formulae are of course very useful, in practice.

\section*{12d. Hyperspherical laws}

Good news, we can now do some interesting probability over the spheres. Our source of inspiration will be Theorem 12.9, dealing with the 2D real case, and as technical tools, we will use the Wallis 3 and 4 formulae, from Theorems 12.18 and 12.20.

\bigskip

In the real case, many things can be said, and as a main result here, we have:

\index{hyperspherical law}
\index{normal law}

\begin{theorem}
The even moments of hyperspherical variables are given by
$$\int_{S^{N-1}_\mathbb R}x_i^kdx=\frac{(N-1)!!k!!}{(N+k-1)!!}$$
and the odd moments vanish. The normalized hyperspherical variables 
$$y_i=\sqrt{N}x_i$$
become standard normal, following $g_1$, and independent with $N\to\infty$.
\end{theorem}

\begin{proof}
This is something standard, based on Theorem 12.18, as follows:

\medskip

(1) To start with, as before in Theorem 12.9, we can certainly denote the points of the sphere $x\in S^{N-1}_\mathbb R$ as $x=(x_1,\ldots,x_N)$, and therefore talk about the $N$ random variables $x_i:S^{N-1}_\mathbb R\to[-1,1]$, with respect to the uniform, mass 1 measure on the sphere.

\medskip

(2) In order to study the individual and joint laws of these variables $x_1,\ldots,x_N$, and their asymptotics, we can use the Wallis 3 formula, from Theorem 12.18, namely:
$$\int_{S^{N-1}_\mathbb R}x_1^{k_1}\ldots x_N^{k_N}\,dx=\frac{(N-1)!!k_1!!\ldots k_N!!}{(N+\Sigma k_i-1)!!}$$

(3) As a first observation, by assuming that all exponents are 0, except for the $i$-th one, we obtain the individual moment formula in the statement, namely:
$$\int_{S^{N-1}_\mathbb R}x_i^kdx=\frac{(N-1)!!k!!}{(N+k-1)!!}$$

(4) Next, observe that by using this, we have the following estimate, with $N\to\infty$:
\begin{eqnarray*}
\int_{S^{N-1}_\mathbb R}x_i^kdx
&=&\frac{(N-1)!!}{(N+k-1)!!}\times k!!\\
&\simeq&N^{-k/2}k!!\\
&=&N^{-k/2}M_k(g_1)
\end{eqnarray*}

Thus the variables $y_i=\sqrt{N}x_i$ become standard normal with $N\to\infty$, as stated.

\medskip

(5) Getting now to the joint moments of $y_i=\sqrt{N}x_i$, with $N\to\infty$ we have:
\begin{eqnarray*}
\int_{S^{N-1}_\mathbb R}y_1^{k_1}\ldots y_N^{k_N}\,dx
&=&N^{\sum k_i/2}\int_{S^{N-1}_\mathbb R}x_1^{k_1}\ldots x_N^{k_N}\,dx\\
&=&N^{\sum k_i/2}\frac{(N-1)!!k_1!!\ldots k_N!!}{(N+\sum k_i-1)!!}\\
&=&N^{\sum k_i/2}\frac{(N-1)!!}{(N+\sum k_i-1)!!}\times k_1!!\ldots k_N!!\\
&\simeq&N^{\sum k_i/2}N^{-\sum k_i/2}k_1!!\ldots k_N!!\\
&=&k_1!!\ldots k_N!!
\end{eqnarray*}

Thus the asymptotic joint moments factorize, and by using the moment characterization of independence we conclude that $y_i$ are asymptotically independent, as stated.
\end{proof}

In the complex case we have a similar result, as follows:

\begin{theorem}
The moments of complex hyperspherical variables are given by
$$\int_{S^{N-1}_\mathbb C}|z_i|^{2k}dz=\frac{(N-1)!k!}{(N+k-1)!}$$
with the other moments vanishing. The normalized complex hyperspherical variables 
$$w_i=\sqrt{N}z_i$$
become standard complex normal, following $G_1$, and independent with $N\to\infty$.
\end{theorem}

\begin{proof}
This is a standard complex remake of Theorem 12.21, as follows:

\medskip

(1) To start with, we can certainly denote the points of the complex sphere $z\in S^{N-1}_\mathbb C$ as $z=(z_1,\ldots,z_N)$, and therefore talk about the $N$ random variables $z_i:S^{N-1}_\mathbb C\to D$, with $D\subset\mathbb C$ being the unit disk, with respect to the uniform measure on the sphere.

\medskip

(2) In order to study the individual and joint laws of these variables $z_1,\ldots,z_N$, and their asymptotics, we can use the Wallis 4 formula, from Theorem 12.20, namely:
$$\int_{S^{N-1}_\mathbb C}|z_1|^{2k_1}\ldots|z_N|^{2k_N}\,dz=\frac{(N-1)!k_1!\ldots k_n!}{(N+\sum k_i-1)!}$$

(3) As a first observation, by assuming that all exponents are 0, except for the $i$-th one, we obtain the individual moment formula in the statement, namely:
$$\int_{S^{N-1}_\mathbb C}|z_i|^{2k}dz=\frac{(N-1)!k!}{(N+k-1)!}$$

(4) Next, observe that by using this, we have the following estimate, with $N\to\infty$:
$$\int_{S^{N-1}_\mathbb C}|z_i|^{2k}dz
=\frac{(N-1)!}{(N+k-1)!}\times k!
\simeq N^{-k}k!
=N^{-k}M_{2k}(G_1)$$

Thus $w_i=\sqrt{N}z_i$ become standard complex normal with $N\to\infty$, as stated.

\medskip

(5) Getting now to the joint moments of $w_i=\sqrt{N}z_i$, with $N\to\infty$ we have:
\begin{eqnarray*}
\int_{S^{N-1}_\mathbb C}|w_1|^{2k_1}\ldots |w_N|^{2k_N}\,dz
&=&N^{\sum k_i}\int_{S^{N-1}_\mathbb C}|z_1|^{2k_1}\ldots|z_N|^{2k_N}\,dz\\
&=&N^{\sum k_i}\frac{(N-1)!k_1!\ldots k_N!}{(N+\sum k_i-1)!}\\
&=&N^{\sum k_i}\frac{(N-1)!}{(N+\sum k_i-1)!}\times k_1!\ldots k_N!\\
&\simeq&N^{\sum k_i}N^{-\sum k_i}k_1!\ldots k_N!\\
&=&k_1!\ldots k_N!
\end{eqnarray*}

Thus the joint moments factorize, and $w_i$ are asymptotically independent.
\end{proof}

Summarizing, we have nice asymptotic results, both in the real and the complex cases. Getting now to the fixed $N\in\mathbb N$ case, we already know from Theorem 12.9 that at $N=2$, in the real case, we obtain arcsine laws. Another case of interest, that we would like to discuss now in some detail, is $N=4$, also in the real case. We have here:

\begin{theorem}
The even moments of the standard coordinates of $S^3_\mathbb R\subset\mathbb R^4$ are
$$\int_{S^3_\mathbb R}x_i^{2r}dx=\frac{1}{4^r}\cdot\frac{1}{r+1}\binom{2r}{r}$$
the odd moments vanish, and the doubled coordinates follow the Wigner law $\gamma_1$:
$$2x_i\sim\frac{1}{2\pi}\sqrt{4-x^2}\,dx$$
Moreover, in the standard picture of the group $SU_2$, which is as follows,
$$SU_2=\left\{\begin{pmatrix}x+iy&z+it\\ -z+it&x-iy\end{pmatrix}\Big|(x,y,z,t)\in S^3_\mathbb R\right\}$$ 
the main character is one of these doubled coordinates, $\chi=2x$, following $\gamma_1$.
\end{theorem}

\begin{proof}
Many things going on here, the idea being as follows:

\medskip

(1) In even dimensions, $N=2n$, the moment formula in Theorem 12.21 reads:
\begin{eqnarray*}
\int_{S^{2n-1}_\mathbb R}x_i^{2r}dx
&=&\frac{(2n-1)!!(2r)!!}{(2n+2r-1)!!}\\
&=&\frac{2^{n-1}(n-1)!\cdot (2r)!/(2^rr!)}{2^{n+r-1}(n+r-1)!}\\
&=&\frac{1}{4^r}\cdot\frac{(n-1)!(2r)!}{r!(n+r-1)!}
\end{eqnarray*}

In particular, at $n=2$ we obtain the moment formula in the statement, namely:
$$\int_{S^3_\mathbb R}x_i^{2r}dx=\frac{1}{4^r}\cdot\frac{(2r)!}{r!(r+1)!}=\frac{1}{4^r}\cdot\frac{1}{r+1}\binom{2r}{r}$$

Now the moments of $2x_i$ being the Catalan numbers, we have $2x_i\sim\gamma_1$, as stated.

\medskip

(2) Regarding now the second assertion, and with my apologies for the conflicting notations for coordinates, but $(x,y,z,t)\in S^3_\mathbb R$ is very standard in the $SU_2$ context, at the beginning, we have the defining formula for the group $SU_2$, namely:
$$SU_2=\left\{U\in M_2(\mathbb C)\Big|U^*=U^{-1},\,\det U=1\right\}$$

In practice, if we pick $U=\binom{a\ b}{c\ d}$ of determinant 1, the equation $U^*=U^{-1}$ reads:
$$\begin{pmatrix}\bar{a}&\bar{c}\\ \bar{b}&\bar{d}\end{pmatrix}
=\begin{pmatrix}d&-b\\ -c&a\end{pmatrix}$$

Thus we must have $d=\bar{a}$, $c=-\bar{b}$, and since with these conditions imposed, the condition $\det U=1$ reads $|a|^2+|b|^2=1$, we are led to the following conclusion:
$$SU_2=\left\{\begin{pmatrix}a&b\\ -\bar{b}&\bar{a}\end{pmatrix}\Big|\ |a|^2+|b|^2=1\right\}$$

But with $a=x+iy$, $b=z+it$ this gives the formula in the statement, namely:
$$SU_2=\left\{\begin{pmatrix}x+iy&z+it\\ -z+it&x-iy\end{pmatrix}\Big|(x,y,z,t)\in S^3_\mathbb R\right\}$$

Summarizing, we have $SU_2\simeq S^3_\mathbb R$, and it is also known that the uniform measure on $SU_2$, in a group-theoretical sense, coincides with the uniform measure on $S^3_\mathbb R$.

\medskip

(3) But with this, done, because if we look at the main character of $SU_2$, in the spirit of the computations done in chapter 4 for some other matrix groups, this character is $\chi=2x$ in the above $SU_2\simeq S^3_\mathbb R$ picture, which by (1) follows the Wigner law $\gamma_1$.
\end{proof}

And we will end this chapter with this. Needless to say, Theorem 12.23 is just a beginning, and many other things can be said, as a continuation of that. More later.

\section*{12e. Exercises}

Welcome to geometry and probability, and as exercises on all this, we have:

\begin{exercise}
Learn more about the work of Fresnel, and his integrals.
\end{exercise}

\begin{exercise}
Look up or find applications of the stereographic projections.
\end{exercise}

\begin{exercise}
Have everything understood, regarding the ellipses, and conics.
\end{exercise}

\begin{exercise}
Review 3D orientation and geography, latitude and longitude.
\end{exercise}

\begin{exercise}
Then do the same, for spherical coordinates in $N$ dimensions.
\end{exercise}

\begin{exercise}
Have a closer look at hyperspherical laws, in even dimensions.
\end{exercise}

\begin{exercise}
Have a closer look at hyperspherical laws, in odd dimensions.
\end{exercise}

\begin{exercise}
Have as well a closer look at what happens in the complex case.
\end{exercise}

As bonus exercise, learn about $SU_2$, $SO_3$ and related topics, as much as you can.

\part{Advanced aspects}

\ \vskip50mm

\begin{center}
{\em Swing low, sweet chariot

Coming for to carry me home

Swing low, sweet chariot

Coming for to carry me home}
\end{center}

\chapter{Orthogonal polynomials}

\section*{13a. Hilbert spaces}

Welcome to this final Part IV, and what is to be done here? Most likely, more specialized laws, and more specialized techniques too. And this will be indeed what we will be talking about, over the next 100 pages. Hang on, tough material to come.

\bigskip

To start with, we need some sort of trick, for truly taking off, from what we did so far in this book. And the trick, which is standard in advanced mathematics, is that of using infinite dimensionality. So, let us formulate the following definition:

\index{Hilbert space}
\index{scalar product}

\begin{definition}
A Hilbert space is a complex vector space $H$, with a scalar product on it $<x,y>$, taken linear at left and antilinear at right,
$$<\lambda x,y>=\lambda<x,y>\quad,\quad <x,\lambda y>=\bar{\lambda}<x,y>$$
which is complete with respect to corresponding norm
$$||x||=\sqrt{<x,x>}$$
in the sense that any sequence $\{x_n\}$ which is a Cauchy sequence, having the property $||x_n-x_m||\to0$ with $n,m\to\infty$, has a limit, $x_n\to x$.
\end{definition}

Here the fact that $||x||=\sqrt{<x,x>}$ is indeed a norm, satisfying $||x+y||\leq||x||+||y||$, follows from Cauchy-Schwarz, $|\!<x,y>\!|\leq||x||\cdot||y||$, which itself follows as in the $H=\mathbb C^N$ case, by looking at $f(t)=||wx+ty||^2$, as explained in chapter 9. At the level of examples of Hilbert spaces, generalizing the basic $H=\mathbb C^N$ example, we have:

\index{square-summable}
\index{Cauchy-Schwarz}

\begin{proposition}
Given an index set $I$, which can be finite or not, the space of square-summable vectors having indices in $I$, namely
$$l^2(I)=\left\{(x_i)_{i\in I}\Big|\sum_i|x_i|^2<\infty\right\}$$
is a Hilbert space, with scalar product as follows:
$$<x,y>=\sum_ix_i\bar{y}_i$$
When $I$ is finite, $I=\{1,\ldots,N\}$, we obtain in this way the usual space $H=\mathbb C^N$.
\end{proposition}

\begin{proof}
This is something very standard, the idea being as follows:

\medskip

(1) We know that $l^2(I)\subset\mathbb C^I$ is the space of vectors satisfying $||x||<\infty$. We want to prove that $l^2(I)$ is a vector space, that $<x,y>$ is a scalar product on it, that $l^2(I)$ is complete with respect to $||.||$, and finally that for $|I|<\infty$ we have $l^2(I)=\mathbb C^{|I|}$.

\medskip

(2) The last assertion, $l^2(I)=\mathbb C^{|I|}$ for $|I|<\infty$, is clear, because in this case the sums are finite, so the condition $||x||<\infty$ is automatic. So, we know at least one thing.

\medskip

(3) Regarding the rest, we have Cauchy-Schwarz, $|<x,y>|\leq||x||\cdot||y||$, proved as usual, by looking at $f(t)=||wx+ty||^2$, and we deduce that we have $||x+y||\leq||x||+||y||$. Thus $l^2(I)$ is indeed a vector space, the other vector space conditions being trivial.

\medskip

(4) Also, the quantity $<x,y>$ is surely a scalar product on this vector space, because all the conditions for a scalar product are trivially satisfied.

\medskip

(5) Finally, the fact that our space $l^2(I)$ is indeed complete with respect to its norm $||.||$ follows in the obvious way, the limit of a Cauchy sequence $\{x_n\}$ being the vector $y=(y_i)$ given by $y_i=\lim_{n\to\infty}x_{ni}$, with all the verifications here being trivial.
\end{proof}

Going now a bit more abstract, we have, more generally, the following result:

\index{square-summable}
\index{square-integrable}

\begin{theorem}
Given an arbitrary space $X$ with a positive measure $\mu$ on it, the space of square-summable complex functions on it, namely
$$L^2(X)=\left\{f:X\to\mathbb C\Big|\int_X|f(x)|^2\,d\mu(x)<\infty\right\}$$
is a Hilbert space, with scalar product as follows:
$$<f,g>=\int_Xf(x)\overline{g(x)}\,d\mu(x)$$
When $X=I$ is discrete, meaning that the measure $\mu$ on it is the counting measure, $\mu(\{x\})=1$ for any $x\in X$, we obtain in this way the previous spaces $l^2(I)$.
\end{theorem}

\begin{proof}
This is something routine, remake of Proposition 13.2, as follows:

\medskip

(1) The proof of the first, and main assertion is something perfectly similar to the proof of Proposition 13.2, by replacing everywhere the sums by integrals. 

\medskip

(2) With the remark that I forgot to say in the statement that the $L^2$ functions are by definition taken up to equality almost everywhere, $f=g$ when $||f-g||=0$.

\medskip

(3) As for the last assertion, when $\mu$ is the counting measure all our integrals here become usual sums, and we recover in this way Proposition 13.2.
\end{proof}

At the level of the general theory now, the idea is that the various things that we know from chapter 9, regarding the geometry of $\mathbb C^N$, and of the linear maps $f:\mathbb C^N\to\mathbb C^N$ too, have extensions to the Hilbert space setting. Leaving aside the linear maps, that we will discuss later, in chapter 14, let us record some geometry basics, as follows:

\begin{proposition}
Given a Hilbert space $H$, with scalar product $<x,y>$ and norm $||x||=\sqrt{<x,x>}$, the norms, or lengths, are subject to the parallelogram identity
$$||x+y||^2+||x-y||^2=2(||x||^2+||y||^2)$$
and the scalar products can be recovered from lengths via the polarization identity:
$$<x,y>
=\frac{||x+y||^2-||x-y||^2+i||x+iy||^2-i||x-iy||^2}{4}$$
Also, we can say that two vectors are orthogonal when their scalar product is zero.
\end{proposition}

\begin{proof}
This is something that we know from chapter 9, for $H=\mathbb C^N$, and the proof in general is identical. To be more precise, the parallelogram rule comes from:
\begin{eqnarray*}
||x+y||^2+||x-y||^2
&=&<x+y,x+y>+<x-y,x-y>\\
&=&2||x||^2+2||y||^2+2Re(<x,y>)-2Re(<x,y>)\\
&=&2(||x||^2+||y||^2)
\end{eqnarray*}

As for the polarization identity, this can be established as follows:
\begin{eqnarray*}
&&||x+y||^2-||x-y||^2+i||x+iy||^2-i||x-iy||^2\\
&=&||x||^2+||y||^2-||x||^2-||y||^2+i||x||^2+i||y||^2-i||x||^2-i||y||^2\\
&&+2Re(<x,y>)+2Re(<x,y>)+2iIm(<x,y>)+2iIm(<x,y>)\\
&=&4<x,y>
\end{eqnarray*}

Finally, the last assertion is a convention, coming from what happens in $H=\mathbb C^N$, which itself was a convention too, coming from what happens in the real space $\mathbb R^N$.
\end{proof}

At a more advanced level now, the idea is that we can talk about orthonormal bases, and the related notion of dimension of a Hilbert space. Let us start with:

\begin{theorem}
Any system of linearly independent vectors $\{f_1,\ldots,f_n\}$ can be turned into an orthogonal system $\{e_1,\ldots,e_n\}$ by using the Gram-Schmidt procedure,
$$e_1=f_1$$
$$e_2=f_2+\alpha_1f_1$$
$$e_3=f_3+\beta_1f_1+\beta_2f_2$$
$$e_4=f_4+\gamma_1f_1+\gamma_2f_2+\gamma_3f_3$$
$$\vdots$$
with the needed scalars $\alpha_i,\beta_i,\gamma_i,\ldots$ being uniquely determined.
\end{theorem}

\begin{proof}
Many things can be said here, depending on how sharp you want to be, with the essentials of what is to be known being as follows:

\medskip

(1) Let us first study the case $n=2$. With $e_1=f_1$ and $e_2=f_2+\alpha_1f_1$ as in the statement, the needed orthogonality condition can be processed as follows:
\begin{eqnarray*}
e_1\perp e_2
&\iff&<f_1,f_2+\alpha_1f_1>=0\\
&\iff&\alpha_1<f_1,f_1>=-<f_1,f_2>\\
&\iff&\alpha_1=-\frac{<f_1,f_2>}{<f_1,f_1>}
\end{eqnarray*}

Thus, we get our result, and with the remark that, alternatively, we can set:
$$e_2=f_2-Proj_{e_1}(f_2)$$

Indeed, with the above formula of $\alpha_1$ in hand, the vector $e_2=f_2+\alpha_1f_1$ that we get is precisely this one. Or, we can simply argue that this latter vector $e_2$ does the job, and with some basic linear algebra telling us that this vector $e_2$ is indeed unique.

\medskip 

(2) At $n=3$ now, with $e_1,e_2$ already constructed, and with $e_3=f_3+\beta_1f_1+\beta_2f_2$ as in the statement, the first orthogonality condition can be processed as follows:
\begin{eqnarray*}
e_1\perp e_3
&\iff&<f_1,f_3+\beta_1f_1+\beta_2f_2>=0\\
&\iff&\beta_1<f_1,f_1>+\beta_2<f_1,f_2>=-<f_1,f_3>
\end{eqnarray*}

As for the second orthogonality condition, this can be now processed as follows:
\begin{eqnarray*}
e_2\perp e_3
&\iff&<f_2,f_3+\beta_1f_1+\beta_2f_2>=0\\
&\iff&\beta_1<f_2,f_1>+\beta_2<f_2,f_2>=-<f_2,f_3>
\end{eqnarray*}

Thus, we are led to the following system, for the parameters $\beta_1,\beta_2$:
$$\beta_1<f_1,f_1>+\beta_2<f_1,f_2>=-<f_1,f_3>$$
$$\beta_1<f_2,f_1>+\beta_2<f_2,f_2>=-<f_2,f_3>$$

Now let us compute the determinant of this system. This is given by:
\begin{eqnarray*}
D
&=&\begin{vmatrix}
<f_1,f_1>&<f_1,f_2>\\
<f_2,f_1>&<f_2,f_2>
\end{vmatrix}\\
&=&<f_1,f_1><f_2,f_2>-<f_1,f_2><f_2,f_1>\\
&=&||f_1||^2||f_2||^2-|<f_1,f_2>|^2
\end{eqnarray*}

But this is exactly the quantity from the Cauchy-Schwarz inequality, so we have $D\geq0$, with equality when $f_1,f_2$ are proportional. Now since $f_1,f_2$ were assumed to be linearly independent, we conclude that we have $D>0$, so our system has indeed solutions.

\medskip

(3) Alternatively, we can say at $n=3$ that with the vectors $e_1,e_2$ being already constructed, we can construct the vector $e_3$ as follows, obviously doing the orthogonality job, and with its uniqueness coming from some standard linear algebra:
$$e_3=f_3-Proj_{e_1}(f_3)-Proj_{e_2}(f_3)$$

(4) Summarizing, we have two possible proofs for our result. Getting now to the general case, as a first proof, which is perhaps the most straightforward, we can set:
$$e_1=f_1$$
$$e_2=f_2-Proj_{e_1}(f_2)$$
$$e_3=f_3-Proj_{e_1}(f_3)-Proj_{e_2}(f_3)$$
$$e_4=f_4-Proj_{e_1}(f_4)-Proj_{e_2}(f_4)-Proj_{e_3}(f_4)$$
$$\vdots$$

Indeed, these vectors do indeed the needed orthogonality job, and their uniqueness is clear too, via some basic linear algebra, that we will leave here as an exercise.

\medskip

(5) Alternatively, by doing some explicit computations, as in (1,2), at step $k+1$ of the orthogonalization algorithm, the system to be solved is as follows:
$$x_1<f_1,f_1>+x_2<f_1,f_2>+\ldots+x_k<f_1,f_k>=-<f_1,f_{k+1}>$$
$$x_1<f_2,f_1>+x_2<f_2,f_2>+\ldots+x_k<f_2,f_k>=-<f_2,f_{k+1}>$$
$$\vdots$$
$$x_1<f_k,f_1>+x_2<f_k,f_2>+\ldots+x_k<f_k,f_k>=-<f_k,f_{k+1}>$$

The determinant of this system is then the Gram determinant of our system of vectors $f_1,\ldots,f_k$, given by the following formula, and with exercise for you to show $D_k>0$:
$$D_k=\begin{vmatrix}
<f_1,f_1>&<f_1,f_2>&\ldots&<f_1,f_k>\\
<f_2,f_1>&<f_2,f_2>&\ldots&<f_2,f_k>\\
&\vdots&\vdots\\
<f_k,f_1>&<f_k,f_2>&\ldots&<f_k,f_k>
\end{vmatrix}$$

But with this, we can work out some explicit formulae for the vectors $e_k$, alternative to those in (4), which remain a bit theoretical. To be more precise, we have:
$$e_k=\frac{1}{D_{k-1}}\begin{vmatrix}
<f_1,f_1>&<f_1,f_2>&\ldots&<f_1,f_k>\\
<f_2,f_1>&<f_2,f_2>&\ldots&<f_2,f_k>\\
&\vdots&\vdots\\
<f_{k-1},f_1>&<f_{k-1},f_2>&\ldots&<f_{k-1},f_k>\\
f_1&f_2&\ldots&f_k
\end{vmatrix}$$

And I will leave some study here to you, as an instructive exercise, and please do better than my students, who usually stop after the first 2-3 steps of the algorithm.
\end{proof}

Getting back now to our Hilbert space questions, we have the following result:

\index{orthonormal basis}
\index{Gram-Schmidt}
\index{separable space}
\index{orthogonal basis}
\index{dimension}
\index{countable dimension}
\index{Zorn Lemma}

\begin{theorem}
Any Hilbert space $H$ has an orthonormal basis $\{e_i\}_{i\in I}$, which is by definition a set of vectors whose span is dense in $H$, and which satisfy
$$<e_i,e_j>=\delta_{ij}$$
with $\delta$ being a Kronecker symbol. The cardinality $|I|$ of the index set, which can be finite, countable, or uncountable, depends only on $H$, and is called dimension of $H$. We have
$$H\simeq l^2(I)$$
in the obvious way, mapping $\sum\lambda_ie_i\to(\lambda_i)$. The Hilbert spaces with $\dim H=|I|$ being countable, such as $l^2(\mathbb N)$, are all isomorphic, and are called separable.
\end{theorem}

\begin{proof}
We have several assertions here, the idea being as follows:

\medskip

(1) In finite dimensions an orthonormal basis $\{e_i\}_{i\in I}$ can be constructed by starting with any vector space basis $\{f_i\}_{i\in I}$, and using the Gram-Schmidt procedure from Theorem 13.5. As for the other assertions, these are all clear, from basic linear algebra.

\medskip

(2) In general, the same method works, namely Gram-Schmidt, with a subtlety coming from the fact that the basis $\{e_i\}_{i\in I}$ will not span in general the whole $H$, but just a dense subspace of it, as it is in fact obvious by looking at the standard basis of $l^2(\mathbb N)$. 

\medskip

(3) And there is a second subtlety as well, coming from the fact that the recurrence procedure needed for Gram-Schmidt must be replaced by some sort of ``transfinite recurrence'', using standard tools from logic, and more specifically the Zorn lemma.
\end{proof}

\section*{13b. Weierstrass, Fourier}

According to Theorem 13.6, there is only one separable Hilbert space, up to isomorphism. In order to further comment on this, we will need the following key result:

\begin{theorem}[Weierstrass]
Any continuous function on a closed interval
$$f:[a,b]\to\mathbb R$$
can be uniformly approximated by polynomials.
\end{theorem}

\begin{proof}
We can assume $[a,b]=[0,1]$. Let us define the family of Bernstein polynomials as follows, and with this reminding us the binomial laws from chapter 3:
$$b_{kn}(x)=\binom{n}{k}x^k(1-x)^{n-k}$$

Our precise claim is that we can get our approximation $f_n\to f$ by setting:
$$f_n(x)=\sum_{k=0}^nf\left(\frac{k}{n}\right)b_{kn}(x)$$

(1) In order to prove this, observe first we have the following formulae, which in probabilistic terms are dealing with the moments of order $0,1,2$ of the binomial laws:
$$\sum_kb_{kn}(x)=1\quad,\quad 
\sum_k\frac{k}{n}\cdot b_{kn}(x)=x\quad,\quad 
\sum_k\left(x-\frac{k}{n}\right)^2b_{kn}(x)=\frac{x(1-x)}{n}$$

(2) In order to estimate the error $|f_n-f|$, we can use the uniform continuity property of $f$. Indeed, pick $\varepsilon>0$, and then $\delta>0$ such that the following happens:
$$|x-y|<\delta\implies|f(x)-f(y)|<\varepsilon$$

Now with this done, we have the following estimate, using the first formula in (1) at the first step, the uniform continuity at the last step, and with $M=\sup|f|$:
\begin{eqnarray*}
&&|f_n(x)-f(x)|\\
&=&\left|\sum_k\left(f\left(\frac{k}{n}\right)-f(x)\right)b_{kn}(x)\right|\\
&\leq&\sum_k\left|f\left(\frac{k}{n}\right)-f(x)\right|b_{kn}(x)\\
&=&\sum_{\left|x-\frac{k}{n}<\delta\right|}\left|f\left(\frac{k}{n}\right)-f(x)\right|b_{kn}(x)
+\sum_{\left|x-\frac{k}{n}\geq\delta\right|}\left|f\left(\frac{k}{n}\right)-f(x)\right|b_{kn}(x)\\
&\leq&\varepsilon+M\sum_{\left|x-\frac{k}{n}\geq\delta\right|}b_{kn}(x)
\end{eqnarray*}

(3) The point now is that the last sum on the right can be estimated by using the Chebycheff inequality, based on the third formula from (1), and we obtain:
\begin{eqnarray*}
\sum_{\left|x-\frac{k}{n}\geq\delta\right|}b_{kn}(x)
&\leq&\sum_k\delta^{-2}\left(x-\frac{k}{n}\right)^2b_{kn}(x)\\
&=&\delta^{-2}\,\frac{x(1-x)}{n}\\
&\leq&\frac{\delta^{-2}}{4n}
\end{eqnarray*}

(4) Now by putting everything together, we obtain the following estimate:
$$|f_n(x)-f(x)|\leq \varepsilon+\frac{\delta^{-2}M}{4n}$$

Thus we have indeed $|f_n-f|\to0$, uniform convergence, as desired.
\end{proof}

Getting back now to our Hilbert space questions, we can formulate:

\index{orthogonal polynomials}
\index{Weierstrass basis}

\begin{theorem}
The following happen, in relation with separability:
\begin{enumerate}
\item The Hilbert space $H=L^2[-1,1]$ is separable, with orthonormal basis coming by applying Gram-Schmidt to the basis $\{x^k\}_{k\in\mathbb N}$, coming from Weierstrass.

\item In fact, any $H=L^2(\mathbb R,\mu)$, with $d\mu(x)=f(x)dx$, is separable, and the same happens in higher dimensions, for $H=L^2(\mathbb R^N,\mu)$, with $d\mu(x)=f(x)dx$.

\item More generally, given a separable abstract measured space $X$, the associated Hilbert space of square-summable functions $H=L^2(X)$ is separable.
\end{enumerate}
\end{theorem}
 
\begin{proof}
Many things can be said here, the idea being as follows:

\medskip

(1) The fact that $H=L^2[-1,1]$ is separable is clear indeed from the Weierstrass density theorem, which provides us with the algebraic basis $g_k=x^k$, which can be orthogonalized by using the Gram-Schmidt procedure, as explained in Theorem 13.6.

\medskip

(2) Regarding now more general spaces, of type $H=L^2(\mathbb R,\mu)$, we can use here the same argument, after modifying if needed our measure $\mu$, in order for the functions $g_k=x^k$ to be indeed square-summable. As for higher dimensions, the situation here is similar, because we can use here the multivariable polynomials $g_k(x)=x_1^{k_1}\ldots x_N^{k_N}$.

\medskip

(3) Finally, the last assertion, regarding the general spaces of type $H=L^2(X)$, which generalizes all this, comes as a consequence of the general measure theory that we learned in chapter 5, and I will leave working out the details here as an instructive exercise.
\end{proof}

As our next result now, on the same topics, namely finding orthonormal bases for separable Hilbert spaces, following Fourier, and unrelated to Weierstrass, we have:

\begin{theorem}
The space of $2\pi$-periodic square-summable functions on $\mathbb R$,
$$L^2(\mathbb R)_{per}=\left\{f:\mathbb R\to\mathbb C\,\Big|\,f(t)=f(t+2\pi),\,\int_{-\pi}^\pi|f(t)|^2dt<\infty\right\}$$
has $\{e^{int}\}_{n\in\mathbb Z}$ as orthonormal basis, with respect to the normalized mass $1$ measure.
\end{theorem}

\begin{proof}
This is something quite tricky, which came as a big surprise at the time of its discovery by Fourier, the idea with all this being as follows:

\medskip

(1) For reasons that will become clear in a moment, consider the following trigonometric polynomials, with $c_k\in\mathbb R$ being chosen for having mass $1$:
$$E_k(t)=c_k\left(\frac{1+\cos t}{2}\right)^k$$

Observe that $c_k$ can be computed explicitly, by using $1+\cos t=2\cos^2(t/2)$ and the Wallis formulae from chapter 12, but in what follows, we will not need this. Our claim, which is the reason for introducing these functions $E_k$, is that for any $\delta>0$ we have:
$$\lim_{k\to\infty}\ \sup_{\delta<|t|<\pi}E_k(t)=0$$

(2) So, let us prove this claim. As mentioned, $c_k$ can be computed explicitly, but in what follows we will only need the following elementary estimate:
\begin{eqnarray*}
1
&=&\frac{c_k}{\pi}\int_0^\pi\left(\frac{1+\cos t}{2}\right)^kdt\\
&>&\frac{c_k}{\pi}\int_0^\pi\left(\frac{1+\cos t}{2}\right)^k\sin tdt\\
&=&\frac{c_k}{\pi}\left[-\frac{2}{k+1}\left(\frac{1+\cos t}{2}\right)^{k+1}\right]_0^\pi\\
&=&\frac{2c_k}{\pi(k+1)}
\end{eqnarray*}

Now since $E_k$ is decreasing on $[0,\pi]$, we obtain from this, for $\delta<|t|<\pi$:
$$E_k(t)<E_k(\delta)<\frac{\pi(k+1)}{2}\left(\frac{1+\cos\delta}{2}\right)^k$$

But this proves our claim, because for $\delta>0$ we have $(1+\cos\delta)/2<1$, as needed.

\medskip

(3) Getting now to what we wanted to do, we must prove that $\{e^{int}\}_{n\in\mathbb Z}$ spans a dense subset of $L^2(\mathbb R)_{per}$. Since $C(\mathbb R)_{per}\subset L^2(\mathbb R)_{per}$ is dense, it is enough to prove that any $f\in C(\mathbb R)_{per}$ can be approximated by trigonometric polynomials $\sum_nc_ne^{int}$. Moreover, since $||.||_2\leq||.||_\infty$, it is enough to prove our approximation with respect to $||.||_\infty$.

\medskip

(4) All in all, it remains to prove that given a function $f\in C(\mathbb R)_{per}$ and a number $\varepsilon>0$, we can always come with a trigonometric polynomial $\sum_nc_ne^{int}$, such that:
$$\left|f(t)-\sum_nc_ne^{int}\right|<\varepsilon\quad,\quad\forall t\in[-\pi,\pi]$$

(5) But for this, we can use the polynomials $E_k$ from (1). Let us set indeed:
$$Q_k(t)=\frac{1}{2\pi}\int_{-\pi}^\pi f(t-s)E_k(s)\,ds$$

As a first observation, with the change of variables $s\to t-s$, we have the following alternative formula, which shows that $Q_k(t)$ are indeed trigonometric polynomials:
$$Q_k(t)=\frac{1}{2\pi}\int_{-\pi}^\pi f(s)E_k(t-s)\,ds$$

(6) Now given $\varepsilon>0$, let us prove that the estimate in (4) holds indeed, with the trigonometric polynomial there being $Q_k(t)$, for $k>>0$ large enough. For this purpose, we use the uniform continuity of $f$, which tells us that we can find $\delta>0$ such that:
$$|s-t|<\delta\implies|f(s)-f(t)|<\varepsilon$$

Indeed, by using this, we have the following estimate, for the error in (4):
\begin{eqnarray*}
|Q_k(t)-f(t)|
&=&\frac{1}{2\pi}\left|\int_{-\pi}^\pi(f(t-s)-f(t))E_k(s)\,ds\right|\\
&\leq&\frac{1}{2\pi}\int_{-\pi}^\pi|f(t-s)-f(t)|E_k(s)\,ds
\end{eqnarray*}

(7) Now let us split the last integral into three parts, according to:
$$[-\pi,\pi]=[-\pi,-\delta]\cup[-\delta,\delta]\cup[\delta,\pi]$$

On the middle part the integrand is $<\varepsilon$, so the middle integral is $<\varepsilon$. As for the other two integrals, on $[-\pi,-\delta]\cup[\delta,\pi]$, we can use here (6), telling us that $E_k(t)\to0$ uniformly, on that domain. Indeed, with $k>>0$ big enough the other two integrals are $<\varepsilon$ too, so we have obtained (4) as desired, with $\varepsilon\to3\varepsilon$, which finishes the proof. 
\end{proof}

Still with me I hope, after all these computations. As a more digest formulation of Theorem 13.9, in relation with our abstract considerations in this chapter, we have:

\begin{theorem}
We have an isomorphism $L^2(\mathbb R)_{per}\simeq l^2(\mathbb Z)$, as follows:
\begin{enumerate}
\item Associated to $f\in L^2(\mathbb R)_{per}$ are its Fourier coefficients, given by:
$$\widehat{f}(n)=\frac{1}{2\pi}\int_{-\pi}^\pi f(t)e^{int}dt$$

\item Associated to $g\in l^2(\mathbb Z)$ is the series $S_g(t)=\sum_{n\in\mathbb Z}g(n)e^{-int}$.
\end{enumerate}
\end{theorem}

\begin{proof}
This is something self-explanatory, based on the general orthonormal basis theory from Theorem 13.6, with the Fourier coefficients of $f\in L^2(\mathbb R)_{per}$ being its coefficients $\widehat{f}(n)=<f,e^{-int}>$ with respect to the basis $\{e^{int}\}$ from Theorem 13.9.
\end{proof}

Finally, no discussion about Fourier series would be complete without:

\begin{theorem}
We have the following formula of Euler,
$$1+\frac{1}{4}+\frac{1}{9}+\frac{1}{16}+\frac{1}{25}+\ldots=\frac{\pi^2}{6}$$
computing the value $\zeta(2)$ of the Riemann zeta function, and solving the Basel problem.
\end{theorem}

\begin{proof}
The nonzero Fourier coefficients of the function $f(t)=t$ on $[-\pi,\pi]$ are:
$$\widehat{f}(n)
=\frac{1}{2\pi}\int_{-\pi}^\pi te^{int}\,dt
=\frac{1}{2\pi}\left[\frac{1-int}{n^2}\,e^{int}\right]_{-\pi}^\pi
=\frac{(-1)^{n+1}}{n}\,i$$

Now by using the isomorphism in Theorem 13.10, at the level of norms, we obtain:
$$\sum_{n\in\mathbb Z^*}\frac{1}{n^2}=\frac{1}{2\pi}\int_{-\pi}^\pi t^2dt=\frac{1}{2\pi}\cdot\frac{2\pi^3}{3}=\frac{\pi^2}{3}$$

And we therefore solved the Basel problem, just like that. Amazing.
\end{proof}

\section*{13c. Orthogonal polynomials}

Let us go back now to Theorem 13.8 and its proof, which was something quite subtle. That material leads us into orthogonal polynomials, which are defined as follows:

\index{orthogonal polynomials}
\index{real measure}

\begin{definition}
The orthogonal polynomials with respect to $d\mu(x)=f(x)dx$ are polynomials $P_k\in\mathbb R[x]$ of degree $k\in\mathbb N$, which are orthogonal inside $H=L^2(\mathbb R,\mu)$:
$$\int_\mathbb RP_k(x)P_l(x)f(x)dx=0\quad,\quad\forall k\neq l$$
Equivalently, these orthogonal polynomials $\{P_k\}_{k\in\mathbb N}$, which are each unique modulo scalars, appear from the Weierstrass basis $\{x^k\}_{k\in\mathbb N}$, by doing Gram-Schmidt.
\end{definition}

Observe that the orthogonal polynomials exist indeed for any real measure $d\mu(x)=f(x)dx$, as explained above. It is possible to be a bit more explicit here, as follows:

\index{moments of measure}

\begin{theorem}
The orthogonal polynomials with respect to $\mu$ are given by
$$P_k=c_k\begin{vmatrix}
M_0&M_1&\ldots&M_k\\
M_1&M_2&\ldots&M_{k+1}\\
\vdots&\vdots&&\vdots\\
M_{k-1}&M_k&\ldots&M_{2k-1}\\
1&x&\ldots&x^k
\end{vmatrix}$$
where $M_k=\int_\mathbb Rx^kd\mu(x)$ are the moments of $\mu$, and $c_k\in\mathbb R^*$ can be any numbers.
\end{theorem}

\begin{proof}
Let us first see what happens at small values of $k\in\mathbb N$. At $k=0$ our formula is as follows, stating that the first polynomial $P_0$ must be a constant, as it should:
$$P_0=c_0|M_0|=c_0$$

At $k=1$ now, again by using $M_0=1$, the formula is as follows:
$$P_1=c_1\begin{vmatrix}M_0&M_1\\ 1&x\end{vmatrix}=c_1(x-M_1)$$

But this is again the good formula, because the degree is 1, and we have:
\begin{eqnarray*}
<1,P_1>
&=&c_1<1,x-M_1>\\
&=&c_1(<1,x>-<1,M_1>)\\
&=&c_1(M_1-M_1)\\
&=&0
\end{eqnarray*}

At $k=2$ now, things get more complicated, with the formula being as follows:
$$P_2=c_2\begin{vmatrix}
M_0&M_1&M_2\\
M_1&M_2&M_3\\
1&x&x^2
\end{vmatrix}$$

However, no need for big computations here, in order to check the orthogonality, because by using the fact that $x^k$ integrates up to $M_k$, we obtain:
$$<1,P_2>=\int_\mathbb RP_2(x)d\mu(x)=c_2\begin{vmatrix}
M_0&M_1&M_2\\
M_1&M_2&M_3\\
M_0&M_1&M_2
\end{vmatrix}=0$$

Similarly, again by using the fact that $x^k$ integrates up to $M_k$, we have as well:
$$<x,P_2>=\int_\mathbb RxP_2(x)d\mu(x)=c_2\begin{vmatrix}
M_0&M_1&M_2\\
M_1&M_2&M_3\\
M_1&M_2&M_3
\end{vmatrix}=0$$

Thus, result proved at $k=0,1,2$, and the proof in general is similar.
\end{proof}

In practice now, all this leads us to a lot of interesting combinatorics, and countless things can be said. For the simplest measured space $X\subset\mathbb R$, which is the interval $[-1,1]$, with its uniform measure, the orthogonal basis problem can be solved as follows:

\index{Legendre polynomials}
\index{Rodrigues formula}
\index{Legendre equation}

\begin{theorem}
The orthonormal polynomials for $L^2[-1,1]$, subject to
$$\int_{-1}^1P_k(x)P_l(x)\,dx=\delta_{kl}$$
and called Legendre polynomials, satisfy the following differential equation,
$$(1-x^2)P_k''(x)-2xP_k'(x)+k(k+1)P_k(x)=0$$
which is the Legendre equation from physics. Moreover, we have the formula
$$(k+1)P_{k+1}(x)=(2k+1)xP_k(x)-kP_{k-1}(x)$$
called Bonnet recurrence formula, as well as the formula
$$P_k(x)=\frac{1}{2^kk!}\cdot\frac{d^k}{dx^k}\left(1-x^2\right)^k$$
called Rodrigues formula for the Legendre polynomials.
\end{theorem}

\begin{proof}
As a first observation, we are not lost somewhere in abstract math, because of the occurrence of the Legendre equation. As for the proof, this goes as follows:

\medskip

(1) The first assertion is clear, because the Gram-Schmidt procedure applied to the Weierstrass basis $\{x^k\}$ can only lead to a certain family of polynomials $\{P_k\}$, with each $P_k$ being of degree $k$, and also unique, if we assume that it has positive leading coefficient, with this $\pm$ choice being needed, as usual, at each step of Gram-Schmidt.

\medskip

(2) In order to have now an idea about these beasts, here are the first few of them, which can be obtained say via a straightforward application of Gram-Schmidt:
\begin{eqnarray*}
P_0&=&1\\
P_1&=&x\\
P_2&=&(3x^2-1)/2\\
P_3&=&(5x^3-3x)/2\\
P_4&=&(35x^4-30x^2+3)/8\\
P_5&=&(63x^5-70x^3+15x)/8
\end{eqnarray*}

(3) Now thinking about what Gram-Schmidt does, this is certainly something by recurrence. And examining the recurrence leads to the Legendre equation, as stated. As for the Bonnet recurrence formula, the story here is similar.

\medskip

(4) Regarding the Rodrigues formula, by uniqueness no need to try to understand where this formula comes from, and we have two choices here, either by verifying that $\{P_k\}$ is orthonormal, or by verifying the Legendre equation. And both methods work.
\end{proof}

The above result is just the tip of the iceberg, and as a continuation, we have:

\index{Jacobi polynomials}
\index{Chebycheff polynomials}

\begin{theorem}
The orthogonal polynomials for $L^2[-1,1]$, with measure
$$d\mu(x)=(1-x)^a(1+x)^bdx$$
called Jacobi polynomials, satisfy as well a degree $2$ equation, namely
$$(1-x^2)P_k''(x)+(b-a-(a+b+2)x)P_k'(x)+k(k+a+b+1)P_k(x)=0$$
as well as an order $2$ recurrence relation, and are given by the following formula:
$$P_k(x)=\frac{(-1)^k}{2^kk!}(1-x)^{-a}(1+x)^{-b}\frac{d^k}{dx^k}\left[(1-x)^a(1+x)^b(1-x^2)^k\right]$$
At $a=b=0$ we recover the Legendre polynomials, and at $a=b=\pm\frac{1}{2}$ we recover the Chebycheff polynomials of the first and second kind, from trigonometry. 
\end{theorem}

\begin{proof}
There are many things going on here, the idea being as follows:

\medskip

(1) To start with, in what regards the precise statement, the order 2 recurrence relation mentioned there is something quite complicated, as follows:
\begin{eqnarray*}
&&2k(k+a+b)(2k+a+b-2)P_k(x)\\
&=&(2k+a+b-1)\left[(2k+a+b)(2k+a+b-2)x+a^2-b^2\right]P_{k-1}(x)\\
&-&2(k+a-1)(k+b-1)(2k+a+b)P_{k-2}(x)
\end{eqnarray*}

(2) Regarding now the proof, the statement itself appears as a generalization of Theorem 13.14, which corresponds to the particular case $a=b=0$, and the proof is quite similar. We will leave learning more about all this as an interesting exercise.

\medskip

(3) For completness, let us record as well a few numerics, as follows:
\begin{eqnarray*}
P_0&=&1\\
P_1&=&(a+1)+(a+b+2)\,\frac{x-1}{2}\\
P_2&=&\frac{(a+1)(a+2)}{2}+(a+2)(a+b+3)\,\frac{x-1}{2}\\
&+&\frac{(a+b+3)(a+b+4)}{2}\left(\frac{x-1}{2}\right)^2
\end{eqnarray*}

(4) Regarding now the main particular cases of the Jacobi polynomials, these are the Gegenbauer polynomials, appearing at $a=b$. However, there is not that much of a simplification when passing from general parameters $a,b$ to equal parameters, $a=b$, so in practice, the main particular cases are those indicated in the statement, namely:

\medskip

-- The Legendre polynomials, that we know well from Theorem 13.14, appearing at the simplest values of the parameters, namely $a=b=0$.

\medskip

-- The Chebycheff polynomials of the first kind $T_k$, which are given by the formula $T_k(\cos t)=\cos(kt)$ from trigonometry, appearing at $a=b=-\frac{1}{2}$.

\medskip

-- The Chebycheff polynomials of the second kind $U_k$, which are given by the formula $U_k(\cos t)\sin t=\sin((k+1)t)$, appearing at $a=b=\frac{1}{2}$.

\medskip

(5) So, this was for the story of the Jacobi polynomials, and their main particular cases, and in practice, we will leave some further learning here as an exercise, coming as a continuation of the further learning of Theorem 13.14, and its details.
\end{proof}

Getting now to other spaces $X\subset\mathbb R$, of particular interest here is the following result, which complements well Theorem 13.14, for the needs of basic quantum mechanics:

\index{Laguerre polynomials}

\begin{theorem}
The orthogonal polynomials for $L^2[0,\infty)$, with scalar product
$$<f,g>=\int_0^\infty f(x)g(x)e^{-x}\,dx$$
are the Laguerre polynomials $\{P_k\}$, satisfying the following differential equation,
$$xP_k''(x)+(1-x)P_k'(x)+kP_k(x)=0$$
as well as the following order $2$ recurrence relation,
$$(k+1)P_{k+1}(x)=(2k+1-x)P_k(x)-kP_{k-1}(x)$$
and which are given by the following formula,
$$P_k(x)=\frac{e^x}{k!}\cdot\frac{d^k}{dx^k}\left(e^{-x}x^k\right)$$
called Rodrigues formula for the Laguerre polynomials.
\end{theorem}

\begin{proof}
The story here is very similar to that of the Legendre and Jacobi polynomials, and many further things can be said here, with exercise for you to learn a bit about all this. Let us record as well a few numeric values, for the Laguerre polynomials:
\begin{eqnarray*}
P_0&=&1\\
P_1&=&1-x\\
P_2&=&(x^2-4x+2)/2\\
P_3&=&(-x^3+9x^2-18x+6)/6\\
P_4&=&(x^4-16x^3+72x^2-96x+24)/24
\end{eqnarray*}

Finally, for the story to be complete, no discussion about the Laguerre polynomials would be complete without a word about their use, in quantum mechanics. And here, as usual, we will leave some exploration of this as an instructive exercise.
\end{proof}

Finally, regarding the space $X=\mathbb R$ itself, we have here the following result:

\index{Hermite polynomials}

\begin{theorem}
The orthogonal polynomials for $L^2(\mathbb R)$, with scalar product
$$<f,g>=\int_0^\infty f(x)g(x)e^{-x^2}\,dx$$
are the Hermite polynomials $\{P_k\}$, satisfying the following differential equation,
$$P_k''(x)-2xP_k'(x)+P_k(x)=0$$
as well as the following order $2$ recurrence relation,
$$P_{k+1}(x)=2xP_k(x)-2kP_{k-1}(x)$$
and which are given by the following formula,
$$P_k(x)=(-1)^k e^{x^2}\cdot\frac{d^k}{dx^k}\big(e^{-x^2}\big)$$
called Rodrigues formula for the Hermite polynomials.
\end{theorem}

\begin{proof}
As before, the story here is quite similar to that of the Legendre and other orthogonal polynomials, and exercise for you to learn a bit about all this. Let us record as well a few numeric values, for the Hermite polynomials:
\begin{eqnarray*}
P_0&=&1\\
P_1&=&2x\\
P_2&=&4x^2-2\\
P_3&=&8x^3-12x\\
P_4&=&16x^4-48x^2+12\\
P_5&=&32x^5-160x^3+120x
\end{eqnarray*}

With of course, exercise for you to deduce all these formulae.
\end{proof}

And with this, good news, end of the story with the orthogonal polynomials, at least at the very introductory level, and this due to the following fact, which is something quite technical, and that we will not attempt to prove, or even explain in detail here:

\begin{fact}
From an abstract point of view, coming from degree $2$ equations, and Rodrigues formulae for the solutions, there are only three types of ``classical'' orthogonal polynomials, namely the Jacobi, Laguerre and Hermite ones, discussed above.
\end{fact}

Finally, as already mentioned, the above results are very useful in the context of basic quantum mechanics, and more specifically, for solving the hydrogen atom, following Heisenberg and Schr\"odinger. Have a look here at Griffiths \cite{gri}, or at my book \cite{ba3}.

\section*{13d. Circular measures}

We would like to end this chapter with the discussion of a question that we were stuck with, at the beginning of chapter 7, namely random walks on basic finite graphs. Let us start with a general comment, making the link between things done in the above, which were many, and the material from chapter 7, which was quite varied too:

\begin{comment}
The following happen, with the laws being up to scalars:
\begin{enumerate}
\item The simplest orthogonal polynomials are the Chebycheff ones, $T_k$ and $U_k$.

\item These appear via $T_k(\cos t)=\cos(kt)$ and $U_k(\cos t)\sin t=\sin((k+1)t)$.

\item They satisfy $P_{k+1}=2xP_k-P_{k-1}$, with $T_0=1,T_1=x$ and $U_0=1,U_1=2x$.

\item $T_k$ are the orthogonal polynomials for $1/\sqrt{1-x^2}$, arcsine law on $[-1,1]$.

\item $U_k$ are the orthogonal polynomials for $\sqrt{1-x^2}$, semicircle law on $[-1,1]$.
\end{enumerate}
\end{comment}

To be more precise, (1) is something that we learned from our study from the previous section, with the Chebycheff polynomials certainly getting the simplicity crown. Next, (2) is the definition of these polynomials, and with this specifying them, as orthogonal polynomials, up to scalars. Then (3) comes from some trigonometry, or from the mathematics in Theorem 13.15, at $a=b=\pm\frac{1}{2}$. And finally (4) and (5) come from Theorem 13.15 at $a=b=-\frac{1}{2},\frac{1}{2}$, with the masses of the measures there being $\pi,\pi/2$.

\bigskip

In relation now with our graph considerations from chapter 7, we have:

\begin{fact}
The characteristic polynomial of the segment graph with $N$ vertices is
$$P_N(x)=U_N\left(\frac{x}{2}\right)$$
and this allows in principle to solve the random walk question, for this graph.
\end{fact}

To be more precise, this is something that we came upon in chapter 7, based on our recurrence formula found there, namely $P_{N+1}=xP_N-P_{N-1}$, with $P_0=1,P_1=x$.

\bigskip

However, it is not about this, which remains something quite complicated, that we would like to talk about. The point is that we have something better, as follows:

\begin{principle}
The random walk question for basic finite graphs is best solved by mapping the measure $\mu$ to the unit circle, via some sort of Fourier transform. 
\end{principle}

So, let us get into this. We will be actually interested in a family of graphs much larger than the segment ones, namely the ADE family, whose basics are as follows:

\index{algebraic manifolds}
\index{quantum field theory}
\index{subfactors}
\index{knots and links}
\index{planar algebras}

\begin{fact}
The ADE graphs classify the following:
\begin{enumerate}
\item Basic Lie groups and algebras.

\item Subgroups of $SU_2$ and of $SO_3$.

\item Singularities of algebraic manifolds.

\item Basic invariants of knots and links.

\item Subfactors and planar algebras of small index.

\item Subgroups of the quantum permutation group $S_4^+$.

\item Basic quantum field theories, and other physics beasts.
\end{enumerate}
\end{fact}

Which sounds quite exciting, good mathematics that we will be learning here. Getting to work now, we first need to know what the ADE graphs exactly are. The A graphs, which are the simplest, are as follows, with the distinguished vertex being denoted $\bullet$, and with $A_n$ having by definition $n\geq2$ vertices, and $\tilde{A}_{2n}$ having $2n\geq2$ vertices:
$$A_n=\bullet-\circ-\circ\cdots\circ-\circ-\circ\hskip18mm 
A_{\infty}=\bullet-\circ-\circ-\circ\cdots\hskip7mm$$
\vskip-3mm
$$\ \ \ \ \ \ \ \tilde{A}_{2n}=
\begin{matrix}
\circ&\!\!\!\!-\circ-\circ\cdots\circ-\circ-&\!\!\!\!\circ\\
|&&\!\!\!\!|\\
\bullet&\!\!\!\!-\circ-\circ-\circ-\circ-&\!\!\!\!\circ\\
\\
\\
\end{matrix}\hskip20mm 
\tilde{A}_\infty=
\begin{matrix}
\circ&\!\!\!\!-\circ-\circ-\circ\cdots\\
|&\\
\bullet&\!\!\!\!-\circ-\circ-\circ\cdots\\
\\
\\
\end{matrix}
\hskip15mm$$
\vskip-7mm

You might probably say, why not stopping here, and doing our unfinished business for the segment and the circle, with whatever new ideas that we might have. Good point, but in answer, these ideas will apply as well, with minimal changes, to the D graphs, which are as follows, with $D_n$ having $n\geq3$ vertices, and $\tilde{D}_n$ having $n+1\geq5$ vertices:
$$D_n=\bullet-\circ-\circ\dots\circ-
\begin{matrix}\ \circ\\
\ |\\
\ \circ \\
\ \\
\  \end{matrix}-\circ\hskip71mm$$
\vskip-7mm
$$\hskip7mm\tilde{D}_n=\bullet-
\begin{matrix}\circ\\
|\\
\circ\\
\ \\
\ \end{matrix}-\circ\dots\circ-
\begin{matrix}\ \circ\\
\ |\\
\ \circ \\
\ \\
\  \end{matrix}-\circ\hskip18mm$$
\vskip-7mm
$$\hskip50mm D_\infty=\bullet-
\begin{matrix}\circ\\
|\\
\circ\\
\ \\
\ \end{matrix}-\circ-\circ\cdots$$
\vskip-7mm

So, these are the famous AD graphs, appearing as explained in Fact 13.22.

\bigskip

As a comment now, the labeling conventions for the AD graphs, while very standard, can be a bit confusing. The first graph in each series is by definition as follows:
$$A_2=\bullet-\circ\hskip13mm 
\tilde{A}_2=\begin{matrix}
\circ\\
||\\
\bullet\\
&\\
&\\
\end{matrix}\hskip13mm 
D_3=\begin{matrix}\ \circ\\
\ |\\
\ \bullet \\
\ \\
\  \end{matrix}-\circ \hskip13mm
\tilde{D}_4=\bullet-\!\!\!\!\!\begin{matrix}
\circ\hskip5mm \circ\\
\backslash\ \,\slash\\
\circ\\
&\\
&\\
\end{matrix}\!\!\!\!\!\!\!\!\!\!-\circ$$
\vskip-7mm

Finally, there are also a number of exceptional ADE graphs. First we have:
$$E_6=\bullet-\circ-
\begin{matrix}\circ\\
|\\
\circ\\
\ \\
\ \end{matrix}-
\circ-\circ\hskip71mm$$
\vskip-13mm
$$E_7=\bullet-\circ-\circ-
\begin{matrix}\circ\\
|\\
\circ\\
\ \
\\
\ \end{matrix}-
\circ-\circ\hskip18mm$$
\vskip-15mm
$$\hskip30mm E_8=\bullet-\circ-\circ-\circ-
\begin{matrix}\circ\\
|\\
\circ\\
\ \\
\ \end{matrix}-
\circ-\circ$$
\vskip-5mm

Then, we have extended versions of the above exceptional graphs, as follows:
$$\tilde{E}_6=\bullet-\circ-\begin{matrix}
\circ\\
|
\\
\circ\\
|&\\
\circ&\!\!\!\!-\ \circ\\
\ \\
\   \\
\ \\
\ \end{matrix}-\circ\hskip71mm$$
\vskip-22mm
$$\tilde{E}_7=\bullet-\circ-\circ-
\begin{matrix}\circ\\
|\\
\circ\\
\ \\
\ \end{matrix}-
\circ-\circ-\circ\hskip18mm$$
\vskip-15mm
$$\hskip30mm \tilde{E}_8=\bullet-\circ-\circ-\circ-\circ-
\begin{matrix}\circ\\
|\\
\circ\\
\ \\
\ \end{matrix}-
\circ-\circ$$
\vskip-5mm

And good news, according to the general ADE theory, that is all. Getting now to work, we have plenty of graphs here, and the problem that we would like to solve is:

\begin{problem}
How to count loops on the ADE graphs?
\end{problem}

In answer now, as a first observation, we know that $A_\infty$ and $\tilde{A}_\infty$ are respectively the graphs that we previously called $\mathbb N$ and $\mathbb Z$, and studied in great detail, in chapter 7. So, based on our previous work for these graphs, where the combinatorics naturally led us into generating series, let us formulate the following definition:

\index{Poincar\'e series}
\index{bipartite graph}

\begin{definition}
The Poincar\'e series of a rooted bipartite graph $X$ is
$$f(z)=\sum_{k=0}^\infty L_{2k}z^k$$
where $L_{2k}$ is the number of $2k$-loops based at the root.
\end{definition}

To be more precise, observe that all the above ADE graphs are indeed bipartite. Now the point is that, for a bipartite graph, the loops based at any point must have even length. Thus, in order to study the loops on the ADE graphs, we are led to $f(z)$.

\bigskip

Before getting into computations, let us introduce as well:

\index{positive spectral measure}

\begin{definition}
The positive spectral measure $\nu$ of a rooted bipartite graph $X$ is the real probability measure having the numbers $L_{2k}$ as moments:
$$\int_\mathbb Rx^kd\nu(x)=L_{2k}$$
Equivalently, we must have the Stieltjes transform formula
$$f(z)=\int_\mathbb R\frac{1}{1-xz}\,d\nu(x)$$
where $f$ is the Poincar\'e series of $X$.
\end{definition}

Here the existence of $\nu$, and the fact that this is indeed a positive measure, meaning a measure supported on $[0,\infty)$, comes from the following simple fact:

\begin{theorem}
The positive spectral measure of a rooted bipartite graph $X$ is given by the following formula, with $d$ being the adjacency matrix of the graph,
$$\nu=law(d^2)$$
and with the probabilistic computation being with respect to the expectation 
$$A\to<A>$$
with $<A>$ being the $(*,*)$-entry of a matrix $A$, where $*$ is the root.
\end{theorem}

\begin{proof}
With the above conventions, we have the following computation:
\begin{eqnarray*}
f(z)
&=&\sum_{k=0}^\infty L_{2k}z^k\\
&=&\sum_{k=0}^\infty\left<d^{2k}\right>z^k\\
&=&\left<\frac{1}{1-d^2z}\right>
\end{eqnarray*}

But this shows that we have $\nu=law(d^2)$, as desired.
\end{proof}

The above result shows that computing $\nu$ might be actually a simpler problem than computing $f$, and in practice, this is indeed the case. So, in what follows we will rather forget about loops and Definition 13.24, and use Definition 13.25 instead, with our computations to follow being based on the concrete interpretation from Theorem 13.26.

\bigskip

However, even with this probabilistic trick in our bag, things are not exactly trivial. Let us introduce as well the following notion, which is something more subtle, coming from the work of Jones on subfactor theory, and related topics \cite{jo1}, \cite{jo2}, \cite{jo3}, \cite{jo4}:

\index{circular measure}
\index{circular spectral measure}

\begin{definition}
The circular measure $\varepsilon$ of a rooted bipartite graph $X$ is given by
$$d\varepsilon(q)=d\nu((q+q^{-1})^2)$$
where $\nu$ is the associated positive spectral measure.
\end{definition}

To be more precise, we know from Theorem 13.26 that the positive measure $\nu$ is the spectral measure of a certain positive matrix, $d^2\geq0$, and it follows from this, and from basic spectral theory, that this measure is supported by the positive reals:
$$supp(\nu)\subset\mathbb R_+$$

But then, with this observation in hand, we can define indeed the circular measure $\varepsilon$ as above, as being the pullback of $\nu$ via the following map:
$$\mathbb R\cup\mathbb T\to\mathbb R_+\quad,\quad 
q\to (q+q^{-1})^2$$

As a basic example for this, to start with, assume that $\nu$ is a discrete measure, supported by $n$ positive numbers $x_1<\ldots<x_n$, with corresponding densities $p_1,\ldots,p_n$:
$$\nu=\sum_{i=1}^n p_i\delta_{x_i}$$

For each $i\in\{1,\ldots,n\}$ the equation $(q+q^{-1})^2=x_i$ has then four solutions, that we can denote $q_i,q_i^{-1},-q_i,-q_i^{-1}$. And with this notation, we have:
$$\varepsilon=\frac{1}{4}\sum_{i=1}^np_i\left(\delta_{q_i}+\delta_{q_i^{-1}}+\delta_{-q_i}+\delta_{-q_i^{-1}}\right)$$

In general, the basic properties of $\varepsilon$ can be summarized as follows:

\begin{theorem}
The circular measure has the following properties:
\begin{enumerate}
\item $\varepsilon$ has equal density at $q,q^{-1},-q,-q^{-1}$.

\item The odd moments of $\varepsilon$ are $0$.

\item The even moments of $\varepsilon$ are half-integers.

\item When $X$ has norm $\leq 2$, $\varepsilon$ is supported by the unit circle.

\item When $X$ is finite, $\varepsilon$ is discrete.

\item If $K$ is a solution of $d=K+K^{-1}$, then $\varepsilon=law(K)$. 
\end{enumerate}
\end{theorem}

\begin{proof}
These results can be deduced from definitions, the idea being that (1-5) are trivial, and that (6) follows from the formula of $\nu$ from Theorem 13.26.
\end{proof}

Getting now to computations, we first have the following result:

\begin{theorem}
The circular measure of the basic index $4$ graph, namely 
$$\begin{matrix}
&\circ&\!\!\!\!-\circ-\circ\cdots\circ-\circ-&\!\!\!\!\circ\cr
\tilde{A}_{2n}=&|&&\!\!\!\!|\cr
&\bullet&\!\!\!\!-\circ-\circ-\circ-\circ-&\!\!\!\!\circ\cr\cr\cr\end{matrix}$$
\vskip-7mm

\noindent is the uniform measure on the $2n$-roots of unity.
\end{theorem}

\begin{proof}
Let us identify the vertices of $X=\tilde{A}_{2n}$ with the group $\{w^k\}$ formed by the $2n$-th roots of unity in the complex plane, where $w=e^{\pi i/n}$. The adjacency matrix of $X$ acts then on the functions $f\in C(X)$ in the following way:
$$df(w^s)=f(w^{s-1})+f(w^{s+1})$$

But this shows that we have $d=K+K^{-1}$, where $K$ is given by:
$$Kf(w^s)=f(w^{s+1})$$

Thus we can use Theorem 13.26 and Theorem 13.28 (6), and we get:
$$\varepsilon=law(K)$$

But this is the uniform measure on the $2n$-roots of unity, as claimed.
\end{proof}

Observe that Theorem 13.29 nukes our previous work on the circle graph. In order to deal now with the other ADE graphs, let us introduce the following densities:
$$\alpha=Re(1-q^2)$$
$$\beta=Re(1-q^4)$$
$$\gamma=Re(1-q^6)$$

We have then the following result, $d_n$ being the uniform measure on the $2n$-th roots of unity, and $d_n'=2d_{2n}-d_n$ being the uniform measure on the odd $4n$-roots of unity:

\index{ADE graph}
\index{circular measure}
\index{cyclotomic measure}

\begin{theorem}
The circular measures of ADE graphs are given by:
\begin{enumerate}
\item $A_{n-1}\to\alpha_n$.

\item $\tilde{A}_{2n}\to d_n$.

\item $D_{n+1}\to\alpha_n'$.

\item $\tilde{D}_{n+2}\to (d_n+d_1')/2$.

\item $E_6\to\alpha_{12}+(d_{12}-d_6-d_4+d_3)/2$.

\item $E_7\to\beta_9'+(d_1'-d_3')/2$.

\item $E_8\to\alpha_{15}'+\gamma_{15}'-(d_5'+d_3')/2$.

\item $\tilde{E}_{n+3}\to (d_n+d_3+d_2-d_1)/2$.
\end{enumerate}
\end{theorem}

\begin{proof}
This is something which can be proved in three steps, as follows:

\medskip

(1) For the simplest graph, namely the circle $\tilde{A}_{2n}$, we already have the result, from Theorem 13.29, with the proof there being something elementary.

\medskip

(2) For the other non-exceptional graphs, that is, of type A and D, the same method works, namely direct loop counting, with some matrix tricks.

\medskip

(3) As for the exceptional graphs, of type E, basically the same method works, but the computations are more complicated, using some root of unity know-how. 

\medskip

So, this was for the idea, and in practice, read Jones \cite{jo1}, \cite{jo2}, \cite{jo3}, \cite{jo4}, and then have a look at my paper with Bisch \cite{bbi}, where all this is explained.
\end{proof}

And we will end the present chapter with this. As a conclusion to what we did, infinite dimensionality, coming from the key notion of Hilbert space, provides a whole new point of view on probability theory. We will keep investigating this, in what follows.

\section*{13e. Exercises}

This was a a quite advanced chapter, and as exercises on this, we have:

\begin{exercise}
Is there a complex analogue of $<x,y>=||x||\cdot||y||\cdot\cos t$?
\end{exercise}

\begin{exercise}
Do some Gram-Schmidt computations, the more the better.
\end{exercise}

\begin{exercise}
Experiment a bit with Weierstrass approximation.
\end{exercise}

\begin{exercise}
Learn about the various types of Fourier transforms available.
\end{exercise}

\begin{exercise}
Clarify all the basics, in relation with orthogonal polynomials.
\end{exercise}

\begin{exercise}
Learn the Legendre, Jacobi, Laguerre and Hermite polynomials.
\end{exercise}

\begin{exercise}
Learn about Chebycheff polynomials, as much as you can.
\end{exercise}

\begin{exercise}
Do the computations for the remaining AD graphs.
\end{exercise}

As bonus exercise, reiterated, read the work of Jones \cite{jo1}, \cite{jo2}, \cite{jo3}, \cite{jo4}.

\chapter{Spectral measures}

\section*{14a. Linear operators}

As a continuation of what we did in the previous chapter, and with the aim of reaching to something really conceptual, advanced and powerful, here will be our principle:

\begin{principle}
Advanced probability and mathematics simplify when assuming that your variables no longer commute, $fg\neq gf$, and with examples for this coming from:
\begin{enumerate}
\item Usual matrices, $A\in M_N(\mathbb C)$, which generically do not commute, indeed.

\item Random matrices, $Z\in M_N(L^\infty(X))$, with $X$ being a probability space.

\item Linear operators $T:H\to H$, with $H$ being a Hilbert space.
\end{enumerate}
\end{principle}

To be more precise, what we have here is some sort of hierarchy, with the matrices in (1), that we certainly love, being at the bottom. Then we have the random matrices in (2), appearing as a smart joint generalization of the matrices in (1), which appear with $X=\{.\}$, and of the usual random variables, which appear at $N=1$. And then we have the operators in (3), that we will get to love too, which generalize the random matrices in (2), appearing via $H=\mathbb C^N\otimes L^2(X)$, via some standard identifications.

\bigskip

As for the principle itself, you will have to trust me here, I mean I'm not talking nonsense, that is promised. Getting started now, let us first talk about operators:

\index{infinite matrix}
\index{linear operator}
\index{bounded operator}
\index{continuous operator}

\begin{theorem}
Given a Hilbert space $H$, consider the linear operators $T:H\to H$, and for each such operator define its norm by the following formula:
$$||T||=\sup_{||x||=1}||Tx||$$
The operators which are bounded, $||T||<\infty$, form then a complex algebra $B(H)$, which is complete with respect to $||.||$. When $H$ comes with a basis $\{e_i\}_{i\in I}$, we have
$$B(H)\subset M_I(\mathbb C)$$
with the correspondence $T\to M$ coming via the usual linear algebra formulae, namely:
$$T(x)=Mx\quad,\quad M_{ij}=<Te_j,e_i>$$
In infinite dimensions, the inclusion $B(H)\subset M_I(\mathbb C)$ is not an equality.
\end{theorem}

\begin{proof}
This is something straightforward, the idea being as follows:

\medskip

(1) The fact that we have indeed an algebra, satisfying the product condition in the statement, follows from the following estimates, which are all elementary:
$$||S+T||\leq||S||+||T||\quad,\quad 
||\lambda T||=|\lambda|\cdot||T||\quad,\quad 
||ST||\leq||S||\cdot||T||$$

(2) Regarding now the completness assertion, if $\{T_n\}\subset B(H)$ is Cauchy then $\{T_nx\}$ is Cauchy for any $x\in H$, so we can define the limit $T=\lim_{n\to\infty}T_n$ by setting:
$$Tx=\lim_{n\to\infty}T_nx$$

Let us first check that the application $x\to Tx$ is linear. We have:
\begin{eqnarray*}
T(x+y)
&=&\lim_{n\to\infty}T_n(x+y)\\
&=&\lim_{n\to\infty}T_n(x)+T_n(y)\\
&=&\lim_{n\to\infty}T_n(x)+\lim_{n\to\infty}T_n(y)\\
&=&T(x)+T(y)
\end{eqnarray*}

Similarly, we have $T(\lambda x)=\lambda T(x)$, and we conclude that $x\to Tx$ is linear.

\medskip

(3) With this done, it remains to prove now that we have $T\in B(H)$, and that $T_n\to T$ in norm. For this purpose, observe that we have:
\begin{eqnarray*}
||T_n-T_m||\leq\varepsilon\ ,\ \forall n,m\geq N
&\implies&||T_nx-T_mx||\leq\varepsilon\ ,\ \forall||x||=1\ ,\ \forall n,m\geq N\\
&\implies&||T_nx-Tx||\leq\varepsilon\ ,\ \forall||x||=1\ ,\ \forall n\geq N\\
&\implies&||T_Nx-Tx||\leq\varepsilon\ ,\ \forall||x||=1\\
&\implies&||T_N-T||\leq\varepsilon
\end{eqnarray*}

But this gives both $T\in B(H)$, and $T_N\to T$ in norm, and we are done.

\medskip

(4) Regarding the embedding, the correspondence $T\to M$ in the statement is indeed linear, and its kernel is $\{0\}$, so we have indeed an embedding as follows, as claimed:
$$B(H)\subset M_I(\mathbb C)$$

In finite dimensions we have an isomorphism, because any matrix $M\in M_N(\mathbb C)$ determines a linear operator $T:\mathbb C^N\to\mathbb C^N$, given by the following formula:
$$<Te_j,e_i>=M_{ij}$$

However, in infinite dimensions we have matrices not producing operators, as for instance the all-one matrix, so the embedding $B(H)\subset M_I(\mathbb C)$ is not an isomorphism.
\end{proof}

As a second and last basic result regarding the operators, we will need as well:

\index{adjoint operator}
\index{adjoint matrix}

\begin{theorem}
Each operator $T\in B(H)$ has an adjoint $T^*\in B(H)$, given by: 
$$<Tx,y>=<x,T^*y>$$
The operation $T\to T^*$ is antilinear, antimultiplicative, involutive, and satisfies:
$$||T||=||T^*||\quad,\quad ||TT^*||=||T||^2$$
When $H$ comes with a basis $\{e_i\}_{i\in I}$, the operation $T\to T^*$ corresponds to
$$(M^*)_{ij}=\overline{M}_{ji}$$ 
at the level of the associated matrices $M\in M_I(\mathbb C)$.
\end{theorem}

\begin{proof}
This is standard too, and can be proved in 3 steps, as follows:

\medskip

(1) The existence of the adjoint operator $T^*$, given by the formula in the statement, comes from the fact that the function $\varphi(x)=<Tx,y>$ being a linear map $H\to\mathbb C$, we must have a formula as follows, for a certain vector $T^*y\in H$:
$$\varphi(x)=<x,T^*y>$$

Moreover, since this vector is unique, $T^*$ is unique too, and we have as well:
$$(S+T)^*=S^*+T^*\quad,\quad
(\lambda T)^*=\bar{\lambda}T^*\quad,\quad 
(ST)^*=T^*S^*\quad,\quad 
(T^*)^*=T$$

Observe also that we have indeed $T^*\in B(H)$, because:
\begin{eqnarray*}
||T||
&=&\sup_{||x||=1}\sup_{||y||=1}<Tx,y>\\
&=&\sup_{||y||=1}\sup_{||x||=1}<x,T^*y>\\
&=&||T^*||
\end{eqnarray*}

(2) Regarding now $||TT^*||=||T||^2$, which is a key formula, observe that we have:
$$||TT^*||
\leq||T||\cdot||T^*||
=||T||^2$$

On the other hand, we have as well the following estimate:
\begin{eqnarray*}
||T||^2
&=&\sup_{||x||=1}|<Tx,Tx>|\\
&=&\sup_{||x||=1}|<x,T^*Tx>|\\
&\leq&||T^*T||
\end{eqnarray*}

By replacing $T\to T^*$ we obtain from this $||T||^2\leq||TT^*||$, as desired.

\medskip

(3) Finally, when $H$ comes with a basis, the formula $<Tx,y>=<x,T^*y>$ applied with $x=e_i$, $y=e_j$ translates into the formula $(M^*)_{ij}=\overline{M}_{ji}$, as desired.
\end{proof}

Let us discuss now the diagonalization problem for the operators $T\in B(H)$, in analogy with the diagonalization problem for the usual matrices $A\in M_N(\mathbb C)$. We first have:

\begin{definition}
The spectrum of an operator $T\in B(H)$ is the set
$$\sigma(T)=\left\{\lambda\in\mathbb C\Big|T-\lambda\not\in B(H)^{-1}\right\}$$
where $B(H)^{-1}\subset B(H)$ is the set of invertible operators.
\end{definition}

As a basic example, in the finite dimensional case, $H=\mathbb C^N$, the spectrum of a usual matrix $A\in M_N(\mathbb C)$ is the collection of its eigenvalues, taken without multiplicities. We will see many other examples. In general, the spectrum has the following properties:

\index{eigenvalue}

\begin{proposition}
The spectrum of $T\in B(H)$ contains the eigenvalue set
$$\varepsilon(T)=\left\{\lambda\in\mathbb C\Big|\ker(T-\lambda)\neq\{0\}\right\}$$
and $\varepsilon(T)\subset\sigma(T)$ is an equality in finite dimensions, but not in infinite dimensions.
\end{proposition}

\begin{proof}
We have several assertions here, the idea being as follows:

\medskip

(1) First of all, the eigenvalue set is indeed the one in the statement, because $Tx=\lambda x$ tells us precisely that $T-\lambda$ must be not injective. The fact that we have $\varepsilon(T)\subset\sigma(T)$ is clear as well, because if $T-\lambda$ is not injective, it is not bijective.

\medskip

(2) In finite dimensions we have $\varepsilon(T)=\sigma(T)$, because $T-\lambda$ is injective if and only if it is bijective, with the boundedness of the inverse being automatic. 

\medskip

(3) In infinite dimensions we can assume $H=l^2(\mathbb N)$, and the shift operator $S(e_i)=e_{i+1}$ is injective but not surjective. Thus $0\in\sigma(T)-\varepsilon(T)$.
\end{proof}

As a second basic result about the spectrum, which is of great use, we have:

\begin{proposition}
The spectrum of an operator $T\in B(H)$ satisfies:
$$\sigma(T)\subset D_0(||T||))$$
In other words, the spectral radius $\rho(T)=\sup_{\lambda\in\sigma(T)}|\lambda|$ satisfies $\rho(T)\leq||T||$. 
\end{proposition}

\begin{proof}
We use the standard fact that for $||S||<1$, we have the following formula:
$$(1-S)^{-1}=1+S+S^2+\ldots$$

Indeed, by using this inversion formula, we have the following computation:
\begin{eqnarray*}
\lambda>||T||
&\implies&\Big|\Big|\frac{T}{\lambda}\Big|\Big|<1\\
&\implies&1-\frac{T}{\lambda}\in B(H)^{-1}\\
&\implies&\lambda-T\in B(H)^{-1}\\
&\implies&\lambda\notin\sigma(T)
\end{eqnarray*}

Thus, we are led to the conclusion in the statement.
\end{proof}

Let us develop now some more general theory. Here is a key result about products:

\index{spectrum of products}

\begin{theorem}
We have the following formula, valid for any operators $S,T$:
$$\sigma(ST)\cup\{0\}=\sigma(TS)\cup\{0\}$$
In finite dimensions we have $\sigma(ST)=\sigma(TS)$, but this fails in infinite dimensions.
\end{theorem}

\begin{proof}
There are several assertions here, the idea being as follows:

\medskip

(1) Let us first prove the main assertion, stating that the sets $\sigma(ST),\sigma(TS)$ coincide outside 0. We first prove that we have the following implication:
$$1\notin\sigma(ST)\implies1\notin\sigma(TS)$$

So, assume that the operator $1-ST$ is invertible, with inverse denoted $R$:
$$R=(1-ST)^{-1}$$

We have then $RST=R-1$, and by using this, we have the following computation:
\begin{eqnarray*}
(1+TRS)(1-TS)
&=&1+TRS-TS-TRSTS\\
&=&1+TRS-TS-TRS+TS\\
&=&1
\end{eqnarray*}

Similarly, $(1-TS)(1+TRS)=1$. Thus $1-TS$ is invertible, proving our claim. Now by multiplying by scalars, we deduce that for any $\lambda\in\mathbb C-\{0\}$ we have, as desired:
$$\lambda\notin\sigma(ST)\implies\lambda\notin\sigma(TS)$$

(2) Regarding now the counterexample to the formula $\sigma(ST)=\sigma(TS)$, in general, let us take $S$ to be the shift on $H=L^2(\mathbb N)$, given by the following formula:
$$S(e_i)=e_{i+1}$$

As for $T$, we can take it to be the adjoint of $S$, which is the following operator:
$$S^*(e_i)=\begin{cases}
e_{i-1}&{\rm if}\ i>0\\
0&{\rm if}\ i=0
\end{cases}$$

Let us compose now these two operators. In one sense, we have:
$$S^*S=1\implies 0\notin\sigma(SS^*)$$

In the other sense, however, the situation is different, as follows:
$$SS^*=Proj(e_0^\perp)\implies 0\in\sigma(SS^*)$$

Thus, the spectra do not match on $0$, and we have our counterexample, as desired.
\end{proof}

As our next result, inspired from the linear algebra from chapter 9, we have:

\index{rational calculus}

\begin{theorem}
We have the ``rational functional calculus'' formula
$$\sigma(f(T))=f(\sigma(T))$$
valid for any rational function $f\in\mathbb C(X)$ having poles outside $\sigma(T)$.
\end{theorem}

\begin{proof}
This can be proved in two steps, as follows:

\medskip

(1) Assume first that our rational function $f\in\mathbb C(X)$ is a usual polynomial $P\in\mathbb C[X]$. We pick a scalar $\lambda\in\mathbb C$, and we decompose the polynomial $P-\lambda$, as follows:
$$P(X)-\lambda=c(X-r_1)\ldots(X-r_n)$$

We have then the following equivalences, which give the result:
\begin{eqnarray*}
\lambda\notin\sigma(P(T))
&\iff&P(T)-\lambda\in B(H)^{-1}\\
&\iff&c(T-r_1)\ldots(T-r_n)\in B(H)^{-1}\\
&\iff&T-r_1,\ldots,T-r_n\in B(H)^{-1}\\
&\iff&r_1,\ldots,r_n\notin\sigma(T)\\
&\iff&\lambda\notin P(\sigma(T))
\end{eqnarray*}

(2) In general now, we pick a scalar $\lambda\in\mathbb C$, we write $f=P/Q$, and we set $F=P-\lambda Q$. By using what we found in (1), for this polynomial $F\in\mathbb C[X]$, we obtain:
\begin{eqnarray*}
\lambda\in\sigma(f(T))
&\iff&F(T)\notin B(H)^{-1}\\
&\iff&0\in\sigma(F(T))\\
&\iff&0\in F(\sigma(T))\\
&\iff&\exists\mu\in\sigma(T),F(\mu)=0\\
&\iff&\lambda\in f(\sigma(T))
\end{eqnarray*}

Thus, we are led to the formula in the statement.
\end{proof}

As a first application of the above methods, we have the following key result:

\begin{theorem}
The following happen:
\begin{enumerate}
\item For a unitary operator, $U^*=U^{-1}$, we have $\sigma(U)\subset\mathbb T$. 

\item For a self-adjoint operator, $T=T^*$, we have $\sigma(T)\subset\mathbb R$.
\end{enumerate}
\end{theorem}

\begin{proof}
This is something quite tricky, based on Theorem 14.8, as follows:

\medskip

(1) Assuming $U^*=U^{-1}$, we have the following norm computation:
$$||U||
=\sqrt{||UU^*||}
=\sqrt{1}
=1$$

With $D$ being the unit disk, we obtain from this, by using Proposition 14.6:
$$\sigma(U)\subset D$$

On the other hand, once again by using $U^*=U^{-1}$, we have as well:
$$||U^{-1}||
=||U^*||
=||U||
=1$$

Thus, as before with $D$ being the unit disk in the complex plane, we have:
$$\sigma(U^{-1})\subset D$$

Now by using Theorem 14.8, we obtain $\sigma(U)
\subset D\cap D^{-1}
=\mathbb T$, as desired.

\medskip

(2) Consider the following rational function, depending on a parameter $r\in\mathbb R$:
$$f(z)=\frac{z+ir}{z-ir}$$

Then for $r>>0$ the operator $f(T)$ is well-defined, and we have:
$$\left(\frac{T+ir}{T-ir}\right)^*
=\frac{T-ir}{T+ir}
=\left(\frac{T+ir}{T-ir}\right)^{-1}$$

Thus $f(T)$ is unitary, and by (1) we have $\sigma(T)\subset f^{-1}(\mathbb T)=\mathbb R$, as desired.
\end{proof}

Next, we have the following key result, which is something quite far-reaching:

\begin{theorem}
The spectral radius of an operator $T\in B(H)$ is given by
$$\rho(T)=\lim_{n\to\infty}||T^n||^{1/n}$$
and in this formula, we can replace the limit by an inf.
\end{theorem}

\begin{proof}
We have several things to be proved, the idea being as follows:

\medskip

(1) Our first claim is that the numbers $u_n=||T^n||^{1/n}$ satisfy:
$$(n+m)u_{n+m}\leq nu_n+mu_m$$

Indeed, we have the following estimate, obtained by using the well-known Young inequality $ab\leq a^p/p+b^q/q$, with exponents $p=(n+m)/n$ and $q=(n+m)/m$:
\begin{eqnarray*}
u_{n+m}
&=&||T^{n+m}||^{1/(n+m)}\\
&\leq&||T^n||^{1/(n+m)}||T^m||^{1/(n+m)}\\
&\leq&||T^n||^{1/n}\cdot\frac{n}{n+m}+||T^m||^{1/m}\cdot\frac{m}{n+m}\\
&=&\frac{nu_n+mu_m}{n+m}
\end{eqnarray*}

(2) Our second claim is that the second assertion holds, namely:
$$\lim_{n\to\infty}||T^n||^{1/n}=\inf_n||T^n||^{1/n}$$

For this purpose, we just need the inequality found in (1). Indeed, fix $m\geq1$, let $n\geq1$, and write $n=lm+r$ with $0\leq r\leq m-1$. By using twice $u_{ab}\leq u_b$, we get:
$$u_n
\leq\frac{1}{n}( lmu_{lm}+ru_r)
\leq u_{m}+\frac{r}{n}\,u_1$$

It follows that we have $\lim\sup_nu_n\leq u_m$, which proves our claim.

\medskip

(3) Summarizing, we are left with proving the main formula, which is as follows, and with the remark that we already know that the sequence on the right converges:
$$\rho(T)=\lim_{n\to\infty}||T^n||^{1/n}$$

In one sense, we can use the polynomial calculus formula $\sigma(T^n)=\sigma(T)^n$. Indeed, this gives the following estimate, valid for any $n$, as desired:
$$\rho(T)
=\sup_{\lambda\in\sigma(T)}|\lambda|
=\sup_{\rho\in\sigma(T)^n}|\rho|^{1/n}
=\sup_{\rho\in\sigma(T^n)}|\rho|^{1/n}
=\rho(T^n)^{1/n}
\leq||T^n||^{1/n}$$

(4) For the reverse inequality, we fix a number $\rho>\rho(T)$, and we want to prove that we have $\rho\geq\lim_{n\to\infty}||T^n||^{1/n}$. By using the Cauchy formula, we have:
\begin{eqnarray*}
\frac{1}{2\pi i}\int_{|z|=\rho}\frac{z^n}{z-T}\,dz
&=&\frac{1}{2\pi i}\int_{|z|=\rho}\sum_{k=0}^\infty z^{n-k-1}T^k\,dz\\
&=&\sum_{k=0}^\infty\frac{1}{2\pi i}\left(\int_{|z|=\rho}z^{n-k-1}dz\right)T^k\\
&=&\sum_{k=0}^\infty\delta_{n,k+1}T^k\\
&=&T^{n-1}
\end{eqnarray*}

By applying the norm we obtain from this formula the following estimate:
$$||T^{n-1}||
\leq\frac{1}{2\pi}\int_{|z|=\rho}\left|\left|\frac{z^n}{z-T}\right|\right|\,dz
\leq\rho^n\cdot\sup_{|z|=\rho}\left|\left|\frac{1}{z-T}\right|\right|$$

Since the sup does not depend on $n$, by taking $n$-th roots, we obtain in the limit:
$$\rho\geq\lim_{n\to\infty}||T^n||^{1/n}$$

Now recall that $\rho$ was by definition an arbitrary number satisfying $\rho>\rho(T)$. Thus, we have obtained the following estimate, valid for any $T\in B(H)$:
$$\rho(T)\geq\lim_{n\to\infty}||T^n||^{1/n}$$

Thus, we are led to the conclusion in the statement.
\end{proof}

In the case of the normal elements, we have the following finer result:

\index{normal operator}

\begin{theorem}
The spectral radius of a normal element,
$$TT^*=T^*T$$
is equal to its norm.
\end{theorem}

\begin{proof}
In the case $T=T^*$ we have $||T^n||=||T||^n$ for any exponent of the form $n=2^k$, by using $||TT^*||=||T||^2$, and by taking $n$-th roots we get, as desired:
$$\rho(T)\geq||T||$$

In the general normal case $TT^*=T^*T$ we have $T^n(T^n)^*=(TT^*)^n$, and by using this, along with the above result for self-adjoints, applied to $TT^*$, we obtain:
\begin{eqnarray*}
\rho(T)
&=&\lim_{n\to\infty}||T^n||^{1/n}\\
&=&\sqrt{\lim_{n\to\infty}||T^n(T^n)^*||^{1/n}}\\
&=&\sqrt{\lim_{n\to\infty}||(TT^*)^n||^{1/n}}\\
&=&\sqrt{\rho(TT^*)}\\
&=&\sqrt{||T||^2}\\
&=&||T||
\end{eqnarray*}

Thus, we are led to the conclusion in the statement.
\end{proof}

\section*{14b. Spectral theorems}

By using Theorem 14.11 we can say a number of things about the normal operators, commonly known as ``spectral theorem for normal operators''. We first have:

\begin{theorem}
Given $T\in B(H)$ normal, we have a morphism of algebras
$$\mathbb C[X]\to B(H)\quad,\quad 
P\to P(T)$$
having the properties $||P(T)||=||P_{|\sigma(T)}||$, and $\sigma(P(T))=P(\sigma(T))$.
\end{theorem}

\begin{proof}
This is an improvement of Theorem 14.8 for polynomials, in the normal case, with the extra assertion being the norm estimate. But the element $P(T)$ being normal, we can apply to it the spectral radius formula for normal elements, and we obtain:
\begin{eqnarray*}
||P(T)||
&=&\rho(P(T))\\
&=&\sup_{\lambda\in\sigma(P(T))}|\lambda|\\
&=&\sup_{\lambda\in P(\sigma(T))}|\lambda|\\
&=&||P_{|\sigma(T)}||
\end{eqnarray*}

Thus, we are led to the conclusions in the statement.
\end{proof}

At a more advanced level now, we have the following result:

\begin{theorem}
Given $T\in B(H)$ normal, we have a morphism of algebras
$$C(\sigma(T))\to B(H)\quad,\quad 
f\to f(T)$$
which is isometric, $||f(T)||=||f||$, and has the property $\sigma(f(T))=f(\sigma(T))$.
\end{theorem}

\begin{proof}
The idea here is to ``complete'' the morphism in Theorem 14.12. Indeed, by Stone-Weierstrass, that morphism has a unique isometric extension, as follows:
$$C(\sigma(T))\to B(H)\quad,\quad  
f\to f(T)$$

It remains to prove $\sigma(f(T))=f(\sigma(T))$, and we can do this by double inclusion:

\medskip

``$\subset$'' Given a continuous function $f\in C(\sigma(T))$, we must prove that we have:
$$\lambda\notin f(\sigma(T))\implies\lambda\notin\sigma(f(T))$$

For this purpose, consider the following function, which is well-defined:
$$\frac{1}{f-\lambda}\in C(\sigma(T))$$

We can therefore apply this function to $T$, and we obtain:
$$\left(\frac{1}{f-\lambda}\right)T=\frac{1}{f(T)-\lambda}$$

In particular $f(T)-\lambda$ is invertible, so  $\lambda\notin\sigma(f(T))$, as desired.

\medskip

``$\supset$'' Given a continuous function $f\in C(\sigma(T))$, we must prove that we have: 
$$\lambda\in f(\sigma(T))\implies\lambda\in\sigma(f(T))$$

But this is the same as proving that we have:
$$\mu\in\sigma(T)\implies f(\mu)\in\sigma(f(T))$$

For this purpose, we approximate our function by polynomials, $P_n\to f$, and we examine the following convergence, which follows from $P_n\to f$:
$$P_n(T)-P_n(\mu)\to f(T)-f(\mu)$$

We know from polynomial functional calculus that we have:
$$P_n(\mu)
\in P_n(\sigma(T))
=\sigma(P_n(T))$$

Thus, the operators $P_n(T)-P_n(\mu)$ are not invertible. On the other hand, we know that the set formed by the invertible operators is open, so its complement is closed. Thus the limit $f(T)-f(\mu)$ is not invertible either, and so $f(\mu)\in\sigma(f(T))$, as desired.
\end{proof}

Even more generally now, we have the following result:

\begin{theorem}
Given $T\in B(H)$ normal, we have a morphism of algebras as follows, with $L^\infty$ standing for abstract measurable functions, or Borel functions,
$$L^\infty(\sigma(T))\to B(H)\quad,\quad 
f\to f(T)$$
which is isometric, $||f(T)||=||f||$, and has the property $\sigma(f(T))=f(\sigma(T))$.
\end{theorem}

\begin{proof}
As before, the idea will be that of ``completing'' what we have:

\medskip

(1) Given a vector $x\in H$, consider the following functional:
$$C(\sigma(T))\to\mathbb C\quad,\quad 
g\to<g(T)x,x>$$

By the Riesz theorem, this functional must be the integration with respect to a certain measure $\mu$ on the space $\sigma(T)$. Thus, we have a formula as follows:
$$<g(T)x,x>=\int_{\sigma(T)}g(z)d\mu(z)$$

Now given an arbitrary Borel function $f\in L^\infty(\sigma(T))$, as in the statement, we can define a number $<f(T)x,x>\in\mathbb C$, by using exactly the same formula, namely:
$$<f(T)x,x>=\int_{\sigma(T)}f(z)d\mu(z)$$

Thus, we have managed to define numbers $<f(T)x,x>\in\mathbb C$, for all vectors $x\in H$, and in addition we can recover these numbers as follows, with $g_n\in C(\sigma(T))$:
$$<f(T)x,x>=\lim_{g_n\to f}<g_n(T)x,x>$$ 

(2) In order to define now numbers $<f(T)x,y>\in\mathbb C$, for all vectors $x,y\in H$, we can use a polarization trick. Indeed, for any operator $S\in B(H)$ we have:
$$<S(x+y),x+y>=<Sx,x>+<Sy,y>+<Sx,y>+<Sy,x>$$

By replacing $y\to iy$, we have as well the following formula:
$$<S(x+iy),x+iy>=<Sx,x>+<Sy,y>-i<Sx,y>+i<Sy,x>$$

By multiplying this formula by $i$, and summing with the first one, we obtain:
\begin{eqnarray*}
<S(x+y),x+y>+i<S(x+iy),x+iy>
&=&(1+i)[<Sx,x>+<Sy,y>]\\
&+&2<Sx,y>
\end{eqnarray*}

(3) But with this, we can finish. Indeed, by combining (1,2), given a Borel function $f\in L^\infty(\sigma(T))$, we can define numbers $<f(T)x,y>\in\mathbb C$ for any $x,y\in H$, and we obtain in this way a certain operator $f(T)\in B(H)$, having all the desired properties.
\end{proof}

Good news, we can now diagonalize the normal operators. Let us start with:

\begin{theorem}
Any self-adjoint operator $T\in B(H)$ can be diagonalized,
$$T=U^*M_fU$$
with $U:H\to L^2(X)$ being a unitary operator from $H$ to a certain $L^2$ space associated to $T$, with $f:X\to\mathbb R$ being a certain function, once again associated to $T$, and with
$$M_f(g)=fg$$
being the usual multiplication operator by $f$, on the Hilbert space $L^2(X)$.
\end{theorem}

\begin{proof}
The construction of $U,f$ can be done in several steps, as follows:

\medskip

(1) We first prove the result in the special case where our operator $T$ has a cyclic vector $x\in H$, with this meaning that the following holds:
$$\overline{span\left(T^kx\Big|n\in\mathbb N\right)}=H$$

For this purpose, let us go back to the proof of Theorem 14.14. We will use the following formula from there, with $\mu$ being the measure on $X=\sigma(T)$ associated to $x$:
$$<g(T)x,x>=\int_{\sigma(T)}g(z)d\mu(z)$$

Our claim is that we can define a unitary $U:H\to L^2(X)$, first on the dense part spanned by the vectors $T^kx$, by the following formula, and then by continuity:
$$U[g(T)x]=g$$

Indeed, the following computation shows that $U$ is well-defined, and isometric:
\begin{eqnarray*}
||g(T)x||^2
&=&<g(T)x,g(T)x>\\
&=&<g(T)^*g(T)x,x>\\
&=&<|g|^2(T)x,x>\\
&=&\int_{\sigma(T)}|g(z)|^2d\mu(z)\\
&=&||g||_2^2
\end{eqnarray*}

We can then extend $U$ by continuity into a unitary $U:H\to L^2(X)$, as claimed. Now observe that we have the following formula:
\begin{eqnarray*}
UTU^*g
&=&U[Tg(T)x]\\
&=&U[(zg)(T)x]\\
&=&zg
\end{eqnarray*} 

Thus our result is proved in the present case, with $U$ as above, and with $f(z)=z$.

\medskip

(2) We discuss now the general case. Our first claim is that $H$ has a decomposition as follows, with each $H_i$ being invariant under $T$, and admitting a cyclic vector $x_i$:
$$H=\bigoplus_iH_i$$

Indeed, this is something elementary, the construction being by recurrence in finite dimensions, in the obvious way, and by using the Zorn lemma in general. Now with this decomposition in hand, we can make a direct sum of the diagonalizations obtained in (1), for each of the restrictions $T_{|H_i}$, and we obtain the formula in the statement.
\end{proof}

We have the following technical generalization of the above result:

\begin{theorem}
Any family of commuting self-adjoint operators $T_i\in B(H)$ can be jointly diagonalized,
$$T_i=U^*M_{f_i}U$$
with $U:H\to L^2(X)$ being a unitary operator from $H$ to a certain $L^2$ space associated to $\{T_i\}$, with $f_i:X\to\mathbb R$ being certain functions, once again associated to $T_i$, and with
$$M_{f_i}(g)=f_ig$$
being the usual multiplication operator by $f_i$, on the Hilbert space $L^2(X)$.
\end{theorem}

\begin{proof}
This is similar to the proof of Theorem 14.15, by suitably modifying the measurable calculus formula, and the measure $\mu$ itself, as to have this formula working for all the operators $T_i$. With this modification done, everything extends.
\end{proof}

We can now discuss the case of the arbitrary normal operators, as follows:

\begin{theorem}
Any normal operator $T\in B(H)$ can be diagonalized,
$$T=U^*M_fU$$
with $U:H\to L^2(X)$ being a unitary operator from $H$ to a certain $L^2$ space associated to $T$, with $f:X\to\mathbb C$ being a certain function, once again associated to $T$, and with
$$M_f(g)=fg$$
being the usual multiplication operator by $f$, on the Hilbert space $L^2(X)$.
\end{theorem}

\begin{proof}
Consider the decomposition of $T$ into its real and imaginary parts:
$$T=\frac{T+T^*}{2}+i\cdot\frac{T-T^*}{2i}$$

We know that the real and imaginary parts are self-adjoint operators. Now since $T$ was assumed to be normal, $TT^*=T^*T$, these real and imaginary parts commute:
$$\left[\frac{T+T^*}{2}\,,\,\frac{T-T^*}{2i}\right]=0$$

Thus Theorem 14.16 applies to these real and imaginary parts, and gives the result.
\end{proof}

This was for our series of diagonalization theorems. There is of course one more result here, regarding the families of commuting normal operators, as follows:

\index{commuting normal operators}
\index{diagonalization}

\begin{theorem}
Any family of commuting normal operators $T_i\in B(H)$ can be jointly diagonalized,
$$T_i=U^*M_{f_i}U$$
with $U:H\to L^2(X)$ being a unitary operator from $H$ to a certain $L^2$ space associated to $\{T_i\}$, with $f_i:X\to\mathbb C$ being certain functions, once again associated to $T_i$, and with
$$M_{f_i}(g)=f_ig$$
being the usual multiplication operator by $f_i$, on the Hilbert space $L^2(X)$.
\end{theorem}

\begin{proof}
This is similar to the proof of Theorem 14.15 and Theorem 14.17, by combining the arguments there. To be more precise, this follows as Theorem 14.15, by using the decomposition trick from the proof of Theorem 14.17.
\end{proof}

With the above diagonalization results in hand, we can now ``fix'' the continuous and measurable functional calculus theorems, with a key complement, as follows:

\index{polynomial calculus}
\index{rational calculus}
\index{holomorphic calculus}
\index{continuous calculus}
\index{measurable calculus}
\index{adjoint operator}

\begin{theorem}
Given a normal operator $T\in B(H)$, the following hold, for both the continuous calculus and the measurable calculus morphisms:
\begin{enumerate}
\item These morphisms are $*$-morphisms.

\item The function $\bar{z}$ gets mapped to $T^*$.

\item The functions $Re(z),Im(z)$ get mapped to $Re(T),Im(T)$.

\item The function $|z|^2$ gets mapped to $TT^*=T^*T$.

\item If $f$ is real, then $f(T)$ is self-adjoint. 
\end{enumerate}
\end{theorem}

\begin{proof}
These assertions are more or less equivalent, with (1) being the main one, which obviously implies everything else. But this assertion (1) follows from the diagonalization result for the normal operators, from Theorem 14.17.
\end{proof}

Getting now to applications, the above results are quite powerful, and many things can be said, in analogy with what we know about usual matrices. Let us record here:

\index{polar decomposition}
\index{partial isometry}

\begin{theorem}
Any bounded operator $T\in B(H)$ can be decomposed as
$$T=U|T|$$
with $U$ being a partial isometry, and with $|T|=\sqrt{T^*T}$.
\end{theorem}

\begin{proof}
The operator $T^*T$ being self-adjoint, and even positive, in the sense that we have $<T^*Tx,x>\geq0$ for any $x\in H$, we can extract its square root $|T|=\sqrt{T^*T}$, by using the spectral theorem. Now observe that we have the following formula:
\begin{eqnarray*}
<|T|x,|T|y>
&=&<x,|T|^2y>\\
&=&<x,T^*Ty>\\
&=&<Tx,Ty>
\end{eqnarray*}

We conclude that the following linear application is well-defined, and isometric:
$$U:Im|T|\to Im(T)\quad,\quad 
|T|x\to Tx$$

Now by continuity we can extend this isometry $U$ into an isometry between certain Hilbert subspaces of $H$, as follows:
$$U:\overline{Im|T|}\to\overline{Im(T)}\quad,\quad 
|T|x\to Tx$$

Moreover, we can further extend $U$ into a partial isometry $U:H\to H$, by setting $Ux=0$, for any $x\in\overline{Im|T|}^\perp$, and with this convention, the result follows. 
\end{proof}

\section*{14c. Operator algebras}

Getting now a bit abstract, we need to talk about operator algebras, as to suitably cover the random matrix algebras $M_N(L^\infty(X))$. Which can be something quite tricky, but fortunately, we have the following quick and clever definition, due to Gelfand:

\begin{definition}
An abstract operator algebra, or $C^*$-algebra, is a complex algebra $A$ having a norm $||.||$ and an involution $*$, subject to the following conditions:
\begin{enumerate}
\item $A$ is closed with respect to the norm.

\item We have $||aa^*||=||a||^2$, for any $a\in A$.
\end{enumerate}
\end{definition}

In other words, what we did here is to axiomatize the abstract properties of the operator algebras $A\subset B(H)$, assumed to be norm closed and stable under taking adjoints, both very natural conditions, without reference to the ambient Hilbert space $H$.

\bigskip

As basic examples now, we have the usual matrix algebras $M_N(\mathbb C)$, with the norm and the involution being the usual matrix norm and involution, given by:
$$||A||=\sup_{||x||=1}||Ax||\quad,\quad 
(A^*)_{ij}=\overline{A}_{ji}$$

Some other basic examples are the algebras $L^\infty(X)$ of essentially bounded functions $f:X\to\mathbb C$ on a measured space $X$, with the usual norm and involution, namely:
$$||f||=\sup_{x\in X}|f(x)|\quad,\quad 
f^*(x)=\overline{f(x)}$$

We can put our two basic classes of examples together, as follows:

\begin{proposition}
The random matrix algebras $A=M_N(L^\infty(X))$ are $C^*$-algebras, with their usual norm and involution, given by:
$$||Z||=\sup_{x\in X}||Z_x||\quad,\quad 
(Z^*)_{ij}=\overline{Z}_{ij}$$
These algebras generalize both the algebras $M_N(\mathbb C)$, and the algebras $L^\infty(X)$.
\end{proposition}

\begin{proof}
The fact that the $C^*$-algebra axioms are satisfied is clear from definitions. As for the last assertion, this follows by taking $X=\{.\}$ and $N=1$, respectively.
\end{proof}

We can in fact say more about the above algebras, as follows:

\begin{proposition}
Any algebra of type $L^\infty(X)$ is an operator algebra, as follows:
$$L^\infty(X)\subset B(L^2(X))\quad,\quad 
f\to(g\to fg)$$
More generally, any random matrix algebra is an operator algebra, as follows,
$$M_N(L^\infty(X))\subset B\left(\mathbb C^N\otimes L^2(X)\right)$$
with the embedding being the above one, tensored with the identity.
\end{proposition}

\begin{proof}
We have two assertions to be proved, the idea being as follows:

\medskip

(1) Given $f\in L^\infty(X)$, consider the following operator, acting on $H=L^2(X)$:
$$T_f(g)=fg$$

Observe that $T_f$ is indeed well-defined, and bounded as well, because:
$$||fg||_2
=\sqrt{\int_X|f(x)|^2|g(x)|^2d\mu(x)}
\leq||f||_\infty||g||_2$$

The application $f\to T_f$ being linear, involutive, continuous, and injective as well, we obtain in this way a $C^*$-algebra embedding $L^\infty(X)\subset B(H)$, as desired.

\medskip

(2) Regarding the second assertion, this is best viewed in the following way:
\begin{eqnarray*}
M_N(L^\infty(X))
&=&M_N(\mathbb C)\otimes L^\infty(X)\\
&\subset&M_N(\mathbb C)\otimes B(L^2(X))\\
&=&B\left(\mathbb C^N\otimes L^2(X)\right)
\end{eqnarray*}

Here we have used (1), and some standard tensor product identifications.
\end{proof}

Our purpose now will be to develop the spectral theory of $C^*$-algebras, and in particular that of the random matrix algebras $A=M_N(L^\infty(X))$ that we are interested in, one of our objectives being that of talking about spectral measures, in the normal case, in analogy with what we know about the usual matrices. Let us start with:

\begin{theorem}
Given an element $a\in A$ of a $C^*$-algebra, define its spectrum as:
$$\sigma(a)=\left\{\lambda\in\mathbb C\Big|a-\lambda\notin A^{-1}\right\}$$
The following spectral theory results hold, exactly as in the $A=B(H)$ case:
\begin{enumerate}
\item We have $\sigma(ab)\cup\{0\}=\sigma(ba)\cup\{0\}$.

\item We have $\sigma(f(a))=f(\sigma(a))$, for any $f\in\mathbb C(X)$ having poles outside $\sigma(a)$.

\item The spectrum $\sigma(a)$ is compact, non-empty, and contained in $D_0(||a||)$.

\item The spectra of unitaries $(u^*=u^{-1})$ and self-adjoints $(a=a^*)$ are on $\mathbb T,\mathbb R$.

\item The spectral radius of normal elements $(aa^*=a^*a)$ is given by $\rho(a)=||a||$.
\end{enumerate}
In addition, assuming $a\in A\subset B$, the spectra of $a$ with respect to $A$ and to $B$ coincide.
\end{theorem}

\begin{proof}
Here the assertions (1-5), which are of course formulated a bit informally, are known to us for the algebra $A=B(H)$, and the proof in general is similar:

\medskip

(1) Assuming that $1-ab$ is invertible, with inverse $c$, we have $abc=cab=c-1$, and it follows that $1-ba$ is invertible too, with inverse $1+bca$. Thus $\sigma(ab),\sigma(ba)$ agree on $1\in\mathbb C$, and by linearity, it follows that $\sigma(ab),\sigma(ba)$ agree on any point $\lambda\in\mathbb C^*$.

\medskip

(2) The formula $\sigma(f(a))=f(\sigma(a))$ is clear for polynomials, $f\in\mathbb C[X]$, by factorizing $f-\lambda$, with $\lambda\in\mathbb C$. Then, the extension to the rational functions is straightforward, because $P(a)/Q(a)-\lambda$ is invertible precisely when $P(a)-\lambda Q(a)$ is.

\medskip

(3) By using $1/(1-b)=1+b+b^2+\ldots$ for $||b||<1$ we obtain that $a-\lambda$ is invertible for $|\lambda|>||a||$, and so $\sigma(a)\subset D_0(||a||)$. It is also clear that $\sigma(a)$ is closed, so what we have is a compact set. Finally, assuming $\sigma(a)=\emptyset$ the function $f(\lambda)=\varphi((a-\lambda)^{-1})$ is well-defined, for any $\varphi\in A^*$, and by Liouville we get $f=0$, contradiction.

\medskip

(4) Assuming $u^*=u^{-1}$ we have $||u||=1$, and so $\sigma(u)\subset D_0(1)$. But with $f(z)=z^{-1}$ we obtain via (2) that we have as well $\sigma(u)\subset f(D_0(1))$, and this gives $\sigma(u)\subset\mathbb T$. As for the result regarding the self-adjoints, this can be obtained from the result for the unitaries, by using (2) with functions of type $f(z)=(z+it)/(z-it)$, with $t\in\mathbb R$.

\medskip

(5) It is routine to check, by integrating quantities of type $z^n/(z-a)$ over circles centered at the origin, and estimating, that the spectral radius is given by $\rho(a)=\lim||a^n||^{1/n}$. But in the self-adjoint case, $a=a^*$, this gives $\rho(a)=||a||$, by using exponents of type $n=2^k$, and then the extension to the general normal case is straightforward.

\medskip 

(6) Regarding now the last assertion, the inclusion $\sigma_B(a)\subset\sigma_A(a)$ is clear. For the converse, assume $a-\lambda\in B^{-1}$, and set $b=(a-\lambda )^*(a-\lambda )$. We have then:
$$\sigma_A(b)-\sigma_B(b)=\left\{\mu\in\mathbb C-\sigma_B(b)\Big|(b-\mu)^{-1}\in B-A\right\}$$

Thus this difference in an open subset of $\mathbb C$. On the other hand $b$ being self-adjoint, its two spectra are both real, and so is their difference. Thus the two spectra of $b$ are equal, and in particular $b$ is invertible in $A$, and so $a-\lambda\in A^{-1}$, as desired.
\end{proof}

We can now prove a key result about the operator algebras, as follows:

\begin{theorem}[Gelfand]
If $X$ is a compact space, the algebra $C(X)$ of continuous functions on it $f:X\to\mathbb C$ is a $C^*$-algebra, with usual norm and involution, namely:
$$||f||=\sup_{x\in X}|f(x)|\quad,\quad 
f^*(x)=\overline{f(x)}$$
Conversely, any commutative $C^*$-algebra is of this form, $A=C(X)$, with 
$$X=\Big\{\chi:A\to\mathbb C\ ,\ {\rm normed\ algebra\ character}\Big\}$$
with topology making continuous the evaluation maps $ev_a:\chi\to\chi(a)$.
\end{theorem}

\begin{proof}
The first assertion is clear from definitions, and regarding the converse:

\medskip

(1) Given a commutative $C^*$-algebra $A$, let us define indeed a space $X$ as in the statement. Then $X$ is compact, and $a\to ev_a$ is a morphism of algebras, as follows:
$$ev:A\to C(X)$$

(2) We first prove that $ev$ is involutive. We use the following formula, which is similar to the $z=Re(z)+iIm(z)$ decomposition formula for usual complex numbers:
$$a=\frac{a+a^*}{2}+i\cdot\frac{a-a^*}{2i}$$

Thus it is enough to prove $ev_{a^*}=ev_a^*$ for the self-adjoint elements $a$. But this is the same as proving that $a=a^*$ implies that $ev_a$ is a real function, which is in turn true, by Theorem 14.24, because $ev_a(\chi)=\chi(a)$ is an element of $\sigma(a)$, contained in $\mathbb R$.

\medskip

(3) Since $A$ is commutative, each element is normal, so $ev$ is isometric:
$$||ev_a||
=\rho(a)
=||a||$$

(4) It remains to prove that $ev$ is surjective. But this follows from the Stone-Weierstrass theorem, because $ev(A)$ is a closed subalgebra of $C(X)$, which separates the points.
\end{proof}

As a main consequence of the Gelfand theorem, we have:

\begin{theorem}
For any normal element $a\in A$ we have an identification as follows:
$$<a>=C(\sigma(a))$$
In addition, given a function $f\in C(\sigma(a))$, we can apply it to $a$, and we have
$$\sigma(f(a))=f(\sigma(a))$$
which generalizes the previous rational calculus formula, in the normal case.
\end{theorem}

\begin{proof}
Since $a$ is normal, the $C^*$-algebra $<a>$ that is generates is commutative, so if we denote by $X$ the space of the characters $\chi:<a>\to\mathbb C$, we have:
$$<a>=C(X)$$

Now since the map $X\to\sigma(a)$ given by evaluation at $a$ is bijective, we obtain:
$$<a>=C(\sigma(a))$$

Thus, we are dealing here with usual functions, and this gives all the assertions.
\end{proof}

As another application of the Gelfand theorem, based on it, we can formulate:

\index{quantum space}
\index{compact quantum space}

\begin{definition}
Given an arbitrary $C^*$-algebra $A$, we agree to write 
$$A=C(X)$$
and call the abstract space $X$ a compact quantum space.
\end{definition}

In other words, we can define the category of compact quantum spaces $X$ as being the category of the $C^*$-algebras $A$, with the arrows reversed. A morphism $f:X\to Y$ corresponds by definition to a morphism $\Phi:C(Y)\to C(X)$, a product of spaces $X\times Y$ corresponds by definition to a product of algebras $C(X)\otimes C(Y)$, and so on.

\bigskip

All this is of course a bit speculative, but let us have some fun with this. We can easily construct algebras with generators and relations, a bit like we construct algebraic manifolds by imposing polynomial relations to the coordinates, and possibilities here are endless, sky is the limit. As the simplest examples of such constructions, we have:

\index{free sphere}
\index{free real sphere}
\index{free complex sphere}
\index{free coordinates}

\begin{definition}
We have free real and complex spheres, defined via
$$C(S^{N-1}_{\mathbb R,+})=C^*\left(x_1,\ldots,x_N\Big|x_i=x_i^*,\sum_ix_i^2=1\right)$$
$$C(S^{N-1}_{\mathbb C,+})=C^*\left(x_1,\ldots,x_N\Big|\sum_ix_ix_i^*=\sum_ix_i^*x_i=1\right)$$
where the symbol $C^*$ stands for universal enveloping $C^*$-algebra.
\end{definition}

Here the fact that these algebras are indeed well-defined comes from the following estimate, which shows that the biggest $C^*$-norms on these $*$-algebras are bounded:
$$||x_i||^2
=||x_ix_i^*||
\leq\Big|\Big|\sum_ix_ix_i^*\Big|\Big|
=1$$

As a first result now, regarding the above free spheres, we have:

\index{liberation}

\begin{theorem}
We have embeddings of compact quantum spaces, as follows,
$$\xymatrix@R=15mm@C=15mm{
S^{N-1}_{\mathbb R,+}\ar[r]&S^{N-1}_{\mathbb C,+}\\
S^{N-1}_\mathbb R\ar[r]\ar[u]&S^{N-1}_\mathbb C\ar[u]
}$$
and the spaces on top appear as liberations of the spaces on the bottom.
\end{theorem}

\begin{proof}
In practice, we must establish the following isomorphisms, where the symbol $C^*_{comm}$ stands for ``universal commutative $C^*$-algebra generated by'':
$$C(S^{N-1}_\mathbb R)=C^*_{comm}\left(x_1,\ldots,x_N\Big|x_i=x_i^*,\sum_ix_i^2=1\right)$$
$$C(S^{N-1}_\mathbb C)=C^*_{comm}\left(x_1,\ldots,x_N\Big|\sum_ix_ix_i^*=\sum_ix_i^*x_i=1\right)$$

It is enough to establish the second isomorphism. So, consider the second universal commutative $C^*$-algebra $A$ constructed above. Since the standard coordinates on $S^{N-1}_\mathbb C$ satisfy the defining relations for $A$, we have a quotient map of as follows:
$$A\to C(S^{N-1}_\mathbb C)$$

Conversely, let us write $A=C(S)$, by using the Gelfand theorem. The variables $x_1,\ldots,x_N$ become in this way true coordinates, providing us with an embedding $S\subset\mathbb C^N$. Also, the quadratic relations take the following form, over the space $S$:
$$\sum_i|x_i|^2=1$$

Thus we have an embedding $S\subset S^{N-1}_\mathbb C$. But this embedding provides us with a quotient map $C(S^{N-1}_\mathbb C)\to A$, as desired, and this gives all the results.
\end{proof}

Quite interesting all this, hope you agree with me. In practice, however, things are a bit more complicated than this, because if we want to do probability over our free spheres, we need uniform measures on them. And these uniform measures must be, by definition, invariant under the action of certain free quantum groups $O_N^+,U_N^+$, still in need to be defined. But we will come to this, with a positive answer, later, in chapter 15.

\section*{14d. Spectral measures}

In order to get now towards noncommutative probability, we first have to develop the theory of positive elements, and linear forms. First, we have the following result:

\begin{proposition}
For an element $a\in A$, the following are equivalent:
\begin{enumerate}
\item $a$ is positive, in the sense that $\sigma(a)\subset[0,\infty)$.

\item $a=b^2$, for some $b\in A$ satisfying $b=b^*$.

\item $a=cc^*$, for some $c\in A$.
\end{enumerate}
\end{proposition}

\begin{proof}
This is something very standard, as follows:

\medskip

$(1)\implies(2)$ Observe first that $\sigma(a)\subset\mathbb R$ implies $a=a^*$. Thus the algebra $<a>$ is commutative, and by using Theorem 14.26, we can set $b=\sqrt{a}$.

\medskip

$(2)\implies(3)$ This is trivial, because we can simply set $c=b$. 

\medskip

$(2)\implies(1)$ This is clear too, because we have:
$$\sigma(a)
=\sigma(b^2)
=\sigma(b)^2\subset\mathbb R^2
=[0,\infty)$$

$(3)\implies(1)$ We can proceed here by contradiction. Indeed, by multiplying $c$ by a suitable element of the algebra $<cc^*>$, we are led to the existence of an element $d\neq0$ satisfying $-dd^*\geq0$. By writing now $d=x+iy$ with $x=x^*,y=y^*$ we have:
$$dd^*+d^*d
=2(x^2+y^2)
\geq0$$

Thus $d^*d\geq0$, which is easily seen to contradict the condition $-dd^*\geq0$.
\end{proof}

We can talk as well about positive linear forms, in a straightforward way, as follows:

\begin{definition}
Consider a linear map $\varphi:A\to\mathbb C$.
\begin{enumerate}
\item $\varphi$ is called positive when $a\geq0\implies\varphi(a)\geq0$.

\item $\varphi$ is called faithful and positive when $a\geq0,a\neq0\implies\varphi(a)>0$.
\end{enumerate}
\end{definition}

In the commutative case, $A=C(X)$, the positive linear forms appear as follows, with $\mu$ being positive, and strictly positive if we want $\varphi$ to be faithful and positive:
$$\varphi(f)=\int_Xf(x)d\mu(x)$$

In general, the positive linear forms can be thought of as being integration functionals with respect to some underlying ``positive measures''. By using this philosophy, we can establish a key theoretical result, further clarifying the axiomatics, as follows:

\index{representation theorem}
\index{GNS theorem}
\index{operator algebra}

\begin{theorem}
Any $C^*$-algebra appears as an operator algebra:
$$A\subset B(H)$$
Moreover, when $A$ is separable, which is usually the case, $H$ can be taken separable.
\end{theorem}

\begin{proof}
This result, called GNS representation theorem after Gelfand-Naimark-Segal, comes as a continuation of the Gelfand theorem, the idea being as follows:

\medskip

(1) To start with, the result holds in the commutative case, $A=C(X)$. Indeed, we can pick a positive measure on $X$, and construct our embedding as follows:
$$C(X)\subset B(L^2(X))\quad,\quad f\to[g\to fg]$$

(2) In general the proof is similar, the idea being that given a $C^*$-algebra $A$ we can construct a Hilbert space $H=L^2(A)$, and then an embedding as above:
$$A\subset B(L^2(A))\quad,\quad a\to[b\to ab]$$

(3) Finally, the last assertion is clear, because when $A$ is separable, meaning that it has a countable algebraic basis, so does the associated Hilbert space $H=L^2(A)$.
\end{proof}

Still talking general theory of the $C^*$-algebras, we have as well the following result:

\index{finite dimensional algebra}
\index{multimatrix algebra}

\begin{theorem}
The finite dimensional $C^*$-algebras are exactly the algebras
$$A=M_{n_1}(\mathbb C)\oplus\ldots\oplus M_{n_k}(\mathbb C)$$
with norm $||(a_1,\ldots,a_k)||=\sup_i||a_i||$, and involution $(a_1,\ldots,a_k)^*=(a_1^*,\ldots,a_k^*)$.
\end{theorem}

\begin{proof}
In one sense this is clear. In the other sense, this is routine, coming by splitting the unit of our algebra $A\subset B(H)$ as a sum of central minimal projections:
$$1=p_1+\ldots+p_k$$

Indeed, when doing so, each of the $*$-algebras $A_i=p_iAp_i$ follows to be a matrix algebra, $A_i\simeq M_{n_i}(\mathbb C)$, and this gives the direct sum decomposition in the statement.
\end{proof}

Let us record as well the quantum space formulation of the above result:

\index{finite quantum space}

\begin{theorem}
The finite quantum spaces are exactly the disjoint unions of type
$$X=M_{n_1}\sqcup\ldots\sqcup M_{n_k}$$
where $M_n$ is the finite quantum space given by $C(M_n)=M_n(\mathbb C)$.
\end{theorem}

\begin{proof}
For a compact quantum space $X$, coming from a $C^*$-algebra $A$ via the formula $A=C(X)$, being finite can only mean that the following number is finite:
$$|X|=\dim_\mathbb CA<\infty$$

But then, by using Theorem 14.33, we are led to the conclusion that we must have:
$$C(X)=M_{n_1}(\mathbb C)\oplus\ldots\oplus M_{n_k}(\mathbb C)$$

Now since direct sums of algebras $A$ correspond to disjoint unions of quantum spaces $X$, via the correspondence $A=C(X)$, this leads to the conclusion in the statement.
\end{proof}

Getting back now to what we wanted to do, quantum probability, let us formulate:

\begin{definition}
Let $A$ be a $C^*$-algebra, given with a positive trace $tr:A\to\mathbb C$.
\begin{enumerate}
\item The elements $a\in A$ are called random variables.

\item The moments of such a variable are the numbers $M_k(a)=tr(a^k)$.

\item The law of such a variable is the functional $\mu_a:P\to tr(P(a))$.
\end{enumerate}
\end{definition}

To be more precise, here the exponent $k=\circ\bullet\bullet\circ\ldots$ is by definition a colored integer, and the powers $a^k$ are defined by the following formulae, and multiplicativity: 
$$a^\emptyset=1\quad,\quad
a^\circ=a\quad,\quad
a^\bullet=a^*$$ 

As for the polynomial $P$, this is a noncommuting $*$-polynomial in one variable: 
$$P\in\mathbb C<X,X^*>$$

Observe that the law is uniquely determined by the moments, because we have:
$$P(X)=\sum_k\lambda_kX^k
\implies\mu_a(P)=\sum_k\lambda_kM_k(a)$$

At the level of the general theory, we have the following key result, extending the various results from linear algebra, regarding the self-adjoint and normal matrices:

\index{normal element}
\index{spectral measure}

\begin{theorem}
Let $A$ be a $C^*$-algebra, with a trace $tr$, and consider an element $a\in A$ which is normal, in the sense that $aa^*=a^*a$.
\begin{enumerate}
\item $\mu_a$ is a complex probability measure, satisfying $supp(\mu_a)\subset\sigma(a)$.

\item In the self-adjoint case, $a=a^*$, this measure $\mu_a$ is real.

\item Assuming that $tr$ is faithful, we have $supp(\mu_a)=\sigma(a)$.
\end{enumerate}
\end{theorem}

\begin{proof}
In the normal case, $aa^*=a^*a$, the Gelfand theorem, or rather the subsequent continuous functional calculus theorem, tells us that we have: 
$$<a>=C(\sigma(a))$$

Thus the functional $f(a)\to tr(f(a))$ can be regarded as an integration functional on the algebra $C(\sigma(a))$, and by the Riesz theorem, this latter functional must come from a probability measure $\mu$ on the spectrum $\sigma(a)$, in the sense that we must have:
$$tr(f(a))=\int_{\sigma(a)}f(z)d\mu(z)$$

We are therefore led to the conclusions in the statement, with the uniqueness assertion coming from the fact that the elements $a^k$, taken as usual with respect to colored integer exponents, $k=\circ\bullet\bullet\circ\ldots$\,, generate the whole $C^*$-algebra $C(\sigma(a))$.
\end{proof}

As a first concrete application now, getting back to the random matrices, we have:

\begin{theorem}
Given a random matrix $Z\in M_N(L^\infty(X))$ which is normal, 
$$ZZ^*=Z^*Z$$
its law, which is by definition the following abstract functional,
$$\mu:\mathbb C<X,X^*>\to\mathbb C\quad,\quad 
P\to\frac{1}{N}\int_Xtr(P(Z))$$
when restricted to the usual polynomials in two variables,
$$\mu:\mathbb C[X,X^*]\to\mathbb C\quad,\quad 
P\to\frac{1}{N}\int_Xtr(P(Z))$$
must come from a probability measure on the spectrum $\sigma(Z)\subset\mathbb C$, as follows:
$$\mu(P)=\int_{\sigma(T)}P(x)d\mu(x)$$
We agree to use the symbol $\mu$ for all these notions.
\end{theorem}

\begin{proof}
This follows indeed from what we know from Theorem 14.36, applied to the normal element $a=Z$, belonging to the $C^*$-algebra $A=M_N(L^\infty(X))$.
\end{proof}

At the level of the basic examples, the situation is as follows:

\begin{theorem}
The basic random matrices $Z\in M_N(L^\infty(X))$ are as follows:
\begin{enumerate}
\item In the case $N=1$ the random matrix is a usual random variable, $f\in L^\infty(X)$, automatically normal, and its law as defined above is the usual law.

\item In the case $X=\{.\}$ the random matrix is a usual scalar matrix, $A\in M_N(\mathbb C)$, and in the diagonalizable case, the law is $\mu=\frac{1}{N}\left(\delta_{\lambda_1}+\ldots+\delta_{\lambda_N}\right)$.
\end{enumerate}
\end{theorem}

\begin{proof}
This is clear indeed, with the first assertion coming from definitions, and the second assertion coming by diagonalizing the matrix, in the obvious way.
\end{proof}

At a more advanced level now, the main problem regarding the random matrices is that of computing the law of various classes of such matrices, coming in series:

\begin{question}
What is the law of random matrices coming in series
$$Z_N\in M_N(L^\infty(X))$$
in the $N>>0$ regime?
\end{question}

The general strategy here, coming from physicists, is that of computing first the asymptotic law $\mu^0$, in the $N\to\infty$ limit, and then looking for the higher order terms as well, as to finally reach to a series in $N^{-1}$ giving the law of $Z_N$, as follows: 
$$\mu_N=\mu^0+N^{-1}\mu^1+N^{-2}\mu^2+\ldots$$

As a basic example here, of particular interest are the matrices having i.i.d. complex normal entries, under the constraint $Z=Z^*$. Here the asymptotic law $\mu^0$ is the Wigner semicircle law on $[-2,2]$. We will discuss this a bit later, in chapter 16.

\section*{14e. Exercises} 

This was a key chapter, aiming high, and as exercises on this, we have:

\begin{exercise}
Learn the Young inequality, that we used in the above.
\end{exercise}

\begin{exercise}
Learn as well the Cauchy formula, that we used in the above.
\end{exercise}

\begin{exercise}
Clarify the abstract details, for the measurable functional calculus.
\end{exercise}

\begin{exercise}
Learn a bit about compact operators too, and singular values.
\end{exercise}

\begin{exercise}
Learn the Liouville theorem, that we used in the above.
\end{exercise}

\begin{exercise}
And learn as well Hahn-Banach, that we certainly used too.
\end{exercise}

\begin{exercise}
Meditate on the integration problem for the free spheres.
\end{exercise}

\begin{exercise}
Learn about group algebras, and abstract group duals.
\end{exercise}

As bonus exercise, and no surprise here, read an operator algebra book.

\chapter{Partitions, cumulants}

\section*{15a. Cumulants, inversion}

What is next? Time to calm down I guess, after all the things that we learned in the previous two chapters. The operator algebra and random matrix philosophy seems to open up an infinite new world to us, and as a preparation for our exploration, I would say that the wisest would be to try to make a link with the material from Parts I-III.

\bigskip

In answer, looking back at what we did in Parts I-III, our flagship technology there was the moment method. So, the idea will be to refine a bit that material, by using the notion of cumulant of Rota, which is cutting edge weaponry, and then look for further improvements and versions of that combinatorics, applying to the quantum world.

\bigskip

Getting started now, following Rota, let us formulate the following key definition:

\index{cumulant}
\index{cumulant-generating function}
\index{Taylor coefficient}
\index{Fourier transform}
\index{generating series}

\begin{definition}
Associated to any real probability measure $\mu=\mu_f$ is the following modification of the logarithm of the Fourier transform $F_\mu(\xi)=E(e^{i\xi f})$,
$$K_\mu(\xi)=\log E(e^{\xi f})$$
called cumulant-generating function. The Taylor coefficients $k_n(\mu)$ of this series, given by
$$K_\mu(\xi)=\sum_{n=1}^\infty k_n(\mu)\,\frac{\xi^n}{n!}$$
are called cumulants of the measure $\mu$. We also use the notations $k_f,K_f$ for these cumulants and their generating series, where $f$ is a variable following the law $\mu$.
\end{definition}

In other words, the cumulants are more or less the coefficients of the logarithm of the Fourier transform $\log F_\mu$, up to some normalizations. To be more precise, we have $K_\mu(\xi)=\log F_\mu(-i\xi)$, so the formula relating $\log F_\mu$ to the cumulants $k_n(\mu)$ is:
$$\log F_\mu(-i\xi)=\sum_{n=1}^\infty k_n(\mu)\,\frac{\xi^n}{n!}$$

We will see in a moment the reasons for the above normalizations, namely change of variables $\xi\to -i\xi$, and Taylor coefficients instead of plain coefficients, the idea being that for simple laws like $p_t,g_t$, we will obtain in this way very simple quantities.

\bigskip

As a first observation, the sequence of cumulants $k_1,k_2,k_3,\ldots$ appears as a modification of the sequence of moments $M_1,M_2,M_3,\ldots\,$, the numerics being as follows:

\index{sequence of moments}
\index{sequence of cumulants}

\begin{proposition}
The sequence of cumulants $k_1,k_2,k_3,\ldots$ appears as a modification of the sequence of moments $M_1,M_2,M_3,\ldots\,$, and uniquely determines $\mu$. We have
$$k_1=M_1$$
$$k_2=-M_1^2+M_2$$
$$k_3=2M_1^3-3M_1M_2+M_3$$
$$k_4=-6M_1^4+12M_1^2M_2-3M_2^2-4M_1M_3+M_4$$
$$\vdots$$
in one sense, and in the other sense we have
$$M_1=k_1$$
$$M_2=k_1^2+k_2$$
$$M_3=k_1^3+3k_1k_2+k_3$$
$$M_4=k_1^4+6k_1^2k_2+3k_2^2+4k_1k_3+k_4$$
$$\vdots$$
with in both cases the correspondence being polynomial, with integer coefficients.
\end{proposition}

\begin{proof}
We know from Definition 15.1 that the cumulants are given by:
$$\log E(e^{\xi f})=\sum_{s=1}^\infty k_s(f)\,\frac{\xi^s}{s!}$$

By exponentiating, we obtain from this the following formula:
$$E(e^{\xi f})=\exp\left(\sum_{s=1}^\infty k_s(f)\,\frac{\xi^s}{s!}\right)$$

Now by looking at the terms of order $1,2,3,4$, this gives the above formulae.
\end{proof}

The interest in the cumulants comes from the fact that $\log F_\mu$, and so the cumulants $k_n(\mu)$ too, linearize the convolution. To be more precise, we have the following result:

\index{linearization of convolution}
\index{additivity of cumulants}

\begin{theorem}
The cumulants have the following properties:
\begin{enumerate}
\item $k_n(cf)=c^nk_n(f)$.

\item $k_1(f+d)=k_1(f)+d$, and $k_n(f+d)=k_n(f)$ for $n>1$.

\item $k_n(f+g)=k_n(f)+k_n(g)$, if $f,g$ are independent.
\end{enumerate}
\end{theorem}

\begin{proof}
Here (1) and (2) are both clear from definitions, because we have:
\begin{eqnarray*}
K_{cf+d}(\xi)
&=&\log E(e^{\xi(cf+d)})\\
&=&\log[e^{\xi d}\cdot E(e^{\xi cf})]\\
&=&\xi d+K_f(c\xi)
\end{eqnarray*}

As for (3), this follows from the fact that the Fourier transform $F_f(\xi)=E(e^{i\xi f})$ satisfies the following formula, whenever $f,g$ are independent random variables:
$$F_{f+g}(\xi)=F_f(\xi)F_g(\xi)$$

Indeed, by applying the logarithm, we obtain the following formula:
$$\log F_{f+g}(\xi)=\log F_f(\xi)+\log F_g(\xi)$$

With the change of variables $\xi\to-i\xi$, we obtain the following formula:
$$K_{f+g}(\xi)=K_f(\xi)+K_g(\xi)$$

Thus, at the level of coefficients, we obtain $k_n(f+g)=k_n(f)+k_n(g)$, as claimed.
\end{proof}

At the level of basic examples now, we have the following result:

\begin{theorem}
The sequence of cumulants $k_1,k_2,k_3,k_4,\ldots$ is as follows:
\begin{enumerate}
\item For $\mu=\delta_c$ the cumulants are $c,0,0,0,\ldots$

\item For $\mu=g_t$ the cumulants are $0,t,0,0,\ldots$

\item For $\mu=p_t$ the cumulants are $t,t,t,t,\ldots$

\item For $\mu=b_t^2$ the cumulants are $0,t,0,t,\ldots$
\end{enumerate}
Also, for the compound Poisson laws the cumulants are $k_n(p_\nu)=M_n(\nu)$.
\end{theorem}

\begin{proof}
We have 5 computations to be done, the idea being as follows:

\medskip

(1) For a Dirac mass, $\mu=\delta_c$, we have the following computation:
$$K_\mu(\xi)
=\log E(e^{c\xi})
=\log(e^{c\xi})
=c\xi$$

But the plain coefficients of this series are the numbers $c,0,0,\ldots\,$, and so the Taylor coefficients of this series are these same numbers $c,0,0,\ldots\,$, as claimed.

\medskip

(2) For a normal law, $\mu=g_t$, we have the following computation:
$$K_\mu(\xi)
=\log F_\mu(-i\xi)
=\log\exp\left[-t(-i\xi)^2/2\right]
=t\xi^2/2$$

But the plain coefficients of this series are the numbers $0,t/2,0,0,\ldots\,$, and so the Taylor coefficients of this series are the numbers $0,t,0,0,\ldots\,$, as claimed.

\medskip

(3) For a Poisson law, $\mu=p_t$, we have the following computation:
$$K_\mu(\xi)
=\log F_\mu(-i\xi)
=\log\exp\left[(e^{i(-i\xi)}-1)t\right]
=(e^\xi-1)t$$

But the plain coefficients of this series are the numbers $t/n!$, and so the Taylor coefficients of this series are the numbers $t,t,t,\ldots\,$, as claimed.

\medskip

(4) For a Bessel law, $\mu=b_t^2$, we have the following computation:
$$K_\mu(\xi)
=\log F_\mu(-i\xi)
=\log\exp\left[\left(\frac{e^\xi+e^{-\xi}}{2}-1\right)t\right]
=\left(\frac{e^\xi+e^{-\xi}}{2}-1\right)t$$

But the plain coefficients of this series are the numbers $(1+(-1)^n)t/n!$, so the Taylor coefficients of this series are the numbers $0,t,0,t,\ldots\,$, as claimed.

\medskip

(5) In order to prove now the last assertion, which generalizes (3,4), by continuity we can assume that our input measure is discrete, $\nu=\sum_i t_i\delta_{z_i}$ with $t_i>0$ and $z_i\in\mathbb R$. By using the Fourier transform formula for $p_\nu$, from chapter 4, we obtain:
\begin{eqnarray*}
K_{p_\nu}(\xi)
&=&\log F_{p_\nu}(-i\xi)\\
&=&\log\exp\left[\sum_it_i(e^{\xi z_i}-1)\right]\\
&=&\sum_it_i\sum_{n\geq1}\frac{(\xi z_i)^n}{n!}\\
&=&\sum_{n\geq1}\frac{\xi^n}{n!}\sum_it_iz_i^n\\
&=&\sum_{n\geq1}\frac{\xi^n}{n!}\,M_n(\nu)
\end{eqnarray*}

Thus, we are led to the conclusion in the statement.
\end{proof}

Getting back to theory now, the sequence of cumulants $k_1,k_2,k_3,\ldots$ appears as a modification of the sequence of moments $M_1,M_2,M_3,\ldots\,$, and understanding the relation between moments and cumulants will be our next task. Let us start with:

\index{M\"obius function}
\index{lattice}
\index{lattice of partitions}

\begin{definition}
The M\"obius function of any lattice, and so of $P$, is given by
$$\mu(\pi,\nu)=\begin{cases}
1&{\rm if}\ \pi=\nu\\
-\sum_{\pi\leq\tau<\nu}\mu(\pi,\tau)&{\rm if}\ \pi<\nu\\
0&{\rm if}\ \pi\not\leq\nu
\end{cases}$$
with the construction being performed by recurrence.
\end{definition}

As an illustration here, for $P(2)=\{||,\sqcap\}$, we have by definition:
$$\mu(||,||)=\mu(\sqcap,\sqcap)=1$$

Also, $||<\sqcap$, with no intermediate partition in between, so we obtain:
$$\mu(||,\sqcap)=-\mu(||,||)=-1$$

Finally, we have $\sqcap\not\leq||$, and so we have as well the following formula:
$$\mu(\sqcap,||)=0$$

Thus, the M\"obius matrix $M_{\pi\nu}=\mu(\pi,\nu)$ of the lattice $P(2)=\{||,\sqcap\}$ is as follows:
$$M=\begin{pmatrix}1&-1\\ 0&1\end{pmatrix}$$

Back to the general case now, the main interest in the M\"obius function comes from the M\"obius inversion formula, which can be formulated as follows:

\index{M\"obius inversion}
\index{M\"obius matrix}

\begin{theorem}
We have the following implication,
$$f(\pi)=\sum_{\nu\leq\pi}g(\nu)
\quad\implies\quad
g(\pi)=\sum_{\nu\leq\pi}\mu(\nu,\pi)f(\nu)$$
valid for any two functions $f,g:P(n)\to\mathbb C$.
\end{theorem}

\begin{proof}
Consider the adjacency matrix of $P$, given by the following formula:
$$A_{\pi\nu}=\begin{cases}
1&{\rm if}\ \pi\leq\nu\\
0&{\rm if}\ \pi\not\leq\nu
\end{cases}$$

Our claim is that the inverse of this matrix is the M\"obius matrix of $P$, given by:
$$M_{\pi\nu}=\mu(\pi,\nu)$$

Indeed, the above matrix $A$ is upper triangular, and when trying to invert it, we are led to the recurrence in Definition 15.5, so to the M\"obius matrix $M$. Thus we have:
$$M=A^{-1}$$

Thus, in practice, we are led to the inversion formula in the statement.
\end{proof}

As an illustration here, for $P(2)$ the formula $M=A^{-1}$ appears as follows:
$$\begin{pmatrix}1&-1\\ 0&1\end{pmatrix}=
\begin{pmatrix}1&1\\ 0&1\end{pmatrix}^{-1}$$

With these ingredients in hand, let us go back to probability. We first have:

\index{cumulant}
\index{classical cumulant}
\index{generalized cumulant}
\index{multiplicativity over blocks}

\begin{definition}
We define quantities $M_\pi(f),k_\pi(f)$, depending on partitions 
$$\pi\in P(k)$$
by starting with $M_n(f),k_n(f)$, and using multiplicativity over the blocks. 
\end{definition}

To be more precise, the convention here is that for the one-block partition $1_n\in P(n)$, the corresponding moment and cumulant are the usual ones, namely:
$$M_{1_n}(f)=M_n(f)\quad,\quad k_{1_n}(f)=k_n(f)$$

Then, for an arbitrary partition $\pi\in P(k)$, we decompose this partition into blocks, having sizes $b_1,\ldots,b_s$, and we set, by multiplicativity over blocks:
$$M_\pi(f)=M_{b_1}(f)\ldots M_{b_s}(f)\quad,\quad k_\pi(f)=k_{b_1}(f)\ldots k_{b_s}(f)$$

With this convention, following Rota and others, we can now formulate a key result, fully clarifying the relation between moments and cumulants, as follows:

\index{moment-cumulant formula}

\begin{theorem}
We have the moment-cumulant formulae
$$M_n(f)=\sum_{\nu\in P(n)}k_\nu(f)\quad,\quad 
k_n(f)=\sum_{\nu\in P(n)}\mu(\nu,1_n)M_\nu(f)$$
or, equivalently, we have the moment-cumulant formulae
$$M_\pi(f)=\sum_{\nu\leq\pi}k_\nu(f)\quad,\quad 
k_\pi(f)=\sum_{\nu\leq\pi}\mu(\nu,\pi)M_\nu(f)$$
where $\mu$ is the M\"obius function of $P(n)$.
\end{theorem}

\begin{proof}
There are several things going on here, the idea being as follows:

\medskip

(1) According to our conventions above, the first set of formulae is equivalent to the second set of formulae. Also, due to the M\"obius inversion formula, in the second set of formulae, the two formulae there are in fact equivalent. Thus, the 4 formulae in the statement are all equivalent. In what follows we will focus on the first 2 formulae.

\medskip

(2) Let us first work out some examples. At $n=1,2,3$ the moment formula gives the following equalities, which are in tune with the findings from Proposition 15.2:
$$M_1=k_|=k_1$$
$$M_2=k_{|\,|}+k_\sqcap=k_1^2+k_2$$
$$M_3=k_{|\,|\,|}+k_{\sqcap|}+k_{\sqcap\hskip-2.8mm{\ }_|}+k_{|\sqcap}+k_{\sqcap\hskip-0.5mm\sqcap}=k_1^3+3k_1k_2+k_3$$

At $n=4$ now, which is a case which is of particular interest for certain considerations to follow, the computation is as follows, again in tune with Proposition 15.2:
\begin{eqnarray*}
M_4
&=&k_{|\,|\,|}+(\underbrace{k_{\sqcap\,|\,|}+\ldots}_{6\ terms})+(\underbrace{k_{\sqcap\,\sqcap}+\ldots}_{3\ terms})+(\underbrace{k_{\sqcap\hskip-0.5mm\sqcap\,|}+\ldots}_{4\ terms})+k_{\sqcap\hskip-0.5mm\sqcap\hskip-0.5mm\sqcap}\\
&=&k_1^4+6k_1^2k_2+3k_2^2+4k_1k_3+k_4
\end{eqnarray*}

As for the cumulant formula, at $n=1,2,3$ this gives the following formulae for the cumulants, again in tune with the findings from Proposition 15.2:
$$k_1=M_|=M_1$$
$$k_2=(-1)M_{|\,|}+M_\sqcap=-M_1^2+M_2$$
$$k_3=2M_{|\,|\,|}+(-1)M_{\sqcap|}+(-1)M_{\sqcap\hskip-2.8mm{\ }_|}+(-1)M_{|\sqcap}+M_{\sqcap\hskip-0.5mm\sqcap}=2M_1^3-3M_1M_2+M_3$$

Finally, at $n=4$, after computing the M\"obius function of $P(4)$, we obtain the following formula for the fourth cumulant, again in tune with Proposition 15.2:
\begin{eqnarray*}
k_4
&=&(-6)M_{|\,|\,|}+2(\underbrace{M_{\sqcap\,|\,|}+\ldots}_{6\ terms})+(-1)(\underbrace{M_{\sqcap\,\sqcap}+\ldots}_{3\ terms})+(-1)(\underbrace{M_{\sqcap\hskip-0.5mm\sqcap\,|}+\ldots}_{4\ terms})+M_{\sqcap\hskip-0.5mm\sqcap\hskip-0.5mm\sqcap}\\
&=&-6M_1^4+12M_1^2M_2-3M_2^2-4M_1M_3+M_4
\end{eqnarray*}

(3) Time now to get to work, and prove the result. As mentioned above, the formulae in the statement are all equivalent, and it is enough to prove the first one, namely:
$$M_n(f)=\sum_{\nu\in P(n)}k_\nu(f)$$

In order to do this, we use the very definition of the cumulants, namely:
$$\log E(e^{\xi f})=\sum_{s=1}^\infty k_s(f)\,\frac{\xi^s}{s!}$$

By exponentiating, we obtain from this the following formula:
$$E(e^{\xi f})=\exp\left(\sum_{s=1}^\infty k_s(f)\,\frac{\xi^s}{s!}\right)$$

(4) Let us first compute the function on the left. This is easily done, as follows:
$$E(e^{\xi f})
=E\left(\sum_{n=0}^\infty\frac{(\xi f)^n}{n!}\right)
=\sum_{n=0}^\infty M_n(f)\,\frac{\xi^n}{n!}$$

(5) Regarding now the function on the right, this is given by:
\begin{eqnarray*}
\exp\left(\sum_{s=1}^\infty k_s(f)\,\frac{\xi^s}{s!}\right)
&=&\sum_{p=0}^\infty\frac{\left(\sum_{s=1}^\infty k_s(f)\,\frac{\xi^s}{s!}\right)^p}{p!}\\
&=&\sum_{p=0}^\infty\frac{1}{p!}\sum_{s_1=1}^\infty k_{s_1}(f)\,\frac{\xi^{s_1}}{s_1!}\ldots\ldots\sum_{s_p=1}^\infty k_{s_p}(f)\,\frac{\xi^{s_p}}{s_p!}\\
&=&\sum_{p=0}^\infty\frac{1}{p!}\sum_{s_1=1}^\infty\ldots\sum_{s_p=1}^\infty k_{s_1}(f)\ldots k_{s_p}(f)\,\frac{\xi^{s_1+\ldots+s_p}}{s_1!\ldots s_p!}
\end{eqnarray*}

But the point now is that all this leads us into partitions. Indeed, we are summing over indices $s_1,\ldots,s_p\in\mathbb N$, which can be thought of as corresponding to a partition of $n=s_1+\ldots+s_p$. So, let us rewrite our sum, as a sum over partitions. For this purpose, recall that the number of partitions $\nu\in P(n)$ having blocks of sizes $s_1,\ldots,s_p$ is:
$$\binom{n}{s_1,\ldots,s_p}=\frac{n!}{p_1!\ldots p_s!}$$

Also, when resumming over partitions, there will be a $p!$ factor as well, coming from the permutations of $s_1,\ldots,s_p$. Thus, our sum can be rewritten as follows:
\begin{eqnarray*}
\exp\left(\sum_{s=1}^\infty k_s(f)\,\frac{\xi^s}{s!}\right)
&=&\sum_{n=0}^\infty\sum_{p=0}^\infty\frac{1}{p!}\sum_{s_1+\ldots+s_p=n}k_{s_1}(f)\ldots k_{s_p}(f)\,\frac{\xi^n}{s_1!\ldots s_p!}\\
&=&\sum_{n=0}^\infty\frac{\xi^n}{n!}\sum_{p=0}^\infty\frac{1}{p!}\sum_{s_1+\ldots+s_p=n}\binom{n}{s_1,\ldots,s_p}k_{s_1}(f)\ldots k_{s_p}(f)\\
&=&\sum_{n=0}^\infty\frac{\xi^n}{n!}\sum_{\nu\in P(n)}k_\nu(f)
\end{eqnarray*}

(6) We are now in position to conclude. According to (3,4,5), we have:
$$\sum_{n=0}^\infty M_n(f)\,\frac{\xi^n}{n!}=\sum_{n=0}^\infty\frac{\xi^n}{n!}\sum_{\nu\in P(n)}k_\nu(f)$$

Thus, we have the following formula, valid for any $n\in\mathbb N$:
$$M_n(f)=\sum_{\nu\in P(n)}k_\nu(f)$$

We are therefore led to the conclusions in the statement.
\end{proof}

The above cumulant technology is quite powerful, and can be applied as well to more specialized measures that we know. We will leave this as a long, pleasant exercise.

\section*{15b. Free probability} 

As an interesting phenomenon, the cumulant technology does not apply well to the measures from chapter 7, namely Wigner, Marchenko-Pastur, arcsine and modified arcsine, and with the reasons for this being something quite subtle, the idea being that these latter measures are ``free'' in nature, and are therefore in need of ``free cumulants''.

\bigskip

In order to discuss this, following Voiculescu and then Speicher, let us start with something quite general, that we met in chapter 14, namely:

\begin{definition}
Let $A$ be a $C^*$-algebra, coming with a trace $tr:A\to\mathbb C$.
\begin{enumerate}
\item The elements $a\in A$ are called noncommutative random variables.

\item The moments of such a variable are the numbers $M_k(a)=tr(a^k)$.

\item The law of such a variable is the functional $\mu_a:P\to tr(P(a))$.
\end{enumerate}
\end{definition}

Generally speaking, this definition is something quite abstract, but there is no other way of doing things, at least at this level of generality. However, in certain special cases, the formalism simplifies, and we recover more familiar objects, as follows:

\index{spectral theorem}
\index{spectral measure}

\begin{proposition}
Assuming that $a\in A$ is normal, $aa^*=a^*a$, its law corresponds to a probability measure on its spectrum $\sigma(a)\subset\mathbb C$, according to the following formula:
$$tr(P(a))=\int_{\sigma(a)}P(x)d\mu(x)$$
When the trace is faithful we have $supp(\mu)=\sigma(a)$. Also, in the particular case where the variable is self-adjoint, $a=a^*$, this law is a real probability measure.
\end{proposition}

\begin{proof}
This is indeed something that we know from chapter 14, the idea being that by Gelfand we have $<a>=C(\sigma(a))$, and the rest comes from the Riesz theorem.
\end{proof}

Getting now where we wanted to get, following Voiculescu \cite{vdn}, let us formulate:

\index{independence}
\index{freeness}

\begin{definition}
Two subalgebras $A,B\subset C$ are called independent when
$$tr(a)=tr(b)=0\implies tr(ab)=0$$
holds for any $a\in A$ and $b\in B$, and free when
$$tr(a_i)=tr(b_i)=0\implies tr(a_1b_1a_2b_2\ldots)=0$$
holds for any $a_i\in A$ and $b_i\in B$.
\end{definition}
 
In order to understand now what is going on, let us discuss some basic models for independence and freeness. We first have the following result, which clarifies things:

\begin{proposition}
Given two algebras $(A,tr)$ and $(B,tr)$, the following hold:
\begin{enumerate}
\item $A,B$ are independent inside their tensor product $A\otimes B$.

\item $A,B$ are free inside their free product $A*B$.
\end{enumerate}
\end{proposition}

\begin{proof}
Both the assertions are clear from definitions, after some standard discussion regarding the tensor product and free product trace. See Voiculescu \cite{vdn}.
\end{proof}

In relation with groups and algebra, we have the following result, with the group $C^*$-algebras involved being some natural completions of the usual group algebras:

\begin{proposition}
We have the following results, valid for group algebras:
\begin{enumerate}
\item $C^*(\Gamma),C^*(\Lambda)$ are independent inside $C^*(\Gamma\times\Lambda)$.

\item $C^*(\Gamma),C^*(\Lambda)$ are free inside $C^*(\Gamma*\Lambda)$.
\end{enumerate}
\end{proposition}

\begin{proof}
This follows from the general results in Proposition 15.12, along with the following two isomorphisms, which are both standard:
$$C^*(\Gamma\times\Lambda)=C^*(\Lambda)\otimes C^*(\Gamma)\quad,\quad 
C^*(\Gamma*\Lambda)=C^*(\Lambda)*C^*(\Gamma)$$

Alternatively, we can prove this directly, by using the fact that each algebra is spanned by the corresponding group elements, and checking the result on group elements.
\end{proof}

In order to study independence and freeness, the main tools are as follows: 

\index{Fourier transform}
\index{Cauchy transform}
\index{R-transform}
\index{convolution}
\index{free convolution}

\begin{theorem}
The convolution is linearized by the log of the Fourier transform,
$$F_f(x)=E(e^{ixf})$$
and the free convolution is linearized by the Voiculescu $R$-transform, which is given by
$$G_\mu(\xi)=\int_\mathbb R\frac{d\mu(t)}{\xi-t}\implies G_\mu\left(R_\mu(\xi)+\frac{1}{\xi}\right)=\xi$$
and so is the inverse of the Cauchy transform, up to a $\xi^{-1}$ factor.
\end{theorem}

\begin{proof}
The first assertion is something that we know well, since chapter 6. Regarding now the second assertion, we need a good model for free convolution, and the best is to use the semigroup algebra of the free semigroup on two generators:
$$A=C^*(\mathbb N*\mathbb N)$$

Indeed, we have some freeness in the semigroup setting, a bit in the same way as for the group algebras $C^*(\Gamma*\Lambda)$, from Proposition 15.13. In addition to this fact, and to what happens in the group algebra case, the following things happen:

\medskip

(1) The variables of type $S^*+f(S)$, with $S\in C^*(\mathbb N)$ being the shift, and with $f\in\mathbb C[X]$ being a polynomial, model in moments all the distributions $\mu:\mathbb C[X]\to\mathbb C$. This is indeed something elementary, which can be checked via a direct algebraic computation.

\medskip

(2) Given $f,g\in\mathbb C[X]$, the variables $S^*+f(S)$ and $T^*+g(T)$, where $S,T\in C^*(\mathbb N*\mathbb N)$ are the shifts corresponding to the generators of $\mathbb N*\mathbb N$, are free, and their sum has the same law as $S^*+(f+g)(S)$. This follows indeed by using a $45^\circ$ argument.

\medskip

(3) But with this in hand, we can see that $\mu\to f$ linearizes the free convolution. We are therefore left with a computation inside $C^*(\mathbb N)$, whose conclusion is that $R_\mu=f$ can be recaptured from $\mu$ via the Cauchy transform $G_\mu$, as stated. See Voiculescu \cite{vdn}.
\end{proof}

As a main application of the above linearization technology, we have:

\index{CLT}
\index{Wigner law}
\index{semicircle law}
\index{normal law}
\index{Gaussian law}

\begin{theorem}[CLT]
Given self-adjoint random variables $x_1,x_2,x_3,\ldots$ assumed to be i.i.d./f.i.d., centered, with variance $t>0$, we have, with $n\to\infty$, in moments,
$$\frac{1}{\sqrt{n}}\sum_{i=1}^nx_i\sim g_t/\gamma_t$$
where the limiting laws $g_t/\gamma_t$ are the following measures,
$$g_t=\frac{1}{\sqrt{2\pi t}}e^{-x^2/2t}dx\quad,\quad 
\gamma_t=\frac{1}{2\pi t}\sqrt{4t^2-x^2}dx$$
which are the normal, or Gaussian, and Wigner semicircle law of parameter $t$.
\end{theorem}

\begin{proof}
In the classical case, we know this from chapter 10. In the free case the proof is similar, based on the linearization technology from Theorem 15.14, as follows:

\medskip

(1) We can assume that we are dealing with the case $t=1$, the proof in general being similar. To start with, the $R$-transform of the variable on the left is given by:
$$R(\xi)
=nR_x\left(\frac{\xi}{\sqrt{n}}\right)
\simeq\xi$$

(2) Regarding now the right term, which is the Wigner semicircle law $\gamma_1$, our first claim is that the Cauchy transform of this law satisfies the following equation:
$$G_{\gamma_1}\left(\xi+\frac{1}{\xi}\right)=\xi$$

Indeed, we know from chapter 7 that the Cauchy transform of $\gamma_1$ is given by:
\begin{eqnarray*}
G_{\gamma_1}(y)
&=&y^{-1}\sum_{k=0}^\infty C_ky^{-2k}\\
&=&y^{-1}\cdot\frac{1-\sqrt{1-4y^{-2}}}{2y^{-2}}\\
&=&\frac{y}{2}\left(1-\sqrt{1-4y^{-2}}\right)\\
&=&\frac{y}{2}-\frac{1}{2}\sqrt{y^2-4}
\end{eqnarray*}

Now with $y=\xi+\xi^{-1}$, this formula becomes, as claimed in the above:
\begin{eqnarray*}
G_{\gamma_1}\left(\xi+\frac{1}{\xi}\right)
&=&\frac{\xi+\xi^{-1}}{2}-\frac{1}{2}\sqrt{\xi^2+\xi^{-2}-2}\\
&=&\frac{\xi+\xi^{-1}}{2}-\frac{\xi^{-1}-\xi}{2}\\
&=&\xi
\end{eqnarray*}

(3) We conclude from the formula found in (2) and from Theorem 15.14 that the $R$-transform of the Wigner semicircle law $\gamma_1$ is given by the following formula:
$$R_{\gamma_1}(\xi)=\xi$$

In fact, this follows as well from the following formula, coming from the random walk considerations from chapter 7, and from the technical details of the $R$-transform:
$$S+S^*\sim\gamma_1$$

Thus, the laws in the statement have the same $R$-transforms, so they are equal.
\end{proof}

Next, we have the following complex version of the CLT:

\index{CCLT}
\index{Complex CLT}
\index{Voiculescu law}
\index{circular law law}
\index{complex normal law}
\index{complex Gaussian law}

\begin{theorem}[CCLT]
Given variables $x_1,x_2,x_3,\ldots$ which are i.i.d./f.i.d., centered, with variance $t>0$, we have, with $n\to\infty$, in moments,
$$\frac{1}{\sqrt{n}}\sum_{i=1}^nx_i\sim G_t/\Gamma_t$$
where $G_t/\Gamma_t$ are the complex normal and Voiculescu circular law of parameter $t$, given by:
$$G_t=law\left(\frac{1}{\sqrt{2}}(a+ib)\right)\quad,\quad 
\Gamma_t=law\left(\frac{1}{\sqrt{2}}(\alpha+i\beta)\right)$$
where $a,b/\alpha,\beta$ are independent/free, each following the law $g_t/\gamma_t$.
\end{theorem}

\begin{proof}
This follows indeed from the CLT, by taking real and imaginary parts of all the variables involved, and for details and more here, including the combinatorics of the Voiculescu circular law $\Gamma_t$, which is quite subtle, we refer again to \cite{vdn}.
\end{proof}

Let us denote by $\boxplus$ the free convolution of real probability measures, given by the fact that $\mu_{a+b}=\mu_a\boxplus\mu_b$, when $a,b$ are free. With this convention, we have the following discrete version of the CLT, called Poisson Limit Theorem (PLT):

\index{PLT}
\index{Poisson limit}
\index{Poisson law}
\index{free Poisson law}
\index{Marchenko-Pastur law}

\begin{theorem}[PLT]
The following Poisson limits converge, for any $t>0$,
$$p_t=\lim_{n\to\infty}\left(\left(1-\frac{t}{n}\right)\delta_0+\frac{t}{n}\delta_1\right)^{*n}\quad,\quad 
\pi_t=\lim_{n\to\infty}\left(\left(1-\frac{t}{n}\right)\delta_0+\frac{t}{n}\delta_1\right)^{\boxplus n}$$
the limiting measures being the Poisson law $p_t$, and the Marchenko-Pastur law $\pi_t$,
$$p_t=\frac{1}{e^t}\sum_{k=0}^\infty\frac{t^k\delta_k}{k!}\quad,\quad 
\pi_t=\max(1-t,0)\delta_0+\frac{\sqrt{4t-(x-1-t)^2}}{2\pi x}\,dx$$
with the Marchenko-Pastur law being also called free Poisson law. 
\end{theorem}

\begin{proof}
In the classical case, we know this from chapter 4. In the free case the proof is similar, based on the linearization technology from Theorem 15.14, as follows:

\medskip

(1) Consider the measure in the statement, under the convolution sign:
$$\eta=\left(1-\frac{t}{n}\right)\delta_0+\frac{t}{n}\delta_1$$

The Cauchy transform of this measure is easy to compute, and is given by:
$$G_\eta(\xi)=\left(1-\frac{t}{n}\right)\frac{1}{\xi}+\frac{t}{n}\cdot\frac{1}{\xi-1}$$

(2) In order to prove the result, we want to compute the following $R$-transform:
$$R
=R_{\eta^{\boxplus n}}(y)
=nR_\eta(y)$$

According to the formula of $G_\eta$, the equation for this function $R$ is as follows:
$$\left(1-\frac{t}{n}\right)\frac{1}{1/y+R/n}+\frac{t}{n}\cdot\frac{1}{1/y+R/n-1}=y$$

(3) By multiplying both sides by $n/y$, this equation can be written as:
$$\frac{t+yR}{1+yR/n}=\frac{t}{1+yR/n-y}$$

With $n\to\infty$ things simplify, and we obtain the following formula:
$$t+yR=\frac{t}{1-y}$$

Thus we have the following formula, for the $R$-transform that we are interested in:
$$R=\frac{t}{1-y}$$

(4) But this gives the result, since $R_{\pi_t}$ is elementary to compute from what we have, by ``doubling'' the results for the Wigner law $\gamma_t$, and is given by the same formula.
\end{proof}

Finally, we have the following ``compound'' generalization of the PLT:

\index{CPLT}
\index{Compound PLT}
\index{compound Poisson law}

\begin{theorem}[CPLT]
Given a compactly supported positive measure $\rho$, of mass $c=mass(\rho)$, the following compound Poisson limits converge,
$$p_\rho=\lim_{n\to\infty}\left(\left(1-\frac{c}{n}\right)\delta_0+\frac{1}{n}\rho\right)^{*n}
\quad,\quad 
\pi_\rho=\lim_{n\to\infty}\left(\left(1-\frac{c}{n}\right)\delta_0+\frac{1}{n}\rho\right)^{\boxplus n}$$
and if we write $\rho=\sum_{i=1}^sc_i\delta_{z_i}$ with $c_i>0$ and $z_i\in\mathbb R$ we have the formula
$$p_\rho/\pi_\rho={\rm law}\left(\sum_{i=1}^sz_i\alpha_i\right)$$
where the variables $\alpha_i$ are Poisson/free Poisson$(c_i)$, independent/free.
\end{theorem}

\begin{proof}
In the classical case, we know this from chapter 4. In the free case the proof is similar, based on the linearization technology from Theorem 15.14, as follows:

\medskip

(1) Consider the measure in the statement, under the convolution sign:
$$\eta_n=\left(1-\frac{c}{n}\right)\delta_0+\frac{1}{n}\rho$$

The Cauchy transform of $\eta_n$ is then given by the following formula:
$$G_{\eta_n}(\xi)=\left(1-\frac{c}{n}\right)\frac{1}{\xi}+\frac{1}{n}\sum_{i=1}^s\frac{c_i}{\xi-z_i}$$

(2) Consider now the $R$-transform of the measure $\eta_n^{\boxplus n}$, which is given by:
$$R_{\eta_n^{\boxplus n}}(y)=nR_{\eta_n}(y)$$

By using the general theory of the $R$-transform, from Theorem 15.14, the above formula of $G_{\eta_n}$ shows that the equation for $R=R_{\eta_n^{\boxplus n}}$ is as follows:
\begin{eqnarray*}
&&\left(1-\frac{c}{n}\right)\frac{1}{1/y+R/n}+\frac{1}{n}\sum_{i=1}^s\frac{c_i}{1/y+R/n-z_i}=y\\
&\implies&\left(1-\frac{c}{n}\right)\frac{1}{1+yR/n}+\frac{1}{n}\sum_{i=1}^s\frac{c_i}{1+yR/n-yz_i}=1
\end{eqnarray*}

(3) Now multiplying by $n$, then rearranging the terms, and letting $n\to\infty$, we get:
\begin{eqnarray*}
\frac{c+yR}{1+yR/n}=\sum_{i=1}^s\frac{c_i}{1+yR/n-yz_i}
&\implies&c+yR_{\pi_\rho}(y)=\sum_{i=1}^s\frac{c_i}{1-yz_i}\\
&\implies&R_{\pi_\rho}(y)=\sum_{i=1}^s\frac{c_iz_i}{1-yz_i}
\end{eqnarray*}

Thus, we are led to the conclusion in the statement.
\end{proof}

So long for limiting results in classical and free probability. To finish with, for our purposes here, we will need the following consequence of Theorem 15.18:

\index{Bessel law}
\index{free Bessel law}

\begin{theorem}
Consider the Bessel and free Bessel laws, depending on parameters $s\in\mathbb N\cup\{\infty\}$ and $t>0$, which are the following compound Poisson/free Poisson laws,
$$b^s_t=p_{t\varepsilon_s}\quad,\quad 
\beta^s_t=\pi_{t\varepsilon_s}$$
with $\varepsilon_s$ being the uniform measure on the $s$-th roots of unity. Then:
\begin{enumerate}
\item At $s=1$ we recover the Poisson laws $p_t,\pi_t$.

\item At $s=2$ we recover the real Bessel laws $b_t,\beta_t$.

\item At $s=\infty$ we have the complex Bessel laws $B_t,\mathfrak B_t$.
\end{enumerate}
\end{theorem}

\begin{proof}
This is something quite self-explanatory, based on Theorem 15.18, and on the material from chapter 4 too. We refer to \cite{bb+} for more on these laws. 
\end{proof}

\section*{15c. Free cumulants}

Getting now to combinatorial aspects, following Speicher \cite{nsp}, let us formulate:

\index{cumulant}
\index{free cumulant}

\begin{definition}
The free cumulants $\kappa_n(a)$ of a variable $a\in A$ are defined by
$$R_a(\xi)=\sum_{n=1}^\infty\kappa_n(a)\xi^{n-1}$$
with $R$ being as usual the Voiculescu $R$-transform.
\end{definition}

As before with the classical cumulants, we have a number of basic examples and illustrations, and a number of basic general results. Let us start with some numerics:

\begin{proposition}
The free cumulants $\kappa_1,\kappa_2,\kappa_3,\ldots$ appear as a modification of the moments $M_1,M_2,M_3,\ldots\,$, and uniquely determine $\mu$. We have
$$\kappa_1=M_1$$
$$\kappa_2=-M_1^2+M_2$$
$$\kappa_3=2M_1^3-3M_1M_2+M_3$$
$$\kappa_4=-5M_1^4+10M_1^2M_2-2M_2^2-4M_1M_3+M_4$$
$$\vdots$$
in one sense, and in the other sense we have
$$M_1=\kappa_1$$
$$M_2=\kappa_1^2+\kappa_2$$
$$M_3=\kappa_1^3+3\kappa_1\kappa_2+\kappa_3$$
$$M_4=\kappa_1^4+6\kappa_1^2\kappa_2+2\kappa_2^2+4\kappa_1\kappa_3+\kappa_4$$
$$\vdots$$
with in both cases the correspondence being polynomial, with integer coefficients.
\end{proposition}

\begin{proof}
This is something quite instructive, a bit more complicated than the similar computations in the classical case, worth discussing in detail, the idea being as follows:

\medskip

(1) We know that the Cauchy transform is the following function:
$$G(\xi)=\sum_{n=0}^\infty\frac{M_n}{\xi^{n+1}}$$

Consider the inverse of this Cauchy transform $G$, with respect to composition:
$$G(K(\xi))=K(G(\xi))=\xi$$

According to Definition 15.20, the free cumulants $\kappa_n$ appear then as follows:
$$K(\xi)=\frac{1}{\xi}+\sum_{n=1}^\infty\kappa_n\xi^{n-1}$$

Thus, we can compute moments in terms of free cumulants, and vice versa, by using either of the inversion formulae $G(K(\xi))=\xi$ and $K(G(\xi))=\xi$.

\medskip

(2) This was for the theory. In practice now, playing with the original inversion formula from Theorem 15.14, namely $G(K(\xi))=\xi$, proves to be something quite complicated, so we will choose to use instead the other inversion formula, namely:
$$K(G(\xi))=\xi$$

Thus, the equation that we want to use is as follows, with $G=G(\xi)$:
$$\frac{1}{G}+\sum_{n=1}^\infty\kappa_nG^{n-1}=\xi$$

(3) With $\xi=z^{-1}$ our equation takes the following form, with $G=G(z^{-1})$:
$$\frac{1}{G}+\sum_{n=1}^\infty\kappa_nG^{n-1}=z^{-1}$$

Now by multiplying by $z$, our equation takes the following form:
$$\frac{z}{G}+z\sum_{n=1}^\infty\kappa_nG^{n-1}=1$$

Equivalently, our equation is as follows, with $G=G(z^{-1})$ as before:
$$\frac{z}{G}+\sum_{n=1}^\infty\kappa_nz^n\left(\frac{G}{z}\right)^{n-1}=1$$

(4) Observe now that we have the following formula:
$$\frac{G}{z}=\frac{G(z^{-1})}{z}=\frac{\sum_{n=0}^\infty M_nz^{n+1}}{z}=\sum_{n=0}^\infty M_nz^n$$

This suggests introducing the following quantity:
$$F=\sum_{n=1}^\infty M_nz^n$$

Indeed, we have then $G/z=1+F$, and our equation becomes:
$$\frac{1}{1+F}+\sum_{n=1}^\infty\kappa_nz_n(1+F)^{n-1}=1$$

(5) By expanding the fraction on the left, our equation becomes:
$$\sum_{n=0}^\infty(-F)^n+\sum_{n=1}^\infty\kappa_nz_n(1+F)^{n-1}=1$$

Moreover, we can cancel the 1 term on both sides, and our equation becomes:
$$\sum_{n=1}^\infty(-F)^n+\sum_{n=1}^\infty\kappa_nz_n(1+F)^{n-1}=0$$

Alternatively, we can write our equation as follows:
$$\sum_{n=1}^\infty\kappa_nz_n(1+F)^{n-1}=-\sum_{n=1}^\infty(-F)^n$$

(6) Good news, this latter equation is something that we are eventually happy with. By remembering that we have $F=\sum_{n=1}^\infty M_nz^n$, our equation looks as follows:
\begin{eqnarray*}
&&\kappa_1z+\kappa_2z^2(1+M_1z+M_2z^2+\ldots)+\kappa_3z^3(1+M_1z+M_2z^2+\ldots)^2+\ldots\\
&=&(M_1z+M_2z^2+\ldots)-(M_1z+M_2z^2+\ldots)^2+(M_1z+M_2z^2+\ldots)^3-\ldots
\end{eqnarray*}

(7) This was for the hard part, carefully fine-tuning our equation, as to have it as simple as possible, before getting to numeric work. The rest is routine. Indeed, by looking at the terms of order $1,2,3,4$ we obtain, instantly or almost, the formulae of $\kappa_1,\kappa_2,\kappa_3,\kappa_4$ in the statement. As for the formulae for $M_1,M_2,M_3,M_4$, these follow from these. 

\medskip

(8) To be more precise, the equations that we get at order $1,2,3,4$ are as follows:
$$\kappa_1=M_1$$
$$\kappa_2=M_2-M_1^2$$
$$\kappa_2M_1+\kappa_3=M_3-2M_1M_2+M_1^3$$
$$\kappa_4+2\kappa_3M_1+\kappa_2M_2=M_4-2M_1M_3-M_2^2+3M_1^2M_2-M_1^4$$

Thus, we are led to the formulae of $\kappa_1,\kappa_2,\kappa_3,\kappa_4$ in the statement, and then to the formulae of $M_1,M_2,M_3,M_4$ in the statement, as desired. 
\end{proof}

Observe the similarity with the formulae for classical cumulants. In fact, we have:

\begin{conclusion}
The first three classical and free cumulants coincide, but
$$k_4=-6M_1^4+12M_1^2M_2-3M_2^2-4M_1M_3+M_4$$
$$\kappa_4=-5M_1^4+10M_1^2M_2-2M_2^2-4M_1M_3+M_4$$
are different, and the same happens at higher order as well.
\end{conclusion}

This is something quite interesting, and we will back later with a conceptual explanation for this, via partitions, the idea being that all this comes from:
$$P(n)=NC(n)\iff n\leq 3$$

But more on this later. At the level of basic general results, we first have:

\begin{theorem}
The free cumulants have the following properties:
\begin{enumerate}
\item $\kappa_n(\lambda a)=\lambda^n\kappa_n(a)$.

\item $\kappa_n(a+b)=\kappa_n(a)+\kappa_n(b)$, if $a,b$ are free.
\end{enumerate}
\end{theorem}

\begin{proof}
This is something very standard, the idea being as follows:

\medskip

(1) We have the following Cauchy transform computation:
\begin{eqnarray*}
G_{\lambda a}(\xi)
&=&\int_\mathbb R\frac{d\mu_{\lambda a}(t)}{\xi-t}\\
&=&\int_\mathbb R\frac{d\mu_a(s)}{\xi-\lambda s}\\
&=&\frac{1}{\lambda}\int_\mathbb R\frac{d\mu_a(s)}{\xi/\lambda-s}\\
&=&\frac{1}{\lambda}\,G_a\left(\frac{\xi}{\lambda}\right)
\end{eqnarray*}

But this gives the following formula, by using the definition of the $R$-transform:
\begin{eqnarray*}
G_{\lambda a}\left(\lambda R_a(\lambda\xi)+\frac{1}{\xi}\right)
&=&\frac{1}{\lambda}\,G_a\left(R_a(\lambda\xi)+\frac{1}{\lambda\xi}\right)\\
&=&\frac{1}{\lambda}\cdot\lambda\xi\\
&=&\xi
\end{eqnarray*}

Thus we have the formula $R_{\lambda a}(\xi)=\lambda R_a(\lambda\xi)$, which gives (1). 

\medskip

(2) This follows from the standard fact, that we know well from before, that the $R$-transform linearizes the free convolution operation.
\end{proof}

Again in analogy with the classical case, at the level of examples, we have:

\begin{theorem}
The sequence of free cumulants $\kappa_1,\kappa_2,\kappa_3,\kappa_4,\ldots$ is as follows:
\begin{enumerate}
\item For $\mu=\delta_c$ the free cumulants are $c,0,0,0,\ldots$

\item For $\mu=\gamma_t$ the free cumulants are $0,t,0,0,\ldots$

\item For $\mu=\pi_t$ the free cumulants are $t,t,t,t,\ldots$

\item For $\mu=\beta_t^2$ the free cumulants are $0,t,0,t,\ldots$
\end{enumerate}
Also, for the compound free Poisson laws the free cumulants are $\kappa_n(\pi_\nu)=M_n(\nu)$.
\end{theorem}

\begin{proof}
The proofs are analogous to those from the classical case, as follows:

\medskip

(1) For $\mu=\delta_c$ we have $G_\mu(\xi)=1/(\xi-c)$, and so $R_\mu(\xi)=c$, as desired.

\medskip

(2) For $\mu=\gamma_t$ we have, as computed before, $R_\mu(\xi)=t\xi$, as desired.

\medskip

(3) For $\mu=\pi_t$ we have, also from before, $R_\mu(\xi)=t/(1-\xi)$, as desired.

\medskip

(4) For $\mu=\beta_t^2$, a free Bessel law, this can be established directly, but the best is to prove directly the last assertion, which generalizes (3,4). With $\nu=\sum_ic_i\delta_{z_i}$ we have:
\begin{eqnarray*}
R_{\pi_\nu}(\xi)
&=&\sum_i\frac{c_iz_i}{1-\xi z_i}\\
&=&\sum_{n\geq1}\xi^{n-1}\sum_ic_iz_i^n\\
&=&\sum_{n\geq 1}\xi^{n-1}\,M_n(\nu)
\end{eqnarray*}

Thus, we are led to the conclusion in the statement.
\end{proof}

As before in the classical case, we can define now generalized free cumulants, $\kappa_\pi(a)$ with $\pi\in P(k)$, by starting with the numeric free cumulants $\kappa_n(a)$, as follows:

\begin{definition}
We define free cumulants $\kappa_\pi(a)$, depending on partitions 
$$\pi\in P(k)$$
by starting with $\kappa_n(a)$, and using multiplicativity over the blocks. 
\end{definition}

To be more precise, the convention here is that for the one-block partition $1_n\in P(n)$, the corresponding free cumulant is the usual one, namely:
$$\kappa_{1_n}(a)=\kappa_n(a)$$

Then, for an arbitrary partition $\pi\in P(k)$, we decompose this partition into blocks, having sizes $b_1,\ldots,b_s$, and we set, by multiplicativity over blocks:
$$\kappa_\pi(a)=\kappa_{b_1}(a)\ldots\kappa_{b_s}(a)$$

With this convention, we have the following result, due to Speicher \cite{nsp}:

\begin{theorem}
We have the moment-cumulant formulae
$$M_n(a)=\sum_{\nu\in NC(n)}\kappa_\nu(a)\quad,\quad 
\kappa_n(a)=\sum_{\nu\in NC(n)}\mu(\nu,1_n)M_\nu(a)$$
or, equivalently, we have the moment-cumulant formulae
$$M_\pi(a)=\sum_{\nu\leq\pi}\kappa_\nu(a)\quad,\quad 
\kappa_\pi(a)=\sum_{\nu\leq\pi}\mu(\nu,\pi)M_\nu(a)$$
where $\mu$ is the M\"obius function of $NC(n)$.
\end{theorem}

\begin{proof}
As before in the classical case, the 4 formulae in the statement are equivalent, via M\"obius inversion. Thus, it is enough to prove one of them, and we will prove the first formula, which in practice is the most useful one. Thus, we must prove that:
$$M_n(a)=\sum_{\nu\in NC(n)}\kappa_\nu(a)$$

(1) In order to prove this formula, let us get back to the construction of the free cumulants, from Definition 15.20. The Cauchy transform of $a$ is the following function:
$$G_a(\xi)=\sum_{n=0}^\infty\frac{M_n(a)}{\xi^{n+1}}$$

Consider the inverse of this Cauchy transform $G_a$, with respect to composition:
$$G_a(K_a(\xi))=K_a(G_a(\xi))=\xi$$

According to Definition 15.20, the free cumulants $\kappa_n(a)$ appear then as follows:
$$K_a(\xi)=\frac{1}{\xi}+\sum_{n=1}^\infty\kappa_n(a)\xi^{n-1}$$

Thus, we can compute moments in terms of free cumulants by using either of the inversion formulae $G_a(K_a(\xi))=\xi$ and $K_a(G_a(\xi))=\xi$. 

\medskip

(2) In practice, as seen in the proof of Proposition 15.21, the best is to use the second inversion formula, $K_a(G_a(\xi))=\xi$, which after some manipulations reads:
\begin{eqnarray*}
&&\kappa_1z+\kappa_2z^2(1+M_1z+M_2z^2+\ldots)+\kappa_3z^3(1+M_1z+M_2z^2+\ldots)^2+\ldots\\
&=&(M_1z+M_2z^2+\ldots)-(M_1z+M_2z^2+\ldots)^2+(M_1z+M_2z^2+\ldots)^3-\ldots
\end{eqnarray*}

But, in case you have fully understood the proof of Proposition 15.21, you know how to exploit this formula at order $n=1,2,3,4$. The same method works in general, and after some computations, this leads to the formula that we want to establish, namely:
$$M_n(a)=\sum_{\nu\in NC(n)}\kappa_\nu(a)$$

(3) We are therefore led to the conclusions in the statement. All this was of course quite brief, and for details here, we refer for instance to Nica-Speicher \cite{nsp}.
\end{proof}

Observe that the above result leads among others to a more conceptual explanation for Conclusion 15.22, with the equalities and non-equalities there coming from:
$$P(n)=NC(n)\iff n\leq3$$

Observe that this is related as well to our kurtosis comments from chapter 10. Nice.

\section*{15d. The standard cube}

As a main application now of the theories of Rota cumulants and Speicher free cumulants, following Bercovici and Pata \cite{bpa}, we can formulate the following simple and bright definition, making the connection between classical and free:

\begin{definition}
We say that a real probability measure
$$m\in\mathcal P(\mathbb R)$$
is the classical version of another measure, called its free version, or its liberation
$$\mu\in\mathcal P(\mathbb R)$$
when the classical cumulants of $m$ coincide with the free cumulants of $\mu$.
\end{definition}

As a first observation, this definition fits with all the classical and free probability theory developed in the above, and notably with the various limiting measures, coming from the main classical and free versions of the CLT, which form a cube, as follows:

\begin{theorem}
In the standard cube of basic probability limiting measures
$$
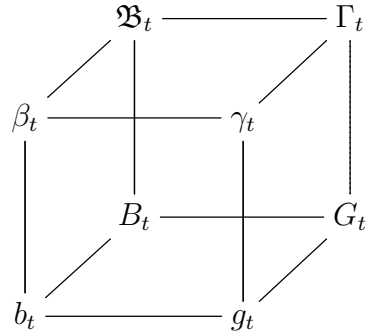
\xymatrix@R=20pt@C=22pt{
&\mathfrak B_t\ar@{-}[rr]\ar@{-}[dd]&&\Gamma_t\ar@{-}[dd]\\
\beta_t\ar@{-}[rr]\ar@{-}[dd]\ar@{-}[ur]&&\gamma_t\ar@{-}[dd]\ar@{-}[ur]\\
&B_t\ar@{-}[rr]\ar@{-}[uu]&&G_t\ar@{.}[uu]\\
b_t\ar@{-}[uu]\ar@{-}[ur]\ar@{-}[rr]&&g_t\ar@{-}[uu]\ar@{-}[ur]
}$$
the upper measures appear as the free versions of the lower measures.
\end{theorem}

\begin{proof}
This follows indeed from our various cumulant formulae found above. Observe that the above cube is in fact something quite subtle, as follows:

\medskip

(1) On the vertical, as stated, we have the classical/free correspondence.

\medskip

(2) On the horizontal we have an interesting discrete/continuous correspondence.

\medskip

(3) And on the depth, we have a real/complex correspondence, interesting as well.
\end{proof}

In order to reach now to an even more advanced theory, depending this time on a parameter $t>0$, still following Bercovici-Pata \cite{bpa}, let us formulate:

\index{convolution semigroup}
\index{free convolution semigroup}
\index{Bercovici-Pata bijection}
\index{classical cumulants}
\index{free cumulants}

\begin{definition}
A convolution semigroup of measures
$$\{m_t\}_{t>0}\quad:\quad m_s*m_t=m_{s+t}$$
is in Bercovici-Pata bijection with a free convolution semigroup of measures
$$\{\mu_t\}_{t>0}\quad:\quad \mu_s\boxplus\mu_t=\mu_{s+t}$$
when the classical cumulants of $m_t$ coincide with the free cumulants of $\mu_t$.
\end{definition}

As before, this fits with all the theory developed so far in this book, and notably with the measures from the standard cube, and we have the following result:

\begin{theorem}
In the standard cube of basic semigroups of measures,
$$
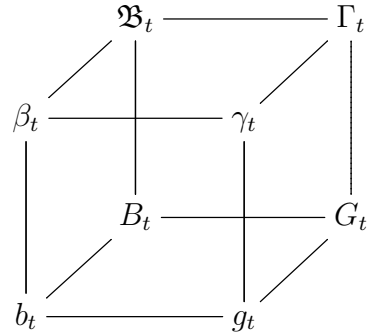
\xymatrix@R=20pt@C=22pt{
&\mathfrak B_t\ar@{-}[rr]\ar@{-}[dd]&&\Gamma_t\ar@{-}[dd]\\
\beta_t\ar@{-}[rr]\ar@{-}[dd]\ar@{-}[ur]&&\gamma_t\ar@{-}[dd]\ar@{-}[ur]\\
&B_t\ar@{-}[rr]\ar@{-}[uu]&&G_t\ar@{.}[uu]\\
b_t\ar@{-}[uu]\ar@{-}[ur]\ar@{-}[rr]&&g_t\ar@{-}[uu]\ar@{-}[ur]
}$$
the upper semigroups are in Bercovici-Pata bijection with the lower semigroups.
\end{theorem}

\begin{proof}
This is a technical improvement of Theorem 15.28, based on the fact that the upper measures in the above diagram form indeed free convolution semigroups, and that the lower measures form indeed classical convolution semigroups, which itself is something very standard, coming from the various results established in the above.
\end{proof}

And with this, end of the story? Not at all, the above standard cube, providing us with some very useful 3D orientation, inside probability theory at large, as briefly explained in the proof of Theorem 15.28, being rather the beginning of something new.

\bigskip

So, let us draw more cubes. To start with, we have the following result:

\begin{theorem}
The moments of the basic limiting laws are given by
$$M_k=\sum_{\pi\in D(k)}t^{|\pi|}$$
with the corresponding categories of partitions $D\subset P$ being as follows,
$$
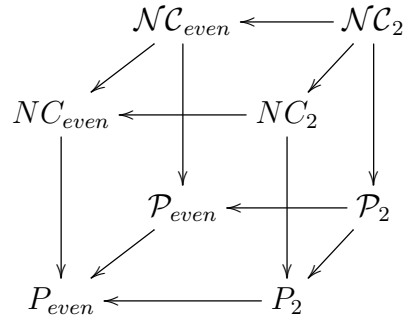
\xymatrix@R=18pt@C=2pt{
&\mathcal{NC}_{even}\ar[dl]\ar[dd]&&\mathcal {NC}_2\ar[dl]\ar[ll]\ar[dd]\\
NC_{even}\ar[dd]&&NC_2\ar[dd]\ar[ll]\\
&\mathcal P_{even}\ar[dl]&&\mathcal P_2\ar[dl]\ar[ll]\\
P_{even}&&P_2\ar[ll]
}$$
and with $|.|$ standing as usual for the number of blocks.
\end{theorem}

\begin{proof}
This is indeed something quite self-explanatory, based on the various moment and cumulant formulae that we found, in this chapter, and previously in this book.
\end{proof}

At a more advanced level now, still talking cubes, we have the following result:

\index{free orthogonal group}
\index{free unitary group}
\index{free rotation}
\index{free reflection group}
\index{Brauer theorem}
\index{easiness}
\index{noncrossing partitions}
\index{quantum reflection group}
\index{standard cube}
\index{liberation}

\begin{theorem}
We have compact groups and quantum groups as follows,
$$\xymatrix@R=16pt@C=16pt{
&K_N^+\ar[rr]&&U_N^+\\
H_N^+\ar[rr]\ar[ur]&&O_N^+\ar[ur]\\
&K_N\ar[rr]\ar[uu]&&U_N\ar[uu]\\
H_N\ar[uu]\ar[ur]\ar[rr]&&O_N\ar[uu]\ar[ur]
}$$
which correspond via Tannakian duality to the following categories of partitions,
$$\xymatrix@R=18pt@C=2pt{
&\mathcal{NC}_{even}\ar[dl]\ar[dd]&&\mathcal {NC}_2\ar[dl]\ar[ll]\ar[dd]\\
NC_{even}\ar[dd]&&NC_2\ar[dd]\ar[ll]\\
&\mathcal P_{even}\ar[dl]&&\mathcal P_2\ar[dl]\ar[ll]\\
P_{even}&&P_2\ar[ll]
}$$
and whose asymptotic laws of truncated characters are the main limiting laws:
$$\xymatrix@R=18pt@C=20pt{
&\mathfrak B_t\ar@{-}[rr]\ar@{-}[dd]&&\Gamma_t\ar@{-}[dd]\\
\beta_t\ar@{-}[rr]\ar@{-}[dd]\ar@{-}[ur]&&\gamma_t\ar@{-}[dd]\ar@{-}[ur]\\
&B_t\ar@{-}[rr]\ar@{-}[uu]&&G_t\ar@{.}[uu]\\
b_t\ar@{-}[uu]\ar@{-}[ur]\ar@{-}[rr]&&g_t\ar@{-}[uu]\ar@{-}[ur]
}$$
Moreover, under very strong axioms for everything, we have uniqueness too.
\end{theorem}

\begin{proof}
There is a long story here, first for formulating the precise statement, which is something non-trivial, that we will not attempt to explain here, and then of course for proving it, and for the whole story here, we refer to \cite{bsp} and related papers.
\end{proof}

Summarizing, standard cube, in its various incarnations, well understood.

\bigskip

Back to the Bercovici-Pata bijection now, as formulated before, there are of course many other instances of it, and with some escaping from our standard cube philosophy. As a main example here, we have the following surprising result, from \cite{bbl}:

\index{normal law}

\begin{theorem}
The normal law $g_1$ is freely infinitely divisible.
\end{theorem}

\begin{proof}
This is something tricky, involving all sorts of not very intuitive computations, and for full details here, we refer here to the original paper \cite{bbl}.
\end{proof}

The above result shows that the normal law $g_1$ should have a ``classical analogue'' in the sense of the Bercovici-Pata bijection. And isn't that puzzling. The problem, however, is that this latter law is difficult to compute, and interpret. See \cite{bbl}.

\bigskip

Still in relation with the Bercovici-Pata bijection, let us also mention that there are many interesting analytic aspects, coming from the combinatorics of the infinitely divisible laws, classical or free. For this, and other analytic aspects, we refer to \cite{bpa}.

\bigskip

Finally, some interesting aspects appear when trying to further build on Theorem 15.32, with noncommutative geometry ideas in mind, with the help of the free spheres constructed in chapter 14, and related objects.  For instance, one interesting thing that happens is that, at the level of the corresponding tori, the classical/free bijection is no longer provided by the Bercovici-Pata correspondence, but rather by the Meixner/free Meixner correspondence. For more on all this, you can have a look at my book \cite{ba3}.

\section*{15e. Exercises}

This was a quite advanced chapter, yes I know, and as exercises here, we have:

\begin{exercise}
Learn more about partitions, and M\"obius inversion.
\end{exercise}

\begin{exercise}
Compute the cumulants of various more specialized measures.
\end{exercise}

\begin{exercise}
Read the proof of the Voiculescu $R$-transform theorem. 
\end{exercise}

\begin{exercise}
Learn more about the Voiculescu circular law $\Gamma_1$.
\end{exercise}

\begin{exercise}
Do the computations, for the low order free cumulants.
\end{exercise}

\begin{exercise}
Then clarify the proof of the Speicher inversion formula.
\end{exercise}

\begin{exercise}
Read more about Bercovici-Pata, notably analytic aspects.
\end{exercise}

\begin{exercise}
Learn about ``easy'' quantum groups, coming from partitions.
\end{exercise}

As bonus exercise, the standard cube is there for everyone to use. Enjoy.

\chapter{Random matrices}

\section*{16a. Wigner matrices}

Time to end this book, and we have kept of course the best for the end. We will be talking here about random matrices, which are the most interesting objects in probability, and perhaps, even in mathematics in general. Or at least, that is the personal opinion of the author of this book, which is myself, and further confirmed by my cats.

\bigskip

We recall that a random matrix algebra is an algebra of type $A=M_N(L^\infty(X))$, and that we are interested in the computation of the laws of the operators $Z\in A$, called random matrices. Regarding the precise classes of random matrices that we are interested in, first we have the complex Gaussian matrices, which are constructed as follows:

\begin{definition}
A complex Gaussian matrix is a random matrix of type
$$Z\in M_N(L^\infty(X))$$
which has i.i.d. complex normal entries.
\end{definition}

We will see that the above matrices have an interesting, and ``central'' combinatorics, among all kinds of random matrices, with the study of the other random matrices being usually obtained as a modification of the study of the Gaussian matrices.

\bigskip

As a somewhat surprising remark, using real normal variables in Definition 16.1, instead of the complex ones appearing there, leads nowhere. The correct real versions of the Gaussian matrices are the Wigner random matrices, constructed as follows: 

\begin{definition}
A Wigner matrix is a random matrix of type
$$Z\in M_N(L^\infty(X))$$
which has i.i.d. complex normal entries, up to the constraint $Z=Z^*$.
\end{definition}

In other words, a Wigner matrix must be as follows, with the diagonal entries being real normal variables, $a_i\sim g_t$, for some $t>0$, the upper diagonal entries being complex normal variables, $b_{ij}\sim G_t$, the lower diagonal entries being the conjugates of the upper diagonal entries, as indicated, and with all the variables $a_i,b_{ij}$ being independent: 
$$Z=\begin{pmatrix}
a_1&b_{12}&\ldots&\ldots&b_{1N}\\
\bar{b}_{12}&a_2&\ddots&&\vdots\\
\vdots&\ddots&\ddots&\ddots&\vdots\\
\vdots&&\ddots&a_{N-1}&b_{N-1,N}\\
\bar{b}_{1N}&\ldots&\ldots&\bar{b}_{N-1,N}&a_N
\end{pmatrix}$$

As a comment here, for many concrete applications the Wigner matrices are in fact the central objects in random matrix theory, and in particular, they are often more important than the Gaussian matrices. In fact, these are the random matrices which were first considered and investigated, a long time ago, by Wigner himself \cite{wig}.

\bigskip

Finally, we will be interested as well in the complex Wishart matrices, which are the positive versions of the above random matrices, constructed as follows: 

\begin{definition}
A complex Wishart matrix is a random matrix of type
$$Z=YY^*\in M_N(L^\infty(X))$$
with $Y$ being a complex Gaussian matrix.
\end{definition}

As before with the Gaussian and Wigner matrices, there are many possible comments that can be made here, of technical or historical nature. First, using in the above real Gaussian variables instead of complex ones leads to a less interesting combinatorics. Also, these matrices were introduced and studied by Marchenko-Pastur in \cite{mpa}, not long after Wigner, and so historically came second. Finally, in what regards their combinatorics and applications, these matrices quite often come first, before both the Gaussian and the Wigner ones, with all this being of course a matter of knowledge and taste.

\bigskip

Summarizing, we have three main types of random matrices, which can be somehow designated as ``complex'', ``real'' and ``positive'', and that we will study in what follows. Let us also mention that there are many other interesting classes of random matrices, usually appearing as modifications of the above. More on these later.

\bigskip

Let us also recall, from chapter 14, what we want to do with these matrices. Given a commutative $C^*$-algebra, $A=C(X)$, the positive linear forms appear as follows, with $\mu$ being positive, and strictly positive if we want $\varphi$ to be faithful and positive:
$$\varphi(f)=\int_Xf(x)d\mu(x)$$

In general, the positive linear forms can be thought of as being integration functionals with respect to some underlying ``positive measures''. Based on this, we can formulate the following very general definition, providing foundations for quantum probability:

\begin{definition}
Let $A$ be a $C^*$-algebra, given with a positive trace $tr:A\to\mathbb C$.
\begin{enumerate}
\item The elements $a\in A$ are called random variables.

\item The moments of such a variable are the numbers $M_k(a)=tr(a^k)$.

\item The law of such a variable is the functional $\mu_a:P\to tr(P(a))$.
\end{enumerate}
\end{definition}

To be more precise, here the exponent $k=\circ\bullet\bullet\circ\ldots$ is by definition a colored integer, and the powers $a^k$ are defined by the following formulae, and multiplicativity: 
$$a^\emptyset=1\quad,\quad
a^\circ=a\quad,\quad
a^\bullet=a^*$$ 

As for the polynomial $P$, this is a noncommuting $*$-polynomial in one variable: 
$$P\in\mathbb C<X,X^*>$$

Observe that the law is uniquely determined by the moments, because we have:
$$P(X)=\sum_k\lambda_kX^k
\implies\mu_a(P)=\sum_k\lambda_kM_k(a)$$

At the level of the general theory, we have the following key result, extending the various results from linear algebra, regarding the self-adjoint and normal matrices:

\index{normal element}
\index{spectral measure}

\begin{theorem}
Let $A$ be a $C^*$-algebra, with a trace $tr$, and consider an element $a\in A$ which is normal, in the sense that $aa^*=a^*a$.
\begin{enumerate}
\item $\mu_a$ is a complex probability measure, satisfying $supp(\mu_a)\subset\sigma(a)$.

\item In the self-adjoint case, $a=a^*$, this measure $\mu_a$ is real.

\item Assuming that $tr$ is faithful, we have $supp(\mu_a)=\sigma(a)$.
\end{enumerate}
\end{theorem}

\begin{proof}
In the normal case, $aa^*=a^*a$, the Gelfand theorem, or rather the subsequent continuous functional calculus theorem, tells us that we have: 
$$<a>=C(\sigma(a))$$

Thus the functional $f(a)\to tr(f(a))$ can be regarded as an integration functional on the algebra $C(\sigma(a))$, and by the Riesz theorem, this latter functional must come from a probability measure $\mu$ on the spectrum $\sigma(a)$, in the sense that we must have:
$$tr(f(a))=\int_{\sigma(a)}f(z)d\mu(z)$$

We are therefore led to the conclusions in the statement, with the uniqueness assertion coming from the fact that the elements $a^k$, taken as usual with respect to colored integer exponents, $k=\circ\bullet\bullet\circ\ldots$\,, generate the whole $C^*$-algebra $C(\sigma(a))$.
\end{proof}

As a first concrete application now of this, by getting back to the random matrices, and to the questions raised in the beginning of this chapter, we have:

\begin{theorem}
Given a random matrix $Z\in M_N(L^\infty(X))$ which is normal, 
$$ZZ^*=Z^*Z$$
its law, which is by definition the following abstract functional,
$$\mu:\mathbb C<X,X^*>\to\mathbb C\quad,\quad 
P\to\frac{1}{N}\int_Xtr(P(Z))$$
when restricted to the usual polynomials in two variables,
$$\mu:\mathbb C[X,X^*]\to\mathbb C\quad,\quad 
P\to\frac{1}{N}\int_Xtr(P(Z))$$
must come from a probability measure on the spectrum $\sigma(Z)\subset\mathbb C$, as follows:
$$\mu(P)=\int_{\sigma(T)}P(x)d\mu(x)$$
We agree to use the symbol $\mu$ for all these notions.
\end{theorem}

\begin{proof}
This follows indeed from what we know from Theorem 16.5, applied to the normal element $a=Z$, belonging to the $C^*$-algebra $A=M_N(L^\infty(X))$. 
\end{proof}

At the level of the basic examples, the situation is as follows:

\begin{theorem}
The basic random matrices $Z\in M_N(L^\infty(X))$ are as follows:
\begin{enumerate}
\item In the case $N=1$ the random matrix is a usual random variable, $f\in L^\infty(X)$, automatically normal, and its law as defined above is the usual law.

\item In the case $X=\{.\}$ the random matrix is a usual scalar matrix, $A\in M_N(\mathbb C)$, and in the normal case, the law is $\mu=\frac{1}{N}\left(\delta_{\lambda_1}+\ldots+\delta_{\lambda_N}\right)$.
\end{enumerate}
\end{theorem}

\begin{proof}
This is something elementary, coming from definitions, as follows:

\medskip

(1) Assume first that we are in the case $N=1$, and so $Z\in L^\infty(X)$. In this case we obtain the usual law of $Z$, because the equation to be satisfied by $\mu$ is:
$$\int_X\varphi(Z)=\int_\mathbb C\varphi(x)d\mu(x)$$

(2) Regarding the second assertion, where we have a matrix $Z\in M_N(\mathbb C)$, which is normal, by changing the basis of $\mathbb C^N$ via a unitary, coming from the spectral theorem, which will not affect our trace computations, we can assume that $Z$ is diagonal:
$$Z\sim
\begin{pmatrix}
\lambda_1\\
&\ddots\\
&&\lambda_N
\end{pmatrix}$$

But for such a diagonal matrix, we have the following formula:
$$tr(Z^k)=\frac{1}{N}(\lambda_1^k+\ldots+\lambda_N^k)$$

More generally, we have the following formula, for any function $\varphi$:
$$tr(\varphi(Z))=\frac{1}{N}(\varphi(\lambda_1)+\ldots+\varphi(\lambda_N))$$

Thus, the law of $Z$ is the average of the Dirac masses at the eigenvalues:
$$\mu=\frac{1}{N}\left(\delta_{\lambda_1}+\ldots+\delta_{\lambda_N}\right)$$

We are therefore led to the conclusions in the statement.
\end{proof}

At a more advanced level now, the main problem regarding the random matrices is that of computing the law of various classes of such matrices, coming in series:

\begin{question}
What is the law of random matrices coming in series
$$Z_N\in M_N(L^\infty(X))$$
in the $N>>0$ regime?
\end{question}

The general strategy here, coming from physicists, is that of computing first the asymptotic law $\mu^0$, in the $N\to\infty$ limit, and then looking for the higher order terms as well, as to finally reach to a series in $N^{-1}$ giving the law of $Z_N$, as follows: 
$$\mu_N=\mu^0+N^{-1}\mu^1+N^{-2}\mu^2+\ldots$$

As a basic example here, of particular interest are the matrices having i.i.d. complex normal entries, under the constraint $Z=Z^*$. Here the asymptotic law $\mu^0$ is the Wigner semicircle law on $[-2,2]$. We will discuss this in a moment, after some preliminaries.

\bigskip

Getting to work now, let us first discuss the Gaussian matrices, from Definition 16.1. As usual, we will use the moment method. We first have the following result:

\begin{theorem}
Given a sequence of Gaussian random matrices
$$Z_N\in M_N(L^\infty(X))$$
having independent $G_t$ variables as entries, for some fixed $t>0$, we have
$$M_k\left(\frac{Z_N}{\sqrt{N}}\right)\simeq t^{|k|/2}|\mathcal{NC}_2(k)|$$
for any colored integer $k=\circ\bullet\bullet\circ\ldots\,$, in the $N\to\infty$ limit.
\end{theorem}

\begin{proof}
This is something standard, which can be done as follows:

\medskip

(1) We fix $N\in\mathbb N$, and we let $Z=Z_N$. Let us first compute the trace of $Z^k$. With $k=k_1\ldots k_s$, and with the convention $(ij)^\circ=ij,(ij)^\bullet=ji$, we have:
\begin{eqnarray*}
Tr(Z^k)
&=&Tr(Z^{k_1}\ldots Z^{k_s})\\
&=&\sum_{i_1=1}^N\ldots\sum_{i_s=1}^N(Z^{k_1})_{i_1i_2}(Z^{k_2})_{i_2i_3}\ldots(Z^{k_s})_{i_si_1}\\
&=&\sum_{i_1=1}^N\ldots\sum_{i_s=1}^N(Z_{(i_1i_2)^{k_1}})^{k_1}(Z_{(i_2i_3)^{k_2}})^{k_2}\ldots(Z_{(i_si_1)^{k_s}})^{k_s}
\end{eqnarray*}

(2) Next, we rescale our variable $Z$ by a $\sqrt{N}$ factor, as in the statement, and we also replace the usual trace by its normalized version, $tr=Tr/N$. Our formula becomes:
$$tr\left(\left(\frac{Z}{\sqrt{N}}\right)^k\right)=\frac{1}{N^{s/2+1}}\sum_{i_1=1}^N\ldots\sum_{i_s=1}^N(Z_{(i_1i_2)^{k_1}})^{k_1}(Z_{(i_2i_3)^{k_2}})^{k_2}\ldots(Z_{(i_si_1)^{k_s}})^{k_s}$$

Thus, the moment that we are interested in is given by:
$$M_k\left(\frac{Z}{\sqrt{N}}\right)=\frac{1}{N^{s/2+1}}\sum_{i_1=1}^N\ldots\sum_{i_s=1}^N\int_X(Z_{(i_1i_2)^{k_1}})^{k_1}(Z_{(i_2i_3)^{k_2}})^{k_2}\ldots(Z_{(i_si_1)^{k_s}})^{k_s}$$

(3) Let us apply now the Wick formula, from chapter 10. We conclude that the moment that we are interested in is given by the following formula:
\begin{eqnarray*}
&&M_k\left(\frac{Z}{\sqrt{N}}\right)\\
&=&\frac{t^{s/2}}{N^{s/2+1}}\sum_{i_1=1}^N\ldots\sum_{i_s=1}^N\#\left\{\pi\in\mathcal P_2(k)\Big|\pi\leq\ker\left((i_1i_2)^{k_1},(i_2i_3)^{k_2},\ldots,(i_si_1)^{k_s}\right)\right\}\\
&=&t^{s/2}\sum_{\pi\in\mathcal P_2(k)}\frac{1}{N^{s/2+1}}\#\left\{i\in\{1,\ldots,N\}^s\Big|\pi\leq\ker\left((i_1i_2)^{k_1},(i_2i_3)^{k_2},\ldots,(i_si_1)^{k_s}\right)\right\}
\end{eqnarray*}

(4) Our claim now is that in the $N\to\infty$ limit the combinatorics of the above sum simplifies, with only the noncrossing partitions contributing to the sum, and with each of them contributing precisely with a 1 factor, so that we will have, as desired:
\begin{eqnarray*}
M_k\left(\frac{Z}{\sqrt{N}}\right)
&=&t^{s/2}\sum_{\pi\in\mathcal P_2(k)}\Big(\delta_{\pi\in NC_2(k)}+O(N^{-1})\Big)\\
&\simeq&t^{s/2}\sum_{\pi\in\mathcal P_2(k)}\delta_{\pi\in NC_2(k)}\\
&=&t^{s/2}|\mathcal{NC}_2(k)|
\end{eqnarray*}

(5) In order to prove this, the first observation is that when $k$ is not uniform, in the sense that it contains a different number of $\circ$, $\bullet$ symbols, we have $\mathcal P_2(k)=\emptyset$, and so:
$$M_k\left(\frac{Z}{\sqrt{N}}\right)=t^{s/2}|\mathcal{NC}_2(k)|=0$$

(6) Thus, we are left with the case where $k$ is uniform. Let us examine first the case where $k$ consists of an alternating sequence of $\circ$ and $\bullet$ symbols, as follows:
$$k=\underbrace{\circ\bullet\circ\bullet\ldots\ldots\circ\bullet}_{2p}$$

In this case it is convenient to relabel our multi-index $i=(i_1,\ldots,i_s)$, with $s=2p$, in the form $(j_1,l_1,j_2,l_2,\ldots,j_p,l_p)$. With this done, our moment formula becomes:
$$M_k\left(\frac{Z}{\sqrt{N}}\right)
=t^p\sum_{\pi\in\mathcal P_2(k)}\frac{1}{N^{p+1}}\#\left\{j,l\in\{1,\ldots,N\}^p\Big|\pi\leq\ker\left(j_1l_1,j_2l_1,j_2l_2,\ldots,j_1l_p\right)\right\}$$

Now observe that, with $k$ being as above, we have an identification $\mathcal P_2(k)\simeq S_p$, obtained in the obvious way. With this done too, our moment formula becomes:
$$M_k\left(\frac{Z}{\sqrt{N}}\right)
=t^p\sum_{\pi\in S_p}\frac{1}{N^{p+1}}\#\left\{j,l\in\{1,\ldots,N\}^p\Big|j_r=j_{\pi(r)+1},l_r=l_{\pi(r)},\forall r\right\}$$

(7) We are now ready to do our asymptotic study, and prove the claim in (4). Let indeed $\gamma\in S_p$ be the full cycle, which is by definition the following permutation:
$$\gamma=(1 \, 2 \, \ldots \, p)$$

In terms of $\gamma$, the conditions $j_r=j_{\pi(r)+1}$ and $l_r=l_{\pi(r)}$ found above read:
$$\gamma\pi\leq\ker j\quad,\quad 
\pi\leq\ker l$$

Counting the number of free parameters in our moment formula, we obtain:
$$M_k\left(\frac{Z}{\sqrt{N}}\right)
=\frac{t^p}{N^{p+1}}\sum_{\pi\in S_p}N^{|\pi|+|\gamma\pi|}
=t^p\sum_{\pi\in S_p}N^{|\pi|+|\gamma\pi|-p-1}$$

(8) The point now is that the last exponent is well-known to be $\leq 0$, with equality precisely when the permutation $\pi\in S_p$ is geodesic, which in practice means that $\pi$ must come from a noncrossing partition. Thus we obtain, in the $N\to\infty$ limit, as desired:
$$M_k\left(\frac{Z}{\sqrt{N}}\right)\simeq t^p|\mathcal{NC}_2(k)|$$

This finishes the proof in the case of the exponents $k$ which are alternating, and the case where $k$ is an arbitrary uniform exponent is similar, by permuting everything.
\end{proof}

We can further improve Theorem 16.9 by using free probability, as follows:

\begin{theorem}
Given a sequence of complex Gaussian matrices
$$Z_N\in M_N(L^\infty(X))$$
having independent $G_t$ variables as entries, with $t>0$, we have
$$\frac{Z_N}{\sqrt{N}}\sim\Gamma_t$$
in the $N\to\infty$ limit, with the limiting measure being Voiculescu's circular law.
\end{theorem}

\begin{proof}
According to Theorem 16.9 that we just proved, the asymptotic moments of the complex Gaussian matrices are given by the following formula:
$$M_k\left(\frac{Z_N}{\sqrt{N}}\right)\simeq t^{|k|/2}|\mathcal{NC}_2(k)|$$

On the other hand, it follows from the free probability theory developed in chapter 15 that an abstract noncommutative variable $a\in A$ is circular, following the Voiculescu circular law $\Gamma_t$, precisely when its moments are given by the following formula:
$$M_k(a)=t^{|k|/2}|\mathcal{NC}_2(k)|$$

Thus, we are led to the conclusion in the statement.
\end{proof}

The above result is of course something quite theoretical, and having it formulated as such is certainly something nice. However, and here comes our point, it is actually possible to use free probability theory in order to go well beyond this, with this time some truly new results on the random matrices. More on this, later in this chapter.

\bigskip
 
Regarding now the Wigner matrices, which are something more understandable, because they are self-adjoint, we first have the following result, about them:

\begin{theorem}
Given a sequence of Wigner random matrices
$$Z_N\in M_N(L^\infty(X))$$
having independent $G_t$ variables as entries, with $t>0$, up to $Z_N=Z_N^*$, we have
$$M_k\left(\frac{Z_N}{\sqrt{N}}\right)\simeq t^{k/2}|NC_2(k)|$$
for any integer $k\in\mathbb N$, in the $N\to\infty$ limit.
\end{theorem}

\begin{proof}
We have two possible proofs here, as follows:

\medskip

(1) The formula in the statement can be certainly established via a direct computation based on the Wick formula, similar to that from the proof of Theorem 16.9.

\medskip

(2) However, the best is to deduce this result from Theorem 16.9 itself. Indeed, we know from there that for Gaussian matrices $Y_N\in M_N(L^\infty(X))$ we have the following formula, valid for any colored integer $K=\circ\bullet\bullet\circ\ldots\,$, in the $N\to\infty$ limit:
$$M_K\left(\frac{Y_N}{\sqrt{N}}\right)\simeq t^{|K|/2}|\mathcal{NC}_2(K)|$$

By doing now some combinatorics, we deduce that we have the following formula for the moments of the matrices $Re(Y_N)$, with respect to usual exponents, $k\in\mathbb N$:
\begin{eqnarray*}
M_k\left(\frac{Re(Y_N)}{\sqrt{N}}\right)
&=&2^{-k}\cdot M_k\left(\frac{Y_N}{\sqrt{N}}+\frac{Y_N^*}{\sqrt{N}}\right)\\
&=&2^{-k}\sum_{|K|=k}M_K\left(\frac{Y_N}{\sqrt{N}}\right)\\
&\simeq&2^{-k}\sum_{|K|=k}t^{k/2}|\mathcal{NC}_2(K)|\\
&=&2^{-k}\cdot t^{k/2}\cdot 2^{k/2}|\mathcal{NC}_2(k)|\\
&=&2^{-k/2}\cdot t^{k/2}|NC_2(k)|
\end{eqnarray*}

Now since the matrices $Z_N=\sqrt{2}Re(Y_N)$ are of Wigner type, this gives the result.
\end{proof}

In order to recapture now the Wigner measure from its moments, we can use:

\begin{proposition}
Given $t>0$, the real measure having as even moments the numbers $M_{2k}=t^kC_k$ and having all odd moments $0$ is the measure
$$\gamma_t=\frac{1}{2\pi t}\sqrt{4t-x^2}dx$$
called Wigner semicircle law on $[-2\sqrt{t},2\sqrt{t}]$.
\end{proposition}

\begin{proof}
This is something that we know well since chapter 7, when first discussing the Wigner law, coming either via Stieltjes inversion, or from a direct integration.
\end{proof}

Now by putting everything together, we obtain the Wigner theorem, as follows:

\begin{theorem}
Given a sequence of Wigner random matrices
$$Z_N\in M_N(L^\infty(X))$$
which by definition have i.i.d. complex normal entries, up to $Z_N=Z_N^*$, we have
$$Z_N\sim\gamma_t$$
in the $N\to\infty$ limit, where $\gamma_t=\frac{1}{2\pi t}\sqrt{4t-x^2}dx$ is the Wigner semicircle law. 
\end{theorem}

\begin{proof}
This follows indeed from all the above, and more specifically, by combining the moment formula from Theorem 16.11 with Proposition 16.12.
\end{proof}

\section*{16b. Wishart matrices}

Let us discuss now the Wishart matrices, which are the positive analogues of the Wigner matrices. Quite surprisingly, the computation here leads to the Catalan numbers, but not in the same way as for the Wigner matrices, the result being as follows:

\begin{theorem}
Given a sequence of complex Wishart matrices
$$W_N=Y_NY_N^*\in M_N(L^\infty(X))$$
with $Y_N$ being $N\times N$ complex Gaussian of parameter $t>0$, we have
$$M_k\left(\frac{W_N}{N}\right)\simeq t^kC_k$$
for any exponent $k\in\mathbb N$, in the $N\to\infty$ limit.
\end{theorem}

\begin{proof}
There are several possible proofs for this result, as follows:

\medskip

(1) A first method is by using the formula that we have in Theorem 16.9, for the Gaussian matrices $Y_N$. Indeed, we know from there that we have the following formula, valid for any colored integer $K=\circ\bullet\bullet\circ\ldots\,$, in the $N\to\infty$ limit:
$$M_K\left(\frac{Y_N}{\sqrt{N}}\right)\simeq t^{|K|/2}|\mathcal{NC}_2(K)|$$

With $K=\circ\bullet\circ\bullet\ldots\,$, alternating word of length $2k$, with $k\in\mathbb N$, this gives:
$$M_k\left(\frac{Y_NY_N^*}{N}\right)\simeq t^k|\mathcal{NC}_2(K)|$$

Thus, in terms of the Wishart matrix $W_N=Y_NY_N^*$ we have, for any $k\in\mathbb N$:
$$M_k\left(\frac{W_N}{N}\right)\simeq t^k|\mathcal{NC}_2(K)|$$

The point now is that, by doing some combinatorics, we have:
$$|\mathcal{NC}_2(K)|=|NC_2(2k)|=C_k$$

Thus, we are led to the formula in the statement.

\medskip

(2) A second method, that we will explain now as well, is by proving the result directly, starting from definitions. The matrix entries of our matrix $W=W_N$ are given by:
$$W_{ij}=\sum_{r=1}^NY_{ir}\bar{Y}_{jr}$$

Thus, the normalized traces of powers of $W$ are given by the following formula:
\begin{eqnarray*}
tr(W^k)
&=&\frac{1}{N}\sum_{i_1=1}^N\ldots\sum_{i_k=1}^NW_{i_1i_2}W_{i_2i_3}\ldots W_{i_ki_1}\\
&=&\frac{1}{N}\sum_{i_1=1}^N\ldots\sum_{i_k=1}^N\sum_{r_1=1}^N\ldots\sum_{r_k=1}^NY_{i_1r_1}\bar{Y}_{i_2r_1}Y_{i_2r_2}\bar{Y}_{i_3r_2}\ldots Y_{i_kr_k}\bar{Y}_{i_1r_k}
\end{eqnarray*}

By rescaling now $W$ by a $1/N$ factor, as in the statement, we obtain:
$$tr\left(\left(\frac{W}{N}\right)^k\right)=\frac{1}{N^{k+1}}\sum_{i_1=1}^N\ldots\sum_{i_k=1}^N\sum_{r_1=1}^N\ldots\sum_{r_k=1}^NY_{i_1r_1}\bar{Y}_{i_2r_1}Y_{i_2r_2}\bar{Y}_{i_3r_2}\ldots Y_{i_kr_k}\bar{Y}_{i_1r_k}$$

By using now the Wick rule, we obtain the following formula for the moments, with $K=\circ\bullet\circ\bullet\ldots\,$, alternating word of lenght $2k$, and with $I=(i_1r_1,i_2r_1,\ldots,i_kr_k,i_1r_k)$:
\begin{eqnarray*}
M_k\left(\frac{W}{N}\right)
&=&\frac{t^k}{N^{k+1}}\sum_{i_1=1}^N\ldots\sum_{i_k=1}^N\sum_{r_1=1}^N\ldots\sum_{r_k=1}^N\#\left\{\pi\in\mathcal P_2(K)\Big|\pi\leq\ker(I)\right\}\\
&=&\frac{t^k}{N^{k+1}}\sum_{\pi\in\mathcal P_2(K)}\#\left\{i,r\in\{1,\ldots,N\}^k\Big|\pi\leq\ker(I)\right\}
\end{eqnarray*}

In order to compute this quantity, we use the standard bijection $\mathcal P_2(K)\simeq S_k$. By identifying the pairings $\pi\in\mathcal P_2(K)$ with their counterparts $\pi\in S_k$, we obtain:
\begin{eqnarray*}
M_k\left(\frac{W}{N}\right)
&=&\frac{t^k}{N^{k+1}}\sum_{\pi\in S_k}\#\left\{i,r\in\{1,\ldots,N\}^k\Big|i_s=i_{\pi(s)+1},r_s=r_{\pi(s)},\forall s\right\}
\end{eqnarray*}

Now let $\gamma\in S_k$ be the full cycle, which is by definition the following permutation:
$$\gamma=(1 \, 2 \, \ldots \, k)$$

The general factor in the product computed above is then 1 precisely when following two conditions are simultaneously satisfied:
$$\gamma\pi\leq\ker i\quad,\quad 
\pi\leq\ker r$$

Counting the number of free parameters in our moment formula, we obtain:
$$M_k\left(\frac{W}{N}\right)
=t^k\sum_{\pi\in S_k}N^{|\pi|+|\gamma\pi|-k-1}$$

The point now is that the last exponent is well-known to be $\leq 0$, with equality precisely when the permutation $\pi\in S_k$ is geodesic, which in practice means that $\pi$ must come from a noncrossing partition. Thus we obtain, in the $N\to\infty$ limit:
$$M_k\left(\frac{W}{N}\right)\simeq t^kC_k$$

Thus, we are led to the conclusion in the statement.
\end{proof}

In order to recapture now the corresponding measure, we can use our knowledge of the Marchenko-Pastur law $\pi_1$, from chapter 7, and we obtain in this way:

\begin{theorem}
Given a sequence of complex Wishart matrices
$$W_N=Y_NY_N^*\in M_N(L^\infty(X))$$
with $Y_N$ being $N\times N$ complex Gaussian of parameter $t>0$, we have
$$\frac{W_N}{tN}\sim\frac{1}{2\pi}\sqrt{4x^{-1}-1}\,dx$$
with $N\to\infty$, with the limiting measure being the Marchenko-Pastur law $\pi_1$.
\end{theorem}

\begin{proof}
This follows indeed from Theorem 16.14, via our knowledge of $\pi_1$.
\end{proof}

As a comment now, while the above result is definitely something interesting at $t=1$, at general $t>0$ this looks more like a ``fake'' generalization of the $t=1$ result, because the law $\pi_1$ stays the same, modulo a trivial rescaling. The reasons behind this phenomenon are quite subtle, and skipping some discussion, the point is that Theorem 16.15 is indeed something ``fake'' at general $t>0$, and the correct generalization of the $t=1$ computation, involving more general classes of complex Wishart matrices, is as follows:

\begin{theorem}
Given a sequence of general complex Wishart matrices
$$W_N=Y_NY_N^*\in M_N(L^\infty(X))$$
with $Y_N$ being $N\times M$ complex Gaussian of parameter $1$, we have
$$\frac{W_N}{N}\sim\max(1-t,0)\delta_0+\frac{\sqrt{4t-(x-1-t)^2}}{2\pi x}\,dx$$
with $M=tN\to\infty$, with the limiting measure being the Marchenko-Pastur law $\pi_t$.
\end{theorem}

\begin{proof}
This is again something which is very standard, as follows:

\medskip

(1) In order to prove the formula in the statement, we can proceed as usual, by using the Wick formula. The matrix entries of our Wishart matrix $W=W_N$ are given by:
$$W_{ij}=\sum_{r=1}^MY_{ir}\bar{Y}_{jr}$$

Thus, the normalized traces of powers of $W$ are given by the following formula:
\begin{eqnarray*}
tr(W^k)
&=&\frac{1}{N}\sum_{i_1=1}^N\ldots\sum_{i_k=1}^NW_{i_1i_2}W_{i_2i_3}\ldots W_{i_ki_1}\\
&=&\frac{1}{N}\sum_{i_1=1}^N\ldots\sum_{i_k=1}^N\sum_{r_1=1}^M\ldots\sum_{r_k=1}^MY_{i_1r_1}\bar{Y}_{i_2r_1}Y_{i_2r_2}\bar{Y}_{i_3r_2}\ldots Y_{i_kr_k}\bar{Y}_{i_1r_k}
\end{eqnarray*}

By rescaling now $W$ by a $1/N$ factor, as in the statement, we obtain:
$$tr\left(\left(\frac{W}{N}\right)^k\right)=\frac{1}{N^{k+1}}\sum_{i_1=1}^N\ldots\sum_{i_k=1}^N\sum_{r_1=1}^M\ldots\sum_{r_k=1}^MY_{i_1r_1}\bar{Y}_{i_2r_1}Y_{i_2r_2}\bar{Y}_{i_3r_2}\ldots Y_{i_kr_k}\bar{Y}_{i_1r_k}$$

(2) By using now the Wick rule, we obtain the following formula for the moments, with $K=\circ\bullet\circ\bullet\ldots\,$, alternating word of lenght $2k$, and $I=(i_1r_1,i_2r_1,\ldots,i_kr_k,i_1r_k)$:
\begin{eqnarray*}
M_k\left(\frac{W}{N}\right)
&=&\frac{1}{N^{k+1}}\sum_{i_1=1}^N\ldots\sum_{i_k=1}^N\sum_{r_1=1}^M\ldots\sum_{r_k=1}^M\#\left\{\pi\in\mathcal P_2(K)\Big|\pi\leq\ker I\right\}\\
&=&\frac{1}{N^{k+1}}\sum_{\pi\in\mathcal P_2(K)}\#\left\{i\in\{1,\ldots,N\}^k,r\in\{1,\ldots,M\}^k\Big|\pi\leq\ker I\right\}
\end{eqnarray*}

(3) In order to compute this quantity, we use the standard bijection $\mathcal P_2(K)\simeq S_k$. By identifying the pairings $\pi\in\mathcal P_2(K)$ with their counterparts $\pi\in S_k$, we obtain:
\begin{eqnarray*}
M_k\left(\frac{W}{N}\right)
&=&\frac{1}{N^{k+1}}\sum_{\pi\in S_k}\#\left\{i\in\{1,\ldots,N\}^k,r\in\{1,\ldots,M\}^k\Big|i_s=i_{\pi(s)+1},r_s=r_{\pi(s)}\right\}
\end{eqnarray*}

Now let $\gamma\in S_k$ be the full cycle, which is by definition the following permutation:
$$\gamma=(1 \, 2 \, \ldots \, k)$$

The general factor in the product computed above is then 1 precisely when following two conditions are simultaneously satisfied:
$$\gamma\pi\leq\ker i\quad,\quad 
\pi\leq\ker r$$

Counting the number of free parameters in our expectation formula, we obtain:
$$M_k\left(\frac{W}{N}\right)
=\frac{1}{N^{k+1}}\sum_{\pi\in S_k}N^{|\gamma\pi|}M^{|\pi|}
=\sum_{\pi\in S_k}N^{|\gamma\pi|-k-1}M^{|\pi|}$$

(4) Now by using the same arguments as in the case $M=N$, from the proof of Theorem 16.14, we conclude that in the $M=tN\to\infty$ limit the permutations $\pi\in S_k$ which matter are those coming from noncrossing partitions, and so that we have:
$$M_k\left(\frac{W}{N}\right)
\simeq\sum_{\pi\in NC(k)}N^{-|\pi|}M^{|\pi|}
=\sum_{\pi\in NC(k)}t^{|\pi|}$$

(5) But these numbers are the moments of the Marchenko-Pastur law $\pi_t$, which in addition has the density given by the formula in the statement, as desired.
\end{proof}

And with this, end of our basic learning regarding the random matrices, the conclusion being that the asymptotic laws of these matrices are the main laws in free probability.

\section*{16c. Shifted semicircles}

Our goal now will be that of explaining a surprising result, due to Aubrun \cite{aub}, stating that when suitably block-transposing the entries of a complex Wishart matrix, we obtain as asymptotic distribution a shifted version of Wigner's semicircle law. 

\bigskip

Following \cite{aub}, and some more recent work in \cite{bne}, where the Aubrun finding was further studied, and better understood, let us start with the following definition:

\begin{definition}
The partial transpose of a complex Wishart matrix $W$ of parameters $(dn,dm)$ is the matrix
$$\widetilde{W}=(id\otimes t)W$$
where $id$ is the identity of $M_d(\mathbb C)$, and $t$ is the transposition of $M_n(\mathbb C)$. 
\end{definition}

In more familiar terms of bases and indices, the standard decomposition $\mathbb C^{dn}=\mathbb C^d\otimes\mathbb C^n$ induces an algebra decomposition $M_{dn}(\mathbb C)=M_d(\mathbb C)\otimes M_n(\mathbb C)$, and with this convention made, the partial transpose matrix $\widetilde{W}$ constructed above has entries as follows:
$$\widetilde{W}_{ia,jb}=W_{ib,ja}$$

Our goal in what follows will be that of computing the law of $\widetilde{W}$, first when $d,n,m$ are fixed, and then in the $d\to\infty$ regime. For this purpose, we will need a number of standard facts regarding the noncrossing partitions. Let us start with:

\begin{proposition}
For a permutation $\sigma\in S_p$, we have the formula
$$|\sigma|+\#\sigma=p$$
where $|\sigma|$ is the number of cycles of $\sigma$, and $\#\sigma$ is the minimal $k\in\mathbb N$ such that $\sigma$ is a product of $k$ transpositions. Also, the following formula defines a distance on $S_p$, 
$$(\sigma,\pi)\to\#(\sigma^{-1}\pi)$$
and the set of permutations $\sigma\in S_p$ which saturate the triangular inequality 
$$\#\sigma+\#(\sigma^{-1}\gamma)=\#\gamma=p-1$$
where $\gamma\in S_p$ is a full cycle, is in bijection with the set $NC(p)$. 
\end{proposition}

\begin{proof}
All this is standard combinatorics, that we will leave as an exercise.
\end{proof}

We use the standard bijection $NC(p)\simeq NC_2(2p)$, denoted $\pi\to\widetilde{\pi}$, obtained by fattening the partitions. We have the following formula,  where $\vee$ is the join operation on $NC_2(2p)$, and $\rho_{12} = (12)(34)\ldots(2p-1,2p)$ is the fattened identity permutation:
$$|\pi|=|\widetilde{\pi}\vee\rho_{12}|$$

Similarly, we have the formula $|\pi\gamma|=|\widetilde{\pi}\vee\rho_{14}|$, where $\rho_{14}$ is the pairing corresponding to the fattening of the inverse full cycle $\gamma^{-1}(i) = i-1$, which pairs an element $2i$ with $2(i-1)-1 = 2i-3$, or, equivalently, an element $i \in \{1, \ldots, 2p\}$ with $i+(-1)^{i+1}3$.

\bigskip

We will need the following well-known result, which is something standard:

\begin{proposition}
The number $||\pi||$ of blocks having even size is given by
$$1+||\pi||=|\pi\gamma|$$
for every noncrossing partition $\pi \in NC(p)$.
\end{proposition}

\begin{proof}
We use a recurrence over the number of blocks of $\pi$. If $\pi$ has just one block, its associated geodesic permutation is $\gamma$ and we have:
$$|\gamma^2|=\begin{cases}
1&(p\ \text{odd})\\ 
2&(p\ \text{even})\\
\end{cases}$$

For the partitions $\pi$ having more than one block, we can assume without loss of generality that $\pi = \hat 1_k \sqcup \pi'$, where $\hat 1_k$ is a contiguous block of size $k$. Recall that the number of blocks of the permutation $\pi\gamma$ is given by the following formula, where $\rho_{14} \in P_2(2p)$ is the pair partition which pairs an element $i$ with $i+(-1)^{i+1}3$:
$$|\pi\gamma|=|\widetilde{\pi}\vee\rho_{14}|$$

If $k$ is an even number, $k=2r$, consider the following partition, which contains the block $(1 \, 4 \, 5 \, 8 \, \ldots 4r-3 \, 4r)$, along with the blocks coming from the elements of the form $4i+2, 4i+3$ from $\{1, \ldots, 4r\}$ and from $\pi'$:
$$\sigma=\widetilde{\hat 1_{2r} \sqcup \pi'}\vee \rho_{14}$$

We can count the blocks of the join of two partitions by drawing them one beneath the other and counting the number of connected components of the curve, without taking into account the possible crossings. We conclude that we have the following formula, where $\rho'_{14}$ is $\rho_{14}$ restricted to the set $\{2k+1, 2k+2 \ldots, 2p\}$:
$$|\widetilde{\pi}\vee\rho_{14}|=1+|\widetilde{\pi'}\vee\rho'_{14}|$$

If $k$ is odd, $k=2r+1$, there is no extra block appearing, so we have:
$$|\widetilde{\pi}\vee\rho_{14}|=|\widetilde{\pi'}\vee\rho'_{14}|$$

Thus, we are led to the conclusion in the statement.
\end{proof}

We can now investigate the block-transposed Wishart matrices, and we have:

\begin{theorem}
For any $p\geq 1$ we have the formula
$$\lim_{d\to\infty}(E\circ tr)\big(m\widetilde{W}\big)^p
=\sum_{\pi\in NC(p)}m^{|\pi|}n^{||\pi||}$$
where $|.|$ and $||.||$ are the number of blocks, and the number of blocks of even size. 
\end{theorem}

\begin{proof}
The matrix elements of the partial transpose matrix are given by:
$$\widetilde{W}_{ia,jb}=W_{ib,ja}=(dm)^{-1}\sum_{k=1}^d\sum_{c=1}^mG_{ib,kc}\bar{G}_{ja,kc}$$

This gives the following formula:
\begin{eqnarray*}
tr(\widetilde{W}^p)
&=&(dn)^{-1}(dm)^{-p}\sum_{i_1,\ldots,i_p=1}^d\sum_{a_1,\ldots,a_p=1}^n\prod_{s=1}^p W_{i_sa_{s+1},i_{s+1}a_s} \\
&=&(dn)^{-1}(dm)^{-p}\sum_{i_1,\ldots,i_p=1}^d\sum_{a_1,\ldots,a_p=1}^n\prod_{s=1}^p \sum_{j_1,\ldots,j_p=1}^d\sum_{b_1,\ldots,b_p=1}^mG_{i_sa_{s+1},j_sb_s}\bar{G}_{i_{s+1}a_s,j_sb_s}
\end{eqnarray*}

After interchanging the product with the last two sums, the average of the general term can be computed by the Wick rule, namely:
$$E\left(\prod_{s=1}^pG_{i_sa_{s+1},j_sb_s}\bar{G}_{i_{s+1}a_s,j_sb_s}\right)
=\sum_{\pi\in S_p}\prod_{s=1}^p\delta_{i_s,i_{\pi(s)+1}}\delta_{a_{s+1},a_{\pi(s)}}\delta_{j_s,j_{\pi(s)}}\delta_{b_s,b_{\pi(s)}}$$

Let $\gamma\in S_p$ be the full cycle $\gamma=(1 \, 2 \, \ldots \, p)^{-1}$. The general factor in the above product is 1 if and only if the following four conditions are simultaneously satisfied:
$$\gamma^{-1}\pi\leq \ker i\quad,\quad
\pi\gamma \leq \ker a\quad,\quad
\pi \leq \ker j\quad,\quad
\pi \leq \ker b$$

Counting the number of free parameters in the above equation, we obtain:
\begin{eqnarray*}
(E\circ tr)(\widetilde{W}^p)
&=&(dn)^{-1}(dm)^{-p}\sum_{\pi\in S_p}d^{|\pi|+|\gamma^{-1}\pi|}m^{|\pi|}n^{|\pi\gamma|}\\
&=&\sum_{\pi\in S_p}d^{|\pi|+|\gamma^{-1}\pi|-p-1}m^{|\pi|-p}n^{|\pi\gamma|-1}
\end{eqnarray*}

The exponent of $d$ in the last expression on the right is:
\begin{eqnarray*}
N(\pi)
&=&|\pi|+|\gamma^{-1}\pi|-p-1\\
&=&p-1-(\#\pi+\#(\gamma^{-1}\pi))\\
&=&p-1-(\#\pi+\#(\pi^{-1}\gamma))
\end{eqnarray*}

As explained in the beginning of this section, this quantity is known to be $\leq 0$, with equality iff $\pi$ is geodesic, hence associated to a noncrossing partition. Thus:
$$(E\circ tr)(\widetilde{W}^p)=(1+O(d^{-1}))m^{-p}n^{-1}\sum_{\pi\in NC(p)}m^{|\pi|} n^{|\pi\gamma|}$$

Together with $|\pi\gamma|=||\pi||+1$, this gives the result.
\end{proof}

We would like now to find an equation for the moment generating function of the asymptotic law of $m\widetilde{W}$. This moment generating function is defined by:
$$F(z)=\lim_{d\to\infty}(E\circ tr)\left(\frac{1}{1-zm\widetilde{W}}\right)$$

We have the following result, regarding this moment generating function:

\begin{theorem}
The moment generating function of $m\widetilde{W}$ satisfies the equation
$$(F-1)(1-z^2F^2)=mzF(1+nzF)$$
in the $d\to\infty$ limit.
\end{theorem}

\begin{proof}
We use the formula in Theorem 16.20. If we denote by $N(p,b,e)$ the number of partitions in $NC(p)$ having $b$ blocks and $e$ even blocks, we have:
\begin{eqnarray*}
F
&=&1+\sum_{p=1}^\infty\sum_{\pi\in NC(p)} z^pm^{|\pi|}n^{||\pi||}\\
&=&1+\sum_{p=1}^\infty\sum_{b=0}^\infty\sum_{e=0}^\infty z^pm^bn^eN(p,b,e)
\end{eqnarray*}

Let us try to find a recurrence formula for the numbers $N(p,b,e)$. If we look at the block containing $1$, this block must have $r\geq 0$ other legs, and we get:
\begin{eqnarray*}
N(p,b,e) 
&=&\sum_{r\in 2\mathbb N}\sum_{p=\Sigma p_i+r+1}\sum_{b=\Sigma b_i+1}\sum_{e=\Sigma e_i}N(p_1,b_1,e_1)\ldots N(p_{r+1},b_{r+1},e_{r+1})\\
&+&\sum_{r\in 2\mathbb N+1}\sum_{p=\Sigma p_i+r+1}\sum_{b=\Sigma b_i+1}\sum_{e=\Sigma e_i+1}N(p_1,b_1,e_1)\ldots N(p_{r+1},b_{r+1},e_{r+1})
\end{eqnarray*}

Here $p_1,\ldots,p_{r+1}$ are the number of points between the legs of the block containing 1, so that we have $p=(p_1+\ldots+p_{r+1})+r+1$, and the whole sum is split over two cases, $r$ even or odd, because the parity of $r$ affects the number of even blocks of our partition. Now by multiplying everything by a $z^pm^bn^e$ factor, and by carefully distributing the various powers of $z,m,b$ on the right, we obtain the following formula:
\begin{eqnarray*}
z^pm^bn^eN(p,b,e)
&=&m\sum_{r\in 2\mathbb N}z^{r+1}\sum_{p=\Sigma p_i+r+1}\sum_{b=\Sigma b_i+1}\sum_{e=\Sigma e_i}\prod_{i=1}^{r+1}z^{p_i}m^{b_i}n^{e_i}N(p_i,b_i,e_i)\\
&+&mn\sum_{r\in 2\mathbb N+1}z^{r+1}\sum_{p=\Sigma p_i+r+1}\sum_{b=\Sigma b_i+1}\sum_{e=\Sigma e_i+1}\prod_{i=1}^{r+1}z^{p_i}m^{b_i}n^{e_i}N(p_i,b_i,e_i)
\end{eqnarray*}

Let us sum now all these equalities, over all $p\geq 1$ and over all $b,e\geq 0$. According to the definition of $F$, at left we obtain $F-1$. As for the two sums appearing on the right, that is, at right of the two $z^{r+1}$ factors, when summing them over all $p\geq 1$ and over all $b,e\geq 0$, we obtain in both cases $F^{r+1}$. So, we have the following formula:
\begin{eqnarray*}
F-1
&=&m\sum_{r\in 2\mathbb N}(zF)^{r+1}+mn\sum_{r\in 2\mathbb N+1}(zF)^{r+1}\\
&=&m\,\frac{zF}{1-z^2F^2}+mn\,\frac{z^2F^2}{1-z^2F^2}\\
&=&mzF\,\frac{1+nzF}{1-z^2F^2}
\end{eqnarray*}

But this gives the formula in the statement, and we are done.
\end{proof}

Our goal now will be that of further processing the formula in Theorem 16.21, as to reach to a formula for the density of the corresponding law. This is something quite tricky, and as a first result here, we can reformulate Theorem 16.21 as follows:

\begin{theorem}
The Cauchy transform of $m\widetilde{W}$ satisfies the equation
$$(\xi G-1)(1-G^2)=mG(1+nG)$$
in the $d\to\infty$ limit. Moreover, this equation simply reads
$$R=\frac{m}{2}\left(\frac{n+1}{1-z}-\frac{n-1}{1+z}\right)$$
with the substitutions $G\to z$ and $\xi\to R+z^{-1}$.
\end{theorem}

\begin{proof}
We have two assertions to be proved, the first one being standard, and the second one being something quite magic, the idea being as follows:

\medskip

(1) Consider the equation of $F$, found in Theorem 16.21, namely:
$$(F-1)(1-z^2F^2)=mzF(1+nzF)$$

With $z\to\xi^{-1}$ and $F\to\xi G$, so that $zF\to G$, we obtain, as desired:
$$(\xi G-1)(1-G^2)=mG(1+nG)$$

(2) Thus, we have our equation for the Cauchy transform, and with this in hand, we can try to go ahead, and use somehow the Stieltjes inversion formula, in order to reach to a formula for the density. This is certainly possible, but our claim is that we can do better, by performing first some clever manipulations on the Cauchy transform.

\medskip

(3) To be more precise, let us look at the equation of the Cauchy transform that we have. With the substitutions $\xi \to K$ and $G\to z$, this equation becomes:
$$(zK-1)(1-z^2)=mz(1+nz)$$

The point now is that with $K\to R+z^{-1}$ this latter equation becomes:
$$zR(1-z^2)=mz(1+nz)$$

But the solution of this latter equation is trivial to compute, given by:
$$R=m\,\frac{1+nz}{1-z^2}=\frac{m}{2}\left(\frac{n+1}{1-z}-\frac{n-1}{1+z}\right)$$

Thus, we are led to the conclusion in the statement.
\end{proof}

Now let us recall from chapter 15 that we have the following key definition:

\begin{definition}
Given a real probability measure $\mu$, its $R$-transform is:
$$G_\mu(\xi)=\int_\mathbb R\frac{d\mu(t)}{\xi-t}\implies 
G_\mu\left(R_\mu(\xi)+\frac{1}{\xi}\right)=\xi$$
That is, the $R$-transform is the inverse of the Cauchy transform, up to a $\xi^{-1}$ factor.
\end{definition}

Getting back now to our questions, we would like to find the probability measure having as $R$-transform the function in Theorem 16.22. But here, we can only expect to find some kind of modification of the Marchenko-Pastur law, so as a first piece of work, let us just compute the $R$-transform of the Marchenko-Pastur law. We have here:

\begin{proposition}
The $R$-transform of the Marchenko-Pastur law $\pi_t$ is
$$R_{\pi_t}(\xi)=\frac{t}{1-\xi}$$ 
for any $t>0$.
\end{proposition}

\begin{proof}
This can be done in two steps, as follows:

\medskip

(1) At $t=1$, we know that the moments of $\pi_1$ are the Catalan numbers, $M_k=C_k$, and we obtain that the Cauchy transform is given by the following formula:
$$G(\xi)=\frac{1}{2}-\frac{1}{2}\sqrt{1-4\xi^{-1}}$$

Now with $R(\xi)=\frac{1}{1-\xi}$ being the function in the statement, at $t=1$, we have:
\begin{eqnarray*}
G\left(R(\xi)+\frac{1}{\xi}\right)
&=&G\left(\frac{1}{1-\xi}+\frac{1}{\xi}\right)\\
&=&G\left(\frac{1}{\xi-\xi^2}\right)\\
&=&\frac{1}{2}-\frac{1}{2}\sqrt{1-4\xi+4\xi^2}\\
&=&\frac{1}{2}-\frac{1}{2}(1-2\xi)\\
&=&\xi
\end{eqnarray*}

Thus, the function $R(\xi)=\frac{1}{1-\xi}$ is indeed the $R$-transform of $\pi_1$, in the above sense.

\medskip

(2) In the general case, $t>0$, the proof is similar, by using the moment formula for $\pi_t$, that we know from the above, from the proof of the Marchenko-Pastur theorem.
\end{proof}

All this is very nice, and we can now further build on Theorem 16.22, as follows:

\begin{theorem}
The $R$-transform of $m\widetilde{W}$ is given by
$$R=R_{\pi_s}-R_{\pi_t}$$
in the $d\to\infty$ limit, where $s=m(n+1)/2$ and $t=m(n-1)/2$.
\end{theorem}

\begin{proof}
We know from Theorem 16.22 that the $R$-transform of $m\widetilde{W}$ is given by:
$$R=\frac{m}{2}\left(\frac{n+1}{1-z}-\frac{n-1}{1+z}\right)$$

By using now the formula in Proposition 16.24, this gives the result.
\end{proof}

We can now recover the original result of Aubrun \cite{aub}, as follows:

\begin{theorem}
For a block-transposed Wishart matrix $\widetilde{W}=(id\otimes t)W$ we have, in the $n=\beta m\to\infty$ limit, with $\beta>0$ fixed, the formula
$$\frac{\widetilde{W}}{d}\sim\gamma_{\beta}^1$$
with $\gamma_\beta^1$ being the shifted version of the semicircle law $\gamma_\beta$, with support centered at $1$.
\end{theorem}

\begin{proof}
This follows from Theorem 16.25. Indeed, in the $n=\beta m\to\infty$ limit, with $\beta>0$ fixed, we are led to the following formula for the Stieltjes transform:
$$f(x)=\frac{\sqrt{4\beta-(1-x)^2}}{2\beta\pi}$$

But this is the density of the shifted semicircle law having support as follows:
$$S=[1-2\sqrt{\beta},1+2\sqrt{\beta}]$$

Thus, we are led to the conclusion in the statement. See \cite{aub}, \cite{bne}.
\end{proof}

As a continuation of this, let us discuss  now more general block modifications. Consider a complex Wishart matrix of parameters $(dn,dm)$. In other words, we start with a $dn\times dm$ matrix $Y$ having independent complex $G_1$ entries, and we set:
$$W=YY^*$$

This matrix has size $dn\times dn$, and is best thought of as being a $d\times d$ array of $n\times n$ matrices. We will be interested here in the study of the block-modified versions of $W$, obtained by applying to the $n\times n$ blocks a given linear map, as follows:
$$\varphi:M_n(\mathbb C)\to M_n(\mathbb C)$$

By reasoning as in the proof of Theorem 16.20, we have the following asymptotic moment formula, extending the usual moment computation for the Wishart matrices:

\index{block-modified matrix}

\begin{theorem}
The asymptotic moments of a block-modified Wishart matrix 
$$\widetilde{W}=(id\otimes\varphi)W$$
with parameters $d,m,n\in\mathbb N$, as above, are given by the formula
$$\lim_{d\to\infty}M_e\left(\frac{\widetilde{W}}{d}\right)=\sum_{\sigma\in NC_p}(mn)^{|\sigma|}(M^\sigma_e\otimes M^\gamma_e)(\Lambda)$$
where $\Lambda\in M_n(\mathbb C)\otimes M_n(\mathbb C)$ is the square matrix associated to $\varphi:M_n(\mathbb C)\to M_n(\mathbb C)$.
\end{theorem}

\begin{proof}
This is something quite standard, coming as usual from the Wick formula, via some computations similar to those done before in this chapter.
\end{proof}

In order to interpret the moment formula above, we can use free probability. Indeed, by using the free cumulant theory from chapter 15, we have the following result:

\begin{proposition}
Given a square matrix $\Lambda\in M_n(\mathbb C)\otimes M_n(\mathbb C)$, having distribution 
$$\rho=law(\Lambda)$$
the moments of the compound free Poisson law $\pi_{mn\rho}$ are given by
$$M_e(\pi_{mn\rho})=\sum_{\sigma\in NC_p}(mn)^{|\sigma|}(M^\sigma_e\otimes M^\sigma_e)(\Lambda)$$
for any choice of the extra parameter $m\in\mathbb N$.
\end{proposition}

\begin{proof}
We can use the fact, that we know from chapter 15, that the free cumulants of $\pi_{mn\rho}$ coincide with the moments of $mn\rho$. Thus, these free cumulants are:
\begin{eqnarray*}
\kappa_e(\pi_{mn\rho})
&=&M_e(mn\rho)\\
&=&mn\cdot M_e(\Lambda)\\
&=&mn\cdot (M^\gamma_e\otimes M^\gamma_e)(\Lambda)
\end{eqnarray*}

But with this cumulant formula in hand, by using now Speicher's free moment-cumulant formula, explained in chapter 15 as well, this gives the result.
\end{proof}

We can see now an obvious similarity with the formula in Theorem 16.27. In order to exploit this similarity, once again by following \cite{bne}, let us introduce:

\index{multiplicative matrix}

\begin{definition}
We call a square matrix $\Lambda\in M_n(\mathbb C)\otimes M_n(\mathbb C)$ multiplicative when
$$(M^\sigma_e\otimes M^\gamma_e)(\Lambda)=(M^\sigma_e\otimes M^\sigma_e)(\Lambda)$$
holds for any $p\in\mathbb N$, any exponents $e_1,\ldots,e_p\in\{1,*\}$, and any $\sigma\in NC_p$.
\end{definition}

Now with the above notion in hand, which is something quite tricky, we can formulate an asymptotic result regarding the block-modified Wishart matrices, as follows:

\begin{theorem}
Consider a block-modified Wishart matrix 
$$\widetilde{W}=(id\otimes\varphi)W$$
and assume that the matrix $\Lambda\in M_n(\mathbb C)\otimes M_n(\mathbb C)$ associated to $\varphi$ is multiplicative. Then
$$\frac{\widetilde{W}}{d}\sim\pi_{mn\rho}$$
holds, in moments, in the $d\to\infty$ limit, where $\rho=law(\Lambda)$.
\end{theorem}

\begin{proof}
By comparing the moment formulae in Theorem 16.27 and in Proposition 16.28, we conclude that the asymptotic formula $\frac{\widetilde{W}}{d}\sim\pi_{mn\rho}$ is equivalent to the following equality, which should hold for any $p\in\mathbb N$, and any exponents $e_1,\ldots,e_p\in\{1,*\}$:
$$\sum_{\sigma\in NC_p}(mn)^{|\sigma|}(M^\sigma_e\otimes M^\gamma_e)(\Lambda)=\sum_{\sigma\in NC_p}(mn)^{|\sigma|}(M^\sigma_e\otimes M^\sigma_e)(\Lambda)$$

Now by assuming that the matrix $\Lambda$ is multiplicative, in the sense of Definition 16.29, these two sums are trivially equal, and this gives the result.
\end{proof}

\section*{16d. Asymptotic freeness}

We would like to end the present chapter, and book, with some key asymptotic freeness results of Voiculescu \cite{vdn}. Let us begin with the Wigner matrices. We have here:

\index{Wigner matrix}
\index{asymptotic freeness}

\begin{theorem}
Given a family of sequences of Wigner matrices, 
$$Z^i_N\in M_N(L^\infty(X))\quad,\quad i\in I$$
with pairwise independent entries, each following the complex normal law $G_t$, with $t>0$, up to the constraint $Z_N^i=(Z_N^i)^*$, the rescaled sequences of matrices
$$\frac{Z^i_N}{\sqrt{N}}\in M_N(L^\infty(X))\quad,\quad i\in I$$
become with $N\to\infty$ semicircular, each following the Wigner law $\gamma_t$, and free.
\end{theorem}

\begin{proof}
This is something quite subtle, the idea being as follows:

\medskip

(1) First of all, we know from Theorem 16.13 that for any $i\in I$ the corresponding sequence of rescaled Wigner matrices becomes semicircular in the $N\to\infty$ limit:
$$\frac{Z_N^i}{\sqrt{N}}\simeq\gamma_t$$ 

(2) Thus, what is new here, and that we have to prove, is the asymptotic freeness assertion. For this purpose we can assume that we are dealing with the case of 2 sequences of matrices, $|I|=2$. So, assume that we have Wigner matrices as follows:
$$Z_N,Z_N'\in M_N(L^\infty(X))$$

We have to prove that these matrices become asymptotically free, with $N\to\infty$.

\medskip

(3) The best is to prove this by using a trick, based on the result in Theorem 16.10. Consider indeed the following random matrix:
$$Y_N=\frac{1}{\sqrt{2}}(Z_N+iZ_N')$$

This is then a complex Gaussian matrix, and so by using Theorem 16.10, we obtain that in the limit $N\to\infty$, we have:
$$\frac{Y_N}{\sqrt{N}}\simeq\Gamma_t$$

(4) Now recall that the circular law $\Gamma_t$ was by definition the law of the following variable, with $a,b$ being semicircular, each following the law $\gamma_t$, and free:
$$c=\frac{1}{\sqrt{2}}(a+ib)$$

We are therefore in the situation where the variable $(Z_N+iZ_N')/\sqrt{N}$, which has asymptotically semicircular real and imaginary parts, converges to the distribution of $a+ib$, equally having semicircular real and imaginary parts, but with these real and imaginary parts being free. Thus $Z_N,Z_N'$ become asymptotically free, as desired.
\end{proof}

Getting now to the complex case, we have a similar result here, as follows:

\index{Gaussian matrix}
\index{asymptotic freeness}

\begin{theorem}
Given a family of sequences of complex Gaussian matrices, 
$$Z^i_N\in M_N(L^\infty(X))\quad,\quad i\in I$$
with pairwise independent entries, each following the complex normal law $G_t$, with $t>0$, the rescaled sequences of matrices
$$\frac{Z^i_N}{\sqrt{N}}\in M_N(L^\infty(X))\quad,\quad i\in I$$
become with $N\to\infty$ circular, each following the Voiculescu law $\Gamma_t$, and free.
\end{theorem}

\begin{proof}
This follows from Theorem 16.31, which applies to the real and imaginary parts of our complex Gaussian matrices, and gives the result.
\end{proof}

We have as well a result for the Wishart matrices, as follows:

\begin{theorem}
Given a family of sequences of complex Wishart matrices, 
$$Z^i_N=Y^i_N(Y^i_N)^*\in M_N(L^\infty(X))\quad,\quad i\in I$$
with each $Y^i_N$ being a $N\times M$ matrix, with entries following the normal law $G_1$, and with all these entries being pairwise independent, the rescaled sequences of matrices
$$\frac{Z^i_N}{N}\in M_N(L^\infty(X))\quad,\quad i\in I$$
become with $M=tN\to\infty$ Marchenko-Pastur, each following the law $\pi_t$, and free.
\end{theorem}

\begin{proof}
Here the first assertion is the Marchenko-Pastur theorem, and the second assertion follows from the freeness result for the Gaussian matrices.
\end{proof}

\section*{16e. Exercises}

Congratulations for having read this book, and no exercises for this final chapter. Instead, you can have a look at the various books and articles referenced below.

\baselineskip=14pt

\printindex


\begin{thebibliography}{99}

\baselineskip=14.3pt

\bibitem{afa}G.W. Anderson and B. Farrell, Asymptotically liberating sequences of random unitary matrices, {\em Adv. Math.} {\bf 255} (2014), 381--413.

\bibitem{agz}G.W. Anderson, A. Guionnet and O. Zeitouni, An introduction to random matrices, Cambridge Univ. Press (2010).

\bibitem{ans}M. Anshelevich, Free Meixner states, {\em Comm. Math. Phys.} {\bf 276} (2007), 863--899.

\bibitem{apo}T.M. Apostol, Introduction to analytic number theory, Springer (1976).

\bibitem{anv}O. Arizmendi, I. Nechita and C. Vargas, On the asymptotic distribution of block-modified random matrices, {\em J. Math. Phys.} {\bf 57} (2016), 1--27.

\bibitem{arn}V.I. Arnold, Mathematical methods of classical mechanics, Springer (1974).

\bibitem{ame}N.W. Ashcroft and N.D. Mermin, Solid state physics, Saunders College Publ. (1976).

\bibitem{ati}M.F. Atiyah, The geometry and physics of knots, Cambridge Univ. Press (1990).

\bibitem{aub}G. Aubrun, Partial transposition of random states and non-centered semicircular distributions, {\em Random Matrices Theory Appl.} {\bf 1} (2012), 125--145.

\bibitem{bds}J. Baik, P. Deift and T. Suidan, Combinatorics and random matrix theory, AMS (2016).

\bibitem{ba1}T. Banica, Discrete random variables (2026).

\bibitem{ba2}T. Banica, Methods of free probability (2024).

\bibitem{ba3}T. Banica, Introduction to modern physics (2025).

\bibitem{bb+}T. Banica, S.T. Belinschi, M. Capitaine and B. Collins, Free Bessel laws, {\em Canad. J. Math.} {\bf 63} (2011), 3--37.

\bibitem{bbi}T. Banica and D. Bisch, Spectral measures of small index principal graphs, {\em Comm. Math. Phys.} {\bf 269} (2007), 259--281.

\bibitem{bne}T. Banica and I. Nechita, Asymptotic eigenvalue distributions of block-transposed Wishart matrices, {\em J. Theoret. Probab.} {\bf 26} (2013), 855--869.

\bibitem{bsp}T. Banica and R. Speicher, Liberation of orthogonal Lie groups, {\em Adv. Math.} {\bf 222} (2009), 1461--1501.

\bibitem{bax}R.J. Baxter, Exactly solved models in statistical mechanics, Academic Press (1982).

\bibitem{bbl}S.T. Belinschi, M. Bo\.zejko, F. Lehner and R. Speicher, The normal distribution is $\boxplus$-infinitely divisible, {\em Adv. Math.} {\bf 226} (2011), 3677--3698.

\bibitem{bzy}I. Bengtsson and K. \.Zyczkowski, Geometry of quantum states, Cambridge Univ. Press (2006).

\bibitem{bpa}H. Bercovici and V. Pata, Stable laws and domains of attraction in free probability theory, {\em Ann. of Math.} {\bf 149} (1999), 1023--1060.

\bibitem{bgv}N. Berline, E. Getzler and M. Vergne, Heat kernels and Dirac operators, Springer (2004).

\bibitem{bcg}P. Biane, M. Capitaine and A. Guionnet, Large deviation bounds for matrix Brownian motion, {\em Invent. Math.} {\bf 152} (2003), 433--459.

\bibitem{blo}G. Blower, Random matrices: high dimensional phenomena, Cambridge Univ. Press (2009).

\bibitem{blu}S.J. Blundell and K.M. Blundell, Concepts in thermal physics, Oxford Univ. Press (2006).

\bibitem{bol}B. Bollob\'as, Random graphs, Cambridge Univ. Press (1985).

\bibitem{bos}A. Bose, Random matrices and non-commutative probability, CRC Press (2021).

\bibitem{csn}B. Collins and P. \'Sniady, Integration with respect to the Haar measure on unitary, orthogonal and symplectic groups, {\em Comm. Math. Phys.} {\bf 264} (2006), 773--795.

\bibitem{con}A. Connes, Noncommutative geometry, Academic Press (1994).

\bibitem{cur}S. Curran, Quantum rotatability, {\em Trans. Amer. Math. Soc.} {\bf 362} (2010), 4831--4851.

\bibitem{dei}P. Deift, Orthogonal polynomials and random matrices: a Riemann-Hilbert approach, AMS (1999).

\bibitem{dif}P. Di Francesco, Meander determinants, {\em Comm. Math. Phys.} {\bf 191} (1998), 543--583.

\bibitem{dms}P. Di Francesco, P. Mathieu and D. S\'en\'echal, Conformal field theory, Springer (1996).

\bibitem{dir}P.A.M. Dirac, Principles of quantum mechanics, Oxford Univ. Press (1930).

\bibitem{dre}M. Dresher, The mathematics of games of strategy, Dover (1981).

\bibitem{du1}R. Durrett, Probability: theory and examples, Cambridge Univ. Press (1990).

\bibitem{du2}R. Durrett, Random graph dynamics, Cambridge Univ. Press (2010).

\bibitem{epu}D.E. Evans and M. Pugh, Spectral measures and generating series for nimrep graphs in subfactor theory, {\em Comm. Math. Phys.} {\bf 295} (2010), 363--413.

\bibitem{eyn}B. Eynard, Counting surfaces, Birkh\"auser (2016).

\bibitem{fel}W. Feller, An introduction to probability theory and its applications, Wiley (1950).

\bibitem{fer}E. Fermi, Thermodynamics, Dover (1937). 

\bibitem{fls}R.P. Feynman, R.B. Leighton and M. Sands, The Feynman lectures on physics, Caltech (1963).

\bibitem{fse}P. Flajolet and R. Sedgewick, Analytic combinatorics, Cambridge Univ. Press (2009).

\bibitem{fsn}M. Fukuda and P. \'Sniady, Partial transpose of random quantum states: exact formulas and meanders, {\em J. Math. Phys.} {\bf 54} (2013), 1--31.

\bibitem{goo}D.L. Goodstein, States of matter, Dover (1975).

\bibitem{glm}P. Graczyk, G. Letac and H. Massam, The complex Wishart distribution and the symmetric group, {\em Ann. Statist.} {\bf 31} (2003), 287--309.

\bibitem{gsw}M.B. Green, J.H. Schwarz and E. Witten, Superstring theory, Cambridge Univ. Press (2012). 

\bibitem{gri}D.J. Griffiths and D.F. Schroeter, Introduction to quantum mechanics, Cambridge Univ. Press (2018).

\bibitem{gkz}A. Guionnet, M. Krishnapur and O. Zeitouni, The single ring theorem, {\em Ann. of Math.} {\bf 174} (2011), 1189--1217.

\bibitem{gur}R. Gurau, Random tensors, Oxford Univ. Press (2016).

\bibitem{haa}U. Haagerup, On Voiculescu's R and S transforms for free non-commuting random variables, {\em Fields Inst. Comm.} {\bf 12} (1997), 127--148.

\bibitem{hth}U. Haagerup and S. Thorbj\o rnsen, Random matrices with complex Gaussian entries, {\em Exposition. Math.} {\bf 21} (2003), 293--337.

\bibitem{hpe}F. Hiai and D. Petz, The semicircle law, free random variables and entropy, AMS (2000).

\bibitem{hjo}R.A. Horn and C.R. Johnson, Matrix analysis, Cambridge Univ. Press (1985).

\bibitem{hua}K. Huang, Introduction to statistical physics, CRC Press (2001).

\bibitem{joh}K. Johansson, Shape fluctuations and random matrices, {\em Comm.  Math. Phys.} {\bf 209} (2000), 437--476.

\bibitem{jo1}V.F.R. Jones, Index for subfactors, {\em Invent. Math.} {\bf 72} (1983), 1--25.

\bibitem{jo2}V.F.R. Jones, On knot invariants related to some statistical mechanical models, {\em Pacific J. Math.} {\bf 137} (1989), 311--334.

\bibitem{jo3}V.F.R. Jones, Planar algebras I (1999).

\bibitem{jo4}V.F.R. Jones, The annular structure of subfactors, {\em Monogr. Enseign. Math.} {\bf 38} (2001), 401--463.

\bibitem{jun}K. Jung, Amenability, tubularity, and embeddings into $R^\omega$, {\em Math. Ann.} {\bf 338} (2007), 241--248.

\bibitem{kad}L.P. Kadanoff, Statistical physics: statics, dynamics and renormalization, World Scientific (2000).

\bibitem{ksp}C. K\"ostler, R. Speicher, A noncommutative de Finetti theorem: invariance under quantum permutations is equivalent to freeness with amalgamation, {\em Comm. Math. Phys.} {\bf 291} (2009), 473--490.

\bibitem{kum}M. Kumar, Quantum: Einstein, Bohr, and the great debate about the nature of reality, Norton (2009). 

\bibitem{lbl}T. Lancaster and K.M. Blundell, Quantum field theory for the gifted amateur, Oxford Univ. Press (2014).

\bibitem{lax}P. Lax, Functional analysis, Wiley (2002).

\bibitem{lig}T.M. Liggett, Interacting particle systems, Springer (1985).

\bibitem{liu}W. Liu, General de Finetti type theorems in noncommutative probability, {\em  Comm. Math. Phys.} {\bf 369} (2019), 837--866.

\bibitem{mpa}V.A. Marchenko and L.A. Pastur, Distribution of eigenvalues in certain sets of random matrices, {\em Mat. Sb.} {\bf 72} (1967), 507--536.

\bibitem{meh}M.L. Mehta, Random matrices, Elsevier (2004).

\bibitem{msp}J.A. Mingo and R. Speicher, Free probability and random matrices, Springer (2017).

\bibitem{nsp}A. Nica and R. Speicher, Lectures on the combinatorics of free probability, Cambridge Univ. Press (2006).

\bibitem{nch}M.A. Nielsen and I.L. Chuang, Quantum computation and quantum information, Cambridge Univ. Press (2000).

\bibitem{pbe}R.K. Pathria and and P.D. Beale, Statistical mechanics, Elsevier (1972).

\bibitem{pbo}M. Potters and J.P. Bouchaud, A first course in random matrix theory, Cambridge Univ. Press (2020).

\bibitem{rud}W. Rudin, Real and complex analysis, McGraw-Hill (1966).

\bibitem{ksc}K. Schm\"udgen, The moment problem, Springer (2017).

\bibitem{sch}D.V. Schroeder, An introduction to thermal physics, Oxford Univ. Press (1999).

\bibitem{sp1}R. Speicher, Multiplicative functions on the lattice of noncrossing partitions and free convolution, {\em Math. Ann.} {\bf 298} (1994), 611--628.

\bibitem{sp2}R. Speicher, Combinatorial theory of the free product with amalgamation and operator-valued free probability theory, {\em Mem. Amer. Math. Soc.} {\bf 132} (1998).

\bibitem{ste}A.M. Steane, Thermodynamics, Oxford Univ. Press (2016).

\bibitem{twe}P. Tarrago and M. Weber, Unitary easy quantum groups: the free case and the group case, {\em Int. Math. Res. Not.} {\bf 18} (2017), 5710--5750.

\bibitem{tli}N.H. Temperley and E.H. Lieb, Relations between the ``percolation'' and ``colouring'' problem and other graph-theoretical problems associated with regular planar lattices: some exact results for the ``percolation'' problem, {\em Proc. Roy. Soc. London} {\bf 322} (1971), 251--280.

\bibitem{twi}C.A. Tracy and H. Widom, Level-spacing distributions and the Airy kernel, {\em Comm. Math. Phys.} {\bf 159} (1994), 151--174.

\bibitem{vo1}D.V. Voiculescu, Addition of certain noncommuting random variables, {\em J. Funct. Anal.} {\bf 66} (1986), 323--346.

\bibitem{vo2}D.V. Voiculescu, Multiplication of certain noncommuting random variables, {\em J. Operator Theory} {\bf 18} (1987), 223--235.

\bibitem{vo3}D.V. Voiculescu, Limit laws for random matrices and free products, {\em Invent. Math.} {\bf 104} (1991), 201--220.

\bibitem{vdn}D.V. Voiculescu, K.J. Dykema and A. Nica, Free random variables, AMS (1992).

\bibitem{von}J. von Neumann, Mathematical foundations of quantum mechanics, Princeton Univ. Press (1955).

\bibitem{vmo}J. von Neumann and O. Morgenstern, Theory of games and economic behavior, Princeton Univ. Press (1944).

\bibitem{wan}S. Wang, Quantum symmetry groups of finite spaces, {\em Comm. Math. Phys.} {\bf 195} (1998), 195--211.

\bibitem{wat}J. Watrous, The theory of quantum information, Cambridge Univ. Press (2018).

\bibitem{swe}S. Weinberg, Lectures on quantum mechanics, Cambridge Univ. Press (2012).

\bibitem{wei}D. Weingarten, Asymptotic behavior of group integrals in the limit of infinite rank, {\em J. Math. Phys.} {\bf 19} (1978), 999--1001.

\bibitem{we1}H. Weyl, The theory of groups and quantum mechanics, Princeton Univ. Press (1931).

\bibitem{we2}H. Weyl, The classical groups: their invariants and representations, Princeton Univ. Press (1939).

\bibitem{wig}E. Wigner, Characteristic vectors of bordered matrices with infinite dimensions, {\em Ann. of Math.} {\bf 62} (1955), 548--564. 

\bibitem{wo1}S.L. Woronowicz, Compact matrix pseudogroups, {\em Comm. Math. Phys.} {\bf 111} (1987), 613--665.

\bibitem{wo2}S.L. Woronowicz, Tannaka-Krein duality for compact matrix pseudogroups. Twisted SU(N) groups, {\em Invent. Math.} {\bf 93} (1988), 35--76.

\bibitem{zin}P. Zinn-Justin, Jucys-Murphy elements and Weingarten matrices, {\em Lett. Math. Phys.} {\bf 91} (2010), 119--127.

\end{thebibliography}
\end{document}